\documentclass[11pt]{article}
\usepackage[margin=1in]{geometry}
\usepackage[T1]{fontenc}
\usepackage{lmodern,microtype,amsmath,amssymb,amsthm,mathtools,mathrsfs}
\usepackage{booktabs,longtable,enumitem,graphicx,tikz,float,xurl,hyperref}
\usetikzlibrary{arrows.meta,positioning,calc,decorations.pathreplacing,decorations.markings,patterns}
\hypersetup{colorlinks=true,linkcolor=blue!50!black,citecolor=blue!50!black,urlcolor=blue!50!black,
 pdftitle={Quantum unique ergodicity and the critical crossover for Gaussian block band matrices}}
\allowdisplaybreaks
\numberwithin{equation}{section}
\newtheorem{theorem}{Theorem}[section]
\newtheorem{proposition}[theorem]{Proposition}
\newtheorem{lemma}[theorem]{Lemma}
\newtheorem{corollary}[theorem]{Corollary}
\theoremstyle{definition}

\theoremstyle{remark}
\newtheorem{remark}[theorem]{Remark}
\newcommand{\Naux}{{N_0}}
\newcommand{\R}{\mathbb R}
\newcommand{\C}{\mathbb C}
\newcommand{\T}{\mathbb T}
\newcommand{\Z}{\mathbb Z}
\newcommand{\E}{\mathbb E}
\newcommand{\Prob}{\mathbb P}
\newcommand{\one}{\mathbf 1}

\DeclareMathOperator{\Tr}{Tr}
\DeclareMathOperator{\Var}{Var}

\DeclareMathOperator{\diag}{diag}
\DeclareMathOperator{\rank}{rank}
\DeclareMathOperator{\GUE}{GUE}
\DeclareMathOperator{\HS}{HS}

\DeclareMathOperator{\Law}{Law}
\newcommand{\norm}[1]{\lVert#1\rVert}

\newcommand{\ip}[2]{\langle#1,#2\rangle}
\newcommand{\av}[1]{\langle#1\rangle}
\title{Quantum unique ergodicity and the critical crossover\\for 1-d Gaussian band matrices}
\author{Promit Ghosal \& Guillaume Remy}

\hypersetup{hidelinks}
\begin{document}
\maketitle
{\renewcommand{\thefootnote}{}%
 \footnotetext{P.~Ghosal: Department of Statistics, University of Chicago, USA.
 Email: \href{mailto:EMAIL1}{\nolinkurl{promit@uchicago.edu}}}%
 \footnotetext{G.~Remy: Axiom Math, Palo Alto, USA.
 }%
}
\begin{abstract}
We consider an $N\times N$ complex 1-d Gaussian random band matrix  %
with band width of order  $W$. The matrix has the periodic block form used
by Yau and Yin \cite{YY}. The  $N$ sites are split into $n=N/W$ blocks of size $W$
arranged on a ring, and the entries are nonzero only inside a block and
between adjacent blocks. The conjectured localization--delocalization
transition occurs at $W\sim\sqrt N$, that is, when the number of blocks
$n$ is comparable to the block size $W$. We prove quantum unique ergodicity
for bulk eigenvectors, at every deterministic rank, whenever $W\gg\sqrt N$,
that is, $W^2/N\to\infty$. This requires no power or logarithmic separation
between $W^2$ and $N$. At the threshold, $N/W^2\to\tau\in(0,\infty)$, we
show that the rescaled bulk eigenvalues and the block-mass profiles of their eigenvectors (probability measures  on the circle that place at the point $k/n$ the
squared norm of the eigenvector restricted to the $k$-th block) converge jointly to a unique law. The limit law is encoded by a random entire function, the limit of a
perturbed characteristic determinant of the random band matrix, whose zeros give the limits
of the rescaled eigenvalues and whose derivatives there give the limits of
the profile integrals (i.e., the integrals of test functions against their
block-mass profiles) of the eigenvectors. The law of this function is
determined by the Dorokhov--Mello--Pereyra--Kumar (DMPK) diffusion \cite{BBDR} of the
unitary class. In this critical law, the
limits of the block mass profiles have strictly positive variance against every
nonconstant smooth test function on the circle. Thus the marked law is not
concentrated on the uniform profile. As $\tau\to0$ this
variance tends to zero and the local spectrum tends to the sine process.
As $\tau\to\infty$ the local spectrum tends to a Poisson process. Both
endpoint limits are taken after the critical limit. Thus the critical
family interpolates between the delocalized regime $W\gg\sqrt N$ and the
localized regime $W\ll\sqrt N$. This confirms, for this model, the prediction of Fyodorov and
Mirlin~\cite{FM91} that at $W^2\asymp N$ the joint law of eigenvalues and
eigenvectors depends only on $N/W^2$ and interpolates between the
delocalized and localized phases. Lastly the full content of this work has been autoformalized in Lean4 by AxiomProver.
\end{abstract}
\setcounter{tocdepth}{2}

\section{Introduction}\label{sec:intro}

 Anderson~\cite{And58} discovered, in his study of electron transport in disordered materials, that disorder can stop a quantum particle from spreading. Since then, the transition between localized and delocalized states has been a central problem of mathematical physics. In Anderson's model the Hamiltonian is
the random Schr\"odinger operator $-\Delta+\lambda V$ on $\Z^d$. Here
$\Delta$ is the graph Laplacian, the values of the potential $V$ are iid with a
nondegenerate common law, and $\lambda>0$ measures the strength of the disorder. %
An
eigenvector is \emph{localized} if its mass is concentrated in a region
much smaller than the system, and \emph{delocalized} if it spreads over
the whole system. In
dimension one, every eigenvector is localized, however weak the
disorder~\cite{CKM87}. In higher dimensions, localization is known when
the disorder is strong~\cite{FS83,AM93}, see~\cite{AWbook} for an
overview. The delocalization expected at weak disorder remains a major
open problem, and the transition between the two regimes has been a
central question of mathematical physics ever since.

Random band matrices are a prominent alternative model for this
transition~\cite{CMI90,FM91,Spe11,Bou18}. A random band matrix
$H=(H_{xy})$ is a random Hermitian matrix whose entries are independent up
to the symmetry $H_{yx}=\overline{H_{xy}}$ and centered. Its variance
profile $S_{xy}=\E|H_{xy}|^2$ is of order $W^{-1}$ when $x$ and $y$ are
at distance at most $W$, and is negligible otherwise. It is normalized so
that $\sum_yS_{xy}=1$. The parameter $W$ is called the \emph{band width}.
When $W$ is of order one, $H$ is a short-range operator, like a random
Schr\"odinger operator. When $W$ is comparable to the dimension $N$ of
the matrix, $H$ is a mean-field matrix of Wigner type~\cite{Wig55}. As $W$
grows, band matrices therefore interpolate between local and mean-field
models, and the tools of random matrix theory become available for a
system with a spatial structure. Heuristically, the band width $W$ plays
the role of the inverse disorder strength $\lambda^{-1}$ of the Anderson model. Here and below, \emph{one dimensional} random band matrices means that the $N$ sites are arranged
along a line or a ring, and that the entry $H_{xy}$ can be nonzero only
when $x$ and $y$ are within distance $W$ along it.

In one dimension, numerical work~\cite{CMI90} and the nonrigorous
supersymmetric $\sigma$-model calculation of Fyodorov and
Mirlin~\cite{FM91} predict that bulk eigenvectors have localization length
of order $W^2$. Since this length cannot exceed $N$, one arrives at the
following conjecture for energies in the bulk of the spectrum.
\begin{enumerate}[label=(\roman*)]
\item If $W^2\gg N$, the eigenvectors are delocalized and the local
eigenvalue statistics are those of the Gaussian unitary ensemble (GUE). In its strongest form, delocalization is characterized by  \emph{quantum unique ergodicity}
(QUE) which is achieved when every eigenvector spreads evenly over the system, in the sense that
the mass it puts on any fixed region is close to the fraction of the
system that this region occupies. 
\item If $W^2\ll N$, the eigenvectors are localized and the local
eigenvalue statistics follow a Poisson point process. %
\end{enumerate}
The same physical analysis~\cite{FM91} also describes the critical window
$W^2\asymp N$. There the statistics of eigenvectors and eigenvalues are
expected to depend on the system only through the ratio $N/W^2$, and to
interpolate between (i) and (ii) as this ratio moves from $0$ to
$\infty$. At the spectral edge in one dimension, the corresponding transition occurs at
$W\asymp N^{5/6}$ as this was established by Sodin~\cite{Sod10} with the
moment method. See this beautiful survey \cite{Bou18} on further discussion about history of this problem.

 Away from the critical window, the transition in dimension one 
is now
well understood. On the delocalized side, a sequence of works established
delocalization, quantum diffusion, quantum unique ergodicity and bulk
universality under conditions of the form $W\ge N^a$ for various exponents
$a>1/2$~\cite{EK11,BE17,BYY20,YangYin,DY24}. Through a dynamical analysis
of the loop hierarchy, Yau and Yin~\cite{YY} proved delocalization,
quantum unique ergodicity and GUE statistics under
$W\ge N^{1/2+\varepsilon}$ for every fixed $\varepsilon>0$. Their
model is a Gaussian block model, and it coincides with the model of this
paper (at one value of the parameter $b$); see \eqref{eq:model} below.
Erd\H{o}s and Riabov~\cite{ER} extended this result to general variance
profiles and entry distributions, in both the real symmetric and the
complex Hermitian symmetry classes. On the localized side, localization
was proved under conditions $W\ll N^a$ for various exponents
$a<1/2$~\cite{Sch09,PSSS,CS22,CPSS24}. More recently, for a class of real
block band matrices, Drogin~\cite{Drogin} proved that 
correlations between two sites at distance $r$ of an eigenvector decay like
$e^{-cr/W^2}$, up to a factor polynomial in $W$. This gives localization
on the scale $W^2$ whenever $W^2\log W\ll N$. Fan and Zou~\cite{FZ}
proved Poisson statistics for a real Gaussian block model when
$W^2\log N/N\to0$.

These results leave two questions open. The first concerns the approach
to the critical window from the delocalized side. The results
of~\cite{YY,ER} require a fixed power separation $W^2\ge N^{1+\varepsilon}$.
The supersymmetric results of M.~and T.~Shcherbina~\cite{SS18,SS19}, which
give the GUE limit of the two-point correlation function, require at least
a logarithmic separation, for instance $W\ge Cn\log^5n$ in~\cite{SS19}.
None of them covers, for example, $W^2=N\log\log W$. The second question
concerns the critical window itself, where much less is known.
Shcherbina~\cite{Sh19} identified a critical interpolation for the second
correlation function of characteristic polynomials at $N\asymp W^2$. These
observables, however, do not determine the joint law of all eigenvalues,
and they do not provide any direct information about the spatial distribution of
eigenvectors. A detailed
comparison with all of these works is given in
Subsection~\ref{subsec:comparison}.

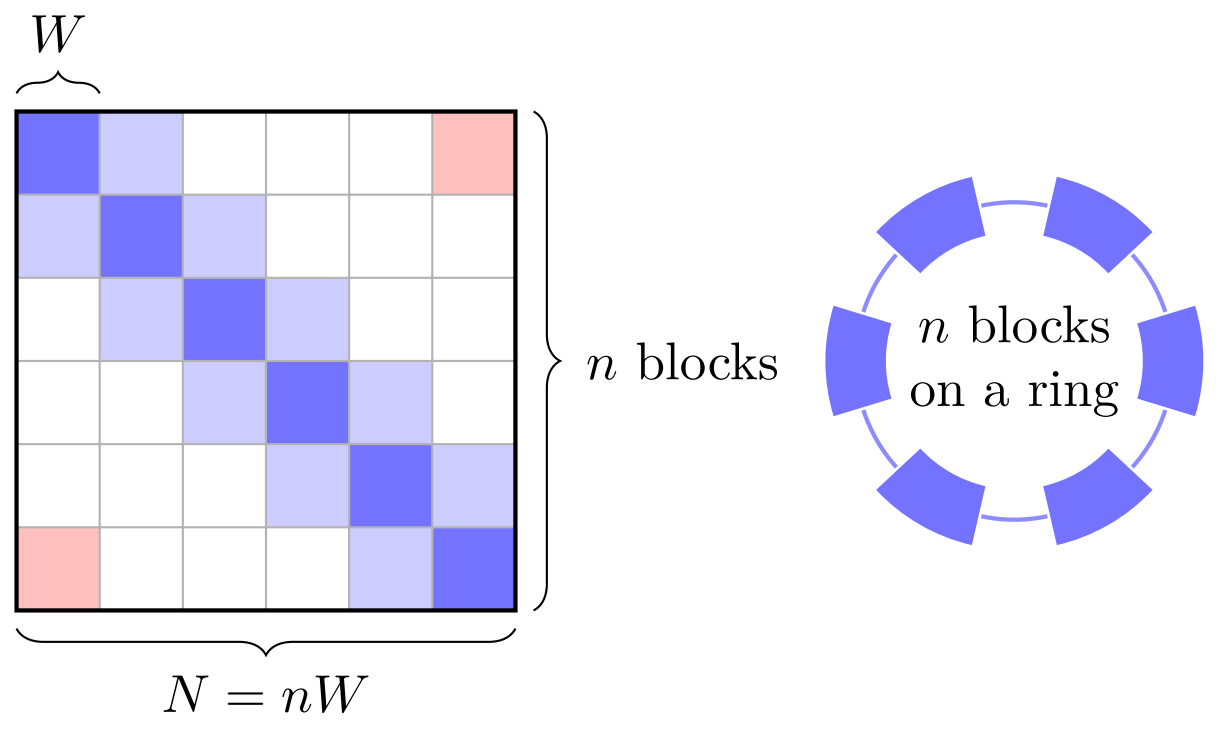
\begin{figure}[ht]
\centering
\begin{tikzpicture}[x=1cm,y=1cm]
 \def\b{0.55}
 \foreach \i in {0,...,5}{
   \foreach \j in {0,...,5}{
     \pgfmathtruncatemacro{\d}{abs(\i-\j)}
     \ifnum\d=0
       \fill[blue!55] (\j*\b,-\i*\b) rectangle ++(\b,-\b);
     \else\ifnum\d=1
       \fill[blue!20] (\j*\b,-\i*\b) rectangle ++(\b,-\b);
     \else\ifnum\d=5
       \fill[red!25] (\j*\b,-\i*\b) rectangle ++(\b,-\b);
     \fi\fi\fi
   }
 }
 \draw[step=\b,black!30,thin] (0,0) grid (6*\b,-6*\b);
 \draw[thick] (0,0) rectangle (6*\b,-6*\b);
 \draw[decorate,decoration={brace,amplitude=4pt}] (0,0.12)--(\b,0.12)
   node[midway,above=4pt,font=\small] {$W$};
 \draw[decorate,decoration={brace,amplitude=5pt,mirror}] (0,-6*\b-0.12)--(6*\b,-6*\b-0.12)
   node[midway,below=5pt,font=\small] {$N=nW$};
 \draw[decorate,decoration={brace,amplitude=5pt}] (6*\b+0.12,0)--(6*\b+0.12,-6*\b)
   node[midway,right=6pt,font=\small,align=left] {$n$ blocks};
 \begin{scope}[shift={(6.6,-1.65)}]
  \foreach \a in {0,60,...,300}{
    \fill[blue!55] (\a-17:1.25) arc (\a-17:\a+17:1.25) -- (\a+17:0.85) arc (\a+17:\a-17:0.85) -- cycle;
  }
  \foreach \a in {30,90,...,330}{
    \draw[thick,blue!45] (\a-12:1.05) arc (\a-12:\a+12:1.05);
  }
  \node[font=\small,align=center] at (0,0) {$n$ blocks\\ on a ring};
 \end{scope}
\end{tikzpicture}
\caption{Left: the matrix written in $n$ blocks of size
$W$, here with $n=6$. Diagonal blocks (dark) are independent GUE matrices,
blocks next to the diagonal (light) are independent Gaussian couplings
between neighbouring blocks, and the two corner blocks (red) couple the
last block to the first. All other entries vanish. Right: the same
geometry seen as $n$ blocks on a ring, each coupled to its two
neighbours.}
\label{fig:block-ring}
\end{figure}

The goal of this paper is to answer both questions for one concrete
model, shown in Figure~\ref{fig:block-ring}. We split the $N$ coordinates
into $n$ consecutive blocks of $W$ coordinates each, so that $N=nW$, and
we arrange the blocks on a circle, so that the last block is a neighbour
of the first. We study a complex Gaussian Hermitian matrix whose entries
are nonzero only inside a block and between two neighbouring blocks. So,
written in blocks, it is a periodic block tridiagonal matrix, and its band
width is of order $W$. 
In this model the
critical scaling $W^2\asymp N$ reads $n\asymp W$, and the natural control
parameter is
\[
 \frac{n}W=\frac{N}{W^2}.
\]
We ask two questions about this model. First, when does every bulk
eigenvector spread evenly over all coordinates, in the sense of quantum
unique ergodicity described in (i) above? Second, when the eigenvectors do
not spread evenly, what is the joint law of the eigenvalues and of the way
their eigenvectors distribute mass over the blocks? Our three main
results, stated precisely in Subsection~\ref{subsec:main}, answer these
questions as follows. The first result answers the first question all the
way up to the threshold. The second and third answer the second question
exactly at the threshold.
\begin{enumerate}[label=(\alph*)]
\item If $n/W\to0$, quantum unique ergodicity holds at every deterministic
bulk rank (Theorem~\ref{thm:que}). We measure the uneven spreading of an
eigenvector by the mean-square deviation $\mathcal Q_\alpha(v)$ of
\eqref{eq:que-statistic}, which compares the mass that the eigenvector
puts on a test region with the size of that region. This deviation is
bounded by a constant multiple of $n/W$. So no power or logarithmic
separation between $W^2$ and $N$ is needed.
\item If $n/W\to\tau\in(0,\infty)$, the rescaled eigenvalues and the
block-mass profiles of their eigenvectors converge jointly, as a marked
point process, to a unique law, which we construct
(Theorem~\ref{thm:critical}).
\item In this critical law the eigenvectors are found not to be spread evenly. In fact, the
variance of their block masses against any nonconstant smooth test
function is strictly positive. This variance tends to zero as
$\tau\downarrow0$. The spectrum tends to the sine process as
$\tau\downarrow0$ and to the Poisson point process as $\tau\to\infty$, where both
limits are taken after the critical limit (Theorem~\ref{thm:endpoints}).
\end{enumerate}
Thus, for this model, the critical family of laws connects the two sides
of the conjecture, as predicted in~\cite{FM91}.

The rest of this section is organized as follows. In
Subsection~\ref{subsec:model} we define the model and the observables we
study. In Subsection~\ref{subsec:main} we state the main results. In
Subsection~\ref{subsec:comparison} we compare them with earlier work.
In Subsection~\ref{subsec:proof-ideas} we explain the main ideas of the
proofs. Further notation is collected in Subsection~\ref{subsec:notation}.

The remaining sections are organized in three parts, which follow the
three parts of Subsection~\ref{subsec:proof-ideas}. Section~\ref{sec:inputs}
fixes notation and collects the results from earlier work that we use, and
Section~\ref{sec:counting} proves the counting and density estimates used
throughout. Part~I, which proves Theorem~\ref{thm:que}, occupies
Sections~\ref{sec:coarea}--\ref{sec:angular-que}. Part~II, which proves
Theorem~\ref{thm:critical}, occupies
Sections~\ref{sec:transmission}--\ref{sec:t22-full-critical}, and
Section~\ref{sec:smoothed-hierarchy} records a consequence of the critical
law. Part~III, which proves Theorem~\ref{thm:endpoints}, uses
Section~\ref{sec:window-transfer} and occupies
Sections~\ref{sec:t22-signed-response}--\ref{sec:t22-positive-log-poisson}.

\subsection{The model and its observables}\label{subsec:model}

In this subsection we define the matrix, the masses of its eigenvectors
on the blocks, and the statistic that measures how evenly an eigenvector
spreads. We work with one specific Gaussian model. It is simple enough to
avoid technical complications that are not central to our results, and
its exact unitary invariance is used throughout the proofs.

Fix $b\in(0,1/2)$ and put $c=1-2b$. For
integers $n\ge4$ and $W\ge1$, let $N=nW$ and
\[
 \Lambda_{n,W}=(\Z/n\Z)\times\{1,\ldots,W\},\qquad
 I_k=\{k\}\times\{1,\ldots,W\}.
\]
The set $I_k$ is the $k$-th \emph{block}. We use the letters $k,\ell$ for
block indices in $\Z/n\Z$ and the letters $x,y$ for sites in
$\Lambda_{n,W}$. Let $H=H^*$ act on $\C^{\Lambda_{n,W}}$. Its
upper-triangular entries are independent centered Gaussian random
variables with variances
\begin{equation}\label{eq:model}
 s_{(k,a),(\ell,d)}:=\E|H_{(k,a),(\ell,d)}|^2
 =\frac{c\one_{\ell=k}+b\one_{\ell=k+1}+b\one_{\ell=k-1}}{W}.
\end{equation}
An off-diagonal entry of variance $s$ has independent real and imaginary
parts, each of variance $s/2$. A diagonal entry is real and has variance
$c/W$. The indices $k,\ell$ in \eqref{eq:model} are read modulo $n$.

In other words, the diagonal blocks are independent GUE matrices with
entry variance $c/W$. The forward bonds, that is, the blocks
$H_{I_kI_{k+1}}$, are independent complex Gaussian matrices with entry
variance $b/W$, and the reverse bonds are their adjoints. Here GUE means
the Gaussian law on Hermitian matrices that is invariant under conjugation
by unitary matrices; its normalization is fixed in
\eqref{eq:gue-normalizer}. Every row of $H$ has total variance $c+2b=1$,
so the empirical eigenvalue distribution of $H$ converges to Wigner's
semicircle law on $[-2,2]$ as $W\to\infty$. The condition $b\in(0,1/2)$
makes both $b$ and $c$ positive, so that both the blocks and the bonds are
random. The choice $b=1/3$ gives equal variances on the three blocks; this
is the model of~\cite{YY}. We denote the law of $H$ by $\Prob_{n,W,b}$. An
expectation $\E$ without a subscript is taken under this law.

Write $\lambda_1<\cdots<\lambda_N$ for the eigenvalues of $H$; the spectrum is simple almost surely. The integer $\alpha$ in $\lambda_\alpha$ is the
\emph{rank} of the eigenvalue. Choose normalized eigenvectors
$\psi_\alpha$. Let $P_k$ be the orthogonal projection onto $\C^{I_k}$, and
define
\begin{equation}\label{eq:mass}
 p_{\alpha k}=\norm{P_k\psi_\alpha}^2,\qquad
 \nu_\alpha=\sum_{k=0}^{n-1}p_{\alpha k}\delta_{k/n}\in\mathcal P(\T),
 \qquad \T=\R/\Z.
\end{equation}
The number $p_{\alpha k}$ is the \emph{mass} of $\psi_\alpha$ in block
$k$, and $\sum_kp_{\alpha k}=1$. The measure $\nu_\alpha$ is the
\emph{mass profile} of $\psi_\alpha$: it places the mass of block $k$ at
the point $k/n$ of the circle $\T$. Here $\mathcal P(\T)$ is the space of
probability measures on $\T$ with the topology of weak convergence. The
phases of the eigenvectors play no role in $p_{\alpha k}$ or $\nu_\alpha$.

For a deterministic real vector
$v=(v_x)_{x\in\Lambda_{n,W}}$, called a \emph{site test}, define
\begin{equation}\label{eq:que-statistic}
 \mathcal Q_\alpha(v)=
 \E\left|\sum_{x\in\Lambda_{n,W}}
       v_x\bigl(|\psi_\alpha(x)|^2-N^{-1}\bigr)\right|^2.
\end{equation}
If $\psi_\alpha$ were spread exactly evenly, then  $|\psi_\alpha(x)|^2=N^{-1}$
 for every $x$ and $\mathcal Q_\alpha(v)=0$. \emph{Quantum unique
ergodicity} (QUE), in the sense used in this paper, means that
$\mathcal Q_\alpha(v)$ tends to zero uniformly over deterministic bulk ranks
$\alpha$ and deterministic tests with $\norm v_\infty\le1$. A rank is in
the \emph{bulk} if $\vartheta\le\alpha/N\le1-\vartheta$ for a fixed
$\vartheta\in(0,1/2)$. ``Uniformly'' means that the suprema over $\alpha$
and $v$ are taken outside the expectation. It does not allow the test $v$
to be chosen after looking at $\psi_\alpha$.

\subsection{Main results}\label{subsec:main}

We now state three main theorems. The first concerns the regime $n/W\to0$. The
second and third concern the regime $n/W\to\tau\in(0,\infty)$.

\paragraph{QUE at a deterministic rank.} The first theorem bounds the QUE
statistic \eqref{eq:que-statistic} by a constant times $n/W$, for every
deterministic bulk rank and every bounded deterministic test.

\begin{theorem}[Prescribed-rank quantum unique ergodicity]\label{thm:que}
Fix $b\in(0,1/2)$, $\vartheta\in(0,1/2)$ and $T<\infty$.
There is $C=C_{b,\vartheta,T}<\infty$ such that, for $4\le n\le TW$,
uniformly over deterministic ranks $\alpha$ and real deterministic site tests $v$,
\begin{equation}\label{eq:que-bound}
 \sup_{\vartheta\le\alpha/N\le1-\vartheta}
 \sup_{\norm v_\infty\le1}
 \mathcal Q_\alpha(v)
 \le C\left[\frac{n}W+\left(\frac{n}W\right)^{3/2}\right]
        +\frac{1}{W+1}.
\end{equation}
Consequently the left side tends to zero whenever $W\to\infty$ and
$n/W\to0$. %
\end{theorem}

Theorem~\ref{thm:que} claims a bound by a constant factor, with no power
$N^\epsilon$ and no logarithmic factor. It therefore gives QUE under the
single condition $n/W\to0$, with no further relation between $n$ and $W$.
It is proved in Section~\ref{sec:angular-que}. 

\begin{remark}[Parameters in Theorem~\ref{thm:que}]\label{rem:que-parameters}
The condition $b\in(0,1/2)$ makes both $b$ and $c=1-2b$ positive, so that
both the blocks and the bonds carry Gaussian randomness. The constant may
deteriorate as $b$ approaches either endpoint. The parameter $\vartheta$
keeps the rank away from the two spectral edges.  The last term $(W+1)^{-1}$ comes from
passing from block tests to site tests; see Part~I of
Subsection~\ref{subsec:proof-ideas}. The theorem does not assert a bound for
the expectation of a supremum over $v$.
\end{remark}

\paragraph{The critical joint law.} To state the second theorem we rescale
the eigenvalues near a fixed energy. Fix $E\in(-2,2)$ and define
\begin{equation}\label{eq:unfolding}
 \rho(E)=\frac{\sqrt{4-E^2}}{2\pi},\qquad
 u_\alpha=N\rho(E)(\lambda_\alpha-E),\qquad
 \Xi_W=\sum_{\alpha=1}^N\delta_{(u_\alpha,\nu_\alpha)}.
\end{equation}
The number $\rho(E)$ is the semicircle density. After the rescaling, the
points $u_\alpha$ have mean spacing of order one near $u=0$. The random
measure $\Xi_W$ is a \emph{marked point process} on $\R\times\mathcal P(\T)$:
each point $u_\alpha$ carries the mark $\nu_\alpha$. A marked point process
is \emph{locally finite} if it has finitely many points in each bounded
spectral interval. Convergence in distribution of such processes is
tested by integration against continuous functions with compact support
in the spectral coordinate. Since $\mathcal P(\T)$ is compact, this is the
usual vague topology on $\R\times\mathcal P(\T)$.

The second theorem gives one limit law for the rescaled eigenvalues and
all their mass profiles together.

\begin{theorem}[Critical joint law]\label{thm:critical}
Fix $b\in(0,1/2)$, $E\in(-2,2)$ and $\tau\in(0,\infty)$. There is a random
marked point process $\Xi_{\tau,b,E}$, whose law depends only on $\tau$, $b$
and $E$, with the following properties. It is a random locally finite
collection of points $(u,\nu)$, where $u\in\R$ is a rescaled eigenvalue and
its mark $\nu\in\mathcal P(\T)$ is a probability measure on the circle,
the limiting mass profile of the corresponding eigenvector. If
$W\to\infty$ and $n/W\to\tau$, along any sequence and for either parity of
$n$, then
\begin{equation}\label{eq:critical-limit}
 \Xi_W\ \Longrightarrow\ \Xi_{\tau,b,E}
\end{equation}
in distribution, as random point measures on $\R\times\mathcal P(\T)$. Here
$\mathcal P(\T)$ carries the topology of weak convergence, and point
measures carry the vague topology, that is, convergence against continuous
test functions with compact support. Almost surely, no two points of
$\Xi_{\tau,b,E}$ have the same rescaled eigenvalue $u$. Moreover,
\begin{equation}\label{eq:critical-intensity}
 \E\sum_{(u,\nu)\in\Xi_{\tau,b,E}}\delta_u=du,
 \qquad
 \E\sum_{(u,\nu)\in\Xi_{\tau,b,E}}\delta_u\otimes\nu=du\,dr.
\end{equation}
The law of $\Xi_{\tau,b,E}$ is invariant under spectral translations and
under simultaneous rotations of all its marks. It is determined by the
construction in Theorem~\ref{thm:entire-law}.
\end{theorem}

Theorem~\ref{thm:critical} claims convergence of the full marked process,
not only of eigenvalue correlations. For every continuous spatial tests
$q_1,\ldots,q_k$, the process marked by
$(\nu(q_1),\ldots,\nu(q_k))$  converges jointly with the spectral locations.
One may order points locally in a bounded interval whose endpoints have
zero limiting probability of being spectral locations. The first identity in \eqref{eq:critical-intensity} says the
limit has one point per unit length on average. The second says that, on
average with respect to spectral intensity, the mark is the uniform measure
$dr$ on $\T$. The proof is
completed in Section~\ref{sec:entire-law}.

The law $\Xi_{\tau,b,E}$ is built from finite approximants. Each
approximant is the law of a finite system of stochastic differential
equations, so it is a probability law. Removing some of the spatial tests
gives the approximant for the remaining tests. The approximants share
common noises and one common matrix closing variable, so they retain
correlations between marks. We do not assert that the marks are
independent of each other, or of the eigenvalues.

\paragraph{Profile variance and the two endpoints.} The third theorem
shows that the critical marked law is not concentrated on uniform profiles,
and identifies the limits
of the spectrum as $\tau\downarrow0$ and $\tau\to\infty$. For a real smooth
function $\phi$ on $\T$, write $\bar\phi=\int_\T\phi(r)\,dr$ and
$\nu(\phi)=\int_\T\phi\,d\nu$. The quantity $|\nu_\alpha(\phi)-\bar\phi|^2$
measures how far the profile $\nu_\alpha$ is from uniform when tested
against $\phi$.

\begin{theorem}[Profile variance and spectral endpoints]\label{thm:endpoints}
Fix $b\in(0,1/2)$ and $E\in(-2,2)$, and consider the family of marked
point processes $\Xi_{\tau,b,E}$ of Theorem~\ref{thm:critical} for all
$\tau>0$. For every $\tau>0$, every $w\in C_c^\infty(\R)$ with $w\ge0$ and
$\int_\R w=1$, and every smooth real nonconstant $\phi$ on $\T$,
\begin{equation}\label{eq:variance-limit}
 \lim_{\substack{W\to\infty\\n/W\to\tau}}
 \E\sum_\alpha w(u_\alpha)
             |\nu_\alpha(\phi)-\bar\phi|^2
 =V_{\tau,b,E}(\phi)>0.
\end{equation}
The limit does not depend on $w$, and
\begin{equation}\label{eq:small-variance}
 \lim_{\tau\downarrow0}V_{\tau,b,E}(\phi)=0.
\end{equation}
Let $\mathcal X_{\tau,b,E}$ be the \emph{spectral part} of $\Xi_{\tau,b,E}$, i.e., 
$
 \mathcal X_{\tau,b,E}=\sum_{(u,\nu)\in\Xi_{\tau,b,E}}\delta_u .$
Then, in
distribution with respect to the vague topology on point measures on $\R$,
\begin{equation}\label{eq:spectral-endpoints}
 \mathcal X_{\tau,b,E}\Longrightarrow\mathcal S
       \quad(\tau\downarrow0),\qquad
 \mathcal X_{\tau,b,E}\Longrightarrow\mathcal P
       \quad(\tau\to\infty).
\end{equation}
Here $\mathcal S$ is the determinantal point process of intensity one
with kernel
\[
 K(x,y)=\frac{\sin(\pi(x-y))}{\pi(x-y)},\qquad K(x,x)=1.
\]
The process $\mathcal P$ is the Poisson process with intensity one, that
is, its counts in disjoint bounded intervals are independent Poisson
variables whose means are the lengths of the intervals.
\end{theorem}

Theorem~\ref{thm:endpoints} claims three things. The limiting profile
variance $V_{\tau,b,E}(\phi)$ is strictly positive for every
nonconstant $\phi$, and, by uniform integrability
(see the proof of the variance assertions in
Section~\ref{sec:t22-signed-response}), it is also the corresponding
expectation under $\Xi_{\tau,b,E}$, so the limiting marked law is not
concentrated on the uniform profile. This averaged conclusion does not assert that every mark
is nonuniform. The
 variance tends to zero as $\tau\downarrow0$, which matches
Theorem~\ref{thm:que}. The spectral part of the critical law tends to the
sine process as $\tau\downarrow0$ and to the Poisson process as
$\tau\to\infty$.

\begin{remark}[Order of limits and the weight $w$]\label{rem:endpoint-parameters}
The limits in \eqref{eq:spectral-endpoints} are taken after the critical
limit \eqref{eq:critical-limit}. They do not assert convergence along an
arbitrary sequence with $n/W\to\infty$, or a joint profile limit at either
endpoint. Moreover, the endpoint statements
concern only the eigenvalues. They do not describe the joint law of
eigenvalues and profiles as $\tau\to0$ or $\tau\to\infty$. The weight $w$ in
\eqref{eq:variance-limit} lets us average the profile fluctuations over the
eigenvalues near $E$ without conditioning on any event. Because
$\int w=1$ and the critical law is invariant under spectral translations
(Theorem~\ref{thm:critical}), this average does not depend on $w$.
\end{remark}

\subsection{Comparison with earlier results}\label{subsec:comparison}

In this subsection we compare Theorems~\ref{thm:que}, \ref{thm:critical}, and \ref{thm:endpoints}
with earlier results. For each work we recall the model, the condition on
the band width, and the statement, and then explain how it relates to
ours. Three features matter throughout: the condition on the band width,
the observable that is controlled, and the order of limits. Throughout, $N=nW$ denotes the total matrix dimension,
as in most of the cited papers. 

\paragraph{Delocalization phase}
Yau and Yin~\cite{YY} consider the Gaussian periodic block model with
\[
 \E|H_{xy}|^2=\frac{1}{3W}\,\one\bigl(|[x]-[y]|\le1\bigr),
\]
where $[x]$ is the block of $x$. This is \eqref{eq:model} with $b=1/3$.
Assume $W\ge N^{1/2+\mathfrak c}$ for a fixed $\mathfrak c>0$. Their
delocalization theorem~\cite[Theorem 2.2]{YY} states that for every fixed
$\kappa,\delta>0$ and $D>0$,
\begin{equation}\label{eq:cmp-yy-deloc}
 \Prob\Bigl(\max_{k:\,|\lambda_k|\le2-\kappa}\norm{\psi_k}_\infty^2
     \le N^{-1+\delta}\Bigr)\ge1-N^{-D}
\end{equation}
for all large $N$. Their quantum unique ergodicity
theorem~\cite[Theorem 2.5]{YY} concerns sets of blocks. Fix
$0<\tau_*<\mathfrak c/2$ and an energy $|E|<2-\kappa$, and let $J_E$ be
the set of ranks $k$ with
$|\lambda_k-E|\le N^{-1-\tau_*}(W^2/N)^{1/3}$. Then, uniformly over sets
$A$ of blocks,
\begin{equation}\label{eq:cmp-yy-que}
 \Prob\Biggl(\max_{k\in J_E}
   \Bigl|\sum_{a\in A}\norm{P_a\psi_k}^2-\frac{|A|W}N\Bigr|
   \ge\frac{|A|W}{N^{1+\tau_*/6}}\Biggr)\le N^{-\tau_*/6}.
\end{equation}
They also prove local laws, quantum diffusion and bulk
universality~\cite[Theorems 2.3, 2.4 and 2.6]{YY}. Erd\H{o}s and Riabov~\cite{ER} treat a much broader class. Their matrices
have a general variance profile $S$ with $\sum_aS_{ab}=1$ and
$S_{ab}\le C/W$, entries with all moments bounded, and either symmetry
class~\cite[Definition 2.1]{ER}; a block profile with three equal
variances is one of their examples. They assume
$W^2\ge N^{1+\zeta}$ for a fixed $\zeta>0$. For every fixed $\xi>0$ and
every deterministic diagonal matrix $A$, they prove that, with very high
probability~\cite[Corollary 3.5]{ER},
\begin{equation}\label{eq:cmp-er}
 \max_{i,j:\,|\lambda_i|,|\lambda_j|\le2-\kappa}
 \Bigl|\ip{\psi_i}{A\psi_j}-\delta_{ij}\frac{\Tr A}N\Bigr|
 \le N^{\xi}\sqrt{\frac{N}{W^2}}\,\norm A.
\end{equation}

The bound \eqref{eq:cmp-er} is the natural point of comparison for
Theorem~\ref{thm:que}. Take $i=j=\alpha$ and $A=\diag(v)$. Then the
quantity inside the absolute value in \eqref{eq:cmp-er} is exactly the
random variable whose second moment is $\mathcal Q_\alpha(v)$ in
\eqref{eq:que-statistic}. Theorem~\ref{thm:que} therefore gives, for
$n\le W$,
\begin{equation}\label{eq:cmp-ours}
 \Bigl(\E\Bigl|\ip{\psi_\alpha}{\diag(v)\psi_\alpha}
   -\frac{\Tr\diag(v)}N\Bigr|^2\Bigr)^{1/2}
 \le C\sqrt{\frac{N}{W^2}}+\frac{1}{\sqrt{W+1}},
 \qquad \norm v_\infty\le1.
\end{equation}
So the two bounds have the same main scale $\sqrt{N/W^2}$. The bounds
\eqref{eq:cmp-yy-que} and \eqref{eq:cmp-er} hold with high probability simultaneously over many
ranks. On the other hand, the bound \eqref{eq:cmp-ours} has a
constant in place of the factor $N^{\xi}$. This difference decides which
band widths are covered. The exponents $\mathfrak c,\tau_*,\zeta,\xi$
in \eqref{eq:cmp-yy-deloc}--\eqref{eq:cmp-er} are fixed and positive, and
they cannot be chosen to depend on $N$. For example, let
$n=\lfloor W/\log\log W\rfloor$. Then
\[
 \frac{N}{W^2}\asymp\frac{1}{\log\log W}\longrightarrow0,
 \qquad\text{but}\qquad
 N^{\xi}\sqrt{\frac{N}{W^2}}\longrightarrow\infty
 \quad\text{for every fixed }\xi>0,
\]
so \eqref{eq:cmp-er} gives no information, and the condition of
\eqref{eq:cmp-yy-deloc} fails. In contrast, the right side of
\eqref{eq:cmp-ours} tends to zero. What is new in
Theorem~\ref{thm:que} is thus the mean-square bound with a constant
factor, and the resulting QUE under the single condition $n/W\to0$, for
the specific Gaussian block ring \eqref{eq:model}.

\paragraph{Two-point correlations above the threshold.}
Shcherbina and Shcherbina~\cite{SS19} study the Gaussian Hermitian block
model on a path of $n$ blocks, each of size $W$, with $N=nW$. The entries
in the blocks $(j,k)$ have variance $J_{jk}$, where
\begin{equation}\label{eq:cmp-ss19-model}
 J=\frac{1}W\bigl(I+\beta\Delta\bigr),\qquad \beta<\frac{1}{4},
\end{equation}
and $\Delta$ is the discrete Laplacian on $\{1,\ldots,n\}$ with Neumann
boundary conditions. Away from the two ends of the path, $J$ puts
variance $(1-2\beta)/W$ inside a block and $\beta/W$ between neighboring
blocks. So in the interior it is the profile \eqref{eq:model} with
$b=\beta$. The differences are the boundary condition and the range
$\beta<1/4$. Let $R_2$ denote the two-point correlation function of the
eigenvalues and $\rho$ the semicircle density. Their main
result~\cite[Theorem 1.1]{SS19} states that, for every $E\in(-2,2)$,
\begin{equation}\label{eq:cmp-ss19}
 \frac{1}{(N\rho(E))^2}\,
 R_2\Bigl(E+\frac{\xi_1}{\rho(E)N},\,E+\frac{\xi_2}{\rho(E)N}\Bigr)
 \longrightarrow
 1-\frac{\sin^2(\pi(\xi_1-\xi_2))}{\pi^2(\xi_1-\xi_2)^2}
\end{equation}
as $n,W\to\infty$ with $W\ge Cn\log^5n$. The proof uses a supersymmetric
integral representation and a transfer operator for ratios of
characteristic polynomials. In our notation the condition reads
$n/W\le(C\log^5n)^{-1}$, so it requires $n/W\to0$ at a logarithmic rate.
Their observable is also different. For instance, \eqref{eq:cmp-ss19} is a statement
about pairs of eigenvalues, and it says nothing about eigenvectors.

A related result is the sigma-model approximation of~\cite{SS18}. There
the variance matrix is
\[
 J=\frac{1}W+\frac{\beta\Delta}{W^2},\qquad \beta>0,
\]
so the coupling between blocks is of order $W^{-2}$ rather than $W^{-1}$.
The limit is taken in two steps. First $W\to\infty$ with $n$ and $\beta$
fixed, this produces a sigma model~\cite[Theorems 1.1--1.2]{SS18}. Then
$n,\beta\to\infty$ with $\beta\ge Cn\log^2n$, and the two-point function
of the sigma model converges to the GUE limit~\eqref{eq:cmp-ss19} for
$|E|<\sqrt2$~\cite[Theorems 1.3--1.4]{SS18}. Both the scaling of the
coupling and the order of limits differ from our setting, in which $b$ is
fixed and $W\to\infty$ together with $n$.

\paragraph{Characteristic polynomials at the threshold.}
Shcherbina~\cite{Sh19} considers the Gaussian Hermitian band matrix whose
variance matrix is
\[
 J=\bigl(-W^2\Delta+I\bigr)^{-1},
\]
with $\Delta$ the Laplacian on $\{1,\ldots,N\}$ with Neumann boundary
conditions. Its entries are exponentially small when $|i-j|\gg W$. The
observable is the normalized second correlation function of
characteristic polynomials,
\[
 \bar F_2(x_1,x_2)=\frac{\E\bigl[\det(x_1-H)\det(x_2-H)\bigr]}{\E\bigl[\det(E-H)^2\bigr]},
 \qquad
 x_{1,2}=E\pm\frac{\xi}{2N\rho(E)}.
\]
Earlier works, recalled in~\cite[Theorem 1.1]{Sh19}, showed that
$\bar F_2(x_1,x_2)$ tends to $\sin(\pi\xi)/(\pi\xi)$ when
$W\ge N^{1/2+\theta}$, and to $1$ when $1\ll W\le\sqrt{N/(C\log N)}$.
At the threshold $N=C_*W^2$, Shcherbina proves~
\cite[Theorem 1.2]{Sh19} 
\begin{displaymath}\label{eq:cmp-sh19}
 \lim_{N\to\infty}\bar F_2\Bigl(E+\frac{\xi}{2N\rho(E)},\,
   E-\frac{\xi}{2N\rho(E)}\Bigr)
 =\bigl(e^{-C^*\Delta_U-i\xi\hat\nu}\,1,\,1\bigr),
 \qquad C^*=\frac{C_*}{(2\pi\rho(E))^2}.
\end{displaymath}
Here $(\cdot,\cdot)$ is the inner product on the two-dimensional sphere,
$\Delta_U$ is a Laplace operator on the sphere,
and $\hat\nu$ is a
multiplication operator.
As $C^*\to0$ and $C^*\to\infty$,  the right side
recovers the two limits above~
\cite[Remark 1.1]{Sh19}.

Of all the earlier results, \eqref{eq:cmp-sh19} is the closest in spirit
to Theorems~\ref{thm:critical} and~\ref{thm:endpoints}. The constant
$C_*=N/W^2$ plays the role of our $\tau=n/W=N/W^2$, and its two limits
mirror the sine and Poisson endpoints in
\eqref{eq:spectral-endpoints}. The difference lies in the observable.
The function $\bar F_2$ is one average of a product of determinants. It
does not determine the joint law of the eigenvalues, and it carries no
information about eigenvectors. Theorem~\ref{thm:critical} identifies
the joint law of all rescaled eigenvalues together with all their mass
profiles.

\paragraph{The continuum Anderson Hamiltonian in one dimension.}
A complete picture of a one-dimensional crossover, including
eigenvectors, is available for a continuum model. Dumaz and
Labb\'e~\cite{DL20,DL23,DL24} study the Anderson Hamiltonian
\[
 \mathcal H_L=-\partial_t^2+dB
\]
on the segment $[0,L]$ with Dirichlet boundary conditions, where $dB$ is
a white noise potential. In this model the disorder is fixed, and the
role of our band width is played by the energy level $E$ at which one
looks. For $1\ll E\ll L$ the eigenfunctions are localized on the scale
$E$. In this regime, let $n(E)$ be the density of states and $U_i$ the
centre of mass of the eigenfunction $\varphi_i$. The pairs
$(Ln(E)(\lambda_i-E),\,U_i/L)$ then converge to a Poisson point process
of intensity $d\lambda\otimes du$~\cite[Theorem 1]{DL24}. So the
transition is expected at $E\asymp L$, and the ratio $L/E$ corresponds to
our ratio $N/W^2=n/W$.

The critical regime is treated in~\cite{DL23}. Fix $\tau>0$, and $B,\Re\mathcal W,\Im\mathcal W$ be independent standard real
 Brownian motions. For
$\lambda\in\R$, define the diffusions $\Theta_\lambda$ and $\rho_\lambda$
on $[0,1]$ by $\Theta_\lambda(0)=\rho_\lambda(0)=0$ and
\begin{equation}\label{eq:cmp-dl-sde}
 d\Theta_\lambda=\frac{\lambda}{2}\,dt+\frac{\sqrt\tau}{2}\,dB
   +\frac{\sqrt\tau}{2\sqrt2}\,\Re\bigl(e^{2i\Theta_\lambda}d\mathcal W\bigr),
 \qquad
 d\rho_\lambda=\frac{\tau}{8}\,dt
   +\frac{\sqrt\tau}{2\sqrt2}\,\Im\bigl(e^{2i\Theta_\lambda}d\mathcal W\bigr).
\end{equation}
The point process
$\mathrm{Sch}_\tau=\{\lambda\in\R:\Theta_\lambda(1)\in\pi\Z\}$ was
introduced by Kritchevski, Valk\'o and Vir\'ag~\cite{KVV12} as the
scaling limit of critical discrete random Schr\"odinger operators. Let
$E=E(L)\sim L/\tau$, and let $\ell_E\in[0,\pi)$ be a deterministic shift
of order one, fixed by the boundary condition. Dumaz and Labb\'e
prove~\cite[Theorem 1]{DL23} that the point process
\begin{equation}\label{eq:cmp-dl}
 \Bigl\{\Bigl(\tfrac L{\sqrt E}(\lambda-E)+2\ell_E,\;
   L\varphi_\lambda^2(Lt)\,dt\Bigr)\Bigr\}
 \;\Longrightarrow\;
 \Bigl\{\Bigl(\lambda,\;
   \frac{e^{2\rho_\lambda(t)}\,dt}{\int_0^1e^{2\rho_\lambda(s)}\,ds}\Bigr):
   \lambda\in\mathrm{Sch}_\tau\Bigr\}
\end{equation}
converges in law on $\R\times\mathcal P([0,1])$, with the vague topology
in the first coordinate. Here the left side runs over the eigenvalues
$\lambda$ of $\mathcal H_L$, and $\varphi_\lambda$ is the normalized
eigenfunction. They also identify the intensity of the limit
explicitly~\cite[Theorem 3]{DL23}. They construct a random operator
$\mathtt{CS}_\tau$ whose eigenvalues and eigenvectors realize the limit,
and they prove convergence to it at the level of operators~\cite[Theorems
2, 4 and 5]{DL23}. In particular, the typical eigenvector has the shape
predicted by Rifkind and Vir\'ag~\cite{RV18}.

Structurally, \eqref{eq:cmp-dl} is the closest earlier result to
Theorem~\ref{thm:critical}. In both cases the limit is a point process of
eigenvalues marked by normalized mass profiles. In both cases it is
indexed by the ratio of the system size to the localization length. The
models and the endpoints differ, however. The operator $\mathcal H_L$ has
a single channel, and the transition is reached by moving the energy. Our
matrix has $W$ channels in each block; the energy stays in the bulk, and
the transition is reached by moving the band width. As $\tau\downarrow0$, the noise in the Dumaz--Labb\'e model disappears, and
its eigenvalues become those of the free Laplacian, without disorder. In
the rescaled spectral variable these are evenly spaced at distance $2\pi$,
forming the deterministic lattice $2\pi\Z$, often called a \emph{picket
fence}. In contrast, in Theorem~\ref{thm:endpoints}
the limit as $\tau\downarrow0$ is the sine process. Heuristically, a
single-channel system behaves like the free Laplacian on scales shorter
than its localization length, while a block band matrix behaves like a
GUE matrix. As $\tau\to\infty$, both families tend to the  Poisson processes,
with intensity $d\lambda/(2\pi)$ in the displayed $\mathrm{Sch}_\tau$
coordinate unlike the unit intensity case in our RBM setting. The precise
stationary comparison could be found in \mbox{%
\cite[Theorem 1.4(2)]{Nakano15}}\hskip0pt%
, for
$\tau=4/\beta$. %

 Finally, in \eqref{eq:cmp-dl} the profile mark lives on a
segment with Dirichlet ends, and it is an explicit functional of the same
Brownian motions that drive $\Theta_\lambda$. In
Theorem~\ref{thm:critical} the marks are block-mass profiles on the
circle $\T$.

\paragraph{Transmission diffusions.}
Bachmann, Butz and De Roeck~\cite{BBDR} study the Anderson model on a
strip with a fixed number $N$ of channels, a model of a disordered
quantum wire. In a limit of weak disorder $\lambda\to0$, with the
transverse hopping also tending to zero and time rescaled by
$\lambda^{-2}$, the transfer matrices converge to the solution of a
linear stochastic differential equation~\cite[Theorem 8]{BBDR}. In the
unitary class, the induced transmission eigenvalues $T_1,\ldots,T_N$
follow the Dorokhov--Mello--Pereyra--Kumar (DMPK)
process~\cite[Corollary 10]{BBDR}. This is the diffusion that describes
how the channel transmission probabilities of a disordered wire evolve
with its length $s$. For the conductance $g_N(s)=\sum_kT_k(s)$ of this
diffusion, they also prove Ohm's law~\cite[Theorem 2]{BBDR}, i.e., for every
$p\ge1$ and $T>0$,
\begin{equation}\label{eq:cmp-ohm}
 \lim_{N\to\infty}\frac{\E\,g_N(s)^p}{N^p}=\frac{1}{(1+s)^p}
 \qquad\text{uniformly in }s\in[0,T].
\end{equation}
The number of channels $N$ is fixed in the scaling limit and is sent to
infinity only afterwards, at bounded time $s$. Our proof also uses a transmission diffusion, but it differs from that of
\cite{BBDR} in three ways. First, we do not start from a weak-disorder
limit. We derive the transmission, together with the spatial sources that
record where the response lies inside the path, directly from the Gaussian
blocks of \eqref{eq:model}. Second, in~\cite{BBDR} the number of channels
is fixed in the limit and sent to infinity only afterwards, at bounded
time. In our model the number of channels is the block size $W$, which
tends to infinity together with the length, and we control all channels
at once, at the critical time scale. Third, a transmission diffusion
describes an open wire. To obtain the spectrum of the ring, we must also
close the wire, and we identify the law of the matrix that does this.

\paragraph{Localized side.}
Drogin~\cite{Drogin} considers real symmetric block tridiagonal matrices
on a path of $n$ blocks of size $W$. The entries of
$\sqrt W H$ are independent up to symmetry and \emph{$M$-regular}, i.e., they
are centered, have variance at most one and fourth moment at most $M$,
and have a $C^2$ density whose logarithm has second derivative bounded by
$M$. The main result in~\cite[Theorem 1]{Drogin} states that for every
$E_0>0$ and all sites $x,y$,
\begin{equation}\label{eq:cmp-drogin}
 \E\sum_{j:\,|\lambda_j|\le E_0}|\psi_j(x)|\,|\psi_j(y)|
 \le 2W^{C}\exp\Bigl(-c\,\frac{|x-y|}{W^2}\Bigr),
\end{equation}
with $C,c>0$ depending only on $M$ and $E_0$. The bound gives
exponential localization on the scale $W^2$, with a polynomial prefactor
in $W$. The right side is small once $|x-y|$ is much larger than
$W^2\log W$. In particular it gives localization when $W^2\log W\ll N$.
It is derived from a fractional moment bound on the
resolvent~\cite[Theorem 2]{Drogin}.

Fan and Zou~\cite{FZ} study the real symmetric Gaussian periodic block
model with three equal variances, that is,
$\E H_{xy}^2=(3W)^{-1}\one(|[x]-[y]|\le1)$. Let $\rho_W$ be the limit of
the density of states as $N\to\infty$ with $W$ fixed. They prove~\cite[Theorem 1.3]{FZ} that for every fixed $E\in(-2,2)$,
\begin{equation}\label{eq:cmp-fz}
 \frac{W^2\log N}N\longrightarrow0 \qquad \text{implies} \qquad \sum_{k=1}^N\delta_{N\rho_W(E)(\lambda_k-E)}
 \Longrightarrow \mathrm{PPP}(d\lambda),
\end{equation}
where $\mathrm{PPP}(d\lambda)$ is the Poisson point process of intensity
one, and the convergence is in distribution in the vague topology. They
leave the complex Hermitian case for later work.

The localized-side results \eqref{eq:cmp-drogin} and \eqref{eq:cmp-fz}
hold along every sequence in their regime. Our Poisson conclusion in
Theorem~\ref{thm:endpoints} is of a different kind. It is the endpoint
$\tau\to\infty$ of the critical complex family, taken after the critical
limit \eqref{eq:critical-limit}. We do not infer from it the
corresponding joint limit in $W$ and $N$ without this order of limits.

\paragraph{Other results that are used in our proof.}
Several previous results about mean-field matrices are used inside our
proof. We state all these results with details in Section~\ref{sec:inputs}. The block-Gaussian
counting estimates of Peled, Schenker, Shamis and Sodin~\cite{PSSS}
control arbitrary short spectral intervals, uniformly in $W$. The local
law of Cipolloni, Erd\H{o}s and Henheik~\cite{CEH} for deformed matrices
allows Hermitian deformations of polynomially large norm. The
exact complex Gaussian kernel of Johansson~
\cite{JohanssonKernel},
with the uniform contour proof given below, identifies the local
spectrum of a newly exposed Gaussian block. The eigenvector moment
flow~\cite{BY,Benigni}, the finite speed of propagation of that
flow~\cite{BHY}, and the direct eigenvector comparison of Benigni and
Lopatto~\cite{BL} supply tools for the conditional local eigenvector
calculation.  The relative-frame and conditional forms  needed here are proved in
this paper with the precise Gaussian laws and normalizations. %

\subsection{Proof ideas}\label{subsec:proof-ideas}

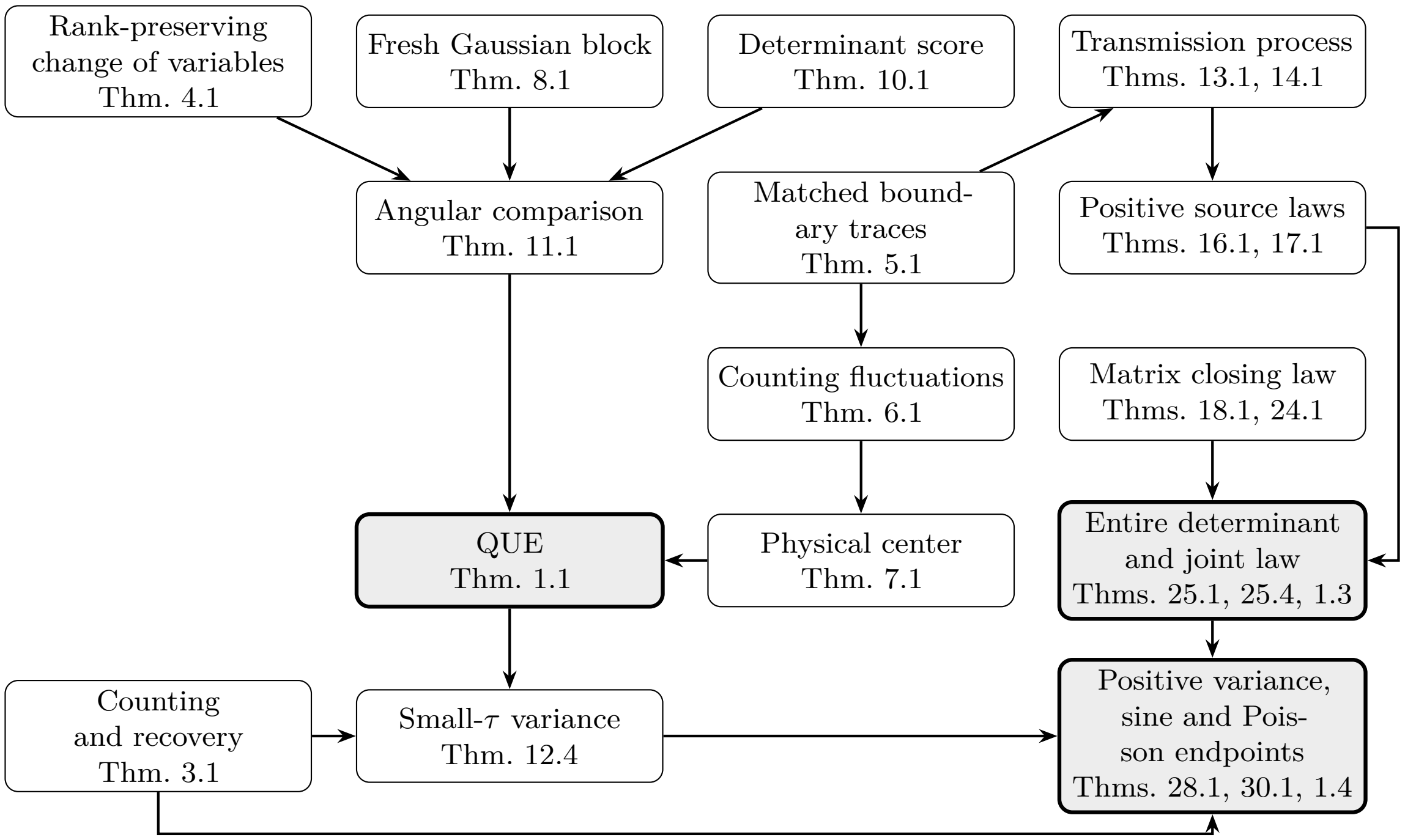
\begin{figure}[ht]
\centering
\begin{tikzpicture}[x=1cm,y=1cm,
 box/.style={draw,rounded corners,align=center,text width=3.1cm,
              inner sep=3pt,font=\footnotesize,minimum height=1cm},
 main/.style={box,fill=black!7,very thick},
 arr/.style={-{Stealth[length=2mm]},thick}]
 \node[box] (coarea) at (0,0) {Rank-preserving change of variables\\Thm.~\ref{thm:coarea}};
 \node[box] (palm) at (3.8,0) {Fresh Gaussian block\\Thm.~\ref{thm:fresh-palm}};
 \node[box] (score) at (7.6,0) {Determinant score\\Thm.~\ref{thm:determinant-score}};
 \node[box] (trans) at (11.4,0) {Transmission process\\Thms.~\ref{thm:transmission-generator}, \ref{thm:transmission-law}};
 \node[box] (angular) at (3.8,-1.8) {Angular comparison\\Thm.~\ref{thm:physical-angular}};
 \node[box] (layers) at (7.6,-1.8) {Matched boundary traces\\Thm.~\ref{thm:matched-traces}};
 \node[box] (source) at (11.4,-1.8) {Positive source laws\\Thms.~\ref{thm:source-process}, \ref{thm:complex-source-law}};
 \node[box] (rankloc) at (7.6,-3.6) {Counting fluctuations\\Thm.~\ref{thm:rank-location}};
 \node[box] (closing) at (11.4,-3.6) {Matrix closing law\\Thms.~\ref{thm:relative-mixed-vectors}, \ref{thm:matrix-closing}};
 \node[main] (que) at (3.8,-5.4) {QUE\\Thm.~\ref{thm:que}};
 \node[box] (center) at (7.6,-5.4) {Physical center\\Thm.~\ref{thm:physical-center}};
 \node[main] (entire) at (11.4,-5.4) {Entire determinant and joint law\\Thms.~\ref{thm:entire-law}, \ref{thm:joint-recovery}, \ref{thm:critical}};
 \node[box] (count) at (0,-7.3) {Counting and recovery\\Thm.~\ref{thm:counting}};
 \node[box] (small) at (3.8,-7.3) {Small-$\tau$ variance\\Thm.~\ref{thm:small-critical}};
 \node[main] (ends) at (11.4,-7.3) {Positive variance, sine and Poisson endpoints\\Thms.~\ref{thm:t22-sine-endpoint}, \ref{thm:t22-Poisson-endpoint}, \ref{thm:endpoints}};
 \draw[arr] (coarea)--(angular);
 \draw[arr] (palm)--(angular);
 \draw[arr] (score)--(angular);
 \draw[arr] (angular)--(que);
 \draw[arr] (layers)--(rankloc);
 \draw[arr] (rankloc)--(center);
 \draw[arr] (center)--(que);
 \draw[arr] (layers)--(trans);
 \draw[arr] (trans)--(source);
 \draw[arr] (source.east)--++(.35,0)|-(entire.east);
 \draw[arr] (closing)--(entire);
 \draw[arr] (que)--(small);
 \draw[arr] (count)--(small);
 \draw[arr] (small)--(ends);
 \draw[arr] (entire)--(ends);
 \draw[arr] (count.south)--++(0,-.45)-|(ends.south);
\end{tikzpicture}
\caption{Flow of the argument. Each box names a main result, and each arrow
points from a result to a later result whose proof uses it. The left part
is the route to Theorem~\ref{thm:que} (Part~I of
Subsection~\ref{subsec:proof-ideas}): the exact change of variables of
Theorem~\ref{thm:coarea}, the conditional law of a fresh Gaussian block and
the determinant score are the inputs of the angular comparison, while the
matched boundary traces and the counting fluctuations lead to the bound on
the physical center. The right column is the route to
Theorem~\ref{thm:critical} (Part~II): the transmission and source processes
and the matrix closing law are separate inputs to the entire determinant.
The bottom row is Part~III. Shaded boxes contain the three main theorems.}
\label{fig:argument}
\end{figure}

In this subsection we explain how the three main theorems are proved. Part~I concerns Theorem~\ref{thm:que}, Part~II concerns Theorem~\ref{thm:critical}, and Part~III concerns Theorem~\ref{thm:endpoints}. In each part we first state the estimate we want. We then reduce it, step by step, to simpler statements. After each step we say what remains to be proved and why the remaining statement is difficult, and then we explain how we prove it. Figure~\ref{fig:argument} %
shows how the main results of the paper fit together.

\paragraph{Part I: quantum unique ergodicity at a fixed rank
(Theorem~\ref{thm:que}).}

\emph{The goal.} Fix a deterministic bulk rank $\alpha$ and a
deterministic site test $v$ with $\norm v_\infty\le1$. We want to prove
\begin{equation}\label{eq:intro-que-goal}
 \mathcal Q_\alpha(v)
 =\E\Bigl|\sum_xv_x\bigl(|\psi_\alpha(x)|^2-N^{-1}\bigr)\Bigr|^2
 \le C\Bigl[\frac{n}W+\Bigl(\frac{n}W\Bigr)^{3/2}\Bigr]+\frac{1}{W+1},\end{equation}
with a constant $C$ that does not depend on $n$, $W$, $\alpha$ or $v$.

\emph{Why the known methods do not suffice.} The results of~\cite{YY,ER}
are proved in the sense of stochastic domination
(Subsection~\ref{subsec:notation}). Such estimates lose a fixed factor $N^\xi$; see \eqref{eq:cmp-er}. This loss is harmless when $W^2\ge N^{1+\zeta}$, but it is fatal when $W^2/N$ grows slower than any
power, for instance when $n=\lfloor W/\log\log W\rfloor$. A second
difficulty is that the eigenvalue $\lambda_\alpha$ of the fixed rank
$\alpha$ is random and depends on every entry of $H$. An estimate at a
fixed energy, or an estimate averaged over many ranks, therefore does not
directly control $\mathcal Q_\alpha(v)$.

\emph{Step 1: from site tests to Fourier modes of the block masses.} We
first remove the fine structure inside the blocks. Conjugating $H$ by
independent Haar unitary matrices inside its blocks does not change its
law or the rank of any eigenvalue. Conditionally on the block masses
$p_{\alpha k}$, the rotated eigenvector inside each block is then a
uniform random vector on a sphere. Computing its second moments, we
obtain the exact identity \eqref{eq:haar-site-bridge}, which reads
\[
 \mathcal Q_\alpha(v)
 =\E\Bigl|\sum_k\bar v_k\Bigl(p_{\alpha k}-\frac{1}n\Bigr)\Bigr|^2
 +\frac{1}{W+1}\,\E\sum_kp_{\alpha k}^2
       \bigl(\overline{v^2}_k-\bar v_k^2\bigr),
\]
where $\bar v_k$ and $\overline{v^2}_k$ are the averages of $v$ and $v^2$
over block $k$. The second term is at most $(W+1)^{-1}$ which provides the last
term in \eqref{eq:intro-que-goal}. It therefore suffices to bound the
first term, which involves only a test $t_k=\bar v_k$ that is constant on
blocks.

Next we expand the block masses in Fourier modes on the ring. For
$1\le m<n$ put
\[
 z_m=e^{2\pi im/n},\qquad
 A_m=\sum_{k=0}^{n-1}z_m^{-k}p_{\alpha k},\qquad
 \gamma_m=\E|A_m|^2,\qquad
 d_m=|1-z_m|.
\]
The law of $H$ is invariant under cyclic shifts of the blocks, so distinct
nonzero modes are orthogonal in expectation. Together with Parseval's
identity, this gives
\begin{equation}\label{eq:intro-parseval}
 \E\Bigl|\sum_kt_k\Bigl(p_{\alpha k}-\frac{1}n\Bigr)\Bigr|^2
 =\sum_{m\ne0}|\widehat t_m|^2\gamma_m,
 \qquad
 \sum_{m}|\widehat t_m|^2\le1,
\end{equation}
where $\widehat t_m=n^{-1}\sum_kt_kz_m^k$; see the proof of
Theorem~\ref{thm:que} in Section~\ref{sec:angular-que}. Thus
\eqref{eq:intro-que-goal} follows from the bound
\begin{equation}\label{eq:intro-mode-goal}
 \gamma_m\le C\Bigl\{\frac{1}{Wd_m}+\frac{1}{\sqrt n\,W^{3/2}d_m^2}\Bigr\}
 \qquad (1\le m<n),
\end{equation}
using $d_m=2\sin(\pi m/n)\ge4/n$.

\emph{What remains, and the challenge.} We must bound one Fourier mode
$\gamma_m$ of the block masses at the fixed rank $\alpha$, without any
factor $N^\xi$. The difficulty is still the random energy
$\lambda_\alpha$.

\emph{Step 2: a size-biased law turns $\gamma_m$ into a statement about
one pair of blocks.} Fix an adjacent pair $S=\{i,i+1\}$ and let
$m_S=p_{\alpha i}+p_{\alpha,i+1}$ be the mass of $\psi_\alpha$ on $S$. We
define the size-biased law
\begin{equation}\label{eq:intro-cancel}
 d\Prob^\sharp_{\alpha,S}=\frac{n}{2}\,m_S\,d\Prob,
\end{equation}
which has total mass one because $\sum_i(p_{\alpha i}+p_{\alpha,i+1})=2$.
Let $v=P_S\psi_\alpha/\sqrt{m_S}$ be the normalized restriction of
$\psi_\alpha$ to $S$, and let $J=P_i-P_{i+1}$. Then
$v^*Jv=(p_{\alpha i}-p_{\alpha,i+1})/m_S$ is the mass difference between
the two blocks of the pair, relative to the mass of the pair. Under
$\Prob^\sharp$ the denominator $m_S$ cancels, and cyclic invariance gives
the exact identity \eqref{eq:signed-pair}:
\begin{equation}\label{eq:intro-signed-pair}
 \E^\sharp_{\alpha,S}\bigl[v^*Jv\;\overline{A_m}\bigr]
 =\frac{1-z_m}{2}\,z_m^i\,\gamma_m.
\end{equation}
So $\gamma_m$ is controlled by the correlation, under $\Prob^\sharp$,
between the mass difference $v^*Jv$ on one pair and the global mode $A_m$.
If $v^*Jv$ were close to a quantity that is small and only weakly
correlated with $A_m$, then the left side would be small, and
\eqref{eq:intro-signed-pair} would give a bound on $\gamma_m$ with the
factor $d_m=|1-z_m|$.

\emph{What remains, and the challenge.} We must understand the law of the
direction $v$ under $\Prob^\sharp$, at the random energy $\lambda_\alpha$.
For this we need an exact formula for the eigenpair at a prescribed rank.

\emph{Step 3: an exact change of variables at the prescribed rank.} We
split off an independent GUE matrix $G$ on the pair $S$, so that
$H_{SS}=V+G$ as in \eqref{eq:pair-splitting}. This is possible because
both the diagonal blocks and the bond inside the pair have positive
variance. We then eliminate the remaining blocks by a Schur complement.
The eigenvalue equation $H\psi_\alpha=\lambda_\alpha\psi_\alpha$ becomes
the equation $A(E)v=0$ on the pair, where $A(E)$ is the Schur complement
of $H-E$ on $S$ and $E=\lambda_\alpha$; see \eqref{eq:schur-data}.
Finally we change variables from $G$ to the eigenpair $(E,v)$. This is
Theorem~\ref{thm:coarea}. The resulting density
\eqref{eq:complete-coarea} contains three factors:
$(i)$ a squared characteristic polynomial of the Schur complement on the
orthogonal complement of $v$; $(ii)$ an eigenvalue-counting indicator,
which uses Sylvester's law of inertia to fix the rank $\alpha$ even when
the exterior has eigenvalues near $E$; and $(iii)$ the Jacobian
$v^*\mathcal D(E)v$. The Jacobian is the derivative in $E$ of the zero
eigenvalue of $A(E)$, and \eqref{eq:clock-reconstruction} shows that it
equals $m_S^{-1}$. Hence
\[
 m_S\cdot v^*\mathcal Dv=1,
\]
and under the size-biased law \eqref{eq:intro-cancel} the Jacobian is
replaced by the constant $n/2$. This exact cancellation is the reason for
working with $\Prob^\sharp$. At fixed exterior data and fixed energy, the
direction $v$ then has an explicit density on the sphere, of the form
\eqref{eq:angular-measure} which is a Gaussian factor times the function
$\mathfrak Z(v)$ of factors $(i)$ and $(ii)$.

\emph{Step 4: comparison of the mass difference with a deterministic
center.} With this density we compare $v^*Jv$ with a center
$j_S(\lambda_\alpha)$ that depends only on the exterior and on the energy.
The center is defined in \eqref{eq:dyson-center} through the self-consistent
equation of the fresh GUE block. Theorem~\ref{thm:physical-angular} states
that for every bounded, phase-invariant Lipschitz function $\phi$ of the
eigenvector,
\begin{equation}\label{eq:intro-angular}
 \bigl|\E^\sharp_{\alpha,S}
   \bigl[(v^*Jv-j_S(\lambda_\alpha))\,\phi(\psi_\alpha)\bigr]\bigr|
 \le\frac{C(M+L_\phi)}W,
\end{equation}
where $M$ is a bound for $\phi$ and $L_\phi$ its Lipschitz constant. We
prove \eqref{eq:intro-angular} by Gaussian integration by parts on the
sphere (Lemma~\ref{lem:angular-stein} and
Proposition~\ref{prop:angular-comparison}). The integration by parts
produces three kinds of terms. The first contains the logarithmic gradient
of $\mathfrak Z(v)$. We compute this gradient exactly and bound it under
$\Prob^\sharp$ in Theorem~\ref{thm:determinant-score}. The bound uses the log-concavity of the eigenvalue law in a fixed
rank sector, which we derive from the Karlin--McGregor representation
of noncolliding Brownian motion and Pr\'ekopa's
theorem~\cite{Prekopa}. The second contains
second moments of $v$. These are given, at fixed energy, by the moment
theorem for a fresh Gaussian block (Theorem~\ref{thm:fresh-palm}) whose proof follows the localization and averaging method of
Benigni~\cite{Benigni}, which is built on the eigenvector moment
flow~\cite{BY} and its finite speed of propagation~\cite{BHY}. The
third comes from rewriting the eigenvector $\psi_\alpha=\mathcal L_Ev/
\norm{\mathcal L_Ev}$ in terms of $v$, and is bounded by
Lemma~\ref{lem:angular-ratio}. The regularity of the exterior needed in
these steps holds outside a polynomially small event
(Lemma~\ref{lem:physical-local-window}), which we obtain from the deformed-matrix local law and eigenvector
estimates of Cipolloni, Erd\H{o}s and Henheik~\cite{CEH} and the
block-Gaussian counting bounds of Peled, Schenker, Shamis and
Sodin~\cite{PSSS}.

We apply \eqref{eq:intro-angular} with $\phi=\overline{A_m}$, which is
bounded by one and has Lipschitz constant at most two. Combined with
\eqref{eq:intro-signed-pair} and the Cauchy--Schwarz inequality, using
$\E^\sharp|A_m|^2=\gamma_m$, this gives
\begin{equation}\label{eq:intro-mode}
 \frac{d_m}{2}\,\gamma_m
 \le\frac{C}W+\bigl(\E^\sharp_{\alpha,S}|j_S(\lambda_\alpha)|^2\bigr)^{1/2}
       \gamma_m^{1/2}.
\end{equation}

\emph{What remains, and the challenge.} We need a bound on the second
moment of the center at the random energy $\lambda_\alpha$. At a fixed
deterministic energy the center is of order $W^{-1}$
(Lemma~\ref{lem:smooth-center}), so its square is of order $W^{-2}$. We
cannot simply insert the random energy $\lambda_\alpha$, because it is
correlated with the center.

\emph{Step 5: the center at the random energy.} We prove in
Theorem~\ref{thm:physical-center} that
\begin{equation}\label{eq:intro-center}
 \E^\sharp_{\alpha,S}|j_S(\lambda_\alpha)|^2\le\frac{C}{\sqrt n\,W^{3/2}}.
\end{equation}
For the proof we split off independent GUE matrices $G_e$ on every edge
$e$ of the ring at once, as in \eqref{eq:simultaneous-splitting}, and we
isolate the scalar part $\xi_e=(2W)^{-1}\Tr G_e$ of each. By the
Hellmann--Feynman formula, the derivative of $\lambda_\alpha$ in $\xi_e$
is exactly the pair mass $p_{\alpha i}+p_{\alpha,i+1}$, that is, the
size-biasing weight. We therefore introduce the antiderivative in energy
\[
 F_e=\int_{e_\alpha}^{\lambda_\alpha}f_e(E)\,dE,
 \qquad f_e=|j_e^{\rm sm}|^2,
\]
where $e_\alpha$ is the classical location of the rank $\alpha$ and
$j_e^{\rm sm}$ is a smooth version of the center
(Lemma~\ref{lem:smooth-center}). Its derivative in $\xi_e$ is the
size-biased quantity we want, so Gaussian integration by parts in the
variables $\xi_e$ gives the exact identity
\eqref{eq:energy-primitive-identity}:
\[
 \E^\sharp_{\alpha,e}f_e(\lambda_\alpha)
 =\frac{1}{2}\,\E\int_{e_\alpha}^{\lambda_\alpha}D(E)\,dE,
 \qquad D(E)=\sigma^{-2}\sum_e\xi_ef_e(E).
\]
The right side is small for two reasons. First, the integration interval
is short. Theorem~\ref{thm:rank-location} shows that
$\norm{\lambda_\alpha-e_\alpha}_{L^p}\le C/(n\sqrt W)$ where its proof uses the counting estimate of~\cite{PSSS} for the last
fresh Gaussian block. Second, $D(E)$ is
small at every fixed energy: two further Gaussian integrations by parts
give its second moment exactly \eqref{eq:divergence-isometry}, and the
cross terms decay with the distance between edges
(Lemma~\ref{lem:smooth-center}), so $\norm{D(E)}_2\le C\sqrt n/W$.
Lemma~\ref{lem:random-energy} combines these two facts and gives
\eqref{eq:intro-center}.

\emph{Step 6: conclusion.} Inserting \eqref{eq:intro-center} into
\eqref{eq:intro-mode} gives
$d_m\gamma_m/2\le C/W+Cn^{-1/4}W^{-3/4}\gamma_m^{1/2}$. Young's
inequality absorbs half of the left side and yields
\eqref{eq:intro-mode-goal}, which is \eqref{eq:mode-bound}. By Step~1 this
proves Theorem~\ref{thm:que}.

\paragraph{Part II: construction of the critical law
(Theorem~\ref{thm:critical}).}

\emph{The goal.} Let $n/W\to\tau\in(0,\infty)$. We want to show that the
marked point process $\Xi_W$ of rescaled eigenvalues $u_\alpha$ and mass
profiles $\nu_\alpha$ near $E$, defined in \eqref{eq:unfolding}, converges
to a unique limit, 
and we want to describe this limit.

\emph{Step 1: tightness.} Tightness is the easy half. The counting
estimate of Theorem~\ref{thm:counting} bounds all factorial moments of the
number of rescaled eigenvalues in an interval, uniformly in $n$ and $W$. This estimate is a direct consequence of the deformed block-Gaussian
counting theorem of~\cite{PSSS}.
It follows that the sequence $\Xi_W$ is tight, so that it has limits along
subsequences. Every such subsequential limit is again a marked point
process of the same kind. Its points are simple, that is, no two of them
sit at the same rescaled eigenvalue. Its marks are probability measures on
the circle, that is, the limiting mass profiles still have total mass one,
and no mass is lost in the limit. The latter holds because the circle is
compact, so a weak limit of probability measures on it is again a
probability measure. The limit has unit spectral intensity. At the critical band width we obtain this by comparing the ring with shorter rings, on which the
local law of Erd\H{o}s and Riabov~\cite{ER} applies
(Theorem~\ref{thm:critical-density}).

\emph{What remains, and the challenge.} We must prove that all
subsequential limits coincide. The marks $\nu_\alpha$ are not determined
by the correlation functions of the eigenvalues, so the methods
of~\cite{Sh19,SS18,SS19} do not apply. We need one family of observables
that sees the eigenvalues and the marks at the same time, and whose limit
we can compute.

\emph{Step 2: one random entire function encodes the eigenvalues and the
marks.} Choose continuous functions $q_0=1,q_1,\ldots,q_k\ge0$ on $\T$ and
put $Q_a=\sum_jq_a(j/n)P_j$. Each $Q_a$ is block diagonal and multiplies
block $j$ by $q_a(j/n)$; in particular $Q_0=I$ and
$\ip{\psi_\alpha}{Q_a\psi_\alpha}=\nu_\alpha(q_a)$. Consider the ratio of
characteristic polynomials
\begin{equation}\label{eq:intro-determinant}
 D_W(\boldsymbol z)=
 \frac{\det(H-E+W^{-2}\sum_{a=0}^kz_aQ_a)}{\det(H-E)},
 \qquad \boldsymbol z=(z_0,\ldots,z_k)\in\C^{k+1}.
\end{equation}
It is an entire function of $\boldsymbol z$. On the first coordinate axis,
$D_W(z_0,0,\ldots,0)$ vanishes exactly when $E-W^{-2}z_0$ is an eigenvalue,
that is, at $z_\alpha=W^2(E-\lambda_\alpha)$. Since $N\approx\tau W^2$,
these zeros are the rescaled eigenvalues up to a constant factor. At a
simple zero, first-order perturbation theory shows that the ratio of
partial derivatives is $\ip{\psi_\alpha}{Q_a\psi_\alpha}$. In the limit we
therefore read off
\begin{equation}\label{eq:intro-extractor}
 u_\alpha=-\tau\rho(E)z_\alpha,
 \qquad
 \nu_\alpha(q_a)=
 \frac{\partial_{z_a}D(z_\alpha,0)}{\partial_{z_0}D(z_\alpha,0)}
\end{equation}
from the limit $D$ of $D_W$ (Theorem~\ref{thm:joint-recovery}). Hence it
suffices to prove that $D_W$ converges in law, as a random entire
function, to a unique limit $D$.

\emph{Step 3: bounded variables that determine the law of $D$.} We do not
work with $D_W$ directly, since its moments may be large. Instead we use
its inverse values on the imaginary cone, $D_W(i\boldsymbol y)^{-1}$ with
$\boldsymbol y\ge0$. Since $H-E$ is Hermitian and
$W^{-2}\sum_ay_aQ_a\succeq0$, these variables have modulus at most one.
Their mixed moments, together with complex conjugates, therefore
determine their joint law, and with it the law of $D$
(Theorem~\ref{thm:entire-law}). To compute these moments we write each
inverse determinant as a Gaussian integral. If $\zeta$ is a complex
Gaussian vector with covariance $W^{-2}\sum_ay_aQ_a$, independent of $H$,
then
\[
 D_W(i\boldsymbol y)^{-1}
 =\E_\zeta\exp\bigl\{-i\,\zeta^*(H-E)^{-1}\zeta\bigr\}.
\]
Using one independent probe vector $\zeta_j$ for each factor, the mixed
moments become expectations of
$\exp\{-i\sum_j\sigma_j\mathcal T_{jj}\}$, where
$\mathcal T=(\zeta_j^*(H-E)^{-1}\zeta_l)_{j,l}$ is the matrix of resolvent
quadratic forms of the probes; see \eqref{eq:determining-bounded-moments}.
Tightness of the entire functions follows from bounds on the trace term
and the Hilbert--Schmidt norm in $\log D_W$, which come from the same
probe representation.

\emph{What remains, and the challenge.} We must find the joint limit law
of the probe matrix $\mathcal T$ of the ring, jointly for several probes.
Two obstacles appear. First, the ring has no end from which we could grow
it block by block. Second, several probes need the joint law of a whole
matrix of boundary entries at the place where the ring is closed, not just
the law of each entry.

\emph{Step 4: cutting the ring.} We cut the ring at an independent buffer
of three consecutive blocks, which we call the collar. What remains is a
path. By Woodbury's identity and the determinant of a periodic closing,
the determinant ratio of the ring is expressed exactly through two
objects, \eqref{eq:complete-closing-determinant}: the resolvent of the
path with two matched leads attached at its ends, and a closing matrix
built from the collar. As a consequence, the probe matrix of the ring is
\begin{equation}\label{eq:intro-closing}
 K^{\rm ring}=\Re K+\eta\,V^*\mathcal C_WV,
 \qquad V^*V=-\Im K/\eta,
\end{equation}
see \eqref{eq:closed-source-matrix}. Here $K$ is the probe matrix of the
path with matched leads, $V$ collects the responses of the probes at the
two ends of the path, $\eta>0$ is the imaginary part of the lead
parameter, and $\mathcal C_W$ is the closing matrix of the collar. We
therefore need the joint law of $(K,V^*V)$ for the path (Steps~5--7), and
the law of $\mathcal C_W$ given the path (Step~8).

\emph{Step 5: the path with matched leads.} At each end of the path we attach the parameter $z_\partial$ of \eqref{eq:matched-parameters}, ($b$ is the bond variance; see \eqref{eq:model})
\[
 z_\partial=a^{-1}+ba,\qquad a=\frac{-E-i\sqrt{4-E^2}}{2}.
\]
It is chosen
so that the limiting trace equation of the path has the same fixed point
at every length. Physically, it plays the role of a semi-infinite lead
without reflection. Since $|a|=1$ and $\Im a<0$, its
imaginary part $\eta=\Im z_\partial$ is bounded below. Together with the
fresh Gaussian bonds, this gives inverse moments that are uniform in the
rest of the path (Lemma~\ref{lem:gaussian-absorption}), and the diagonal
traces of the resolvent concentrate at this fixed point with error
$O(W^{-1})$ (Theorem~\ref{thm:matched-traces}).

We then grow the path one layer at a time by the exact Schur recursion
\eqref{eq:two-ended-recursion}. If $X_r$ is the Schur complement of a path of $r$ blocks, adding one fresh block at each end
gives the exact recursion \eqref{eq:two-ended-recursion},
\[
 X_{r+2}=B^*(A-E-X_r)^{-1}B,
\]
where $A$ consists of fresh GUE blocks and $B$ of fresh Gaussian bonds. Integrating over each new Gaussian layer,
conditionally on the previous ones, gives the drift and the covariance of
the change. On the time scale $t=k/W$, the transmission eigenvalues of
the path then follow a diffusion. We identify its generator in
Theorem~\ref{thm:transmission-generator} to be the DMPK (Dorokhov--Mello--Pereyra--Kumar) generator of the
unitary class, with the constraint $0\preceq T\preceq I$ retained. A DMPK limit for transfer matrices of a weakly disordered strip with
a fixed number of channels was obtained by~\cite{BBDR}, and here we derive it directly from the Gaussian
blocks, with the number of channels growing. The
generator alone does not determine the law, because the entrance at time
zero has infinite conductance. We identify the law through products of
determinant ratios that have a positive finite-dimensional dual process,
and we exclude mass hidden at zero transmission
(Theorem~\ref{thm:transmission-law}).

\emph{What remains, and the challenge.} The transmission describes how
the path transmits, but not where inside the path the response to each
boundary input is located. The probe matrix needs exactly this
information, jointly for all the weights $Q_a$.

\emph{Step 6: the sources.} For a weight $q$ we track the interior mass
matrix $Y^q$ (associated to any path of blocks), which records how much of the response to each boundary
input lies in each block, weighted by $q$.  We call these matrices
sources. These are given by the formula: 
\[
 Y^q_j=\mathcal U_j^*\diag(q_1I_W,\ldots,q_jI_W)\,\mathcal U_j
\]
where $\mathcal U_j$ denotes the map of an incoming boundary vector to a path of $j$ blocks. Theorem~\ref{thm:full-source-generator} identifies the generator
of the sources jointly with the transmission. It involves the current
$J=\mathbb S^* P \mathbb S-P$ (with $\mathbb S$ being the scattering matrix of the path and $P$ projects onto the left boundary outputs), the difference of two orthogonal projections that encode the
transmission eigenvalues. This generator acts on trace observables; it
does not by itself define a probability law. We construct such a law by
finite systems of stochastic differential equations
\eqref{eq:finite-source-sde} in a finite-dimensional space. Their
solutions preserve the pair of projections
(Proposition~\ref{prop:finite-gram-law}), and the sources are Gram
matrices of the solution, hence positive semidefinite. A finite system
keeps only finitely many transmission eigenvalues. We remove this cutoff
with a coupling bound that does not depend on the dimension
(Section~\ref{sec:positive-source-law}); the result is
Theorem~\ref{thm:source-process}.

\emph{Step 7: complex probe matrices.} The sources of Step~6 give only the
imaginary part of the resolvent quadratic forms of the probes, but the
probe matrix also needs their real part. Let $G$ be the resolvent of the
path with its two matched leads, let $\mathcal J$ collect the probe
vectors, and let $V_{\rm pr}$ be their boundary responses. We add one
equation for the full complex matrix, driven by the same noise $dA$ as
the Gram process:
\[
 \mathcal K=\mathcal J^*G\mathcal J,\qquad
 d\mathcal K=-2i\eta\,V_{\rm pr}^*\,dA\,V_{\rm pr},\qquad
 \Im\mathcal K=-\eta\,V_{\rm pr}^*V_{\rm pr}.
\]
The last identity is the imaginary-part identity of the resolvent. It ties
$\mathcal K$ to the sources and fixes all normalizations
(Theorem~\ref{thm:complex-source-law}). Finally we replace the global
Gaussian probes $\zeta_a$, which are spread over the whole path, by many
independent probes on single blocks, born at random times with density
$q_a$. A Lindeberg replacement with error uniform in $W$ justifies this
step (Lemma~\ref{lem:t22-global-Gaussian-source}).

\emph{What remains, and the challenge.} We must identify the closing
matrix $\mathcal C_W$ of the collar, jointly in all directions needed by
the probes. Its Schur complement can be large, and the laws of single
entries do not determine the joint law.

\emph{Step 8: the closing matrix.} The closing matrix is a matrix-valued
Herglotz function of the collar. Its imaginary part is positive in the
upper half-plane. As a function of the microscopic spectral variable $z$
near $E$, it has the form
\[
 \mathcal C_W(z)=C_W+\sum_k\frac{u_ku_k^*}{\lambda_k-z},
 \qquad C_W=C_W^*.
\]
The poles $\lambda_k$ are the rescaled eigenvalues of the collar, and the
residues are built from the overlaps $u_k$ of their eigenvectors with the
boundary directions. By Lemma~\ref{lem:t22-Herglotz-tail}, three facts fix
the limit. First, the poles converge to a sine process
(Theorem~\ref{thm:external-sine}, from the exact kernel of the Gaussian
ensemble with external source~\cite{JohanssonKernel}), and the overlaps
converge to independent Gaussian vectors
(Theorem~\ref{thm:relative-mixed-vectors}, which adapts the eigenvector
moment flow~\cite{BY,Benigni}, its finite speed of
propagation~\cite{BHY}, and the comparison of Benigni and
Lopatto~\cite{BL}). Second, no mass escapes to infinity, i.e.,
$\mathbb E\,\Im\mathcal C_W(iY)\to I_r$ for every $Y>0$
(Theorem~\ref{thm:t22-exact-boundary-mean}). Third, single boundary
entries are Cauchy, which fixes the real constant $C_W$
(Corollary~\ref{cor:t22-physical-one-channel-Cauchy}). This extends the
scalar argument of Aizenman and Warzel~\cite{AW} to matrices and gives
Theorem~\ref{thm:matrix-closing}: the closing matrix converges to
\[
 \mathsf C_r(0)=\frac1\pi\lim_{R\to\infty}
 \sum_{|x_j|\le R}\frac{g_jg_j^*}{x_j},
\]
where $\{x_j\}$ is one sine process and the $g_j$ are independent
Gaussian vectors; see \eqref{eq:canonical-cauchy}. It does not depend
on the path. All probes share this one matrix, which keeps their
correlations.

\emph{Step 9: assembling the law.} Inserting the limit of the closing
matrix into \eqref{eq:intro-closing}, we obtain the limiting probe matrix
$\mathcal T$ in \eqref{eq:limiting-closed-source}. Its law is determined
by the source process of Steps~5--7 and the closing matrix of Step~8. The
probes have no mass on the three collar blocks. We restore this mass at a
cost that vanishes, because the counting estimates bound the resolvent (\cite{PSSS}),
and because the limiting marks have no atoms
(Theorem~\ref{thm:small-critical}). This gives Theorem~\ref{thm:entire-law}.
Step~2 then recovers all rescaled eigenvalues and all their profiles
jointly (Theorem~\ref{thm:joint-recovery}), which proves
Theorem~\ref{thm:critical}.

\paragraph{Part III: the profile variance and the two endpoints
(Theorem~\ref{thm:endpoints}).}

\emph{The goals.} We want to prove three statements about the critical law
$\Xi_{\tau,b,E}$: the profile variance $V_{\tau,b,E}(\phi)$ in
\eqref{eq:variance-limit} is strictly positive and does not depend on the
spectral weight $w$; it tends to zero as $\tau\downarrow0$; and the
spectrum tends to the sine process as $\tau\downarrow0$ and to the Poisson
process as $\tau\to\infty$, as in \eqref{eq:spectral-endpoints}.

\emph{The challenge.} The law $\Xi_{\tau,b,E}$ is given through the
construction of Part~II, not by explicit formulas. Its properties must be
read from that construction, or proved directly for the matrix with
estimates that are uniform in $W$. There is also a technical point. The
variance statistic is a sum over all eigenvalues in the support of $w$,
whose number is not bounded, so weak convergence alone does not give
convergence of its expectation.

\emph{Step 1: passing the variance to the limit, and independence of $w$.}
The counting estimates of Theorem~\ref{thm:counting} (from~\cite{PSSS}) bound all moments of
the number of points in a window. This gives the uniform integrability
needed to pass the variance statistic to the limit. The limit law is
invariant under spectral translations (Proposition~\ref{prop:stationarity}).
Hence the measure $u\mapsto\E\sum w(u)|\nu(\phi)-\bar\phi|^2$ is a
constant multiple of Lebesgue measure, and the variance does not depend on
$w$ once $\int w=1$.

\emph{Step 2: strict positivity of the variance.} The limit law is also
invariant under rotations of the circle, so the variance is diagonal in
Fourier modes:
\[
 V_{\tau,b,E}(\phi)=\sum_{m\ne0}|\widehat\phi_m|^2v_m,
 \qquad
 v_m=\E\sum_{(u,\nu)}w(u)\,|\nu(e^{2\pi im\cdot})|^2 .
\]
Every nonconstant $\phi$ has a nonzero coefficient $\widehat\phi_m$, so it
suffices to prove $v_m>0$ for every $m\ne0$, i.e.\ that the profiles
fluctuate in every mode. We detect this by scattering. We express the
one-point statistics of the ring through the path of Part~II with matched
leads at its ends (Theorem~\ref{thm:t22-physical-lead-profile}). A wave
entering from one lead is partly reflected and partly transmitted; we
observe the transmission through bounded transforms
$F_{L,R}=\det(I-yT_{L,R})/\det(I+xT_{L,R})$ of the matrices
$T_L=I-S_{11}^*S_{11}$ and $T_R=I-S_{22}^*S_{22}$. A weak absorber is a
small imaginary potential $i\kappa q(s)$ inside the path, which removes
intensity at the position $s$ at a rate proportional to $q(s)$. To first
order in $\kappa$, it changes $F_L+F_R$ by an integral of $q$ against the
intensity profile of the waves crossing the path
(Lemma~\ref{lem:t22-pass-derivative}). For the signed weight
$q_m(s)=\cos(2\pi ms/\tau)$, defined by linearity from its positive and
negative parts, this response vanishes for a flat profile, since $q_m$
has mean zero, and Theorem~\ref{thm:t22-response-variance} bounds its
square by $C\,v_m$. It remains to show that it is not zero. In the limit
it is an explicit functional $J_{\tau,m}$ of the transmission diffusion
absorbed at its corner. We show $J_{\tau,m}>0$ by conditioning on the
absorption time. Uder this conditioning, the diffusion is then a bridge whose occupation is concentrated near its initial end, where $q_m\ge1/2$
(Lemmas~\ref{lem:t22-lifetime-bridge}
and~\ref{lem:t22-bridge-occupation},
Proposition~\ref{prop:t22-nonzero-response}).

\emph{Step 3: the variance vanishes as $\tau\downarrow0$.} Here we use
Part~I. The mode bound \eqref{eq:mode-bound} holds at every fixed rank,
while $v_m$ is a statistic in a fixed microscopic energy window. To pass
from ranks to energies, we write $H=H_0+tI$ with $t=N^{-1}\Tr H$. The
scalar $t$ is a Gaussian variable independent of $H_0$, its standard
deviation is about $\sqrt n$ eigenvalue spacings, and it shifts all
eigenvalues without changing any eigenvector. Averaging over $t$ therefore
turns the window statistic into an average of prescribed-rank statistics
(Theorem~\ref{thm:window-transfer}). With \eqref{eq:mode-bound} this gives
$v_m\le C[\tau/|m|+\tau^{3/2}/m^2]^{1-\eta}$, and hence
$V_{\tau,b,E}(\phi)\le C\tau^{1-\eta}$ (Theorem~\ref{thm:small-critical}).

\emph{Step 4: the sine endpoint.} We add a scalar multiple of the identity
on one adjacent pair $S$ and set $H(t)=H+t(2W\rho)^{-1}P_S$. Two facts make
this flow useful. First, the values of $t$ at which $E$ is an eigenvalue of
$H(t)$ are, up to the factor $-2W\rho$, the eigenvalues of the Schur
complement on $S$, which is a fresh GUE matrix plus an independent
deformation. Its local eigenvalues converge to the sine process
(Lemma~\ref{lem:t22-crossing-sine}, based on Theorem~\ref{thm:external-sine}), i.e.\ on the exact finite kernel of~\cite{JohanssonKernel}.
Second, by the Hellmann--Feynman formula each rescaled eigenvalue moves with
speed:
$$\dot u_\alpha=\frac{n}{2}\norm{P_S\psi_\alpha}^2.$$
If this speed is close
to one, then each eigenvalue crosses $E$ at the time $t\approx-u_\alpha$,
so the eigenvalues and the crossing values are paired within a small
error. Speed close to one means that the mass of $\psi_\alpha$ on $S$ is
close to its average $2/n$. We prove this local mass flatness, with error
$C\sqrt\tau$, in Theorem~\ref{thm:t22-local-mass-flat}, and we control it
along the whole flow in Lemma~\ref{lem:t22-flow-admission}. The pairing
argument then gives Theorem~\ref{thm:t22-sine-endpoint}, with error
$C\tau^{1/8}$.

\emph{Step 5: the Poisson endpoint.} We compare the original ring with $m$
independent smaller rings of the same total dimension. To build the
comparison, we place $m$ collars around the ring. Each collar block has two
Gaussian bonds, which can be paired in two ways: one pairing gives the
original ring, and the other gives $m$ independent rings. All entries keep
their original laws (Subsection on same-volume rewiring in
Section~\ref{sec:t22-positive-log-poisson}).

The eigenvalue statistic $\Tr f(\Lambda(H-E))$, with $\Lambda=N\rho(E)$, is
not directly a determinant. We therefore approximate it by the positive
logarithmic statistic \eqref{eq:positive-log-statistic}, which is a
smoothed log-determinant with a small absorber $i\delta P_C$ on the blocks
away from the collars (Lemma~\ref{lem:t22-finite-log-recovery}). Changing
the pairing at one collar changes this log-determinant by an amount
controlled by the change of the scattering matrix of the rest of the
system (Lemma~\ref{lem:t22-unitary-log-comparison}). That change is in
turn bounded by the response of the absorber, measured by a fractional
trace of the clock of the blocks in $C$, by
\eqref{eq:remote-scattering-source}. Since these blocks are at distance at
least $R$ from every collar, this response decays like $e^{-cR}$
(Theorem~\ref{thm:t22-full-remote-Schatten}). Summing over the collars
gives \eqref{eq:same-volume-log-error}. For the independent rings, the
Laplace functional of the sum of their statistics factorizes. Each ring
carries a fraction $a_i$ of the volume and contributes few eigenvalues near
$E$, so the product tends to the Poisson Laplace functional, with error
$C\sum_ia_i^2$ by \eqref{eq:independent-ring-laplace}. We take
$m\asymp\sqrt\tau$ rings and buffers of length $R=A\log(2+\tau)$. We first
let $W\to\infty$ for this fixed partition, then $\tau\to\infty$, and
finally $\delta\downarrow0$. All errors vanish, which proves
Theorem~\ref{thm:t22-Poisson-endpoint}.

\section{Notation and inputs}\label{sec:inputs}

The goal of this section is to fix notation and to collect, in theorem
form, the results from earlier works that our proofs use, together with
the direct proofs and comparison lemmas needed at their application
scope. Imported results are identified by their citations; the complex
fixed-variance kernel theorem and the relative profile and admissibility
bridges have explicit proofs here or in the indicated later sections. The
specializations of these results to a newly exposed Gaussian block are
proved later and are kept separate from the arguments that compare that
block with the original ring.

The inputs enter as follows. Lemma~\ref{lem:gaussian-identities} gives the
Gaussian moment and integration-by-parts identities used throughout.
Theorem~\ref{thm:external-count}, from~\cite{PSSS}, gives every count of
eigenvalues in short intervals; it is used in Section~\ref{sec:counting}.
Theorems~\ref{thm:external-local} and~\ref{thm:external-eth},
from~\cite{CEH}, give local laws for a Gaussian matrix with a large
deterministic deformation; they localize eigenvalue paths. Theorem~\ref{thm:external-sine}, proved directly from the complex
Gaussian kernel, identifies the local spectrum of a fresh Gaussian block.
Theorems~\ref{thm:external-emf} and~\ref{thm:external-finite-speed},
from~\cite{Benigni,BHY}, control the conditional moments of eigenvectors;
they are adapted in Sections~\ref{sec:fresh-gaussian}
and~\ref{sec:relative-vectors}. Theorem~\ref{thm:external-band-local},
from~\cite{ER}, gives a local law on a shorter ring and is used only for
the mean density in Section~\ref{sec:counting}.

\subsection{Matrices, norms and Gaussian normalization}\label{subsec:notation}

In this subsection, we fix the conventions needed to compare estimates from different parts
of the proof. Lemma~\ref{lem:gaussian-identities} records the Gaussian
identities in our normalization. We specify matrix norms and dimensions,
explain the range of moment estimates, and then give the Gaussian laws
and their basic calculations.

\paragraph{Matrices and norms.} For a matrix $A$ on $\C^d$, let
\[
 \av A=d^{-1}\Tr A,\qquad
 \norm A_{\HS}=(\Tr A^*A)^{1/2},\qquad
 \Re A=\frac{A+A^*}{2},\qquad \Im A=\frac{A-A^*}{2i}.
\]
The notation $\norm A$ means the operator norm. For Hermitian matrices,
$A\preceq B$ means $v^*Av\le v^*Bv$ for every vector $v$.
The full physical ring has dimension $N=nW$. Auxiliary dimensions
$N_0$ or $d$ are defined locally: in particular a fresh pair has
dimension $2W$, a short ring has dimension $mW$, and the finite source
SDE has an independently chosen projection rank.
 We use the unnormalized trace $\Tr$ unless angle brackets are displayed.
The Schatten expression is $\norm A_p^p=\Tr|A|^p$, also when $0<p<1$; in
that range it is not a norm. The notation $\norm X_{L^p}=(\E|X|^p)^{1/p}$
always refers to the specified probability law and to a fixed order $p$,
unless a growing range of $p$ is stated.

\paragraph{Convention on moments and widths.}\label{par:moment-width}
Many auxiliary estimates involve inverse factors with absorption, their
derivatives, or propagated clocks. We state each such estimate only for
$W\ge W_0(p,\ell)$, where $p$ is the fixed moment order and $\ell$ is the
longest fixed word length used in it. The constants may also depend on
the compact parameter sets and the time horizon stated with the estimate.
When we say that such an estimate is ``uniform in $W$'', we mean uniform
for all $W$ above this threshold. Moments of growing order appear only in
explicitly stated ranges $p\le aW$, and each product or derivative
reduces $a$ by a fixed factor. We never claim inverse moments of
arbitrary order at a fixed small width. Exact identities and estimates
for bounded statistics are not affected by this convention. They hold in
the finite-dimensional generality in which they are stated; in
particular, Theorem~\ref{thm:que} holds for all $W$.

 \paragraph{Gaussian laws.} Let $\GUE_d(t)$ denote the Gaussian law on
$d\times d$ Hermitian matrices with density
\begin{equation}\label{eq:gue-normalizer}
 Z_d(t)^{-1}e^{-\Tr A^2/(2t)}\,dA,
 \qquad Z_d(t)=(2\pi t)^{d/2}(\pi t)^{d(d-1)/2},
\end{equation}
where $dA$ is Lebesgue measure on the real diagonal and the complex upper
triangle. Under $\GUE_d(t)$ a diagonal entry has variance $t$ and an
off-diagonal entry has $\E|A_{ij}|^2=t$. The normalized GUE is
$\GUE_d(1/d)$. A \emph{standard circular Gaussian vector} $g\in\C^d$ has
density $\pi^{-d}e^{-\norm g^2}$, so $\E gg^*=I$ and $\E gg^T=0$.

The next lemma records the Gaussian identities used throughout the paper.

\begin{lemma}[Gaussian moments and integration by parts]\label{lem:gaussian-identities}
Let $g\in\C^d$ be a standard circular Gaussian vector and let $A,B$ be
deterministic $d\times d$ matrices. Then
\begin{equation}\label{eq:gaussian-wick}
 \E(g^*Ag)=\Tr A,\qquad
 \E[(g^*Ag)(g^*Bg)]=\Tr A\Tr B+\Tr(AB).
\end{equation}
For a real Gaussian coordinate $X$ of variance $t$,
$\E[Xf(X)]=t\E f'(X)$. For a circular Gaussian coordinate $Z$ of variance
$t$,
\begin{equation}\label{eq:gaussian-ibp}
 \E[Zf]=t\E\partial_{\bar Z}f,
 \qquad \E[\bar Zf]=t\E\partial_Zf.
\end{equation}
These identities hold whenever $f$ is locally weakly differentiable and
the following functions are integrable with respect to the relevant
Gaussian measure: $f$ itself, its weak derivatives, and the products of
$f$ with the coordinates that appear in the identities. For a single real
variable it is enough that $f$ is locally absolutely continuous. For a
complex variable $z=x+iy$, the derivatives are the Wirtinger derivatives
$\partial_z=\frac{1}{2}(\partial_x-i\partial_y)$ and
$\partial_{\bar z}=\frac{1}{2}(\partial_x+i\partial_y)$, built from the two
real weak derivatives. In the multivariate weak-derivative formulation, all coordinates other than the differentiated coordinate have positive variance. The differentiated real coordinate may have variance zero, in which case both sides vanish.
\end{lemma}

The first identity in \eqref{eq:gaussian-wick} is the mean of a quadratic
form, and the second is its second moment; the term $\Tr(AB)$ is the only
correction to independence. The identities \eqref{eq:gaussian-ibp} are
used for every Gaussian integration by parts in the paper.

\begin{proof}
We expand the entries. For \eqref{eq:gaussian-wick}, write
$g^*Ag=\sum_{i,j}\bar g_iA_{ij}g_j$ and use the fourth moment of
the circular coordinates,
\[
 \E\bar g_i g_j\bar g_k g_\ell=\delta_{ij}\delta_{k\ell}
 +\delta_{i\ell}\delta_{kj}.
\]
The term $\delta_{ij}\delta_{k\ell}$ gives $\Tr A\Tr B$, and the term
$\delta_{i\ell}\delta_{kj}$ gives $\Tr(AB)$. The first identity follows
in the same way from $\E\bar g_ig_j=\delta_{ij}$. The real identity
$\E[Xf(X)]=t\E f'(X)$ follows by integrating by parts against the density
$(2\pi t)^{-1/2}e^{-x^2/(2t)}$. Writing $Z=X_1+iX_2$ with independent real
$X_1,X_2$ of variance $t/2$ and applying the real identity to each part
gives \eqref{eq:gaussian-ibp}, since
$\partial_{\bar Z}=\frac{1}{2}(\partial_{X_1}+i\partial_{X_2})$.
To extend the identities to such $f$, we first cut $f$ off at a large
radius. Let $\chi_R$ be a smooth function that equals one on the ball of
radius $R$, vanishes outside the ball of radius $2R$, and satisfies
$|\nabla\chi_R|\le C/R$. The product $\chi_Rf$ has compact support, and on
that support the Gaussian density and its derivatives are bounded. So we
can approximate $\chi_Rf$ by smooth functions in the local Sobolev norm,
and integration by parts holds for $\chi_Rf$. It remains to let
$R\to\infty$. The derivative of the cutoff produces one extra term,
\[
 \Bigl|\int f\,\nabla\chi_R\,d\gamma\Bigr|
 \le\frac CR\int_{|x|>R}|f|\,d\gamma\longrightarrow0 .
\]
All other terms converge to the corresponding terms for $f$, because
$f$, its weak derivatives, and the coordinate products are integrable
with respect to $\gamma$. \end{proof}

\begin{remark}[Gradient inequalities]\label{rem:poincare}
We also use the Gaussian Poincar\'e inequality and its gradient inequality
for a fixed moment order $p$, with the same covariance of real coordinates.
Matrix derivatives are always taken before a weighted resolvent is
estimated.
\end{remark}

\paragraph{Stochastic domination.} We write $X\prec Y$ if for every fixed
$\epsilon,D>0$,
\begin{equation}\label{eq:stochastic-domination}
 \Prob(X>N^\epsilon Y)\le N^{-D}
\end{equation}
for all sufficiently large $N$, uniformly over the deterministic
parameters stated with the estimate.

\begin{remark}[The fixed loss in $\prec$]\label{rem:domination-loss}
In \eqref{eq:stochastic-domination} the exponent $\epsilon$ is fixed before
$N\to\infty$. A stochastic-domination estimate therefore does not imply an
estimate with a constant factor. This is one reason why the proof of
Theorem~\ref{thm:que} uses the exact formulas of Section~\ref{sec:coarea}.
\end{remark}

\subsection{Deformed block Gaussian matrix counting}\label{subsec:external-count}

We need bounds on how many eigenvalues can lie in a short interval,
even after the surrounding matrix has been fixed. Theorem~\ref{thm:external-count} gives these bounds from the independent
Gaussian diagonal blocks. We state the counting theorem, derive a version
that retains a spatial projection, and use it to control logarithms of
determinants with full or partial absorption.

\begin{theorem}[Peled--Schenker--Shamis--Sodin]\label{thm:external-count}
Let $A_0$ be a deterministic Hermitian matrix of dimension
$d=d_1+\cdots+d_k$, and let $G_j\sim\GUE_{d_j}(1/d_j)$ be independent.
Set $A=A_0+\bigoplus_{j=1}^kG_j$. There is an absolute $C$ such that
for every real interval $I$ and every integer $r\ge1$,
\begin{equation}\label{eq:external-count}
 \E\prod_{j=0}^{r-1}\bigl(\Tr\one_I(A)-j\bigr)
       \le(Cd|I|)^r.
\end{equation}
\end{theorem}

The left side of \eqref{eq:external-count} is the $r$-th factorial moment of the number $\Tr\one_I(A)$ of eigenvalues of $A$ in $I$. The bound is of
order $(d|I|)^r$, where $d|I|$ is the number of eigenvalues expected in $I$
for a density of order $d$, and it holds for every interval, however short.
This is the unitary case of \cite[Theorems 2--3]{PSSS} with the
normalization \eqref{eq:gue-normalizer}.

\begin{remark}[The deformation]\label{rem:count-deformation}
The deformation $A_0$ need not have bounded norm. In our ring, conditioning
on all bonds turns $H$ into a deterministic matrix plus independent GUE
diagonal blocks, so the hypotheses of Theorem~\ref{thm:external-count} are
preserved after this conditioning.
\end{remark}

\begin{lemma}[Projected Wegner bound with an independent absorber]
\label{lem:t26-projected-wegner}
Let $H=D+\bigoplus_{j=1}^k A_j$ act on $\C^{kW}$, where the blocks
$A_j\sim\GUE_W(c/W)$ are independent and $c>0$. Let $D$ be Hermitian and
$\Gamma$ positive semidefinite, jointly independent of the blocks; they
may be random and of any norm. For real $e$, $\eta>0$ and a union $B$ of
blocks, put $G=(H-e+i\Gamma+i\eta I)^{-1}$. Then
\begin{equation}\label{eq:t26-projected-wegner}
 \E[-\Im\Tr P_BG]\le\pi K_c W|B|.
\end{equation}
By Fatou's lemma the bound also holds for the boundary values
$\eta\downarrow0$. For a single block $j$, it suffices that $A_j$ is jointly independent of $(D,\Gamma,(A_i)_{i\ne j})$.
\end{lemma}
\begin{proof}
We prove the bound for one block $j$ and then sum over $j\in B$. We
condition on $(D,\Gamma)$ and on all blocks except $A_j$. The Schur
complement of $M=H-e+i\Gamma+i\eta I$ onto block $j$ then has the form
$A_j+F+iY+i\eta I$, where $F=F^*$ and $Y\succeq0$ are fixed by the
conditioning. To see that its imaginary part is at least $\eta I$, put
$w=(v,-M_{22}^{-1}M_{21}v)$; then $w^*Mw=v^*M_{\rm Schur}v$, and
$\Im w^*Mw\ge\eta|w|^2\ge\eta|v|^2$.

For each value of $A_j$, the function
\[
 \zeta\longmapsto-\Im\Tr(A_j+F+\zeta Y+i\eta I)^{-1}
\]
is bounded, nonnegative and harmonic in the upper half-plane. The
Poisson formula at $\zeta=i$ therefore writes it as an average over a
real background:
\[
 -\Im\Tr(A_j+F+iY+i\eta I)^{-1}
 =\int_\R\frac{-\Im\Tr(A_j+F+tY+i\eta I)^{-1}}{\pi(1+t^2)}\,dt.
\]
This is \cite[Lemma 3.1]{PSSS}, with the sign conjugated. For real $t$,
the one-block case of Theorem~\ref{thm:external-count}, after scaling by
$\sqrt c$, bounds the mean spectral density of $A_j+F+tY$ by $K_cW$,
uniformly in $F+tY$. Since $\int\eta\,dx/(x^2+\eta^2)=\pi$, the
conditional expectation of each integrand is at most $\pi K_cW$. Tonelli's
theorem and removing the conditioning give the one-block bound. Summing
over $j\in B$ gives \eqref{eq:t26-projected-wegner}. Since $-\Im G\succeq0$,
Fatou's lemma applies as $\eta\downarrow0$.
\end{proof}

\begin{lemma}[Uniform moments of the singular log statistic]
\label{lem:t26-singular-log-moments}
Let $H=D+\bigoplus_{j=1}^k A_j$ be as in
Lemma~\ref{lem:t26-projected-wegner}, with $d=kW\le N$. For every fixed
$r>0$ and $K<\infty$,
\[
 \left\|\frac{1}{2}\sum_{\lambda\in\operatorname{Spec}H}
                 \log\left(1+\frac{a^2}{(\lambda-e)^2}\right)
                  \right\|_{L^r}\le C_{r,K,c},
\]
uniformly in real $e$, in $D$, in $W,k,N$, and in $0\le aN\le K$. The
bound holds conditionally on any independent Hermitian $D$, hence also
after averaging over such $D$.
\end{lemma}
\begin{proof}
We rescale the spectrum, putting $u=N(\lambda-e)$ and
$\mathcal N(I)=\#\{u\in I\}$. Since $d/N\le1$,
Theorem~\ref{thm:external-count} bounds the factorial moments,
$\E(\mathcal N(I))_j\le(K_c|I|)^j$. For an integer $r\ge1$ we expand
$\mathcal N^r$ into falling factorials with Stirling numbers, and obtain
\[
 \|\mathcal N([-v,v])\|_r\le C_r v^{1/r}\quad(0<v\le1),
 \qquad \|\mathcal N([-v,v])\|_r\le C_r v\quad(v\ge1).
\]
The first-moment bound shows that almost surely no point lies at $u=0$.

We split the sum into dyadic shells. Since $\$a/|\lambda-e|=aN/|u|$, the
log summand is at most $C_K(1+j)$ on $2^{-j-1}<|u|\le2^{-j}$, and at most
$C_K2^{-2j}$ on $2^j<|u|\le2^{j+1}$. By Minkowski's inequality, the
contributions of the small and the large shells with index $j\ge m$ are
bounded by
\[
 C_{r,K}\sum_{j\ge m}(1+j)2^{-j/r},\qquad
 C_{r,K}\sum_{j\ge m}2^{-j}.
\]
Both series converge, and their tails vanish as $m\to\infty$, uniformly in
all the parameters. This also lets us remove the cutoffs near the pole
and at infinite energy. For noninteger $r$ we use the next integer
moment. The only probabilistic input is the uniform count theorem for
deformed blocks.
\end{proof}

\begin{lemma}[Partial absorbers and the Cauchy formula for the log]
\label{lem:t26-passive-log}
Let $H$ be as in Lemma~\ref{lem:t26-singular-log-moments}. Let $Y\succeq0$
be independent of the Gaussian blocks, jointly with $D$, and let
$0\preceq Q\preceq I$ be fixed. Put
\[
 B_a(H+iY;Q)=\log\left|
       \frac{\det(H-e+iY+iaQ)}{\det(H-e+iY)}\right|.
\]
For every fixed real $e$ and $a\ge0$, almost surely,
\begin{equation}\label{eq:t26-passive-log-formula}
 B_a(H+iY;Q)=\int_\R B_a(H+tY;Q)\frac{dt}{\pi(1+t^2)}\ge0.
\end{equation}
For every Hermitian $A$ with $A-e$ invertible,
\begin{equation}\label{eq:t26-partial-log}
 0\le B_a(A;Q)\le B_a(A;I)
 =\frac{1}{2}\sum_\lambda\log(1+a^2/(\lambda-e)^2).
\end{equation}
Consequently $\|B_a(H+iY;Q)\|_r\le C_{r,K,c}$ when $aN\le K$,
uniformly in $D,Y,Q$ and $d\le N$.
\end{lemma}
\begin{proof}
 We first add a regulator $i\eta I$ with $\eta>0$. Along
any increasing absorber $\Gamma(s)$,
\[
 \frac{d}{ds}\log|\det(A-e+i\Gamma(s)+i\eta I)|
 =-\Im\Tr[\Gamma'(s)(A-e+i\Gamma(s)+i\eta I)^{-1}]\ge0,
\]
because $-\Im G=G^*(\Gamma(s)+\eta I)G\succeq0$ and $\Gamma'(s)\succeq0$.
We apply this from $0$ to $aQ$ and then from $aQ$ to $aI$, and let
$\eta\downarrow0$. This proves \eqref{eq:t26-partial-log}. We never need
$Q$ and $A$ to be simultaneously diagonalizable. Fix $0\le s\le a$. The function
\[
 \zeta\longmapsto-\Im\Tr[Q(H-e+\zeta Y+isQ+i\eta I)^{-1}]
\]
is harmonic, nonnegative and bounded by $d/\eta$, so the Poisson formula
applies. We integrate first in $s$ and then in the Cauchy variable $t$;
Tonelli's theorem holds by nonnegativity. The fundamental theorem of
calculus then gives \eqref{eq:t26-passive-log-formula} with $i\eta I$
added in all four determinants. For $0<\eta\le a$, the increment with real
background is bounded by its full-absorber version, in which each
eigenvalue of $H+tY$ contributes at most
\[
 \frac{1}{2}\log\frac{(\lambda-e)^2+(\eta+a)^2}{(\lambda-e)^2+\eta^2}
 \le\frac{1}{2}\log\left(1+\frac{4a^2}{(\lambda-e)^2}\right).
\]
Conditionally on $(D,Y)$ and the bonds, every fixed moment of this sum is
bounded uniformly in $t$ by Lemma~\ref{lem:t26-singular-log-moments}. The
Cauchy measure has mass one, so dominated convergence in the joint
Gaussian and Cauchy measure passes to the integral without regulator, in
$L^r$. Almost surely, the determinant with real background is nonzero for
each fixed $t$, by the count theorem, and $H-e$ is invertible; by
positivity, $H-e+iY$ and $H-e+iY+iaQ$ are invertible too. So the left
side converges pointwise, which proves the identity; the case $a=0$ is
trivial. Finally, Jensen's inequality for the Cauchy probability measure
and \eqref{eq:t26-partial-log} give the moment bound.
\end{proof}

 \subsection{A local law for a large deformation}\label{subsec:external-local}

We compare the resolvent of a Gaussian matrix with a deterministic
approximation, even when the added deterministic matrix is large.
Theorems~\ref{thm:external-local} and~\ref{thm:external-eth} control the
resolvent and eigenvectors in a regular spectral window. We define the
approximation through its matrix Dyson equation, state the estimates,
and explain the regular-window consequence used below.

Let $F=F^*\in\C^{\Naux \times \Naux }$ be deterministic, $s>0$, and $A=F+G$ with
$G\sim\GUE_\Naux (s/\Naux )$. For $z\in\C\setminus\R$, let $M(z)$ be the solution of
\begin{equation}\label{eq:mde}
 M(z)=(F-z-s\av{M(z)})^{-1},\qquad
 \Im z\,\Im M(z)\succ0.
\end{equation}
Equation \eqref{eq:mde} is the matrix Dyson equation of $A$; its solution
$M(z)$ is the deterministic approximation of the resolvent of $A$.
Write $m=\av M$, $\rho_s=\pi^{-1}|\Im m|$, and
$\Gamma(z)=\Im M(z)/\Im m(z)$. The ratio $\Gamma$ is positive in either
half-plane and has normalized trace one. The next theorem compares the resolvent of $A$ with $M$, in quadratic
forms and in normalized traces.

\begin{theorem}[Cipolloni--Erd\H{o}s--Henheik]\label{thm:external-local}
Suppose $s^{-1/2}\norm F\le C \Naux ^A$ for fixed $A,C$.
Let $G_A(z)=(A-z)^{-1}$. For every fixed $\epsilon>0$, uniformly over
$\Naux |\Im z|\rho_s(z)\ge \Naux ^\epsilon$ and
$|z|\le\sqrt s\,\Naux ^{A+1/\epsilon}$, deterministic vectors $x,y$ and
deterministic matrices $B$,
\begin{align}
 |x^*(G_A(z)-M(z))y|
 &\prec\sqrt{\frac{\rho_s(z)}{\Naux |\Im z|}}
             \sqrt{(x^*\Gamma(z)x)(y^*\Gamma(z)y)},
             \label{eq:external-isotropic}\\
 |\av{(G_A(z)-M(z))B}|
 &\prec\frac{\sqrt{\av{\Gamma(z)BB^*}}}{\Naux |\Im z|}.
             \label{eq:external-averaged}
\end{align}
\end{theorem}

The bound \eqref{eq:external-isotropic} is the isotropic law and
\eqref{eq:external-averaged} is the averaged law. Both keep the relative
matrix $\Gamma(z)$ on the right side. Replacing $\Gamma$ by an operator-norm
bound would lose polynomial factors for the frames used later.

These are \cite[Theorem 2.3, equations (2.8)--(2.9)]{CEH}, with
$\lambda=\sqrt s$ and normalized Hilbert--Schmidt norm in the source.
Their equation (2.8) is the averaged bound and (2.9) is the isotropic bound.
GUE has circular off-diagonal coordinates, so its parameter
$\E\chi_{{od}}^2$ is zero and satisfies Assumption 2.1 there.

\begin{remark}[Rigidity on a window]\label{rem:ceh-rigidity}
When $s$ and the free-convolution density $\rho_s$ stay bounded above and
below on a fixed window, the averaged law \eqref{eq:external-averaged} also
gives bulk eigenvalue rigidity on that window.
\end{remark}

The next result is the consequence of the eigenvector theorem
of~\cite{CEH} for a shrinking window around a regular energy $E_0$. It
compares quadratic forms of eigenvectors with a deterministic quantity.
Recall from \eqref{eq:mde} that $M_s$ is the deterministic approximation
of the resolvent, $m_s=\av{M_s}$ its normalized trace,
$\rho_s=\pi^{-1}|\Im m_s|$ the corresponding density, and
$\Gamma=\Im M_s/\Im m_s$ the relative covariance matrix. We also use the
\emph{fluctuation height} $\eta_f(x)$, defined by
\[
 \Naux \eta_f(x)\rho_s(x+i\eta_f(x))=1 .
\]
It is the height above the energy $x$ at which a window of the spectrum
contains, on average, about one eigenvalue; in the bulk it is of order
$\Naux^{-1}$.
 
\begin{theorem}[Regular-window consequence of the CEH eigenvector theorem]
\label{thm:external-eth}
Work in the setting of Theorem~\ref{thm:external-local}. Let the variance
$s$ range over a fixed compact subset of $(0,\infty)$, and assume the
following regularity on a fixed neighborhood of $E_0$, i.e., 
the boundary values $m_s$ and $M_s$ have bounded derivatives, $\Im m_s$
is bounded below, $\norm{M_s}$ is bounded, and $|1-s\av{M_s^2}|$ is
bounded below. Fix a sufficiently small $\beta>0$ and put
$h=\Naux^{-1+\beta}$. Let $J$ be a deterministic Hermitian matrix with
$\norm J\le1$. Then, for every fixed $0<\xi<\beta$,
\begin{equation}\label{eq:external-eth}
 \bigl|u_i^*Ju_j-\delta_{ij}\av{\Gamma(\lambda_i+i\Naux^\xi
                  \eta_f(\lambda_i))J}\bigr|\prec \Naux^{-1/2},
\end{equation}
uniformly over all eigenvalues $\lambda_i,\lambda_j$ in the window
$[E_0-h,E_0+h]$. The probability estimate refers to one event, on which
the bound holds simultaneously for every pair of indices whose
eigenvalues lie in the window. The regularity constants and the compact
range of $s$ are fixed before the dimension $\Naux$ grows.
\end{theorem}
 
Take two eigenvectors $u_i$ and $u_j$ whose eigenvalues lie in a
window of length $2\Naux^{-1+\beta}$ around $E_0$, which contains only a
small power of $\Naux$ eigenvalues. Their overlap through $J$ is then
deterministic up to an error of order $\Naux^{-1/2}$, apart from the small
loss built into $\prec$. For two different eigenvectors the overlap is
close to zero. For the same eigenvector it is close to a normalized trace
of $J$ against the relative covariance $\Gamma$, evaluated slightly above
the real axis, at height $\Naux^\xi$ times the fluctuation height. The
choice of that height does not matter, because $\Gamma$ varies smoothly
near a regular energy; the proof uses exactly this fact.
 
\begin{proof}
\emph{Plan of the proof.} The estimate is a special case
of~\cite[Theorem 2.5, (2.12)--(2.15)]{CEH}, restricted to a short window
around a regular energy. Two things must be checked. First, the
hypotheses and the error term of that theorem take a simple form in our
regular window (Step~1). Second, the theorem of~\cite{CEH} evaluates the
deterministic centre at a height that is tied to the size of the
stochastic-domination loss, whereas in \eqref{eq:external-eth} the height
exponent $\xi$ is fixed in advance. We remove this mismatch by showing
that changing the height moves the centre by much less than the error we
allow (Step~2). Step~3 collects the probability bounds.
 
\emph{Step 1: the theorem of~\cite{CEH} in a regular window.} That theorem
contains an indicator,
\[
 \theta_{ij}=\one_{\{s\Re\av{M_iM_j^*}>1/2\}},
\]
where $M_i$ and $M_j$ denote the deterministic approximation evaluated
at the two eigenvalues. It selects pairs of eigenvalues that are close
enough for their eigenvectors to interact. In our window it equals one
for large $\Naux$, because both eigenvalues are within $h$ of $E_0$ and
$s\av{|M_s(E_0)|^2}=1$ at a regular energy. The theorem is stated for the
regularized observable
\[
 J^\circ=J-\av{M_iJM_j^*}\av{M_iM_j^*}^{-1}I,
\]
from which the part of $J$ proportional to the identity, in the
appropriate sense, has been removed. Its error term is governed by a
control function built from these matrices. In our setting all the
matrices in this control function are bounded. The only denominator,
$1-s\av{M_iM_j}$, stays away from zero by the assumed stability and by
continuity of $M_s$ near $E_0$. Hence the control function is bounded,
and the theorem of~\cite{CEH} gives an error of order $\Naux^{-1/2}$ up to
the stochastic-domination loss, with the centre evaluated at an
auxiliary height.
 
\emph{Step 2: moving the centre to the prescribed height.} Recall that
$X\prec Y$ means that, for every fixed loss $\delta>0$ and every $D>0$,
the event $\{X>\Naux^\delta Y\}$ has probability at most $\Naux^{-D}$ for
large $\Naux$. Fix $0<\xi<\beta<1/2$ and a loss $\delta>0$, and put
$\xi'=\min(\xi,\delta/4)$. Write
\[
 a_\xi(x)=\av{\Gamma(x+i\Naux^\xi\eta_f(x))J}
\]
for the centre at height exponent $\xi$. Both $a_\xi(x)$ and $a_{\xi'}(x)$
are close to the boundary value of the same function at the real energy
$x$. Since this function has a bounded derivative in the regular class,
and both heights are of order $\Naux^{\xi-1}$ and $\Naux^{\xi'-1}$, we get,
uniformly in the regular class,
\[
 |a_\xi(x)-a_{\xi'}(x)|\le C(\Naux^{\xi-1}+\Naux^{\xi'-1})
 <\Naux^{\delta-1/2}-\Naux^{\delta/2-1/2}
\]
for all sufficiently large $\Naux$. The last inequality holds because
$\xi<1/2$, so the left side is smaller than any fixed power above
$\Naux^{-1/2}$.
 
This inequality says that the two centres are much closer to each other
than the allowed error. Consequently, if the overlap misses the centre
$a_\xi$ by more than $\Naux^{\delta-1/2}$, then it misses the centre
$a_{\xi'}$ by more than $\Naux^{\delta/2-1/2}$. In other words, the bad
event at centre $a_\xi$ with loss $\delta$ is contained in the bad event
at centre $a_{\xi'}$ with loss $\delta/2$.
 
\emph{Step 3: the probability bound.} At the height exponent
$\xi'<\delta/2$, the auxiliary-height estimate of Step~1 applies with loss
$\delta/2$. We use it with tail exponent $D+20$. A union bound over at
most $\Naux^{20}$ pairs of grid points then gives a single event of
probability at least $1-\Naux^{-D}$ on which the bound holds for all pairs
of eigenvalues in the window. By Step~2, the same event controls the
centre $a_\xi$ with loss $\delta$. Since $\delta$ and $D$ were arbitrary,
this proves \eqref{eq:external-eth} with the order of quantifiers stated
in the theorem. We do not claim the fixed-$\xi$ statement for
deformations outside the regular class.
\end{proof}

\begin{remark}[What the check gives]\label{rem:eth-loss}
The check keeps a fixed stochastic-domination loss; see
Remark~\ref{rem:domination-loss}. It does not give an eigenvector bound with
a constant factor.
\end{remark}

\subsection{Fixed-energy complex Gaussian statistics}\label{subsec:external-sine}

In this subsection we prove that a deformed complex GUE matrix has sine
statistics at a fixed energy in the bulk, with an explicit error
$O(\Naux ^{-1/2})$ recorded in Theorem~\ref{thm:external-sine}. In every
application below, the Gaussian part is a full fresh block whose variance
is bounded below independently of the dimension. For this case we work
directly with the exact finite kernel of the Gaussian ensemble with
external source, and evaluate it by a steepest descent argument. In
particular, we do not need the small-time complex extension stated in a
remark of Landon--Sosoe--Yau~\cite{LSY}.

\begin{theorem}[Fixed-variance complex Gaussian bulk statistics]
\label{thm:external-sine}
Let $F=F^*\in\C^{\Naux \times \Naux }$ with $\|F\|\le \Naux ^A$, and let
$H=F+\sqrt t\,G_\Naux $, where $G_\Naux $ is a normalized complex GUE matrix and
$0<t_-\le t\le t_+<\infty$. Fix $E$ in a compact real interval. Suppose
that the solution of \eqref{eq:mde} at $E$ with positive imaginary part
satisfies $\Im m_{F,t}(E)\ge v_*>0$, and put
$\rho=\pi^{-1}\Im m_{F,t}(E)$. Then, after conjugating the correlation
kernel by a real scalar factor, we have, uniformly for $\xi,\zeta$ in
compact real sets,
\begin{equation}\label{eq:t26-complex-sine-kernel}
 \frac{e^{-h_\Naux (\xi)}K_\Naux (E+\xi/(\Naux \rho),E+\zeta/(\Naux \rho))
             e^{h_\Naux (\zeta)}}{\Naux \rho}
 =\frac{\sin\pi(\xi-\zeta)}{\pi(\xi-\zeta)}+O(\Naux ^{-1/2}),
\end{equation}
with the diagonal value taken by continuity. The constants depend only
on $A,t_-,t_+,v_*$ and the compact sets. Hence every fixed factorial
correlation, and the local point process, converge to those of the
unit-intensity sine process; in particular, the normalized mean density
of points at $E$ is $\rho+O(\Naux ^{-1/2})$. We assume neither a limiting
empirical measure of $F$ nor a bound on $\|F\|$ independent of $\Naux $.
\end{theorem}

\begin{proof}
The eigenvalues of a Gaussian matrix with a fixed
matrix added to it form a determinantal process. All their correlation
functions are determinants built from the
correlation kernel. For this ensemble the kernel has an exact formula as a
double contour integral(Step~1). We show that near the energy $E$, after rescaling by the
local density, this kernel converges to the sine kernel. The method is
steepest descent. We deform the two curves so that they pass through the
saddle points of the exponent, the points where its derivative vanishes, and
so that the integrand is largest there and decays quickly away from them, as
in Figure~\ref{fig:sine-contours} (Steps~2--5). Moving one curve across the
other picks up a residue, which is exactly the
sine kernel (Step~6), while the remaining double integral is small
(Step~7). Step~8 passes from the kernel to the point process.

\emph{Step 1:} By unitary invariance we may take $F$
diagonal, with eigenvalues $d_j$. We first fix the normalization. For
distinct increasing $d_j$, the symmetrized density of the eigenvalues of
$H$ is
\[
 \frac{1}{\Naux !\Delta(d)}
 \det[p_{t/\Naux }(x_i-d_j)]_{i,j=1}^\Naux \,\Delta(x),\qquad
 p_s(x)=(2\pi s)^{-1/2}e^{-x^2/(2s)}.
\]
This is the Gaussian density with external source, that is, the
Vandermonde transform of the kernel of independent Brownian motions. Its
total mass is one, with no unknown constant: by Andreief's identity the
integral equals $\det[\int p_{t/\Naux }(x-d_i)x^{j-1}dx]/\Delta(d)$, and this is
one because Gaussian convolution keeps monic polynomials monic. The
factor $\Naux !$ belongs to the density on $\R^\Naux $ and is absent on the ordered
chamber. The resulting biorthogonal kernel is the exact formula
\cite[(2.18)]{JohanssonKernel} with $S=t/\Naux $:
\begin{equation}\label{eq:t26-johansson-kernel}
 K_\Naux (u,v)=\frac{\Naux e^{\Naux (v^2-u^2)/(2t)}}{(2\pi i)^2t}
 \int_\gamma dz\int_{\Gamma_L}dw\,
 \frac{e^{\Naux (w^2-2vw-z^2+2uz)/(2t)}}{w-z}
            \prod_{j=1}^\Naux \frac{w-d_j}{z-d_j}.
\end{equation}
Here $\gamma$ encircles all $d_j$ counterclockwise, and $\Gamma_L=L+i\R$
is oriented upward and lies strictly to the right of $\gamma$; this fixes
both signs. With this normalization the factorial densities are
$\det[K_\Naux (x_i,x_j)]$, including the factor $\Naux !$. Repeated $d_j$ follow by
continuity of the Gaussian law and of \eqref{eq:t26-johansson-kernel}
on contours away from the real poles.

\emph{Step 2:} We put
\[
 z_0=x_0+iy_0=E+tm_{F,t}(E),\qquad
 S(z)=\frac{(z-E)^2}{2t}+\Naux ^{-1}\sum_j\log(z-d_j).
\]
The equation for $m_{F,t}$ gives $S'(z_0)=0$, $y_0=t\pi\rho\ge t_-v_*$, and
$\Naux ^{-1}\sum_j|z_0-d_j|^{-2}=1/t$. By Cauchy--Schwarz,
$|m_{F,t}(E)|\le t^{-1/2}$, so $x_0$ stays in a fixed compact set and
$y_0\le\sqrt t$. Globally we only use $\Re S$; near the two nonreal
saddles $z_0,\bar z_0$ we choose analytic branches of the logarithms.

\emph{Step 3:} For real $x$, let $y(x)>0$ solve
\[
 \Naux ^{-1}\sum_j((x-d_j)^2+y(x)^2)^{-1}=1/t
\]
when the left side at $y=0$ exceeds $1/t$, and set $y(x)=0$ otherwise. The
graphs of $\pm y$, with the shared real gaps cancelled, bound a union of
intervals containing all $d_j$. After an arbitrarily small smoothing at
the real endpoints, we use them as $\gamma$. They lie in
$[\min d_j-\sqrt t,\max d_j+\sqrt t]+i[-\sqrt t,\sqrt t]$, and their total
length is polynomial in $\Naux $. Indeed, for fixed $y>0$, clearing
denominators turns the level equation into a polynomial equation in $x$
of degree at most $2\Naux $, so the coarea formula gives
$\int|y'(x)|dx\le2\Naux \sqrt t$, while the horizontal length is $O(\Naux ^A+1)$.
The finitely many critical points and endpoints do not affect this.

\emph{Step 4:} With
$a_j=x-d_j$, put $B=\langle(a_j^2+y^2)^{-2}\rangle$,
$A_1=\langle a_j(a_j^2+y^2)^{-2}\rangle$ and
\[
 \psi(x)=x+t\langle a_j/(a_j^2+y(x)^2)\rangle.
\]
Where $y>0$, differentiation gives
\[
 y'=-A_1/(yB),\qquad
 \psi'=2t(y^2B+A_1^2/B)>0,\qquad
 \frac{d}{dx}\Re S(x+iy(x))=(\psi(x)-E)/t.
\]
On a real gap, $\psi'=1-t\langle(x-d_j)^{-2}\rangle\ge0$, so by continuity
$\psi$ is increasing everywhere. Near $x_0$ all derivatives are bounded,
$y\ge y_0/2$, and Jensen's inequality gives $B\ge t^{-2}$, so
$\psi'\ge c>0$. Hence there is a fixed $\delta>0$ with
\begin{equation}\label{eq:t26-contour-horizontal}
 \Re S(x\mathbin\pm iy(x))-\Re S(z_0)
 \ge c\min\{|x-x_0|^2,\delta|x-x_0|\}.
\end{equation}
For $|x-x_0|>\delta$ this follows by integrating the monotone $\psi$, since
$\psi-E$ already has a fixed sign and size at $x_0\pm\delta$. The formula
for $y'$ and Cauchy--Schwarz also give $|y'|\le C$ near the saddles.

\begin{figure}[ht]
\centering
\begin{tikzpicture}[x=1cm,y=1cm,>=Stealth]
 \draw[->,black!60] (-4.6,0)--(5.7,0) node[below] {$\Re z$};
 \draw[very thick,postaction={decorate},
   decoration={markings,mark=at position .30 with {\arrow{>}},
                        mark=at position .80 with {\arrow{>}}}]
   (2.8,0) .. controls (2.5,0.9) and (1.4,1.35) .. (0,1.4)
           .. controls (-1.6,1.45) and (-3.4,1.0) .. (-3.6,0)
           .. controls (-3.4,-1.0) and (-1.6,-1.45) .. (0,-1.4)
           .. controls (1.4,-1.35) and (2.5,-0.9) .. (2.8,0);
 \draw[very thick,postaction={decorate},
   decoration={markings,mark=at position .25 with {\arrow{>}}}]
   (4.4,0) .. controls (4.35,.45) and (3.45,.45) .. (3.4,0)
           .. controls (3.45,-.45) and (4.35,-.45) .. (4.4,0);
 \foreach \x in {-3.2,-2.6,-2.3,-1.7,-1.1,-0.8,-0.2,0.5,0.9,1.6,2.1,2.5,3.7,4.1}
   \fill (\x,0) circle (1.6pt);
 \node[below,font=\small] at (-2.0,-.05) {$d_j$};
 \draw[line width=3.2pt,red!70!black,opacity=.30]
   (0.6,1.38) .. controls (1.5,1.33) and (2.5,0.9) .. (2.8,0)
              .. controls (2.5,-0.9) and (1.5,-1.33) .. (0.6,-1.38);
 \draw[thick,blue!70!black,->] (0.6,-2.6)--(0.6,2.7)
   node[above,font=\small] {$\Gamma=x_0+i\R$};
 \draw[line width=2.2pt,red!70!black] (0.6,-1.38)--(0.6,1.38);
 \fill[red!70!black] (0.6,1.38) circle (2.3pt) node[above left,font=\small] {$z_0$};
 \fill[red!70!black] (0.6,-1.38) circle (2.3pt) node[below left,font=\small] {$\bar z_0$};
 \draw[thick,black!55,dotted,->] (5.1,-2.6)--(5.1,2.7)
   node[above,font=\small,black!70] {$\Gamma_L$};
 \node[font=\small] at (-2.3,1.65) {$\gamma$: graph of $y(x)$};
 \node[font=\small] at (-2.3,-1.7) {graph of $-y(x)$};
 \node[font=\small,red!70!black,anchor=west] at (3.0,1.15) {residue arc};
 \draw[red!70!black,->,thin] (2.95,1.1)--(2.45,.75);
 \node[font=\small,red!70!black,anchor=east] at (0.45,-0.55) {segment};
\end{tikzpicture}
\caption{The contours in the proof of Theorem~\ref{thm:external-sine}. The
dots are the eigenvalues $d_j$ of the deformation $F$. The contour $\gamma$
(black, counterclockwise) is the union of the graphs of $\pm y(x)$, where
$y(x)$ solves the level equation of Step~3; it may have several
components, each enclosing a group of the $d_j$. The vertical contour
starts as $\Gamma_L$ (dotted), to the right of every component, and is
moved left to $\Gamma=x_0+i\R$ (blue), which crosses $\gamma$ exactly at
the saddle points $z_0$ and $\bar z_0$. Along $\gamma$ the function
$\Re S$ is smallest at the saddles, by \eqref{eq:t26-contour-horizontal},
and along $\Gamma$ it is largest there, by Step~5. This is the steepest
descent geometry that makes the remaining double integral of order
$\Naux ^{-1/2}$ in Step~7. Moving the vertical contour across the pole $w=z$
leaves a residue integrated over the part of $\gamma$ to the right of
$\Gamma$ (shaded). Components entirely on that side, such as the small
one on the right, contribute nothing. On the main component the arc from
$\bar z_0$ to $z_0$ may be replaced by the straight segment between them,
which lies on $\Gamma$ (solid red); its integral is the sine kernel of
Step~6.}
\label{fig:sine-contours}
\end{figure}
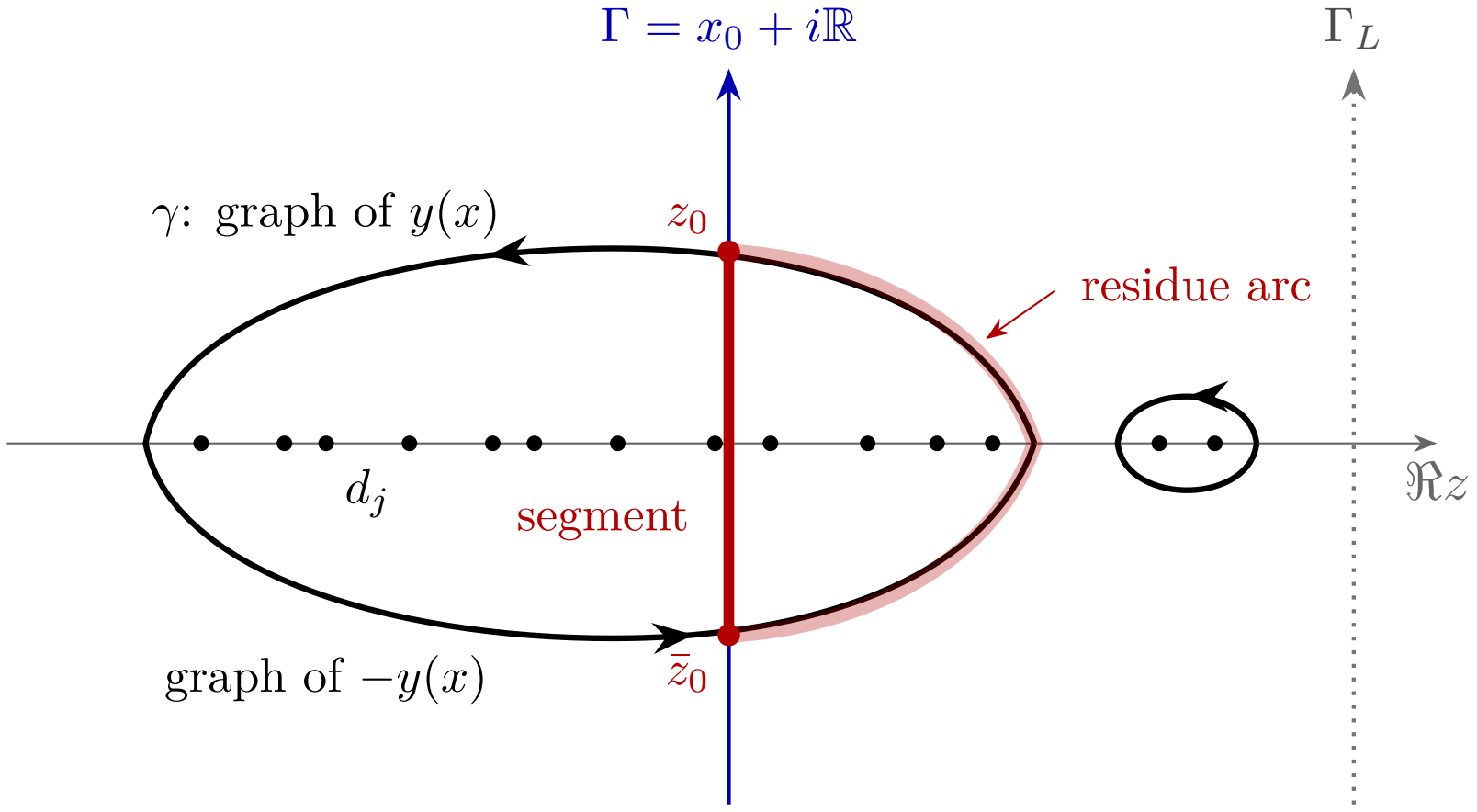

\emph{Step 5:} We move the
vertical contour to $\Gamma=x_0+i\R$. Along it,
\[
 \frac{d}{dv}\Re S(x_0+iv)
 =v\left[\left\langle((x_0-d_j)^2+v^2)^{-1}\right\rangle-1/t\right].
\]
The bracket decreases strictly for $v>0$ and vanishes at $y_0$, so the two
global maxima are at $v=\pm y_0$, with second derivatives
$-2y_0^2\langle|z_0-d_j|^{-4}\rangle\le-c$. So $\Re S$ drops quadratically
near each maximum and by a fixed amount away from them. For large $v$,
\[
 \Re S(x_0+iv)-\Re S(z_0)
 \le-\frac{v^2-y_0^2}{2t}
       +\tfrac12\log(1+v^2/y_0^2),
\]
which gives a Gaussian tail. None of these bounds depends on the
eigenvalues $d_j$ far from $E$.

\emph{Step 6:} Moving the upward
contour to the left crosses the pole $w=z$. The residue is integrated
over the part of $\gamma$ to the right of $\Gamma$; components lying
entirely on that side contribute zero. On the component containing
$x_0$, this arc runs from $\bar z_0$ to $z_0$, and its integral equals the
integral over the straight segment. For $u=E+\xi/(\Naux \rho)$ and
$v=E+\zeta/(\Naux \rho)$, the normalized residue is
\[
 \frac{e^{\Naux (v^2-u^2)/(2t)}}{2\pi it\rho}
 \int_{\bar z_0}^{z_0}e^{(\xi-\zeta)z/(t\rho)}dz
 =e^{h_\Naux (\xi)-h_\Naux (\zeta)}
                     \frac{\sin\pi(\xi-\zeta)}{\pi(\xi-\zeta)},
\]
with $h_\Naux (\xi)=\xi(x_0-E)/(t\rho)-\xi^2/(2\Naux t\rho^2)$. The orientation
gives the diagonal value $+1$.

\emph{Step 7:} Where the contours
meet, they cross transversally: the graph has horizontal tangent
component one and bounded slope, while $\Gamma$ is vertical. In local
coordinates $a=x-x_0$ and $b=v\mp y_0$ their distance is at least
$c\sqrt{a^2+b^2}$, so the normalized integrand is bounded locally by
\[
 C\frac{e^{-c\Naux (a^2+b^2)}}{\sqrt{a^2+b^2}}\,da\,db,
\]
whose integral is $O(\Naux ^{-1/2})$. Pairs of opposite saddle neighbourhoods
have a denominator bounded away from zero and give $O(\Naux ^{-1})$. Away from
the saddles, the strict drops of $\Re S$, the polynomial length of
$\gamma$ and the Gaussian tail on $\Gamma$ give an exponentially small
contribution. The microscopic factors are at most $e^{C|x-x_0|}$ and are
absorbed by \eqref{eq:t26-contour-horizontal}. Since the transverse
singularity is integrable in two variables, removing small
neighbourhoods and letting their radius tend to zero also justifies the
contour crossing. This proves \eqref{eq:t26-complex-sine-kernel}.

\emph{Step 8:} The conjugation cancels in every correlation
determinant, so uniform kernel convergence gives all fixed factorial
correlations on compact sets. The counting bound of
Theorem~\ref{thm:external-count} gives uniform exponential moments of
local counts, so these correlations determine the subsequential limits
and give tightness. The diagonal of the kernel gives the density, and all
estimates are uniform in the stated class.
\end{proof}

\subsection{The conditional eigenvector moment estimate}\label{subsec:external-emf}

We need to understand eigenvector overlaps after the eigenvalue path
has been fixed. Theorem~\ref{thm:external-emf} gives their conditional
moments for a real Gaussian evolution. We state its assumptions and
covariance profile, then compare the profiles at positive and zero
height. Section~\ref{sec:fresh-gaussian} supplies the Hermitian argument
and retains the spectral weight required by our application.

\paragraph{Setting.} Let $D$ be a deterministic diagonal
$(\eta_*,r_\Naux )$ -regular matrix in the sense of \cite[Definition 1.1]{Benigni},
and assume
\[
 \Naux ^{-1}\le\eta_*\le \Naux ^{-\epsilon_0},\qquad
 \Naux ^{\epsilon_0}\eta_*\le r_\Naux \le \Naux ^{-\epsilon_0}.
\]
Regularity means that $\Im\av{(D-x-i\eta)^{-1}}$ lies between fixed
positive constants for $|x-E_0|\le r_\Naux $ and $\eta_*\le\eta\le10$. Let
$G_\Naux ^{\rm O}$ be normalized GOE: its off-diagonal variance is $\Naux ^{-1}$ and
its diagonal variance is $2\Naux ^{-1}$. Let $D+\sqrt t\,G_\Naux ^{\rm O}$ be evolved
for a further Gaussian time $\upsilon$. The scale conditions are
\[
 \Naux ^\omega\eta_*\le t\le \Naux ^{-\omega}r_\Naux ,\qquad
 \Naux ^{-1+a_0}\le\upsilon\le \Naux ^{-a_0}t,
\]
for fixed positive slack $\omega,a_0$. For a direction $q$, write
 \[
 \widetilde\sigma_t^2(q,x,h)=
 \sum_a\frac{t|q_a|^2}{(D_a-x-t\Re m_{D,t}(x+ih))^2+
                         (t\Im m_{D,t}(x+ih))^2}.
\]
For a particle configuration $\xi=(\xi_i)$ of fixed total size
$r$, its covariance product is
$\Sigma_t(q,\xi,h)=\prod_i\widetilde\sigma_t^2(q,\gamma_{i,t},h)^{\xi_i}$.
This is the profile in \mbox{%
\cite[equation (1.3)]{Benigni}}\hskip0pt%
. The positive
smoothing height is used to evaluate $m_{D,t}$.

The next theorem says that these conditional moments are, to leading
order, the moments of independent Gaussian overlaps with variances
$\widetilde\sigma_t^2$.

\begin{theorem}[Benigni, real conditional moment-flow estimate]\label{thm:external-emf}
Let $q$ be a deterministic real unit vector. For normalized real
eigenvectors $u_i$, let $f_\upsilon(\xi)$ be the conditional
expectation, given the complete eigenvalue path, of
$\prod_i|\sqrt \Naux q^Tu_i|^{2\xi_i}/(2\xi_i-1)!!$,
with $(-1)!!=1$. On the trajectory event of
\cite[Lemma 3.2]{Benigni}, these moments satisfy
\begin{equation}\label{eq:external-emf}
 f_\upsilon(\xi)
 =\Sigma_t(q,\xi,\Naux ^{-\epsilon}\upsilon)
 \left[1+O_r\left(
 \Naux ^{p_r\epsilon}\{(\Naux \upsilon)^{-1/2}
                          +(\upsilon/t)^{1/3}\}\right)\right]
\end{equation}
for configurations in a fixed interior part of the regular energy
window; $p_r<\infty$ depends only on $r$.
\end{theorem}

The relative error in \eqref{eq:external-emf} is small when
$\Naux ^{-1}\ll\upsilon\ll t$. The estimate is conditional on the full
eigenvalue path, and it holds on an event that is described next.

This is \cite[Proposition 3.4, (3.15), conditioning (1.20)]{Benigni}.
The hypotheses on good trajectories are: $(i)$ local rigidity at precision $\Naux ^\epsilon(\Naux ^{-1}+\upsilon)$; $(ii)$ the averaged and directional resolvent laws on $\Naux ^\vartheta/\Naux \le\eta\le \Naux ^{-\vartheta}t$ in that interior window, with respective relative errors $\Naux ^\epsilon/(\Naux \eta)+\upsilon/t$ and
$\Naux ^{2\epsilon}/\sqrt{\Naux \eta}+\upsilon/t$; and  $(iii)$ finite propagation of the
retained short-range moment flow. The cited lemma supplies these with
arbitrarily high fixed polynomial probability. Its short-range distance is
an index distance, chosen between $\Naux \upsilon$ and $\Naux t$ with fixed powers of
slack.

\begin{remark}[The trajectory event]\label{rem:emf-event}
The event in~\cite{Benigni} includes eigenvector resolvent estimates.
 The disintegration that we 
need for a spectral weight is part of the proof of
Theorem~\ref{thm:fresh-palm}.
\end{remark}

\begin{lemma}[Relative comparison of the two covariance profiles]
\label{lem:t26-positive-height-profile}
Let $D=D^*$ have arbitrary finite dimension, let $t>0$, $h\ge0$, and
let $m=u+iv$ with $v>0$. Define the positive operators
\[
 \widetilde\Gamma=t\bigl((D-x-tu)^2+t^2v^2I\bigr)^{-1},\qquad
 \Gamma_h=\frac{\Im(D-x-ih-tm)^{-1}}v,
 \qquad \epsilon_h=\frac{h}{tv}.
\]
Then, without a bound on $\norm D$ or a lower bound on a tested direction,
\begin{equation}
\label{eq:t26-profile-loewner}
 (1+\epsilon_h)^{-1}\widetilde\Gamma
 \preceq\Gamma_h\preceq(1+\epsilon_h)\widetilde\Gamma.
\end{equation}
In particular if $v\ge c_*>0$, then uniformly over all nonzero vectors $q$,
\[
 \left|\frac{q^*\Gamma_hq}{q^*\widetilde\Gamma q}-1\right|
 \le\frac{h}{tc_*}.
\]
Products of $r$ such quadratic forms have relative ratio between
$(1+\epsilon)^{-r}$ and $(1+\epsilon)^r$ when every
$\epsilon_h\le\epsilon$. Thus at fixed $r$ the relative error is
$O_r(h/t)$ on the regular window, uniformly in dimension and direction,
when $h/t\to0$. This comparison advances the covariance normalization in
Theorem~\ref{thm:fresh-palm}; it supplies no conditional moment-flow estimate.
\end{lemma}

\begin{proof}
Diagonalize $D$ and put $a=D_j-x-tu$, $b=tv>0$. The ratio of the
corresponding two positive diagonal entries is
\[
 (1+\epsilon_h)\frac{a^2+b^2}{a^2+(1+\epsilon_h)^2b^2}.
\]
Since $a^2+b^2\le a^2+(1+\epsilon_h)^2b^2
\le(1+\epsilon_h)^2(a^2+b^2)$, this ratio lies between
$(1+\epsilon_h)^{-1}$ and $1+\epsilon_h$. Functional calculus gives
\eqref{eq:t26-profile-loewner}; taking quadratic forms and multiplying
proves the remaining assertions. At $h=0$ the two profiles agree exactly.
For $h>0$ they are generally different.
\end{proof}

 \subsection{Finite propagation of the conditional moment flow}\label{subsec:external-finite-speed}

We need a quantitative limit on how far conditional eigenvector
moments can travel during a short Gaussian evolution.
 Theorem~\ref{thm:external-finite-speed} shows that the flow with short
jumps has a negligible effect at sufficiently distant labels. We state
the distance and time scales, give the estimate, and explain how its
single scale differs from our later two-scale localization argument.

Write $a_\Naux \ll b_\Naux $ when $\Naux ^\epsilon a_\Naux \le b_\Naux $ for some fixed $\epsilon>0$.

\begin{theorem}[Bourgade--Huang--Yau]\label{thm:external-finite-speed}
Suppose that $\norm{H_0}\le \Naux ^A$ and
$c\le\Im\av{(H_0-x-i\eta)^{-1}}\le C$ for
$|x-E_0|\le r$, $\eta_*\le\eta\le1$, where
$\Naux ^{-1}\ll\eta_*\ll r\le1$.
Consider real Dyson Brownian motion with
$\eta_*\ll t_0\le t\le t_0+\ell/\Naux \ll r$.
Fix a particle number $k$ and a fixed interior energy window.
Let $U_{\rm sh}(s',s)$ be the conditional moment-flow propagator
retaining only jumps of index distance at most $\ell$.
For configurations $\xi=(x_1,\ldots,x_k)$ and
$\zeta=(y_1,\ldots,y_k)$ in increasing order, let
$d_{\rm in}(\xi,\zeta)$ be the largest, over $a\le k$, of the
number of classical eigenvalue locations in the interior window
between $x_a$ and $y_a$.
Put $\psi=\Naux ^\epsilon$, with a fixed sufficiently small $\epsilon>0$.
Outside an event of arbitrarily high polynomially small probability,
uniformly for $t_0\le s'\le s\le t$ and bounded $h$ whose support
has $d_{\rm in}$-distance at least $\psi\ell$ from $\xi$,
\begin{equation}\label{eq:external-finite-speed}
 |U_{\rm sh}(s',s)h(\xi)|
       \le \norm h_\infty \Naux ^k e^{-c\psi}.
\end{equation}
\end{theorem}

The right side of \eqref{eq:external-finite-speed} is smaller than every
negative power of $\Naux $, because $\psi=\Naux ^\epsilon$. So a function supported
at index distance at least $\Naux ^\epsilon\ell$ from $\xi$ has a negligible
effect at $\xi$.

This is \cite[Assumption 1.3, Lemma 3.2 and Corollary 3.3]{BHY}.
The proof uses local rigidity and a smoothed moving eigenvalue weight.

\begin{remark}[Modifications proved later]\label{rem:bhy-modifications}
Section~\ref{sec:relative-vectors} repeats the weighted energy estimate
of~\cite{BHY} with smoothing length $L>\ell$, the complex one-particle
generator and nonnegative Ornstein--Uhlenbeck killing. Those modifications
are proved there; they are not part of the statement imported here.
\end{remark}

\subsection{A local law on a shorter band ring}\label{subsec:external-band}

We use a shorter ring to determine the mean spectral density of the
original matrix. Theorem~\ref{thm:external-band-local} controls individual
resolvent entries on that shorter ring, where the bandwidth has  a fixed
power margin. We state the specialization of~\cite{ER}, describe its
spatial decay, and record the scope needed for the density comparison.

\begin{theorem}[Erd\H{o}s--Riabov, short-ring specialization]
\label{thm:external-band-local}
Let $K$ have the ring law of this paper with $m$ blocks, $m\geq 4$, and dimension
$\Naux =mW$. Fix $\zeta_0>0$ and suppose $W^2\ge \Naux ^{1+\zeta_0}$.
For a fixed bulk interval, fixed $0<\delta_0<\zeta_0$, and
$\eta\ge \Naux ^{-1+\delta_0}$, put $\ell=\min(\Naux ,W/\sqrt\eta)$. Also require $|e+i\eta|\le C_0$ for a fixed domain constant $C_0$. 
For every fixed $D\ge6$,
\begin{equation}\label{eq:short-ring-local}
 |G^K_{xy}(e+i\eta)-m_{\rm sc}(e+i\eta)\delta_{xy}|
 \prec\left\{\frac{1}{\ell\eta}
       \left(1+\frac{|x-y|_\Naux}{\ell}\right)^{-D}\right\}^{1/2}.
\end{equation}
\end{theorem}

The bound \eqref{eq:short-ring-local} compares every resolvent entry of $K$
with the semicircle Stieltjes transform $m_{\rm sc}$. The error has polynomial spatial decay of any fixed order $D\ge6$ on the scale $\ell$..

This is an internal specialization of \cite[Theorem 3.3, (3.17)]{ER}, using the fixed-constant polynomial admissibility check below; the displayed estimate is not a separate native input. Here $|x-y|_\Naux $ is circular site distance; changing an
orbital position within its block changes this distance by at most $W$.

\begin{remark}[Use on a short ring only]\label{rem:er-short}
The condition $W^2\ge \Naux ^{1+\zeta_0}$  holds on a ring with $m$ blocks only
when $m$ is much smaller than $W$; it does not hold on the critical ring
of Theorem~\ref{thm:critical}. No local law for the full critical ring is
inferred from this specialization. The verification of the variance
profile, including both parities and arbitrary fixed $b\in(0,1/2)$, is
given in Section~\ref{sec:counting}.
\end{remark}

\section{Counting, mean density and spectral recovery}
\label{sec:counting}
This section collects the spectral facts about the original ring that
every later argument uses. We need three things. Short spectral intervals
must rarely contain many eigenvalues, the mean density of eigenvalues must be
the semicircle density even at critical band width, and the critical limit
must have no preferred spectral position or block. These facts serve all
three parts of the proof.
Theorem~\ref{thm:counting} bounds the factorial moments of the number of
eigenvalues in a short interval, uniformly over deterministic deformations.
In words, the eigenvalues of the ring do not cluster, and this gives
tightness of the marked spectrum. Theorem~\ref{thm:critical-density} shows
that the mean density is the semicircle density up to an error $W^{-1/32}$,
uniformly for $n$ up to $TW$, so that every critical limit has unit spectral
intensity. Proposition~\ref{prop:stationarity} shows that every critical
limit is invariant under spectral translations and under rotations of the
circle, with mean intensity $du\,dr$.
We prove Theorem~\ref{thm:counting} at the start of the section. In
Subsection~\ref{subsec:s3-smoothing-and-uniform-integrability} we show that smoothed spectral observables recover the
original ones, with bounded moments. In Subsection~\ref{subsec:s3-mean-density-at-critical-bandwidth} we obtain the
mean density by comparison with shorter rings, on which the local law of
Erd\H{o}s and Riabov applies. In Subsection~\ref{subsec:s3-translation-and-rotation-invariance} we prove the two
invariances: the spectral one from the independent scalar Gaussian part of
the matrix, and the spatial one from cyclic symmetry.

\begin{theorem}\label{thm:counting}
Put $c=1-2b$ and $K_b=C/\sqrt c$, where $C$ is the constant in
Theorem~\ref{thm:external-count}. For every $n\ge4$, $W\ge1$, real
interval $J$, integer $r\ge1$, and deterministic Hermitian matrix $A$,
\begin{equation}\label{eq:factorial-count}
 \E\bigl(\Tr\one_J(H+A)\bigr)_r\le(K_bN|J|)^r,
 \qquad (x)_r=x(x-1)\cdots(x-r+1).
\end{equation}
Fix $|E|<2$, write $\rho=\rho(E)$ and $\mathsf K=K_b/\rho$,
and let $\mathcal N(I)=\#\{\alpha:u_\alpha\in I\}$. Then
\begin{align}
 \E e^{t\mathcal N([-R,R])}
 &\le \exp\{2R\mathsf K(e^t-1)\}, &&t\ge0,
 \label{eq:count-exponential}\\
 \E\sum_{\substack{\alpha\ne\beta\\|u_\alpha|,|u_\beta|\le R}}
 \one_{\{|u_\alpha-u_\beta|<a\}}
 &\le8\mathsf K^2Ra+16\mathsf K^2a^2. \label{eq:close-pairs}
\end{align}
The marked processes are tight in the space of locally finite measures
on $\R\times\mathcal P(\T)$. Every subsequential limit has simple
spectral locations and probability-measure marks.
\end{theorem}

There is no restriction on $n/W$ in this theorem. In particular, the
counting estimates remain available before any limiting law is identified.
The deterministic deformation in \eqref{eq:factorial-count} may have
arbitrarily large norm.

\subsection{Smoothing and uniform integrability}\label{subsec:s3-smoothing-and-uniform-integrability}

 Smoothing makes spectral observables easier to handle, but we must
recover the original eigenvalues without losing mass in distant tails.
Lemma~\ref{lem:smoothing} bounds this recovery error and all fixed
moments of the smoothed observables. We define the two versions, then
apply the counting bounds separately to the smoothing error and to
products of the smoothed densities.

 Let $p_\delta(t)=\delta/[\pi(t^2+\delta^2)]$. For a spatial test $q$ set
\[
 D_W(u,\delta;q)=
 \sum_\alpha p_\delta(u-u_\alpha)\nu_\alpha(q),\qquad
 \nu(q)=\int_\T q\,d\nu.
\]
For $F\in C_c(\R\times\T)$ define
\[
 \zeta_W(F)=\sum_\alpha\int_\T F(u_\alpha,r)\nu_\alpha(dr),
 \qquad
 \zeta_{W,\delta}(F)=\sum_\alpha\int_\R\int_\T
 p_\delta(u-u_\alpha)F(u,r)\nu_\alpha(dr)\,du.
\]

\begin{lemma}\label{lem:smoothing}
Suppose $F$ is supported in $[-R,R]\times\T$ and extended by zero.
Define its modulus in the spectral coordinate by
\[
 \omega_F(h)=\sup_{|u-v|\le h,r}|F(u,r)-F(v,r)|.
\]
For $0<h\le1$,
\begin{equation}\label{eq:uniform-unsmoothing}
 \E|\zeta_W(F)-\zeta_{W,\delta}(F)|
 \le\mathsf K\left[
 2(R+1)\omega_F(h)+\frac{8R\norm F_\infty\delta}{\pi h}\right].
\end{equation}
For fixed $\delta>0$ and every integer $p\ge1$,
\begin{equation}\label{eq:smoothed-moments}
 \sup_{u_1,\ldots,u_p;\,\norm{q_j}_\infty\le1}
 \E\prod_{j=1}^p|D_W(u_j,\delta;q_j)|
 \le p!C_{b,E,\delta}^{\,p}.
\end{equation}
\end{lemma}

The bound pays the distant spectral tails as well as the local smoothing error. Its fixed-order moments justify subsequent expectation limits.

\begin{proof}
The first factorial estimate gives the measure inequality
$\E\sum_\alpha\delta_{u_\alpha}\le\mathsf K\,du$: first on
intervals, and then on all Borel sets by monotone approximation.
Write
$g(t)=\int_\R\sup_r|F(u+t,r)-F(u,r)|\,du$.
For $|t|\le h\le1$ this is at most $2(R+1)\omega_F(h)$; for
every $t$ it is at most $4R\norm F_\infty$. Tonelli's theorem bounds
the left side of \eqref{eq:uniform-unsmoothing} by
$\mathsf K\int p_\delta(t)g(t)\,dt$. Now use
$\int_{|t|>h}p_\delta(t)\,dt\le2\delta/(\pi h)$.
Thus distant eigenvalues, as well as those in the support of $F$, are
included in the estimate.

For the moment bound partition $\R$ into
$I_j=u+[j,j+1)$, and put
$a_j=\sup_{v\in I_j}p_\delta(u-v)$ and $A_\delta=\sum_j a_j$.
The last sum is finite and independent of $u$. Since
$|D_W(u,\delta;q)|\le\sum_j a_j\mathcal N(I_j)$, Jensen's inequality
and \eqref{eq:count-exponential}, applied to unit intervals, give
\[
 \E\exp\{|D_W(u,\delta;q)|/A_\delta\}
 \le\sum_j\frac{a_j}{A_\delta}\E e^{\mathcal N(I_j)}
 \le e^{\mathsf K(e-1)}.
\]
The inequality $x^p\le p!e^x$ bounds each $p$th moment.
H\"older's inequality with exponent $p$ then bounds the product in
\eqref{eq:smoothed-moments}. The same argument at order $2p$
proves uniform integrability of all the products at order $p$.
\end{proof}

\begin{proof}[Proof of Theorem~\ref{thm:counting}]
Condition on every off-diagonal block. The remaining diagonal blocks
are independent $\sqrt c$ times normalized GUE. Divide the matrix and
the interval by $\sqrt c$, apply Theorem~\ref{thm:external-count},
and integrate over the conditioned blocks. This proves
\eqref{eq:factorial-count} with its full deformation uniformity.
The binomial expansion
\[
 e^{tN}=\sum_{r=0}^{N}\frac{(N)_r}{r!}(e^t-1)^r
\]
has nonnegative terms. Its expectation proves
\eqref{eq:count-exponential}. To prove \eqref{eq:close-pairs},
divide the spectral line into intervals $[ja,(j+1)a)$.
If this interval contains the smaller member of a close pair, both
members belong to $[ja,(j+2)a]$. At most $2R/a+4$ such intervals
are needed; each contributes at most $(2\mathsf Ka)^2$ ordered pairs.

The space $\mathcal P(\T)$ is compact. The local exponential count
bound therefore gives tightness of the marked point measures. The
close-pair bound excludes collisions in any subsequential limit: first
restrict to a compact spectral interval, then let $a\downarrow0$.
The marks remain probability measures because their ambient space is
closed and compact.
\end{proof}

\subsection{Mean density at critical bandwidth}\label{subsec:s3-mean-density-at-critical-bandwidth}

  We want the mean density at critical bandwidth, where a local law for
the full ring is unavailable at the required scale.
Theorem~\ref{thm:critical-density} obtains the semicircle density by
comparison with shorter rings. We first check the hypotheses of their
local law, then estimate the cost of cutting the original ring and
remove the remaining energy smoothing.

\begin{lemma}[Complete short-ring profile check]\label{lem:t26-er-admissibility}
Fix $b\in(0,1/2)$ and the bulk-domain constants of
\cite[Definition 2.4]{ER}. Consider a ring of $m$ blocks of size $W$, of
dimension $\Naux=mW$. Let $J_W=W^{-1}\one\one^*$ be the matrix that
averages over one block, and let $C$ be the $m\times m$ block variance
matrix, with entries $C_{kl}=c\delta_{kl}+b\one_{\{l=k\pm1\}}$, indices
read modulo $m$. Then the variance profile $S=C\otimes J_W$ of the short
ring is admissible in the sense of \cite[Definitions 2.1, 2.4 and 2.5]{ER},
uniformly in $m\ge4$ of either parity. More precisely, for any fixed
$D\ge6$ one may take the control function
\begin{equation}\label{eq:t26-er-control}
 U_\eta(x,y)=\frac{1}{\ell\eta}
       \left(1+\frac{|x-y|_\Naux}{\ell}\right)^{-D},
 \qquad \ell=\min(\Naux,W/\sqrt\eta),
\end{equation}
where $|x-y|_\Naux$ is the circular distance between sites. With this
choice, clauses (i)--(vii) of \cite[Definition 2.5]{ER} hold with fixed
constants that depend only on $b$, $D$ and the bulk-domain constants, and
not on $m$ or $W$. No logarithmic loss is allowed in any clause.
\end{lemma}

The length $\ell$ is the distance over which the resolvent spreads at
spectral height $\eta$. It is the diffusive length $W/\sqrt\eta$, capped
by the size of the ring. The control function has size $1/(\ell\eta)$ for
sites within distance $\ell$ of each other, and decays polynomially, with
a power $D$ that can be chosen as large as needed, at larger distances.

\begin{proof}
We verify, one after the other, the conditions
in \cite[Definitions 2.1, 2.4 and 2.5]{ER}. The general conditions on the
profile, its spectrum and the entry distribution are checked directly
(Step~1). The main condition, clause (i), asks that two deterministic
matrices built from the profile are bounded by the control function; we
compute them through the Green function of a lazy random walk on the ring
of blocks, which has an explicit formula (Step~2). Clauses (ii) and (iii)
are elementary properties of the control function (Step~3). Clauses (iv)
and (v) ask that the control function behaves well under products and
sums over an intermediate site; these follow from the triangle inequality
for the circular distance and elementary row sums (Step~4). Clause (vi)
is a regularity estimate in position (Step~5), and clause (vii) concerns
the constant Fourier mode of the blocks (Step~6). Step~7 applies the local
law. Throughout we write $d(x,y)=|x-y|_\Naux$.

\emph{Step 1: the general conditions.} The matrix $S$ is symmetric and
nonnegative. Each of its rows sums to $c+2b=1$, so it is stochastic, and
its entries are at most $1/W$. Its nonzero eigenvalues are
$1-4b\sin^2(\pi j/m)$, $0\le j<m$, which lie in $[1-4b,1]$. So the
spectrum of $S$ stays a fixed distance above $-1$ for every fixed
$b<1/2$, and no assumption on the parity of $m$ is needed. The normalized
entries are independent centered Gaussians, so all their moments are
bounded. The complex off-diagonal entries have zero second moment, and a
real diagonal is allowed by the covariance convention of
\cite[Section 5.1.1]{ER}. The band-width condition
$W^2\ge \Naux^{1+\zeta_0}$ is not part of admissibility; it is imposed
only in the short-ring theorem. Finally, the support condition of
\cite[Example 2.3]{ER} is stated with an inequality in the reverse
direction: finite range means vanishing at block distance \emph{greater}
than $K_1$. Our range is one, and we verify admissibility directly.

All bounds obtained below in clauses (i) and (vi) decay exponentially in
the distance. Since $e^{-ar}\le C_{a,D}(1+r)^{-D}$ for $r\ge0$, every such
exponential bound is also bounded by a constant multiple of the
polynomial control function \eqref{eq:t26-er-control}.

\emph{Step 2: clause (i).} Clause (i) concerns the matrices
\[
 \Theta(z_1,z_2)=\bigl((I-aC)^{-1}-I\bigr)\otimes J_W
\]
and its analogue $\Xi$, where $a$ is a product of the semicircle
Stieltjes transforms at the two spectral parameters. Here $(I-aC)^{-1}$
is the Green function of a lazy random walk on the ring of blocks, killed
at rate $1-a$. For $0<a<1$ it has the explicit form
\begin{equation}\label{eq:killed-walk}
 (I-aC)^{-1}_{0d}
 =\frac{q^d+q^{m-d}}{\sqrt{(1-ac)^2-4a^2b^2}(1-q^m)},\qquad
 q=\frac{1-ac-\sqrt{(1-ac)^2-4a^2b^2}}{2ab},
\end{equation}
where $d$ is the block distance and the branch is fixed by $q(0)=0$. As
real $a$ increases to one, both $1-q$ and the square root are of order
$\sqrt{1-a}$. For complex $|a|<1$, expanding in powers of $a$ by the
hitting times of the walk gives $|q(a)|\le q(|a|)$, and the power series
of the Green function is bounded entrywise by its value at $|a|$.

For $\Theta$ we take $a=m(z_1)\overline{m(z_2)}$. With
$\eta=\min_j|\Im z_j|$ we have $1-|a|\ge c\eta$, so \eqref{eq:killed-walk}
gives $|\Theta_{xy}|\le CU_\eta(x,y)$. Concretely, above the scale
$\eta\asymp m^{-2}$ the bound is
$C(W\sqrt\eta)^{-1}e^{-cd(x,y)\sqrt\eta/W}$, and below it the bound is
$C/(\Naux\eta)$. Moving a site within its block changes distances by at
most $2W$, which is harmless since $W\le C\ell$.

For $\Xi$ we take $a=m(z_1)m(z_2)$, with both parameters in one
half-plane. The phases of the two Stieltjes transforms avoid $0$ and
$\pi$ in the bulk, so their product avoids the positive real direction.
Together with $1-4b>-1$ this gives
$\inf_{\lambda\in[1-4b,1]}|1-a\lambda|\ge c_b>0$. Compactness of the
closed bulk domain then gives $|q(a)|\le q_b<1$, and therefore
$|\Xi_{xy}|\le CW^{-1}e^{-cd(x,y)/W}\le CU_1(x,y)$. This proves (i).

\emph{Step 3: clauses (ii) and (iii).} Put $F_s(r)=(1+r/s)^{-D/2}$, so that
$U_\eta^{1/2}=(\ell\eta)^{-1/2}F_\ell(d)$. We use two elementary facts
about the circular distance: it satisfies the triangle inequality, and a
ball of radius $R$ contains at most $2R+1$ sites. Consequently, for every
fixed $\sigma>1$,
\[
 \sum_a(1+d(x,a)/\ell)^{-\sigma}\le C_\sigma\ell,
 \qquad
 \sum_{d(a,y)<R}(1+d(a,y)/L_0)^{-\sigma}
       \le\min(2R+1,C_\sigma L_0).
\]
The first sum shows that the row sum of $U_\eta$ is at most $C/\eta$, and
clearly $\max U_\eta\le(\ell\eta)^{-1}$. For a lower bound, $\ell\eta\le
C_0\Naux$ and $1+d/\ell\le\Naux$, so $\min U_\eta\ge
C_0^{-1}\Naux^{-D-1}$; this is the required fixed polynomial lower bound,
after increasing the width threshold. For $\eta\le(W/\Naux)^2$ we have
$\ell=\Naux$, and since $d\le\Naux/2$, $U_\eta\asymp(\Naux\eta)^{-1}$.
Finally, if $\eta_1\le\eta_2$ then $\ell_1\ge\ell_2$ and
$\ell_1\eta_1\le\ell_2\eta_2$, so $U_\eta$ increases as $\eta$ decreases.
These facts are clauses (ii) and (iii).

\emph{Step 4: clauses (iv) and (v).} These clauses compare control
functions at two scales $L_0=\ell(\eta_1)\ge l=\ell(\eta_2)$. With
$h\in\{L_0,l\}$, we prove the three estimates
\begin{align}
 \sup_a F_l(d(x,a))^2F_{L_0}(d(a,y))^2
 &\le F_{L_0}(d(x,y))^2,\label{eq:t26-er-product}\\
 \sum_aF_l(d(x,a))F_{L_0}(d(a,y))
 &\le ClF_{L_0}(d(x,y)),\label{eq:t26-er-conv}\\
 \sum_a\frac{(d(x,a)+W)\wedge L_0}{L_0}
        F_l(d(x,a))F_h(d(a,y))
 &\le C\frac{l^2}{L_0}F_{L_0}(d(x,y)),\label{eq:t26-er-weighted}
\end{align}
with constants independent of the scales. The first is a product bound,
the second a convolution bound, and the third a convolution bound with a
weight that is small when the intermediate site is close to $x$.

The product bound \eqref{eq:t26-er-product} follows from the triangle
inequality in the form
$(1+d(x,a)/l)(1+d(a,y)/L_0)\ge1+d(x,y)/L_0$, using $l\le L_0$.
For the two sums put $s=D/2\ge3$ and $\delta=d(x,y)$, and distinguish two
cases.
If $\delta\le2l$, then $F_{L_0}(\delta)$ is bounded below by a constant,
and the convolution bound follows from the row sum of Step~3. For the
weighted bound we use
\[
 \frac{(d(x,a)+W)\wedge L_0}{L_0}
 \le\min\{1,C(l/L_0)(1+d(x,a)/l)\}
\]
together with the row sum with exponent $s-1>1$.

If $\delta>2l$, we split the sum at $d(a,y)=\delta/2$. On the part where
$d(a,y)\ge\delta/2$, we extract $F_h(d(a,y))\le2^sF_{L_0}(\delta)$ and use
the same row sums. On the part where $d(a,y)<\delta/2$, the site $a$ is
far from $x$. We put $r=L_0/l$, $\widehat x=\delta/l$ and $t=\widehat x/r$.
Then $F_l(d(x,a))\le(1+\widehat x/2)^{-s}$, the weight is at most
$\min\{1,C(1+\widehat x)/r\}$, and the remaining sum over $a$ is at most
$Cl\min(\widehat x,r)$ for $h=L_0$ and at most $Cl$ for $h=l$. Dividing by
the right sides of \eqref{eq:t26-er-conv} and \eqref{eq:t26-er-weighted},
we must bound the resulting ratios. For $t\le1$ they are at most constant
multiples of $\widehat x(1+\widehat x/2)^{-s}$,
$(1+\widehat x)\widehat x(1+\widehat x/2)^{-s}$ or
$(1+\widehat x)(1+\widehat x/2)^{-s}$. For $t>1$ they are at most
$C4^sr^{1-s}$ or $C4^sr^{2-s}$. All of these are bounded because $s\ge3$.

This proves the three estimates. Multiplying them by the amplitudes
$(\ell\eta)^{-1/2}$ that we left out gives (2.23)--(2.25) of \cite{ER}.
In the weighted case the remaining factor is exactly $l^2/L_0$, so no
ratio of scales is lost. These are clauses (iv) and (v).

\emph{Step 5: clause (vi), localized regularity.} This clause asks that
$\Theta$ changes little when one site is moved by a short distance. Assume
$\eta\ge c_2m^{-2}$. Then $|1-q^m|\ge c(c_2,b)>0$, and
\[
 \frac{1-q}{\sqrt{(1-ac)^2-4a^2b^2}}
       =\frac{q}{ab(1+q)}
\]
is uniformly bounded for $|a|\le1$, with its removable value at $a=0$.
Hence the difference of the Green function at neighbouring blocks is at
most $Ce^{-cd\sqrt\eta}/W$ after tensoring with $J_W$. For a displacement
of $r\le\ell/W$ blocks we add up the differences along the way, at the
cost of a fixed exponential factor, since intermediate distances differ
from the endpoint distances by at most $r$. For larger $r$ we use the
bound without differences. The term $-I\otimes J_W$ lives on one block
and obeys the same bound with the allowance $W$. Returning to site
distances, and using the remark at the end of Step~1, we obtain
\[
 |\Theta_{xy}-\Theta_{xz}|
 \le C\frac{(d(y,z)+W)\wedge\ell}{\ell}
              \{U_\eta(x,y)+U_\eta(x,z)\},
\]
including the transition region $\ell=\Naux$ under the lower bound on
$\eta$. This is clause (vi).

\emph{Step 6: clause (vii), the constant block Fourier mode.} This clause
asks that each row of $\Theta$ is close to its average. In the Fourier
basis of the ring of blocks, $C$ has eigenvalues
$\widehat C_j=1-4b\sin^2(\pi j/m)$, and
\[
 \max_{x,y}\left|\Theta_{xy}-\frac{1}\Naux \sum_z\Theta_{xz}\right|
 \le\frac{1}{mW}\sum_{j=1}^{m-1}
       \left|\frac{a\widehat C_j}{1-a\widehat C_j}\right|
 \le\frac{C_b}{mW}\sum_{j=1}^{m-1}
       \{1+\sin^{-2}(\pi j/m)\}
 \le C_b\frac{\Naux }{W^2}.
\]
For $\widehat C_j\ge0$ we used $|1-a\widehat C_j|\ge1-\widehat C_j$, and
for $\widehat C_j<0$ the fixed gap above $-1$ from Step~1. This proves
(vii) on the whole spectral domain.

\emph{Step 7: conclusion.} All admissibility constants are now fixed and
independent of $m$ and $W$. Applying \cite[Theorem 3.3, (3.17)]{ER} with
$u=e_x$, $v=e_y$ and the control function \eqref{eq:t26-er-control} gives
\eqref{eq:short-ring-local}.
\end{proof}

\begin{theorem}\label{thm:critical-density}
The normalized expected empirical measure has a smooth density
$\rho_{n,W}$. For fixed $b,\kappa,T$, uniformly for
$4\lceil W^{15/16}\rceil\le n\le TW$ and $|e|\le2-\kappa$,
\begin{equation}\label{eq:critical-density}
 |\rho_{n,W}(e)-\rho(e)|\le C_{b,\kappa,T}W^{-1/32}.
\end{equation}
Consequently every critical subsequential marked limit has unit
spectral intensity.
\end{theorem}

Only the mean density is identified here at the critical scale. The shorter-ring local law is used within its fixed power range, then transferred to the original ring.

\begin{proof}
We first control unsmoothing. Split every diagonal Gaussian block into
two independent halves. From the second half extract its common scalar
  $t=N^{-1}\Tr A^{(2)}$, of variance
  $\sigma^2=c/(2NW)$. Gaussian orthogonality makes $t$ independent of
the remaining matrix $H_*$. The first half of the diagonal randomness
still gives an expected density bounded by $\sqrt2K_b$.
Convolution with the density of $t$ therefore shows

\begin{equation}\label{eq:density-smooth}
 \norm{\rho_{n,W}}_\infty\le K_b,\qquad
 \norm{\rho_{n,W}'}_\infty\le C_b\sqrt{NW},\qquad
 \norm{\rho_{n,W}-p_\eta*\rho_{n,W}}_\infty
 \le C_b\eta\sqrt{NW}\log(1/\eta)+C_b\eta.
\end{equation}

For the last inequality integrate the derivative bound over $|t|\le1$
against $p_\eta(t)$ and the uniform density bound over $|t|>1$.

Choose $m=\lceil W^{15/16}\rceil$,   $\Naux =mW$  and $\eta=W^{-7/4}$.
Around a site $x$ retain $m$ consecutive blocks with $x$ at distance
at least $m/3$ blocks from either endpoint. Close this segment with
an independent Gaussian bond to obtain a matrix $K$ with exactly the
short-ring law. For large $W$,

\[
 W^2\ge \Naux ^{1+1/128},\qquad
 \eta\ge \Naux ^{-1+1/256},\qquad
 \ell=W^{15/8},\quad \ell\eta=W^{1/8},\quad
 \Naux /\ell\asymp W^{1/16}.
\]

Theorem~\ref{thm:external-band-local} gives
$|G^K_{xx}-m_{\rm sc}|\prec W^{-1/16}$ and, at either endpoint block,
$\max_{y\in\partial K}|G^K_{xy}|\prec W^{-A}$ for any fixed $A$. Indeed $d(x,y)\ge\Naux/4$ for large $W$, and the corrected control is at most $CW^{-1/16-D/32}$. Given $A$, choose $D\ge\max(6,32A)$. The density theorem and its unit-intensity conclusion are unchanged.

On the enlarged probability space set
$H_0=K\oplus H_{K^cK^c}$ and $V=H-H_0$. The row of $G^{H_0}$
at $x$ vanishes outside $K$, and only its $2W$ endpoint columns
meet $V$. The resolvent identity therefore gives
\[
 |G^H_{xx}-G^K_{xx}|
 \le\sqrt{2W}\max_{y\in\partial K}|G^K_{xy}|\norm V\,\eta^{-1}.
\]
There are only finitely many Gaussian blocks in $V$, so $\norm V\prec1$.
On the exceptional event use the deterministic bound $2/\eta$.
The expectation of this difference is smaller than every fixed power
of $W^{-1}$. Spending a fixed stochastic-domination exponent in the
diagonal estimate yields
$|\E G^H_{xx}-m_{\rm sc}|\le CW^{-1/32}$, uniformly in $x$.
Averaging the imaginary part gives this bound for
$p_\eta*\rho_{n,W}-\pi^{-1}\Im m_{\rm sc}$.
The semicircle smoothing error is $O_\kappa(\eta)$; by
\eqref{eq:density-smooth} the remaining error is
$O_{b,T}(W^{-1/4}\log W)$. This proves \eqref{eq:critical-density}.

For $f\in C_c(\R)$,

\[
 \E\sum_\alpha f(u_\alpha)
 =\int f(u)\frac{\rho_{n,W}(E+u/(N\rho(E)))}{\rho(E)}\,du
 \longrightarrow\int f(u)\,du.
\]

The exponential count bound passes this expectation through weak
convergence, proving the intensity assertion.
\end{proof}

\subsection{Translation and rotation invariance}\label{subsec:s3-translation-and-rotation-invariance}

 We show that no spectral origin or spatial block is distinguished in
the limiting marked law. Proposition~\ref{prop:stationarity} gives
spectral translation invariance, spatial rotation invariance, and mean
intensity of the limiting process. We obtain the first from an independent scalar
Gaussian shift and the others from cyclic symmetry and the mean-density
estimate.

 \begin{proposition}\label{prop:stationarity}
Every critical subsequential marked law is invariant under spectral
translations and simultaneous rotations of its circle marks. Its mean
spatial intensity is $du\,dr$.
\end{proposition}

The limiting spectral translation symmetry comes from an exact scalar Gaussian decomposition. Block rotation symmetry simultaneously fixes the mean spatial mark.

\begin{proof}
Let   $t=N^{-1}\Tr H$  and $H_0=H-tI$. Then
  $t\sim N(0,c/(NW))$  is independent of $H_0$. The marked spectrum of
$H$ is that of $H_0$ translated by the independent Gaussian
  $N\rho(E)t$, whose variance is $c\rho(E)^2n$.
If $\theta_s$ translates the spectral coordinate by $s$, contraction
of total variation under a measurable map gives
\begin{equation}\label{eq:stationarity-tv}
 d_{\rm TV}(\Law(\Xi_W),\Law(\theta_s\Xi_W))
 \le \min\left(1,\frac{|s|}{\rho(E)\sqrt{2\pi cn}}\right).
\end{equation}
Indeed the distance between $N(0,\sigma^2)$ and $N(s,\sigma^2)$
equals $2\Phi(|s|/(2\sigma))-1\le |s|/(\sqrt{2\pi}\sigma)$.
At critical bandwidth $n\to\infty$, so the limiting law is stationary.

Cyclic block invariance preserves each ordered eigenvalue and gives
\[
 \E\sum_\alpha f(u_\alpha)p_{\alpha k}
 =\frac{1}n\E\sum_\alpha f(u_\alpha).
\]
The uniform grid measure on $\T$ converges to Lebesgue measure; use
Theorem~\ref{thm:critical-density} and the count bounds to pass to
the mean marked intensity. Approximate any fixed circle rotation by a
multiple of $1/n$. Continuity of rotation on $\mathcal P(\T)$ then
passes the finite cyclic invariance to the limit.
\end{proof}

\part{} 

\section{An exact change of variables at a prescribed rank}
\label{sec:coarea}
The goal of this section is an exact formula for the law of the eigenpair
$(\lambda_\alpha,\psi_\alpha)$ at a deterministic rank $\alpha$. It is the first
step of Part~I of Subsection~\ref{subsec:proof-ideas}. It replaces estimates
at a fixed energy, which do not see the rank, by a formula that keeps both
the law of the energy and the index of the selected eigenvalue.
Theorem~\ref{thm:coarea} says the following. We fix every entry of the
matrix except a fresh GUE part on two adjacent blocks. Then the eigenvalue of
rank $\alpha$, together with the direction $v$ of its eigenvector on the pair,
has an explicit joint density. This density consists of a squared
characteristic polynomial, an eigenvalue-counting indicator that fixes the
rank, and a Jacobian equal to the inverse of the mass of the eigenvector on
the pair. Under the law size-biased by this mass, the pair law
\eqref{eq:pair-law}, the Jacobian cancels exactly.
In Subsection~\ref{subsec:coarea-statement} we expose the Gaussian pair, define its Schur
complement and state the formula. In Subsection~\ref{subsec:jacobian} we compute the
Jacobian of eigenpair coordinates, which gives the squared determinant and the
normalizing constant. In Subsection~\ref{subsec:supported-source} we show how a local inverse and
a restored Gaussian row can be used inside the formula without dropping the
weight that selects the rank. In Subsection~\ref{subsec:coarea-proof} we prove the theorem.

The formula is obtained in three moves. We expose an independent Gaussian
matrix on two adjacent blocks. We eliminate all other blocks by a Schur
complement, so that the eigenvalue equation becomes an equation on the
pair. We then change variables from the exposed Gaussian matrix to the
selected eigenpair. The Jacobian of this change contains the speed at
which an eigenvalue of the Schur complement crosses zero, and this speed
is exactly the inverse mass of the exposed pair. This is why the QUE
argument uses the pair-mass law \eqref{eq:pair-law}, under which the
crossing speed cancels.

The section is organized as follows. Theorem~\ref{thm:coarea} is stated in
Subsection~\ref{subsec:coarea-statement} and proved in
Subsection~\ref{subsec:coarea-proof}. Lemma~\ref{lem:eigenpair-jacobian}
in Subsection~\ref{subsec:jacobian} computes the Jacobian of the eigenpair
coordinates of a Hermitian matrix; it gives the constant and the squared
determinant in Theorem~\ref{thm:coarea}. Lemmas~\ref{lem:supported-source}
and~\ref{lem:row-restoration} in Subsection~\ref{subsec:supported-source}
are not used in the proof of Theorem~\ref{thm:coarea}. They show how to
use a local inverse inside the formula \eqref{eq:complete-coarea} without
dropping the weight that selects the rank. The pair-mass law
\eqref{eq:pair-law} is the law $\Prob^\sharp$ of
Theorems~\ref{thm:physical-center} and~\ref{thm:physical-angular}, and the
formula \eqref{eq:complete-coarea} is used again in
Section~\ref{sec:determinant-score}.

\subsection{The formula}\label{subsec:coarea-statement}

\begin{figure}[H]
\centering
\begin{tikzpicture}[x=1cm,y=1cm]
 \begin{scope}
  \foreach \k in {0,...,9}{\draw[black!40,line width=1.2pt] (36*\k:1.5)--(36*\k+36:1.5);}
  \foreach \k in {0,...,9}{
   \ifnum\k<2
    \node[circle,draw,fill=black!25,minimum size=7mm,inner sep=0pt] at (36*\k+72:1.5) {};
   \else
    \node[circle,draw,fill=white,minimum size=7mm,inner sep=0pt] at (36*\k+72:1.5) {};
   \fi}
  \node[font=\small] at (90:.45) {pair $S$};
  \node[font=\small] at (0,-.35) {exterior $R$};
  \node[font=\scriptsize] at (72:2.15) {$i$};
  \node[font=\scriptsize] at (108:2.15) {$i+1$};
 \end{scope}
 \begin{scope}[shift={(5.2,-1.4)}]
  \draw (0,0) rectangle (3,3);
  \fill[black!25] (0,2.1) rectangle (.9,3);
  \draw (.9,0)--(.9,3); \draw (0,2.1)--(3,2.1);
  \node at (.45,2.55) {$H_{SS}$};
  \node at (1.95,2.55) {$C^*$};
  \node at (.45,1.05) {$C$};
  \node at (1.95,1.05) {$H_R$};
  \node[font=\small,anchor=west] at (3.2,2.55) {$H_{SS}=V+G$};
  \node[font=\small,anchor=west] at (3.2,1.05) {$\omega=(H_R,C,V)$};
 \end{scope}
\end{tikzpicture}
\caption{The splitting used in Theorem~\ref{thm:coarea}. Left: the ring of
blocks, with the exposed pair $S=\{i,i+1\}$ shaded and the remaining blocks
$R$ white. Right: the matrix $H$ in the decomposition $S\oplus R$. The
shaded corner $H_{SS}$ is split as $V+G$ by \eqref{eq:pair-splitting}, with
$G$ a fresh GUE matrix independent of $\omega=(H_R,C,V)$. In the proof,
$\omega$ is held fixed and only $G$ is changed into eigenpair coordinates;
eliminating the unshaded part gives the Schur complement $A$ of
\eqref{eq:schur-data}.}
\label{fig:coarea-split}
\end{figure}
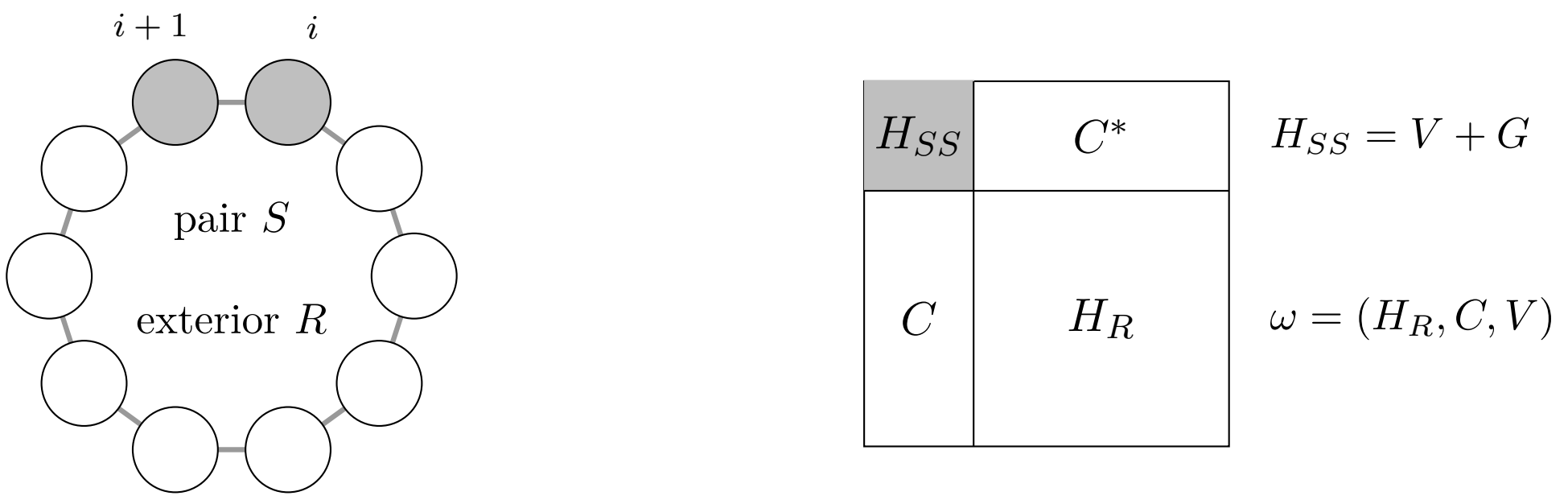

  We write the selected eigenpair law in coordinates on two adjacent
blocks. Theorem~\ref{thm:coarea} gives an exact density that retains both
the selected rank and the mass needed to reconstruct the full vector.
We expose the Gaussian pair, define its Schur complement and
reconstruction map, and then state the formula illustrated in
 Figure~\ref{fig:coarea-split}.

\paragraph{The exposed pair.} Fix an adjacent pair of blocks
$S=\{i,i+1\}$, write $R=S^c$ for the remaining blocks and $d=2W$ for the
dimension of the pair, and choose $0<s\le\min(c,b)$. Put $t=s/d$. There is
an exact independent splitting
\begin{equation}\label{eq:pair-splitting}
 H_{SS}=V+G,\qquad G\sim\GUE_d(t),
\end{equation}
where $G$ is independent of $V$ and of all entries of $H$ outside $H_{SS}$.
All remaining variances are nonnegative. Write $C=H_{RS}$ for the
coupling between $R$ and $S$, and $\omega=(H_R,C,V)$ for all randomness
except $G$.

\paragraph{The Schur complement.} For $E\notin\operatorname{Spec}H_R$ set
\begin{equation}\label{eq:schur-data}
\begin{aligned}
 F(E)&=V-C^*(H_R-E)^{-1}C,& T(E)&=(H_R-E)^{-1}C,\\
 \mathcal D(E)&=I+T(E)^*T(E),& A&=F(E)+G-E.
\end{aligned}
\end{equation}
The matrix $A=A(E)$ is the Schur complement of $H-E$ on the pair $S$. The
\emph{reconstruction map} is $\mathcal L_Ev=(v,-T(E)v)$ in the
decomposition $S\oplus R$; it extends a vector on the pair to a vector on
the whole ring. Let $N_-$ count strictly negative eigenvalues, and put
$N_R(E)=N_-(H_R-E)$.

The next theorem gives the law of the selected eigenpair, written as an
integral over $\omega$, the energy $E$ and the direction $v$ of the
eigenvector on the pair. It also gives the same law under the pair-mass
weight.

\begin{theorem}[Rank-preserving eigenpair formula]\label{thm:coarea}
Let $E=\lambda_\alpha$, $m=p_{\alpha i}+p_{\alpha,i+1}$ and
$v=P_S\psi_\alpha/\sqrt m$. Then
\begin{equation}\label{eq:clock-reconstruction}
 Av=0,\qquad m^{-1}=v^*\mathcal Dv,\qquad
 \psi_\alpha=\frac{\mathcal L_Ev}{\norm{\mathcal L_Ev}}.
\end{equation}
Let $\nu_d$ be normalized projective Haar measure and $U_v$ an
orthonormal frame of $v^\perp$. Define
\[
 F_v=U_v^*FU_v,\quad B=F_v+G_{d-1},\quad M=B-E,\qquad
 G_{d-1}\sim\GUE_{d-1}(t),
\]
and
\begin{equation}\label{eq:coarea-regression}
 \mathcal R(E,v,F)=
 \exp\left\{-\frac{(E-v^*Fv)^2+
                    2\norm{(I-vv^*)Fv}^2}{2t}\right\},\qquad
 c_{d,t}=\frac{1}{(d-1)!\,t^{d-1}\sqrt{2\pi t}}.
\end{equation}
For every nonnegative measurable observable $\Phi$ of the
reconstructed matrix and its selected eigenpair,
\begin{align}
 \E_\alpha\Phi
 =\int d\Prob(\omega)\,dE\,d\nu_d(v)\,
 c_{d,t}\mathcal R(E,v,F)(v^*\mathcal Dv)
 \E_{G_{d-1}}\!\left[
 \det(M)^2\one_{\{N_-(M)=\alpha-N_R(E)-1\}}\Phi\right].
 \label{eq:complete-coarea}
\end{align}
This measure has total mass one. The pair-mass law
\begin{equation}\label{eq:pair-law}
 d\Prob^\sharp_{\alpha,S}
 =\frac{n}{2}(p_{\alpha i}+p_{\alpha,i+1})\,d\Prob
\end{equation}
also has total mass one; its formula is
\eqref{eq:complete-coarea} with $v^*\mathcal Dv$ replaced by $n/2$.
\end{theorem}

Theorem~\ref{thm:coarea} claims an exact identity, with no error term, for
the law of the eigenpair at the rank $\alpha$. The density in
\eqref{eq:complete-coarea} has four factors: the Gaussian regression factor
$c_{d,t}\mathcal R$, the crossing speed $v^*\mathcal Dv$, the squared
determinant $\det(M)^2$, and the inertia indicator that fixes the rank.
Under the pair-mass law \eqref{eq:pair-law} the crossing speed is replaced
by the constant $n/2$.  All later expectations $\E^\sharp_{\alpha,S}$ and
probabilities $\Prob^\sharp_{\alpha,S}$, also with an edge $e$ in place of
$S$, are taken under this size-biased physical law, with density
$(n/2)(p_{\alpha i}+p_{\alpha,i+1})$, also where \eqref{eq:pair-law} is
not cited again. In particular $d\Prob^\sharp_{\alpha,S}\le(n/2)\,d\Prob$.

The inertia indicator fixes the original rank even when the exterior has
eigenvalues near $E$. The density includes the complete Gaussian
regression and determinant factors. In particular, conditioning the
physical row in this formula does not make that row Gaussian.

\begin{remark}[The choice of $s$ and $t$]\label{rem:coarea-parameters}
The bound $s\le\min(c,b)$ is what makes the splitting
\eqref{eq:pair-splitting} exact. A diagonal block of the pair has entry
variance $c/W$ and the bond inside the pair has entry variance $b/W$,
while $G\sim\GUE_d(t)$ removes variance $t=s/(2W)\le\min(c,b)/W$ from each
entry of $H_{SS}$. Hence $V$ has nonnegative variances. Any positive $s$ in this range gives
the formula. Section~\ref{sec:clock-center} uses the same kind of
splitting simultaneously on every edge, with $s=\min(c,b)$.
\end{remark}

\subsection{The Jacobian of eigenpair coordinates}\label{subsec:jacobian}

In this subsection we compute the Jacobian of the map from a Hermitian
matrix to one of its eigenpairs and the remaining matrix on the orthogonal
complement. It gives the squared determinant and the constant $c_{d,t}$ in
Theorem~\ref{thm:coarea}, and it fixes the normalization of projective
measure used later.

\begin{lemma}\label{lem:eigenpair-jacobian}
For a Hermitian $d\times d$ matrix $K$, a selected eigenpair may be
written as $K=Evv^*+U_vBU_v^*$. Its Lebesgue change of variables,
summed over the $d$ eigenpairs, is
\begin{equation}\label{eq:eigenpair-jacobian}
 dK=\frac{\pi^{d-1}}{(d-1)!}
       \det(B-E)^2\,dE\,dB\,d\nu_d(v).
\end{equation}
Here the equality means the integral identity for a sum over selected
eigenpairs; it does not count that sum twice and $dK,dB$ use the diagonal real coordinates and the real and imaginary parts of strict-upper entries, without a $\sqrt2$ rescaling; $\nu_d$ has mass one. 
\end{lemma}

The lemma says that the volume element of $K$ factors into the eigenvalue
$E$, the remaining matrix $B$ on $v^\perp$, and the direction $v$, with
density $\det(B-E)^2$. The measure $\nu_d$ on eigenvector directions has
total mass one, and the constant $\pi^{d-1}/(d-1)!$ is the volume of the
space of directions.

\begin{remark}[Normalization]\label{rem:jacobian-normalization}
The convention that $\nu_d$ has total mass one fixes the constant in the
physical formula \eqref{eq:complete-coarea}. With another normalization of
projective measure, the constant $c_{d,t}$ changes by the same factor.
\end{remark}

\begin{proof}
The challenge is to identify both the density and the constant. We
compute the differential at one point, and then the volume of the space
of directions.

\emph{The density.} By unitary invariance it is enough to take $v=e_1$.
Use complex coordinates $z\in\C^{d-1}$ for the eigenvector direction near
$e_1$. Differentiation of the off-diagonal block of $K$ gives
\[
 dK_{\rm off}=(E-B)\,dz.
\]
The real Jacobian of this complex linear map is $|\det(E-B)|^2$. The real
diagonal coordinate along $v$ is $dE$, and the remaining Hermitian
coordinates are $dB$. The differential is triangular in these variables,
so
\[
 dK=|\det(E-B)|^2\,dE\,dB\,d\mathrm{vol}(v),
\]
where $d\mathrm{vol}$ is the invariant volume on $\C P^{d-1}$ in this
coordinate convention. Since $B-E$ is Hermitian,
$|\det(E-B)|^2=\det(B-E)^2$.

\emph{The constant.} The invariant volume of $\C P^{d-1}$ for this
coordinate convention is $\pi^{d-1}/(d-1)!$. To see the constant directly,
the projective chart $v=(1,z)/\sqrt{1+\norm z^2}$ has volume density
\[
 (1+\norm z^2)^{-d}\,d^{2d-2}z,
\]
whose integral is $\pi^{d-1}/(d-1)!$ by polar integration.
Dividing that volume by its total mass gives $\nu_d$ and proves
\eqref{eq:eigenpair-jacobian}.
\end{proof}

\subsection{A supported source and a restored Gaussian row}\label{subsec:supported-source}

In this subsection we show how a local inverse can be used inside the
formula \eqref{eq:complete-coarea} without dropping the weight that selects
the rank. Lemma~\ref{lem:supported-source} writes a quadratic form of the
pseudoinverse of $H-E$ as a quadratic form of the pseudoinverse of the
Schur complement $A$. Lemma~\ref{lem:row-restoration} then bounds moments of
that quadratic form by integrating out one row and column of $M$.

Put $z=p_{\alpha i}/m$ and $J=P_{i+1}-P_i$ on the pair. The number $z$ is
the fraction of the pair mass that lies in block $i$.
Here $J$ is the negative of the contrast $J_e=P_i-P_{i+1}$ used
in the scalar-center and angular sections. This choice fixes the sign
of the displayed supported vector $f$; its norm and pseudoinverse
quadratic form are unchanged when $f$ is replaced by $-f$.

The next lemma introduces a vector supported on the pair and orthogonal to
$\psi_\alpha$, and computes its norm and its quadratic form.

\begin{lemma}\label{lem:supported-source}
The vector
\[
 f=(p_{\alpha i}p_{\alpha,i+1})^{1/4}
       \left(\frac{P_{i+1}\psi_\alpha}{p_{\alpha,i+1}}
                    -\frac{P_i\psi_\alpha}{p_{\alpha i}}\right)
 =\frac{(I-vv^*)Jv}{2[z(1-z)]^{3/4}}
\]
is supported on $S$ and orthogonal to $\psi_\alpha$. It satisfies
\begin{equation}\label{eq:supported-inverse}
 \norm f^2=[z(1-z)]^{-1/2},\qquad
 f^*(H-E)^+f=f^*A^+f.
\end{equation}
The superscript $+$ denotes the inverse on the orthogonal complement
of the kernel.
\end{lemma}

The lemma claims that the quadratic form of the pseudoinverse of $H-E$
at the supported vector $f$ equals the quadratic form of the pseudoinverse
of the Schur complement $A$, although the two kernels and their orthogonal
complements differ. This lets the later estimates work with the
$d\times d$ matrix $A$ only.

\begin{remark}[Why $f$ is orthogonal to $\psi_\alpha$]\label{rem:supported-orthogonal}
Orthogonality of $f$ to $\psi_\alpha$ removes the arbitrary kernel component
of the solution of $(H-E)y=f$. With the powers $1/4$ and $3/4$ in the
definition of $f$, the norm $\norm f^2=[z(1-z)]^{-1/2}$ depends only on the
ratio $z$.
\end{remark}

\begin{proof}
We prove the properties of $f$ first and the identity of quadratic forms
second.

\emph{Properties of $f$.} The vector $f$ is supported on $S$ because
$P_i\psi_\alpha$ and $P_{i+1}\psi_\alpha$ are. Its inner product with
$\psi_\alpha$ is
\[
 \psi_\alpha^*f=(p_{\alpha i}p_{\alpha,i+1})^{1/4}
 \left(\frac{p_{\alpha,i+1}}{p_{\alpha,i+1}}
       -\frac{p_{\alpha i}}{p_{\alpha i}}\right)=0 .
\]
Since $v=P_S\psi_\alpha/\sqrt m$ and $z=p_{\alpha i}/m$,
\[
 v^*Jv=\frac{p_{\alpha,i+1}-p_{\alpha i}}{m}=1-2z.
\]
Direct subtraction of the component of $Jv$ parallel to $v$ then proves
the expression of $f$ through $(I-vv^*)Jv$. For the norm,
\[
 \norm f^2=(p_{\alpha i}p_{\alpha,i+1})^{1/2}
 \left(\frac{1}{p_{\alpha,i+1}}+\frac{1}{p_{\alpha i}}\right)
 =\frac{m}{(p_{\alpha i}p_{\alpha,i+1})^{1/2}}
 =[z(1-z)]^{-1/2},
\]
since $p_{\alpha i}p_{\alpha,i+1}=m^2z(1-z)$.

\emph{The quadratic form.} Solve $(H-E)y=f$ with $y\perp\psi_\alpha$; this
$y$ is $(H-E)^+f$. The exterior equation, as in
\eqref{eq:block-equations} below, gives $y_R=-Ty_S$. The interior equation
then gives
\[
 Ay_S=f,\qquad\text{so}\qquad y_S=A^+f+av
\]
for some scalar $a$, since the kernel of $A$ is spanned by $v$. The
constraint $y\perp\psi_\alpha$ fixes $a$, but $f^*v=0$ removes this term
from $f^*y$. Since $f$ is supported on $S$, $f^*(H-E)^+f=f^*y=f^*y_S=f^*A^+f$.
This proves \eqref{eq:supported-inverse}.
\end{proof}

\begin{figure}[H]
\centering
\begin{tikzpicture}[x=1cm,y=1cm]
 \draw (0,0) rectangle (3,3);
 \fill[black!20] (0,2.4) rectangle (3,3);
 \fill[black!20] (0,0) rectangle (.6,2.4);
 \draw (.6,0)--(.6,3); \draw (0,2.4)--(3,2.4);
 \node at (.3,2.7) {$x$};
 \node at (1.8,2.7) {$y^*$};
 \node at (.3,1.2) {$y$};
 \node at (1.8,1.2) {$B_0$};
 \node[font=\small,anchor=west] at (3.3,2.7) {row along $w=U_v^*f/\norm f$};
 \node[font=\small,anchor=west] at (3.3,1.2) {$B_0$ fixed during the row integration};
 \node[font=\small,anchor=west] at (3.3,.5) {$h=x-y^*B_0^{-1}y$};
\end{tikzpicture}
\caption{The split \eqref{eq:row-split} used in
Lemma~\ref{lem:row-restoration}. The shaded row and column are the
coordinates $(x,y)$ of $M$ along the direction $w$ of the supported source
$f$; the unshaded block $B_0$ is the restriction of $M$ to the orthogonal
complement of $w$. In the proof, $(x,y)$ is integrated with $B_0$ fixed. The
scalar $h$ is the Schur complement of $B_0$ in $M$; it determines both
$\det M$ and the quadratic form $w^*M^{-1}w$.}
\label{fig:row-split}
\end{figure}
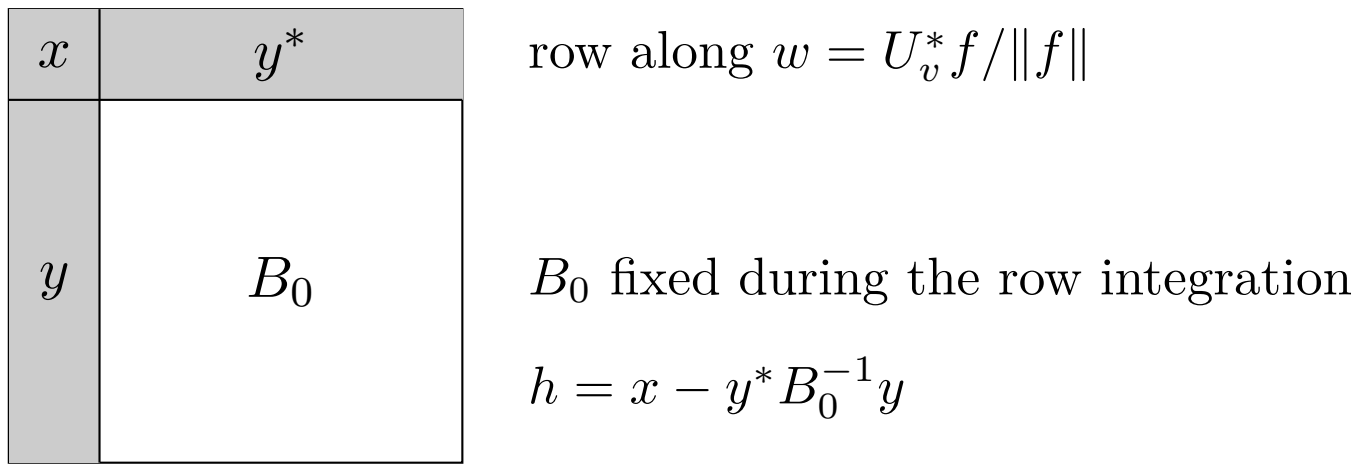

Fix $E,v,\omega$ in \eqref{eq:complete-coarea}, and split $M$ along the unit
vector $w=U_v^*f/\norm f$ in $\C^{d-1}$, as in
Figure~\ref{fig:row-split}:
\begin{equation}\label{eq:row-split}
 M=\begin{pmatrix}x&y^*\\y&B_0\end{pmatrix},\qquad
 x\sim N(\mu_0,t),\qquad y=y_0+\sqrt t\,g.
\end{equation}
Here $g$ is a standard circular Gaussian vector. These are independent
Gaussian \emph{prior} coordinates: their law is the law of the
corresponding entries of $M=F_v+G_{d-1}-E$ before the weight in
\eqref{eq:complete-coarea} is applied. Conditional on $B_0$ define
\begin{align}
 q^2={}&(\mu_0-y_0^*B_0^{-1}y_0-t\Tr B_0^{-1})^2+t
       +t^2\Tr B_0^{-2}+2t\norm{B_0^{-1}y_0}^2 .
 \label{eq:row-spread}
\end{align}

The next lemma bounds the moments of the quadratic form of
Lemma~\ref{lem:supported-source} by an expression in which the row $(x,y)$
has been integrated out.

\begin{lemma}\label{lem:row-restoration}
For $0<r\le2$ and any nonnegative observable $\Psi$ fixed during
the $(x,y)$ integration,
\begin{equation}\label{eq:row-restoration}
 \E_\alpha[\Psi|f^*(H-E)^+f|^r]
 \le\sum_{|\beta-\alpha|\le1}
       \E_\beta[\Psi\norm f^{2r}q^{-r}].
\end{equation}
The same statement holds under the respective pair-mass laws.
On the right, the observable is reconstructed using the same
$(E,v,\omega,B_0)$ coordinates.
\end{lemma}

The lemma replaces a negative power $|h|^{-r}$ of a random scalar by the
deterministic quantity $q^{-r}$, at the cost of summing over the three
ranks $\alpha-1,\alpha,\alpha+1$. Deleting a row is an integration device
within the same full density: its Gaussian regression, determinant and
inertia factors are all restored on the right side.

\begin{remark}[The range of $r$]\label{rem:row-range}
The condition $r\le2$ is what allows the concavity step in the proof: the
function $x\mapsto x^{(2-r)/2}$ is concave on $[0,\infty)$ exactly when
$r\le2$. For $r=2$ the bound uses no concavity.
\end{remark}

\begin{proof}
The challenge is that the left side of \eqref{eq:row-restoration} contains
a negative power of the random scalar $h=x-y^*B_0^{-1}y$, while the weight
$\det(M)^2$ and the inertia indicator also depend on $h$. We first express
all three quantities through $h$, then integrate $(x,y)$ with $B_0$ fixed,
and finally account for the inertia indicator.

\emph{Expressions through $h$.} Schur's formula for the split
\eqref{eq:row-split} gives
\[
 \det M=\det B_0\,h,\qquad w^*M^{-1}w=h^{-1}.
\]
By Lemma~\ref{lem:supported-source}, $f$ is supported on $S$ and
orthogonal to $v$, so $f=\norm f\,U_vw$, and the inverse of $A$ on
$v^\perp$ is $U_vM^{-1}U_v^*$. Hence
$f^*(H-E)^+f=\norm f^2h^{-1}$, and in \eqref{eq:complete-coarea}
\[
 \det(M)^2\,|f^*(H-E)^+f|^r=\det(B_0)^2\norm f^{2r}|h|^{2-r}.
\]

\emph{The row integral.} The Gaussian Wick identity \eqref{eq:gaussian-wick}
gives $\E_{x,y}h^2=q^2$, with precisely the four terms in
\eqref{eq:row-spread}: the square of the mean of $h$, the variance $t$ of
$x$, the variance $t^2\Tr B_0^{-2}$ of the quadratic term $t\,g^*B_0^{-1}g$,
and the variance $2t\norm{B_0^{-1}y_0}^2$ of the linear term. By concavity,
\[
 \E_{x,y}|h|^{2-r}\le(\E_{x,y}h^2)^{1-r/2}
                   =q^{-r}\E_{x,y}h^2 .
\]
The original row integral on the left of \eqref{eq:row-restoration},
weighted by the determinant, is bounded by this expression after dropping
its sign restriction on $h$.

\emph{The inertia indicator.} The sign of $h$ decides whether $M$ has one
more negative eigenvalue than $B_0$. Restoring the two signs of $h$
therefore changes $N_-(M)$ by at most one, and uses only the three ranks
$\beta$ with $|\beta-\alpha|\le1$ displayed in \eqref{eq:row-restoration}.
Both the crossing derivative and its replacement $n/2$ under the pair-mass
law are fixed in this row integration. All remaining factors, including
$\det B_0^2$, the regression density and their normalizers, are identical
on both sides. Integrating them proves the assertion.
\end{proof}

\subsection{Proof of Theorem~\ref{thm:coarea}}\label{subsec:coarea-proof}

In this subsection we prove Theorem~\ref{thm:coarea}. The challenge is to
change variables from the Gaussian matrix $G$ to the selected eigenpair of
$H$ while keeping two things exactly: the law of the energy $\lambda_\alpha$,
which depends on all of $H$, and the index $\alpha$ of the selected
eigenvalue. The tools are Lemma~\ref{lem:eigenpair-jacobian}, the
normalization \eqref{eq:gue-normalizer} of GUE, and Schur congruence. We
prove four claims.

\begin{figure}[H]
\centering
\begin{tikzpicture}[x=1.4cm,y=1cm]
 \begin{scope}
  \clip (-.1,-1.25) rectangle (4.2,2.1);
  \foreach \c/\s in {1.7/.45,.3/.35,-.5/.6}{
   \draw[black!55,thick,domain=0:4.1,samples=40] plot (\x,{\c-\s*\x-0.06*\x*\x});}
  \draw[very thick,domain=0:4.1,samples=40] plot (\x,{1.0-.8*\x-0.06*\x*\x});
 \end{scope}
 \draw[->] (-.1,0)--(4.35,0) node[below left] {$E$};
 \fill (1.12,0) circle (.06);
 \node[below,font=\small] at (1.12,-.05) {$\lambda_\alpha$};
 \draw[->] (2.3,1.4) node[above,font=\small] {slope $-v^*\mathcal Dv$}--(1.2,.08);
 \node[font=\small,anchor=west] at (4.75,1.5) {eigenvalue branches of $A(E)$};
 \node[font=\small,anchor=west] at (4.75,.9) {zeros: eigenvalues of $H$};
 \node[font=\small,anchor=west] at (4.75,.3) {below $0$ at $E=\lambda_\alpha$:};
 \node[font=\small,anchor=west] at (4.75,-.2) {$N_-(M)=\alpha-1-N_R(E)$};
\end{tikzpicture}
\caption{A schematic picture for Claim~3 in the proof of
Theorem~\ref{thm:coarea}; the curves are not computed. The horizontal axis
is the physical energy $E$, and each curve is an eigenvalue of the Schur
complement $A(E)$ of \eqref{eq:schur-data}. Each curve decreases, because
$\partial_EA=-\mathcal D\preceq-I$. A zero of a curve is an eigenvalue of
$H$. The thick curve vanishes at $E=\lambda_\alpha$; its slope there is
$-v^*\mathcal Dv$, and this is the crossing factor in
\eqref{eq:complete-coarea}. The remaining curves that lie below zero at
$\lambda_\alpha$ are the negative eigenvalues of $M$; together with the
$N_R(E)$ negative eigenvalues of $H_R-E$ they count the $\alpha-1$
eigenvalues of $H$ below $\lambda_\alpha$.}
\label{fig:crossing}
\end{figure}
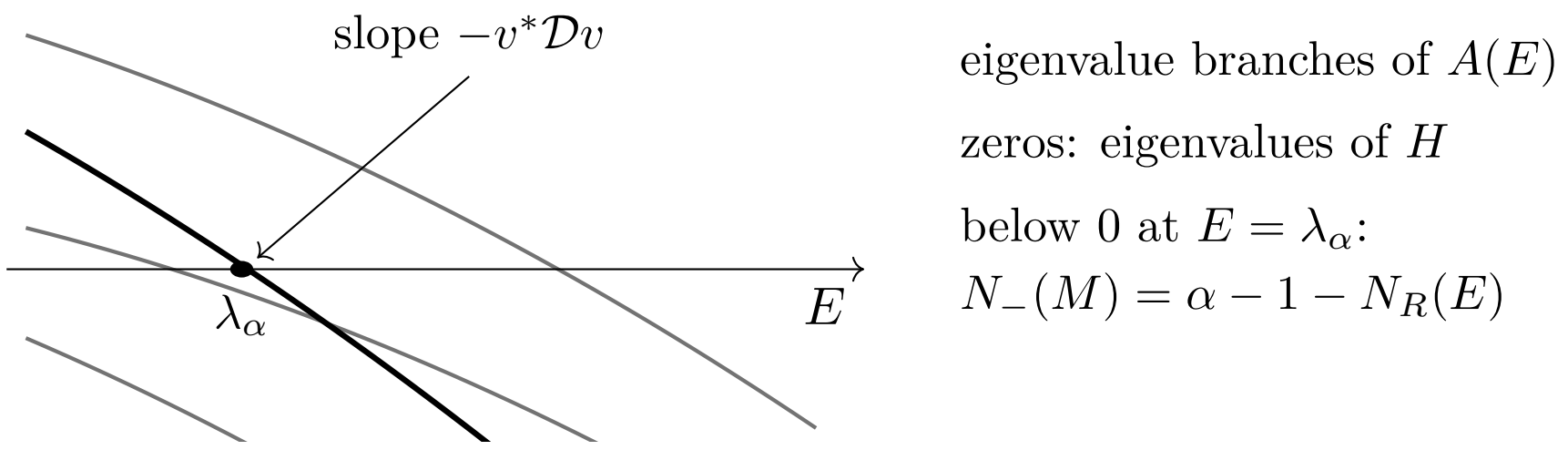

\begin{enumerate}[label=\textbf{Claim \arabic*.},leftmargin=*,align=left]
\item The identities \eqref{eq:clock-reconstruction} hold.
\item At fixed $(E,\omega)$, the GUE density of $G$ in eigenpair coordinates
gives the factor $c_{d,t}\mathcal R(E,v,F)\det(M)^2$ and the law
$\GUE_{d-1}(t)$ of $G_{d-1}$.
\item The passage from the eigenvalue of $F(E)+G$ to the physical energy $E$
produces the factor $v^*\mathcal Dv$, and the selected rank $\alpha$
corresponds to the indicator $\one_{\{N_-(M)=\alpha-N_R(E)-1\}}$.
\item The pair-mass law \eqref{eq:pair-law} has total mass one, and under it
the factor $v^*\mathcal Dv$ is replaced by $n/2$.
\end{enumerate}
Claims~2 and~3 together give \eqref{eq:complete-coarea}; taking $\Phi=1$
shows that this measure has total mass one. Claim~4 gives the statement
about \eqref{eq:pair-law}. Figure~\ref{fig:crossing} illustrates Claim~3.

\begin{proof}[Proof of Claim 1]
Write $\psi_S=P_S\psi_\alpha$ and $\psi_R=P_R\psi_\alpha$. The two block
equations of $(H-E)\psi_\alpha=0$ in the decomposition $S\oplus R$ are
\begin{equation}\label{eq:block-equations}
 C\psi_S+(H_R-E)\psi_R=0,\qquad
 (V+G-E)\psi_S+C^*\psi_R=0 .
\end{equation}
The first equation gives $\psi_R=-T\psi_S$. Substituting it into the second
gives $A\psi_S=0$, hence $Av=0$. Since $\psi_S=\sqrt m\,v$, we get
$\psi_\alpha=\sqrt m\,\mathcal L_Ev$. Taking squared norms,
\[
 1=\norm{\psi_\alpha}^2=m\bigl(\norm v^2+\norm{Tv}^2\bigr)=m\,v^*\mathcal Dv.
\]
This proves $m^{-1}=v^*\mathcal Dv$ and the reconstruction formula in
\eqref{eq:clock-reconstruction}.
\end{proof}

\begin{proof}[Proof of Claim 2]
Fix $E$ and $\omega$. By Claim~1, $Av=0$ means that $K=F(E)+G$ has the
eigenpair $(E,v)$, so we may write $K=Evv^*+U_vBU_v^*$ as in
Lemma~\ref{lem:eigenpair-jacobian}. The density of $G=K-F$ is
$Z_d(t)^{-1}e^{-\norm{K-F}_{\HS}^2/(2t)}$, and the Gaussian norm separates
exactly:
\[
 \norm{Evv^*+U_vBU_v^*-F}_{\HS}^2
 =(E-v^*Fv)^2+2\norm{(I-vv^*)Fv}^2
                         +\norm{B-F_v}_{\HS}^2.
\]
The first two terms give $\mathcal R(E,v,F)$. The last term, together with
the normalizer $Z_{d-1}(t)$, is the density of $B=F_v+G_{d-1}$ with
$G_{d-1}\sim\GUE_{d-1}(t)$. By Lemma~\ref{lem:eigenpair-jacobian} the volume
element contributes $\det(B-E)^2=\det(M)^2$ and the constant
$\pi^{d-1}/(d-1)!$. Combining that constant with $Z_{d-1}(t)/Z_d(t)$ from
\eqref{eq:gue-normalizer},
\[
 \frac{\pi^{d-1}}{(d-1)!}\cdot\frac{Z_{d-1}(t)}{Z_d(t)}
 =\frac{\pi^{d-1}}{(d-1)!}\cdot\frac{1}{\sqrt{2\pi t}\,(\pi t)^{d-1}}
 =c_{d,t}.
\]
This proves Claim~2.
\end{proof}

\begin{proof}[Proof of Claim 3]
We first compute the crossing factor and then the index.

\emph{The crossing factor.} As the physical energy varies,
\[
 F'(E)=-C^*(H_R-E)^{-2}C=-T^*T,\qquad\text{hence}\qquad
 \partial_E A=-\mathcal D.
\]
A simple zero of a local eigenvalue of $A(E)$, with eigenvector $v$,
therefore has absolute crossing derivative $|v^*\partial_EAv|=v^*\mathcal Dv$.
This is the factor in \eqref{eq:complete-coarea} that converts the
eigenvalue coordinate of Claim~2 into the physical energy; see the thick
curve in Figure~\ref{fig:crossing}.

\emph{The index.} Schur congruence gives the decomposition of inertia of
$H-E$ into that of $H_R-E$ and that of $A(E)$:
\[
 N_-(H-E)=N_-(H_R-E)+N_-(A(E)).
\]
At $E=\lambda_\alpha$ the left side is $\alpha-1$. The kernel of $A(E)$ is
spanned by $v$, and on $v^\perp$ the matrix $A(E)$ is $U_vMU_v^*$. Hence at
the selected simple zero, precisely $\alpha-1-N_R(E)$ eigenvalues in
$v^\perp$ are negative, which is the indicator in \eqref{eq:complete-coarea}.

\emph{ The global change of variables:} The following argument supplies the
global area formula across every exterior pole, rather than extending a
local identity by an unspecified limit. Fix exterior data with simple
eigenvalues   $\mu_1<\cdots<\mu_{N-d}$  and nonzero vectors
$C^*\varphi_j$, where $\varphi_j$ is an eigenvector of $H_R$.
These conditions hold almost surely: the exterior is a connected block
path with invertible Gaussian bonds; an eigenvector cannot vanish at an
endpoint, by its block recurrence. Conditional on that path, its Gaussian
bond to $S$ is invertible almost surely. Simplicity follows from the
nonzero discriminant polynomial and the continuous diagonal variables.

Decompose $G=G^0+\zeta I_d$, with $G^0$ traceless and
$\zeta\sim N(0,t/d)$ independent. On each of the   $N-d+1$  open intervals
between the $\mu_j$ (including the two unbounded intervals), let
$a_j(E)$,  $1\le j\le d$, be the increasingly ordered eigenvalues of
$F(E)+G^0-E$. They are locally Lipschitz and satisfy
\[
 a_j(E_2)\le a_j(E_1)-(E_2-E_1)\qquad(E_2>E_1)
\]
by the min--max principle and $\partial_E A\preceq-I$. Hence each map
$E\mapsto-a_j(E)$ is strictly increasing, locally absolutely continuous,
and has derivative $v_j^*\mathcal Dv_j$ at a simple branch point. Apply
the one-dimensional area formula on a compact subinterval to the equation
$\zeta=-a_j(E)$, with its Gaussian scalar density and a nonnegative
observable. Exhaust each pole interval by compact subintervals, sum over
all $j$ and all pole intervals, and apply monotone convergence. This
accounts for every crossing exactly once, with precisely the derivative
$v_j^*\mathcal Dv_j$. No integrability estimate at a pole or at infinity
is needed for this nonnegative identity.

At almost every fixed $(E,G^0)$ the branches are simple; this follows
also by the discriminant polynomial and Fubini. Equivalently the preceding
area formula reads as the surface measure at eigenvalue zero of $A(E)$,
multiplied by $v^*\mathcal Dv$. Disintegrating that surface measure in
the remaining matrix coordinates is exactly the fixed-$E$ eigenpair
Jacobian of Claim~2. Tonelli therefore gives the full integral
\eqref{eq:complete-coarea}, not an identity confined to one branch.

No physical eigenvalue is lost at a pole. For fixed $\omega$ as above,
$\det(H-\mu_j)$ is a nonzero polynomial in $G$: along $G=\zeta I_d$
its leading nonzero coefficient, of degree $d-1$, is
$-\norm{C^*\varphi_j}^2\prod_{k\ne j}(\mu_k-\mu_j)$.
Thus $\operatorname{Spec}H\cap\operatorname{Spec}H_R$ is empty almost
surely. The full discriminant is also a nonzero polynomial, so the full
spectrum is simple. A zero pair mass would give a common eigenvalue and
is consequently excluded. The inertia selection above assigns each
remaining crossing to its unique original rank. This proves Claim~3.
\end{proof}

\begin{proof}[Proof of Claim 4]
The law of $H$ is invariant under the cyclic shift of blocks
$k\mapsto k+1$, which maps the rank-$\alpha$ eigenvector to the
rank-$\alpha$ eigenvector of the shifted matrix. Hence $\E p_{\alpha k}$
does not depend on $k$, and since $\sum_kp_{\alpha k}=1$,
\[
 \E p_{\alpha k}=\frac{1}n,\qquad \E m=\frac{2}n.
\]
Thus \eqref{eq:pair-law} is normalized. By Claim~1,
$m\,(v^*\mathcal Dv)=1$, so the weight $\frac{n}{2}m$ in \eqref{eq:pair-law}
multiplies the density in \eqref{eq:complete-coarea} by $\frac{n}{2}m$ and
removes exactly the crossing derivative, leaving the constant $n/2$.
No conditional density in $E$ or in the rank sector is divided out.
This proves Claim~4.
\end{proof}

\begin{proof}[Proof of Theorem~\ref{thm:coarea}]
Claim~1 gives \eqref{eq:clock-reconstruction}. Claim~2 gives the density
of $G$ in eigenpair coordinates at fixed energy, and Claim~3 converts its
eigenvalue coordinate into the physical energy and restricts to the rank
$\alpha$. Integrating over $\omega$ proves \eqref{eq:complete-coarea}, and
$\Phi=1$ shows that the measure has total mass one. Claim~4 proves the
statement about the pair-mass law \eqref{eq:pair-law}.
\end{proof}

\section{Gaussian layers and matched boundary traces}
\label{sec:layers}
Many of our arguments cut the ring into a path and grow the path one
block at a time. This section shows that, with matched ports at its ends, the
normalized trace of the boundary resolvent stays close to the matched value
$a$ at every length, however large the old Schur complement may be. This
control is used both in Part~I and in Part~II.
Theorem~\ref{thm:matched-traces} states that, for every length up to $TW$
and every deterministic bulk energy, the normalized port traces are within
$C\sqrt p/W$ of $a$ in $L^p$. It gives the same control for the one-block
complement inside the ring.
In Subsection~\ref{subsec:s5-uniform-inverse-moments} we prove inverse moments from the fresh Gaussian
block and its coupling to an absorbing boundary, uniformly in the old Schur
matrix (Lemma~\ref{lem:gaussian-absorption}). In Subsection~\ref{subsec:s5-gaussian-contractions-with-matrix-weight} we show
that covariances between a normalized trace and a matrix entry are of order
$W^{-2}$, and that trace fluctuations are of order $W^{-1}$
(Lemma~\ref{lem:trace-matrix}). In Subsection~\ref{subsec:s5-the-mean-equations-and-their-stable-bran} we derive the mean
equations of one layer, show that one layer contracts a small error
(Lemma~\ref{lem:layer-contraction}), and iterate along the path to prove the
theorem.

Fix a compact interval $I\Subset(-2,2)$ and put
\begin{equation}\label{eq:matched-parameters}
\begin{aligned}
 a&=\frac{-E-i\sqrt{4-E^2}}{2},& \eta&=\Im z,\\
 z=z_\partial&=a^{-1}+ba=-E-(1-b)a,&
 z_0&=a^{-1}+2ba=-E-ca.
\end{aligned}
\end{equation}
Thus $|a|=1$, $\Im a<0$, and both $z,z_0$ have positive imaginary
parts bounded below on $I$.
Let $X_r$ be the Schur self-energy of a path of $r$ physical blocks,
coupled by independent Ginibre bonds to two distinct ports.
The two diagonal traces of $R_r=(z-X_r)^{-1}$ are divided by $W$.
Growing the path from its center gives the exact recursion
\begin{equation}\label{eq:two-ended-recursion}
 X_{r+2}=B^*(A-E-X_r)^{-1}B,\qquad
 A=\diag(A_1,A_2),\quad B=\diag(B_1,B_2).
\end{equation}
Here $A_i$ are fresh GUE blocks of entry variance $c/W$ and $B_i$
are fresh circular Gaussian blocks of entry variance $b/W$.

\begin{theorem}\label{thm:matched-traces}
For fixed $b,I,T$ there are constants $\epsilon_0,C>0$ such that
for $2\le p\le\epsilon_0W$, $1\le r\le TW$ and deterministic $E\in I$,
\begin{equation}\label{eq:matched-trace-moments}
 \left\|\max_{i=1,2}
 \left|W^{-1}\Tr(R_r)_{ii}-a\right|\right\|_{L^p}
 \le\frac{C\sqrt p}{W}.
\end{equation}
For the one-block complement $\Sigma=C^*(H_R-E)^{-1}C$ in the ring,
\begin{equation}\label{eq:common-port-trace}
 \left\|W^{-1}\Tr(z_0-\Sigma)^{-1}-a\right\|_{L^p}
 \le\frac{C\sqrt p}{W}.
\end{equation}
The analogous residual-pair traces (which are equal to $W^{-1}\Tr P_a(K(E)-sa)^{-1}$ on the two exposed blocks, with $K(E)=F(E)-E$) after subtracting an independent
$\GUE_{2W}(s/(2W))$ obey the same bound, with constants depending
also on fixed $0 < s < 2 \min(1 - 2b, b)$. The upper-half-plane version of these traces is their complex conjugate.
\end{theorem}

The bounds include the one- and two-block cores and therefore both
parities and bounded ring lengths. The energy is deterministic and its
supremum stays outside the moment; selected energies will be treated
by a separate energy-transport argument.

\subsection{Uniform inverse moments}\label{subsec:s5-uniform-inverse-moments}

We   need an inverse estimate that remains valid after any old
Hermitian Schur matrix has been fixed. Lemma~\ref{lem:gaussian-absorption}
gives bounded inverse moments from the fresh Gaussian block and its
coupling to an absorbing boundary. We separate the nearly singular
subspace, control its size by eigenvalue counting, and use the Gaussian
coupling to control the inverse there.

\begin{lemma}\label{lem:gaussian-absorption}
Let $X=X^*\in\C^{2W\times2W}$ be independent of $A=\diag(A_1,A_2)$.
Let either $B=\diag(B_1,B_2)$ or $B=(B_1^T,B_2^T)^T$.
For $\Im z\ge\eta_0>0$ and $|z|\le M$, define
\[
 \mathcal G=\begin{pmatrix}A-E-X&B\\B^*&zI\end{pmatrix}^{-1}.
\]
There are $\epsilon_1,C>0$, depending only on $b,\eta_0,M$, such that
\begin{equation}\label{eq:absorber-moments}
 \left(\E\norm{\mathcal G}^{\,p}\right)^{1/p}\le C,
 \qquad 1\le p\le\epsilon_1W.
\end{equation}
The assertion is uniform in Hermitian $X$ and real $E$. It also holds
for one physical block and one or two independent ports.
\end{lemma}

Independent Gaussian coupling regularizes the required inverse even at a real physical energy. Uniformity in the old Hermitian matrix is essential for repeated Schur elimination.

\begin{proof}
Put $D=A-E-X$. Let $P$ project onto the $k=\lfloor W/8\rfloor$
eigenvectors of $D$ nearest zero, and let $\delta\le1$ be the
truncated absolute value of the next eigenvalue.
Theorem~\ref{thm:external-count}, after scaling by $\sqrt c$, gives
\[
 \Prob(\delta<t)\le
 \frac{(C_bWt)^{k+1}}{(k+1)!}\le(C_b't)^{k+1}
 \quad(0<t<t_0).
\]
Conditional on $D$, let $d=s_{\min}(B^*|_{\operatorname{Ran}P})\wedge1$
and $K=\norm B$. For a fixed unit vector $v=(v_1,v_2)$, in the
separate-port case
\[
 \norm{B^*v}^2\ \stackrel{\rm law}=
 \frac{b}W\bigl(\norm{v_1}^2G_1+\norm{v_2}^2G_2\bigr),
\]
where $G_1,G_2$ are independent gamma variables of shape $W$.
At least one coefficient is at least $1/2$. In the common-port case
$B^*v$ is directly a circular Gaussian vector with covariance $bI/W$.
The gamma density therefore gives in either case
$\Prob(\norm{B^*v}<t)\le(C_bt^2)^W$.
On $K\le t^{-1/2}$, a net of mesh $t^{3/2}/2$ on the unit sphere
of $\operatorname{Ran}P$ gives
\[
 \Prob(d<t\mid D)
 \le Ce^{-c_bW/t}+(Ct^{-3/2})^{2k}(C_bt^2)^W
 \le2(C_b't)^W
\]
for small fixed $t$. This bound is uniform in the subspace $P$.
Integrating the two lower tails, with a cutoff $t_1<(2C_b')^{-1}$,
shows
\[
 \|\delta^{-1}\|_{L^q}+
 \|d^{-1}\|_{L^q(\,\cdot\,\mid D)}\le C_b,
 \qquad q\le\epsilon_2W.
\]
For example the lower-tail integral is bounded by
$q(C_b')^{k+1}t_1^{k+1-q}/(k+1-q)$, and the denominator remains
a fixed fraction of $k+1$. Gaussian norm tails give
$\|K\|_{L^q}\le C_b(1+\sqrt{q/W})$.

Write $-z^{-1}=\alpha+i\gamma$, where
$\gamma\ge\eta_0/M^2$. If $\norm v=1$ and
$\varepsilon=\norm{(D+(\alpha+i\gamma)BB^*)v}\le1$, the imaginary
quadratic form gives $\norm{B^*v}\le\sqrt{\varepsilon/\gamma}$.
Consequently
\[
 \norm{Dv}\le C(1+K)\sqrt\varepsilon,\qquad
 \norm{(I-P)v}\le C(1+K)\sqrt\varepsilon/\delta.
\]
Using $d\norm{Pv}\le\norm{B^*Pv}$ and
$1\le\norm{Pv}+\norm{(I-P)v}$ yields
$1\le C(1+K)^2\sqrt\varepsilon/(d\delta)$.
The resulting inverse bound also holds when $\varepsilon>1$:
\begin{equation}\label{eq:absorber-deterministic}
 \norm{(D-z^{-1}BB^*)^{-1}}
       \le C(1+K)^4d^{-2}\delta^{-2}.
\end{equation}
H\"older's inequality and the conditional negative moments prove
the required bound for the top Schur inverse. The block inverse
formula proves it for $\mathcal G$, after reducing $\epsilon_1$
by a fixed factor. For either finite core the same argument conditions
on its internal bond and uses its independent exterior ports.
\end{proof}

\subsection{Gaussian contractions with matrix weights}\label{subsec:s5-gaussian-contractions-with-matrix-weight}

  We need covariance estimates that still retain a matrix weight outside
 a normalized trace  . Lemma~\ref{lem:trace-matrix} gives a covariance of
order $W^{-2}$ and trace fluctuations of order $W^{-1}$. We write the
covariance by Gaussian interpolation, differentiate the inverse
matrices, and bound the resulting products using the inverse moments.

\begin{lemma}\label{lem:trace-matrix}
Let $f$ be a trace divided by $W$ of a fixed word in
$\mathcal G,\mathcal G^*$ and bounded deterministic matrices, and let
$U$ be a fixed matrix word of the same type. For each fixed moment
order and word length, 
for $W\ge W_0(p,\ell)$, \begin{equation}\label{eq:trace-matrix-cov}
 \norm{\E(fU)-(\E f)(\E U)}\le CW^{-2}.
\end{equation}
A normalized block trace satisfies
\begin{equation}\label{eq:fresh-trace-fluctuation}
 \|f-\E f\|_{L^p}\le C\sqrt p/W,\qquad 2\le p\le\epsilon_3W,
\end{equation}
when $f$ is a single resolvent trace. All these expectations are
conditional on the independent Hermitian deformation.
\end{lemma}

The contraction retains a matrix test outside the normalized trace. It gives the extra inverse-width factor needed for the conditional mean equation.

\begin{proof}
In real Gaussian coordinates $x$, let
$x_t=e^{-t}x+\sqrt{1-e^{-2t}}\,x'$ with $x'$ independent.
Differentiation of this interpolation gives the bilinear covariance
identity
\[
 \operatorname{Cov}(f(x),U(x))
 =\int_0^\infty e^{-t}
 \E\bigl[\nabla f(x)\cdot\Sigma\nabla U(x_t)\bigr]\,dt.
\]
Here $\Sigma$ is the real-coordinate covariance, and no conjugation
is inserted in this bilinear identity. A term in $df$ is
$W^{-1}\Tr K\,d\mathcal H$; a term in $dU$ is
$L\,d\mathcal H\,N$, with signs included in the factors.
The entry covariance
\[
 \E(d\mathcal H_{ij}d\mathcal H_{k\ell})
   =W^{-1}\sigma_{ij}\delta_{i\ell}\delta_{jk}
\]
contracts these terms to
$W^{-2}L\mathcal V(K)N$.
For one physical block and one port,
$\mathcal V(K)=cP_1KP_1+bP_1KP_2+bP_2KP_1$.
For separate or common ports it is the finite sum over the actual
diagonal blocks and actual bonds, with coefficients $c,b$.
Its operator norm is bounded independently of $W$.
H\"older and Lemma~\ref{lem:gaussian-absorption} bound the product
of the remaining factors, at both $x$ and $x_t$. Integrating $e^{-t}$
proves \eqref{eq:trace-matrix-cov}.

For a single trace, $d\mathcal G=-\mathcal G(d\mathcal H)\mathcal G$
gives a Hilbert--Schmidt gradient bounded by
$CW^{-1/2}\norm{\mathcal G}^2$. The dimension-free Gaussian
$L^p$ gradient inequality, with covariance $O(W^{-1})$, gives
\eqref{eq:fresh-trace-fluctuation}. The growing inverse moments
allow $p$ up to a sufficiently small multiple of $W$.
Smooth cutoff approximation is justified by the same inverse moments.
\end{proof}

\subsection{The mean equations and their stable branch}\label{subsec:s5-the-mean-equations-and-their-stable-bran}

 We identify the scalar value toward which the boundary traces move.
Lemma~\ref{lem:layer-contraction} shows that a fresh layer reduces a
small incoming trace error, up to $O(W^{-2})$. We derive the conditional
mean equations, check stability near the matched solution, and combine
entrance, fluctuation, and exit estimates to prove
Theorem~\ref{thm:matched-traces}.

 Write $m_i=\E W^{-1}\Tr\mathcal G_{ii}$ for the physical blocks
and $u_i$ for the port means. The deterministic comparison equations
for separate ports are
\begin{equation}\label{eq:two-port-dyson}
 m_i=W^{-1}\Tr
 [-E-X-c\diag(m_1,m_2)-b\diag(u_1,u_2)]^{-1}_{ii},
 \qquad u_i=(z-bm_i)^{-1}.
\end{equation}
In this display the same letters on the right denote the unknown
deterministic solution; the actual means satisfy these equations with
an $O(W^{-2})$ residual. For a common port replace the last diagonal
term by $bu_0I$ and set $u_0=(z_0-b(m_1+m_2))^{-1}$.

\begin{lemma}\label{lem:layer-contraction}
Put $d(X)=\max_i|W^{-1}\Tr[(z-X)^{-1}]_{ii}-a|$.
There are $d_*,C>0$ and $q_*<1$, depending only on $b,I$, such that
on $d(X)\le d_*$ the actual outgoing means satisfy
\[
 \max_i|\E W^{-1}\Tr(\mathcal G_{\rm ports})_{ii}-a|
 \le q_*d(X)+CW^{-2}
\]
for separate ports, and the common-port error is at most
$Cd(X)+CW^{-2}$.
\end{lemma}

The matched scalar branch is stable on a fixed neighborhood. The entrance and moment estimates allow this local stability to control the full allowed path length.

\begin{proof}
Integration by parts in an $A$ entry gives
$\E(A_i\mathcal G_{ij})=-c\E(g_i\mathcal G_{ij})$.
For a bond it gives the corresponding term
$-b\E(g_{\rm port}\mathcal G_{ij})$.
For example the top and bottom equations for one block are
\[
 (D-c\bar g_1-b\bar g_2)\E\mathcal G_{11}
 =I+c\operatorname{Cov}(g_1,\mathcal G_{11})
       +b\operatorname{Cov}(g_2,\mathcal G_{11}),
\quad
 (z-b\bar g_1)\E\mathcal G_{22}
 =I+b\operatorname{Cov}(g_1,\mathcal G_{22}).
\]
These signs follow directly from the minus sign in the inverse
derivative. The supported multi-block equations are the same entry
identities with a sum over the incident bonds.
Lemma~\ref{lem:trace-matrix} bounds their covariance terms.
The port equation gives $(z-b\bar m_i)\bar u_i=1+O(W^{-2})$,
so $|\bar u_i|$ is bounded below. The exact identity
\[
 -\Im\mathcal G=(\Im z)\mathcal G^*P_{\rm ports}\mathcal G
\]
then bounds $-\Im\bar u_i$ below. Thus the inverse with the top
mean shifts is uniformly bounded and the normalized equations have
residual $O(W^{-2})$.

At $m_i=u_i=a$, the top inverse equals $R=(z-X)^{-1}$.
Set $Q_{ij}=W^{-1}\Tr R_{ij}R_{ji}$. Normality of $R$ and the
imaginary resolvent identity give
\[
 \max_i\sum_j|Q_{ij}|
 \le(1-b)^{-1}+O(d(X)).
\]
The linearized separate-port equations are
$v=d+(c+b^2a^2)Qv$, $w=ba^2v$. Since
$|c+b^2a^2|<(1-b)^2$ in the bulk, their outgoing coefficient obeys
\begin{equation}\label{eq:strict-contraction}
 q_0(E)=\frac{b(1-b)}{1-b-|c+b^2a^2|}<1.
\end{equation}
The gap is uniform on $I$. The implicit-function theorem gives a
negative-imaginary solution near $a$ and outgoing error
$q_0d+Cd^2$. For the common port the linearized top coefficient is
$Q(cI+b^2a^2\one\one^T)$, of infinity norm at most
$(c+2b^2)/(1-b)+O(d)<1$. This gives its bounded coefficient.

We must also identify the branch selected by the actual means.
For a top inverse $M$ put
$A_{ij}=W^{-1}\Tr M_{ij}M_{ij}^*$. With $v_i=-\Im m_i$ and
$w_i=-\Im u_i$, its two imaginary identities give
\[
 v_i=\sum_j A_{ij}y_j=\sum_j A_{ji}y_j,\qquad
 y_j=cv_j+bw_j
\]
(replace $w_j$ by $w_0$ in the common-port case).
For the approximate means the identities have $O(W^{-2})$ errors.
The port $w$ is bounded below as above. On $d(X)\le d_*$ the top
$v_i$ is also bounded below: the inverse identity comparing $M$
with $R$ gives a fixed lower bound on each normalized squared block
column of $M$ from that of $R$, which equals
$-\Im(W^{-1}\Tr R_{ii})/\eta$. All these quantities are bounded above.

Subtract the approximate and local equations and weight each component
by the geometric mean of the two negative imaginary parts.
Cauchy--Schwarz bounds a top coefficient by
$\sqrt{A^{(1)}_{ij}A^{(2)}_{ji}}$. The displayed identities make
the weighted row sums at most one, up to $O(W^{-2})$.
Each port row loses a fixed fraction to $\Im z>0$; every top row
either loses a fixed fraction already or sends a fixed fraction to a
port, since $bw/(cv+bw)$ is bounded below.
After two iterations all row sums are at most $1-\epsilon$.
The weights have fixed upper and lower bounds. Thus the difference
equations have a uniformly bounded inverse and their solution
difference is $O(W^{-2})$. This proves that the local branch is the
one approximating the actual means, and proves the lemma.
\end{proof}

\begin{proof}[Proof of Theorem~\ref{thm:matched-traces}]
 A path of $r$ physical blocks is attached
at its two ends to two ports. Recall that a port is an extra $W$-dimensional
coordinate block with the complex number $z$ on its diagonal. It
replaces the rest of the system, which we have cut away, by a
semi-infinite lead. After eliminating the path by a Schur complement, we
are left with the $2W\times2W$ matrix $z-X_r$ on the two ports, and
$R_r=(z-X_r)^{-1}$ is the resolvent seen from the ports. 

\emph{Plan of the proof.} The recursion \eqref{eq:two-ended-recursion}
builds a path of $r+2$ blocks from a path of $r$ blocks by adding one
fresh block at each end, as in Figure~\ref{fig:matched-growth}. We
therefore argue by induction on $r$ in steps of two. First we prove the
claim for the two shortest paths, $r=1$ and $r=2$, which start the odd
and the even lengths. Then we show that one step of the recursion
contracts the error. Finally we sum the errors over at most $TW$ steps,
and we control the rare event that the error ever becomes large.

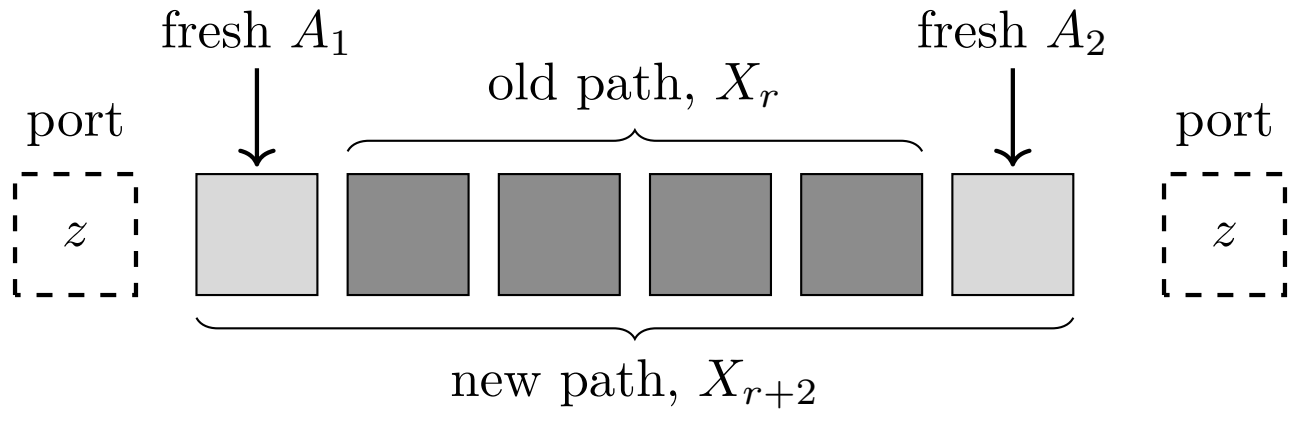
\begin{figure}[ht]
\centering
\begin{tikzpicture}[x=1cm,y=1cm]
  \foreach \x/\f in {0/black!15,1/black!45,2/black!45,3/black!45,4/black!45,5/black!15}
   \draw[fill=\f] (\x,0) rectangle ++(.8,.8);
 \draw[thick,dashed] (-1.2,0) rectangle ++(.8,.8);
 \draw[thick,dashed] (6.4,0) rectangle ++(.8,.8);
 \node[font=\small] at (-.8,.4) {$z$};
 \node[font=\small] at (6.8,.4) {$z$};
 \node[font=\small,above] at (-.8,.85) {port};
 \node[font=\small,above] at (6.8,.85) {port};
 \draw[decorate,decoration={brace,amplitude=4pt}] (1,.95)--(4.8,.95)
   node[midway,above=4pt,font=\small] {old path, $X_r$};
 \draw[decorate,decoration={brace,amplitude=4pt,mirror}] (0,-.15)--(5.8,-.15)
   node[midway,below=4pt,font=\small] {new path, $X_{r+2}$};
 \draw[->,thick] (.4,1.5)--(.4,.85);
 \draw[->,thick] (5.4,1.5)--(5.4,.85);
 \node[font=\small,above] at (.4,1.45) {fresh $A_1$};
 \node[font=\small,above] at (5.4,1.45) {fresh $A_2$};
\end{tikzpicture}
\caption{One step of the recursion \eqref{eq:two-ended-recursion}, drawn
for illustration. The old path (dark blocks) is seen from its ends only
through its Schur complement $X_r$. One fresh GUE block (light) is added at
each end, coupled by fresh Gaussian bonds, and the ports (dashed) are
moved outward. The new Schur complement $X_{r+2}$ depends on the old path only
through $X_r$. Since each step adds two blocks, odd and even lengths are
started separately, from the one-block and the two-block cores.}
\label{fig:matched-growth}
\end{figure}

\emph{Step 1:} Here we seek to derive equation of mean for the shortest path. We start with
a single physical block between the two ports. We order its coordinates
as physical block, left port, right port. For two physical blocks we
order them as the two blocks followed by their ports. The entry
variances, divided by $W^{-1}$, are recorded in the matrices
\[
 S_3=\begin{pmatrix}c&b&b\\b&0&0\\b&0&0\end{pmatrix},\qquad
 S_4=\begin{pmatrix}c&b&b&0\\b&c&0&b\\b&0&0&0\\0&b&0&0\end{pmatrix}.
\]
For example, the first row of $S_3$ says that the physical block has
internal variance $c$ and is coupled with variance $b$ to each port,
while the port rows have no internal randomness. The diagonal entries of
the matrix whose inverse we take are
$\mathbf h=(-E,z,z)$ or $\mathbf h=(-E,-E,z,z)$, respectively.

Let $\bar m_i$ denote the expected normalized trace of the $i$-th
diagonal block of the inverse. Gaussian integration by parts in every
Gaussian entry gives
\begin{equation}\label{eq:mt-vector-dyson}
 \Bigl(\mathbf h_i-\sum_j(S_\ell)_{ij}\bar m_j\Bigr)\bar m_i=1+e_i,
 \qquad |e_i|\le CW^{-2},\quad \ell=3,4 .
\end{equation}
This is the vector version of the semicircle equation. Each block feels
the others through their mean traces, weighted by the variances. Two
facts make it a scalar equation for each $i$. First, the law is
invariant under unitary rotations inside each block, so each mean
diagonal block is a multiple of the identity. Second, the error $e_i$
collects covariances between a normalized trace and a matrix entry.
Lemma~\ref{lem:trace-matrix} shows that each such covariance is
$O(W^{-2})$. That lemma needs inverse moments, and for the two-block core
we obtain them by conditioning on the bond between the two physical
blocks.

\emph{Step 2:} Here we show that the exact equations are solved by $a$, and stably. Without
the errors, \eqref{eq:mt-vector-dyson} is solved by the constant vector
$\bar m_i=a$. On a physical row, the variances sum to $c+2b=1$, and the
equation becomes $(-E-a)a=1$. On a
port row, the only variance is $b$, and the equation becomes
$(z-ba)a=1$. We now ask whether a small error $e$ moves the solution only a little.
Linearizing \eqref{eq:mt-vector-dyson} at $\bar m=a$ produces the
operator $I-a^2S_\ell$, so we must show that it has a bounded inverse.
The rows of $S_\ell$ sum to at most one, so $S_\ell$ itself has norm one
and is not a contraction. Its square, however, is:
\begin{equation}\label{eq:mt-square-norm}
 \|S_3^2\|_\infty=\max\{1-2b(1-b),b\},\qquad
 \|S_4^2\|_\infty=\max\{1-b(1-b),b\}<1 .
\end{equation}
The reason is that mass leaks out at the ports. A port row sums to
$b<1$, and every physical row reaches a port in at most two steps, so
after two steps each row has lost a fixed fraction. For example, the
first row of $S_3^2$ sums to $c+2b^2=1-2b(1-b)$. Since $|a|=1$, we can
expand in powers of $a^4S_\ell^2$:
\[
 (I-a^2S_\ell)^{-1}=(I+a^2S_\ell)\sum_{k\ge0}(a^4S_\ell^2)^k,
 \qquad
 \|(I-a^2S_\ell)^{-1}\|_\infty\le\frac{2}{b(1-b)}.
\]

There is one more issue. The equation \eqref{eq:mt-vector-dyson} is
quadratic, so it has a second solution, and the expansion above is only
useful if the actual means are near $a$ rather than near that other
solution. Lemma~\ref{lem:layer-contraction} settles this with an
imaginary-part identity: the port rows lose a fixed fraction of
$\Im$ to the boundary parameter $\Im z>0$, and this selects $a$. With the
bounded inverse, the residual $e=O(W^{-2})$ becomes an error $O(W^{-2})$
in the means. Finally, a normalized trace depends on the $W^2$ Gaussian
entries in a smooth way, with gradient of order $W^{-1}$. The Gaussian
Poincar\'e inequality therefore bounds its fluctuation around the mean
by $C\sqrt p/W$ in $L^p$. This proves the theorem for $r=1$ and $r=2$.

\emph{Step 3:} Now we show that one step of the recursion contracts the error and those could be summed over many steps. We measure
the distance of a path from the matched value by
\[
 d(X)=\max_{i=1,2}\bigl|W^{-1}\Tr[(z-X)^{-1}]_{ii}-a\bigr| ,
\]
as in Lemma~\ref{lem:layer-contraction}. Suppose that $d(X_r)\le d_*$
for the small constant $d_*$ of that lemma. The lemma repeats Steps~1
and~2 for one new layer, with the old path entering only through $X_r$.
It shows that the mean of the new distance is at most
$q_*d(X_r)+CW^{-2}$, for a fixed $q_*<1$. Adding the fluctuation from the
fresh Gaussian layer, as in Step~2, we get
\begin{equation}\label{eq:mt-one-step}
 d(X_{r+2})\le q_*d(X_r)+\frac C{W^2}+|\xi_{r+2}|,\qquad
 \|\xi_{r+2}\|_{p\mid X_r}\le\frac{C\sqrt p}W,
\end{equation}
where the norm is conditional on the old path. So each step shrinks the
old error by the factor $q_*$ and adds a new error of size $\sqrt p/W$. Iterating \eqref{eq:mt-one-step}
gives a sum in which the error added $k$ steps ago carries the factor
$q_*^k$. By Minkowski's inequality and the geometric series,
\[
 \|d(X_r)\|_{L^p}
 \le\sum_{k\ge0}q_*^k\Bigl(\frac C{W^2}+\frac{C\sqrt p}W\Bigr)
 \le\frac{C\sqrt p}W,
\]
with a constant that does not grow with $r$. This is the uniformity in
the length.

\emph{Step 4:} The goal of this step is to control the rare event that the error becomes large. The argument
of Step~4 used \eqref{eq:mt-one-step}, which is valid only while
$d\le d_*$. We therefore stop the recursion at the first step at which
the distance exceeds $d_*$, and we show that this almost never happens.
Since the old error was small, a first exit requires a fresh error
$|\xi|$ of fixed size. Take $q=\lfloor\epsilon_3W\rfloor$, the largest
moment allowed in Lemma~\ref{lem:trace-matrix}. Markov's inequality with
the $q$-th moment gives $\Prob(|\xi|\ge c)\le(C\sqrt q/W)^q\le(C/\sqrt W)^q$.
A union bound over the at most $TW$ steps gives
\begin{equation}\label{eq:scalar-exit}
 \Prob(\text{exit by }TW)
 \le C_TW(C/\sqrt W)^q\le e^{-cW\log W}.
\end{equation}
On the exit event we only need a crude bound. The normalized traces are
always at most $1/\Im z$, a fixed constant, because $\Im(z-X)\succeq\Im z$
for Hermitian $X$. So for $p\le q/8$ the exit event contributes at most
$O(W^{-1})$ to the $L^p$ norm. Together with Step~3, this proves
\eqref{eq:matched-trace-moments}.

\emph{Step 5:} Now we extend the bound to the full ring. The same scheme proves
the other two bounds. For the ring with one block removed, we grow the
middle $n-3$ physical blocks with separate ports, as above. The two
remaining end blocks, together with their bonds to the removed block,
form one last layer in which both ends meet at a common port. The common
port has the parameter $z_0=a^{-1}+2ba$, since it is coupled to two
neighbours. Lemma~\ref{lem:layer-contraction} gives the mean bound for
this layer, and the Poincar\'e inequality the fluctuation bound. This
proves \eqref{eq:common-port-trace}.

\emph{Step 6:} It remains to prove the last assertion
of the theorem, which concerns the matrix obtained after a fresh Gaussian
part has been removed from an adjacent pair. We subtract an independent
$\GUE_{2W}(s/(2W))$ from the pair. What remains of the pair is a residual
Gaussian layer, with diagonal variance $(c-s/2)/W$ and internal-bond
variance $(b-s/2)/W$. We write
\[
 \alpha=c-\frac s2,\qquad \beta=b-\frac s2,\qquad y=-\Im a>0,
\]
where $\alpha,\beta>0$ because $0<s<2\min(c,b)$. The test matrix of this
layer is
\[
 \mathcal M=(V-E-X-sa)^{-1},
\]
where $V$ is the residual Gaussian layer and $X$ is the Hermitian Schur
complement of the rest of the ring, seen from the pair. The row sum of
its mean variance is $\alpha+\beta=1-b-s$. At the matched value $a$, its
deterministic inverse is again $(z-X)^{-1}$, and its stability gap, which
plays the role of \eqref{eq:mt-square-norm}, is at least $s/(1-b)$. We
make this precise in five parts.

\emph{(a) A deterministic bound on the test matrix.} This layer is used
only once, at the very end, so it does not need the contraction factor
$q_*<1$ of Step~3. Instead, its test matrix is bounded deterministically.
Since $V$ and $X$ are Hermitian,
\[
 \Im(V-E-X-sa)=sy\,I,
 \qquad\text{hence}\qquad
 \norm{\mathcal M}\le(sy)^{-1}
\]
for every Hermitian $X$ and every realization of the Gaussian layer.

\emph{(b) Fluctuations of the block traces.} For the two blocks
$k=1,2$ of the pair, let
\[
 g_k=W^{-1}\Tr P_k\mathcal M .
\]
Its derivative in a direction $\Delta$ of the Gaussian entries is
$-W^{-1}\Tr P_k\mathcal M\Delta\mathcal M$. By the Cauchy--Schwarz
inequality for the Hilbert--Schmidt pairing and part~(a),
\[
 \bigl|W^{-1}\Tr P_k\mathcal M\Delta\mathcal M\bigr|
 \le W^{-1}\norm{\mathcal MP_k\mathcal M}_{\HS}\norm\Delta_{\HS}
 \le W^{-1/2}(sy)^{-2}\norm\Delta_{\HS}.
\]
The entry variances are at most $1/W$, so the Gaussian $L^p$ Poincar\'e
inequality (Remark~\ref{rem:poincare}) gives, conditionally on $X$ and
for every $X$,
\begin{equation}\label{eq:mt-residual-fluct}
 \|g_k-\E g_k\|_{L^p}\le C\frac{\sqrt p}{W}.
\end{equation}

\emph{(c) The mean equation and its matched solution.} As in the proof of
Lemma~\ref{lem:layer-contraction}, Gaussian integration by parts in the
entries of the layer, together with the fluctuation bound
\eqref{eq:mt-residual-fluct}, gives for the means $\bar m_k=\E g_k$
\begin{equation}\label{eq:mt-residual-mean}
 \bar m_k=W^{-1}\Tr P_kM+e_k,\qquad |e_k|\le CW^{-2},
\end{equation}
where
\[
 M=\Bigl(-E-X-\sum_{l=1}^2\sigma_lP_l\Bigr)^{-1},\qquad
 \sigma_l=sa+\alpha\bar m_l+\beta\bar m_{3-l}.
\]
Here $\sigma_l$ is the mean self-energy on block $l$: the fresh part
contributes $sa$, the diagonal block contributes $\alpha\bar m_l$, and the
bond to the other block contributes $\beta\bar m_{3-l}$. Since
$\Im\bar m_l\le0$, we have $\Im(-\sigma_l)\ge sy$, and therefore again
$\norm M\le(sy)^{-1}$.

The matched value $a$ solves the limiting equation. Indeed,
\[
 \alpha+\beta+s=1-b,\qquad -E-(1-b)a=a^{-1}+ba=z,
\]
so if $\bar m_1=\bar m_2=a$, then $\sigma_l=(1-b)a$ and $M$ becomes
\[
 R=(z-X)^{-1}.
\]
The distance of $R$ from the matched value is measured by
\[
 |W^{-1}\Tr P_kR-a|\le d=d(X).
\]

\emph{(d) Stability: the means are close to the matched value.} We compare
the actual means with $a$. Put
\[
 \delta_k=\bar m_k-a,\qquad v_k=-\Im\bar m_k,
\]
\[
 A_{kl}=W^{-1}\norm{P_kMP_l}_{\HS}^2,\qquad
 B_{kl}=W^{-1}\norm{P_kRP_l}_{\HS}^2 .
\]
The difference of the self-energies is
$\sigma_l-(1-b)a=\alpha\delta_l+\beta\delta_{3-l}$, so the resolvent
identity reads
\[
 M-R=M\Bigl(\sum_l\sigma_lP_l-(1-b)a\Bigr)R
    =M\Bigl(\sum_l(\alpha\delta_l+\beta\delta_{3-l})P_l\Bigr)R .
\]
Taking the normalized trace on block $k$, using \eqref{eq:mt-residual-mean}
and the Cauchy--Schwarz inequality gives
\begin{equation}\label{eq:mt-residual-diff}
 |\delta_k|\le d+CW^{-2}
 +\sum_lA_{kl}^{1/2}B_{lk}^{1/2}\bigl(\alpha|\delta_l|+\beta|\delta_{3-l}|\bigr).
\end{equation}
The imaginary parts of the two resolvents give the identities
\begin{align}
 -\Im W^{-1}\Tr P_kM&=\sum_lA_{kl}\,(sy+\alpha v_l+\beta v_{3-l}),
                                     \label{eq:mt-residual-imM}\\
 \sum_lB_{lk}&=\frac{-\Im W^{-1}\Tr P_kR}{(1-b)y}\le\frac{y+d}{(1-b)y},
                                     \label{eq:mt-residual-imR}
\end{align}
where in the second we used that $R$ is normal and that $\Im z=(1-b)y$.

The imaginary parts $v_k$ are bounded below. Indeed, the comparison
\[
 P_kR=P_kM\Bigl(I-\Bigl(\sum_l\sigma_lP_l-(1-b)a\Bigr)R\Bigr)
\]
together with \eqref{eq:mt-residual-imM} gives
\[
 v_k\ge y_0>0\qquad\text{whenever } d\le y/2 \text{ and } W \text{ is large}.
\]

We now close the estimate in a weighted norm. Put
\[
 \mathcal T=\max_k\frac{|\delta_k|}{(v_k+W^{-2})^{1/2}} .
\]
Apply the Cauchy--Schwarz inequality twice to the sum in
\eqref{eq:mt-residual-diff}, and use \eqref{eq:mt-residual-imM},
\eqref{eq:mt-residual-imR} and $A_{kl}\le(sy)^{-2}$. This bounds the sum by
\[
 \mathcal T\,(v_k+CW^{-2})^{1/2}\,\theta,
\]
where
\[
 \theta^2=\Bigl(1-\frac{s}{1-b}\Bigr)\Bigl(1+\frac dy\Bigr)
 \le1-\frac{s}{2(1-b)}
 \qquad\text{for } d\le d_{\rm res}=\frac{sy}{2(1-b)} .
\]
The factor $1-s/(1-b)$ is exactly the stability gap mentioned at the
beginning of this step. Dividing \eqref{eq:mt-residual-diff} by
$(v_k+W^{-2})^{1/2}\ge y_0^{1/2}$ gives
\[
 \mathcal T\le\theta_0\mathcal T+C(d+W^{-2})
\]
with a fixed $\theta_0<1$. Hence
\begin{equation}\label{eq:mt-residual-mean-bound}
 |\E g_k-a|\le C\bigl(d(X)+W^{-2}\bigr)
 \qquad\text{on the event } d(X)\le d_{\rm res}.
\end{equation}

\emph{(e) Conclusion.} The matrix $X$ is the two-port Schur complement of
the exterior path of $n-2$ blocks, attached to the pair as two separate
ports. So \eqref{eq:matched-trace-moments} with $r=n-2$ gives
\[
 \|d(X)\|_{L^p}\le C\frac{\sqrt p}{W}.
\]
On the event $d(X)\le d_{\rm res}$, we combine
\eqref{eq:mt-residual-fluct} and \eqref{eq:mt-residual-mean-bound}. On the
complementary event $d(X)>d_{\rm res}$, we use the deterministic bound
$|g_k|\le(sy)^{-1}$ from part~(a). By Markov's inequality and
\eqref{eq:matched-trace-moments} at order $2p$ (after halving
$\epsilon_0$),
\[
 \Prob\bigl(d(X)>d_{\rm res}\bigr)^{1/p}
 \le\frac{\|d(X)\|_{L^{2p}}^2}{d_{\rm res}^2}
 \le C\frac{p}{W^2},
\]
so this event contributes at most $Cp/W^2\le C\sqrt p/W$ to the $L^p$
norm, for $p\le\epsilon_0W$. Altogether,
\[
 \|g_k-a\|_{L^p}
 \le\|g_k-\E g_k\|_{L^p}+\|\E g_k-a\|_{L^p}
 \le C\frac{\sqrt p}{W},
\]
which is the residual-pair assertion.
\end{proof}

\section{Counting fluctuations and energy transport}
\label{sec:rank-location}
An estimate proved at a deterministic energy cannot be applied directly
at the random eigenvalue $\lambda_\alpha$ of a fixed rank. This section locates
$\lambda_\alpha$ precisely, so that the random energy can be replaced by its
classical location at a small cost. It is used in Part~I, for the physical
center and the angular comparison.
Theorem~\ref{thm:rank-location} bounds the deviation of the eigenvalue
counting function from its classical value by $C_p\sqrt W$ in $L^p$, and its
mean deviation by $C\{n/W+W^{2/5}\log(2W)\}$. As a consequence, the
eigenvalue of a bulk rank lies within $O((n\sqrt W)^{-1})$ of its classical
location, with no extra power or logarithmic factor.
In Subsection~\ref{subsec:s6-averaging-the-final-gaussian-block} we replace the sharp eigenvalue count of the last
Gaussian block by a smooth boundary phase (Lemma~\ref{lem:fresh-count}). In
Subsection~\ref{subsec:s6-the-compensated-two-port-phase} we show that adding this boundary phase cancels the drift
that would otherwise accumulate along the path
(Lemma~\ref{lem:compensated-phase}). In Subsection~\ref{subsec:s6-integration-to-a-random-energy} we bound an
integral up to a random energy by the typical displacement times a
fixed-energy moment bound (Lemma~\ref{lem:random-energy}).

Let $N_H(E)=\#\{\lambda_j<E\}$, let $F_{\rm sc}$ be the semicircle
distribution function, and put
$e_\alpha=F_{\rm sc}^{-1}(\alpha/(nW))$.

\begin{theorem}\label{thm:rank-location}
For fixed $b,T$, a compact bulk interval $I$, and every fixed $p<\infty$,
uniformly for $4\le n\le TW$ and deterministic $E\in I$,
\begin{align}
 \|N_H(E)-nWF_{\rm sc}(E)\|_{L^p}&\le C_p\sqrt W,
 \label{eq:physical-count}\\
 |\E N_H(E)-nWF_{\rm sc}(E)|
 &\le C\{n/W+W^{2/5}\log(2W)\}.
 \label{eq:physical-count-mean}
\end{align}
For every fixed $\vartheta\in(0,1/2)$,
\begin{equation}\label{eq:rank-location}
 \sup_{\vartheta\le\alpha/(nW)\le1-\vartheta}
 \|\lambda_\alpha-e_\alpha\|_{L^p}
 \le \frac{C_{b,T,\vartheta,p}}{n\sqrt W}.
\end{equation}
These eigenvalue-only statements also hold under each normalized
block-mass or pair-mass law.
\end{theorem}

The location scale is much larger than an individual spacing at
critical bandwidth. Its role is to pay an energy integral of a
signed response; no microscopic rigidity claim is made here.

\subsection{Averaging the final Gaussian block}\label{subsec:s6-averaging-the-final-gaussian-block}

 We replace the sharp eigenvalue count of the last block by a smooth
boundary phase. Lemma~\ref{lem:fresh-count} bounds the comparison and the
remaining count fluctuations uniformly in the old deformation. We
interpolate between the Gaussian block and the matched boundary value,
then control smoothing and fluctuations by counting and Gaussian
concentration.

 For Hermitian $F$, define
$\Theta_z(F)=\pi^{-1}\sum_j\arg(z-\lambda_j(F))$, using
$0<\arg(z-\lambda_j)<\pi$. We retain the negative-imaginary
parameter $a$ from \eqref{eq:matched-parameters}.

\begin{lemma}\label{lem:fresh-count}
Let $A\sim\GUE_W(c/W)$ be independent of deterministic $F=F^*$,
put $z_0=-E-ca$ and
$d=|W^{-1}\Tr(z_0-F)^{-1}-a|$. For $d\le d_*$,
\begin{equation}\label{eq:fresh-count}
 \left|\E_A N_-(A-E-F)-\Theta_{z_0}(F)
                       -\frac{cW}{2\pi}\Im a^2\right|
 \le C\{Wd^2+W^{2/5}\log(2W)\}.
\end{equation}
For every $F$ the conditional centered count has $L^p$ norm at
most $C_p\sqrt W$. Constants are uniform in $\norm F$.
\end{lemma}

The fresh Gaussian comparison is uniform in the old deformation. It separates the sharp count from a conditional mean whose phase can be matched along the chain.

\begin{proof}
Set $\varepsilon=W^{-3/5}$ and
\[
 H_t=\sqrt t\,A-E-F-(1-t)ca+i\varepsilon,\quad
 g_t=W^{-1}\Tr H_t^{-1},\quad
 y_t=(1-t)c(-\Im a)+\varepsilon.
\]
The upper-half-plane logarithm has the exact heat derivative
\begin{equation}\label{eq:fresh-log-heat}
 \frac{d}{dt}\E\Tr\log H_t
 =cW\E(ag_t-g_t^2/2)
 =\frac{cWa^2}{2}-\frac{cW}{2}\E(g_t-a)^2.
\end{equation}
For clarity, the Gaussian second derivative of $\Tr\log H$ in
a direction $K$ is $-\Tr(H^{-1}KH^{-1}K)$. Contracting two GUE
increments gives $-(c/W)(\Tr H^{-1})^2$; the heat factor is $1/2$.
The deterministic derivative contributes $ca\Tr H^{-1}$.
These are exactly the two terms in \eqref{eq:fresh-log-heat}.

For $t\ge1/2$, the Gaussian counting bound gives
$\E\Tr|H_t^{-1}|^2\le CW/y_t$.
Gaussian Poincar\'e, using
$\norm{H_t^{-1}KH_t^{-1}}_{\HS}
\le\norm K\,y_t^{-1}\norm{H_t^{-1}}_{\HS}$, yields
\begin{equation}\label{eq:fresh-count-variance}
 \E|g_t-\E g_t|^2\le C/(W^2y_t^3).
\end{equation}
The same bound holds for a bounded deterministic matrix insertion.
For $t\le1/2$ the deterministic height is bounded below, so the
bound follows directly.

The scalar equation
$m_t=W^{-1}\Tr[-E-F-(1-t)ca+i\varepsilon-tcm_t]^{-1}$
has a solution within $C(d+\varepsilon)$ of $a$, uniformly in $t$.
Indeed at $\varepsilon=0$, writing $R=(z_0-F)^{-1}$,
\[
 \Re(1-tcW^{-1}\Tr R^2)
 =1-tcW^{-1}\Tr|R|^2+
             2tcW^{-1}\Tr(\Im R)^2
 \ge1-t+2tc(\Im a)^2-Cd.
\]
We used $W^{-1}\Tr|R|^2=1/c+O(d)$ and Jensen's inequality
for the last trace. The right side is uniformly positive for small
$d$. All higher derivatives are bounded at the positive test height.
The implicit-function theorem and the scalar imaginary-part
uniqueness argument identify this local branch.

To compare $\bar g_t=\E g_t$ with it, diagonalize $F$ and use
Gaussian integration by parts:
\[
 \bar g_t-\av{M(\bar g_t)}
 =tc\,\E[(g_t-\bar g_t)\av{M(\bar g_t)H_t^{-1}}],
 \quad
 M(u)=[-E-F-(1-t)ca+i\varepsilon-tcu]^{-1}.
\]
While $\bar g_t$ is near $a$, $M(\bar g_t)$ is uniformly bounded.
Apply \eqref{eq:fresh-count-variance} to both traces after centering
the second. The residual is $C/(W^2y_t^3)$.
Stability and continuation from $t=0$ give
$|\bar g_t-a|\le C(d+\varepsilon+W^{-2}y_t^{-3})$.
The continuation stays in the stable neighborhood because
$W^{-2}\varepsilon^{-3}=W^{-1/5}\to0$.
Consequently
\[
 \E|g_t-a|^2
 \le C(d^2+\varepsilon^2+W^{-2}y_t^{-3}).
\]
Integrating \eqref{eq:fresh-log-heat} costs at most
$C(Wd^2+W\varepsilon^2+(W\varepsilon^2)^{-1})$.
The initial phase changes by $O(W\varepsilon)$ when removing its
extra height. At $t=1$, replacing the smoothed phase by the sharp
negative count costs $CW\varepsilon\log(2/\varepsilon)$ in expectation:
the scalar kernel error is $C\min(1,\varepsilon/|x|)$; integrate
against the expected density on $|x|\le1$ and use $C\varepsilon$
per eigenvalue outside. Substitution of $\varepsilon=W^{-3/5}$
proves \eqref{eq:fresh-count}.

Finally expose the independent upper-triangular row groups of $A$.
Replacing one group changes the matrix by rank at most two and
therefore changes its negative count by at most two. The bounded
martingale-difference inequality gives
$\|N_--\E_A N_-\|_{L^p}\le C_p\sqrt W$.
It also gives variance at most $2W$ by Efron--Stein.
\end{proof}

\subsection{The compensated two-port phase}\label{subsec:s6-the-compensated-two-port-phase}

 The eigenvalue count alone carries a boundary contribution that would
accumulate along the path. Lemma~\ref{lem:compensated-phase} shows that
adding the boundary phase cancels this drift and leaves controlled
fluctuations. We compute the complete logarithmic increment, center its
conditional mean, and sum the resulting fluctuations before closing
the path.

 Write
\[
 \ell=\pi^{-1}\Im\log a^{-1},\qquad
 h=\pi^{-1}\Im a^2,\qquad F_{\rm sc}(E)=\ell+h/2.
\]
For a physical path of $r$ blocks with its two independent ports put
$\Psi_r=N_r(E)+\pi^{-1}\Im\Tr\log(z-X_r)$.

\begin{lemma}\label{lem:compensated-phase}
For every fixed $p$ and $1\le r\le TW$,
\begin{equation}\label{eq:compensated-phase}
 \|\Psi_r-\E\Psi_r\|_{L^p}\le C_p\sqrt{r+1},\qquad
 |\E\Psi_r-rWF_{\rm sc}(E)-W(2\ell+bh)|
 \le C(r+1)/W.
\end{equation}
\end{lemma}

The phase increment is centered using its complete conditional mean. This cancellation prevents a drift proportional to the full number of blocks.

\begin{proof}
We calculate a complete logarithm so that all real Schur poles cancel.
Given $X$, interpolate one separate-port layer by
\[
 \mathcal H_t=
 \begin{pmatrix}
 \sqrt t\,A-E-X-(1-t)(1-b)a&\sqrt t\,B\\
 \sqrt t\,B^*&(a^{-1}+tba)I
 \end{pmatrix}.
\]
For a common port use $B=(B_1^T,B_2^T)^T$ and bottom
parameter $a^{-1}+2tba$. At $t=0$ the top block is $z-X$.
For $t\le1/2$ both deterministic heights are positive and bounded
below. For $t\ge1/2$, the proof of
Lemma~\ref{lem:gaussian-absorption} applies at variances $tc,tb$;
the extra nonnegative top height can only improve its imaginary
quadratic-form estimate. Thus all inverse moments and covariance
bounds are uniform in $t$.

The mean equations along this interpolation have the constant
comparison value $a$. Their linearized top coefficients are
$t(c+tb^2a^2)Q$ for separate ports and
$tQ(cI+tb^2a^2\one\one^T)$ for the common port.
Their absolute row sums are bounded respectively by
$(c+b^2)/(1-b)+Cd<1$ and $(c+2b^2)/(1-b)+Cd<1$.
The same stable-branch comparison as in
Lemma~\ref{lem:layer-contraction} gives trace errors
$O_{L^2}(d+W^{-1})$ uniformly in $t$.

For separate ports the heat derivative of the expected logarithm is
\[
 W\E\left[
 (1-b)a(g_1+g_2)+ba(u_1+u_2)
 -\frac{c}{2}(g_1^2+g_2^2)-b(g_1u_1+g_2u_2)\right].
\]
For the common port the corresponding expression is
\[
 W\E\left[
 (1-b)a(g_1+g_2)+2bau
 -\frac{c}{2}(g_1^2+g_2^2)-b(g_1+g_2)u\right].
\]
Both have constant term $Wa^2$. Every linear term in the
deviations from $a$ vanishes because $c+2b=1$; their expected
quadratic remainder is $O(Wd^2+W^{-1})$.
The initial extra port logarithm is $2W\log a^{-1}$ for separate
ports and $W\log a^{-1}$ for a common port. Schur's determinant
identity, with continuous logarithms, therefore gives
\begin{align*}
 \E(\Psi_{r+2}-\Psi_r\mid X_r)
 &=2WF_{\rm sc}(E)+O(Wd(X_r)^2+W^{-1}),\\
 \E(\Xi_{r+2}-\Psi_r\mid X_r)
 &=W(\ell+h)+O(Wd(X_r)^2+W^{-1}),
\end{align*}
where $\Xi$ denotes the common-port compensated count.
The integer phase in each identity is fixed by Schur inertia,
not by a choice of a real logarithm modulo $2\pi$.

The conditionally centered increments have bounded fixed $L^p$
norms. Indeed their logarithmic-phase gradient is the imaginary
part of $\Tr\mathcal H_1^{-1}d\mathcal H_1$, of
Hilbert--Schmidt norm at most $C\sqrt W\norm{\mathcal H_1^{-1}}$;
Gaussian covariance supplies $W^{-1}$. Its continuous phase is
locally Sobolev across a Schur pole: the jump of the real negative
count is exactly the opposite of the jump of the boundary logarithm.
One can equivalently use the full, invertible absorbing matrix to
define this phase, then apply smooth cutoffs to its inverse.
The inverse moments justify removing those cutoffs.

By Theorem~\ref{thm:matched-traces}, the predictable error has
$L^p$ norm $C_p/W$ per step. On $d>d_*$ the count and phase
increments are bounded by $CW$, and arbitrarily high fixed
moments of $d$ pay this event at the same rate. The martingale
moment inequality bounds the sum of centered increments by
$C_p\sqrt{r+1}$, while the accumulated predictable error is
$C_pr/W\le C_{p,T}\sqrt{r+1}$.

It remains to fix the starting constant. In the $r=1,2$ core,
interpolate each Gaussian entry by $\sqrt t$ and each deterministic
diagonal by $a^{-1}+t s_i a$, where $s_i$ is its variance row sum.
The constant mean is $a$ and the finite core stability was proved
in Theorem~\ref{thm:matched-traces}. The heat derivative has constant
$(W/2)a^2\sum_i s_i=(W/2)a^2(r+2b)$ and an $O(W^{-1})$
remainder. At $t=0$ there are $(r+2)W$ copies of
$\log a^{-1}$. Thus
$$\E\Psi_r=W[(r+2)\ell+(r+2b)h/2]+O(W^{-1}).$$
The fluctuation of the core is bounded as for the increments. For $r=1,2$
let $\mathcal H$ be the matrix of the core together with its two ports,
with $-E$ on the physical blocks and $z$ on the ports. By Schur's
determinant identity, with the continuous logarithms used above, $\Psi_r$
is $\pi^{-1}\Im\log\det\mathcal H$ up to a deterministic constant, locally
Sobolev across Schur poles as for the increments. Its gradient in the Gaussian coordinates of the core, including the internal
bond when $r=2$, is $\pi^{-1}\Im\Tr\mathcal H^{-1}d\mathcal H$. All these
coordinates have variance at most $1/W$, so the weighted squared gradient
is at most $CW^{-1}\norm{\mathcal H^{-1}}_{\HS}^2\le C\norm{\mathcal H^{-1}}^2$.
Lemma~\ref{lem:gaussian-absorption} bounds the moments of
$\norm{\mathcal H^{-1}}$: for $r=1$ directly, with one physical block and
two ports, and for $r=2$ conditionally on the internal bond, which enters
as the Hermitian matrix $X$, uniformly in it. The Gaussian $L^p$
Poincar\'e inequality therefore gives $\|\Psi_r-\E\Psi_r\|_{L^p}\le C_p$
for $r=1,2$; this is the initial term of the martingale bound.
This initializes the recursion and proves \eqref{eq:compensated-phase}.
\end{proof}

\begin{proof}[Proof of Theorem~\ref{thm:rank-location}]
Delete block zero. Grow the middle $n-3$ blocks with separate
ports, then append its two end blocks and their actual bonds to
block zero as the common-port layer above. The result satisfies
\[
 \|\Xi_{n-1}-\E\Xi_{n-1}\|_{L^p}\le C_p\sqrt n,\qquad
 \E\Xi_{n-1}
 =(n-1)WF_{\rm sc}+W(\ell+bh)+O(n/W).
\]
The common-port trace has residual $O_{L^q}(W^{-1})$ for every
fixed $q$. Lemma~\ref{lem:fresh-count} therefore applies to the
remaining independent $A_0$ with $F=\Sigma$; outside its regular
set, the deterministic $CW$ count/phase bound and high moments
of the residual pay the error. Schur inertia gives exactly
$N_H=N_{n-1}+N_-(A_0-E-\Sigma)$.
The deterministic constants satisfy
\[
 (n-1)F_{\rm sc}+\ell+bh+ch/2=nF_{\rm sc},
 \qquad b+c/2=1/2.
\]
The fresh centered count costs $C_p\sqrt W$ and its mean error
costs $C_pW^{2/5}\log(2W)$. Since $n\le TW$, this proves
\eqref{eq:physical-count} and \eqref{eq:physical-count-mean}.

Monotonicity at deterministic energies $e_\alpha\pm t$, where the
semicircle density has a fixed positive lower bound, gives
$\Prob(|\lambda_\alpha-e_\alpha|>t)
\le C_q(n\sqrt W\,t)^{-q}$ for $t\ge C/(n\sqrt W)$ within
a compact bulk interval. Integrate with $q>p$.
For the remaining tails use this bound at the endpoints of a
slightly larger bulk interval and H\"older with $q>2p$.
All fixed matrix-norm moments are uniformly bounded: the block
norm tails and
$\norm H\le\max_i\norm{A_i}+2\max_i\norm{B_i}$
give this for $n\le TW$. This proves \eqref{eq:rank-location}.
Finally a function of the full spectrum is invariant under cyclic
block shifts, so $\E[np_kf(\operatorname{Spec}H)]=\E f$;
the same equality holds for the pair bias.
\end{proof}

\subsection{Integration to a random energy}\label{subsec:s6-integration-to-a-random-energy}

  We need to integrate an estimate up to a random energy that may depend
on the integrand. Lemma~\ref{lem:random-energy} bounds this integral by
the typical displacement times a fixed-energy moment bound. We split
the possible displacement into successively larger intervals and use
the tail of the terminal energy to sum their contributions.

\begin{lemma}\label{lem:random-energy}
Suppose $\|\Lambda-e\|_{L^q}\le C_q\omega$ and
$\sup_E\|D(E)\|_{L^{2p}}\le K$, with $p\geq 1$ and $q>2p$.
Then
\begin{equation}\label{eq:random-energy}
 \left\|\int_e^\Lambda D(E)\,dE\right\|_{L^p}
 \le C_{p,q}\omega K.
\end{equation}
For an expectation of the absolute integral (the case $p=1$) it suffices to use
$\sup_E\|D(E)\|_{L^2}\le K$ and $q>2$.
\end{lemma}

The integrand and the terminal energy may be dependent. Fixed moments and a dyadic tail decomposition are sufficient for the stated integral bound.

\begin{proof}
Write the oriented integral as an integral over the line with an
indicator that $\Lambda$ reaches the integration point. Inside
distance $\omega$, Minkowski's inequality, valid for $p\ge1$, gives at most $2\omega K$.
On the annulus of radii $2^j\omega,2^{j+1}\omega$, H\"older and
Markov give an $L^p$ integrand bound
$C_qK2^{-jq/(2p)}$. The annulus length is at most
$2^{j+2}\omega$. The resulting geometric sum converges since
$q>2p$. For the expectation of the absolute integral,
Cauchy--Schwarz gives the annular factor $2^{-jq/2}$ instead;
its length-weighted sum converges for $q>2$.
No independence is used.
\end{proof}
\section{Weighted derivatives and the physical center}
\label{sec:clock-center}
In Part~I, the angular comparison compares the mass difference on a pair
with a deterministic center. This section proves that the center, evaluated
at the random eigenvalue and under the pair law, is small. It provides the
second of the two estimates that give quantum unique ergodicity.
Theorem~\ref{thm:physical-center} states that
$\E^\sharp_{\alpha,e}|j_e(\lambda_\alpha)|^2\le C/(\sqrt n\,W^{3/2})$ at every
deterministic bulk rank. In words, the center at the random energy is small
on average, with a bound that is sufficient for the Fourier argument of
Section~\ref{sec:angular-que}.
In Subsection~\ref{subsec:s7-a-weighted-ward-identity} we propagate a positive weight through one layer
without bounding its largest eigenvalue (Lemma~\ref{lem:weighted-layer}). In
Subsection~\ref{subsec:s7-full-and-tagged-clocks} we show that the full clock grows at most linearly while a
local tag has bounded moments (Lemma~\ref{lem:clock-moments}). In
Subsection~\ref{subsec:s7-differentiating-the-trace-equations} we differentiate the trace equations when a source is
inserted (Lemmas~\ref{lem:t26-differentiated-covariance}
and~\ref{lem:tagged-trace}). In Subsection~\ref{subsec:s7-a-smooth-center-through-real-schur-poles} we construct a smooth
center that remains usable through real Schur poles
(Lemma~\ref{lem:smooth-center}). In Subsection~\ref{subsec:s7-the-gaussian-energy-primitive} we integrate in energy
up to the selected eigenvalue and use Gaussian integration by parts to prove
the theorem.

Use a simultaneous independent splitting
\begin{equation}\label{eq:simultaneous-splitting}
 H=V+\sum_eG_e,\qquad
 G_e\sim\GUE_{2W}(s/(2W)),\qquad s=\min(c,b),
\end{equation}
where $G_e$ is supported on the two blocks of edge $e$.
This is possible because a diagonal block receives total variance
$s/W\le c/W$ and a bond receives $s/(2W)\le b/W$.
All residual entries of $V$ and all $G_e$ are independent.
For a fixed edge $e$, put $H^{(e)}=H-G_e$.
The law of $H$ is exactly the original law.

At a fixed real $E$, let $K_e(E)$ be the Schur complement of
$H^{(e)}-E$ on $S_e$, and put   $\Naux =2W$, $\mu=\bar a$.
On the regular domain let $\zeta_e$ be the upper-half-plane root of
  $\zeta=s\Naux ^{-1}\Tr(K_e-\zeta)^{-1}$. Define

\begin{equation}\label{eq:dyson-center}
 j_e(E)=\frac{s}{\Naux }\Tr J_e|(K_e-\zeta_e)^{-1}|^2,\qquad
 J_e=P_i-P_{i+1}.
\end{equation}

The value on the exceptional domain will be a bounded, equivariant
extension. The regular domain is specified below.

\begin{theorem}[The center at the selected eigenvalue]\label{thm:physical-center}
Fix $b$, $T$ and the bulk margin $\vartheta$, and let $4\le n\le TW$. Let
$\alpha$ be a deterministic bulk rank, and let $e$ be an edge of the
ring, that is, an adjacent pair of blocks. Consider the regular center
$j_e$ of \eqref{eq:dyson-center}, extended by zero outside its regular
domain. Then, under the exact pair law \eqref{eq:pair-law}, the center
evaluated at the selected eigenvalue $\lambda_\alpha$ is small in mean
square:
\begin{equation}\label{eq:physical-center}
 \E^\sharp_{\alpha,e}|j_e(\lambda_\alpha)|^2
 \le\frac{C_{b,T,\vartheta}}{\sqrt n\,W^{3/2}} .
\end{equation}

The regular domain can be described explicitly. Let $I$ be a fixed
compact interval in $(-2,2)$ that contains a neighbourhood of all
classical locations $e_\alpha$ with $\vartheta\le\alpha/(nW)\le1-\vartheta$,
and fix $\epsilon_{\rm reg}>0$. With the trace errors $\epsilon_{e,i}$ of
Lemma~\ref{lem:tagged-trace}, put
\[
 \epsilon_e^2=\sum_{i\in e}|\epsilon_{e,i}|^2,
 \qquad
 D_e(E)=\bigl\{E\in I,\ \epsilon_e(E)\le\epsilon_{\rm reg}\bigr\}.
\]
On this domain the following hold, with fixed constants.
\begin{enumerate}[label=(\alph*)]
\item The deformation $F=K_e(E)+E$ satisfies the regularity condition
\eqref{eq:fresh-regularity}, with $c_0=s$.
\item The Dyson root is $\zeta_e=sm_s(E)$, where $m_s$ is the Stieltjes
transform of the Dyson solution formed with this $F$. It has a fixed
positive imaginary part and a bounded modulus.
\item The matrix $s|K_e-\zeta_e|^{-2}$ is the covariance $\Gamma$ used
both in Theorem~\ref{thm:fresh-palm} and in
\eqref{eq:angular-reference}.
\end{enumerate}
Finally, the selected eigenvalue lies outside the regular domain only
with very small probability: for every fixed $q$,
\[
 \Prob^\sharp_{\alpha,e}\bigl(D_e(\lambda_\alpha)^c\bigr)\le C_qW^{1-q}.
\]
\end{theorem}

Let us decode the theorem in simple words. The center is the deterministic quantity with which the
angular comparison of Section~\ref{sec:angular-que} compares the mass
difference of the eigenvector on the two blocks of the pair. At a fixed
deterministic energy it is of order $W^{-1}$, so its square is of order
$W^{-2}$. The theorem shows that, even when it is evaluated at the random
selected eigenvalue and under the pair law, which weights the matrix law
by the eigenvector's mass on the pair, its mean square is at most of
order $n^{-1/2}W^{-3/2}$. This is weaker than $W^{-2}$ when $n\ll W$, but
it is enough for the Fourier argument in the proof of
Theorem~\ref{thm:que}. The second part of the statement says that the
regular domain on which the center is defined is the set of bulk
energies where the trace errors of the pair are small. There the fresh
Gaussian block is regular, its Dyson root and covariance are the ones
used in the moment theorem and in the angular measure, and the selected
eigenvalue falls outside this domain only with probability smaller than
any fixed power of $W$.

\subsection{A weighted Ward identity}\label{subsec:s7-a-weighted-ward-identity}

 We need to propagate a positive weight without assuming that its
largest eigenvalue is bounded. Lemma~\ref{lem:weighted-layer} controls
the mean and fluctuations of its outgoing trace in terms of its incoming
mass. We retain the weight inside the inverse identities and Gaussian
contractions, so that the estimate can be repeated along the path.

 For one separate-port step write
\[
 \mathcal G=\begin{pmatrix}A-E-X&B\\B^*&zI\end{pmatrix}^{-1},
 \quad K=z-X,\quad L=K\mathcal G_{12}.
\]
If $S\succeq0$ is a source arriving from the old path, let
$Y=R^*SR$, $R=K^{-1}$ and $u=W^{-1}\Tr Y$.

\begin{lemma}\label{lem:weighted-layer}
For an independent Hermitian input $X$, the weighted factors
$K\mathcal G_{11}$, $\mathcal G_{11}K^*$ and $L$ have uniformly
bounded fixed moments for $W\ge W_0(p)$.
On the scalar domain $d(X)\le d_*$,
\begin{equation}\label{eq:weighted-ward}
 \norm{\E(LL^*\mid X)-I}\le C(d(X)+W^{-2}).
\end{equation}
For $Y\succeq0$ fixed with the past,
\begin{align}
 \E[W^{-1}\Tr L^*YL\mid X,Y]
 &\le(1+Cd+CW^{-2})u,\label{eq:clock-mean}\\
 \|W^{-1}\Tr L^*YL-\E[W^{-1}\Tr L^*YL\mid X,Y]\|_{p\mid X,Y}
 &\le\frac{C_p}{W}\sqrt{W^{-1}\Tr Y^2}
 \le C_pu/\sqrt W .\label{eq:clock-noise}
\end{align}
Without the scalar-domain restriction the outgoing trace has
conditional $L^p$ norm at most $C_pu$.
\end{lemma}

The estimate keeps the old positive weights inside the Gaussian contraction. It can therefore be iterated without an operator-norm bound on the old clock.

\begin{proof}
The top inverse is
$\mathcal G_{11}=(A-E-X-z^{-1}BB^*)^{-1}$. Thus
\begin{align*}
 K\mathcal G_{11}
 &=I+[(z+E)I-A+z^{-1}BB^*]\mathcal G_{11},\\
 \mathcal G_{11}K^*
 &=I+\mathcal G_{11}[(\bar z+E)I-A+z^{-1}BB^*].
\end{align*}
These identities and Lemma~\ref{lem:gaussian-absorption} give the
weighted moments. The top-right block equation similarly gives
$K\mathcal G_{12}=[(z+E)I-A]\mathcal G_{12}-B\mathcal G_{22}$.

The exact Ward identity is
$\mathcal G_{12}\mathcal G_{12}^*=-\eta^{-1}\Im\mathcal G_{11}$.
Let
$M=[-E-X-\diag(cm_i+bu_i)]^{-1}$
be the deterministic top mean inverse from
\eqref{eq:two-port-dyson}. At the matched solution
$KM=I+O(d)$, and
$-\Im(cm_i+bu_i)=\eta+O(d)$.
Consequently
$-\eta^{-1}K(\Im M)K^*=I+O(d)$.
It remains to estimate the weighted random mean error.

In the exact top integration-by-parts equation keep
$\mathcal G_{11}K^*$ inside each covariance:
\[
 K(\E\mathcal G_{11}-\bar M)K^*
 =K\bar M\sum_j\sigma_j
       \operatorname{Cov}(g_j,\mathcal G_{11}K^*).
\]
Here $\bar M$ uses the actual scalar means and the finite sum runs
over the incident physical and port coordinates. The contraction in
Lemma~\ref{lem:trace-matrix} still gives $W^{-2}$; its rightmost
factor remains the weighted inverse
$\mathcal G_{11}K^*$ throughout differentiation.
Its moments were just proved, as were those of $K\bar M$ on this
domain. Replacing the scalar means by their deterministic values
also costs $O(W^{-2})$ after these weights. Taking the imaginary
part in the Ward identity proves \eqref{eq:weighted-ward}.
Trace pairing with $Y$ proves \eqref{eq:clock-mean}.

For the fluctuation, differentiate
$W^{-1}\Tr YLL^*$ before estimating. Every term contains $Y$
between bounded-moment weighted inverse factors; cyclicity and
Hilbert--Schmidt Cauchy--Schwarz give gradient norm at most
$CW^{-1}\norm Y_{\HS}$ times such factors.
The Gaussian covariance contributes $W^{-1/2}$ in the gradient
inequality. This is
$CW^{-3/2}\norm Y_{\HS}
=CW^{-1}\sqrt{W^{-1}\Tr Y^2}$, proving
\eqref{eq:clock-noise}. Finally
$W^{-1}\Tr L^*YL\le\norm L^2u$ proves the unrestricted moment
bound. Positivity gives $\Tr Y^2\le(\Tr Y)^2$, which proves its
last displayed inequality.
\end{proof}

\subsection{Full and tagged clocks}\label{subsec:s7-full-and-tagged-clocks}

 We distinguish the response to energy on every block from the response
to a perturbation on only one or two blocks. Lemma~\ref{lem:clock-moments}
shows that the full response grows at most linearly, while a local tag
has bounded trace moments. We derive their common positive recursion,
then use the weighted one-layer estimate and control the exit event.

 For a nonnegative scalar perturbation on a physical block,
differentiate its exact Schur complement. If $S=-\dot X\succeq0$,
the outgoing source is
\begin{equation}\label{eq:tagged-clock}
 Y_+=\mathcal G_{12}^*P\mathcal G_{12}+L^*YL,
\end{equation}
where $P$ is the projection onto the newly tagged blocks.
The signs follow from
$\dot{(A-E-X)}=P+S$ and
$\dot X_+=-B^*(A-E-X)^{-1}(P+S)(A-E-X)^{-1}B$.
For the full physical energy derivative, $S=X'\succeq0$ and
$P=I$ at every step. For an edge tag, $P=0$ except at at most
two steps.

\begin{lemma}\label{lem:clock-moments}
For every fixed $p<\infty$, $W\ge W_0(p)$  and $r\le TW$, the full clock has
$\|W^{-1}\Tr Y_r\|_{L^p}\le C_p(r+1)$.
For a tag on at most two physical blocks,
\begin{equation}\label{eq:tag-clock-moments}
 \sup_{r\le TW}\|W^{-1}\Tr Y_r\|_{L^p}\le C_p.
\end{equation}
These are full prior moments, including the scalar exit event.
\end{lemma}

All moment orders here are fixed unless the growing range is explicitly stated. The homogeneous tagged clock has no accumulated additive source, unlike the full clock.

\begin{proof}
Write $u=W^{-1}\Tr Y$ and $\delta=d(X)$.
At an untagged step Lemma~\ref{lem:weighted-layer} gives conditional
mean at most $(1+C\delta+C/W^2)u$ and centered $q$th moment at
most $C_qu/\sqrt W$. At an injection replace $u$ by $1+u$
and add a constant to the mean. Expanding $u_+^p$ around its mean
for integer $p\ge2$, the linear centered term vanishes and all
other fluctuation terms cost $C_pu^p/W$. Thus an untagged step obeys
$\E(u_+^p\mid X,Y)\le(1+C_p\delta+C_p/W)u^p$.

The scalar recursion is
$\delta_+\le q_*\delta+C/W^2+|\xi_+|$,
$\|\xi_+\|_{2\mid X}\le C/W$. Its dependence on the clock is
retained:
\[
 \E(|\xi_+|u_+^p\mid X,Y)
 \le \|\xi_+\|_{2\mid X}\|u_+^p\|_{2\mid X,Y}
 \le C_pu^p/W.
\]
Choose $K_p(1-q_*)>2C_p$ and then reduce the fixed scalar radius.
The polynomial inequality
$(1+K_pq_*\delta)(1+C_p\delta+C_p/W^2)
\le(1+K_p\delta)(1+C'_p/W^2)$
holds there. Multiplication by the linear weight $1+K_p\delta_+$
and the preceding correlated estimate give
\[
 M_+\le(1+C_p/W)M,\qquad
 M=\E[\one_{\rm regular}(1+K_p\delta)u^p].
\]
At an injection it gives $M_+\le C_p(1+M)$.
There are at most two injections, so iteration for $r\le TW$
proves a uniform bound. The finite core injection moments follow
from Lemma~\ref{lem:gaussian-absorption}.
For the full clock the same argument instead gives
$M_+\le(1+C_p/W)M+C_pM^{1-1/p}+C_p$.
Its $p$th root grows by at most a fixed constant per step, up to
the factor $1+C_p/W$; discrete iteration yields $C_p(r+1)^p$.

The stopped scalar process has exit probability
$e^{-cW\log W}$ by \eqref{eq:scalar-exit}.
Without stopping, the fresh weighted norm bounds give crude
$2p$th moments at most $C_p^{r+1}$, with the finitely many
injections included. H\"older bounds the omitted part by
$C_p^{TW+1}e^{-cW\log W/2}$, which tends to zero.
This removes the stop. Real $p\ge1$ follows by monotonicity.
\end{proof}

\subsection{Differentiating the trace equations}\label{subsec:s7-differentiating-the-trace-equations}

  We must control how a boundary trace changes when a source is inserted
in the old path. Lemma~\ref{lem:t26-differentiated-covariance} bounds the
differentiated Gaussian covariance, and Lemma~\ref{lem:tagged-trace}
propagates the resulting trace estimates. We differentiate the exact
mean equations, retain every source factor, and then iterate; merely
differentiating  an error bound  would not justify   the argument.

\begin{lemma}[Differentiated conditional covariance]
\label{lem:t26-differentiated-covariance}
Condition on the old two-port path and its positive source $S$, before
exposing the next Gaussian layer. Put $K=z-X$, $Y=K^{-*}SK^{-1}$,
$u=W^{-1}\Tr Y$, and let $\widehat S$ embed $S$ in the old-input
block of the augmented layer matrix. Differentiate
$X(h)=X-hS+O(h^2)$ at zero, keeping $E,z$ and every fresh Gaussian
coordinate fixed. For a normalized block trace $g_i=W^{-1}\Tr P_i\mathcal G$,
\begin{align}
 \|\partial_h\operatorname{Cov}(g_i,\mathcal G)\|
                           &\le CW^{-2}\|Y\|\le Cu/W,
                                      \label{eq:t26-covariance-derivative}\\
 \|\dot g_i-\E(\dot g_i\mid X,S)\|_{p\mid X,S}
                           &\le C_pW^{-3/2}\|Y\|_{\rm HS}
                            \le C_pu/\sqrt W.
                                      \label{eq:t26-derivative-noise}
\end{align}
The constants are uniform in the old Hermitian deformation and source.
The corresponding assertions with a bounded deterministic insertion
added to $\widehat S$ have $1+u$ in place of $u$.
\end{lemma}

\begin{remark} In Lemma~\ref{lem:layer-contraction}, we
derived the mean equations of one fresh layer by Gaussian integration
by parts. Their error terms are covariances $\operatorname{Cov}(g_i,\mathcal G)$
between a normalized trace and the layer resolvent, and
Lemma~\ref{lem:trace-matrix} shows that these are $O(W^{-2})$. To follow
a source through the path, we now differentiate the mean equations in
the strength $h$ of the source. So we need two facts about the
derivatives. First, the derivative of each covariance error must still
be small. Second, the derivative of each trace must concentrate around
its conditional mean.
 The main point is the size in which we measure the source. The old
source $S$ can have a huge operator norm, because it comes from a long
path that we do not control. What the new layer actually sees, however,
is $S$ transported through the old port resolvent $K^{-1}$, that is,
the matrix $Y=K^{-*}SK^{-1}$, and its normalized trace $u$. Both bounds
of the lemma are proportional to $u$. This is what allows us to iterate
the lemma along the whole path without any bound on $\|S\|$.
\end{remark}

\begin{proof}
The source enters the fresh layer only through the
old input matrix, so differentiating in its strength simply inserts one copy
of the source between two resolvents (Step~1). We write the covariance between a trace and the
resolvent as an exact integral over a Gaussian interpolation, which couples
the layer with an independent copy of itself (Step~2). Differentiating this
formula gives four terms, one for each resolvent in which the source can be
inserted, as in Figure~\ref{fig:dc-insertions} (Step~3). The key point is
that the source never appears alone. It is always sandwiched between weighted
resolvent factors with bounded moments, so each insertion costs only the size
of the transported source and not its possibly huge norm (Step~4). Step~5
bounds the fluctuation of the derivative by the Gaussian Poincar\'e
inequality, which controls the variance of a function of Gaussian variables
by the size of its gradient, and Step~6 justifies the differentiations.

\emph{Step 1:} We demonstrate first how source enters the layer. The layer resolvent
$\mathcal G$ depends on the source only through the old input
$X(h)=X-hS+O(h^2)$, which sits in the old-input block. Since
$\partial_h\mathcal G^{-1}=\widehat S$ at $h=0$, the usual formula for
the derivative of an inverse gives
\begin{equation}\label{eq:dc-resolvent-derivative}
 \dot{\mathcal G}=-\mathcal G\,\widehat S\,\mathcal G.
\end{equation}
Thus differentiating in $h$ replaces one factor $\mathcal G$ by
$-\mathcal G\widehat S\mathcal G$. In other words, it inserts the source
$\widehat S$ once, between two resolvents. All fresh Gaussian
coordinates, the energy $E$ and the port parameter $z$ are held fixed.

\emph{Step 2:} Now we seek derive an an exact formula for the covariance. To differentiate a
covariance we first need a formula for it. We use the Gaussian
interpolation from the proof of Lemma~\ref{lem:trace-matrix}. In real
coordinates $x$ of the fresh layer, we put
$x_v=e^{-v}x+\sqrt{1-e^{-2v}}\,x'$, where $x'$ is an independent copy.
This is a second sample of the layer that is correlated with the first
by the factor $e^{-v}$. We write $\mathcal G_0$ and $\mathcal G_v$ for the
resolvents at $x$ and $x_v$, with the same old $X$ and $S$. The
interpolation identity gives
\begin{equation}\label{eq:dc-covariance-formula}
 \operatorname{Cov}(g_i,\mathcal G)
 =W^{-2}\int_0^\infty e^{-v}
   \E\bigl[\mathcal G_v\,\mathcal V(\mathcal G_0P_i\mathcal G_0)\,
                            \mathcal G_v\bigm| X,S\bigr]\,dv.
\end{equation}
The gradient of the trace $g_i$ produces the inner factor
$\mathcal G_0P_i\mathcal G_0$. The gradient of the matrix $\mathcal G$
produces the two outer factors $\mathcal G_v$. The factor $W^{-2}$ comes
from the entry variances, which are of order $W^{-1}$, and from the
normalization of the trace. The map $\mathcal V$ records which entries of the layer are correlated:
\[
 \mathcal V(T)=c\sum_{a=1}^2P_aTP_a
       +b\sum_{a=1}^2(P_aTQ_a+Q_aTP_a).
\]
Here $P_a$ projects onto the $a$-th new physical block and $Q_a$ onto its
port. The first sum comes from the GUE blocks, with variance $c$, and the
second from the bonds between a block and its port, with variance $b$.
For a common port we set $Q_1=Q_2=Q$. Only these entries are random. The
real diagonal of each GUE block is included through the entry covariance
$\E(d\mathcal H_{ij}d\mathcal H_{kl})=W^{-1}\sigma_{ij}\delta_{il}\delta_{jk}$.

\emph{Step 3:} We now
differentiate \eqref{eq:dc-covariance-formula} in $h$. The integrand is a
product of four resolvents: two inside $\mathcal V$ and two outside it.
By the product rule and \eqref{eq:dc-resolvent-derivative}, the
derivative is minus $W^{-2}\int_0^\infty e^{-v}\E[\,\cdot\mid X,S]\,dv$
applied to the sum of the following four terms:
\begin{align*}
 \mathcal G_v\widehat S\mathcal G_v
         \mathcal V(\mathcal G_0P_i\mathcal G_0)\mathcal G_v,\qquad 
 & \quad \mathcal G_v\mathcal V(
         \mathcal G_0\widehat S\mathcal G_0P_i\mathcal G_0)
                                                        \mathcal G_v,\\
 \mathcal G_v\mathcal V(
         \mathcal G_0P_i\mathcal G_0\widehat S\mathcal G_0)
                                                        \mathcal G_v,\quad
 & \quad\mathcal G_v\mathcal V(\mathcal G_0P_i\mathcal G_0)
                    \mathcal G_v\widehat S\mathcal G_v.
\end{align*}
Each term contains exactly one source insertion, at one of the four
resolvents, as in Figure~\ref{fig:dc-insertions}. Each insertion stays
at the sample of the resolvent it came from: an outer factor carries
$\mathcal G_v\widehat S\mathcal G_v$, and an inner factor carries
$\mathcal G_0\widehat S\mathcal G_0$. We never move a source from one
sample to the other, since the two samples are different matrices.

\begin{figure}[H]
\centering
\begin{tikzpicture}[x=1cm,y=1cm,
  box/.style={draw,thick,minimum width=1.05cm,minimum height=.75cm,font=\small}]
 \node[box] (a) at (0,0) {$\mathcal G_v$};
 \node[font=\small] at (1.05,0) {$\mathcal V\Bigl($};
 \node[box] (b) at (2.1,0) {$\mathcal G_0$};
 \node[font=\small] at (3.05,0) {$P_i$};
 \node[box] (c) at (4.0,0) {$\mathcal G_0$};
 \node[font=\small] at (4.95,0) {$\Bigr)$};
 \node[box] (d) at (5.9,0) {$\mathcal G_v$};
 \foreach \n/\k in {a/1,b/2,c/3,d/4}{
   \draw[->,thick,red!70!black] (\n.north)++(0,.75) -- (\n.north)
     node[pos=0,above,font=\small,red!70!black] {term \k};
 }
 \draw[decorate,decoration={brace,amplitude=4pt,mirror}] (1.55,-.55)--(4.55,-.55)
   node[midway,below=4pt,font=\small] {inner factors, sample $x$};
 \node[font=\small,below] at (0,-.6) {sample $x_v$};
 \node[font=\small,below] at (5.9,-.6) {sample $x_v$};
\end{tikzpicture}
\caption{The four terms in Step~3 of the proof of
Lemma~\ref{lem:t26-differentiated-covariance}. The integrand of the
covariance formula \eqref{eq:dc-covariance-formula} is a product of four
resolvents. Differentiating in the source strength inserts
$\widehat S$ at exactly one of them, replacing $\mathcal G$ by
$-\mathcal G\widehat S\mathcal G$ at the same sample. The outer factors
belong to the correlated sample $x_v$ and the inner factors to $x$.
Each of the four terms is bounded in Step~4 by the same argument, so
the count of terms, and not their position, is what matters for the
constant.}
\label{fig:dc-insertions}
\end{figure}
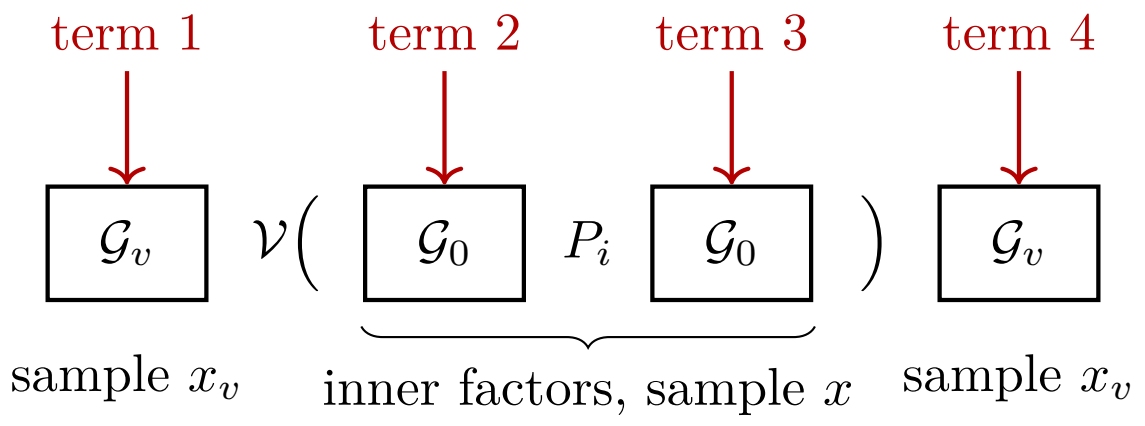

\emph{Step 4:} Now we show that each source insertion costs only $\|Y\|$. A single
insertion $\mathcal G_a\widehat S\mathcal G_a$ may look dangerous,
because $\|S\|$ is not controlled. The point is that $S$ never appears
alone. Since $\widehat S$ lives in the old-input block, which we label
$1$, and since $S=K^*YK$ by the definition of $Y$, we have the exact
factorization
\begin{equation}\label{eq:dc-factorization}
 \mathcal G_a\widehat S\mathcal G_a
   =\bigl((\mathcal G_a)_{\cdot1}K^*\bigr)\,Y\,
    \bigl(K(\mathcal G_a)_{1\cdot}\bigr).
\end{equation}
So the large matrix $K$ is always multiplied by a block of the
resolvent. These weighted factors are well behaved: by
Lemma~\ref{lem:weighted-layer}, $K\mathcal G_{11}$ and $\mathcal G_{11}K^*$
have all fixed moments bounded, uniformly in the old path. The port
blocks of $\mathcal G$ are obtained from the top block by the block
inverse equations, which multiply it by $-z^{-1}B$ or by its adjoint,
where $B$ is a fresh Gaussian bond. Hence the port factors also have
bounded moments.

Every term of Step~3 is therefore a product of factors with bounded
moments and one matrix $Y$. By the conditional H\"older inequality,
which applies to the two correlated samples as well, each term has norm
at most $C\|Y\|$. Including the prefactor $W^{-2}$ and using
$\|Y\|\le\Tr Y=Wu$, which holds because $Y\succeq0$, we obtain
\[
 \|\partial_h\operatorname{Cov}(g_i,\mathcal G)\|
 \le CW^{-2}\|Y\|\le\frac{Cu}W.
\]
This is \eqref{eq:t26-covariance-derivative}. Note that the full factor
$W^{-2}$ of the covariance has survived the differentiation.

\emph{Step 5:} Next we upper bound the fluctuation of the derivative. We now turn to the
second bound. By \eqref{eq:dc-resolvent-derivative}, the derivative of
the trace is
\[
 \dot g_i=-W^{-1}\Tr P_i\mathcal G\widehat S\mathcal G.
\]
To bound its fluctuation in the fresh layer, we use the Gaussian
Poincar\'e inequality, which controls a function by its gradient. Its
differential in the fresh Gaussian entries is
\[
 d\dot g_i=W^{-1}\Tr\big[
  (\mathcal G\widehat S\mathcal G P_i\mathcal G
   +\mathcal G P_i\mathcal G\widehat S\mathcal G)
  \,d\mathcal H\big].
\]
Each of the two terms contains one source insertion. We factor it as in
\eqref{eq:dc-factorization} and use $\|AYB\|_{\rm HS}\le\|A\|\,\|Y\|_{\rm HS}\,\|B\|$.
This bounds the Hilbert--Schmidt norm of the gradient, in every fixed
conditional moment, by $C_pW^{-1}\|Y\|_{\rm HS}$. The entries of the
fresh layer have variance of order $W^{-1}$, which contributes a further
factor $W^{-1/2}$. Hence
\[
 \|\dot g_i-\E(\dot g_i\mid X,S)\|_{p\mid X,S}
 \le C_pW^{-3/2}\|Y\|_{\rm HS}
 \le\frac{C_pu}{\sqrt W},
\]
where the last step uses $\|Y\|_{\rm HS}\le\Tr Y=Wu$, again because
$Y\succeq0$. This is \eqref{eq:t26-derivative-noise}.

\emph{Step 6:} No we turn to final bound. All
differentiations above are first carried out with cutoffs on the
Gaussian coordinates and on the inverses, so that every expression is
smooth and bounded. The uniform moments of the weighted factors then
let us remove the cutoffs by dominated convergence. For fixed $X$ and
$S$ and $h$ in a compact interval, the input $X(h)$ stays Hermitian, so
the same inverse moments hold along the whole interval.

Finally, suppose that a bounded deterministic matrix is added to
$\widehat S$. Its operator norm is bounded and its Hilbert--Schmidt norm
is $O(\sqrt W)$. Repeating Steps~4 and~5 with these norms in place of
$\|Y\|$ and $\|Y\|_{\rm HS}$ gives contributions $O(W^{-2})$ and
$O(W^{-1})$, which are smaller than $C/W$ and $C/\sqrt W$. So the bounds
hold with $1+u$ in place of $u$.
\end{proof}

 \begin{lemma}\label{lem:tagged-trace}
Let $\xi_g=(2W)^{-1}\Tr G_g$ in \eqref{eq:simultaneous-splitting}.
For $g\ne e$ put
\[
 R_e^0=(K_e-s\mu)^{-1},\quad
 \epsilon_{e,i}=W^{-1}\Tr(R_e^0)_{ii}-\mu,\quad
 B_{eg}=\partial_{\xi_g}K_e,\quad
 \lambda_{eg}=W^{-1}\Tr(R_e^0)^*B_{eg}R_e^0.
\]
Then, for every fixed $p$ and deterministic $E\in I$,
\begin{align}
 \max_i\|\epsilon_{e,i}\|_{L^p}&\le C_p/W,&
 \|\lambda_{eg}\|_{L^p}&\le C_p,\label{eq:tagged-base}\\
 \max_i\|\partial_{\xi_g}\epsilon_{e,i}\|_{L^p}
 &\le C_p(\rho^{d(e,g)}+W^{-1/2}),\qquad \rho<1 .
 \label{eq:tagged-decay}
\end{align}
The derivative in $\xi_e$ is zero. The full energy derivative obeys
\begin{equation}\label{eq:full-energy-trace}
 \sup_{E\in I}\max_i\|\epsilon'_{e,i}(E)\|_{L^p}
 \le C_p n/\sqrt W.
\end{equation}
\end{lemma}

The spatial derivative decays from the tagged location until it reaches its fluctuation scale. This is the estimate needed to sum Gaussian divergence cross terms.

\begin{proof}
At an untagged separate-port step let $d_i$ be the derivative of
$W^{-1}\Tr R_{ii}$. Since $E,z,a$ are then fixed, differentiation
of \eqref{eq:two-port-dyson} at the matched point gives exactly
\[
 v=d+(c+b^2a^2)Qv,\qquad w=ba^2v.
\]
There is no constant term. Its coefficient is the strict contraction
\eqref{eq:strict-contraction}. With actual scalar error $\delta$,
the coefficient changes by $O(\delta)$.
Write the incoming source as $S=K^*YK$.
Lemma~\ref{lem:t26-differentiated-covariance} expands the complete
covariance derivative, including the two interpolation samples. To see
how it enters the mean equation, retain its exact form
\[
 (D-\mathcal S(\bar g))\E\mathcal G
       =I+\sum_j\sigma_j\operatorname{Cov}(g_j,\mathcal G),
\]
with the appropriate block projections on each term. Here $D$ is the
augmented deterministic layer matrix, $\bar g$ are its actual scalar
means, and $\mathcal S$ is the finite scalar self-energy map.
Differentiate this identity before estimating. The new source is
$\dot D=\widehat S$, the scalar feedback is
$\mathcal S(\dot{\bar g})$, and the derivative of its covariance
right side is exactly the four-term formula above.
If $\bar M=(D-\mathcal S(\bar g))^{-1}$, then
\[
 \dot{\bar M}=-\bar M
        (\widehat S-\mathcal S(\dot{\bar g}))\bar M.
\]
Differentiating the bounded multipliers of the covariance residual
therefore adds at most
$CW^{-2}(\|Y\|+\max_j|\dot{\bar g}_j|)$.
The fixed finite scalar stability gap absorbs the last term.
Together with \eqref{eq:t26-covariance-derivative}, the differentiated
 residual is at most $Cu/W$. This is a derivative of the exact mean
identity, not of its $O(W^{-2})$ estimate.

Changing $M$ to $R$ costs $C\delta u$ in normalized traces,
because $KM-I=O(\delta)$ and positivity permits trace pairing
with $Y$. The Gaussian gradient, now with one $Y$ insertion,
has conditional fluctuation at most
$C_pW^{-1}\sqrt{W^{-1}\Tr Y^2}\le C_pu/\sqrt W$.
Consequently
\[
 |\E(d_+\mid X,S)|\le q_*|d|+C\delta u+Cu/W,\qquad
 \|d_+-\E(d_+\mid X,S)\|_{p\mid X,S}\le C_pu/\sqrt W.
\]
The conditioning used in this recursion is the sigma field generated
by the already exposed \emph{physical} blocks and bonds of $H^{(e)}$.
The matrices $X,S,Y$ and their derivatives are measurable there.
The next physical layer, including its endpoint bonds, is independent
of this field; its source projection is a deterministic choice of the
tagged block, and is never a conditioning variable. Although the
auxiliary scalar $\xi_g$ occurs in two blocks, pointwise
$\partial_{\xi_g}H^{(e)}=P_g$ for $g\ne e$. Its derivative is
therefore the same deterministic directional derivative under the
marginal independent-entry law of $H^{(e)}$. We do not condition on
$\xi_g$ or reveal one part of its auxiliary decomposition while
using freshness of the other. Exposing the original physical entries
of $H^{(e)}$ is what makes this filtration valid.

 At a tagged step there is instead a bounded source term.
Theorem~\ref{thm:matched-traces}, Lemma~\ref{lem:clock-moments}
and H\"older give $\|\delta u\|_{L^p}\le C_p/W$.
Iterating from the last tagged layer gives
$C_p(q_*^\ell+W^{-1/2})$ after $\ell$ untagged layers.

For the energy derivative let $s_r=W^{-1}\Tr R_r^*X_r'R_r$.
The exact outgoing derivative is
\[
 f'_{+,i}=W^{-1}\Tr[
 \mathcal G_{21}(I+X')\mathcal G_{12}
                         -z'\mathcal G_{22}^{\,2}]_{ii}.
\]
The bounded $I,z'$ insertions add constants to the preceding
estimates. At the matched solution their constant forcing cancels
after subtracting $a'$, since
\[
 a'=\frac{a^2}{1-a^2},\quad
 z'=(b-a^{-2})a',\quad 1+(1-b)a'+z'=0.
\]
Thus the same contraction gives outgoing error
$q_*|f'-a'|+C\delta(1+s_r)+C(1+s_r)/W$
and noise $C_p(1+s_r)/\sqrt W$.
The full clock bound proves
$\|f'_r-a'\|_{L^p}\le C_p(r+1)/\sqrt W$.
The finite-core derivatives follow by differentiating their stable
finite mean equations. Exits are removed exactly as in
Lemma~\ref{lem:clock-moments}.

The final residual pair has a fixed stability gap, as in
Theorem~\ref{thm:matched-traces}; it needs a bounded coefficient,
not a further contraction. Its source is precisely
\[
 B_{eg}=(P_g)_{SS}
 +H_{SR}(H_{RR}-E)^{-1}(P_g)_{RR}
                   (H_{RR}-E)^{-1}H_{RS}\succeq0 .
\]
A source inside the selected pair is a bounded direct injection.
Otherwise it enters the two-ended path at most twice, and at
least $d(e,g)-2$ untagged layers remain after its last injection.
This proves \eqref{eq:tagged-decay}; its clock is the final trace
pairing with $Y$ and therefore has bounded moments.
The base residual is Theorem~\ref{thm:matched-traces}.
The final energy derivative gives \eqref{eq:full-energy-trace}.
\end{proof}

\subsection{A smooth center through real Schur poles}\label{subsec:s7-a-smooth-center-through-real-schur-poles}

  Real Schur poles make the eliminated matrix singular, but the
absorbing resolvent still gives a usable center.
Lemma~\ref{lem:smooth-center} constructs a smooth version with controlled
Gaussian derivatives and negligible replacement error at the selected
rank. We use the full compressed inverse, continue its regular root,
and introduce cutoffs whose errors are paid under the pair law.

\begin{lemma}\label{lem:smooth-center}
There is a globally smooth real function $j_e^{\rm sm}(E)$ of the
auxiliary scalar coordinates, supported in a fixed bulk interval,
which agrees with \eqref{eq:dyson-center} on a regular domain and obeys the following: for each fixed $p\ge1$ there are constants $C_p,W_0$, with the other fixed parameters as in the standing width convention, such that  for $W>W_0$, we have 
\begin{equation}\label{eq:smooth-center}
 \|j_e^{\rm sm}(E)\|_{L^p}\le C_p/W,\quad
 \|\partial_{\xi_g}j_e^{\rm sm}(E)\|_{L^p}
 \le C_p(\rho^{d(e,g)}+W^{-1/2}),\quad
 \partial_{\xi_e}j_e^{\rm sm}=0 .
\end{equation}
At $\lambda_\alpha$, bounded pair-law observables change by
$O_D(W^{-D})$ when the regular center is replaced by this function.
\end{lemma}

The auxiliary center is smooth in the Gaussian coordinates even through an exterior pole. Its replacement error is paid at the original selected rank.

\begin{proof}
Use the full compressed inverse
\[
 \mathcal R_e(\zeta)=
 P_{S_e}(H^{(e)}-E-\zeta P_{S_e})^{-1}P_{S_e}\big|_{S_e},
 \qquad\Im\zeta>0.
\]
It equals $(K_e-\zeta)^{-1}$ away from Schur poles.
Almost surely every physical exterior bond is invertible.
A null vector of the full matrix must vanish on both selected
blocks by its imaginary quadratic form. Their equations force the
neighboring path endpoints to vanish, and successive invertible
bonds force the entire vector to vanish. Hence this full matrix is
invertible for every real $E$ and every real scalar variation.
Its compression is smooth and has norm at most $(\Im\zeta)^{-1}$.

At $\zeta_0=s\mu$, let $\epsilon_e^2=\sum_i|\epsilon_{e,i}|^2$.
The residual certifies stability. With $R_0=\mathcal R_e(\zeta_0)$,
$y_0=\Im\zeta_0$ and   $\Delta_0=1-s\Naux ^{-1}\Tr R_0^2$,

\[
 s\Naux ^{-1}\Tr|R_0|^2=1+O(\epsilon_e),\qquad
 \Re\Delta_0
 =1-s\Naux ^{-1}\Tr|R_0|^2+
       2sy_0^2\Naux ^{-1}\Tr|R_0|^4
 \ge2y_0^2/s-C\epsilon_e.
\]

Jensen's inequality was used in the last step.
For sufficiently small fixed $\epsilon_e<r_0$ the implicit-function
theorem gives the unique nearby Dyson root, with fixed positive
height and displacement $O(\epsilon_e)$.
At that root $\Re\Delta\ge2(\Im\zeta)^2/s$.

Choose a smooth cutoff $\chi(\epsilon_e^2)$ equal to one for
$\epsilon_e\le r_0/4$ and zero for $\epsilon_e\ge r_0/2$.
Set
  $j_e^{\rm sm}=\chi\,s\Naux ^{-1}\Tr J_e|\mathcal R_e(\zeta_e)|^2$
 there and zero elsewhere, with a fixed smooth energy cutoff when
needed. This is globally smooth, bounded by one, and independent
of its own $\xi_e$ at fixed $E$.
At the root, with $y=\Im\zeta_e$,
  $j_e=(s/y)\Im[\Naux ^{-1}\Tr J_e\mathcal R_e(\zeta_e)]$.
The contrast at $\zeta_0$ is a residual difference, hence
$|j_e|\le C\epsilon_e$. This proves the first bound in
\eqref{eq:smooth-center}.

For the derivative, implicit differentiation gives

\[
 \partial_{\xi_g}\zeta_e
 =-\frac{s\Naux ^{-1}\Tr[\mathcal R_e(\zeta_e)
                      B_{eg}\mathcal R_e(\zeta_e)]}{\Delta}.
\]

Resolvent identities on the fixed-height disk give, for
$J=I,J_e,P_i,P_{i+1}$,

\[
 \left|\Naux ^{-1}\Tr J[
 \mathcal R_e(\zeta)B_{eg}\mathcal R_e(\zeta)-R_0B_{eg}R_0]\right|
 \le C|\zeta-\zeta_0|\lambda_{eg}.
\]

Indeed write $\mathcal R_e(\zeta)=R_0[I-(\zeta-\zeta_0)R_0]^{-1}$;
the extra multiplier differs from $I$ by $O(|\zeta-\zeta_0|)$.
Positivity of $B_{eg}$ and Hilbert--Schmidt Cauchy--Schwarz give
the displayed bound without $\norm{B_{eg}}$.
Differentiating the contrast formula and the cutoff therefore yields
\[
 |\partial_{\xi_g}j_e^{\rm sm}|
 \le C\left(\sum_i|\partial_{\xi_g}\epsilon_{e,i}|
                         +\epsilon_e\lambda_{eg}\right).
\]
Lemma~\ref{lem:tagged-trace} and H\"older prove
\eqref{eq:smooth-center}, off the real Schur poles of $K_e$, where
$\epsilon_{e,i}$, $B_{eg}$ and $\lambda_{eg}$ are defined. At a fixed energy
these poles occur only on a null set of the Gaussian variables, so the
$L^p$ norms in \eqref{eq:smooth-center} are unaffected.

Finally Lemma~\ref{lem:random-energy},
\eqref{eq:rank-location} and \eqref{eq:full-energy-trace} give
$\|\epsilon_e(\lambda_\alpha)\|_{L^q}\le C_q/W$ for every fixed
$q$. It suffices to prove this for $q\ge1$; smaller positive moments follow by monotonicity of $L^q$ norms on the probability space. Use a slightly larger bulk interval and pay clipping by the
rank-location tails. Thus certificate failure has probability
$O_q(W^{-q})$. Multiplication by the pair density, at most
$n\le TW$, still makes its bounded contribution smaller than every
fixed power. We now specify the regular domain of
Theorem~\ref{thm:physical-center}. Let $\epsilon_*$ be the constant of
Lemma~\ref{lem:regularity-certificate} for $c_0=s$ and the bulk interval
$I$. Fix $\epsilon_{\rm reg}\le\min(r_0/4,\epsilon_*)$, with $r_0$ the
radius of the cutoff $\chi$ above (not the window radius of that lemma),
so small that $C\epsilon_{\rm reg}/s$, with $C$ the displacement constant
of the Dyson root above, is less than the radius of the Rouch\'e disk in
the proof of Lemma~\ref{lem:regularity-certificate}. Put
$D_e(E)=\{E\in I,\ \epsilon_e(E)\le\epsilon_{\rm reg}\}$; it depends only
on $E$ and $H^{(e)}$. Since $H^{(e)}_{S_eS_e}$ is the matrix $V$ of
\eqref{eq:pair-splitting} for $G=G_e$, $K_e(E)=F(E)-E$, and $\mu=\bar a$
is the point $\mu$ of that lemma. Hence, for $F=K_e(E)+E$ and $c_0=s$,
its residual is
$\av{(K_e-s\mu)^{-1}}-\mu=(\epsilon_{e,i}+\epsilon_{e,i+1})/2$, of modulus
at most $\epsilon_e$, and on $D_e(E)$ it gives \eqref{eq:fresh-regularity}
with $c_0=s$ and fixed constants, without any norm bound. The root
satisfies $|\zeta_e/s-\mu|\le C\epsilon_e/s$, so $\zeta_e/s$ is the unique
zero in the Rouch\'e disk at $z=E$; this zero is continuous in $z$ and is
the Dyson solution for $\Im z>0$, hence equals $m_s(E)$, and
$\zeta_e=sm_s(E)$ has fixed positive height and bounded modulus. So
$s|K_e-\zeta_e|^{-2}$ is the covariance both of
Theorem~\ref{thm:fresh-palm} and of \eqref{eq:angular-reference}. On
$D_e(E)$ also $\chi(\epsilon_e^2)=1$, so, with the energy cutoff equal to
one on $I$, $j_e^{\rm sm}=j_e$ there. Finally, the moment bound above,
Markov's inequality, the rank-location tails of
Theorem~\ref{thm:rank-location} for the exit of $\lambda_\alpha$ from $I$,
and $d\Prob^\sharp_{\alpha,e}\le(n/2)\,d\Prob\le TW\,d\Prob$ give
$\Prob^\sharp_{\alpha,e}(D_e(\lambda_\alpha)^c)\le C_qW^{1-q}$. The pole
and norm restrictions of the angular comparison are paid separately in
Step~3 of the proof of Theorem~\ref{thm:physical-angular}.
\end{proof}

\subsection{The Gaussian energy primitive}\label{subsec:s7-the-gaussian-energy-primitive}

 We now finish the estimate for the physical center in
Theorem~\ref{thm:physical-center}. Integrating its square from a fixed
energy to the selected eigenvalue converts the pair-mass weight into an
exact Gaussian derivative. We apply integration by parts, bound the
resulting divergence and its cross terms, and use the random-energy
integral estimate to conclude.

 \begin{proof}[Proof of Theorem~\ref{thm:physical-center}]
Decompose $G_e=G_e^0+\xi_eI_{S_e}$, where
  $\sigma^2=\Var\xi_e=s/\Naux ^2$.
The scalars $\xi_e$ are independent of one another and of all
traceless parts and $V$. Let $\vartheta(E)$ be a smooth bulk
cutoff equal to one near every allowed $e_\alpha$, and define
\[
 f_e(E)=\vartheta(E)|j_e^{\rm sm}(E)|^2,\qquad
 D(E)=\sigma^{-2}\sum_e\xi_ef_e(E).
\]
Each $f_e$ is independent of its own scalar. If
$\Lambda=\lambda_\alpha(H)$, Hellmann--Feynman gives almost everywhere
$\partial_{\xi_e}\Lambda=p_i+p_{i+1}=q_e$.
The bounded energy primitive
$F_e=\int_{e_\alpha}^{\Lambda}f_e(E)\,dE$ therefore satisfies
$\partial_{\xi_e}F_e=q_ef_e(\Lambda)$.
Gaussian integration by parts and cyclic invariance give the exact
identity
\begin{equation}\label{eq:energy-primitive-identity}
 \E^\sharp_{\alpha,e}f_e(\Lambda)
 =\frac{1}{2}\E\int_{e_\alpha}^{\Lambda}D(E)\,dE.
\end{equation}
The factor $1/2$ follows from
$\sum_e\E[q_ef_e(\Lambda)]=2\E^\sharp_{\alpha,e}f_e(\Lambda)$.
The original ordered eigenvalue is used throughout, so no derivative
of an inertia indicator appears.

At each deterministic energy, two integrations by parts give
\begin{equation}\label{eq:divergence-isometry}
 \E D(E)^2
 =\sigma^{-2}\sum_e\E f_e^2
   +\sum_{e\ne g}\E[
      (\partial_{\xi_g}f_e)(\partial_{\xi_e}f_g)].
\end{equation}
In a cross term the first integration differentiates $f_g$ in
$\xi_e$; the second differentiates $f_e$ in $\xi_g$.
The possible second derivative vanishes since $f_g$ is independent
of $\xi_g$. Approximation by functions smooth and truncated only
in the other scalar coordinates preserves this independence.
The first-derivative moment bounds justify passage to the limit
in the Gaussian Sobolev norm, and hence justify
\eqref{eq:divergence-isometry} without second-derivative assumptions.

By \eqref{eq:smooth-center} and H\"older,
$\E f_e^2\le C/W^4$ and
$\|\partial_{\xi_g}f_e\|_2
\le CW^{-1}(\rho^{d(e,g)}+W^{-1/2})$.
At most two edges occur at each positive cyclic distance. Therefore
\begin{equation}\label{eq:divergence-bound}
 \E D(E)^2\le C(n/W^2+n^2/W^3),\qquad
 \sup_E\|D(E)\|_2\le C_T\sqrt n/W.
\end{equation}
Apply the expectation version of Lemma~\ref{lem:random-energy}
to \eqref{eq:energy-primitive-identity}, using
$\omega=(n\sqrt W)^{-1}$ from Theorem~\ref{thm:rank-location}.
It bounds the nonnegative left side by
$C\omega\sqrt n/W=C/(\sqrt n\,W^{3/2})$.
The energy cutoff and the smooth-center replacement cost arbitrary
fixed negative powers under the pair law, by rank-location moments
and Lemma~\ref{lem:smooth-center}. Choose those powers first and
absorb the errors. This proves \eqref{eq:physical-center}.
\end{proof}
\section{Conditional eigenvectors of a fresh Gaussian block}
\label{sec:fresh-gaussian}
When a block is exposed as a fresh GUE matrix, we need the moments of
its eigenvector components, conditionally on everything else and at the
eigenvalue selected by the rank. This section shows that these moments are
Gaussian, with the correct covariance. It is used in the angular comparison
of Part~I and in the closing law of Part~II.
Theorem~\ref{thm:fresh-palm} states that the $r$-th moment of
$|\sqrt\Naux q^*u_\ell|^2$, weighted by the density that selects the
eigenvalue, equals $r!(q^*\Gamma q)^r$ times the sector mass, up to a relative
error $\Naux^{-\kappa_r}$ and an additive error $\Naux^{-K}$. The weight is
kept unnormalized, so rare rank sectors cause no loss from division.
In Subsection~\ref{subsec:s8-one-bounded-trace-certifies-regularity} we show that one bounded trace error certifies a
whole regular window, without bounding the norm of the deformation
(Lemma~\ref{lem:regularity-certificate}). In Subsection~\ref{subsec:s8-a-relative-local-law-at-coarse-height} we prove a
local law relative to the covariance of each direction
(Lemma~\ref{lem:coarse-relative}). In Subsection~\ref{subsec:s8-fixed-variance-hermitian-moment-flow} we analyze the
conditional eigenvector moment flow, following the localization and averaging
method of~\cite{Benigni} (Lemma~\ref{lem:t26-conditional-moments}). In
Subsection~\ref{subsec:s8-the-exact-scalar-palm-weight} we insert the exact scalar weight that selects the
eigenvalue and complete the proof.

Let   $F=F^*\in\C^{\Naux \times \Naux }$, $\norm F\le \Naux ^A$, and let
$0<c_0<1$ be fixed. Write
$M_s(z)=(F-z-sm_s(z))^{-1}$ with $m_s=\av{M_s}$.
Assume that for fixed $r_0,s_0,v_0>0$, $s_0<c_0/2$,
\begin{equation}\label{eq:fresh-regularity}
 v_0\le\Im m_s(x+i\eta)\le v_0^{-1},\quad
 s\in[c_0-s_0,c_0],\quad |x-E|\le r_0,\quad 0\le\eta\le1 .
\end{equation}
Put $\Gamma=c_0|F-E-c_0m_{c_0}(E)|^{-2}$, so   $\Tr\Gamma=\Naux $.
Split $A_\Naux \sim\GUE_\Naux (c_0/\Naux )$ as $A_\Naux =A_\Naux ^0+gI$,
where $g\sim N(0,c_0/\Naux ^2)$  is independent of the trace-free part.
Let $(\nu_\ell^0,u_\ell)$ be the ordered eigenpairs of   $F+A_\Naux ^0$  and
define

\begin{equation}\label{eq:palm-scalar-weight}
 w_{\ell,E}=\varphi_{\sqrt{c_0}/\Naux }(E-\nu_\ell^0),\qquad
 Z_{\ell,E}=\E w_{\ell,E}.
\end{equation}

Here $\varphi_\sigma$ is the centered real Gaussian density of
variance $\sigma^2$.

\begin{theorem}\label{thm:fresh-palm}
For each fixed integer $r\ge1$ there is $\kappa_r>0$ such that,
for every fixed $K>0$, every deterministic complex vector $q$ and
every ordered rank $\ell$,

\begin{equation}\label{eq:fresh-palm-moments}
 \left|\E[w_{\ell,E}|\sqrt \Naux q^*u_\ell|^{2r}]
       -r!(q^*\Gamma q)^r Z_{\ell,E}\right|
 \le C_r\{\Naux ^{-\kappa_r}Z_{\ell,E}+C_{r,K}\Naux ^{-K}\}
                         (q^*\Gamma q)^r.
\end{equation}

The constants are uniform over \eqref{eq:fresh-regularity} and
the fixed polynomial norm bound. No lower bound on $Z_{\ell,E}$
is assumed.
\end{theorem}

The measure in \eqref{eq:palm-scalar-weight} is exactly the
eigenpair density at energy $E$; it uses the original scalar
Gaussian variance. Precisely, let $(\lambda_\ell,u_\ell)=(\nu_\ell^0+g,u_\ell)$
be the ordered rank-$\ell$ eigenpair of $F+A_\Naux$, and let $f\ge0$ be
measurable and phase invariant on the unit sphere. Then, for every $E$,
\[
 \E[w_{\ell,E}f(u_\ell)]
 =\int c_{\Naux,t}\mathcal R(E,v,F)\,
   \E_{G_{\Naux-1}}\bigl[\det(M)^2\one_{\{N_-(M)=\ell-1\}}\bigr]
   f(v)\,d\nu_\Naux(v),\qquad t=c_0/\Naux,
\]
with the notation of Theorem~\ref{thm:coarea} for the constant
deformation $F$. Indeed, as $g$ is independent of $(\nu^0,u)$, both sides integrated over
$E$ in a Borel set $B$ equal $\E[\one_{\{\lambda_\ell\in B\}}f(u_\ell)]$
(the right side by Claim~2 in the proof of Theorem~\ref{thm:coarea} and
Lemma~\ref{lem:eigenpair-jacobian}), so they agree for almost every $E$.
For bounded $f$ both are continuous in $E$, by dominated convergence:
$\varphi_{\sqrt{c_0}/\Naux}$ is bounded and continuous, $\mathcal R\le1$
is continuous, $\det(M)^2$ is dominated by a polynomial in
$\norm{G_{\Naux-1}}$ locally uniformly in $E$, and the indicator is
continuous at $E$ since $E\notin\operatorname{Spec}(F_v+G_{\Naux-1})$
almost surely. Hence they
agree for every $E$, and monotone convergence removes the boundedness of
$f$. The proof of Theorem~\ref{thm:physical-angular} needs this pointwise
form, since there $F=F(E)$. The order $r$ is fixed before   $\Naux \to\infty$.
We do not use a growing-order moment-flow assertion.

\subsection{One bounded trace certifies regularity}\label{subsec:s8-one-bounded-trace-certifies-regularity}

 We need a simple way to check that a deformed Gaussian block has a
regular bulk window. Lemma~\ref{lem:regularity-certificate} obtains a
whole fixed window and a covariance comparison from one bounded trace
error, without bounding the deformation norm. We prove stability of
the scalar equation and continue its solution around the reference
point.

 \begin{lemma}\label{lem:regularity-certificate}
Fix $0<c_0<1$ and $E$ in a compact semicircle bulk interval.
With $\mu=(-E+i\sqrt{4-E^2})/2$, put
$\chi=\av{(F-E-c_0\mu)^{-1}}-\mu$.
There are fixed $\epsilon_*,r_0,s_0,v_0>0$ such that
$|\chi|\le\epsilon_*$ implies \eqref{eq:fresh-regularity}, and
\begin{equation}\label{eq:certificate}
 |m_{c_0}(E)-\mu|\le C|\chi|,\qquad
 (1-C|\chi|)\Gamma_0\preceq\Gamma\preceq(1+C|\chi|)\Gamma_0,
 \quad \Gamma_0=c_0|F-E-c_0\mu|^{-2}.
\end{equation}
This lemma has no norm assumption on $F$.
\end{lemma}

A single trace value at a bounded imaginary shift gives a fixed regularity rectangle. The absence of a deformation norm bound lets this certificate survive nearby Schur poles.

\begin{proof}
Let $R=(F-E-c_0\mu)^{-1}$ and $v=\Im\mu$.
The imaginary resolvent identity and Jensen give
\[
 \av{|R|^2}=\frac{v+\Im\chi}{c_0v},\qquad
 \Re(1-c_0\av{R^2})
 =-\Im\chi/v+2c_0\av{(\Im R)^2}
 \ge-|\chi|/v+2c_0(v+\Im\chi)^2 .
\]
The last expression is bounded below uniformly for small
$\epsilon_*$. All derivatives of
$\mathcal F(s,z,m)=m-\av{(F-z-sm)^{-1}}$
are bounded on a fixed neighborhood with
$s\ge c_0/2$ and $\Im(z+sm)\ge c_0v/4$.
Choose a small disk around $\mu$ so that its derivative varies
by less than one quarter of this lower bound. If
$|\chi|,|s-c_0|,|z-E|$ are sufficiently small, the linear term
dominates the residual and Taylor remainder on its boundary.
Rouch\'e's theorem gives a unique nearby zero and displacement
$C(|\chi|+|s-c_0|+|z-E|)$.
For $\Im z>0$ it is the unique Dyson solution, since
$s\av{|M_s|^2}=s\Im m_s/(\Im z+s\Im m_s)<1$
and subtraction of two equations followed by Cauchy--Schwarz
proves uniqueness.

The resulting density has a fixed positive lower bound on
$|x-E|\le2r_0$. Its Poisson representation extends this to
$|x-E|\le r_0$, $0\le\eta\le1$; explicitly the lower bound is
at least $(v/\pi)\arctan r_0$ after reducing $v$ by a fixed
factor. The imaginary equation also gives $|m_s|^2\le1/s$,
which provides the upper bound.
Finally diagonalize $F$. Changing $\mu$ to $m_{c_0}(E)$ changes
each denominator by $O(|\chi|)$, while its modulus is bounded
below by $c_0v$. The squared reciprocal moduli therefore differ
relatively by $O(|\chi|)$, proving \eqref{eq:certificate}.
\end{proof}

\subsection{A relative local law at coarse height}\label{subsec:s8-a-relative-local-law-at-coarse-height}

  A small absolute error can be useless in a direction with very small
variance. Lemma~\ref{lem:coarse-relative} instead measures the resolvent
error relative to that direction's deterministic covariance. We first
control the normalized trace by its stable mean equation, then use
Gaussian gradient bounds to obtain the directional estimate.

\begin{lemma}\label{lem:coarse-relative}
Under \eqref{eq:fresh-regularity}, for a polynomial-size deterministic
family of vectors $q$, with probability at least   $1-\Naux ^{-D}$  for
every fixed $D$, simultaneously on
$|x-E|\le r_0/2$,   $\Naux ^{-1/100}\le\eta\le1$,

\begin{equation}\label{eq:coarse-relative}
 |\av{G_s(z)}-m_s(z)|\le \Naux ^{-1/3},\qquad
 |q^*(G_s(z)-M_s(z))q|
       \le \Naux ^{-1/3}q^*\Im M_s(z)q,
\end{equation}

where   $G_s(z)=(F+\sqrt s\,G_\Naux -z)^{-1}$.
\end{lemma}

The error is measured relative to the deterministic covariance in each direction. This retains information even when that covariance is polynomially small.

\begin{proof}
Let   $m_\Naux =\av{G_s}$ and $\bar m=\E m_\Naux $.
Gaussian gradient bounds give
  $\|m_\Naux -\bar m\|_{L^p}\le C_p\Naux ^{-1}\eta^{-2}$.
The integration-by-parts equation, after diagonalizing $F$, yields
  $|\bar m-\av{(F-z-s\bar m)^{-1}}|\le C\Naux ^{-1}\eta^{-4}$.
The scalar stability follows from
\[
 \Re(1-s\av{M_s^2})
 =\frac{\eta}{\eta+s\Im m_s}
        +2s\av{(\Im M_s)^2}
 \ge2(c_0-s_0)v_0^2 .
\]
At height one the scalar map is a contraction with constant at most
$s<1$. Continue its approximate solution down in $\eta$ using this
stability; the error   $C\Naux ^{-1}\eta^{-4}$  remains inside the fixed
stable neighborhood. This proves the trace estimate.

For the directional estimate retain the weights in the identities

\[
 G_sM_s^{-1}=I-G_s(\sqrt s\,G_\Naux +sm_s),\qquad
 G_s^*M_s^{-*}=I-G_s^*(\sqrt s\,G_\Naux +s\bar m_s).
\]

They imply
$\|\norm{G_sq}\|_{L^p}+\|\norm{G_s^*q}\|_{L^p}
\le C_p\eta^{-1}\norm{M_sq}$.
Poincar\'e therefore gives centered error
  $C_p\Naux ^{-1/2}\eta^{-2}\norm{M_sq}^2$.
Testing the exact mean equation

\[
 \E G_s-M_s
 =sM_s\{(\bar m-m_s)\E G_s+\E[(m_\Naux -\bar m)G_s]\}
\]

on $q$ gives mean error at most
  $C\Naux ^{-1}\eta^{-5}\norm{M_sq}^2$.
Since $q^*\Im M_sq=(\eta+s\Im m_s)\norm{M_sq}^2$,
both are relative errors and are smaller than   $\Naux ^{-1/3}$  on
the stated domain after fixed-moment Markov bounds.
Finally   $q^*\Im M_sq\ge C\Naux ^{-2\max(A,1)}\norm q^2$.
A sufficiently fine polynomial mesh controls interpolation in
$z$ even in these small directions; higher fixed moments pay
its size and that of the vector family.
\end{proof}

\subsection{Fixed-variance Hermitian moment flow}\label{subsec:s8-fixed-variance-hermitian-moment-flow}

We   prove the conditional Gaussian moment comparison for the complex
matrix used here. Lemma~\ref{lem:t26-conditional-moments} gives a
quantitative error on typical eigenvalue paths. Following the
 localization and averaging method   of~\cite{Benigni}
, we derive the
complex generator, restrict its propagation, compare the local
covariance averages, and control the exceptional paths.

Fix a terminal variance $c=c_0$, put   $\tau=\Naux ^{-\theta}$  with
$\theta=1/100$, and write $H_u=F+G_{c-\tau+u}$ for
$0\le u\le\tau$. The Gaussian matrix at time zero and its subsequent
increments are independent. For nonzero $q$, divide $q$ by
$\sqrt{q^*\Gamma q}$; throughout this subsection $q^*\Gamma q=1$.
In particular   $\|q\|\le C\Naux ^{\max(A,1)}$. Set

\[
 X_i(u)=|\sqrt \Naux q^*u_i(u)|^2,\qquad
 F_r(u,\xi)=\E\left[\prod_i\frac{X_i(u)^{\xi_i}}{\xi_i!}
           \,\middle|\,\mathcal L\right],\qquad
 \mathcal L=\sigma(\lambda_i(u):i\le \Naux ,\ 0\le u\le\tau),
\]

where $|\xi|=r$ is fixed. The independent phases of the eigenvectors
are immaterial to these radial observables. The angular Brownian noises
are independent of the eigenvalue noises, by the covariance argument in
\eqref{eq:t26-vector-dbm}. The conditional moments therefore solve

\begin{equation}\label{eq:t26-complex-moment-generator}
 \partial_uF_r=\mathcal B_uF_r,\qquad
 (\mathcal B_uf)(\xi)=\sum_{i\ne j}c_{ij}(u)\xi_i(1+\xi_j)
                          [f(\xi^{i,j})-f(\xi)],\qquad
 c_{ij}=\frac{1}{\Naux (\lambda_i-\lambda_j)^2}.
\end{equation}

For clarity this normalization follows directly from the fields
$X_{ij},Y_{ij}$ of Section~\ref{sec:relative-vectors}. On the two
radial variables $x=X_i,y=X_j$, their generator is
$c_{ij}\{xy(\partial_x-\partial_y)^2+(y-x)(\partial_x-\partial_y)\}$.
Applied to $x^ay^b/(a!b!)$, the transfer coefficients are
$a(b+1)$ and $b(a+1)$. This also proves detailed balance with
\emph{counting measure} on configurations: the reverse coefficient at
$\xi^{i,j}$ is $(\xi_j+1)\xi_i$. The normalization $2^rr!$ used
for variance-two complex Gaussians in \mbox{%
\cite[(3.4)]{BY} }\hskip0pt%
differs by the
constant $2^r$ on the fixed particle sector. Our variance-one circular
normalization is $r!$.

In teh following lemma, we follow the eigenvectors of $H_u$ along
the Gaussian flow for a short time $\tau$, conditionally on the whole
eigenvalue path $\mathcal L$. For a configuration $\xi$ of $r$ particles
placed on eigenvalue labels, $F_r(u,\xi)$ is the conditional moment of
the corresponding eigenvector overlaps with the fixed direction $q$,
normalized so that the Gaussian value is one. These moments evolve by
the generator \eqref{eq:t26-complex-moment-generator}. A particle at
label $i$ jumps to label $j$ with a rate built from the coefficient
$c_{ij}$, which is large when $\lambda_i$ and $\lambda_j$ are close. The
following lemma shows that after time $\tau$ every moment in the bulk is within
  $\Naux ^{-\kappa_r}$  of its Gaussian value one, for all eigenvalue paths
except a set of tiny probability.

 However, there are severe challenges in going forward. The jump rates are random and depend on the
whole spectrum, so we cannot solve the flow. We use instead a
maximum-principle argument: at a configuration where $F_r$ is largest,
the flow pulls the maximum down towards one at rate $1/\eta$. To make
this work we must control three things. Jumps over long distances do
not see the local spectrum, so we must show that they are negligible.
The flow must not carry mass in from far away, which is a statement of
finite speed of propagation. And the pull towards one must be computed
from the local law, using the moments with one particle fewer. This is
why the proof is an induction on $r$.

\begin{lemma}[Conditional moments at fixed Gaussian variance]
\label{lem:t26-conditional-moments}
Under \eqref{eq:fresh-regularity}, for each fixed $r$ there is
$\kappa_r>0$ such that outside a measurable set of full eigenvalue
paths of probability at most   $\Naux ^{-D}$, for any prescribed fixed $D$,

\begin{equation}\label{eq:t26-conditional-moment-result}
 |F_r(\tau,r\delta_k)-1|\le C_r\Naux ^{-\kappa_r}
 \qquad (|\lambda_k(\tau)-E|\le \Naux ^{-1+\epsilon}).
\end{equation}

The same estimate holds for each fixed particle configuration whose
labels all lie in the shrinking bulk windows used in the proof.
Constants and exceptional powers are uniform in the deformation and
in the normalized direction. The good path set may depend on that
direction. The exponent $D$ is chosen after $r,A$.
\end{lemma}

\begin{proof}
All bounds below concern the conditional moments in
\eqref{eq:t26-complex-moment-generator} themselves, without any
modification. We proceed in five steps: the scales and the good events,
short-range propagation, long jumps and flattening, the removal identity,
and the maximum principle with the induction.

\emph{Step 1: the scales.} We choose $\epsilon>0$ small, depending on
$r$, and use four scales,

\begin{equation}\label{eq:t26-moment-scales}
 \eta=\Naux ^{-\epsilon}\tau,\qquad
 \ell=\lfloor \Naux ^{1-\theta/3}\rfloor,\qquad
 w=\Naux ^{1-\theta/6},\qquad R=\Naux ^\epsilon\ell.
\end{equation}

The spectral height $\eta$ is the time scale of the relaxation, and it
is a little shorter than the total time $\tau$. The rank distance $\ell$
separates short jumps from long jumps. The rank distance $w$ is the
radius of a window around the target configuration, used to flatten the
initial data. The rank distance $R$ is how far a short jump process can
travel in time $\tau$, up to a tiny error. These scales are ordered as

\[
 \Naux \sqrt\tau\ll\ell\ll w\ll \Naux ,\qquad \Naux \eta\ll \Naux \tau\ll\ell.
\]

Figure~\ref{fig:cm-scales} shows them. For the finitely many nested
windows in the induction on $r$ we may use fixed multiples of $w$.

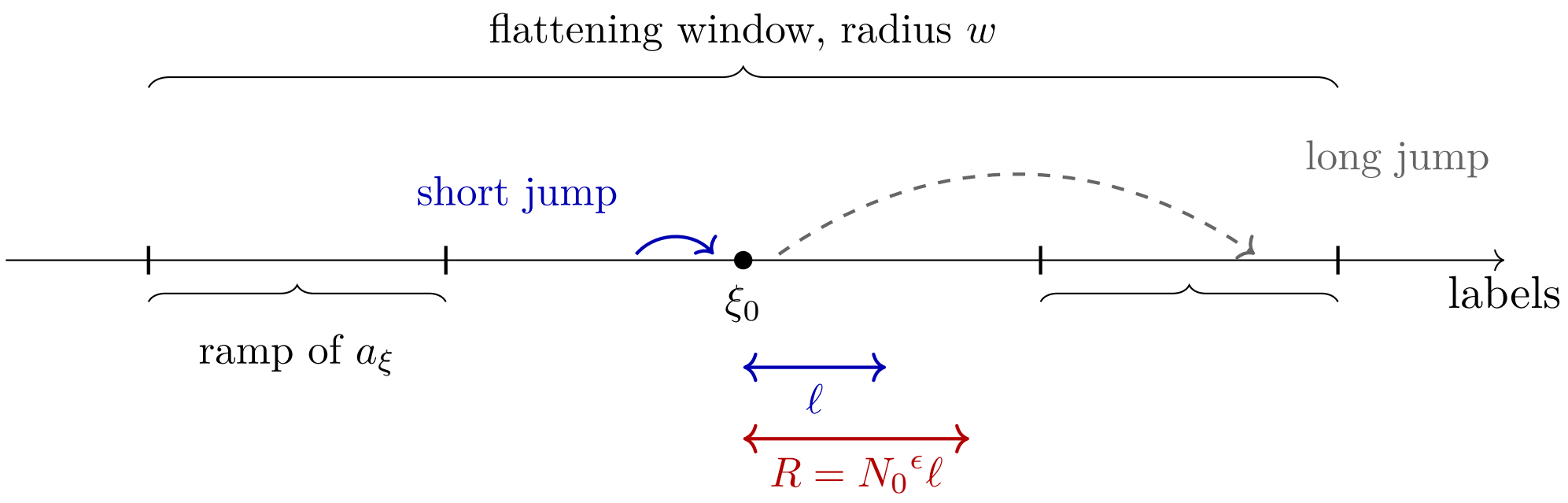
\begin{figure}[ht]
\centering
  \begin{tikzpicture}[x=1cm,y=1cm]
 \draw[->] (-6.2,0)--(6.4,0) node[below] {labels};
 \fill (0,0) circle (2.2pt) node[below=2pt,font=\small] {$\xi_0$};
 \draw[thick] (-5,-.12)--(-5,.12) (5,-.12)--(5,.12);
 \draw[thick] (-2.5,-.12)--(-2.5,.12) (2.5,-.12)--(2.5,.12);
 \draw[decorate,decoration={brace,amplitude=5pt}] (-5,1.45)--(5,1.45)
   node[midway,above=5pt,font=\small] {flattening window, radius $w$};
 \draw[decorate,decoration={brace,amplitude=4pt,mirror}] (-2.5,-.35)--(-5,-.35)
   node[midway,below=4pt,font=\small] {ramp of $a_\xi$};
 \draw[decorate,decoration={brace,amplitude=4pt}] (2.5,-.35)--(5,-.35);
 \draw[<->,thick,blue!70!black] (0,-.9)--(1.2,-.9)
   node[midway,below,font=\small,blue!70!black] {$\ell$};
 \draw[<->,thick,red!70!black] (0,-1.5)--(1.9,-1.5)
   node[midway,below,font=\small,red!70!black] {$R=\Naux ^\epsilon\ell$};
 \draw[->,thick,blue!70!black] (-.9,.05) to[bend left=50] (-.25,.05);
 \draw[->,thick,black!60,dashed] (.3,.05) to[bend left=35] (4.3,.05);
 \node[font=\small,blue!70!black] at (-1.9,.55) {short jump};
 \node[font=\small,black!60] at (5.5,.85) {long jump};
\end{tikzpicture}
 \caption{The scales in the proof of
Lemma~\ref{lem:t26-conditional-moments}, measured in eigenvalue labels. Jumps shorter than $\ell$ are kept in the short-range
flow $\mathcal S$, and in time $\tau$ that flow moves particles by at most
  $R=\Naux ^\epsilon\ell$, up to a tiny error (Step~2). Longer jumps are rare
and are treated as a perturbation (Step~3). Around the target
configuration $\xi_0$ the initial data are flattened: the weight $a_\xi$
equals one at distance at most $w/2$, decreases linearly to zero between
$w/2$ and $w$, and vanishes beyond. Since $\ell\ll w$, the short flow
cannot carry information from outside the window to $\xi_0$, and the
averaging error is measured by $\ell/w$.}
\label{fig:cm-scales}
\end{figure}

\emph{The good events.} The maximum principle needs the local law
along the whole flow. For each fixed $u$, Theorem~\ref{thm:external-local} applies to
$H_u=F+G_{c-\tau+u}$ with $s=c-\tau+u\in[c_0-s_0,c_0]$,
deterministic direction $q$ and averaged test $I$.
Fresh regularity gives $\Naux h\rho_s\ge\Naux^{4\epsilon}v_0/\pi$.
Diagonalizing $F$ and using $|m_s|\le s^{-1/2}$ gives
\[
 |f-x-ih-sm_s|\asymp1+|f-E|
       \asymp|f-E-c_0m_{c_0}(E)|.
\]
The corresponding relative covariances are comparable, so
$q^*\Gamma(z)q\le C$ when $q^*\Gamma q=1$.
For the overlap bound use $q^*\Im G_uq=q^*\Im M_sq+q^*\Im(G_u-M_s)q$,
not a bound on $|q^*M_sq|$. Choose the domination loss below
$\epsilon$ and the domain exponent below $4\epsilon$.
Polynomial meshes in $u,x,h$, resolvent Lipschitz bounds and the
Brownian modulus give a joint event $\mathcal A$, of
probability at least   $1-\Naux ^{-D}$  for any $D$, on which, for
  $\Naux ^{-1+4\epsilon}\le h\le1$  in the fixed bulk window,

\begin{align}
 |\av{G_u(x+ih)}-m_{c-\tau+u}(x+ih)|
                    &\le \Naux ^\epsilon/(\Naux h),\notag\\
 |q^*(G_u-M_{c-\tau+u})q|
                    &\le C \Naux ^\epsilon/\sqrt{\Naux h},\notag\\
 X_i(u)&\le C \Naux ^{5\epsilon}
                 \quad (|\lambda_i(u)-E|\le r_0/4).
                                                     \label{eq:t26-moment-local-law}
\end{align}

The first line controls the eigenvalue density, the second the
eigenvector overlaps in the direction $q$, and the third says that no
single overlap is large. The third follows from the second by the
positive spectral decomposition of $\Im G_u$ at height
  $\Naux ^{-1+4\epsilon}$.

We also need the deterministic approximations to vary slowly. The
regularity assumption gives $\|M_s\|\le C$ and
$\Re(1-s\av{M_s^2})\ge2s(\Im m_s)^2\ge c>0$. By implicit
differentiation, $m_s$ and the relative quadratic forms of
$sM_s^*M_s$ have bounded derivatives. Hence, in the shrinking windows,
\[
 q^*\Im M_{c-\tau+u}(x+ih)q
   =\Im m_{c-\tau+u}(x+ih)
             [1+O(|x-E|+\tau+h)].
\]
All these statements are relative to $q$. They remain useful even when
the original covariance in the direction $q$ is polynomially small.

A second event $\mathcal R$, also of probability at least   $1-\Naux ^{-D}$,
collects the spectral facts used for propagation: the global upper
count and the bulk lower count of \eqref{eq:t26-long-gaps}, and
$\sup_{u\le\tau}\|H_u-H_0\|\le C\sqrt\tau$. The Gaussian modulus of
continuity and Lipschitz bounds for the resolvent fill in the time
meshes. No moments of inverse eigenvalue gaps are used.

The moments are conditional on $\mathcal L$, while $\mathcal A$ also
involves eigenvectors. We therefore disintegrate $\mathcal A$ over
$\mathcal L$. Starting from the failure probability   $\Naux ^{-2D}$, we discard
the paths on which   $\Prob(\mathcal A^c\mid\mathcal L)>\Naux ^{-D}$. By Markov's
inequality this costs probability   $\Naux ^{-D}$. Every raw moment is bounded
deterministically by   $\Naux ^{C(r,A)}$. So on the retained paths, the local
conditional moments are   $O(\Naux ^{C_r\epsilon})$, and every use of a local
law inside a conditional moment costs an additional error
  $\Naux ^{-D+C(r,A)}$. In particular we never condition the moment equation on
$\mathcal A$.

\emph{Step 2: short-range propagation.} Let $\mathcal S$ keep only the
jumps with $|i-j|<\ell$ in \eqref{eq:t26-complex-moment-generator}, and
let $U_{s,u}$ be its propagator from time $s$ to time $u$. Because it is
the propagator of a Markov jump process, it is positive, preserves
constants, and contracts the supremum norm. We now show that it moves
mass by at most $R$ labels, up to a tiny error. We need this only after
integration in the time variables of Duhamel's formula, and at the final
time. The argument below proves it with arbitrarily small polynomial
error on good spectral paths.

We start with one particle. Then the short flow is a unitary rotation of
the eigenvectors without Ornstein--Uhlenbeck killing, and the mean
squared entries of the rotation matrix $M_{s,u}$ solve the one-particle
equation. We use the weight argument of
Lemma~\ref{lem:t26-relative-localization}, with smoothing distance
  $d=\ell/\Naux $  and a weight centred at $\lambda_k(s)$ for a fixed bulk label
$k$. The initial centre is measurable at time $s$. By Weyl's
inequality, eigenvalues move by at most $C\sqrt\tau=o(d)$, so the same
label stays in the central part of the weight. By the upper count, a
rank distance larger than $R$ at time $s$ costs the weight factor
  $e^{-cR/\ell}=e^{-c\Naux ^\epsilon}$. The exact cancellation of the
eigenvalue drift gives the Gronwall exponent   $C\tau(\Naux /\ell)\log \Naux =o(1)$.
Therefore the expected row mass
$\E_{\rm aux}\sum_{|j-k|>R}|M_{s,u}(k,j)|^2$, restricted to $\mathcal R$,
is exponentially small, for every deterministic pair $s,u$.

The column estimate follows from the same forward calculation, with the
opposite weight that increases away from the initial central label. The
initial weighted mass is one, and reaching a far row costs
  $e^{c\Naux ^\epsilon}$. This also controls transitions from the local
support to configurations outside the bulk window. Neither estimate
needs to reverse the roles of $s$ and $u$.

For $r$ particles we lift the one-particle rotation to the symmetric
power. On the normalized occupation basis of   $\operatorname{Sym}^r\C^\Naux $,
the unitary $\operatorname{Sym}^rM$ has the matrix entries
\[
 \frac{\operatorname{per}M[\xi,\zeta]}{\sqrt{\xi!\zeta!}}.
\]
Averaging the squared modulus over the short Brownian rotations gives
$U_{s,u}(\xi,\zeta)$. Indeed, differentiating the creation and
annihilation coefficients $\sqrt{\xi_i(1+\xi_j)}$ gives exactly
\eqref{eq:t26-complex-moment-generator}, restricted to short jumps.
Now suppose that the ordered matching distance between $\xi$ and $\zeta$
exceeds $rR$. Then every term of the permanent contains an entry whose
rank displacement exceeds $R$. Cauchy--Schwarz over the at most $r!$
terms, then summation over ordered destination labels using the unit
row sums of $|M_{ij}|^2$, gives
\[
 \sum_{d(\xi,\zeta)>rR}U_{s,u}(\xi,\zeta)
 \le C_r\sum_{i\in\operatorname{supp}\xi}\xi_i
                     \E_{\rm aux}\sum_{|j-i|>R}|M_{s,u}(i,j)|^2 .
\]
So the multiparticle flow inherits the finite speed of the
one-particle flow.

At fixed $r$ there are only polynomially many configurations. By
Tonelli's theorem in one or two time variables and Markov's inequality
over $\mathcal L$, there is a good path set on which all tail errors used
below have time integral at most   $\Naux ^{-D}$, and the errors at the final
time are at most   $\Naux ^{-D}$. Sums over configurations and fixed polynomial
coefficients are paid for by increasing $D$ first. We write $e_D(u)$ for
the resulting nonnegative error, with   $\int_0^\tau e_D(u)du\le \Naux ^{-D}$.
We never use a uniform lower bound on eigenvalue gaps, or time reversal.

\emph{Step 3: long jumps and flattening.} We now prepare the maximum
principle. Fix a target configuration $\xi_0$ in the shrinking window.
To make the initial data harmless far from $\xi_0$, we replace them by
one there. We put
\[
 a_\xi=\frac2w\int_{w/2}^w1_{\{d(\xi,\xi_0)\le v\}}dv,\qquad
 \operatorname{Av}f=a_\xi f+(1-a_\xi),\qquad
 g_u=U_{0,u}\operatorname{Av}F_r(0).
\]
So $a_\xi$ is the ramp of Figure~\ref{fig:cm-scales}, and $g_u$ evolves
the flattened data by the short flow. The ramp satisfies the Lipschitz
bound $|a_\xi-a_\zeta|\le2d(\xi,\zeta)/w$. All configurations in fixed
multiples of the window of radius $w$ have energies   $E+O_r(w/\Naux )$, by the
upper and lower counts at time zero and the displacement $\sqrt\tau$.
Consequently $g_0$ is bounded by   $C\Naux ^{C_r\epsilon}$  globally, and $g_u$
has the same bound.

Next we show that the long jumps, which $g_u$ ignores, change $F_r$ very
little. For a local occupied site $i$ put   $h_\ell=\ell/\Naux $. Positivity and
\eqref{eq:t26-moment-local-law} give

\[
 \sum_{|j-i|\ge\ell}c_{ij}F_1(u,\delta_j)
 \le\frac{C}{h_\ell}\,
       \E[q^*\Im G_u(\lambda_i+ih_\ell)q\mid\mathcal L]
 \le C \Naux ^{C\epsilon}\Naux /\ell.
\]

The idea is that a jump over distance at least $\ell$ has rate at most
of order $h_\ell^{-1}$ times a resolvent kernel at height $h_\ell$, and the
local law bounds that kernel. In the multiparticle case we first remove
the particle at $i$. On $\mathcal A$ the product of the other $r-1$
powers is at most   $\Naux ^{C_r\epsilon}$, and the complement of $\mathcal A$
is paid for as in Step~1. There are at most $r$ occupied destinations,
and they obey the same rough bound. The unweighted tail is at most
  $C\Naux /\ell$. Thus, locally,
  $|(\mathcal B-\mathcal S)F_r(u,\xi)|\le C\Naux ^{C_r\epsilon}\Naux /\ell$, plus a
negligible error, where $\mathcal B$ is the full generator.

Duhamel's formula and the propagation bound of Step~2 then give
  $|F_r(u)-U_{0,u}F_r(0)|\le C\Naux ^{C_r\epsilon}\Naux \tau/\ell+e_D(u)$  on the whole
local window. At the target, and wherever $a_\xi=1$, the flattening
changes the data only beyond distance $w/2\gg R$, so it adds only a
propagation tail. Hence

\begin{equation}\label{eq:t26-long-short-moments}
 |F_r(u,\xi)-g_u(\xi)|
 \le C \Naux ^{C_r\epsilon}\Naux \tau/\ell+e_D(u),
\end{equation}

with a terminal error   $\Naux ^{-D}$  at $u=\tau$. The flattening error
vanishes at time zero in the inner window and is paid for by
propagation. Globally, the long coefficients are polynomially bounded
on $\mathcal R$, so the exterior errors are absorbed by the same tail
estimates.

The maximum principle will be applied to $g_u$, so we also need to know
that flattening almost commutes with the flow. For each flattening
radius $v$, propagation makes the commutator $U\operatorname{Av}-\operatorname{Av}U$
negligible unless $|d(\xi,\xi_0)-v|\le rR$. These radii form a set of
length at most $2rR$ inside an interval of length $w/2$. On that set we
use the rough local moment bound and the supremum contraction.
Integrating over $v$ gives   $C\Naux ^{C_r\epsilon}\ell/w$, with the factor
  $\Naux ^\epsilon$  absorbed into $C_r\epsilon$, plus $e_D(u)$. The same
argument applies after summing the moving destination against a positive
resolvent kernel, since the scalar law bounds its total mass.

\emph{Step 4: the removal identity and the averaged moments.} The
maximum principle needs the value of the flow at a maximum. That value is
a resolvent-weighted average of moments with one particle moved, and we
now compute it. If $i$ is occupied, define
$P_i=\prod_hX_h^{\xi_h-1_{h=i}}/(\xi_h-1_{h=i})!$, the product of the
overlaps of all particles except one at $i$. Moving that particle to an
unoccupied label $j$ gives the exact identity
\begin{equation}\label{eq:t26-complex-removal}
 F_r(u,\xi^{i,j})=\E[P_iX_j\mid\mathcal L].
\end{equation}
For occupied $j\ne i$ the right side is $(1+\xi_j)F_r(u,\xi^{i,j})$, and
for $j=i$ it is $\xi_iF_r(u,\xi)$. So the factorials are accounted for
exactly, and adding or removing all occupied destinations in a sum
weighted by a resolvent kernel at height $\eta$ costs at most
  $C\Naux ^{C_r\epsilon}/(\Naux \eta)$.

Summing \eqref{eq:t26-complex-removal} over $j$ against the resolvent
kernel produces $q^*\Im G_u q$, to which the directional local law
applies. We multiply the local law by the nonnegative $P_i$ before
conditioning. On $\mathcal A$ we use the rough bound already stated, and
on its complement the deterministic polynomial bound. This gives

\[
\begin{aligned}
 \frac{1}\Naux \sum_j\frac{\eta F_r(u,\xi^{i,j})}{(\lambda_i-\lambda_j)^2+\eta^2}
 &=\Im m_{c-\tau+u}(\lambda_i+i\eta)
                         F_{r-1}(u,\xi-\delta_i)\\
 &\quad+O\!\left(\Naux ^{C_r\epsilon}
       [(\Naux \eta)^{-1/2}+w/\Naux +\tau+(\Naux \eta)^{-1}]\right).
\end{aligned}
\]

At occupied destinations the sum means the restored sum, with the error
just bounded. The scalar factor $\Im m$ is kept. This is the key
identity: an average of $r$-particle moments equals an $(r-1)$-particle
moment, which is close to one by the induction hypothesis.

To use this for $g_u$ instead of $F_r$, we apply the commutator bound
and \eqref{eq:t26-long-short-moments}. Pulling the ramp $a_\xi$ out of
the resolvent sum costs at most   $C\Naux ^{C_r\epsilon}(\Naux \eta/w)\log \Naux $.
Indeed, the coefficient difference is at most $\min(1,2|j-i|/w)$.
Dyadic shells between   $\Naux \eta$  and $w$ give   $C\Naux \eta/w$  each, by the local
upper count and the rough moment bound, and the shells beyond $w$ are
controlled by the tail of the positive resolvent kernel. A Lipschitz
bound   $1/\Naux $  for $a_\xi$ would be wrong here, which is why we use the
shell argument. The result is

\begin{equation}\label{eq:t26-averaged-complex-moments}
 \frac{1}\Naux \sum_j\frac{\eta g_u(\xi^{i,j})}{(\lambda_i-\lambda_j)^2+\eta^2}
 =\Im m_{c-\tau+u}(\lambda_i+i\eta)
       [a_\xi F_{r-1}(u,\xi-\delta_i)+(1-a_\xi)]
       +O(\Naux ^{C_r\epsilon}\mathfrak e)+e_D(u),
\end{equation}

where all the errors collected so far are

\[
 \mathfrak e=\ell/w+\Naux \tau/\ell+(\Naux \eta)^{-1/2}
               +w/\Naux +\tau+(\Naux \eta/w)\log \Naux +(\Naux \eta)^{-1}.
\]

The same argument works for all fixed nested windows. Only their fixed
constants depend on $r$.

\emph{Step 5: the maximum principle and the induction on $r$.} For
$r=0$ the moment is one. The induction is carried out on the single flow
$H_u=F+G_{c-\tau+u}$, $0\le u\le\tau$, conditionally on its fixed
eigenvalue path $\mathcal L$, and only the start time of the damping moves.
Put $t_k=1-2^{-(k+1)}$. The assertion at level $k$ is that, outside a set
of paths of probability at most $\Naux^{-D}$,
$|F_k(u,\xi)-1|\le C_k\Naux^{-\kappa_k}$ for all $u\in[t_k\tau,\tau]$ and all
$k$-particle configurations $\xi$ in the $k$-th nested window. Along the
whole flow the variance $c-\tau+u$ lies in $[c_0-\tau,c_0]\subset
[c_0-s_0,c_0]$ for large $\Naux$, so Steps~1--4 apply at every $u$, and the
moments $F_{r-1}(u,\cdot)$ in \eqref{eq:t26-averaged-complex-moments} are
those of the same flow, conditioned on the same $\mathcal L$, to which the
assertion at level $r-1$ applies. Suppose that assertion holds and put
$s=t_{r-1}\tau$. On $[s,\tau]$ the right side of
\eqref{eq:t26-averaged-complex-moments} therefore equals
  $\Im m\,(1+O(\Naux ^{-\kappa_{r-1}}+\Naux ^{C_r\epsilon}\mathfrak e))$, since
$F_{r-1}\approx1$ by induction.

Now look at a configuration where $g_u$ is largest. There every jump
lowers or keeps the value, so the jumps to occupied destinations, which
are nonpositive, can be discarded. On the other jumps we replace $c_{ij}$
by the smaller kernel   $[\Naux ((\lambda_i-\lambda_j)^2+\eta^2)]^{-1}$.
Extending the short sum to all destinations costs
  $C\Naux ^{C_r\epsilon}\Naux \eta/\ell$  after multiplication by $\eta$. The scalar
local law and the removal of the finitely many occupied sites then give
a damping coefficient at least $c_r/\eta$. As a result, the upper Dini
derivative of the positive excess of the maximum satisfies

\[
 D^+\big(\max_\xi g_u(\xi)-1\big)_+
 \le-\frac{c_r}{\eta}\big(\max_\xi g_u(\xi)-1\big)_+
   +\frac{C\Naux ^{C_r\epsilon}}{\eta}
        [\mathfrak e+\Naux \eta/\ell+\Naux ^{-\kappa_{r-1}}]
   +\frac{e_D(u)}{\eta}
\]

whenever the maximizing configuration lies in the local window. The same
inequality holds for the deficit of the minimum. Outside a slightly
larger window, propagation bounds $|g_u-1|$ by $e_D(u)$.

This integrated error is enough. We discard the times where
  $e_D(u)>\Naux ^{-D/2}$, and their total length is at most   $\Naux ^{-D/2}$. At those
times the maximum cannot increase and the minimum cannot decrease,
because the generator is that of a Markov process. At the remaining
times we use the damping inequality, or the smaller exterior bound.
We use the damping inequality only on $[s,\tau]$; on $[0,s]$ the
maximum principle alone keeps the excess below its initial polynomial
bound. For $u\ge t_r\tau$ the damping time is at least
$u-s\ge(t_r-t_{r-1})\tau=2^{-(r+1)}\tau$, and choosing $D$ large pays every
integrated error. Since $\tau/\eta=\Naux^\epsilon$, the initial polynomial
excess is suppressed by a factor $e^{-c_r2^{-(r+1)}\Naux^\epsilon}$, which is
smaller than every power of $\Naux$ for fixed $r$.

It remains to check that every error is a negative power of   $\Naux $. The
main terms are

\[
 \ell/w=\Naux ^{-\theta/6},\quad \Naux \tau/\ell=\Naux ^{-2\theta/3},
 \quad w/\Naux =\Naux ^{-\theta/6},\quad \Naux \eta/\ell=\Naux ^{-2\theta/3-\epsilon},
\]

and the other terms in $\mathfrak e$ also decrease by fixed powers. We
choose $\epsilon<\theta/(100C_r)$, and if necessary we reduce $\kappa_r$
below $\kappa_{r-1}$. The induction then gives
$|g_u(\xi_0)-1|\le C\Naux^{-\kappa_r}$ for $u\in[t_r\tau,\tau]$, and
\eqref{eq:t26-long-short-moments} gives the same bound for $F_r$. The
finitely many enlarged windows stay inside the fixed regular rectangle.
All local-law events and propagation sums were uniform in the ordered
labels, and the polynomial unions over those labels are paid for by the
arbitrary exceptional exponent. This proves the lemma.
\end{proof}

\subsection{The exact scalar Palm weight}\label{subsec:s8-the-exact-scalar-palm-weight}

 The conditional moment estimate must be integrated with the actual
weight that selects an eigenvalue at the prescribed energy. This
completes Theorem~\ref{thm:fresh-palm} without changing the Gaussian
variance. We express the weight through the original scalar Gaussian,
localize its contribution in energy, and pay for exceptional paths
before integrating.

 \begin{proof}[Proof of Theorem~\ref{thm:fresh-palm}]
Use the original variance-$c_0$ GUE and the short flow of
Lemma~\ref{lem:t26-conditional-moments}; no auxiliary scalar variance
 is introduced. If $\nu_i$ are its terminal eigenvalues, its scalar
part and Palm weight are exactly

\[
 g=\frac{\sum_i\nu_i-\Tr F}{\Naux },\qquad
 w_{\ell,E}=\varphi_{\sqrt{c_0}/\Naux }(E-\nu_\ell+g).
\]

Thus $w_{\ell,E}$ is a nonnegative $\mathcal L$-measurable variable
bounded by   $C_{c_0}\Naux $. Except on an event of arbitrary polynomial
decay,   $|g|\le\tfrac14\Naux ^{-1+\epsilon}$. If
$|\nu_\ell-E|>\Naux ^{-1+\epsilon}$  there, the weight is bounded by
  $C\Naux \exp(-c\Naux ^{2\epsilon})$.

 On the remaining good paths the lemma, with $r$ particles at
$\ell$, gives

\[
 \left|\E[|\sqrt \Naux q^*u_\ell|^{2r}\mid\mathcal L]
                     -r!(q^*\Gamma q)^r\right|
 \le C_r\Naux ^{-\kappa_r}(q^*\Gamma q)^r.
\]

Multiply by the actual weight before integrating. The good contribution
is at most

\[ C_r\Naux ^{-\kappa_r}(q^*\Gamma q)^r\E w_{\ell,E}. \]

On bad paths use
  $|\sqrt \Naux q^*u_\ell|^{2r}\le \Naux ^r\|q\|^{2r}$,
$q^*\Gamma q\ge C\Naux ^{-2\max(A,1)}\|q\|^2$ and $w\le C\Naux $.
Taking the exceptional exponent after $r,K,A$  bounds this contribution
by   $C_{r,K}\Naux ^{-K}(q^*\Gamma q)^r$ ; the Gaussian scalar tail and the
exponentially small off-window weight obey the same bound.
The scalar decomposition is orthogonal in GUE, so the eigenvectors
and $w$ have exactly the trace-free law in
\eqref{eq:palm-scalar-weight}. This proves
\eqref{eq:fresh-palm-moments} for every ordered rank, including ranks
with arbitrarily small $Z_{\ell,E}$. No error was divided by that mass.
The zero direction follows directly, and homogeneity restores every
complex $q$.
\end{proof}
\section{Local estimates at the selected energy}
\label{sec:physical-local-window}
The fresh-Gaussian estimates of Section~\ref{sec:fresh-gaussian} assume a
deformation with a regular spectral window. The deformation produced by the
ring may have an exterior eigenvalue, a pole, close to the selected energy.
This section shows that the estimates still hold at the energy selected by
the full matrix. It is used in Part~I.
Lemma~\ref{lem:physical-local-window} gives the counting and eigenvector
bounds at the selected energy, even close to an exterior pole.
Lemma~\ref{lem:t26-pair-pole-exclusion} shows that, under the pair law, the
selected eigenvalue comes within $\varepsilon$ of an exterior eigenvalue only
with probability $O(n\Naux^2\varepsilon)$.
We first control the proximity of poles under the law of the selected
rank. In Step 1-4 below, we show that a short window contains at most one
exterior eigenvalue with high probability, and that removing it leaves a
deformation of controlled size. In Step 5-7 below, we carry the local
estimates through that single remaining pole, which completes
Lemma~\ref{lem:physical-local-window}.

Here $S=\{i,i+1\}$, $d=2W$, $s=\min(b,1-2b)$, $t=s/d$ and
$A(E)=F(E)+G_d(t)-E$. Let $I$ be a fixed compact bulk interval
containing the classical locations of the ranks under consideration. Here $m_E$ is the normalized Dyson solution for $K_0=F(E)-E$, with normalized trace $d^{-1}\Tr$ and $d=2W$; real arguments mean its regular boundary value.
For a Hermitian observable $B$, set
\[
 j_B(E,\lambda)=d^{-1}\Tr B\Gamma_E(\lambda),\qquad
 \Gamma_E(\lambda)=s|F(E)-E-\lambda-sm_E(\lambda)|^{-2}.
\]
The regularity constants below are fixed; probability estimates are
uniform for $4\le n\le TW$, both parities and deterministic bulk ranks.

\begin{lemma}\label{lem:physical-local-window}
 Fix $\gamma_0=1/2$.  Choose a fixed $0<\delta<1$ and a derivative order $k$ from the fixed bulk and regularity parameters, before $\gamma,\epsilon,D$ and the spectral-test bound $K_0$. For fixed $0<\gamma<\gamma_0$, $\epsilon>0$
and $D>0$, outside an event of
original probability at most $C_DW^{-D}$, the following statements
hold at every selected bulk energy.
\begin{enumerate}[label=(\roman*)]
\item For $\eta=W^{-1+\gamma}$,
\begin{equation}\label{eq:physical-local-count}
 N_{A(E)}([-C\eta,C\eta])\ge cW\eta,\qquad
 \sum_{\lambda\in\operatorname{Spec}A(E)}
       \frac{\eta^2}{\lambda^2+\eta^2}\le CW\eta.
\end{equation}
\item If $B$ is either (a) a measurable energy-independent Hermitian contraction
of the exterior, or (b) $f(F(E)-E)$, where $f$ is a fixed real-valued
$C^\infty$ spectral test with $|f|\le1$, all derivatives of orders
$0\le j\le k$ bounded by a fixed $K_0$, and finite limits at
$\pm\infty$,
then, for $|\lambda|\le\delta$ and every normalized eigenvector
$u_\lambda$,
\begin{equation}\label{eq:physical-window-eth}
 |u_\lambda^*Bu_\lambda-j_B(E,\lambda)|
       \le CW^{-1/2+\epsilon}.
\end{equation}
The assertion includes $B=P_i-P_{i+1}$ and the bounded smooth spectral functions in class (b). For these functions the limiting compression
at an infinite rank-one parameter is used.
\end{enumerate}
The same probability statements hold under the pair law, after
increasing the original failure exponent by one.
\end{lemma}

Let us explain in simple words what the above lemma is claiming.  At the selected energy level by the whole ring, the matrix on the pair obtained by eliminating the rest of the ring behaves like a regular deformed Gaussian matrix, even though its deformation can have a pole near that energy. Part~(i) says that a window of width a little larger than one eigenvalue spacing contains about the expected number of eigenvalues. Part~(ii) says that each eigenvector in a bulk window sees a bounded observable essentially through its deterministic centre $j_B(E,\lambda)$, the normalized trace of $B$ against the covariance $\Gamma_E(\lambda)$, with an error of order $W^{-1/2}$ up to a small power. Here \emph{original probability} refers to the unweighted law $\Prob$ of the ring. The pair law has density at most $CW$ with respect to it, which costs one power of $W$.

In the following lemma, we build the necessary tools for this. We split the ring into the exposed pair $S$
of two adjacent blocks, of dimension $d=2W$, and the exterior $R$, which
is everything else. When we eliminate the exterior by a Schur
complement, the result involves the resolvent $(H_R-\lambda)^{-1}$ of
the exterior. This resolvent has poles at the eigenvalues of $H_R$. So
the Schur complement is well behaved at $\lambda_\alpha$ only if
$\lambda_\alpha$ stays away from $\operatorname{Spec}H_R$. The lemma
says that, under the size-biased pair law $\Prob^\sharp_{\alpha,S}$, the
eigenvalue of rank $\alpha$ comes within $\varepsilon$ of an exterior
eigenvalue only with probability $O(n d^2\varepsilon)$. The bound holds
for every exterior, with no assumption on its norm or on the spacing of
its eigenvalues.

\begin{lemma}[Ordered-rank pole exclusion under the physical pair law]
\label{lem:t26-pair-pole-exclusion}
Use the exact original Gaussian ring and splitting
\eqref{eq:pair-splitting}, with $d=2W$ and fixed $s>0$. For every rank
  $1\le\alpha\le N$, every exterior-measurable Borel set
$B(\omega)\subset\R$ of finite Lebesgue measure, and $G=G^0+\zeta I_S$,
where $G^0$ is traceless, one has
\begin{equation}\label{eq:t26-weighted-spectral-averaging}
 \E\!\left[m_S\one_{\{\lambda_\alpha\in B(\omega)\}}
                       \mid\omega,G^0\right]
 \le\frac{d}{\sqrt{2\pi s}}\,|B(\omega)|.
\end{equation}
Here the conditional expectation integrates only over
$\zeta\sim N(0,s/d^2)$, and $m_S=\norm{P_S\psi_\alpha}^2$ is the actual
mass of the eigenvector of ordered rank $\alpha$ on the pair.
Consequently, uniformly over all exterior data, both parities, and all
$n,W$, for every $\varepsilon>0$,
\begin{equation}\label{eq:t26-pole-probability}
 \Prob^\sharp_{\alpha,S}
 (\operatorname{dist}(\lambda_\alpha,\operatorname{Spec}H_R)
                      <\varepsilon)
 \le\frac{nd(d+1)}{\sqrt{2\pi s}}\varepsilon.
\end{equation}
In particular this is $O_{s,T}(W^{3-M})$ for $\varepsilon=W^{-M}$ and
$4\le n\le TW$. The bound involves the unnormalized pair weight, with no
denominator given by the mass of a rank sector. The lemma is a step in
the physical application of Theorem~\ref{thm:fresh-palm}, and it does
not assert that theorem.
\end{lemma}

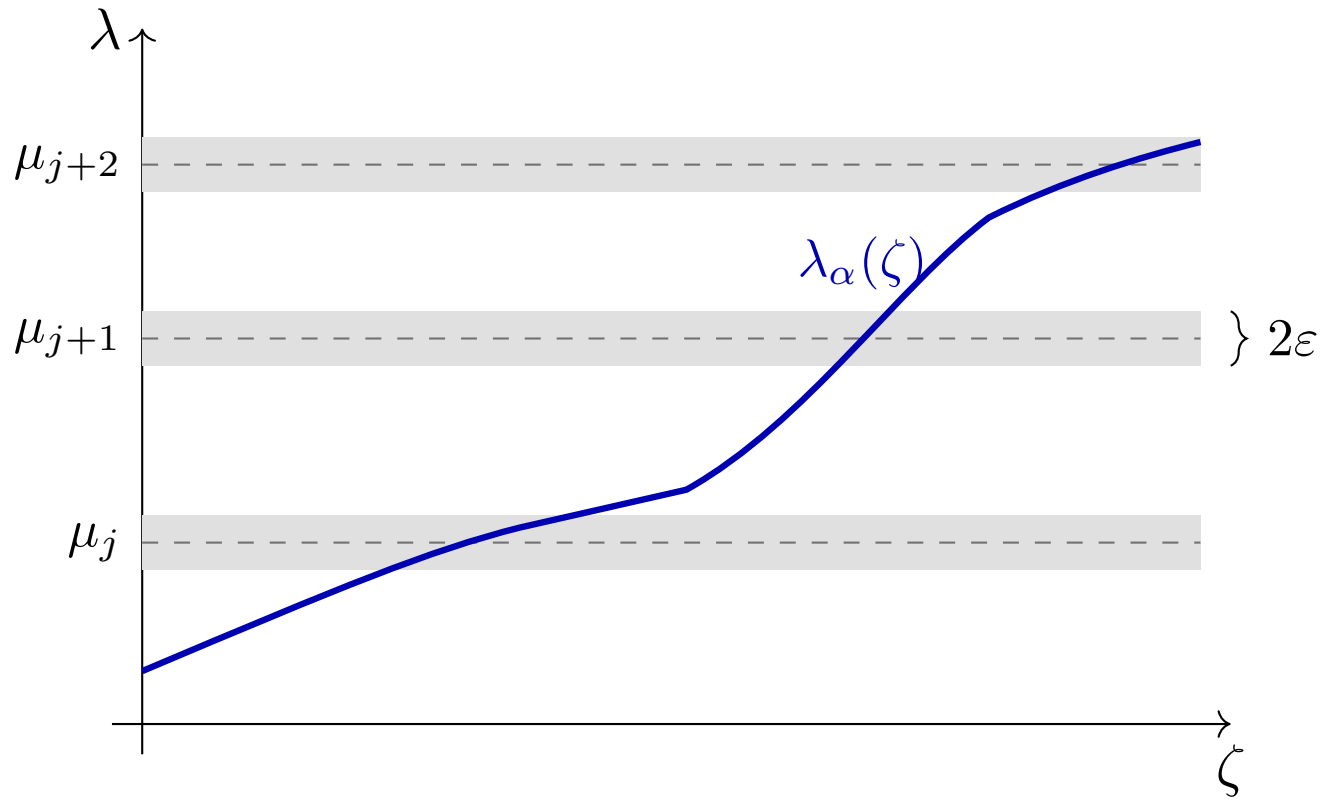
\begin{figure}[ht]
\centering
\begin{tikzpicture}[x=1cm,y=1cm]
 \draw[->] (-.2,0)--(7.2,0) node[below] {$\zeta$};
 \draw[->] (0,-.2)--(0,4.6) node[left] {$\lambda$};
 \foreach \y in {1.2,2.55,3.7}{
   \fill[black!12] (0,\y-.18) rectangle (7,\y+.18);
   \draw[dashed,black!55] (0,\y)--(7,\y);
 }
 \node[font=\small,left] at (0,1.2) {$\mu_{j}$};
 \node[font=\small,left] at (0,2.55) {$\mu_{j+1}$};
 \node[font=\small,left] at (0,3.7) {$\mu_{j+2}$};
 \draw[very thick,blue!70!black]
   (0,.35) .. controls (1.3,.9) and (1.9,1.15) .. (2.5,1.3)
           -- (3.6,1.55)
           .. controls (4.4,2.0) and (5.0,2.9) .. (5.6,3.35)
           .. controls (6.1,3.6) and (6.6,3.75) .. (7,3.85);
 \node[font=\small,blue!70!black,anchor=west] at (4.2,3.05) {$\lambda_\alpha(\zeta)$};
 \draw[decorate,decoration={brace,amplitude=3pt,mirror}] (7.2,2.37)--(7.2,2.73)
   node[midway,right=3pt,font=\small] {$2\varepsilon$};
\end{tikzpicture}
\caption{The mechanism in the proof of
Lemma~\ref{lem:t26-pair-pole-exclusion}. The exterior eigenvalues
$\mu_j$ (dashed) do not depend on the scalar part $\zeta$ of the fresh
pair block. The shaded bands are their $\varepsilon$-neighbourhoods, and
the set $B$ is the union of the bands whose labels are allowed by
interlacing. The ordered eigenvalue $\lambda_\alpha(\zeta)$ (solid) is
nondecreasing with slope $m_S\le1$. Where the slope is small,
$\lambda_\alpha$ may stay in a band for a long $\zeta$-interval, but
there the weight $m_S$ is also small. The area formula
\eqref{eq:pe-area} says that the weighted time spent in $B$ is at most
$|B|$, whatever the shape of the curve.}
\label{fig:pe-monotone}
\end{figure}

\begin{proof}
We fix everything except the scalar part of the
fresh Gaussian block on the pair, a single real Gaussian variable that shifts
the whole pair by the same amount. As this variable increases, the selected
eigenvalue moves up, and its speed equals the eigenvector's mass on the pair
(Step~1). This speed is exactly the weight that defines the pair law, so the
weighted time that the eigenvalue spends in a set of energies is at most the
length of that set, by a one-dimensional change of variables (Step~2).
Figure~\ref{fig:pe-monotone} shows the mechanism. It remains to choose the
set. By interlacing, which says that the eigenvalues of a matrix and of a
submatrix alternate, only a bounded number of exterior eigenvalues, about the
pair dimension, can be close to the selected one (Step~3). Step~4 adds up the
resulting bound.

\emph{Step 1: the eigenvalue as a function of $\zeta$.} We condition on
$(\omega,G^0)$, so the only remaining randomness is $\zeta$. Then
\[
 H(\zeta)=H(0)+\zeta P_S.
\]
Since $0\preceq P_S\preceq I$, the min--max principle shows that each
ordered eigenvalue $\lambda_\alpha(\zeta)$ is nondecreasing and
$1$-Lipschitz in $\zeta$. In particular it is differentiable at almost
every $\zeta$, and there the Hellmann--Feynman formula gives
\begin{equation}\label{eq:pe-hellmann-feynman}
 \lambda_\alpha'(\zeta)=\psi_\alpha^*P_S\psi_\alpha=m_S.
\end{equation}
The degenerate points, where $\lambda_\alpha$ is a multiple eigenvalue,
can be discarded. Alternatively, the area formula below applies to the
almost-everywhere derivative of the ordered Lipschitz function.

\emph{Step 2: the weighted time in $B$.} For a monotone absolutely
continuous real function $f$, the one-dimensional area formula gives
\begin{equation}\label{eq:pe-area}
 \int_\R f'(x)\one_{\{f(x)\in B\}}\,dx
 =\int_B\#\{x:f(x)=y\}\,dy\le |B| .
\end{equation}
The multiplicity $\#\{x:f(x)=y\}$ is at most one, except at the values
taken on nontrivial intervals where $f$ is constant. There are at most
countably many such values, so they have Lebesgue measure zero. We apply
\eqref{eq:pe-area} to $f=\lambda_\alpha$, whose derivative is $m_S$ by
\eqref{eq:pe-hellmann-feynman}. The set $B(\omega)$ is
exterior-measurable, so it is fixed under the conditioning. The density
of $\zeta\sim N(0,s/d^2)$ is at most $(2\pi s/d^2)^{-1/2}=d/\sqrt{2\pi s}$.
Hence
\[
 \E\!\left[m_S\one_{\{\lambda_\alpha\in B(\omega)\}}
                       \mid\omega,G^0\right]
 =\int_\R\lambda_\alpha'(\zeta)\one_{\{\lambda_\alpha(\zeta)\in B(\omega)\}}
   \,p(\zeta)\,d\zeta
 \le\frac{d}{\sqrt{2\pi s}}\,|B(\omega)|,
\]
where $p$ is the density of $\zeta$. This proves
\eqref{eq:t26-weighted-spectral-averaging}.

\emph{Step 3: which exterior eigenvalues can be close.} To apply Step~2
we need a set $B$ of finite measure. The whole exterior spectrum has up
to   $N-d$  points, which would give a bound proportional to   $N$. Instead,
we use the rank $\alpha$ to locate $\lambda_\alpha$ among the exterior
eigenvalues. We write these as   $\mu_1\le\cdots\le\mu_{N-d}$  and set

\[
 I_\alpha=\{j:1\le j\le N-d,\ \alpha-d\le j\le\alpha\}.
\]

The exterior matrix $H_R$ is obtained from $H$ by deleting the $d$
coordinates of the pair. Cauchy interlacing therefore gives
\[
 \mu_{\alpha-d}\le\lambda_\alpha\le\mu_\alpha,
\]
whenever the corresponding index exists. Let $N_R(\lambda)$ be the
number of exterior eigenvalues below $\lambda$. Away from the null event
that $H$ and $H_R$ have a common eigenvalue, the inequalities are strict,
and so
\[
 \alpha-d\le N_R(\lambda_\alpha)\le\alpha-1 .
\]
The exterior eigenvalues nearest to $\lambda_\alpha$ are the ones with
indices $N_R(\lambda_\alpha)$ and $N_R(\lambda_\alpha)+1$, omitting an
index that does not exist at an end of the spectrum. Both lie in
$I_\alpha$. Hence, if $\lambda_\alpha$ is within $\varepsilon$ of some
exterior eigenvalue, it is within $\varepsilon$ of one with index in
$I_\alpha$.

\emph{Step 4: conclusion.} We take for $B(\omega)$ the union of the
$\varepsilon$-neighbourhoods of the $\mu_j$ with $j\in I_\alpha$. This
set is exterior-measurable, and since $|I_\alpha|\le d+1$, its length is
at most $2(d+1)\varepsilon$. By Step~3,
\[
 \{\operatorname{dist}(\lambda_\alpha,\operatorname{Spec}H_R)<\varepsilon\}
 \subset\{\lambda_\alpha\in B(\omega)\}.
\]
We apply \eqref{eq:t26-weighted-spectral-averaging} to this set and
integrate over the conditioning. The size-biased law has density
$\frac{n}{2}\,m_S$ with respect to $\Prob$, by \eqref{eq:intro-cancel}, so
\[
 \Prob^\sharp_{\alpha,S}
 (\operatorname{dist}(\lambda_\alpha,\operatorname{Spec}H_R)<\varepsilon)
 \le\frac n2\cdot\frac d{\sqrt{2\pi s}}\cdot2(d+1)\varepsilon
 =\frac{nd(d+1)}{\sqrt{2\pi s}}\varepsilon.
\]
This is \eqref{eq:t26-pole-probability}. For $\varepsilon=W^{-M}$,
$d=2W$ and $n\le TW$, the right side is at most
$C_{s,T}W\cdot W^2\cdot W^{-M}=O_{s,T}(W^{3-M})$. The argument made no
assumption on the norm of the exterior or on the spacing of its
eigenvalues.
\end{proof}

The estimates allow every fixed $0<\gamma<\gamma_0$ and $\epsilon>0$.
They do not assert a local law at height $W^{-1}$ or uniformity in a
power tending to zero. Their role is to pay rare events in a later
argument whose final bound has no power loss.

\begin{proof}[Proof of Lemma~\ref{lem:physical-local-window}]
Before the steps, we fix some
parameters and record a few facts used throughout the proof. First, we
choose the radius $\delta$ of the spectral window in part~(ii) smaller
than the radius of the energy neighbourhood on which the regularity
certificate of Lemma~\ref{lem:regularity-certificate} applies. This
choice is made before the loss exponent $\epsilon$ and the tail exponent
$D$, so that neither of them affects it. Second, if no exterior
eigenvalue lies within distance $\varepsilon_0$ of the energy, then the
truncation of Step~2 changes nothing, that is, $F=\widehat F$. This is
the case $\rho=0$ of the rank-one decomposition in Step~5, so it needs no
separate treatment. Third, the hypotheses of the eigenvector theorem,
Theorem~\ref{thm:external-eth}, follow from the regularity condition
\eqref{eq:fresh-regularity}. Indeed, the imaginary part of the Stieltjes
transform is bounded below there, which gives the norm bound
$\|M_s\|\le(sv_0)^{-1}$ on the deterministic approximation and the
stability bound $\Re(1-s\av{M_s^2})\ge2s(\Im m_s)^2$. Differentiating the
Dyson equation implicitly then shows that $m_s$ and $M_s$ have bounded
derivatives. The eigenvector theorem does not by itself bound the number
of eigenvalues in the spectral window; that bound comes from the local
law. Fourth, when we compress the matrix in Step~7, its variance becomes
$s'=s(1-1/d)$. This value stays in a fixed compact subset of
$(0,\infty)$, so the eigenvector theorem applies to the compressed matrix
with the same constants. Fifth, the value of the semicircle Stieltjes
transform in \eqref{eq:uniform-pole-certificate} is
$m_{\rm sc}(E)=\bar a$, the complex conjugate of the matched value $a$ of
\eqref{eq:matched-parameters}. Finally, a word on notation. The constant
$L_0$ used in Step~4 to compare Poisson kernels at two heights is a fixed
number. It is unrelated to the growing cutoff $L_0=W^{1/4}$ used later,
in the proof of Proposition~\ref{prop:physical-real-loop}. Likewise, the
radius of the regularity certificate and the cutoff radius in the smooth
centre of Lemma~\ref{lem:smooth-center} play different roles.

The lemma concerns the Schur complement $A(E)=F(E)+G_d(t)-E$ of the pair
$S$, where
\[
 F(E)=V-C^*(H_R-E)^{-1}C .
\]
The fresh Gaussian part $G_d(t)$ is harmless. The difficulty is the
deformation $F(E)$. It contains the exterior resolvent, which has a pole
at every eigenvalue of $H_R$. So $F(E)$ is not bounded, and it is not even
defined when $E$ is an exterior eigenvalue. The estimates of
Sections~\ref{sec:inputs} and~\ref{sec:fresh-gaussian} require a
deformation of polynomial size with a regular spectral window. They
therefore cannot be applied to $F(E)$ directly. A second difficulty is
that the statement must hold at the selected energy $E=\lambda_\alpha$,
which is random and depends on the exterior.

\emph{Plan of the proof.} We handle the poles in two stages, as in
Figures~\ref{fig:plw-truncation} and~\ref{fig:plw-pole}. In the first stage
(Steps~1--4) we show that, with overwhelming probability, every short
window around a bulk energy contains at most one exterior eigenvalue. We
then replace the exterior resolvent by a truncated version, which is
bounded and Lipschitz in the energy and differs from $F(E)$ by rank at
most one. For the truncated deformation the regularity certificate and
the local law apply on an energy grid. Since a rank-one change is harmless
for counting, this proves part~(i). In the second stage (Steps~5--7) we
treat the remaining single pole. We write the singular part as one
rank-one term with a coefficient $\rho$. While $\rho$ is at most a fixed
power of $W$, we use the eigenvector estimates on a grid in $(E,\rho)$.
When $\rho$ is larger, the matrix is close to its compression orthogonal
to the singular direction, and we apply the same estimates there. Finally
(Step~8) we pass from the deterministic grids to the selected energy and
from the original law to the pair law. All these estimates are
simultaneous over energies, so we need no independence between the
events at different grid points.

\begin{figure}[ht]
\centering
\begin{tikzpicture}[x=1cm,y=1cm]
 \begin{scope}
  \draw[->] (-.2,0)--(7.4,0) node[below,font=\small] {energy};
  \draw[thick] (.6,-.12)--(.6,.12) node[below=5pt,font=\small] {};
  \draw[thick] (6.6,-.12)--(6.6,.12);
  \draw[decorate,decoration={brace,amplitude=4pt,mirror}] (.6,-.45)--(6.6,-.45)
    node[midway,below=4pt,font=\small] {bulk interval $I$};
  \foreach \x in {.9,1.55,2.05,2.9,3.95,4.7,5.25,6.2}
    \fill (\x,0) circle (2pt);
  \fill[red!15] (3.55,-.25) rectangle (4.35,.25);
  \fill (3.95,0) circle (2pt);
  \draw[red!70!black,thick] (3.55,-.25)--(3.55,.25) (4.35,-.25)--(4.35,.25);
  \node[font=\small,red!70!black] at (3.95,.55) {$|x-E|\le\varepsilon_0$};
  \draw[->,thin] (3.95,1.45) node[above,font=\small,align=center]
    {at most one exterior\\ eigenvalue (Step~1)} -- (3.95,.85);
 \end{scope}
 \begin{scope}[shift={(10.4,0)}]
  \draw[->] (-1.9,0)--(2.1,0) node[below,font=\small] {$x$};
  \draw[->] (0,-2.7)--(0,2.8);
  \draw[dashed,black!60,domain=.38:1.9,samples=60] plot (\x,{1/\x});
  \draw[dashed,black!60,domain=-1.9:-.38,samples=60] plot (\x,{1/\x});
  \draw[very thick,blue!70!black,domain=.6:1.9,samples=60] plot (\x,{1/\x});
  \draw[very thick,blue!70!black,domain=-1.9:-.6,samples=60] plot (\x,{1/\x});
  \draw[very thick,blue!70!black] (-.6,-1.667)--(.6,1.667);
  \draw[thin] (.6,-.08)--(.6,.08) node[below=4pt,font=\small] {$\varepsilon_0$};
  \node[font=\small,blue!70!black,anchor=east] at (-.35,1.0) {$g_{\varepsilon_0}$};
  \node[font=\small,black!60,anchor=west] at (.45,2.55) {$x^{-1}$};
 \end{scope}
\end{tikzpicture}
\caption{The truncation in Steps~1--2, drawn schematically. Left: the dots
are the eigenvalues of the exterior $H_R$. With overwhelming probability,
every window of half-length $\varepsilon_0=W^{-A_0}$ around a bulk energy
(shaded) contains at most one of them. Right: the truncated function
$g_{\varepsilon_0}$ (solid) agrees with $x^{-1}$ (dashed) outside
$[-\varepsilon_0,\varepsilon_0]$ and is linear inside. It is bounded by
$\varepsilon_0^{-1}$ and Lipschitz with constant $\varepsilon_0^{-2}$.
Replacing $x^{-1}$ by $g_{\varepsilon_0}$ in the exterior resolvent changes
only the contribution of the single exterior eigenvalue in the window, and
hence changes the deformation $F(E)$ by rank at most one.}
\label{fig:plw-truncation}
\end{figure}
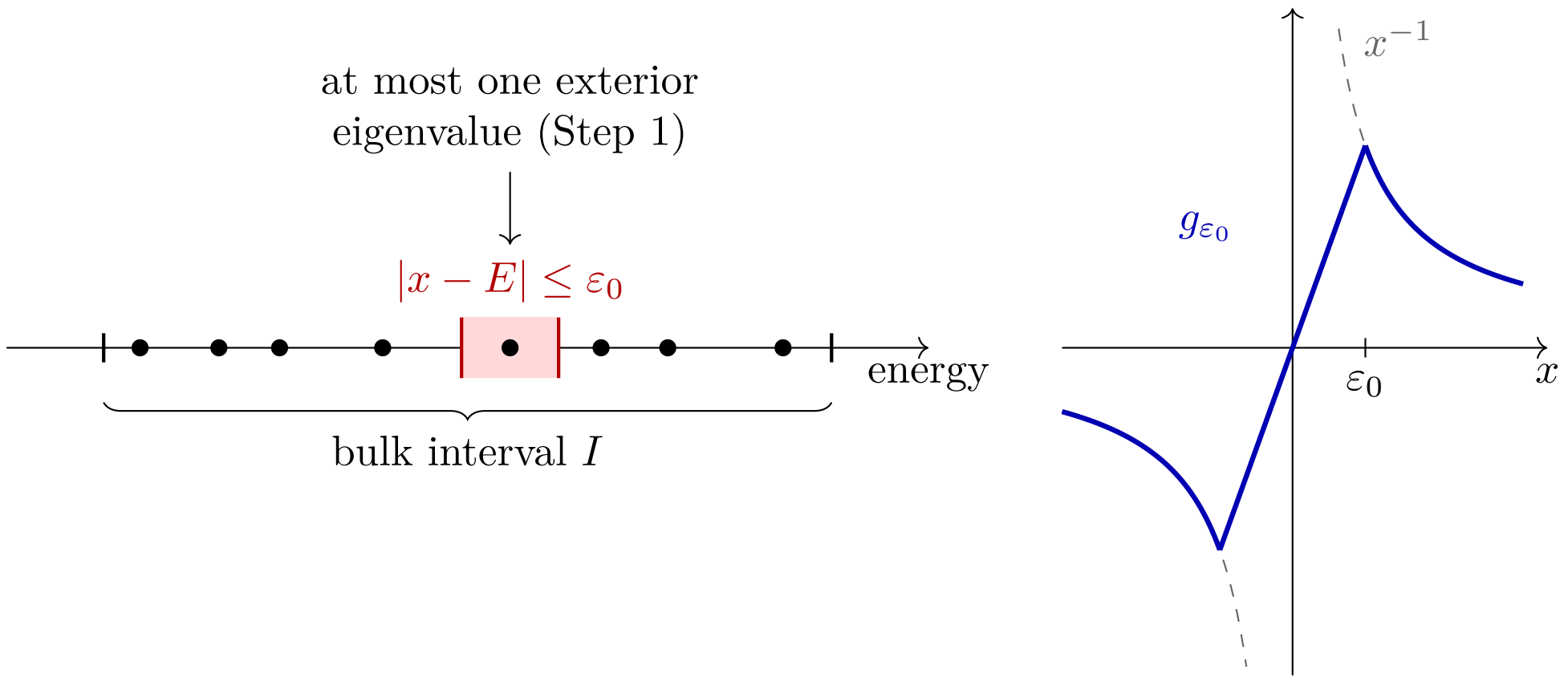

\emph{Step 1: at most one exterior eigenvalue near any bulk energy.} We
first bound the expected number of exterior eigenvalues in an interval
$J$. Conditional on all off-diagonal entries of $H_R$, its diagonal entries
are independent and have densities bounded by $C_b\sqrt W$. Spectral
averaging in one diagonal coordinate therefore gives
\[
 \E N_{H_R}(J)\le C_bN_R\sqrt W|J|,\qquad N_R=(n-2)W .
\]
To exclude two eigenvalues in one window we need the second factorial
moment. Let $H_R^{(a)}$ be the minor with coordinate $a$ removed. Cauchy
interlacing gives $N_{H_R^{(a)}}(J)\ge N_{H_R}(J)-1$. We multiply this
inequality by $(\one_J(H_R))_{aa}$ and sum over $a$. Since
$\sum_a(\one_J(H_R))_{aa}=N_{H_R}(J)$, this gives
\[
 N_{H_R}(J)\bigl(N_{H_R}(J)-1\bigr)
 \le\sum_a(\one_J(H_R))_{aa}\,N_{H_R^{(a)}}(J).
\]
The minor count does not depend on the removed diagonal variable. So we
may first average $(\one_J(H_R))_{aa}$ over that variable, using the
one-coordinate estimate, and then average the minor count, using the
first-moment bound. This proves
\begin{equation}\label{eq:exterior-crowding}
 \E [N_{H_R}(J)(N_{H_R}(J)-1)]
       \le C_b N_R(N_R-1)W|J|^2 .
\end{equation}

Now set $\varepsilon_0=W^{-A_0}$. We cover $I$ by $O(\varepsilon_0^{-1})$
intervals of length $4\varepsilon_0$, so that every window
$[E-\varepsilon_0,E+\varepsilon_0]$ with $E\in I$ lies in one of them.
Markov's inequality for \eqref{eq:exterior-crowding} and a union bound,
with $N_R\le TW^2$, give
\begin{equation}\label{eq:plw-one-pole}
 \Prob\Bigl\{\sup_{E\in I}N_{H_R}([E-\varepsilon_0,E+\varepsilon_0])
                          \ge2\Bigr\}\le CW^{5-A_0}.
\end{equation}
We choose $A_0>D+6$, so this probability is smaller than the required
failure power.

\emph{Step 2: a truncated deformation.} We now remove the possible single
pole by truncating the function $x\mapsto x^{-1}$ at distance
$\varepsilon_0$, as in Figure~\ref{fig:plw-truncation}. We put
\[
 g_{\varepsilon_0}(x)=
 \begin{cases}x/\varepsilon_0^2,&|x|\le\varepsilon_0,\\
 x^{-1},&|x|>\varepsilon_0,\end{cases}
 \qquad \widehat F(E)=V-C^*g_{\varepsilon_0}(H_R-E)C .
\]
The function $g_{\varepsilon_0}$ is continuous, bounded by
$\varepsilon_0^{-1}$, and Lipschitz with constant $\varepsilon_0^{-2}$.
On the complement of the event in \eqref{eq:plw-one-pole}, it differs from
$x^{-1}$ at the spectrum of $H_R-E$ in at most one eigenvalue. Hence
\[
 \rank\bigl(F(E)-\widehat F(E)\bigr)\le1
\]
at every nonpole energy $E\in I$. On the Gaussian norm event for $V$ and
$C$, the bounds on $g_{\varepsilon_0}$ give
\begin{equation}\label{eq:plw-truncated-norm}
 \norm{\widehat F(E)}\le CW^{A_0},\qquad
 \norm{\widehat F(E)-\widehat F(e)}
                         \le CW^{2A_0}|E-e| .
\end{equation}
So the truncated deformation has polynomial size and is Lipschitz in the
energy, which is what the regularity certificate and the local law need.

\emph{Step 3: a uniform regularity certificate.} We show that the
truncated deformation passes the scalar test of
Lemma~\ref{lem:regularity-certificate}, simultaneously for all $E\in I$.
At a fixed positive imaginary height, a rank-one change of a Hermitian
matrix changes its resolvent by a rank-one matrix of trace norm at most
twice the inverse height. So the normalized traces of the resolvents with
$F(E)$ and with $\widehat F(E)$ differ by $O(W^{-1})$. For the true residual
pair, Theorem~\ref{thm:matched-traces} gives these traces at each
deterministic energy, with every fixed moment of size $W^{-1}$. We apply
it on an energy grid of spacing $W^{-M}$ with $M>2A_0+6$. High fixed
moments and a union bound control the grid. Between grid points we use
the Lipschitz bound \eqref{eq:plw-truncated-norm}. Together this gives
\begin{equation}\label{eq:uniform-pole-certificate}
 \max_{a=i,i+1}\left|
 W^{-1}\Tr P_a(\widehat F(E)-E-sm_{\rm sc}(E))^{-1}
                      -m_{\rm sc}(E)\right|
       \le CW^{-1+\epsilon/8}
\end{equation}
simultaneously for $E\in I$, with an arbitrary fixed polynomial failure
power. Lemma~\ref{lem:regularity-certificate} now gives uniform
regularity of the truncated deformation on a fixed energy and variance
rectangle.

\begin{figure}[ht]
\centering
\begin{tikzpicture}[x=1.25cm,y=1cm]
 \draw[->] (-3.3,0)--(3.5,0) node[below,font=\small] {$x$};
 \draw[->] (0,-.1)--(0,3.1);
 \fill[blue!10] (-1.2,0) rectangle (1.2,2.95);
 \draw[very thick,domain=-3.2:3.2,samples=160]
   plot (\x,{2.8*0.09/(\x*\x+0.09)});
 \draw[thick,dashed,red!70!black,domain=-3.2:3.2,samples=160]
   plot (\x,{2.8*0.3*1.2/(\x*\x+1.44)});
 \draw[thin] (1.2,-.08)--(1.2,.08) node[below=4pt,font=\small] {$L_0\eta$};
 \draw[thin] (-1.2,-.08)--(-1.2,.08) node[below=4pt,font=\small] {$-L_0\eta$};
 \node[font=\small,anchor=west] at (.35,2.75) {height $\eta$};
 \node[font=\small,anchor=west,red!70!black] at (1.55,.95) {height $L_0\eta$};
 \node[font=\small,align=center,blue!70!black] at (0,-1.05)
   {eigenvalues here carry at least half of the trace at height $\eta$};
\end{tikzpicture}
\caption{The counting argument of Step~4, drawn schematically. The solid
and dashed curves are scalar Poisson kernels at heights $\eta$ and
$L_0\eta$, normalized as in the normalized imaginary trace. Outside the
shaded region $|x|\le L_0\eta$, the ratio of the first kernel to the
second is at most $2/L_0$. The local law bounds the trace at height
$L_0\eta$ above, so for $L_0$ large the eigenvalues outside the shaded
region carry less than half of the trace at height $\eta$. The local law
also bounds that trace below, so many eigenvalues lie in the shaded region.}
\label{fig:plw-two-heights}
\end{figure}
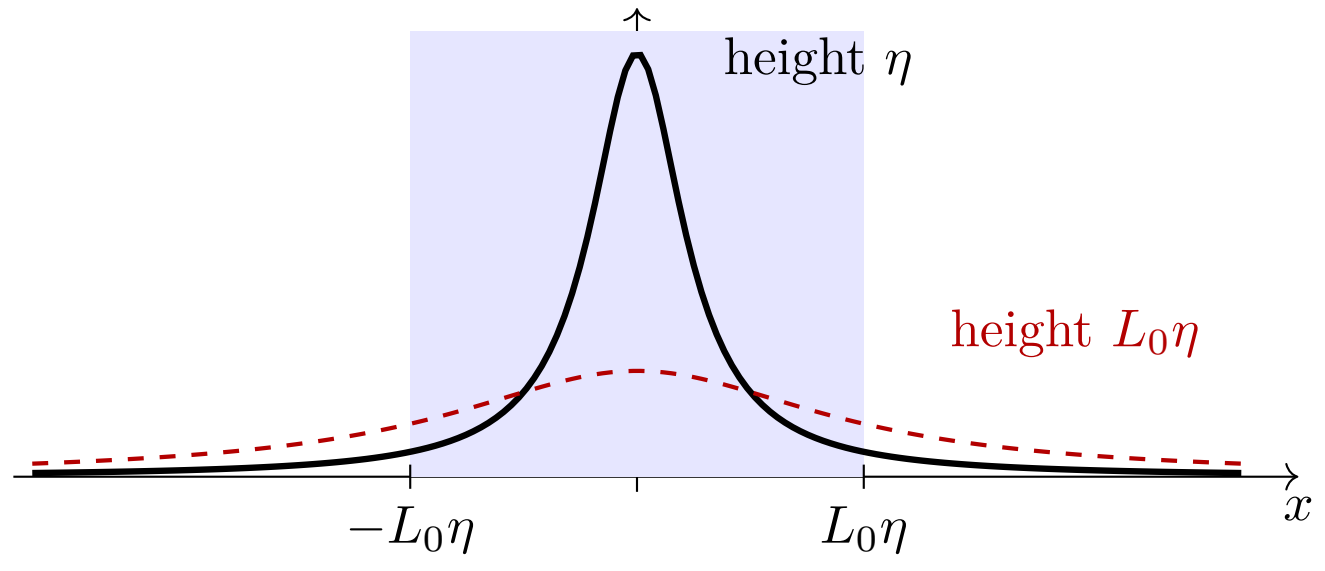

\emph{Step 4: the counting bounds of part~(i).} We apply the local law,
Theorem~\ref{thm:external-local}, to the truncated matrix
$\widehat F(E)+G_d(t)-E$. We use it at the two heights $\eta=W^{-1+\gamma}$
and $L_0\eta$, where $L_0$ is a large fixed constant. Since
$0<\gamma<\gamma_0$, both heights tend to zero within the fixed regularity
rectangle of Step~3, and $W\eta=W^\gamma\to\infty$. Hence the normalized
imaginary traces at the two heights are at least $c_0$ and at most $C_0$,
respectively.

The argument is the one in Figure~\ref{fig:plw-two-heights}. For
$|x|\ge L_0\eta$, the ratio of the scalar Poisson kernels at heights $\eta$
and $L_0\eta$ is at most $2/L_0$. We choose $L_0$ with $2C_0/L_0<c_0/2$.
Then the eigenvalues with $|x|\ge L_0\eta$ contribute less than $c_0d/2$ to
the unnormalized trace at height $\eta$. The remaining contribution, at
least $c_0d/2$, comes from $|x|\le L_0\eta$, where each eigenvalue
contributes at most $\eta^{-1}$. So there are at least $c_0d\eta/2$
eigenvalues there. This gives the lower bound in
\eqref{eq:physical-local-count}. The local law at height $\eta$ also gives
the upper bound in \eqref{eq:physical-local-count}, since the Poisson sum
is $\eta$ times the imaginary trace. Its error is $O(W^\varepsilon)$, which
is smaller than $W\eta$ when $\varepsilon<\gamma$.

These bounds hold on the energy grid of Step~3, and the Lipschitz bound
\eqref{eq:plw-truncated-norm} transfers them to all $E\in I$. Finally we
return from $\widehat F(E)$ to $F(E)$. The two Schur complements differ by
rank at most one. By rank-one interlacing, the count in an interval changes
by at most two, and the Poisson sum changes by $O(1)$. Both changes are
negligible against $W\eta$. This proves part~(i).

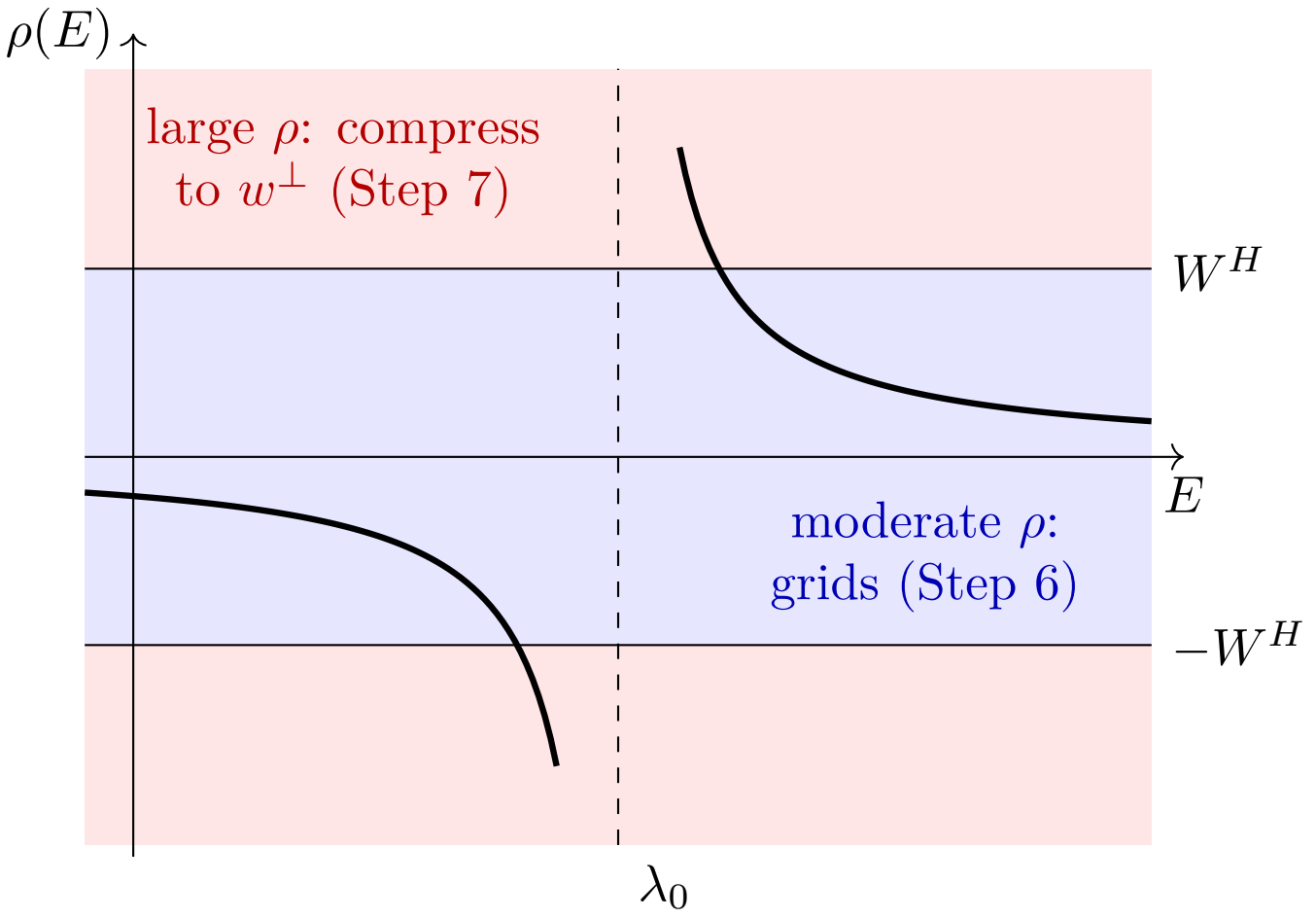
\begin{figure}[ht]
\centering
\begin{tikzpicture}[x=1.1cm,y=.8cm]
 \fill[blue!10] (-.3,-1.6) rectangle (6.3,1.6);
 \fill[red!10] (-.3,1.6) rectangle (6.3,3.3);
 \fill[red!10] (-.3,-3.3) rectangle (6.3,-1.6);
 \draw[->] (-.3,0)--(6.5,0) node[below,font=\small] {$E$};
 \draw[->] (0,-3.4)--(0,3.6) node[left,font=\small] {$\rho(E)$};
 \draw[dashed] (3,-3.3)--(3,3.3);
 \node[font=\small,below right] at (3,-3.3) {$\lambda_0$};
 \draw[very thick,domain=-.3:2.62,samples=80] plot (\x,{-1/(3-\x)});
 \draw[very thick,domain=3.38:6.3,samples=80] plot (\x,{1/(\x-3)});
 \draw[thin] (-.3,1.6)--(6.3,1.6) node[right,font=\small] {$W^H$};
 \draw[thin] (-.3,-1.6)--(6.3,-1.6) node[right,font=\small] {$-W^H$};
 \node[font=\small,align=center,blue!70!black] at (4.9,-.85)
   {moderate $\rho$:\\ grids (Step~6)};
 \node[font=\small,align=center,red!70!black] at (1.3,2.5)
   {large $\rho$: compress\\ to $w^\perp$ (Step~7)};
\end{tikzpicture}
\caption{The two regimes of Steps~5--7, drawn schematically. Near the
single exterior pole $\lambda_0$, the deformation is $F_{\rm reg}(E)+\rho(E)ww^*$
with a coefficient $\rho(E)$ that blows up at $\lambda_0$ (solid curves).
While $|\rho|\le W^H$ (blue band), all matrices have polynomial size and we
use the eigenvector estimates on grids in $E$ and $\rho$. When
$|\rho|>W^H$ (red), the matrix almost decouples along $w$, and its
eigenvectors near the energy are close to those of its compression to
$w^\perp$.}
\label{fig:plw-pole}
\end{figure}

\emph{Step 5: the single pole as a rank-one term.} We now turn to the
eigenvector estimate~(ii), where truncation is not enough: we need the
eigenvectors of the actual Schur complement. On the complement of the
event in \eqref{eq:plw-one-pole}, a bulk energy has at most one exterior
eigenvalue $\lambda_0$ within distance $\varepsilon_0$. Near it, with
$\psi_{\lambda_0}$ the corresponding eigenvector of $H_R$, the exact
deformation is
\[
 F(E)=F_{\rm reg}(E)+\rho(E)ww^*,\qquad
 \rho(E)=-\frac{\norm{C^*\psi_{\lambda_0}}^2}{\lambda_0-E},
 \qquad w=\frac{C^*\psi_{\lambda_0}}{\norm{C^*\psi_{\lambda_0}}}.
\]
If the numerator vanishes, there is no pole. Otherwise the unit vector $w$
is fixed by the exterior. The regular part $F_{\rm reg}(E)$ has polynomial
norm and Lipschitz constant on this pole interval, exactly as in
\eqref{eq:plw-truncated-norm}. Moreover, for every real $\rho$, the matrix
$F_{\rm reg}(E)+\rho ww^*$ differs from $\widehat F(E)$ by rank at most one.
Hence the certificate \eqref{eq:uniform-pole-certificate} and its
regularity consequence hold for these matrices too, uniformly in $\rho$,
including arbitrarily large $\rho$. We now treat $\rho$ as a free real
parameter and distinguish the two regimes of Figure~\ref{fig:plw-pole}.

\emph{Step 6: moderate coefficients.} First let $|\rho|\le W^H$ for a fixed
$H$, chosen in Step~7. Then every matrix involved has polynomial norm. We
condition on the exterior and use polynomial grids in $E$, in $\rho$ and in
the window center. At each grid point, Theorem~\ref{thm:external-eth}, on a
window of radius $W^{-1+\beta}$, bounds every entry
$u_\lambda^*(B-j_B)u_\mu$ of the eigenvectors in the window by
$W^{-1/2+\epsilon'}$. Here $\epsilon'$ is its stochastic-domination loss,
which we take smaller than $\beta$. Let $P$ be the spectral projection onto
the window. Then the scalar local law gives $\rank P\le CW^\beta$. Bounding
the operator norm by the rank times the largest entry gives
\[
 \rank P\le CW^\beta,\qquad
 \norm{P(B-j_B)P}\le CW^{-1/2+2\beta}.
\]
The variation of the deterministic center $j_B$ across the window costs
$O(W^{-1+\beta})$, which is smaller. We choose $\beta<\epsilon/8$. We
increase the failure exponent so that it pays for all grids and for all
$O(nW)$ possible pole labels. This uses only a union bound, so no
independence between grid events is needed.

It remains to pass from grid points to every actual eigenvector. Let $r_W$
be the operator distance from the actual matrix to the nearest grid
matrix. An eigenvector of the actual matrix with eigenvalue $\lambda$ is
then an approximate eigenvector of the grid matrix, with error $r_W$. By
the spectral theorem, its component outside the grid window has norm at
most $2r_WW^{1-\beta}$. We choose the grids so fine that this is
$O(W^{-10})$. This transfers the operator bound to every actual
eigenvector, and proves \eqref{eq:physical-window-eth} in this regime
without any assumption on eigenvalue gaps.

\emph{Step 7: large coefficients and compression.} Now let $|\rho|>W^H$. The
matrix then almost decouples along $w$, so we compare it with its
compression. Let $Q_w=I-ww^*$ and compress to $w^\perp$. The compressed
fresh noise is exactly $\GUE_{d-1}(t)$, with total variance
$s'=s(1-1/d)$. Indeed, complete $w$ to an orthonormal basis $(w,U_w)$ of
$\C^d$ fixed by the exterior. As $G\sim\GUE_d(t)$ is independent of the
exterior, by unitary invariance its matrix in this basis is again
$\GUE_d(t)$ conditionally on the exterior, and $U_w^*GU_w$, its principal
block of size $d-1$, is $\GUE_{d-1}(t)$. In the normalization
$\GUE_\Naux(s/\Naux)$ of \eqref{eq:mde}, used in
Theorems~\ref{thm:external-local} and~\ref{thm:external-eth}, the
uncompressed and compressed noises thus have parameters $s$ in dimension
$d$ and $s'$ in dimension $d-1$, with $t=s/d=s'/(d-1)$. 
The fixed-height certificate and each observable center
change only by $O(W^{-1})$. Indeed, the embedded difference of resolvents
has rank at most two and bounded norm, while $|s-s'|=O(W^{-1})$. So the
estimates of Steps~3--6 apply to the compressed deformation.

To compare the eigenvectors, we write the full matrix in the decomposition
$\C w\oplus w^\perp$ as
\[
 \begin{pmatrix}\rho+a&z^*\\z&B_0\end{pmatrix},
 \qquad |a|+\norm z+\norm{B_0}\le K=CW^{A_0}.
\]
Let $(x,y)$ be a normalized eigenvector with eigenvalue $\lambda$ and
$|\lambda|\le2\delta$. The first row of the eigenvalue equation gives
$(\rho+a-\lambda)x=-z^*y$, and the second gives $(B_0-\lambda)y=-zx$. Since
$|\rho|$ is much larger than $K$, this yields
\[
 |x|\le 2K/|\rho|,\qquad
 \norm{(B_0-\lambda)y}\le2K^2/|\rho| .
\]
We choose $H$ so large that both errors are $O(W^{-10})$, even after the
division by the window radius in Step~6. Then $y$ is an approximate
eigenvector of the compressed matrix, and the compressed window estimate
applies to it.

For the smoothed traces, Schur inversion bounds the trace norm of the
difference from the embedded compressed resolvent by
$C|\rho|^{-1}(1+K/\eta)^2$, once $|\rho|\ge2(K+K^2/\eta+1)$. This also tends
to zero uniformly. Smooth bounded spectral functions have the stated
compression limit, by their resolvent representation with polynomial
parameter bounds. This proves \eqref{eq:physical-window-eth} in the regime
of large coefficients. In particular, the assertion for spectral functions
of $F(E)-E$ uses the limiting compression at an infinite rank-one
parameter, as stated in the lemma.

\emph{Step 8: the selected energy and the pair law.} So far all estimates
hold simultaneously for every $E\in I$ that is not an exterior eigenvalue,
outside an event of the stated original probability. Actual roots avoid
exterior poles almost surely, by the polynomial argument in the proof of
Claim~3 of Theorem~\ref{thm:coarea}. Theorem~\ref{thm:rank-location} pays
for the event that the selected energy leaves the bulk interval. Since the
estimates are simultaneous in $E$, they hold at the selected energy.
Finally, the pair law has density
$d\Prob^\sharp/d\Prob=\tfrac n2m_S\le n/2\le CW$. So increasing the
original failure exponent by one transfers every exceptional event to the
pair law. This completes the proof of Lemma~\ref{lem:physical-local-window}.
\end{proof}

\section{The complete determinant score}
\label{sec:determinant-score}
The density of the eigenvector direction in Theorem~\ref{thm:coarea}
contains a determinant factor. To apply Gaussian integration by parts on the
sphere, we need the logarithmic gradient of this factor, the score, to be
small on average. This section proves that bound. It is used in the angular
comparison of Part~I.
Theorem~\ref{thm:determinant-score} gives an exact formula for the score
and bounds its moments of order $r<2$ under the pair law. Three propositions
control its pieces. Proposition~\ref{prop:physical-real-loop} controls the
real inverse trace through a loop identity.
Proposition~\ref{prop:negative-spread} gives a variance inequality from the
log-concavity of the spectral law. Proposition~\ref{prop:physical-minor-tail}
controls the tail of the remaining minor.
We first justify integration by parts across the singular inertia
boundary (Lemma~\ref{lem:t26-sector-ibp}). In Subsection~\ref{subsec:s10-whole-row-inverse-moments} we prove
inverse moments of whole Gaussian rows (Lemmas~\ref{lem:gaussian-polynomial-gradient}
and~\ref{lem:physical-column}). In Subsection~\ref{subsec:s10-a-real-loop-with-a-positive-source} we control the real
loop with a positive source. In Subsection~\ref{subsec:s10-convexity-of-the-complete-spectral-law} we prove convexity of the
complete spectral law (Lemma~\ref{lem:indexed-spectral-convexity}). In
Subsection~\ref{subsec:s10-restoring-the-deleted-row} we integrate the deleted row while keeping its signs, and
complete the proof.

Keep $d=2W$, $t=s/d$, $Q=I-vv^*$ and $K_0=F(E)-E$.
On $v^\perp$ put
\[
 K=QK_0Q|_{v^\perp},\quad M=K+G_{d-1}(t),\quad R=M^{-1},
 \quad U=\Tr R^2.
\]
Embed $R$ by zero on $v$. For the complete outer data
$\theta=(E,v,H_R,C,V)$ and the prescribed sector
$\ell=\alpha-N_-(H_R-E)-1$, write
\begin{equation}\label{eq:indexed-complement}
 \mathfrak Z_\ell(v)=\E_G[\det(M)^2\one_{N_-(M)=\ell}],
 \qquad
 \av{f}_\theta=
 \frac{\E_G[f(M)\det(M)^2\one_{N_-(M)=\ell}]}{\mathfrak Z_\ell(v)},\quad \overline R=\av R_\theta.
\end{equation}
Only attainable sectors are used. All Gaussian priors in this
definition have mass one.

\begin{theorem}\label{thm:determinant-score}
At every fixed exterior and real nonpole energy,
\begin{equation}\label{eq:exact-determinant-score}
 \nabla_v\log\mathfrak Z_\ell(v)=-4\overline RQK_0v,
\end{equation}
where the horizontal real gradient is defined by
$Df[h]=\Re(h^*\nabla f)$ for $h\perp v$.
For every fixed $0<r<2$,
\begin{equation}\label{eq:physical-determinant-score}
 \E^\sharp_{\alpha,S}\one_{\{E\in I\}}
 \left[\norm{\overline R}^{\,r}
       +\norm{\nabla_v\log\mathfrak Z_\ell(v)}^{\,r}\right]
 \le C_{b,T,\vartheta,I,r}.
\end{equation}
The estimate is uniform for $4\le n\le TW$, both parities,
adjacent pairs and deterministic bulk ranks.
\end{theorem}

No lower bound on a conditional sector probability is required.
In particular a bounded angular vector field of size $W^{-1}$ has
determinant-score contribution $O(W^{-1})$ under the physical law.  The change of variables of
Theorem~\ref{thm:coarea} produces the weight
$1_{\{N_-(M)=\ell\}}\det(M)^2$. Here $N_-(M)$ is the number of negative
eigenvalues of $M$, and the indicator fixes the rank of the eigenvalue we
follow. The normalizing constant of this weight is the sector mass
\[
 Z_\ell(K)=\E\bigl[1_{\{N_-(M)=\ell\}}P^2\bigr].
\]
Later we integrate by parts in the Gaussian coordinates of $G_m(t)$, and
we differentiate $Z_\ell(K)$ in parameters of $K$. Both operations are
delicate, because the indicator is discontinuous. Following lemma resolves this issue.

\begin{lemma}[First derivatives across the singular sector boundary]
\label{lem:t26-sector-ibp}
Fix finite $m$, $t>0$, a Hermitian $K$ and an attainable inertia sector
$\ell$. Put $M=K+G_m(t)$, $P=\det M$ and $C=\operatorname{adj}M$. In every
real Gaussian coordinate, both $1_{\{N_-(M)=\ell\}}P^2$ and
$1_{\{N_-(M)=\ell\}}PC_{ab}$ have weak derivatives equal to their
ordinary derivatives off $\{P=0\}$, with no singular boundary measure.
These functions and their first derivatives are Gaussian integrable.
Consequently, first Gaussian integration by parts with the insertions
$1$ or $R_{ab}=(M^{-1})_{ab}$, and first differentiation of the complete
sector mass in any smooth finite-dimensional parameter of $K$, are valid
at zero smoothing height.
\end{lemma}

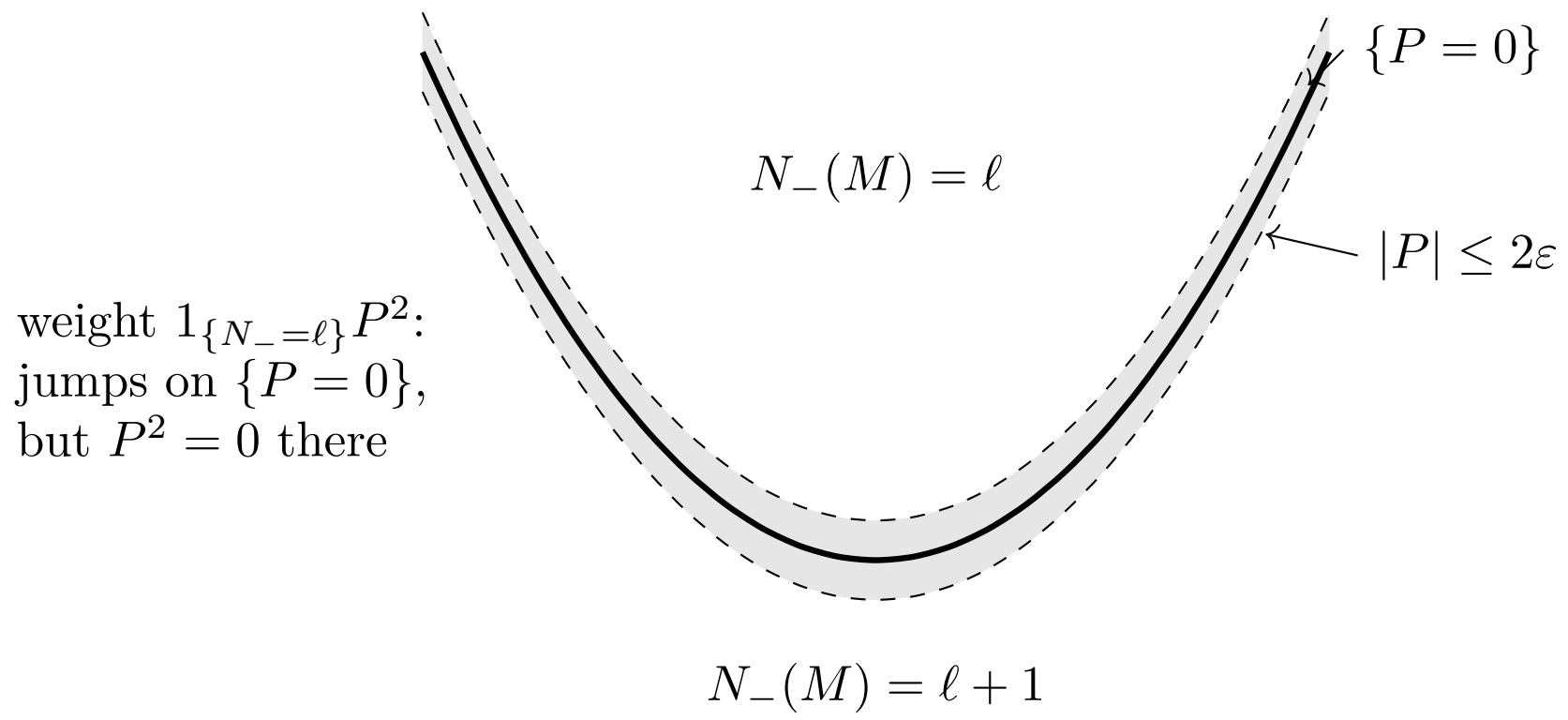
\begin{figure}[ht]
\centering
\begin{tikzpicture}[x=1cm,y=1cm]
 \fill[black!10] plot[smooth,domain=-3.2:3.2] (\x,{.35*\x*\x-1.1+.28})
   -- plot[smooth,domain=3.2:-3.2] (\x,{.35*\x*\x-1.1-.28}) -- cycle;
 \draw[very thick] plot[smooth,domain=-3.2:3.2] (\x,{.35*\x*\x-1.1});
 \draw[dashed] plot[smooth,domain=-3.2:3.2] (\x,{.35*\x*\x-1.1+.28});
 \draw[dashed] plot[smooth,domain=-3.2:3.2] (\x,{.35*\x*\x-1.1-.28});
 \node[font=\small] at (0,1.6) {$N_-(M)=\ell$};
 \node[font=\small] at (0,-2.0) {$N_-(M)=\ell+1$};
 \node[font=\small,anchor=west] at (3.3,2.5) {$\{P=0\}$};
 \draw[->,thin] (3.3,2.5)--(3.05,2.25);
 \node[font=\small,anchor=west] at (3.4,1.05) {$|P|\le2\varepsilon$};
 \draw[->,thin] (3.4,1.05)--(2.75,1.2);
 \node[font=\small,anchor=west,align=left] at (-6.2,.2)
   {weight $1_{\{N_-=\ell\}}P^2$:\\ jumps on $\{P=0\}$,\\ but $P^2=0$ there};
\end{tikzpicture}
\caption{The singular sector boundary in the proof of
Lemma~\ref{lem:t26-sector-ibp}, drawn schematically in two of the real
Gaussian coordinates. The hypersurface $\{P=0\}$ (solid) is where an
eigenvalue of $M$ crosses zero, and it separates two inertia sectors. The
indicator $1_{\{N_-(M)=\ell\}}$ jumps across it, but the factors $P^2$
and $PC_{ab}$ vanish there. The proof first inserts a smooth cutoff
$\chi(P/\varepsilon)$, which vanishes on $\{|P|\le\varepsilon\}$.
Differentiating it produces terms supported in
$\{\varepsilon\le|P|\le2\varepsilon\}$, inside the shaded band
$\{|P|\le2\varepsilon\}$ between the dashed curves. These terms vanish
as $\varepsilon\to0$, because the band shrinks to a set of measure zero.}
\label{fig:sib-sector}
\end{figure}

\begin{proof}
The weight that we differentiate is the squared
determinant multiplied by an indicator that fixes the number of negative
eigenvalues. The indicator jumps where an eigenvalue crosses zero, which is
exactly where the determinant vanishes, and the idea is that the squared
determinant kills the jump, so that no extra boundary term appears. To make
this precise we cut out a thin neighbourhood of the zero set of the
determinant and a far region (Step~1), rewrite the inverse insertion so that
the function becomes a polynomial (Step~2), and show that the terms produced
by the cutoff vanish as the neighbourhood shrinks (Step~3).
Figure~\ref{fig:sib-sector} shows this neighbourhood. Step~4 treats
differentiation in a parameter of the deterministic matrix, and Step~5
explains why no lower bound on the mass of the sector is needed.

\emph{Step 1:} We first remove the singular surface and
infinity. Choose an even smooth function $\chi$ that vanishes on $[-1,1]$
and equals one outside $[-2,2]$. We insert the factor $\chi(P/\varepsilon)$,
which vanishes on the band $\{|P|\le\varepsilon\}$ of
Figure~\ref{fig:sib-sector}, and an ordinary cutoff of radius $B$ in the
Gaussian coordinates. Away from $\{P=0\}$ the inertia of $M$ is locally
constant, because no eigenvalue can cross zero there. So on the support
of $\chi(P/\varepsilon)$ the indicator is locally constant, and every
derivative is an ordinary one.

\emph{Step 2:} For the insertion $R_{ab}$ we
first remove the pole. Since $R=C/P$, we have the polynomial identity
\[
 P^2R_{ab}=PC_{ab},
\]
and we use it before differentiating. So the function to be
differentiated is $1_{\{N_-(M)=\ell\}}PC_{ab}$, a polynomial times an
indicator, and not a function with a pole.

\emph{Step 3:} No we show the cutoff terms vanish. Differentiating the cutoff
$\chi(P/\varepsilon)$ produces the extra factor
$\chi'(P/\varepsilon)\,\partial P/\varepsilon$. It is supported where
$\varepsilon\le|P|\le2\varepsilon$, and there $|P|/\varepsilon\le2$. Hence
the extra terms for the two functions are bounded by
\[
 C|P\,\partial P|1_{\{|P|\le2\varepsilon\}},\qquad
 C|C_{ab}\partial P|1_{\{|P|\le2\varepsilon\}},
\]
respectively. The determinant is a nonzero polynomial in all the
Gaussian coordinates, so its zero set has Lebesgue measure zero. As
$\varepsilon\to0$, both bounds therefore tend to zero almost everywhere.
They also have a polynomial majorant that does not depend on
$\varepsilon$, and polynomials are integrable against the Gaussian
density. The ordinary derivatives are polynomials restricted to the
sector as well. By dominated convergence we first remove $\varepsilon$,
and then the Gaussian tails remove the cutoff $B$. The same argument
applies to every real and imaginary coordinate of the matrix.

Differentiating an inverse produces products of inverse entries, which
could in principle create a new pole. They do not, because
\[
 P^2R_{ab}R_{cd}=C_{ab}C_{cd}
\]
is again a polynomial. This proves the integration by parts with the
insertions $1$ and $R_{ab}$.

\emph{Step 4:} We now justify derivative in a paramaeter of $K$. Now let $K$ depend on a
smooth finite-dimensional parameter. On a compact neighbourhood of the
parameter, the cutoff argument of Steps~1--3 applies, with uniformly
bounded polynomial coefficients. There is also a direct argument.
Locally $P$ is Lipschitz and vanishes on $\{P=0\}$, so
$|P|\le C\operatorname{dist}(M,\{P=0\})$ and
\[
 P^2=O\bigl(\operatorname{dist}(M,\{P=0\})^2\bigr).
\]
Hence $1_{\{N_-=\ell\}}P^2$ is continuously differentiable across
$\{P=0\}$, with zero first derivative there. By Jacobi's formula,
$\partial P=P\Tr(R\,\partial M)$, and $\partial M=\partial K$. It follows
that
\[
 \partial Z_\ell(K)
 =2\E\bigl[1_{\{N_-=\ell\}}P^2\Tr(R\,\partial K)\bigr].
\]

\emph{Step 5: normalization and the physical law.} All these identities
are proved before we divide by $Z_\ell(K)$. To normalize, we only need
$Z_\ell(K)>0$, which holds because the sector is attainable. We do not
need a uniform lower bound on $Z_\ell(K)$. After normalization,
Lemma~\ref{lem:physical-column}, with any fixed exponent between two and
three, gives the integrability of the first and second inverse
insertions under the physical law. So we may condition the identities on
the outer data and then integrate them against the full pair law. The
polynomial argument above is not used to infer any uniform conditional
inverse moment, or any lower bound on small sectors.
\end{proof}

 \subsection{Whole-row inverse moments}\label{subsec:s10-whole-row-inverse-moments}

  We need integrability before using real inverse matrices in the
weighted density. Lemma~\ref{lem:gaussian-polynomial-gradient} gives a
dimension-independent bound for a fixed-degree polynomial, which
Lemma~\ref{lem:physical-column} applies to a complete Gaussian row.
We prove the polynomial bound by Gaussian rotation and then insert
the row variables without discarding their determinant weight.

\begin{lemma}\label{lem:gaussian-polynomial-gradient}
For a nonzero real polynomial $P$ of degree at most $k$ on a
standard Gaussian space of arbitrary finite dimension, and
$2\le p<3$,
\begin{equation}\label{eq:gaussian-polynomial-gradient}
 \E\norm{\nabla P}^{p}|P|^{2-p}\le C_{k,p}\E P^2.
\end{equation}
The constant is independent of dimension and coefficients.
\end{lemma}

The degree is fixed while the Gaussian dimension may grow. This uniformity is what makes the estimate useful for a complete matrix row.

\begin{proof}
For a trigonometric polynomial $f$ of degree at most $k$, put
$Q(z)=z^kf(\theta)$, $z=e^{i\theta}$. Its degree is at most $2k$.
Factoring its roots gives
$|f'|/|f|\le k+\sum_\zeta|e^{i\theta}-\zeta|^{-1}$.
For $Q=(z-\zeta)Q_\zeta$, the coefficient recurrence, read backward
if $|\zeta|\le1$ and forward otherwise, gives
$\norm{Q_\zeta}_2\le2k\norm Q_2$ and
$\sup_{|z|=1}|Q_\zeta(z)|\le(2k)^{3/2}\norm Q_2$.
Thus, with normalized circle measure,
\[
 \int |f'|^p|f|^{2-p}
 \le C_{k,p}\norm f_2^2
       \left(1+\sum_\zeta
                  \int|e^{i\theta}-\zeta|^{-(p-2)}\,d\theta\right)
 \le C_{k,p}\int|f|^2 .
\]
The last integral is bounded uniformly in $\zeta$ because $p-2<1$.
Let $X,Y$ be independent standard Gaussian vectors and apply this
to $f(\theta)=P(X\cos\theta+Y\sin\theta)$. For each $\theta$,
the derivative is $\nabla P(X_\theta)\cdot Y_\theta$, with
$X_\theta,Y_\theta$ independent standard Gaussians. Averaging gives
$\E|N(0,1)|^p$ times the left side of
\eqref{eq:gaussian-polynomial-gradient}; the right side averages
to $C_{k,p}\E P^2$. Tonelli justifies averaging before finiteness
is established.
\end{proof}

\begin{lemma}\label{lem:physical-column}
For every outer-measurable unit vector $w\in v^\perp$ and
$2\le p<3$,
\begin{equation}\label{eq:physical-column}
 \E^\sharp_{\alpha,S}\norm{Rw}^{p}\le C_pW^{p/2},
 \qquad
 \E^\sharp_{\alpha,S}\norm R_{\HS}^{p}\le C_pW^p.
\end{equation}
These estimates hold at every rank and every actual energy.
\end{lemma}

The row is integrated with its restored physical weight. These moments justify the real loop and bound exceptional spectral configurations.

\begin{proof}
Condition on the minor obtained by deleting $w$. The determinant
$P=\det M$ is a real polynomial of degree at most two in the
entire remaining Gaussian row. In whitened real coordinates,
\[
 \norm{\nabla\log|P|}^2
     =t|R_{ww}|^2+2t\sum_{j\ne w}|R_{jw}|^2
     \ge t\norm{Rw}^2.
\]
Apply Lemma~\ref{lem:gaussian-polynomial-gradient}.
Drop the selected inertia only in the nonnegative numerator.
A whole Hermitian row variation has at most one positive and one
negative eigenvalue. On an attainable minor sector the restored
row therefore belongs to at most the three neighboring full sectors.
In \eqref{eq:complete-coarea} the regression, energy, vector and
pair weight $n/2$ are fixed on this row. Integrating the unchanged
base restores three pair laws, each of mass at most one.
This proves the first inequality. Sum over an orthonormal basis
and use
$\norm R_{\HS}^p\le(d-1)^{p/2-1}\sum_a\norm{Re_a}^p$
to prove the second.
\end{proof}

\subsection{A real loop with a positive source}\label{subsec:s10-a-real-loop-with-a-positive-source}

  We control the real inverse trace of the selected matrix, where
individual small eigenvalues can make an absolute bound too costly.
Proposition~\ref{prop:physical-real-loop} gives the needed integrated
estimate. We first compare a centered real trace with its regularized
version, use the local spectral estimates, and then apply the Gaussian
loop identity under the physical pair law.

\begin{lemma}[Centered trace unsmoothing]
\label{lem:t26-score-unsmoothing}
Let $A=K_0+G$ be Hermitian of dimension $d$, let $R=A^+$ be its
Moore--Penrose inverse and $R_\eta=(A-i\eta)^{-1}$, with
$0<\eta\le\delta<1$. Assume $\|G\|\le C_0$, $|g(0)|\le C_0$, and that $g$
is real and $C_0$-Lipschitz on $[-\delta,\delta]$. Suppose that for all
nonzero eigenvalues $\lambda$ in that interval, with eigenvectors
$u_\lambda$,
\[
 |u_\lambda^*K_0u_\lambda-g(\lambda)|\le a_0,
 \qquad
 \sum_\lambda\frac{\eta^2}{\lambda^2+\eta^2}\le C_0d\eta.
\]
Then, with a constant depending only on $C_0$ and $\delta$,
\begin{equation}\label{eq:t26-score-unsmoothing}
 |\Tr(K_0-g(0))(R-\Re R_\eta)|
 \le C a_0\sqrt{d\eta}\,\|R\|_{\rm HS}+Cd\eta+Cd\eta^2.
\end{equation}
\end{lemma}

\begin{figure}[ht]
\centering
\begin{tikzpicture}[x=1.15cm,y=.75cm]
 \draw[->] (-4.2,0)--(4.4,0) node[below] {$\lambda$};
 \draw[->] (0,-3.4)--(0,3.6);
 \draw[thick,domain=0.3:4,samples=80] plot (\x,{1/\x});
 \draw[thick,domain=-4:-0.3,samples=80] plot (\x,{1/\x});
 \draw[thick,dashed,blue!70!black,domain=-4:4,samples=160]
   plot (\x,{\x/(\x*\x+0.36)});
 \fill[red!70!black,opacity=.12] (-1.1,-3.4) rectangle (1.1,3.4);
 \draw[<->] (-.6,-.9)--(.6,-.9) node[midway,below,font=\small] {$2\eta$};
 \node[font=\small,anchor=west] at (.4,3.1) {$1/\lambda$};
 \node[font=\small,anchor=west,blue!70!black] at (2.1,1.05) {$\lambda/(\lambda^2+\eta^2)$};
 \node[font=\small,red!70!black,anchor=west] at (1.15,-2.6) {$|\lambda|\lesssim\eta$};
\end{tikzpicture}
\caption{The two scalar functions compared in the proof of
Lemma~\ref{lem:t26-score-unsmoothing}, drawn schematically. On an
eigenvalue $\lambda$ of $A$, the unsmoothed inverse $R$ acts as
$1/\lambda$ (solid), and $\Re R_\eta$ acts as $\lambda/(\lambda^2+\eta^2)$
(dashed), which is bounded by $1/(2\eta)$. Their difference is
$\eta^2/(\lambda(\lambda^2+\eta^2))$. For $|\lambda|\gg\eta$ it is of order
$\eta^2/|\lambda|^3$, so the two curves almost coincide. It is large only
in the shaded region $|\lambda|\lesssim\eta$. The proof therefore needs
small diagonal coefficients only for eigenvalues near zero, and these are
provided by the hypothesis $u_\lambda^*K_0u_\lambda\approx g(\lambda)\approx
g(0)$.}
\label{fig:su-kernels}
\end{figure}
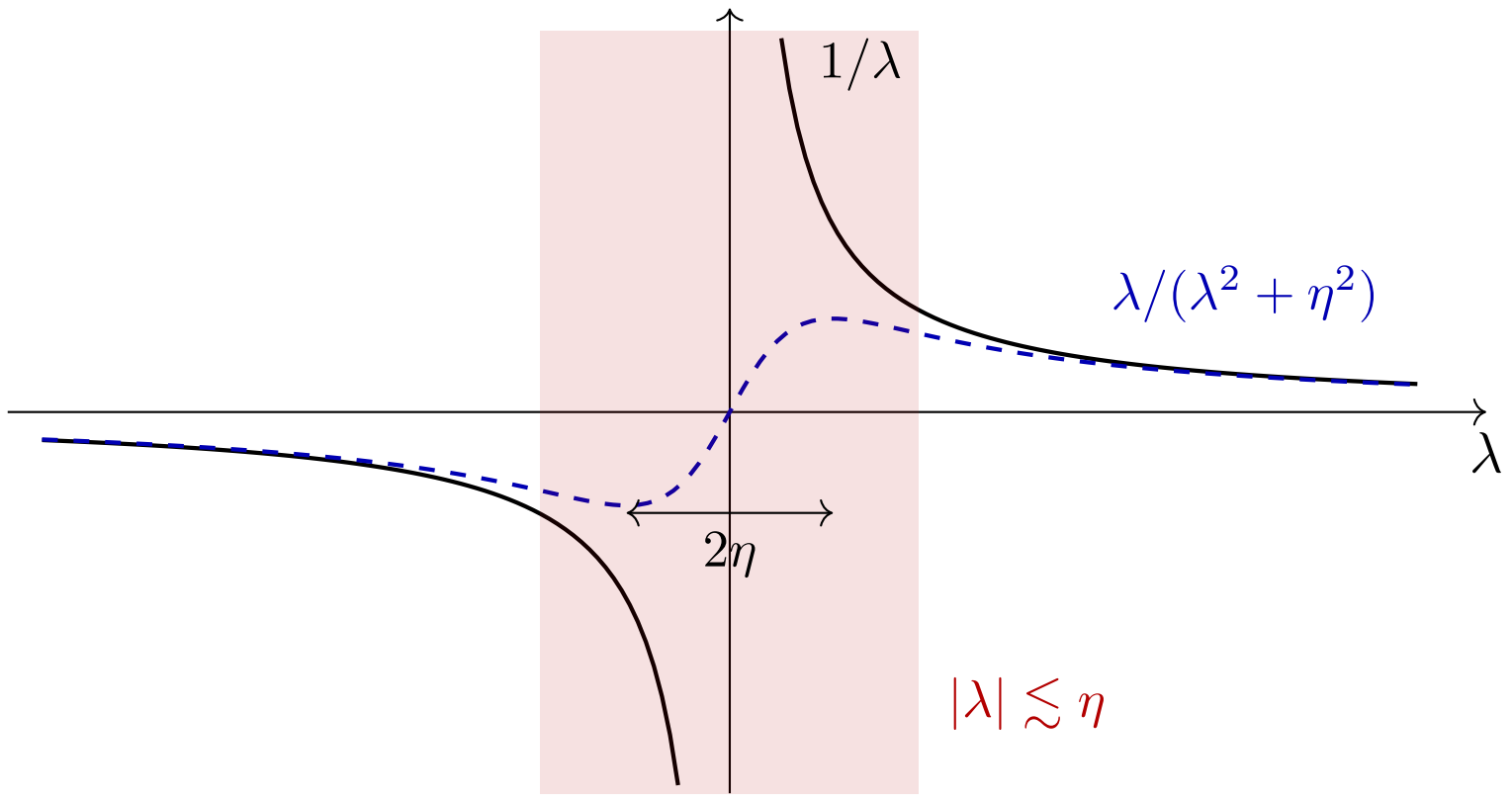

\begin{proof}
We compare the exact inverse with a regularized
inverse, the inverse of the matrix shifted by a small imaginary multiple of
the identity, which is never singular. In the eigenbasis both inverses are
diagonal, so the difference splits into one term per eigenvalue (Step~1). As
Figure~\ref{fig:su-kernels} shows, the two inverses differ substantially only
for eigenvalues close to zero. For these eigenvalues we use the hypothesis
that the corresponding diagonal coefficients are close to a smooth function
(Step~2); eigenvalues away from zero are handled by the eigenvalue equation
itself (Step~3). %

\emph{Step 1: reduction to eigenvalues.} We work in the eigenbasis
$(u_\lambda)$ of $A$, in which both $R$ and $\Re R_\eta$ are diagonal.
The trace in \eqref{eq:t26-score-unsmoothing} is therefore the sum over
eigenvalues of the diagonal coefficient
$u_\lambda^*K_0u_\lambda-g(0)$ times the scalar difference of the two
inverses. For a nonzero eigenvalue this difference is
\[
 \lambda^{-1}-\Re(\lambda-i\eta)^{-1}
 =\lambda^{-1}-\frac{\lambda}{\lambda^2+\eta^2}
 =\frac{\eta^2}{\lambda(\lambda^2+\eta^2)}.
\]
At a zero eigenvalue both real operators in the difference vanish, since
the Moore--Penrose inverse is zero on the kernel and
$\Re(-i\eta)^{-1}=0$. So zero eigenvalues do not contribute.

\emph{Step 2: eigenvalues near zero.} For $0<|\lambda|\le\delta$ we split
the centered diagonal coefficient as
\[
 u_\lambda^*K_0u_\lambda-g(0)
 =\bigl[u_\lambda^*K_0u_\lambda-g(\lambda)\bigr]+\bigl[g(\lambda)-g(0)\bigr].
\]
The first bracket is at most $a_0$ by hypothesis. By the Cauchy--Schwarz
inequality, its total contribution is at most
\[
 a_0\left(\sum_{\lambda\ne0}\lambda^{-2}\right)^{1/2}
       \left(\sum_\lambda\frac{\eta^4}{(\lambda^2+\eta^2)^2}
       \right)^{1/2}
 \le C a_0\|R\|_{\rm HS}\sqrt{d\eta}.
\]
Here the first factor equals $\|R\|_{\rm HS}$, because the eigenvalues of
$R$ are the inverses of the nonzero $\lambda$. For the second factor we
use $\eta^4/(\lambda^2+\eta^2)^2\le\eta^2/(\lambda^2+\eta^2)$ and the
hypothesis $\sum_\lambda\eta^2/(\lambda^2+\eta^2)\le C_0d\eta$.

The second bracket is the Lipschitz part. It is at most $C_0|\lambda|$,
and this factor $|\lambda|$ cancels the $\lambda$ in the denominator of
the scalar difference. Its contribution is therefore at most
\[
 C\sum_\lambda\frac{\eta^2}{\lambda^2+\eta^2}\le Cd\eta.
\]

\emph{Step 3: eigenvalues away from zero.} For $|\lambda|>\delta$ we have
no information on $u_\lambda^*K_0u_\lambda$ from the hypothesis. We use
the eigenvalue equation instead. Since $K_0=A-G$, we have
$u_\lambda^*K_0u_\lambda=\lambda-u_\lambda^*Gu_\lambda$, and so
\[
 |u_\lambda^*K_0u_\lambda-g(0)|\le|\lambda|+2C_0 .
\]
On the other hand, the scalar difference is at most $\eta^2/|\lambda|^3$.
Their product is at most $C_\delta\eta^2$ for each of the at most $d$
eigenvalues, so these contribute at most $C_\delta d\eta^2$. Adding the
three contributions proves \eqref{eq:t26-score-unsmoothing}. At no point
did we use $\|K_0\|$.

Steps 1--3 prove the stated deterministic lemma. The physical-law application belongs to Proposition~\ref{prop:physical-real-loop}, whose proof supplies the inverse-column and exceptional-event estimates.

\end{proof}

 \begin{proposition}\label{prop:physical-real-loop}
There is an outer-measurable error $e_\theta$ such that
\begin{align}
 2t\av U_\theta
 -t\av{(\Tr R+dE/2)^2}_\theta
      &=sd(1-E^2/4)-1+e_\theta,
                                      \label{eq:physical-real-loop}\\
 \E^\sharp_{\alpha,S}\one_{\{E\in I\}}|e_\theta|
      &\le CW^{49/50}.                 \label{eq:real-loop-error}
\end{align}
Consequently, for a fixed $c_b>0$,
\begin{equation}\label{eq:spread-good-event}
 \mathcal G_\alpha=\{E\in I,\ \av U_\theta\ge c_bW^2\},
 \qquad
 \Prob^\sharp_{\alpha,S}(\{E\in I\}\setminus\mathcal G_\alpha)
                  \le CW^{-1/50}.
\end{equation}
\end{proposition}

The square in \eqref{eq:physical-real-loop} has a useful sign.
The order-$W$ source on its right side must be retained; setting it
to zero would destroy the lower spread bound.

\begin{proof}
We first establish the physical trace estimate
\begin{equation}\label{eq:physical-weighted-trace}
 \E^\sharp\one_{\{E\in I\}}
 \left|\Tr G_dR-sE\Tr R-(sd-1)\right|\le CW^{49/50}.
\end{equation}
Put $\eta=W^{-3/5}$ and $L_0=W^{1/4}$.
This height corresponds to $\gamma=2/5<\gamma_0$ in
Lemma~\ref{lem:physical-local-window}; its eigenvector part also allows
the independently chosen sufficiently small window exponent $\beta$.
 For a fixed smooth cutoff let
$K_{0,L_0}=K_0\chi(K_0/L_0)$.
Apply Lemma~\ref{lem:physical-local-window} to $K_{0,L_0}/L_0$.
If $A u_\lambda=\lambda u_\lambda$ and $|\lambda|\le\delta$,
then $\norm{K_0u_\lambda}\le|\lambda|+\norm{G_d}\le C$
on the Gaussian norm event. Hence
\[
 |u_\lambda^*(K_0-K_{0,L_0})u_\lambda|\le C/L_0.
\]
The same tail bound holds for the covariance trace because
$\Gamma_E(\lambda)$ is a squared resolvent with a shift in a
fixed compact upper-half-plane set. Taking imaginary parts of
$\Tr K_0(K_0-\lambda-sm_E)^{-1}
 =d+(\lambda+sm_E)dm_E$ gives the exact center
\[
 g(\lambda)=d^{-1}\Tr K_0\Gamma_E(\lambda)
      =\lambda+2s\Re m_E(\lambda),\quad
 |g'(\lambda)|\le C,\quad
 g(0)=-sE+O(W^{-1+\varepsilon}).
\]
Therefore, simultaneously at the relevant energies and eigenvalues,
\begin{equation}\label{eq:unbounded-observable-eth}
 |u_\lambda^*K_0u_\lambda-g(\lambda)|
       \le C(L_0W^{-1/2+\varepsilon}+L_0^{-1})
       \le CW^{-1/4+\varepsilon}.
\end{equation}

At deterministic exterior data let
$R_\eta=(K_0+G_d-i\eta)^{-1}$. Under the unweighted fresh prior,
Gaussian integration by parts gives
\[
 \E_G\Tr G_dR_\eta=-t\E_G(\Tr R_\eta)^2.
\]
The averaged local law gives
$\norm{\Tr R_\eta-dm_E(i\eta)}_{L^p}\le C_{p,\varepsilon}
 W^\varepsilon/\eta$. Gaussian $L^p$ Poincar\'e gives
\[
 \norm{\Tr G_dR_\eta-\E_G\Tr G_dR_\eta}_{L^p}
 \le C_p\sqrt t\,
 \norm{\norm{R_\eta-R_\eta G_dR_\eta}_{\HS}}_{L^p}
 \le C_p\eta^{-3/2}.
\]
Indeed $\norm{R_\eta}_{\HS}^2
=\eta^{-1}\Im\Tr R_\eta\le Cd/\eta$ on the local-law event;
the deterministic inverse bound and arbitrarily high failure powers
pay its complement. Thus, on the energy and rank-one grids from
the preceding section,
\[
 \Tr G_dR_\eta=-sd\,m_E(i\eta)^2
                      +O(W^{9/10+\varepsilon}).
\]
The infinite-parameter compression used there also applies, since
$G_d$ has bounded norm and the resolvent difference has vanishing
trace norm. Fixed high moments pay the finite grids. These are
simultaneous prior estimates, so one additional failure power transfers
them to the pair law.

The identity $K_0R_\eta=I+i\eta R_\eta-G_dR_\eta$ now gives
\[
 \Tr(K_0-g(0))\Re R_\eta
    =d(1-s|m_E(0)|^2)
                +O(W^{9/10+2\varepsilon}+d\eta).
\]
Here regularity permits replacing $m_E(i\eta)$ by $m_E(0)$,
and the certificate gives
$m_E(0)=m_{\rm sc}(E)+O(W^{-1+\varepsilon})$.
 Apply Lemma~\ref{lem:t26-score-unsmoothing} with
 \eqref{eq:unbounded-observable-eth} and the Poisson bound:
\[
 |\Tr(K_0-g(0))(R-\Re R_\eta)|
 \le CW^{-1/4+\varepsilon}\sqrt{W\eta}\norm R_{\HS}
                                  +CW\eta+CW\eta^2.
\]
For $|\lambda|>\delta$ this follows from
$|u_\lambda^*K_0u_\lambda|\le|\lambda|+\norm{G_d}$;
the zero eigenvalue contributes zero because $R$ is embedded by zero.
Lemma~\ref{lem:physical-column} bounds the first term's expectation
by $CW^{19/20+\varepsilon}$. Replacing $g(0)$ by $-sE$ costs
at most $CW^{1/2+\varepsilon}$, since
$\E^\sharp|\Tr R|\le\sqrt d\,\E^\sharp\norm R_{\HS}\le CW^{3/2}$.
On all exceptional events use
$\Tr K_0R=d-1-\Tr G_dR$ to remove the unbounded deformation.
Column moments with any fixed exponent above two, Gaussian norm
moments and H\"older pay their arbitrarily small polynomial probability.
As $|m_{\rm sc}(E)|=1$, sufficiently small fixed $\varepsilon$
proves \eqref{eq:physical-weighted-trace}.

On each complete sector in \eqref{eq:indexed-complement}, Gaussian
integration by parts gives exactly
\begin{equation}\label{eq:indexed-gaussian-loop}
 \av{\Tr G_dR}_\theta
      =2t\av{\Tr R^2}_\theta-t\av{(\Tr R)^2}_\theta.
\end{equation}
To check the contractions, $\partial R=-R(\partial G)R$
contributes $-t(\Tr R)^2$; differentiating $\det(M)^2$
contributes $2t\Tr R^2$. Lemma~\ref{lem:t26-sector-ibp} proves the zero-height integration
by parts with its actual determinant-square sector and inverse
insertion; its cutoff derivatives tend to zero before normalization.
Condition \eqref{eq:physical-weighted-trace} on $\theta$, then
complete the square using $s=td$ to obtain
\eqref{eq:physical-real-loop}--\eqref{eq:real-loop-error}.
The square is nonnegative and $\inf_{E\in I}(1-E^2/4)>0$.
Markov's inequality proves \eqref{eq:spread-good-event}.
\end{proof}

\subsection{Convexity of the complete spectral law}\label{subsec:s10-convexity-of-the-complete-spectral-law}

  We use the shape of the complete eigenvalue density to control its
fluctuations. Lemma~\ref{lem:indexed-spectral-convexity} proves a weighted
variance inequality in each fixed-inertia sector, for arbitrary
deformation. We compute the curvature of the ordered density, apply the
resulting variance bound, and use it to control the negative-spread
term.

\begin{lemma}\label{lem:indexed-spectral-convexity}
Under the normalized measure
$Z_\ell(K)^{-1}\det(M)^2\one_{N_-(M)=\ell}\,d\Prob(K+G_m(t))$,
the ordered eigenvalues have a negative logarithmic density
$\mathcal U$ on the convex fixed-inertia chamber satisfying
\begin{equation}\label{eq:spectral-curvature}
 \nabla^2\mathcal U\succeq t^{-1}I+
 2\diag(\lambda_j^{-2})+
 \sum_{i<j}\frac{(e_i-e_j)(e_i-e_j)^T}{(\lambda_i-\lambda_j)^2}.
\end{equation}
For every smooth real $f$ of finite weighted energy,
\begin{equation}\label{eq:weighted-spectral-variance}
 \Var f\le
 \E\sum_j\frac{|\partial_jf|^2}{t^{-1}+2\lambda_j^{-2}}.
\end{equation}
\end{lemma}

The fixed-inertia chamber is retained in the spectral density. Convexity therefore applies without a lower bound on the probability of that rank sector.

\begin{proof}
Assume first $k_1<\cdots<k_m$ are the eigenvalues of $K$.
Let $\mathcal C=\{\lambda_1<\cdots<\lambda_m\}$ and
$\Delta(\lambda)=\prod_{i<j}(\lambda_j-\lambda_i)$.
The eigenvalue transition density of Hermitian Brownian motion
is the Vandermonde transform of independent Brownian motions killed
at the boundary of $\mathcal C$:
\[
 p_t^{\rm D}(k,\lambda)=
 \frac{\Delta(\lambda)}{\Delta(k)}p_t^{\mathcal C}(k,\lambda).
\]
See also 
\cite[Section 3.3, (3.10) and Remark 1]{KT}.
This follows by antisymmetrizing the independent Gaussian
heat kernel and taking its Vandermonde transform, or by checking its
initial condition and the generator
$\frac{1}{2}\sum_j\partial_j^2+
\sum_j\sum_{i\ne j}(k_j-k_i)^{-1}\partial_j$. Here, we provide details of the first argument. By unitary invariance of
$G_m(t)$ take $K=\diag(k)$. Weyl's integration formula gives the
ordered eigenvalues of $K+G_m(t)$ a density proportional to
$\Delta(\lambda)^2\int_{U(m)}e^{-\norm{U\Lambda U^*-K}_{\HS}^2/(2t)}\,dU$,
$\Lambda=\diag(\lambda)$, and since
$\norm{U\Lambda U^*-K}_{\HS}^2=|\lambda|^2+|k|^2-2\Tr(U\Lambda U^*K)$, the
Harish-Chandra--Itzykson--Zuber integral
\[
 \int_{U(m)}e^{\Tr(U\Lambda U^*K)/t}\,dU
 =\prod_{p=1}^{m-1}p!\;
 \frac{t^{m(m-1)/2}\det[e^{\lambda_ik_j/t}]_{i,j=1}^m}{\Delta(\lambda)\Delta(k)}
\]
turns it into the external-source density
$\Delta(\lambda)\Delta(k)^{-1}\det[\varphi_t(\lambda_i-k_j)]$ of Step~1 of
the proof of Theorem~\ref{thm:external-sine} (with $t/\Naux$ there
replaced by $t$), $\varphi_t$ the centered Gaussian density of variance
$t$; the constant is one because the heat flow preserves the harmonic
polynomial $\Delta$. Finally
$\det[\varphi_t(\lambda_i-k_j)]=p_t^{\mathcal C}(k,\lambda)$ is the
Karlin--McGregor formula. In the expansion of the determinant, exchanging
the future paths of the two coordinates that meet at the first collision
is a measure-preserving involution between the colliding paths of the
terms $\sigma$ and $\sigma\circ(ij)$, which carry opposite signs; the
remaining paths stay in $\mathcal C$ and, as $k,\lambda\in\mathcal C$,
occur only for $\sigma={\rm id}$.
The killed kernel factors into a free Gaussian kernel and the
probability $S_t(k,\lambda)$ that a bridge stays in $\mathcal C$.
Consequently the exact density on
$\mathcal C_\ell=\mathcal C\cap\{\lambda_\ell<0<\lambda_{\ell+1}\}$
is
\begin{equation}\label{eq:complete-spectral-density}
 \frac{1}{Z_\ell(K)}
 \prod_j\lambda_j^2\,
 \frac{\Delta(\lambda)}{\Delta(k)}
 (2\pi t)^{-m/2}e^{-|\lambda-k|^2/(2t)}
 S_t(k,\lambda).
\end{equation}
The same $Z_\ell(K)$ as in the matrix integral appears here.

On a finite time mesh the centered bridge is Gaussian, and
membership of its affine endpoint interpolation in the convex
chamber is a convex constraint jointly in endpoints and bridge
coordinates. Pr\'ekopa's marginal theorem~\cite{Prekopa} makes the probability
log-concave in the endpoints. Decreasing dense meshes give $S_t$;
a Brownian bridge touches a chamber wall without crossing it with
probability zero. Thus $-\log S_t$ is convex. Taking the Hessian
of minus the logarithm in \eqref{eq:complete-spectral-density}
gives \eqref{eq:spectral-curvature}.

For completeness, on a smooth convex bounded exhaustion solve
$-\mathscr Lu=f-\E f$, where
$\mathscr L=\Delta-\nabla\mathcal U\cdot\nabla$, with Neumann
boundary condition. Integration by parts and the Bochner identity
give
\[
 \Var f=\E\nabla f\cdot\nabla u,\qquad
 \E\nabla u^T(\nabla^2\mathcal U)\nabla u
                         \le\E(\mathscr Lu)^2=\Var f.
\]
The boundary term is nonnegative by convexity. Cauchy--Schwarz
in the Hessian metric proves the Brascamp--Lieb inequality.
Smooth convex approximations and exhaustion pass it to the chamber,
and dropping positive Hessian terms gives
\eqref{eq:weighted-spectral-variance}.
Weighted Sobolev approximation permits every finite-energy test. Here we explain how.  Finite weighted energy $\mathcal E(f)$, the right side of
\eqref{eq:weighted-spectral-variance}, forces $f\in L^2$: for smooth
$\varphi_N$ with $|\varphi_N'|\le1$, $|\varphi_N(x)|\le\min(|x|,N)$ and
$\varphi_N(x)=x$ for $|x|\le N-1$, the bounded tests $\varphi_N\circ f$
have variance at most $\mathcal E(f)$, so if $\Prob(|f|\le A)\ge1/2$,
Chebyshev gives $|\E\varphi_N(f)|\le A+2\mathcal E(f)^{1/2}$ for
$N>A+1$ and $\E\varphi_N(f)^2$ is bounded in $N$. Fatou's lemma gives
$\E f^2<\infty$, and dominated convergence
$\Var f=\lim_N\Var\varphi_N(f)\le\mathcal E(f)$. 

For repeated eigenvalues of $K$ we use the confluent form. Since $\det[e^{\lambda_ik_j/t}]$ is antisymmetric in $k$,
the function
$\widetilde S_t(k,\lambda)=e^{-\lambda\cdot k/t}\det[e^{\lambda_ik_j/t}]/\Delta(k)$,
which equals $S_t(k,\lambda)/\Delta(k)$ for distinct $k$, extends to an
entire function of $(\lambda,k)$, and both sides of the
Harish-Chandra--Itzykson--Zuber identity are continuous in $k$. Hence
\eqref{eq:complete-spectral-density} holds for every $K$, with
$S_t(k,\lambda)/\Delta(k)$ read as $\widetilde S_t(k,\lambda)$. By Weyl's
formula the ordered eigenvalue density of $K+G_m(t)$ is strictly positive
on $\mathcal C$, so $\widetilde S_t(k,\cdot)$ is smooth and strictly
positive there. For distinct $k^{(n)}\to k$,
$\widetilde S_t(k^{(n)},\cdot)\to\widetilde S_t(k,\cdot)$ locally uniformly
on $\mathcal C$, so $-\log\widetilde S_t(k,\cdot)$ is a pointwise limit of
the convex functions $-\log S_t(k^{(n)},\cdot)+\log\Delta(k^{(n)})$, and it
is convex. The Hessian computation above then gives
\eqref{eq:spectral-curvature} and \eqref{eq:weighted-spectral-variance}
for every $K$.
\end{proof}

\begin{proposition}\label{prop:negative-spread}
For $U=\Tr M^{-2}$ under this same normalized law,
\begin{equation}\label{eq:negative-spread}
 \E U\,\E U^{-1}\le64e^8,\qquad
 \E U^{-p}\le C_p(\E U)^{-p}\quad(p>0).
\end{equation}
All constants are independent of $K,m,t,\ell$.
\end{proposition}

The inverse spread is controlled by the mean spread under the same normalized law. Constants do not depend on the deformation or on the selected inertia sector.

\begin{proof}
Dropping the Gaussian term in \eqref{eq:weighted-spectral-variance}
and using $\sum_j\lambda_j^{-4}\le U^2$ gives
\[
 \Var\sqrt U\le\tfrac12\E U,\qquad
 \Var U^{-1/2}\le\tfrac12\E U^{-1},\qquad
 \Var\log U\le2.
\]
These tests are legitimate: the determinant square makes $U$
integrable, $U^{-1}\le\Tr M^2/m^2$, and $\log U$ has finite
second moment at both zero and infinity. Put $A=\E U$.
Paley--Zygmund applied to $\sqrt U$ gives
$\Prob(U\ge A/8)\ge1/8$. Chebyshev for $\log U$ then implies
$A\le8e^4e^{\E\log U}$. Repeating this with $U^{-1}$ gives
$\E U^{-1}\le8e^4e^{-\E\log U}$, proving the first bound.

On the fixed sign chamber put $x_j=|\lambda_j|$ and
$f(x)=(\sum_jx_j^{-2})^{-1/2}$. Direct differentiation gives
\[
 D^2f[a,a]=3\Big(\sum_jx_j^{-2}\Big)^{-5/2}
 \left[\Big(\sum_ja_jx_j^{-3}\Big)^2
       -\Big(\sum_jx_j^{-2}\Big)\sum_ja_j^2x_j^{-4}\right]\le0.
\]
The sign change is linear on this chamber; its density remains
log-concave. Pr\'ekopa applied to the convex hypograph of $f$
makes $\Prob(f\ge r)$ log-concave for $r\ge0$.
If $a=\E f$, Markov gives $\Prob(f\ge2a)\le1/2$, whence
$\Prob(f\ge r)\le2\exp[-(\log2)r/(2a)]$.
As $a^2\le\E U^{-1}\le64e^8/A$, integration gives
$\E f^{2p}\le C_pA^{-p}$, which is the second bound.
\end{proof}

\subsection{Restoring the deleted row}\label{subsec:s10-restoring-the-deleted-row}

 We finish the determinant estimate by integrating the missing Gaussian
row while retaining its signs. Lemma~\ref{lem:signed-row-spread} combines
the potentially singular inverse terms before averaging, and
Proposition~\ref{prop:physical-minor-tail} controls the remaining tail.
We compute the row integral exactly, use the physical local estimates,
and then complete Theorem~\ref{thm:determinant-score}.

 Choose a measurable unit eigenvector $w(\theta)$ of $K$. In its
basis the fresh row has zero off-diagonal mean:
\[
 M=\begin{pmatrix}\mu+x&y^*\\y&B\end{pmatrix},\quad
 x\sim N(0,t),\quad y\sim{\rm CN}(0,tI).
\]
For fixed invertible $B$ write
$C_B=B^{-1}$, $V_B=\Tr C_B^2$,
$a=\mu+x-y^*C_By$, $\delta_B=\mu-t\Tr C_B$ and
$q_B^2=\delta_B^2+t+t^2V_B=\E_{x,y}a^2$.

\begin{lemma}\label{lem:signed-row-spread}
With $U=\Tr M^{-2}$,
\begin{align}
 \E_{x,y}a^2U
 &=V_Bq_B^2+(1+tV_B)^2
             +2t\delta_B\Tr C_B^3-t^2\Tr C_B^4 \notag\\
 &=2V_Bq_B^2+1+tV_B
             -\Tr[C_B^2(\delta_BI-tC_B)^2]
 \le(2V_B+t^{-1})q_B^2.                 \label{eq:signed-row-spread}
\end{align}
\end{lemma}

The cubic and quartic inverse terms are combined before integration
over $B$. Their separate absolute integrals need not be finite.

\begin{proof}
Block inversion gives
$a^2U=a^2V_B+2a\,y^*C_B^3y+(1+y^*C_B^2y)^2$.
The complex Gaussian contractions are
\[
 \E[a\,y^*C_B^3y]
       =t\delta_B\Tr C_B^3-t^2\Tr C_B^4,\qquad
 \E(1+y^*C_B^2y)^2=(1+tV_B)^2+t^2\Tr C_B^4.
\]
Summing and completing the displayed matrix square proves the
identities. Its trace is nonnegative and
$1+tV_B\le q_B^2/t$, proving the inequality.
\end{proof}

\begin{proposition}\label{prop:physical-minor-tail}
For every fixed $p>0$ and $A\ge0$,
\begin{equation}\label{eq:physical-minor-tail}
 \Prob^\sharp_\alpha(\mathcal G_\alpha,V_B\le A)
        \le C_p\left(\frac{A+W}{W^2}\right)^{p/(p+1)}.
\end{equation}
For every $r<2$,
\begin{equation}\label{eq:physical-inverse-spread}
 \E^\sharp_\alpha\one_{\{E\in I\}}q_B^{-r}\le C_r.
\end{equation}
\end{proposition}

The estimate refers to the row-restored physical measure. Treating the smaller matrix under its own unweighted law would change the required tail.

\begin{proof}
At fixed outer data and minor, the selected row sector can occur
only if $N_-(B)\in\{\ell-1,\ell\}$. On $\{V_B\le A\}$,
drop the selected sign restriction in the nonnegative numerator
$\det(M)^2U$ and apply \eqref{eq:signed-row-spread}.
Restore $\det(B)^2q_B^2$, the complete row mass, on these minor
sectors. Its two possible Schur signs give full sectors among
$\ell-1,\ell,\ell+1$. The unchanged pair coarea base integrates
their masses to at most three. Therefore
\[
 \E^\sharp_\alpha\one_{\mathcal G_\alpha,V_B\le A}U
                                      \le3(2A+t^{-1}).
\]
Proposition~\ref{prop:negative-spread}, integrated on
$\mathcal G_\alpha$, gives
$\E^\sharp_\alpha\one_{\mathcal G_\alpha}U^{-p}\le C_pW^{-2p}$.
H\"older with exponents $(p+1)/p,p+1$ proves
\eqref{eq:physical-minor-tail}. The outer event was dropped after
restoration; it was not reassigned to a neighboring rank.

Let $Y=q_B^{-1}\le t^{-1/2}\le C_b\sqrt W$.
If $Y>u\ge1$, then $V_B<C_bW^2/u^2$, so
\[
 \Prob^\sharp_\alpha(\mathcal G_\alpha,Y>u)
                         \le C_pu^{-2p/(p+1)}
       \quad(1\le u\le C_b\sqrt W).
\]
For fixed $r<2$, choose $p$ with $r<2p/(p+1)$ and integrate.
To pay $\mathcal G_\alpha^c$, first note that
Lemma~\ref{lem:physical-local-window} and codimension-two
interlacing give $V_B\ge cW/\eta$ at
$\eta=W^{-1+\gamma}$, with $0<\gamma<\min(\gamma_0,\xi)$,
 outside an event of any fixed failure
power. Thus $q_B^2\ge c/(W\eta)$ there and the diagonal floor
$q_B^2\ge t$ pays the complement:
\[
 \E^\sharp_\alpha q_B^{-2}\le C_\xi W^\xi
                      \quad\hbox{for every fixed }\xi>0.
\]
Using \eqref{eq:spread-good-event}, H\"older bounds the bad
contribution to $\E^\sharp Y^r$ by
$CW^{\xi r/2-(1-r/2)/50}$. Choose $\xi$ sufficiently small.
This proves \eqref{eq:physical-inverse-spread}.
\end{proof}

\begin{proof}[Proof of Theorem~\ref{thm:determinant-score}]
For $0<r\le2$, Jensen gives
$\E_{x,y}|a|^{2-r}\le q_B^{2-r}$.
Since $w^*Rw=a^{-1}$, the selected numerator of
$|w^*Rw|^r$ is bounded by
$\det(B)^2q_B^{2-r}$ on the two attainable minor sectors.
Restore the complete row mass as above, now with insertion
$q_B^{-r}$. Equation~\eqref{eq:physical-inverse-spread} at the
three neighboring ranks gives
$\E^\sharp_\alpha\one_{\{E\in I\}}|w^*Rw|^r\le C_r$.
All ranks involved remain in a fixed smaller bulk interval.

Stabilizer invariance of \eqref{eq:indexed-complement} makes
$\overline R$ commute with $K$ and be scalar on each eigenspace
of $K$. Its maximizing eigenvector is an allowed measurable choice
of $w(\theta)$. Conditional Jensen for $r\ge1$, and then the
first moment for $r<1$, prove the bound for $\norm{\overline R}$.

For a horizontal tangent $h$, choose a local parallel frame
$U_v$ of $v^\perp$ with derivative $\dot U_v=-vh^*U_v$.
Then
\[
 \dot K=-U_v^*h(v^*K_0U_v)-U_v^*K_0v(h^*U_v).
\]
The first derivative of the full determinant-square sector
integral is $2\mathfrak Z_\ell\Tr(\overline R\dot K)$.
 Lemma~\ref{lem:t26-sector-ibp} justifies this differentiation
through every inertia boundary, including multiple singular strata.
Substitution gives
$D\log\mathfrak Z_\ell[h]
=-4\Re(h^*\overline RQK_0v)$, proving the exact score.
At the selected root $QK_0v=-QG_dv$, hence
$\norm{\nabla\log\mathfrak Z_\ell}
\le4\norm{\overline R}\norm{G_d}$.
The pair law has every fixed Gaussian norm moment bounded:
below a fixed norm threshold this is immediate, and above it
the density bound $d\Prob^\sharp\le CW\,d\Prob$ is absorbed by
the Gaussian exponential tail. For $r<2$, choose $a>1$ with
$ra<2$ and apply H\"older. This proves
\eqref{eq:physical-determinant-score}.
\end{proof}
\section{The angular comparison and quantum unique ergodicity}
\label{sec:angular-que}
This section completes Part~I. It turns the estimates for an exposed pair
into uniform spreading of the whole eigenvector, and proves
Theorem~\ref{thm:que}.
Theorem~\ref{thm:physical-angular} states that, under the pair law, the
mass difference between the two blocks of a pair is within $C(M+L_\phi)/W$ of
the center, when tested against any bounded phase-invariant Lipschitz
function of the eigenvector. Proposition~\ref{prop:angular-comparison} is the
abstract version of this comparison on the sphere. Together with
Theorem~\ref{thm:physical-center}, it bounds every Fourier mode of the block
masses, and Theorem~\ref{thm:que} follows.
In Subsection~\ref{subsec:s11-the-finite-angular-measure} we isolate the angular measure and control the cost
of reconstructing a full unit vector from its restriction to the pair
(Lemma~\ref{lem:angular-ratio}). In Subsection~\ref{subsec:s11-two-exact-integration-by-parts-identitie} we derive two exact
integration-by-parts identities (Lemma~\ref{lem:angular-stein}), combine them
in Proposition~\ref{prop:angular-comparison}, and prove
Theorem~\ref{thm:physical-angular}. In Subsection~\ref{subsec:s11-fourier-modes-and-deterministic-site-tes} we sum the Fourier
estimates over deterministic site tests and prove Theorem~\ref{thm:que}.

\begin{theorem}\label{thm:physical-angular}
Fix $b,T,\vartheta$, a deterministic bulk rank $\alpha$ and an
adjacent pair $S=\{i,i+1\}$. Under \eqref{eq:pair-law}, let
$v=P_S\psi_\alpha/\sqrt{p_i+p_{i+1}}$ and $J=P_i-P_{i+1}$.
For every bounded phase-invariant Lipschitz function $\phi$ on the
full unit sphere, with bound $M$ and Lipschitz constant $L_\phi$,
\begin{equation}\label{eq:physical-angular}
 \left|\E^\sharp_{\alpha,S}
 [(v^*Jv-j_S(\lambda_\alpha))\phi(\psi_\alpha)]\right|
 \le C_{b,T,\vartheta}(M+L_\phi)/W.
\end{equation}
Here the contrast has the same sign as $J_e$ in
\eqref{eq:dyson-center}: $v^*Jv=(p_i-p_{i+1})/(p_i+p_{i+1})$.
It is the negative of the local convention in
Subsection~\ref{subsec:supported-source}; only sign-invariant
quadratic-form estimates are transferred from that subsection.
The center is \eqref{eq:dyson-center} on the regular domain, with
the bounded extension used in Theorem~\ref{thm:physical-center}.
\end{theorem}

The constant is uniform over the deterministic choice of the test
function. The estimate is signed and retains its dependence on the
selected full vector, which is why it gives more than an adjacent
mass-variance bound.

\subsection{The finite angular measure}\label{subsec:s11-the-finite-angular-measure}

 We isolate the directional measure appearing in the pair formula.
Lemma~\ref{lem:angular-ratio} controls the derivative cost of
reconstructing a full unit vector, even when the reconstruction map is
poorly conditioned. We define the reference covariance and finite
weighted measure, then use a tangent vector field to estimate the
normalizing ratio.

 For this subsection let   $\Naux \ge2$, $t=s/\Naux $  and $K_0=K_0^*$.
Suppose $\zeta=sm$ solves
  $m=\Naux ^{-1}\Tr(K_0-sm)^{-1}$  with $\Im\zeta=y>0$.
Assume $s,s^{-1},|\zeta|,y^{-1}$ are bounded by fixed constants.
Put
\begin{equation}\label{eq:angular-reference}
 a_*=2\Re\zeta,\quad c_*=|\zeta|^2,\quad
 \Gamma=s(K_0^2-a_*K_0+c_*I)^{-1}.
\end{equation}
The scalar equation and its real and imaginary parts give

\begin{equation}\label{eq:angular-reference-identities}
 \Tr\Gamma=\Naux ,\quad \Tr K_0\Gamma=\Naux a_*,\quad
 K_0^2\Gamma=sI+a_*K_0\Gamma-c_*\Gamma.
\end{equation}

The operators $\Gamma,K_0\Gamma,K_0^2\Gamma$ are uniformly
bounded; no lower eigenvalue bound on $\Gamma$ is assumed.

Write $a(v)=v^*K_0v$ and $b_2(v)=v^*K_0^2v$.
Consider a finite positive measure on projective space,

\begin{equation}\label{eq:angular-measure}
 d\nu(v)=C\exp\{-b_2(v)/t+a(v)^2/(2t)\}
                   \mathfrak Z(v)\,d\nu_\Naux (v),
\end{equation}

where $\mathfrak Z>0$ is smooth and phase invariant.
Let $Z=\nu(1)$ and $H(v)=\nabla\log\mathfrak Z(v)$, using the
real metric with $Df[h]=\Re(h^*\nabla f)$.
Suppose for some $r\ge0$ that

\begin{equation}\label{eq:angular-first-moment}
 \int vv^*\,d\nu(v)\preceq C_0(Z+r)\Gamma/\Naux  .
\end{equation}

Choose $\chi=\chi(b_2)$, nonnegative, smooth, equal to one on $b_2\le R_0$,
supported on $b_2\le R_1$, with fixed $R_0,R_1$ and derivative bounds.
All estimates below keep the mass $Z$ undivided.

\begin{lemma}\label{lem:angular-ratio}
Let $\mathcal L$ be injective and $D=\mathcal L^*\mathcal L$.
If $A$ is a Hermitian contraction, $I$, or $K_0$, then
\begin{equation}\label{eq:angular-ratio}
 \nu\left(\chi\frac{\norm{\mathcal L\Gamma Av}}{\norm{\mathcal Lv}}\right)
 \le C\{Z+r+\nu(\chi\norm H)\}.
\end{equation}
The constant is independent of every singular value of $\mathcal L$.
\end{lemma}

The estimate remains valid for an ill-conditioned positive denominator. This is needed because a physical reconstruction clock need not have a bounded condition number.

\begin{proof}
Work on the complex unit sphere with its real Euclidean metric;
all functions are phase invariant. Put $h=\norm{\mathcal Lv}$,
$q=\norm{\mathcal L\Gamma Av}$ and $T_D=\Tr(D\Gamma)>0$,
and extend $h,q$ homogeneously with degree one.
The tangent field
$Y=\Gamma\nabla h-v\Re(v^*\Gamma\nabla h)$, with
$\nabla h=Dv/h$, has divergence

\[
 \operatorname{div}_{\mathbb S}Y
 =\frac{2T_D}{h}-\frac{v^*D\Gamma Dv}{h^3}
        -(2\Naux -1)\Re(v^*\Gamma Dv)/h
 \ge T_D/h-(2\Naux -1)\norm{\mathcal L\Gamma v}.
\]

Here $D^{1/2}\Gamma D^{1/2}\preceq T_DI$ proves the first
inequality, and Cauchy--Schwarz proves the last term's bound.
Integrate the divergence of $\chi qY$ against
\eqref{eq:angular-measure}.

We bound its four types of terms. First the radial term is at most
  $2\Naux \nu(\chi q\norm{\mathcal L\Gamma v})$.
For each allowed $A$,
$\Gamma A\Gamma A^*\Gamma\preceq C\Gamma$.
Equation~\eqref{eq:angular-first-moment} and Cauchy--Schwarz
therefore bound this term by $C(Z+r)T_D$.

Second, the regression gradient satisfies the exact identity
\[
 \Gamma\nabla\log(\nu/\mathfrak Z)
 =-\frac{2}t\{sv+(a_*-a)\Gamma K_0v
                    -(c_*+b_2-a^2)\Gamma v\}.
\]
This follows by differentiating $-b_2/t+a^2/(2t)$ tangentially
and using \eqref{eq:angular-reference-identities}.
On the cutoff its scalar coefficients are bounded.
Its pairing with $qY$ is bounded by
  $C\Naux q(\norm{\mathcal Lv}+\norm{\mathcal L\Gamma K_0v}
                      +\norm{\mathcal L\Gamma v})$.
The first-moment bound again gives $C(Z+r)T_D$.
The determinant gradient costs at most
$\chi q\norm{\mathcal L\Gamma H}\le C\chi T_D\norm H$,
since both $q$ and $\norm{\mathcal L\Gamma}$ are bounded by
$C\sqrt{T_D}$.

Third, put $B=A^*\Gamma D\Gamma A$.
Away from $q=0$, $\nabla q=Bv/q$ and
$\nabla q^*\Gamma\nabla q\le\Tr(B\Gamma)\le CT_D$.
Also $\nabla h^*\Gamma\nabla h\le T_D$.
Thus $|Dq[Y]|\le CT_D+q\norm{\mathcal L\Gamma v}$,
whose integral is at most   $CT_D[Z+(Z+r)/\Naux ]$.
Regularization of $q$ at zero justifies the calculation.

Finally
$\Gamma\nabla\chi=2\chi'(b_2)(\Gamma K_0^2v-b_2\Gamma v)$.
The same quadratic estimates apply to these bounded spectral
multipliers; in particular
$K_0^4\Gamma^3\preceq C\Gamma$.
The cutoff term is at most   $C(Z+r)T_D/\Naux $.
Divide the sum of all bounds by $T_D$ to prove
\eqref{eq:angular-ratio}.
\end{proof}

\subsection{Two exact integration-by-parts identities}\label{subsec:s11-two-exact-integration-by-parts-identitie}

 We compare the signed pair contrast with its covariance center.
Lemma~\ref{lem:angular-stein} supplies two exact identities, and
Proposition~\ref{prop:angular-comparison} combines them into the angular
estimate. We choose bounded tangent fields, eliminate the longitudinal
term between the identities, and integrate the cutoff and determinant
errors under the actual pair law.

 For Hermitian $\norm J\le1$ put
$j=v^*Jv$,   $j_*=\Naux ^{-1}\Tr\Gamma J$  and
$u_k=v^*K_0^k\Gamma v$ for $k=0,1,2$.
For $A=J,K_0$ define
\[
 r_A=\frac{\Re(v^*\Gamma Av)}{u_0},\quad
 W_A=\Gamma(A-r_AI)v,\quad
 \Delta_A=Dr_A[\Gamma v],\quad
 b_J=\Re(v^*K_0\Gamma Jv)-r_Ju_1.
\]
The fields $W_A$ are norm-tangent; their possible vertical phase
components do not affect phase-invariant functions. Let
$V_\Gamma=u_2-u_1^2/u_0$ and $g=1-V_\Gamma/s$.
Equation~\eqref{eq:angular-reference-identities} gives
\begin{equation}\label{eq:positive-longitudinal}
 g=\frac{u_0}{s}\{(u_1/u_0-a_*/2)^2+y^2\}>0 .
\end{equation}
On $b_2\le R_1$, Jensen gives
$u_0\ge s/(b_2-a_*a+c_*)$, so $u_0,g$ are bounded below.
All the displayed coefficients, their derivatives along $W_J,W_{K_0}$,
and these fields' norms are bounded on the cutoff.

\begin{lemma}\label{lem:angular-stein}
For smooth phase-invariant $F$,

\begin{align}
 2\Naux \nu[g(a-a_*)F]
 &=\nu[W_{K_0}F+F\Re(H^*W_{K_0})-F\Delta_{K_0}],
 \label{eq:longitudinal-stein}\\
 2\Naux \nu[(j-j_*)F]
 &=\nu[W_JF+F\Re(H^*W_J)-F\Delta_J]
           +\frac{2\Naux }{s}\nu[(a-a_*)b_JF].
 \label{eq:contrast-stein}
\end{align}

\end{lemma}

These identities include the derivative of the complete angular density. They allow the determinant score to be estimated only after physical averaging.

\begin{proof}
Extend $r_A$ homogeneously with degree zero. Since $W_A$ is
tangent, its real ambient divergence equals its sphere divergence:
  $\operatorname{div}W_A=2\Re\Tr(\Gamma A)-2\Naux r_A-\Delta_A$.
For $A=K_0$, $r_A=u_1/u_0$ and $W_{K_0}a=2V_\Gamma$.
Substituting the gradient in the proof of
Lemma~\ref{lem:angular-ratio} and
\eqref{eq:angular-reference-identities} shows that divergence plus
regression derivative equals
  $-2\Naux g(a-a_*)-\Delta_{K_0}$.
For $A=J$ the same substitution gives
  $-2\Naux (j-j_*)+(2\Naux /s)(a-a_*)b_J-\Delta_J$.
Integrate the divergence of $FW_A$ and retain the derivative
$\Re(H^*W_A)$ of the full determinant factor. This gives the two
identities with the stated signs.
\end{proof}

The next proposition eliminates the longitudinal term by a second
integration by parts, using the positive coefficient $g$.

\begin{proposition}\label{prop:angular-comparison}
If $\phi$ is phase invariant on the unit sphere of
$\operatorname{Ran}\mathcal L$, with bound $M$ and Lipschitz
constant $L_\phi$ for the chordal Euclidean metric on that unit sphere, then, for
$\Phi(v)=\phi(\mathcal Lv/\norm{\mathcal Lv})$,

\begin{equation}\label{eq:angular-comparison}
 |\nu[\chi(j-j_*)\Phi]|
 \le\frac{C(M+L_\phi)}\Naux \{Z+r+\nu(\chi\norm H)\}.
\end{equation}

\end{proposition}

The error depends on the stated first moment and integrated score. A pointwise comparison of the whole angular density is not required.

\begin{proof}
Use \eqref{eq:contrast-stein} with $F=\chi\Phi$.
Direct differentiation gives
\[
 |W_A\Phi|\le2L_\phi
 \frac{\norm{\mathcal L\Gamma Av}
                  +|r_A|\norm{\mathcal L\Gamma v}}{\norm{\mathcal Lv}}.
\]
Lemma~\ref{lem:angular-ratio} bounds these derivatives.
All cutoff and $\Delta_A$ terms are bounded on their support,
and the determinant term costs $CM\nu(\chi\norm H)$.
For the last term in \eqref{eq:contrast-stein}, apply
\eqref{eq:longitudinal-stein} with
$F=\chi b_J\Phi/(sg)$. Both $b_J/(sg)$ and its derivative along
$W_{K_0}$ are bounded by \eqref{eq:positive-longitudinal}.
The same ratio bound therefore applies. Dividing by   $2\Naux $  proves
\eqref{eq:angular-comparison}. Smooth approximation extends the
argument to the stated Lipschitz class.
\end{proof}

\begin{proof}[Proof of Theorem~\ref{thm:physical-angular}]

\emph{Step 1: the angular form.} We fix the exterior data and the energy
$E$. Theorem~\ref{thm:coarea} expresses the law of the direction $v$ of
the eigenvector on the pair. This law has exactly the angular form
\eqref{eq:angular-measure}, with

\[
 K_0=F-E,\qquad \Naux =2W,
\]

and with the determinant factor

\[
 \mathfrak Z(v)=
 \E_{G_{\Naux -1}}\!\left[
 \det(U_v^*K_0U_v+G)^2
 \one_{\{N_-=\alpha-N_R(E)-1\}}\right].
\]

Here the indicator fixes the inertia of the matrix on $v^\perp$, which is
what keeps the rank equal to $\alpha$. Under the size-biased pair law,
the angular measure is multiplied by the factor $n/2$.

\emph{Step 2: the first moment on a good domain.} The first input of
Proposition~\ref{prop:angular-comparison} is the first moment bound
\eqref{eq:angular-first-moment}. We obtain it from the regular
fresh-Gaussian moment estimate, Theorem~\ref{thm:fresh-palm}, with
  $r=(n/2)C_K\Naux ^{-K}$  and without dividing by the angular mass $Z$. That theorem,
however, needs a deformation of polynomial size, while $F(E)$ is
unbounded near the poles of the exterior. So we first pass to a
polynomial deformation domain. We restrict the outer variables to
\[
 \operatorname{dist}(E,\operatorname{Spec}H_R)\ge W^{-M_0},
 \qquad \norm V+\norm C\le W.
\]
On this domain
\[
 \norm{F(E)}\le W+W^{M_0+2}.
\]
These restrictions depend only on the exterior data and on $E$. So they
are constant during every angular integration by parts, which acts only
on $v$. The regularity certificate required by
Theorem~\ref{thm:fresh-palm} is also a restriction on the outer data.

On this domain we check that the hypotheses of
Theorem~\ref{thm:fresh-palm} are met exactly. For each such $(\omega,E)$,
the angular measure, with its coarea constant, is exactly the eigenpair
density at $E$ in the ordered fresh rank $\alpha-N_R(E)$. The scalar
Gaussian in \eqref{eq:palm-scalar-weight} has variance   $s/\Naux ^2$. This is
precisely the variance of the scalar part of the same fresh
  $\GUE_\Naux (s/\Naux )$. We therefore apply Theorem~\ref{thm:fresh-palm} with
moment order $r=1$ and an arbitrary vector $q$. (This $r$ is the order of
the moment in that theorem, not the remainder $r$ of
\eqref{eq:angular-first-moment}.) For moment order one, its bound reads
  $\E[w|\sqrt \Naux q^*u|^2]\le C(Z+C_K\Naux ^{-K})\,q^*\Gamma q$. Since $q$ is arbitrary, this is a Loewner
inequality for the fresh eigenpair density, with $Z$ the sector mass
$Z_{\ell,E}$ of Theorem~\ref{thm:fresh-palm}. By Step~1 the angular
measure under the size-biased pair law, and its mass, are $n/2$ times
this density and $Z_{\ell,E}$. Multiplying by $n/2$ gives
\eqref{eq:angular-first-moment} with the unnormalized mass $Z$ and
remainder $r=(n/2)C_K\Naux^{-K}$, which is at most $C_K\Naux^{-K}$ after
replacing $K$ by $K+1$, since $n\le TW=T\Naux/2$.

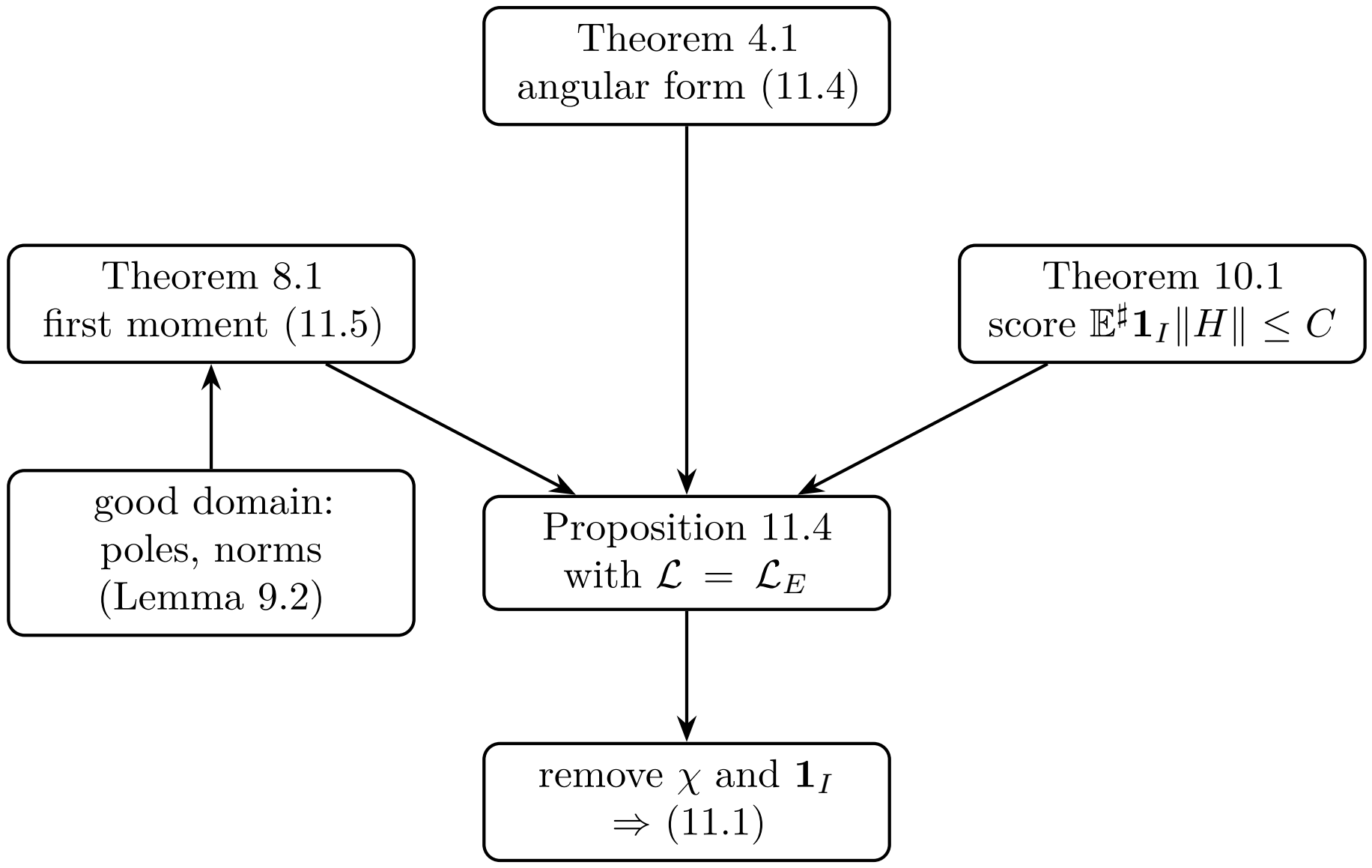
\begin{figure}[ht]
\centering
\begin{tikzpicture}[
  box/.style={draw,thick,rounded corners,align=center,font=\small,
              inner sep=4pt,text width=3.3cm},
  arr/.style={->,thick,>=Stealth}]
 \node[box] (coarea) at (0,0) {Theorem~\ref{thm:coarea}\\ angular form \eqref{eq:angular-measure}};
 \node[box] (palm) at (-4.2,-2.1) {Theorem~\ref{thm:fresh-palm}\\ first moment \eqref{eq:angular-first-moment}};
 \node[box] (score) at (4.2,-2.1) {Theorem~\ref{thm:determinant-score}\\ score $\E^\sharp\one_I\norm H\le C$};
 \node[box] (dom) at (-4.2,-4.3) {good domain:\\ poles, norms\\ (Lemma~\ref{lem:t26-pair-pole-exclusion})};
 \node[box] (prop) at (0,-4.3) {Proposition~\ref{prop:angular-comparison}\\ with $\mathcal L=\mathcal L_E$};
 \node[box] (end) at (0,-6.5) {remove $\chi$ and $\one_I$\\ $\Rightarrow$ \eqref{eq:physical-angular}};
 \draw[arr] (coarea)--(prop);
 \draw[arr] (palm)--(prop);
 \draw[arr] (score)--(prop);
 \draw[arr] (dom)--(palm);
 \draw[arr] (prop)--(end);
\end{tikzpicture}
\caption{The structure of the proof of Theorem~\ref{thm:physical-angular}.
The change of variables of Theorem~\ref{thm:coarea} puts the law of the
direction $v$ in the angular form required by
Proposition~\ref{prop:angular-comparison}. The proposition needs a first
moment bound, which comes from Theorem~\ref{thm:fresh-palm} on a good
domain of polynomial size, and a score bound, which comes from
Theorem~\ref{thm:determinant-score}. The good domain excludes energies
close to the exterior spectrum, whose probability is controlled by
Lemma~\ref{lem:t26-pair-pole-exclusion}, and large exterior norms. The
proposition gives the comparison with the cutoff $\chi$ and on the bulk
interval $I$, and the last step removes both.}
\label{fig:pa-chain}
\end{figure}

\emph{Step 3: the excluded events.} It remains to control the outer
configurations that we excluded in Step~2. There are three of them. The
first is the pole set, where $E$ is within $W^{-M_0}$ of the exterior
spectrum. Lemma~\ref{lem:t26-pair-pole-exclusion}, applied with $M=M_0$,
bounds its probability under the actual pair law by $C_TW^{3-M_0}$. The
second is the event of a large exterior norm. Gaussian norm tails bound
its probability by an arbitrarily high negative power of $W$, even after
multiplication by the factor $n/2$ of the pair law. The third is the complement of the regular domain
$D_S(E)$ of Theorem~\ref{thm:physical-center} for $e=S$, to which Step~2
restricts. On $D_S(E)$, $K_0=F(E)-E=K_S(E)$ satisfies
\eqref{eq:fresh-regularity} with $c_0=s$, the root $\zeta=sm_s(E)$
satisfies the hypotheses of
Subsection~\ref{subsec:s11-the-finite-angular-measure}, and the covariance
$\Gamma$ there is that of Theorem~\ref{thm:fresh-palm}. Under
Theorem~\ref{thm:coarea} the energy is $E=\lambda_\alpha$, so
Theorem~\ref{thm:physical-center} bounds the pair-law probability of this
event by $C_qW^{1-q}$ for every fixed $q$. Each of the three events costs
at most $2M$ times its probability, since $|v^*Jv|\le1$ and the center is
bounded by one.

The order of the choices matters. We first fix $M_0$ large enough for the
desired error in replacing a bounded observable. This fixes the
polynomial norm exponent $M_0+2$. Then we choose the exceptional power
$K$ in Theorem~\ref{thm:fresh-palm} for this fixed exponent. Throughout,
the argument keeps the original dimension, variance, rank, trace
coordinate, determinant and Palm normalization. In particular, it does
not assume a compression theorem with determinant weights at an infinite
pole. The separate eigenvector and counting compression of
Lemma~\ref{lem:physical-local-window} is used only within its stated
scope.

\emph{Step 4: the score.} The second input of
Proposition~\ref{prop:angular-comparison} is the integrated score.
Theorem~\ref{thm:determinant-score} computes the score exactly, and
gives
\[
 \E^\sharp\one_I\norm H\le C,
\]
where $I$ is the fixed bulk energy interval.

\emph{Step 5: the comparison.} We now apply
Proposition~\ref{prop:angular-comparison} with the exact reconstruction
map $\mathcal L=\mathcal L_E$ of \eqref{eq:clock-reconstruction}. Since
$\psi_\alpha=\mathcal L_Ev/\norm{\mathcal L_Ev}$, the function $\Phi$ of the
proposition is $\phi(\psi_\alpha)$, and the center $j_*$ is
$j_S(\lambda_\alpha)$. The proposition bounds the cutoff comparison by
  $C(M+L_\phi)\Naux ^{-1}$  times $Z+r+\nu(\chi\norm H)$. We integrate this over
the exterior and the energy. The $Z$ term then contributes at most one,
because the full law is normalized. The score term is bounded by $C$, by
Step~4. The additive remainder $r$ contributes $O(nW^{-K})$ over the
fixed energy interval, before the factor   $\Naux ^{-1}$. We choose $K$ first
and large enough, so this is negligible. Since   $\Naux =2W$, this gives the
right side of \eqref{eq:physical-angular}, still with the cutoff $\chi$
and restricted to $I$.

\emph{Step 6: removing the cutoffs.} The cutoff $\chi=\chi(b_2)$ equals
one when $b_2(v)=\norm{K_0v}^2\le R_0$. At the actual eigenvector the
eigenvalue equation gives

\[
 K_0v=-G_\Naux v,
\]

so   $b_2(v)\le\norm{G_\Naux }^2$. We take $R_0$ larger than a fixed threshold
for the Gaussian norm. Then $\chi\ne1$ only when   $\norm{G_\Naux }$  exceeds that
threshold. By the exponential Gaussian norm tail and the bound
$d\Prob^\sharp\le(n/2)d\Prob$, this event has a negligible probability
under the pair law. The tails of the rank location
(Theorem~\ref{thm:rank-location}) control the event that $\lambda_\alpha$
leaves the fixed bulk interval $I$. All observables discarded in these
steps are bounded, so each removal costs only the probability of the
removed event. This proves \eqref{eq:physical-angular}.
\end{proof}

\subsection{Fourier modes and deterministic site tests}\label{subsec:s11-fourier-modes-and-deterministic-site-tes}

 We complete Theorem~\ref{thm:que} by turning neighboring-block
comparisons into a bound for every deterministic site test. Cyclic
symmetry gives an exact identity for each Fourier mode of the block
masses. We apply the angular and center estimates to these modes, sum
their contributions, and use unitary invariance within each block.

 \begin{proof}[Proof of Theorem~\ref{thm:que}]
For $1\le m<n$ set
$z_m=e^{2\pi im/n}$,
$A_m=\sum_kz_m^{-k}p_{\alpha k}$,
$\gamma_m=\E|A_m|^2$ and $d_m=|1-z_m|$.
Cyclic invariance gives the exact pair identities
\begin{equation}\label{eq:signed-pair}
 \E^\sharp_{\alpha,S}|A_m|^2=\gamma_m,\qquad
 \E^\sharp_{\alpha,S}
 \left[\frac{p_i-p_{i+1}}{p_i+p_{i+1}}\overline{A_m}\right]
 =\frac{1-z_m}{2}z_m^i\gamma_m.
\end{equation}
Indeed the pair density cancels the contrast denominator.
The cyclic covariance
$\E[p_i\overline{A_m}]=z_m^i\gamma_m/n$
follows by summing its translates against $z_m^{-i}$.
The first identity follows because $|A_m|^2$ is invariant under
cyclic shifts and $\sum_i(p_i+p_{i+1})=2$.

The function $\overline{A_m}$ of the full unit vector is bounded
by one and has Lipschitz constant at most two. Apply
Theorem~\ref{thm:physical-angular}, then
Theorem~\ref{thm:physical-center} and Cauchy--Schwarz:
\[
 \frac{d_m}{2}\gamma_m
 \le C/W+\frac{C}{n^{1/4}W^{3/4}}\sqrt{\gamma_m}.
\]
Young's inequality absorbs $d_m\gamma_m/4$ and gives
\begin{equation}\label{eq:mode-bound}
 \gamma_m\le C\left\{\frac{1}{Wd_m}
              +\frac{1}{\sqrt n\,W^{3/2}d_m^2}\right\}
 \le C\{n/W+(n/W)^{3/2}\},
\end{equation}
where $d_m=2\sin(\pi m/n)\ge4/n$.

For a deterministic real block test $t$, $\norm t_\infty\le1$,
its normalized Fourier coefficients
$\widehat t_m=n^{-1}\sum_kt_kz_m^k$ satisfy
$\sum_m|\widehat t_m|^2=n^{-1}\sum_kt_k^2\le1$.
Cyclic invariance makes distinct nonzero Fourier modes orthogonal
in expectation. Hence
\[
 \E\left|\sum_kt_k(p_{\alpha k}-1/n)\right|^2
 =\sum_{m\ne0}|\widehat t_m|^2\gamma_m
 \le C\{n/W+(n/W)^{3/2}\}.
\]
Both deterministic suprema can now be taken outside expectation.

Finally conjugate the original matrix by independent Haar unitaries
inside its blocks. This preserves its law and the eigenvalue rank.
Conditional on the block masses, the rotated directions are independent
Haar vectors in $\C^W$. Their moments are
$\E|u_a|^2=1/W$ and
$\E|u_a|^2|u_c|^2=(1+\delta_{ac})/[W(W+1)]$.
For $\bar v_k=W^{-1}\sum_av_{ka}$ and
$\overline{v^2}_k=W^{-1}\sum_av_{ka}^2$, conditional variance gives
\begin{equation}\label{eq:haar-site-bridge}
 \mathcal Q_\alpha(v)
 =\E\left|\sum_k\bar v_k(p_{\alpha k}-1/n)\right|^2
 +\frac{1}{W+1}\E\sum_kp_{\alpha k}^2
                     (\overline{v^2}_k-\bar v_k^2).
\end{equation}
Since $\sum_kp_{\alpha k}^2\le1$, the second term is at most
$(W+1)^{-1}$. This proves the stated quantitative bound, including
bounded $n$, and therefore QUE whenever $n/W\to0$.
This argument uses only finitely many fixed auxiliary moment orders.
For the remaining widths $W<W_0$, the defining statistic is at most
$4$, while $n/W\ge4/W_0$. Enlarging $C$ by $W_0$ therefore proves
the same bound at these widths, without an inverse-moment claim there.
 \end{proof}

\section{From prescribed ranks to the physical spectral window}
\label{sec:window-transfer}
Theorem~\ref{thm:que} is a statement at a fixed rank, while the critical
law concerns a fixed microscopic energy window. This section transfers rank
estimates to the window. It then uses them for Part~III, to show that the
profile variance vanishes as $\tau\to0$ and that the marks have no atoms.
Proposition~\ref{prop:refined-count} sharpens the fluctuation of the
eigenvalue count to $a_{n,W}=\sqrt n+n/W+W^{2/5}\log(2W)$.
Theorem~\ref{thm:window-transfer} bounds the window expectation of a bounded,
shift-invariant nonnegative observable by its rankwise bound, without
conditioning on an eigenvalue. Theorem~\ref{thm:small-critical} shows that the
variance of the critical marks is at most $C\tau^{1-\eta}$, and that the marks
have no atoms.
In Subsection~\ref{subsec:s12-an-improved-fresh-count-fluctuation-boun} we prove a sharper count bound for a fresh block
(Lemma~\ref{lem:sharper-fresh-count}) and transfer it to the ring. In
Subsection~\ref{subsec:s12-scalar-gaussian-smoothing} we use the independent scalar Gaussian part of the matrix
to pass from ranks to energy windows. In Subsection~\ref{subsec:s12-diffuse-marks-and-small-critical-varianc} we apply this to
the Fourier modes of the block masses.

 Constants may depend on fixed $b,T$, a compact bulk interval, and a fixed moment order.

\subsection{An improved fresh-count fluctuation bound}\label{subsec:s12-an-improved-fresh-count-fluctuation-boun}
 We sharpen the count fluctuation estimate enough to transfer rank
bounds to a microscopic window. Lemma~\ref{lem:sharper-fresh-count}
gives an order-$W^{1/3}$ fluctuation bound for a fresh deformed block;
Proposition~\ref{prop:refined-count} transfers the improvement to the
ring. We balance the error from smoothing the count against its Gaussian
gradient, then repeat the compensated-phase comparison.

 \begin{lemma}\label{lem:sharper-fresh-count}
Let $A$ be $W\times W$ GUE with entry variance $c/W$ and let $F$ be any
independent Hermitian deformation, with no bound on its norm. At every
real deterministic $e$, conditionally on $F$, for every fixed $p\ge1$,
\[
 \|N_-(A+F-e)-\mathbb E_A N_-(A+F-e)\|_{L^p(A)}
 \le C_{c,p}W^{1/3}.
\]
\end{lemma}

Smoothing balances a count error with a Gaussian gradient error. The bound is uniform over all Hermitian deformations, so it applies after conditioning on an exterior.

\begin{proof}
The single-block case of Theorem~\ref{thm:external-count} gives
$\mathbb E_A(N(J))_r\le(C_cW|J|)^r$, uniformly in $F,e$.
For an interval of length $2\eta$ with $W\eta\ge1$, the conversion
from factorial to ordinary moments consequently gives
$\|N(J)\|_q\le C_{c,q}W\eta$ for every fixed $q\ge1$.
Choose a smooth approximation $f_\eta$ to $\mathbf1_{(-\infty,0)}$,
equal to it off $[-\eta,\eta]$, taking values in $[0,1]$, and with
$|f'_\eta|\le C/\eta$. Put $X_\eta=\operatorname{Tr}f_\eta(A+F-e)$.
Then
\[
 |N_- -X_\eta|\le N([-\eta,\eta]),\qquad
 \|\nabla_A X_\eta\|_{\rm HS}^2
 =\operatorname{Tr}|f'_\eta(A+F-e)|^2
 \le C\eta^{-2}N([-\eta,\eta]).
\]
For $p\ge2$, Gaussian $L^p$ Poincare in the real Hilbert space of
Hermitian matrices, whose covariance is $(c/W)I$, gives
\[
 \|X_\eta-\mathbb E_A X_\eta\|_p
 \le C_p\sqrt{c/W}\,\|\nabla_A X_\eta\|_p
 \le C_{c,p}\eta^{-1/2}.
\]
(The inequality follows, for example, by integrating the Gaussian
Ornstein--Uhlenbeck gradient contraction; its constant depends only on
$p$.) Thus the sharp centered count has norm at most
$C_{c,p}(W\eta+\eta^{-1/2})$. Take $\eta=W^{-2/3}$.
The result for $1\le p<2$ follows from the $p=2$ bound.
\end{proof}

\begin{proposition}[Critical-scale refinement of the physical count bound]
\label{prop:refined-count}
Put
\[
 a_{n,W}=\sqrt n+\frac{n}W+W^{2/5}\log(2W).
\]
At every deterministic compact-bulk energy and fixed $p<\infty$,

\[
 \|N_H(e)-NF_{\rm sc}(e)\|_p\le C_pa_{n,W}.
\]

For every deterministic rank in a fixed bulk margin, with
  $e_\alpha=F_{\rm sc}^{-1}(\alpha/N)$,

\[
 \|\lambda_\alpha-e_\alpha\|_p\le C_p a_{n,W}/N.
\]

In particular, for every fixed $\tau\in(0,T]$, along $n/W\to\tau$,
$a_{n,W}/\sqrt n\to1$. The constants above do not require a lower
bound on $n/W$.
\end{proposition}

At a fixed critical ratio the leading bound is of order the square root of the block count. This improvement is used only to transfer rank bounds into a microscopic energy window.

\begin{proof}
The compensated-phase calculation in Section~\ref{sec:rank-location}
gives
\[
 \|\Xi_{n-1}-\E\Xi_{n-1}\|_p\le C_p\sqrt n,\qquad
 \E\Xi_{n-1}=(n-1)WF_{\rm sc}(e)+W(\ell+bh)+O(n/W).
\]
The remaining independent diagonal block contributes the sharp count
$N_-(A_0-e-\Sigma)$. Lemma~\ref{lem:fresh-count}, with its
original scalar regularity and exceptional-set payment, gives an
$L^p$ predictable error at most
$C_p(W^{-1}+W^{2/5}\log(2W))$. Replace only its final
$C_p\sqrt W$ conditional fluctuation bound by
Lemma~\ref{lem:sharper-fresh-count}; no independence of the compensated exterior
and its deformation is added. The deterministic constants still cancel
by $b+c/2=1/2$, exactly as in that proof. This proves the first estimate,
since $W^{1/3}\le C W^{2/5}\log(2W)$.

The second assertion follows by monotonicity:
on a fixed bulk neighborhood, $|\lambda_\alpha-e_\alpha|>t$ forces
a count deviation at $e_\alpha+t$ or $e_\alpha-t$ of size at least
  $cNt-1$. Use the first assertion with a moment order larger than $p$
and integrate its tail. For the part outside this neighborhood, use
higher count moments at its endpoints and the uniform fixed
moments of $\|H\|$. Taking the auxiliary moment order greater than $2p$
makes this contribution   $O((a_{n,W}/N)^p)$. This also pays the physical
energy tails. Finally $n/W\to\tau>0$ implies the asserted ratio limit.
\end{proof}

\subsection{Scalar Gaussian smoothing}\label{subsec:s12-scalar-gaussian-smoothing}
 We use the scalar part of the original Gaussian matrix to move from
ranks to energy windows. Theorem~\ref{thm:window-transfer} bounds the
window expectation of any bounded, shift-invariant nonnegative
observable from its rankwise bound. We integrate the independent scalar
first, localize the relevant ranks by the refined count estimate, and
sum the rank tails.

 Let $\rho=\rho(E)>0$, $w\ge0$ be continuous and compactly supported,
with $\int w=1$. For every rank, let $A_\alpha$ be a nonnegative
observable satisfying $0\le A_\alpha\le1$ and invariant under adding a
scalar matrix to $H$. It may be correlated with the entire traceless
matrix, the rank, and all spectral gaps. Suppose that for ranks in a
fixed bulk margin containing   $NF_{\rm sc}(E)$,
\[
 \mathbb E A_\alpha\le\epsilon,\qquad 0<\epsilon\le1.
\]
Define its microscopic physical window expectation by

\[
 \mathcal V_{n,W}=\mathbb E\sum_\alpha
       w(N\rho(\lambda_\alpha-E))A_\alpha.
\]

\begin{theorem}[Window transfer without conditioning the energy]
\label{thm:window-transfer}
For every fixed $q>1$,
\[
 \mathcal V_{n,W}
 \le C_q\frac{a_{n,W}}{\sqrt n}\epsilon^{1-1/q}+r_{n,W,q},
\]
where $r_{n,W,q}\to0$ along every fixed critical sequence $n/W\to\tau>0$.
The remainder can be chosen independently of the observables and
$\epsilon$. Constants depend on $b,T,E,w$, the bulk margin and $q$.
In particular
\[
 \limsup_{n/W\to\tau}\mathcal V_{n,W}
 \le C_q\bigl(\limsup_{n/W\to\tau}\epsilon\bigr)^{1-1/q}.
\]
The last formula is interpreted by first applying the finite inequality
with any fixed larger bound on the limsup. No Gaussianity of a physical
current and no rank-conditioned energy-density lower bound is assumed.
\end{theorem}

The transfer retains the dependence between eigenvectors and their selected energies. Its power loss is harmless for vanishing critical variance and is separate from the quantitative QUE theorem.

\begin{proof}
Write exactly $H=H_0+tI$,   $t=N^{-1}\operatorname{Tr}H$. Gaussian
orthogonality gives

\[
 H_0\perp t,\qquad t\sim N(0,\sigma^2),\qquad
 \sigma^2=c/(NW)=c/(nW^2).
\]

Ordered eigenvalues satisfy $\lambda_\alpha=\lambda^0_\alpha+t$;
the invariant $A_\alpha$ is a function of $H_0$ alone. This is the
original Gaussian law disintegrated exactly. If

\[
 h(x)=\mathbb E_t w(N\rho(x+t-E)),
\]

then, for $\operatorname{supp}w\subset[-M,M]$,

\[
 0\le h(x)\le\frac{C}{N\sigma}
 \exp\left[-\frac{((|x-E|-M/(N\rho))_+)^2}{2\sigma^2}\right]
 \le C/\sqrt n.
\]

Here we used $\int w=1$ in the scalar convolution, and
  $N\sigma=\sqrt{cn}$. The bulk-rank estimate of
Proposition~\ref{prop:refined-count}, and the moments of $t$, give
  $\|\lambda^0_\alpha-e_\alpha\|_q\le C_q a_{n,W}/N$.
Since the semicircle density is globally bounded above,
  $|e_\alpha-E|\ge c_E|\alpha-NF_{\rm sc}(E)|/N$.
Writing $m=|\alpha-NF_{\rm sc}(E)|$, Markov's inequality and the
Gaussian bound for $h$ therefore give
\[
 \mathbb E[A_\alpha h(\lambda^0_\alpha)]
 \le\frac{C}{\sqrt n}
 \min\left(\epsilon,\frac{C_q a_{n,W}^q}{(1+m)^q}\right)
\]
for every rank in the fixed bulk margin. For bounded $m$ enlarge $C_q$;
for large $m$ split according to
$|\lambda^0_\alpha-e_\alpha|>|e_\alpha-E|/2$ and use
$e^{-cm^2/n}\le C_q(a_{n,W}/m)^q$ on the complement.
No independence of $A_\alpha$ and $\lambda^0_\alpha$ is used.

The sum over such ranks is bounded using the elementary integral
estimate, valid for $a\ge1$ and $0<\epsilon\le1$,
\[
 \sum_{m\ge0}\min(\epsilon,C_q(a/(1+m))^q)
 \le C'_q a\epsilon^{1-1/q}.
\]
For ranks outside the margin, monotonicity reduces a near-$E$ event
to the deviation of one of its two boundary ranks, whose classical
locations are a fixed distance from $E$. Their higher moment bounds,
and the same scalar Gaussian tail, yield

\[
 \sum_{\alpha\ \mathrm{outside}}
 \mathbb E h(\lambda^0_\alpha)
 \le C_Q\frac{N}{\sqrt n}(a_{n,W}/N)^Q
       +CNn^{-1/2}e^{-cL^2/n}.
\]

Choose a fixed sufficiently large $Q$; along any fixed critical
sequence this tends to zero. This bound defines $r_{n,W,q}$ and is
independent of the observables and $\epsilon$. Summing proves the result.
\end{proof}

\subsection{Diffuse marks and small critical variance}\label{subsec:s12-diffuse-marks-and-small-critical-varianc}
 We determine what the rank estimates imply for critical profile
marks. Theorem~\ref{thm:small-critical} shows that their variance vanishes
as the critical ratio tends to zero and that the marks have no atoms.
We transfer each Fourier bound to the physical window, use spatial
symmetry to diagonalize the covariance, and then control concentration
on short arcs.

 Let $\mathfrak L_{\tau,b,E}$ denote the nonempty family of marked
subsequential laws given by Theorem~\ref{thm:counting} and Proposition~\ref{prop:stationarity} at a fixed critical ratio $\tau>0$.
For $\mathsf P\in\mathfrak L_{\tau,b,E}$ define
\[
 V_{\mathsf P}(w,\phi)=\mathbb E_{\mathsf P}\sum_\alpha
 w(u_\alpha)|\nu_\alpha(\phi)-\bar\phi|^2.
\]
These expectations are finite. Theorem~\ref{thm:counting} and Proposition~\ref{prop:stationarity}
implies $V_{\mathsf P}(w,\phi)=V_{\mathsf P}(\phi)$ whenever
$\int w=1$: the weighted mark intensity is a translation-invariant
locally finite measure on the spectral line. Spatial rotation invariance
similarly diagonalizes its profile covariance. Put
\[
 v_m(\mathsf P)=\mathbb E_{\mathsf P}\sum_\alpha w(u_\alpha)
 \left|\int_{\mathbb T}e^{-2\pi imr}\nu_\alpha(dr)\right|^2,
 \qquad m\ne0.
\]
Then $0\le v_m\le1$, $v_{-m}=v_m$, and for smooth real $\phi$,
\[
 V_{\mathsf P}(\phi)=\sum_{m\ne0}|\widehat\phi_m|^2v_m(\mathsf P).
\]
The absolute summability of smooth Fourier coefficients justifies
expansion, and rotation invariance removes cross terms.

\begin{theorem}[Vanishing small-critical variance and diffuse critical marks]
\label{thm:small-critical}
For every $0<\eta<1$, uniformly over $0<\tau\le1$ and
$\mathsf P\in\mathfrak L_{\tau,b,E}$,
\[
 V_{\mathsf P}(\phi)\le
 C_{b,E,\eta}\|\phi-\bar\phi\|_\infty^2\tau^{1-\eta}.
\]
More precisely, for each fixed $T<\infty$, every $0<\tau\le T$,
all such $\mathsf P$, and each integer $m\ne0$,
\[
 v_m(\mathsf P)\le
 C_{b,E,T,\eta}
 \min\left(1,
   \left[\frac{\tau}{|m|}+\frac{\tau^{3/2}}{m^2}\right]^{1-\eta}\right).
\]
For each fixed $\tau$, almost surely every profile mark of the limiting
point process belongs to $H^{-s}(\mathbb T)$ for every $s>0$ and has
no atoms. These are statements about the actual marked subsequences;
no identification or uniqueness of their critical law is used.
\end{theorem}

The mode estimate controls high spatial frequencies as well as the small-length limit. It also pays the mass of the finite collar without presupposing critical-law uniqueness.

\begin{proof}
Apply Theorem~\ref{thm:window-transfer} to
$A_\alpha=|q\cdot p_\alpha|^2$ for deterministic grid-mean-zero
$\|q\|_\infty\le1$. Theorem~\ref{thm:que} gives $\epsilon\le C\tau_W$ on $n\le W$;
replace it by $\min(1,C\tau_W)$. Choose $q$ in the window bound with $1/q\le\eta$. For a fixed smooth spatial
test, subtract its grid mean first. The grid mean tends to $\bar\phi$,
and the original factorial count bound pays the resulting bounded
window error. Uniform integrability of compact-window counts allows
the expectation to pass to every marked subsequential limit. Scaling
the test proves the first bound, after increasing its constant on a
fixed non-small range of $\tau$.

For the mode bound use the quantitative Fourier inequality in the proof of Theorem~\ref{thm:que}:
\[
 \mathbb E|A_{\alpha,m}|^2
 \le C_T\left[\frac{1}{Wd_m}
       +\frac{1}{\sqrt n W^{3/2}d_m^2}\right],
 \qquad d_m=|1-e^{2\pi im/n}|.
\]
For each fixed integer $m\ne0$ and all large $n$,
$d_m\ge4|m|/n$. Use Theorem~\ref{thm:window-transfer} with this bound for $\epsilon$
and pass to the critical subsequence. All constants are independent
of $m$, so the conclusion holds for every integer $m$ under the
same limiting law.

Choose $0<\eta<2s$. The mode bound, for fixed $\tau$, implies
\[
 \mathbb E_{\mathsf P}\sum_\alpha w(u_\alpha)
 \|\nu_\alpha\|_{H^{-s}}^2
 =1+\sum_{m\ne0}(1+4\pi^2m^2)^{-s}v_m(\mathsf P)<\infty.
\]
Nonnegative Tonelli and an increasing family of spectral-window tests
imply that all marks have finite $H^{-s}$ norm, almost surely. Intersect
over positive rational $s$. Any atom of mass $a>0$ gives infinite
$H^{-s}$ norm when $s<1/2$: equivalently use Wiener's identity
$\lim_{M\to\infty}(2M+1)^{-1}\sum_{|m|\le M}|\widehat\nu_m|^2
=\sum_x\nu(\{x\})^2$, and the Cesaro decay of $v_m$ with Fatou.
The latter argument directly proves absence of atoms in every spectral
window, so does not need a separate regularity convention for measures.
\end{proof}

 \part{}
\section{The physical transmission process}
\label{sec:transmission}
This section begins Part~II. We cut the ring open into a path with
matched ports, and we follow how its transmission changes as the path grows
one block at a time. Time is the number of physical blocks divided by $W$.
At a fixed bulk energy $E$ we use the matched parameter $z=z_b(E)$ of
\eqref{eq:matched-parameters}.
Theorem~\ref{thm:transmission-generator} states that the power sums
$\Tr T^k$ of the transmission matrix are jointly tight, with continuous
limits that solve an explicit martingale problem. The generator is of
Dorokhov--Mello--Pereyra--Kumar type, and the transmission eigenvalues stay
between zero and one.
In Subsection~\ref{subsec:s13-the-scalar-derivative-and-a-contracting-} we compute the derivative of the scalar mean map
and identify its stable fixed point. In Subsection~\ref{subsec:s13-passive-inverse-bounds-and-the-mass-mean} we bound the
passive inverse factors (Lemma~\ref{lem:passive-bounds}) and the weighted
response that must relax (Lemma~\ref{lem:passive-source-mixing}). In
Subsection~\ref{subsec:s13-the-two-replica-contraction} we compute the covariance of two copies of the
transported response (Lemma~\ref{lem:passive-covariance}).

Write the Hermitian two-port Schur matrix as
$\mathbb X=(X_{ij})_{i,j=1}^2$, and eliminate its far port:
\[
 X=X_{11}+X_{12}(z-X_{22})^{-1}X_{21},\quad
 Y=-\Im X\succeq0,\quad K=z-X,\quad R=K^{-1}.
\]
Set
\begin{equation}\label{eq:passive-coordinates}
 Q=R^*YR,\quad S=I-2i\eta R,\quad
 T=4\eta Q=I-S^*S,\quad
 q=\Tr Q,\quad g=\Tr T,\quad D_0=\frac{2}{b(4-E^2)}.
\end{equation}
In particular $0\preceq T\preceq I$. Schur inversion also gives
$g=4\eta^2\norm{[(z-\mathbb X)^{-1}]_{21}}_{\HS}^2$,
so $g$ is the total transmission of this same physical path.
All vectors used in conditional statements are measurable before
the fresh layer.

\begin{theorem}\label{thm:transmission-generator}
The processes $P_k(t)=\Tr T_{\lfloor Wt\rfloor}^k$, $k\ge1$,
are jointly tight on compact positive-time intervals, with continuous
limits. Each limit solves the martingale problem with
\begin{align}
 A_k(P)&=-\frac{D_0k}{2}P_1P_k+
 \frac{D_0k}{2}\sum_{j=1}^{k-1}
             (P_jP_{k-j}-P_{j+1}P_{k-j}),\label{eq:transmission-drift}\\
 B_{k\ell}(P)&=D_0k\ell(P_{k+\ell}-P_{k+\ell+1}).
                                             \label{eq:transmission-bracket}
\end{align}
For every fixed $p\ge1$,
\begin{equation}\label{eq:conductance-moments}
 \E g_r^p\le C_{p,T}\left(1+\frac{W^p}{(r+1)^p}\right),
 \qquad \E P_1(t)\le\frac{2}{D_0t}.
\end{equation}
The constant in the first bound is length-uniform, i.e., for $p\ge1$ there is
$C_p$, uniform for $E$ in a compact bulk set, and then for every $T$ a
$W_0(T)$, such that \eqref{eq:conductance-moments} holds with $C_p$ in
place of $C_{p,T}$ for $W\ge W_0(T)$ and $1\le r\le TW$. In particular
$\E g_r^p\le C_p(1+\varepsilon^{-p})$ for $\varepsilon W\le r\le TW$.
\end{theorem}

The $-P_{k+\ell+1}$ term retains the restriction $T\preceq I$.
A covariance computed with a Hermitian one-port input and then
used without correction would miss this term.

 In Theorem~\ref{thm:matched-traces} we grew
the path from its centre, adding one block at each end, and showed that
the normalized port trace stays within $O(W^{-1})$ of the matched value
$a$. Here we need the same statement for a path grown from one end
only. We grow at the near end, inserting each fresh block between the old path and the measuring port, so the
lengths $r=1,2,\ldots$ form a filtration, and $R_r$ is the resolvent at
the near port. The following lemma has two parts. The first is the same $O(W^{-1})$
bound, uniformly over all lengths up to $TW$. The second is that it
still holds if every physical block carries a small deterministic
absorption, that is, a positive imaginary part of size at most $K/W^2$.

\begin{lemma}[All-length zeroth scalar estimate in passive coordinates]
\label{lem:t26-passive-initial-scalar}
At deterministic $E\in I$, for the actual two-port Gaussian path and
every fixed $p,T$, the near-end filtration process satisfies
\[
 \sup_{1\le r\le TW}\|\delta_{0,r}\|_p\le C_{p,T}/W,
 \qquad \delta_{0,r}=|W^{-1}\Tr R_r-a|.
\]
For any fixed scalar-domain radius $d_*>0$ and $D>0$, the probability
of leaving this domain at any $r\le TW$ is $O_D(W^{-D})$. Without absorption, more
strongly, for every fixed $d_1>0$ there are $c>0$ and $W_0$ such that, for $W\ge W_0$,
\begin{equation}\label{eq:passive-superexp-exit}
 \Prob\Bigl(\max_{1\le r\le TW}\delta_{0,r}>d_1\Bigr)\le e^{-cW\log W}.
\end{equation}
The same
fixed-moment conclusions hold if the physical blocks carry
deterministic nonnegative absorptions at most $K/W^2$. Constants may
depend on $K,p,T,b,I$, and $W$ is sufficiently large for the chosen
fixed moments. %
Furthermore, without absorption the first bound is length-uniform, i.e., for every
$p$ there is $C_p=C_p(b,I)$, and then for every $T$ a $W_0(T)$, with
$\sup_{1\le r\le TW}\|\delta_{0,r}\|_p\le C_p/W$ for $W\ge W_0(T)$.
\end{lemma}

\begin{figure}[ht]
\centering
\begin{tikzpicture}[x=1cm,y=1cm]
 \begin{scope}
  \node[font=\small,anchor=west] at (-1.3,1.25) {centre growth};
  \foreach \x/\f in {0/black!15,1/black!45,2/black!45,3/black!45,4/black!15}
    \draw[fill=\f] (\x,0) rectangle ++(.8,.8);
  \draw[thick,dashed] (-1.1,0) rectangle ++(.8,.8);
  \draw[thick,dashed] (5.1,0) rectangle ++(.8,.8);
  \draw[->,thick] (.4,1.0)--(.4,.85);
  \draw[->,thick] (4.4,1.0)--(4.4,.85);
 \end{scope}
 \begin{scope}[yshift=-2.1cm]
  \node[font=\small,anchor=west] at (-1.3,1.25) {near-end growth};
  \foreach \x/\f in {0/black!45,1/black!45,2/black!45,3/black!45,4/black!15}
    \draw[fill=\f] (\x,0) rectangle ++(.8,.8);
  \draw[very thick] (-1.1,0) rectangle ++(.8,.8);
  \draw[thick,dashed] (5.1,0) rectangle ++(.8,.8);
  \node[font=\small] at (-.7,.4) {$R_r$};
  \draw[->,thick] (4.4,1.0)--(4.4,.85);
  \node[font=\small,below] at (-.7,-.05) {near port};
  \node[font=\small,below] at (5.5,-.05) {far port};
 \end{scope}
\end{tikzpicture}
\caption{Two ways to build the same path, used in the proof of
Lemma~\ref{lem:t26-passive-initial-scalar}. Dark blocks are the old
path and light blocks are the fresh ones. Above, as in
\eqref{eq:two-ended-recursion}, one block is added at each end. Below,
the near port (solid) is kept fixed and each fresh block is inserted at the near end, between the old path and the measuring port. At each fixed length the resulting path consists of the same
independent Gaussian blocks and bonds, so the two constructions have the
same law. Bounds on the law at a single length, such as those of
Theorem~\ref{thm:matched-traces}, therefore transfer from one
construction to the other.}
\label{fig:ps-growth}
\end{figure}
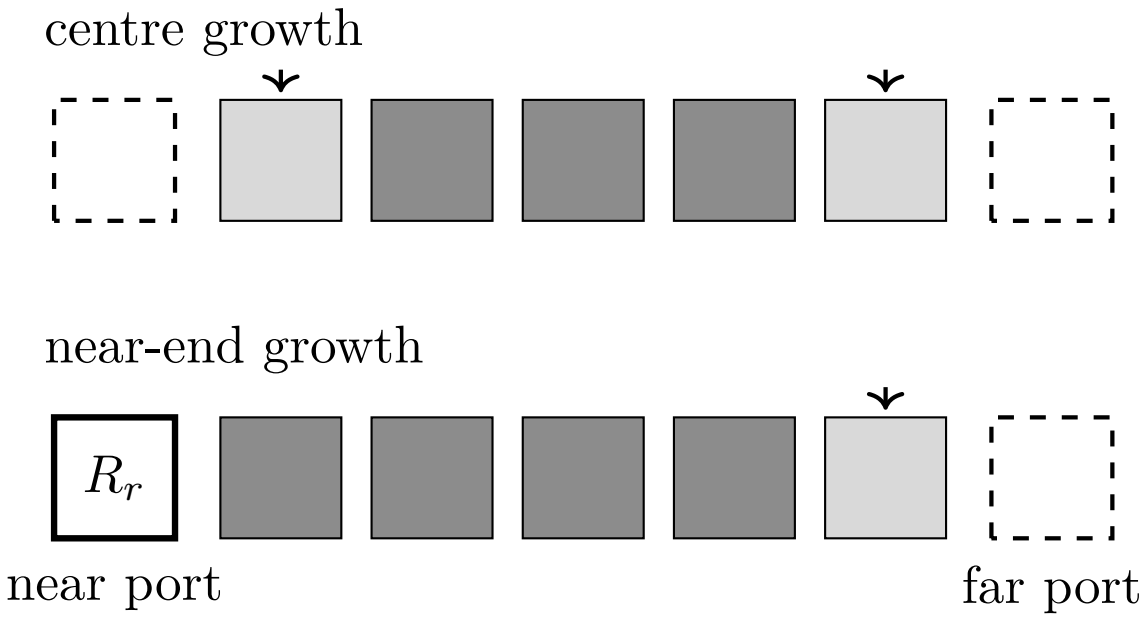

\begin{proof} Without absorption, the path grown from one end has,
at each length, the same law as the path grown from its center, as in
Figure~\ref{fig:ps-growth}, so the bound follows from
Theorem~\ref{thm:matched-traces} (Step~1). With a small absorption on the
blocks, we go through the proof of that theorem again and check that each
ingredient survives: the fixed starting core (Step~2), one extension step,
including a bound on resolvent entries that does not need symmetry (Step~3),
and the summation along the path (Step~4).

\emph{Step 1: the bound without absorption.} We first identify $R_r$
with the object studied in Theorem~\ref{thm:matched-traces}. Let
$\mathbb X_r$ be the Schur complement term of the path of length $r$ at
its two ports. Eliminating the far port by the Schur complement formula
gives exactly
\[
 R_r=\bigl[(z-\mathbb X_r)^{-1}\bigr]_{11},
\]
the block of the two-port resolvent at the near port. Growing at the
near end produces, at each length, the same marginal law of the original
path as growing from the centre in \eqref{eq:two-ended-recursion}. The
all-length two-port estimate of Theorem~\ref{thm:matched-traces}
therefore applies at every $r$, and gives the first claim. Here the
supremum over $r$ is taken outside the $L^p$ norm, so only the law at
each single length matters.

For the exit statement we use Markov's inequality with a large fixed
moment $p$, and then a union bound over the at most $TW$ lengths:
\[
 \Prob\Bigl(\max_{r\le TW}\delta_{0,r}>d_*\Bigr)
 \le\sum_{r\le TW}\frac{\|\delta_{0,r}\|_p^p}{d_*^p}
 \le TW\Bigl(\frac{C_{p,T}}{d_*W}\Bigr)^p.
\]
Taking $p$ large makes this $O_D(W^{-D})$. Only the marginal bounds
enter, and we need no independence between the events at different
lengths.

For \eqref{eq:passive-superexp-exit} we take instead the moment
$p=\lfloor\epsilon_0W\rfloor$, the largest allowed in
Theorem~\ref{thm:matched-traces}. Its bound $C\sqrt p/W$ at each length,
Markov's inequality and the same union bound give
\[
 \Prob\Bigl(\max_{r\le TW}\delta_{0,r}>d_1\Bigr)
 \le TW\Bigl(\frac{C\sqrt{\epsilon_0}}{d_1\sqrt W}\Bigr)^{\lfloor\epsilon_0W\rfloor}
 \le e^{-cW\log W}
\]
for $W$ large.
\emph{Step 2: absorption in a fixed core.} We now add a physical
diagonal absorption $iD$ with $0\preceq D\preceq KW^{-2}I$, and follow
the proof of Theorem~\ref{thm:matched-traces} step by step. That proof
starts with a fixed core. By the resolvent identity, adding $iD$ changes
the core resolvent by a product of two resolvents and $iD$. The core
inverse moments, and $\|D\|\le KW^{-2}$, therefore give an
$O_{L^p}(W^{-2})$ change of the port trace. The inverse moments remain
valid with the extra positive imaginary part, by the Poisson argument
of Lemma~\ref{lem:passive-bounds}.

\emph{Step 3: absorption in one extension step.} The proof of
Theorem~\ref{thm:matched-traces} then extends the path by the
integration-by-parts equations of one layer. With absorption, these
equations acquire the additional term $iD$, of norm at most $KW^{-2}$. So
their residual and their outgoing error increase by $O(W^{-2})$.

The stability argument also needs a bound on the block entries of the
resolvent, and it must now hold for a dissipative old Schur term
$\mathbb X$, meaning that $Y=-\Im\mathbb X\succeq0$. Then the resolvent
need not be normal, so we derive the bound directly. Write $R=(z-\mathbb X)^{-1}$ and
$\eta=\Im z$. Since $R^{-*}-R^{-1}=-2i(\eta I+Y)$, we have the two
identities
\[
 -\Im R=R^*(\eta I+Y)R=R(\eta I+Y)R^* .
\]
Both right sides are at least $\eta R^*R$ and $\eta RR^*$, respectively.
Taking the $(i,i)$ block of each, we get
\[
 \eta\sum_j R_{ji}^*R_{ji}\preceq-\Im R_{ii},\qquad
 \eta\sum_j R_{ij}R_{ij}^*\preceq-\Im R_{ii},
\]
which control the block columns and the block rows of $R$. By the
Cauchy--Schwarz inequality, first for the Hilbert--Schmidt pairing and
then for the sum over $j$, this gives
\[
 \sum_j|W^{-1}\Tr R_{ij}R_{ji}|
 \le-\Im(W^{-1}\Tr R_{ii})/\eta.
\]
This is the same local bound that was used in \eqref{eq:strict-contraction}.

The imaginary-part comparison that selects the correct branch gains a
nonnegative absorption term. So its substochastic row bounds remain
valid. Hence the stopped scalar recursion of
Theorem~\ref{thm:matched-traces} keeps its form:
\[
 d_+\le q_*d+CW^{-2}+|\xi_+|,\qquad \|\xi_+\|_p\le C_p/W,\qquad q_*<1 .
\]

\emph{Step 4: conclusion.} The rest of the proof of
Theorem~\ref{thm:matched-traces} goes through unchanged. The geometric
sum over the recursion gives the $O(W^{-1})$ bound. The two perturbed
cores of Step~2 start both parities. The deterministic port bound
$\|R\|\le\eta^{-1}$ controls the exit event, exactly as in that proof,
and this removes the stops. This proves the assertion with absorption. Taking $D=0$ in Steps~2--4 gives the length-uniform bound
without absorption: the core bound and the one-step constants
$q_*,C,C_p$ do not involve $T$, and the geometric sum gives $C_p/W$ at
every length. The total length enters only through the union bound
over at most $TW$ lengths in the exit estimate, absorbed by
$W\ge W_0(T)$. The same holds for the length-uniform clauses below:
their one-layer constants are chosen first, and $T$ enters only through
exit probabilities over at most $TW$ layers.
\end{proof}

 \subsection{The scalar derivative and a contracting source}\label{subsec:s13-the-scalar-derivative-and-a-contracting-}

 We identify which boundary response loses memory of its incoming
value. The derivative of the scalar mean map has a stable fixed point,
and the source estimates below quantify its contraction. We compute
the derivative map, choose its contracting solution, and control the
finite-width scalar and first-derivative errors after the entrance
interval.

 Use $c=1-2b$, $d_*=1/(1-b)$, $k_-=c+b^2a^2$.
The scalar mean map at a boundary parameter $w$ is
$v=f(\zeta)$, $u=(w-bv)^{-1}$,
$\zeta=-E-cv-bu$. Differentiation at $w=z$ gives
\begin{equation}\label{eq:first-jet-map-paper}
 u'=-u^2\frac{1+cf'(\zeta)}{1+(c+b^2u^2)f'(\zeta)},
 \qquad
 \zeta'=\frac{bu^2}{1+(c+b^2u^2)f'(\zeta)}.
\end{equation}
At the matched value $u=v=a$, the first derivative update is
$s\mapsto-a^2(1+cs)/(1+k_-s)$. Its derivative has modulus
$b^2/|1+k_-s|^2<1$ on a fixed neighborhood of
$|s|\le d_*$. Indeed
$|k_-|\le c+b^2=(1-b)^2$ and
$b/(1-|k_-|d_*)\le1$, with strict inequality uniformly on compact
bulk sets because $a^2\ne1$ there and $b,c>0$.
There is a unique stationary derivative in this disk, characterized by
\begin{equation}\label{eq:stationary-jet}
 k_-s_*^2+(1+ca^2)s_*+a^2=0,\quad |s_*|\le d_*,
 \qquad \theta_*=\frac{ba^2}{1+k_-s_*},\quad|\theta_*|<1.
\end{equation}
The denominator stays bounded away from zero.
The first-jet contraction together with
Theorem~\ref{thm:matched-traces} gives
\begin{equation}\label{eq:passive-jet-moments}
 \norm{\delta_0}_{L^p}+\norm{\delta_1}_{L^p}\le C_p/W,\qquad
 \delta_0=|W^{-1}\Tr R-a|,\quad
 \delta_1=|\partial_z(W^{-1}\Tr R)-s_*|,
\end{equation}
after $C_p\log W$ layers. Here the derivative holds the far absorber
fixed. To justify the same argument for nonnormal $X$, note that
$W^{-1}\Tr(w-X)^{-1}$ maps the upper half-plane into the lower
half-plane and has mass one at infinity. Its Herglotz representation
therefore has the same derivative bounds as a probability Stieltjes
transform. Cauchy's formula on a fixed small disk about $z$ gives
fresh derivative errors $O_{L^p}(W^{-1})$ and conditional means
$O(W^{-2})$. Iteration of the contraction proves
\eqref{eq:passive-jet-moments}; fixed high moments pay exits over
$O(W)$ layers. Higher fixed derivatives follow from the triangular
coefficient $\theta_*^{r+1}$ for the highest incoming derivative. Likewise \eqref{eq:passive-jet-moments} is length-uniform: for
$\varepsilon>0$ and $p$ there is $C_{p,\varepsilon}$, uniform for $E$ in
a compact bulk set, and then for every $T$ a $W_0(T)$, such that it
holds with $C_{p,\varepsilon}$ at every length
$\varepsilon W\le r\le TW$ when $W\ge W_0(T)$. Indeed the contraction
and the fresh errors are one-layer bounds, $\delta_0$ is controlled by
the length-uniform clause of Lemma~\ref{lem:t26-passive-initial-scalar},
and $C_p\log W\le\varepsilon W$ for large $W$.

\subsection{Passive inverse bounds and the mass mean}\label{subsec:s13-passive-inverse-bounds-and-the-mass-mean}

 We propagate transmission while preserving its positivity and upper
bound. Lemma~\ref{lem:passive-bounds} controls the inverse factors, and
Lemma~\ref{lem:passive-source-mixing} controls the weighted response that
must relax. We derive the exact update, bound its conditional noise,
and use the matched mean equations to identify the mass drift.

 For independent fresh blocks $A,B$ of variances $c/W,b/W$, put
\[
 \mathcal G=
 \begin{pmatrix}A-E-X&B\\B^*&zI\end{pmatrix}^{-1},
 \quad X_+=B^*(A-E-X)^{-1}B,\quad L=K\mathcal G_{12}.
\]
Inverse identities give exactly
\begin{equation}\label{eq:passive-update}
 -\Im X_+=B^*(A-E-X)^{-*}Y(A-E-X)^{-1}B,\qquad Q_+=L^*QL.
\end{equation}

\begin{lemma}[Bounds for one fresh layer with a dissipative input]
\label{lem:passive-bounds}
Let the old path enter the fresh layer through an independent matrix $X$
that is dissipative, that is, $-\Im X\succeq0$. With the layer inverse
$\mathcal G$, the weight $K=z-X$, and $L=K\mathcal G_{12}$, the following
hold for $W\ge W_0(p)$, with constants uniform in $X$.
\begin{enumerate}[label=(\roman*)]
\item Every fixed operator moment of $\mathcal G$,
$K\mathcal G_{11}$ and $L$ is bounded.
\item The covariance between a normalized trace of a
fixed word in these inverse factors and another fixed word is
$O(W^{-2})$.
\item For deterministic vectors $u,v$, the
centered bilinear forms $u^*LL^*v$ and
$u^*L\bigl(R_+-(2i\eta)^{-1}I\bigr)L^*v$, where $R_+=\mathcal G_{22}$, have
conditional $L^p$ norm at most $C_pW^{-1/2}\norm u\norm v$.
\end{enumerate}
\end{lemma}

The dissipation of the old input is what makes these bounds possible: it
gives bounded inverse factors and an upper bound on the transmission.
All the bounds hold conditionally on the old path, before the fresh
Gaussian layer is integrated.

\begin{proof}
 For a Hermitian input, the moment bounds were
proved in Lemma~\ref{lem:gaussian-absorption}. We reduce the dissipative
case to the Hermitian one by an analytic continuation in the strength of
the absorbing part, using the Poisson representation of a bounded
analytic function (Step~1). The weighted factors are then reduced to the
unweighted inverse by exact block identities (Step~2). For the
covariance and bilinear bounds we need the Gaussian calculus of
Lemma~\ref{lem:trace-matrix}, which requires smooth bounded functions.
We therefore add a small extra absorption, which makes every inverse
bounded and smooth (Step~3), and check that all the words we need keep
their weights in the right order (Step~4). We then prove the covariance
bound (Step~5) and the bilinear bound (Step~6) at positive absorption,
uniformly in it, and remove the absorption at the end.

\emph{Step 1: the unweighted inverse.} Write $X=X_h-iY$ with $X_h$
Hermitian and $Y\succeq0$. For a complex parameter $\zeta$ with
$\Im\zeta\ge0$ consider
\[
 F(\zeta)=\bigl(A-E-X_h+\zeta Y-z^{-1}BB^*\bigr)^{-1},
\]
so that $F(i)$ is the top block $\mathcal G_{11}$ at the actual input
$X$. Since $\Im(-z^{-1})=\eta/|z|^2$, for every vector $x$
\[
 \Im x^*\bigl(A-E-X_h+\zeta Y-z^{-1}BB^*\bigr)x
 \ge\frac{\eta}{|z|^2}\norm{B^*x}^2
 \ge\frac{\eta}{|z|^2\norm{B^{-1}}^2}\norm x^2 .
\]
Hence, for almost every fresh sample, $F$ is holomorphic in the open
upper half-plane, continuous and bounded on the closed half-plane,
uniformly in $\zeta$ and $Y$. A bounded analytic function is the Poisson
integral of its boundary values, so
$F(i)=\int_\R F(t)\,dt/(\pi(1+t^2))$. The triangle inequality and
Jensen's inequality then give
\[
 \norm{F(i)}^p\le\int_\R\norm{F(t)}^p\,\frac{dt}{\pi(1+t^2)},
 \qquad p\ge1 .
\]
For real $t$, $F(t)$ is the block $\mathcal G_{11}$ at the Hermitian
input $X_h-tY$. By Tonelli's theorem we may take the Gaussian
expectation inside the integral, and Lemma~\ref{lem:gaussian-absorption},
which is uniform in the Hermitian input, bounds it at each $t$. This
proves the bound on the moments of $\mathcal G_{11}$, uniformly in $Y$,
for $1\le p\le\epsilon_1W$.

\emph{Step 2: the full block and the weighted factors.} Schur inversion
expresses the other blocks of $\mathcal G$ through $\mathcal G_{11}$ and
the fresh bond, which gives the bound for the full block. For the
weighted factors we use the two identities
\[
 K\mathcal G_{11}
 =I+[(z+E)I-A+z^{-1}BB^*]\mathcal G_{11},
 \qquad
 \mathcal G_{12}=-z^{-1}\mathcal G_{11}B .
\]
The second gives $L=K\mathcal G_{12}=-z^{-1}K\mathcal G_{11}B$. So the
weighted factors are products of $\mathcal G_{11}$ with polynomials in
$A$ and $B$, and H\"older's inequality together with the Gaussian moments
of $\norm A$ and $\norm B$ proves their moment bounds. The adjoints have
the same order. This proves (i).

\emph{Step 3: regularization by an extra absorption.} For the covariance
and bilinear bounds we regularize, so that the Gaussian calculus of
Lemma~\ref{lem:trace-matrix} applies without any cutoff. For
$\varepsilon>0$ put $X_\varepsilon=X-i\varepsilon I$, which is again
dissipative, and write $\mathcal G_\varepsilon$,
$K_\varepsilon=z-X_\varepsilon$ and
$L_\varepsilon=K_\varepsilon\mathcal G_{\varepsilon,12}$ for the
corresponding objects. In the imaginary part of the full layer matrix the
bond blocks cancel, because they enter as a Hermitian pair:
\[
 \Im x^*\mathcal G_\varepsilon^{-1}x
 =x_1^*(Y+\varepsilon I)x_1+\eta\norm{x_2}^2
 \ge\min(\varepsilon,\eta)\norm x^2,
 \qquad x=(x_1,x_2),
\]
for every Hermitian $A$ and every $B$. Hence
$\norm{\mathcal G_\varepsilon}\le\min(\varepsilon,\eta)^{-1}$ everywhere,
$\mathcal G_\varepsilon$ is a smooth function of the real Gaussian
coordinates, and
\[
 d\mathcal G_\varepsilon=-\mathcal G_\varepsilon(d\mathcal H)\mathcal G_\varepsilon,
\]
where $d\mathcal H$ is the Hermitian differential of
$\bigl(\begin{smallmatrix}A&B\\B^*&0\end{smallmatrix}\bigr)$. The moment
bounds of Steps~1 and~2 apply to the dissipative input $X_\varepsilon$, so
the quantity
\[
 \Lambda_\varepsilon=1+\norm{\mathcal G_\varepsilon}
 +\norm{K_\varepsilon\mathcal G_{\varepsilon,11}}+\norm{L_\varepsilon}
\]
has every fixed moment bounded, uniformly in $\varepsilon>0$, in $X$ and in
$W\ge W_0(p)$. We write $\Lambda$ for the same quantity at a general
dissipative input.

\emph{Step 4: the structure of the words.} Let
$\iota_1,\iota_2:\C^W\to\C^{2W}$ be the embeddings of the physical block
and of the port. Every letter of the words we use has the form
$P\mathcal GQ$, with
\[
 P\in\{\iota_1^*,\iota_2^*,K\iota_1^*\},\qquad Q\in\{\iota_1,\iota_2\};
\]
for example $L=K\iota_1^*\mathcal G\iota_2$. The differential of such a
letter is $-(P\mathcal G)(d\mathcal H)(\mathcal GQ)$, and that of its
adjoint is $-(\mathcal GQ)^*(d\mathcal H)(P\mathcal G)^*$. So the
possibly large weight $K$ enters only through the row
\[
 K\iota_1^*\mathcal G=[K\mathcal G_{11},\,L],
\]
every outside weighted inverse stays in the displayed order, and
$P\mathcal G$, $\mathcal GQ$ and every letter have norm at most
$\Lambda$. Consequently, a word of length $\ell$ whose deterministic
insertions have norm at most $k_0$ has norm at most
$k_0^{\ell+1}\Lambda^\ell$, and its differential is a sum of $\ell$ terms
$a_m(d\mathcal H)c_m$ with
\[
 \norm{a_m}\norm{c_m}\le k_0^{\ell+1}\Lambda^{\ell+1}.
\]

\emph{Step 5: the covariance bound (ii).} Fix $\varepsilon>0$. Let
$f_\varepsilon$ be a normalized trace of a word of length $\ell_f$, and
$U_\varepsilon$ a word of length $\ell_U$, both at the input
$X_\varepsilon$. By Step~3 they are bounded smooth functions with bounded
derivatives, so the Gaussian interpolation identity in the proof of
Lemma~\ref{lem:trace-matrix} applies. It writes the covariance as an
integral of the contraction of $df_\varepsilon$ with $dU_\varepsilon$, at
the sample $x$ and at an interpolated sample $x_t$. By Step~4 this
contraction is
\[
 W^{-2}\sum_{m,m'}a'_{m'}\,\mathcal V(c_ma_m)\,c'_{m'},
 \qquad \norm{\mathcal V(N)}\le C\norm N,
\]
where $\mathcal V$ is the variance map of the layer. The factors are
evaluated at $x$ and at $x_t$, which have the same law; the
arithmetic--geometric mean inequality separates them, and we obtain
\[
 \norm{\E(f_\varepsilon U_\varepsilon)-\E f_\varepsilon\,\E U_\varepsilon}
 \le C\ell_f\ell_U\,k_0^{n}\,W^{-2}\,\E\Lambda_\varepsilon^{n},
 \qquad n=\ell_f+\ell_U+2 .
\]
The right side is bounded uniformly in $\varepsilon$, by Step~3.

It remains to let $\varepsilon\downarrow0$. Outside a null set the layer
matrix at $X$ is invertible, so each word converges as
$\varepsilon\downarrow0$. The words and their products are dominated by
powers of $\Lambda_\varepsilon$, which are bounded in $L^2$ uniformly in
$\varepsilon$, so the family is uniformly integrable. Vitali's theorem
therefore passes the $W^{-2}$ covariance bound to the limit. This proves
(ii).

\emph{Step 6: the bilinear bound (iii) and the trace fluctuation.} Let
$\psi_\varepsilon$ be either $u^*L_\varepsilon L_\varepsilon^*v$ or
$u^*L_\varepsilon(\mathcal G_{\varepsilon,22}-c_0I)L_\varepsilon^*v$, with
$c_0=(2i\eta)^{-1}$; Schur inversion gives $\mathcal G_{22}=R_+$. By
Step~4, each term of $d\psi_\varepsilon$ has the form $x^*(d\mathcal H)y$
with
\[
 \norm x\norm y\le\Lambda_\varepsilon^3(\Lambda_\varepsilon+|c_0|)\norm u\norm v .
\]
Since the entry variances are $O(W^{-1})$, the Gaussian gradient,
weighted by the covariance of the entries, satisfies
\[
 |\nabla\psi_\varepsilon|
 \le CW^{-1/2}\Lambda_\varepsilon^3(\Lambda_\varepsilon+|c_0|)\norm u\norm v .
\]
The dimension-free Gaussian $L^p$ gradient inequality
(Remark~\ref{rem:poincare}) and H\"older's inequality then give
\[
 \norm{\psi_\varepsilon-\E\psi_\varepsilon}_{L^p}\le C_pW^{-1/2}\norm u\norm v,
\]
uniformly in $\varepsilon$. As $\varepsilon\downarrow0$,
$\psi_\varepsilon\to\psi$ almost surely, and
$|\psi_\varepsilon|\le\Lambda_\varepsilon^2(\Lambda_\varepsilon+|c_0|)\norm u\norm v$
is bounded in $L^2$ uniformly in $\varepsilon$. Vitali's theorem gives
$\E\psi_\varepsilon\to\E\psi$, and Fatou's lemma then gives the bilinear
estimate for $\psi$. This proves (iii).

The same argument also controls the normalized port trace. From
$\norm M_{\HS}^2\le2W\norm M^2$ we get
\[
 \bigl|\nabla W^{-1}\Tr\mathcal G_{\varepsilon,22}\bigr|
 \le CW^{-1}\norm{\mathcal G_\varepsilon}^2,
\]
and hence
\[
 \norm{W^{-1}\Tr R_+-\E W^{-1}\Tr R_+}_{L^p}\le C_p/W
\]
for each fixed $p$ and every dissipative $X$.

Finally, a derivative in the boundary parameter inserts bounded
projections between inverse factors, so the whole argument also applies
on a fixed disk of boundary parameters.
\end{proof}

On a small fixed scalar domain let
\[
 M=(z-X-\Omega)^{-1},\quad m=W^{-1}\Tr M,\quad
 v=(z-bm)^{-1},\quad \Omega=c(m-a)+b(v-a).
\]
These equations select the stable mean branch, with
$m,v=a+O(\delta_0)$. Put $h=|v|^2$, $k_h=c+b^2h$,
$A_M=W^{-1}\Tr MM^*$ and $y_M=W^{-1}\Tr M^*YM$.
The advanced/retarded trace equations for a bottom source are
$t_1=A_M(ct_1+bt_2)+O(W^{-2})$ and
$t_2=h(1+bt_1)+O(W^{-2})$. Consequently
\begin{equation}\label{eq:passive-mean}
 \E(LL^*\mid X)=\frac{bh}{1-k_hA_M}(KM)(KM)^*
                                      +O_{\rm op}(W^{-2}).
\end{equation}
The outside $K,K^*$ are kept inside the covariance calculation;
one never multiplies an unweighted error by an unbounded $K$.
The imaginary parts of the scalar equations give
\[
 -\Im m=(\eta-\Im\Omega)A_M+y_M,\quad
 -\Im v=h(\eta+b(-\Im m)),\quad
 \eta-\Im\Omega=c(-\Im m)+b(-\Im v).
\]
Eliminating $m,v$ yields the exact identity
\[
 (\eta-\Im\Omega)(1-k_hA_M)=bh\eta+k_hy_M.
\]
Let $c_0=(2i\eta)^{-1}$, $f=R-c_0I$ and
$\alpha=(1-b)^2/(b\eta)$. Expanding
$M=R(I-\Omega R)^{-1}$ gives
$y_M=(q/W)(1+O(\delta_0))$ with a relative error, by
Hilbert--Schmidt sandwiching with $Y^{1/2}R$.
Thus
\begin{equation}\label{eq:passive-cancelled-mean}
 \E(LL^*\mid X)=
 \frac{I+\Omega f+\bar\Omega f^*}{1+\alpha q/W}
 +O_{\rm op}(\delta_0^2+\delta_0q/W+W^{-2}).
\end{equation}
The linear term from $1-\Im\Omega/\eta$ is precisely the
$-c_0$ in $f$.

\begin{lemma}\label{lem:passive-source-mixing}
For $0<\epsilon<T$ and every fixed $p$,
\begin{equation}\label{eq:passive-source-mixing}
 \sup_{\epsilon W\le r\le TW}
 \left(\norm{\Tr(Q_rf_r)}_{L^p}
 +\norm{\norm{Q_r^{1/2}f_rQ_r^{1/2}}_{\HS}}_{L^p}\right)
 \le C_{p,\epsilon,T}W^{-1/2}.
\end{equation}
The conductance moment bound in \eqref{eq:conductance-moments}
holds already at finite $W$.
\end{lemma}

The scalar mean contracts while the positive source weight is retained. Summing the resulting errors is possible at the critical number of layers.

\begin{proof}
Positivity gives $q_r\le W/(4\eta)$. Equations
\eqref{eq:passive-mean}--\eqref{eq:passive-cancelled-mean} and the
bilinear fluctuation bound give
\[
 \E(q_+\mid X)\le\frac{q}{1+\alpha q/W}
                   +C\delta_0q+CW^{-2}q,\qquad
 \norm{q_+-\E(q_+\mid X)}_{L^p\mid X}\le C_pq/\sqrt W.
\]
Expand a fixed integer power around its conditional mean.
The centered term starts at order two; high fixed moments pay
the remainder. Young's inequalities
$\delta_0q^p\le\varepsilon q^{p+1}/W+
C_{p,\varepsilon}(W\delta_0)^{p+1}/W$ and
$q^p\le\varepsilon q^{p+1}+C_{p,\varepsilon}$ yield
\[
 x_{r+1}\le x_r-c_px_r^{1+1/p}/W+C_p/W,\qquad x_r=\E q_r^p.
\]
Comparison with $C[1+(W/(r+1))^p]$, starting above the deterministic
initial bound, proves the claim.  The constants $c_p,C_p$ of this recursion come from one fresh layer at a deterministic passive input and from the moments of $W\delta_0$, which are length-uniform by Lemma~\ref{lem:t26-passive-initial-scalar}. Hence $C$ does not depend on $T$, which gives the length-uniform clause after \eqref{eq:conductance-moments}. Reducing $c_p$ makes the comparison
map increasing on the whole permitted range; interpolation covers
noninteger $p$.

Hold $X$ fixed and differentiate only the near-port parameter.
The exact identities
\[
 \partial_zL=RL-LR_+,\qquad
 Lf_+L^*=fLL^*-\partial_z(LL^*)
\]
retain the nonnormal order. In this differentiation the full
definition is $\Omega(w)=w+E+cm(w)+bv(w)$.
At the fixed jets $\partial_z\Omega=1-\theta_*$.
Differentiate \eqref{eq:passive-mean} and its scalar identity:
\[
 \E(Lf_+L^*\mid X)=\frac{\theta_*}{1+\alpha q/W}f
             +O_{\rm op}(\delta_0+\delta_1+q/W+W^{-2}).
\]
For example $|\partial_zy_M|\le Cq/W$ follows from
$M=R(I-\Omega R)^{-1}$; it does not require a bound on $Y$.
The coefficient has modulus less than one. Pairing this equation
with $Q$, and then with two factors $Q^{1/2}$, gives geometric
recursions for the two quantities in
\eqref{eq:passive-source-mixing}. Their centered forcing is
$C_pq/\sqrt W$ and their mean forcing is
$Cq(\delta_0+\delta_1+q/W+W^{-2})$.
Start at $\epsilon W/2$, where the just-proved moments are bounded,
and use \eqref{eq:passive-jet-moments}. Iteration proves the result;
arbitrary fixed failure powers and $Q\preceq(4\eta)^{-1}I$
pay jet exits.
\end{proof}

\subsection{The two-replica contraction}\label{subsec:s13-the-two-replica-contraction}

  We compute the noise that remains after the mean transmission has
been removed. Lemma~\ref{lem:passive-covariance} identifies its bilinear
covariance, including the terms that retain the old matrix response.
We contract two copies of a fresh layer in their exact matrix order,
control the residual terms, and use the resulting coefficients to
finish the transmission generator.

 Set $k=c+b^2$, $\sigma_+=k/b$ and $\sigma_-=k_-/(1+k_-s_*)$. Then $\sigma_+/(2\eta^2)=D_0$.

\begin{lemma}\label{lem:passive-covariance}
For $U=LL^*$, $U^0=U-\E(U\mid X)$ and four old vectors, the bilinear
covariance has crossed tensor
\begin{equation}\label{eq:passive-covariance}
 \frac{D_0}{2W}\{I\otimes(I-T)+(I-T)\otimes I\}
 +\frac{\sigma_-}{W}f\otimes f
 +\frac{\bar\sigma_-}{W}f^*\otimes f^*+\mathcal R.
\end{equation}
Here $A\otimes B$ denotes $(v_1^*Av_4)(v_3^*Bv_2)$, and
$|\mathcal R|\le C[W^{-2}+W^{-1}(\delta_0+\delta_1+q/W)]
\prod_j\norm{v_j}$.
\end{lemma}

\begin{remark} The matrix $U=LL^*$ is the transport of one
fresh layer, and its centered part $U^0$ is the noise that drives the
transmission diffusion. The lemma computes the covariance of two
bilinear forms of $U^0$, $v_1^*U^0v_2$ and $v_3^*U^0v_4$, up to errors of
order $W^{-2}$. The answer \eqref{eq:passive-covariance} has three parts.
The first, with $D_0$, is the opposite-sign term. It contains $I-T$
rather than $I$: this is the correction from the upper constraint
$T\preceq I$ on the transmission. The other two, with $\sigma_-$ and
$\bar\sigma_-$, are the same-sign terms. These coefficients become the
diffusion coefficients of the limiting generator.

\end{remark}

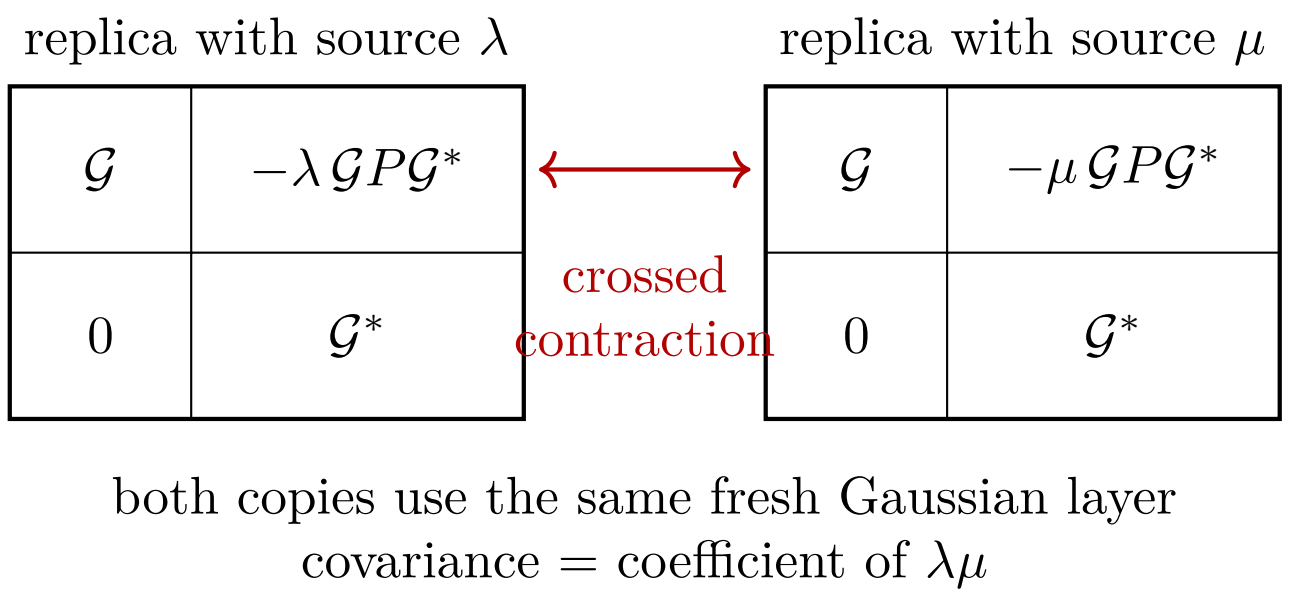
\begin{figure}[ht]
\centering
\begin{tikzpicture}[x=1cm,y=1cm]
 \foreach \xs/\src in {0/\lambda,5.0/\mu}{
  \begin{scope}[xshift=\xs cm]
   \draw[thick] (0,0) rectangle (3.4,2.2);
   \draw[thin] (1.2,0)--(1.2,2.2) (0,1.1)--(3.4,1.1);
   \node[font=\small] at (.6,1.65) {$\mathcal G$};
   \node[font=\small] at (2.3,.55) {$\mathcal G^*$};
   \node[font=\small] at (2.3,1.65) {$-\src\,\mathcal GP\mathcal G^*$};
   \node[font=\small] at (.6,.55) {$0$};
   \node[font=\small] at (1.7,2.5) {replica with source $\src$};
  \end{scope}
 }
 \draw[<->,thick,red!70!black] (3.5,1.65) -- (4.9,1.65);
 \node[font=\small,red!70!black,align=center] at (4.2,.75) {crossed\\ contraction};
 \node[font=\small,align=center] at (4.2,-.75)
   {both copies use the same fresh Gaussian layer\\
    covariance $=$ coefficient of $\lambda\mu$};
\end{tikzpicture}
\caption{The replica construction in the proof of
Lemma~\ref{lem:passive-covariance}. Each copy is the inverse of a
triangular block matrix whose diagonal blocks are $\mathcal G$ and
$\mathcal G^*$. A formal source on the bottom off-diagonal block makes
the top off-diagonal entry equal to $-\lambda\mathcal GP_{\rm bottom}\mathcal G^*$,
a product of the form that enters the bilinear forms of $U$. Both copies are built from
the same fresh Gaussian layer. Gaussian integration by parts couples
them through the crossed contraction (red), and the covariance of the
two bilinear forms is read off as the coefficient of $\lambda\mu$.}
\label{fig:pc-replica}
\end{figure}

\begin{proof}The covariance involves products of a resolvent
$\mathcal G$ and its adjoint $\mathcal G^*$. We compute it with a replica
trick, as in Figure~\ref{fig:pc-replica}. We build a triangular block
matrix whose diagonal blocks are $\mathcal G$ and $\mathcal G^*$ and whose
off-diagonal block carries a formal source $\lambda$. Then the product
$\mathcal GP\mathcal G^*$ appears as the coefficient of $\lambda$. Taking two
such copies, with sources $\lambda$ and $\mu$, and extracting the
coefficient of $\lambda\mu$ produces exactly the covariance we need.
Gaussian integration by parts turns this into a closed linear system,
the ladder, which we solve explicitly. The only subtle point is that the
old input is not Hermitian, so resolvents do not commute with
their adjoints. We reduce only scalar trace words, where the error of
this reduction is controlled.

\emph{Step 1: the replica construction and its errors.} We introduce two
triangular replica inverses. Each has diagonal blocks $\mathcal G$ and
$\mathcal G^*$, and a formal source on the bottom off-diagonal block,
equal to $\lambda$ in the first copy and to $\mu$ in the second. Their
off-diagonal entries are then $-\lambda\mathcal GP_{\rm bottom}\mathcal G^*$,
and its version with $\mu$. Both copies use the same fresh Gaussian
layer. We apply the Gaussian integration by parts \eqref{eq:gaussian-ibp}
to the covariance of the top entries, and we extract the coefficient of
$\lambda\mu$ before bounding any errors. The crossed contraction has the
source $cP_{\rm top}+bP_{\rm bottom}$, which weights the two blocks of the
layer by their variances.

Two kinds of error arise. A centered normalized self-energy costs
$W^{-2}$, by Lemma~\ref{lem:passive-bounds}. Replacing the two remaining
bilinear forms by their means costs $W^{-1}$, times the crossed
prefactor $W^{-1}$. Here $W^{-1}$ bounds the conditional covariance of bilinear forms
$u^*\mathcal Wv$ of general fixed inverse words $\mathcal W$ (products of
a fixed number of blocks of $\mathcal G$, $K\mathcal G_{11}$, $L$, their
adjoints and bounded deterministic matrices), by
$CW^{-1}\prod_j\norm{v_j}$. This follows as the bilinear bound in the
proof of Lemma~\ref{lem:passive-bounds}: after the regularization
$X-i\varepsilon I$, removed by Vitali's theorem as there,
$\partial\mathcal G=-\mathcal G(\partial\mathcal G^{-1})\mathcal G$ gives
$\Var(u^*\mathcal Wv\mid X)\le CW^{-1}\norm u^2\norm v^2$ by the
Gaussian Poincar\'e inequality (Remark~\ref{rem:poincare}), H\"older and
the operator moments of that lemma, and Cauchy--Schwarz gives the
covariance. All derivatives keep their outside factors
$K\mathcal G_{1j}$, and the adjoints of these, attached. So none of these
errors contains $\norm K$.

\emph{Step 2: the ladder for the remaining contractions.} The terms that
remain form a finite linear system. We write $\mathsf M_p^{(1)}$ and
$\mathsf M_p^{(2)}$ for the mean matrices of the two replicas on the
physical blocks $p=1,2$. At the stationary scalar point they are, with
source $\lambda$ or $\mu$,
\[
 \mathsf M_1=\begin{pmatrix}R&-\lambda RR^*\\0&R^*\end{pmatrix},
 \qquad
 \mathsf M_2=\begin{pmatrix}a&-\lambda d_*\\0&\bar a\end{pmatrix}.
\]
For each matrix unit $E_{rs}$ in replica space, the mixed ladder is
\begin{align*}
 \mathsf K_1&=\mathsf M_1^{(2)}
 [cE_{rs}I+c\av{\mathsf K_1}I+b\av{\mathsf K_2}I]
 \mathsf M_1^{(1)},\\
 \mathsf K_2&=\mathsf M_2^{(2)}
 [bE_{rs}I+b\av{\mathsf K_1}I]\mathsf M_2^{(1)}.
\end{align*}
Here the angle brackets mean $W^{-1}\Tr$, taken entrywise in replica
space. The structure is the same as in the two-resolvent system of
Lemma~\ref{lem:full-source-covariance}. The physical block feels itself
with weight $c$ and the port with weight $b$, and the port feels the
physical block with weight $b$. The denominators that arise are $b$,
$1+k_-s_*$ and their conjugates. They and their derivatives in the
sources are uniformly bounded.

\emph{Step 3: reducing the scalar trace words.} The old input is
passive, so $R$ and $R^*$ need not commute. We may therefore reduce only
scalar trace words, as if the input were Hermitian, and we must control
the error. With $\widetilde Q=RYR^*$ we have exactly
\[
 R^*R=c_0(R^*-R)-Q/\eta,\qquad
 RR^*=c_0(R^*-R)-\widetilde Q/\eta.
\]
For a Hermitian input the two positive errors $Q$ and $\widetilde Q$ would
vanish. Here they both have trace $q$, while every other resolvent is
bounded by $\eta^{-1}$. So commuting factors and using partial fractions
costs $O(q/W)$ in every normalized scalar word we use. Together with the
two jet errors $\delta_0$ and $\delta_1$, this gives
\[
 \av{RR^*}=d_*,\quad
 B_*:=\av{R^2R^*}=c_0(d_*+s_*),\quad
 C_*:=\av{R^2(R^*)^2}
                  =\frac{s_*+\bar s_*+2d_*}{4\eta^2},
\]
up to $O(\delta_0+\delta_1+q/W)$.

\emph{Step 4: solving the ladder.} Solving the replica ladder for each
matrix unit gives the following table. Its middle column is the top
$(1,2)$ entry after the outside weights, and its last column is the
other crossed product:
\[
 \begin{array}{c|c|c}
 E_{11}&\lambda(rI-\sigma_-R)&-\mu R\\
 E_{22}&\mu(\bar rI-\bar\sigma_-R^*)&-\lambda R^*\\
 E_{12}&\sigma_+I&\lambda\mu R^*R\\
 E_{21}&\lambda\mu(r_{12}I-rR-\bar rR^*
                                      +\sigma_+R^*R)&I
 \end{array}
\]
To verify the table explicitly, we record the scalar unknowns. They
satisfy
\begin{align*}
 \ell_*&=\frac{-\sigma_+B_*+bd_*as_*}{1+k_-s_*},
 \qquad r=k_-\ell_*-bd_*a,\\
 C_0&=bd_*^2-b^2d_*(a\ell_*+\bar a\bar\ell_*),
 \qquad r_{12}=k\ell_{12}+C_0,\\
 \ell_{12}&=\frac{d_*C_0-rB_*-\bar r\bar B_*+\sigma_+C_*}{b}.
\end{align*}
We obtain these equations successively, first from the diagonal replica
entries and then from the off-diagonal ones. For $E_{21}$ the induced
traces are $\mu\ell_*$, $\lambda\bar\ell_*$ and $\lambda\mu\ell_{12}$. Before
the outside weights are attached, its final top word is
\[
 \lambda\mu\bigl[r_{12}RR^*-rR^2R^*-\bar rR(R^*)^2+\sigma_+RR^*RR^*\bigr],
\]
and this gives the last row of the table.

\emph{Step 5: the scalar simplifications.} The two scalars $r$ and $r_{12}$
simplify to
\begin{equation}\label{eq:replica-scalars}
 r=c_0(\sigma_--\sigma_+),\qquad
 r_{12}=\frac{2\sigma_+-\sigma_--\bar\sigma_-}{4\eta^2}.
\end{equation}
For the first identity, we substitute $B_*=c_0(d_*+s_*)$ and obtain
\[
 r=\frac{-k_-\sigma_+B_*-bd_*a}{1+k_-s_*}.
\]
We then use the two relations
\[
 c_0(1-a^2)=d_*a,\qquad k_-(1+\sigma_+d_*)=\sigma_+-b(1-a^2).
\]
For the second identity, we substitute the first identity and its
conjugate into $C_0$ and $\ell_{12}$, and we replace $B_*$ and $C_*$ by
their values from Step~3. We then use $kd_*=1-b$, $k/b=\sigma_+$, and the
quadratic equation \eqref{eq:stationary-jet} and its conjugate.
Collecting the coefficients of $s_*$ and $\bar s_*$ leaves the second
identity. All these expressions are continuous at $k_-=0$.

\emph{Step 6: summation.} Summing the four crossed products of the table
gives
\begin{align*}
 &\sigma_+(I\otimes R^*R+R^*R\otimes I)
 -r(R\otimes I+I\otimes R)
 -\bar r(R^*\otimes I+I\otimes R^*)\\
 &\hspace{15mm}
 +\sigma_-R\otimes R+\bar\sigma_-R^*\otimes R^*
 +r_{12}I\otimes I.
\end{align*}
We now insert \eqref{eq:replica-scalars}, and we use the identity
\[
 S^*S=I+2i\eta(R^*-R)+4\eta^2R^*R,
\]
which follows from $S=I-2i\eta R$. The result is exactly the tensor in
\eqref{eq:passive-covariance}. We have not commuted any external
factors. The errors of the scalar reductions in Step~3 carry the crossed
prefactor $W^{-1}$, which gives the stated remainder.
\end{proof}

\begin{proof}[Proof of Theorem~\ref{thm:transmission-generator}]
Pair \eqref{eq:passive-cancelled-mean} with $Q$.
Its linear scalar error is
$2\Re(\Omega\Tr Qf)/(1+\alpha q/W)$, whose expected size is
$O(W^{-3/2})$ on compact positive times by
\eqref{eq:passive-jet-moments} and
\eqref{eq:passive-source-mixing}. The other paired errors are
$O(W^{-2})$ in expectation. Expanding the denominator costs
$Cq^3/W^2$. Therefore
\[
 g_{r+1}-g_r=-\frac{D_0}{2W}g_r^2+\Delta M_r+\mathcal R_r,
 \qquad \sum_{\epsilon W\le r\le TW}\E|\mathcal R_r|\to0,
\]
where $\norm{\Delta M_r}_{L^p\mid\mathcal F_r}\le C_pg_r/\sqrt W$;
we used $\alpha/(4\eta)=D_0/2$.

The eigenvalues of $Q_+$ are those of $Q^{1/2}LL^*Q^{1/2}$.
Expand traces about this Hermitian representative. Its centered
Hilbert--Schmidt $L^p$ increment is at most $C_pq/\sqrt W$,
by the four-vector bound in an eigenbasis of $Q$.
Every fixed third-order remainder is therefore summable on
positive-time intervals. Terms containing $f$ are bounded by
\eqref{eq:passive-source-mixing} and fixed powers of $q$.
Their first-order products with $\Omega$ are summable, as are
the two same-sign covariance terms. The remaining tensor, in
an eigenbasis $T=\diag(t_i)$, is
\[
 \E\,dT_{ii}=-\tfrac12D_0P_1t_i\,dt,\qquad
 dT_{ij}\,dT_{kl}
 =\tfrac12D_0t_it_j(2-t_i-t_j)\delta_{il}\delta_{jk}\,dt.
\]
The second differential of $\Tr T^k$ gives
\eqref{eq:transmission-drift}; multiplying the first differentials
of two powers gives \eqref{eq:transmission-bracket}.
For example $A_1=-D_0P_1^2/2$,
$A_2=D_0P_1^2-2D_0P_1P_2$.
Since $0\le P_k\le P_1$, the conductance moments give tightness,
vanishing maximum jumps and uniform integrability of all finite
drift and product martingales. Diagonal extraction proves the
joint assertion for all powers.

Finally $y(t)=\E P_1(t)$ obeys
$y'=-(D_0/2)\E P_1^2\le-(D_0/2)y^2$.
For $0<s<t$, $y(t)\le2/[D_0(t-s)]$.
Taking consistent limits on positive-time intervals and then
$s\downarrow0$ proves the last bound in
\eqref{eq:conductance-moments}.
\end{proof}
\section{A positive construction of the transmission law}
\label{sec:transmission-law}
A generator does not by itself determine a law, especially with
infinitely many channels and an entrance of infinite conductance. This
section constructs the unique law of the transmission process.
Theorem~\ref{thm:transmission-law} states that the transmission processes
converge on compact positive-time intervals to a unique continuous process.
Its determining transforms are expectations of a finite-dimensional dual
diffusion, which is killed when it reaches a corner.
In Subsection~\ref{subsec:s14-the-physical-entrance-and-microscopic-fo} we show that the conductance of the physical path
starts arbitrarily large on the critical time scale, and we derive the
microscopic forcing. In Subsection~\ref{subsec:s14-the-finite-dual-and-its-corner} we introduce the finite dual
diffusion and specify its behaviour at the corner
(Lemma~\ref{lem:t26-dual-boundary}). In Subsection~\ref{subsec:s14-determination-and-absence-of-small-chann} we show that the
dual determines the law and prove the theorem. In Subsection~\ref{subsec:s14-fractional-bounds-for-homogeneous-transm} we rule
out mass hidden in arbitrarily small channels, using fractional powers of the
transmission.

Let $\mu_W(t)=\sum_i t_i(t)\delta_{t_i(t)}$ be the transmission
trace measure. For ordered parameters $0<x_1<\cdots<x_p$ and
$0<y_1<\cdots<y_p<1$, define
\begin{equation}\label{eq:transmission-ratios}
 F_{\mathbf x,\mathbf y}(\mu)=
 \exp\!\left\{\int_{[0,1]}
 \frac{\sum_i[\log(1-y_it)-\log(1+x_it)]}{t}\,d\mu(t)\right\}.
\end{equation}
The integrand at zero is $-\sum_i(x_i+y_i)$.
For a finite matrix this is
$\prod_i\det(I-y_iT)/\det(I+x_iT)$ and lies in $[0,1]$.

\begin{theorem}\label{thm:transmission-law}
The actual transmission trace-measure processes converge on compact
positive-time intervals to a unique continuous process.
Its determining transforms are
\begin{equation}\label{eq:transmission-dual-law}
 \E F_{\mathbf x,\mathbf y}(\mu(t))=
 \E_{\mathbf x,\mathbf y}
 \left[e^{-\int_0^t h(t-s)R_s\,ds}\one_{\{\zeta\le t\}}\right].
\end{equation}
Here $h(s)=4\pi\rho(E)\kappa(s)$ if the physical blocks carry
absorption $\kappa(s)/W^2$, and $h=0$ otherwise.
The positive dual, its total coordinate $R_s$ and extinction time
$\zeta$ are constructed below. The following additional assertions
concern only homogeneous transport, $\kappa\equiv0$. For $0<p<1$,
\begin{equation}\label{eq:fractional-transmission}
 \limsup_{W\to\infty}\E\Tr T_{\lfloor Wt\rfloor}^{p}
                  \le C_p(1+t^{-1}).
\end{equation}
The limit has no trace mass at zero. For each $0<q<1$,
$\limsup_W\E\Tr T_{\lfloor Wt\rfloor}^q\le C_qe^{-c_qt}$
when $t\ge2$.
\end{theorem}

The parameters in this theorem describe transmission, not physical
eigenvalues. Its finite positive construction will be used as initial
data for the source evolution before the ring is closed.

Persistent absorption retains the bounded transform
\eqref{eq:transmission-dual-law}, whose integrand also determines
possible trace mass at zero. No fractional bound, exclusion of that
mass, or exponential decay is asserted for the forced law. Indeed one
positive injection gives $\E\Tr T_+^p\ge c_{p,\kappa}W^{1-2p}$
for $p<1/2$: restrict the injection matrix to a fixed operator-norm
cutoff retaining order-$W$ trace, and use $x^p\ge K^{p-1}x$ there.
Moreover the forced first-trace equation
$(\E g)'=h-(D_0/2)\E g^2$ precludes persistent exponential decay.
These observations do not assert an atom at zero in the limiting
forced trace measure.

 \subsection{The physical entrance and microscopic forcing}\label{subsec:s14-the-physical-entrance-and-microscopic-fo}

 We identify how the physical path enters the transmission limit.
The estimate below shows that its conductance starts arbitrarily large
on the critical time scale. We control a bounded reciprocal of the
conductance, then treat a small positive source to derive the forcing
term used in the finite-dimensional description.

 For any fixed $K>0$, the actual path obeys
\begin{equation}\label{eq:transmission-entrance}
 \limsup_W\Prob\{\min_{1\le r\le tW}g_r\le K\}
                         \le C(1+K)t,\qquad 0<t\le1 .
\end{equation}
To prove this use $q=g/(4\eta)$ and the bounded reciprocal
$r_K(q)=\min(q^{-1},K^{-1})$.
Use the all-length zeroth estimate of
Lemma~\ref{lem:t26-passive-initial-scalar}. Taking the trace against
$Q$ in \eqref{eq:passive-cancelled-mean} and diagonalizing $Q$ in the
bilinear estimate of Lemma~\ref{lem:passive-bounds} gives, conditionally
on the entire old-path filtration $\mathcal F$,
\[
 \mu=\E(q_+\mid\mathcal F)
 =\frac{q}{1+\alpha q/W}+O((\delta_0+W^{-2})q),\qquad
 \|q_+-\mu\|_{p\mid\mathcal F}\le C_pq/\sqrt W.
\]
These estimates use only $\delta_0$, with
$0<q\le W/(4\eta)$; the positive-time source-mixing and the
post-logarithmic first-jet estimates are not needed. Shrink the fixed
scalar domain to have $\mu\ge c_0q$. For $q\ge K$ they imply
 \[
 \E(r_K(q_+)\mid\mathcal F)
 \le q^{-1}+C/W+C\delta_0/q+C/(Wq)+C_{A,K}W^{-A}.
\]
Indeed Taylor expansion on $|q_+-\mu|\le\mu/2$ costs
$C(q_+-\mu)^2/q^3$; higher fixed centered moments pay the
complement and the removed linear term. Stop at the first crossing
below $K$. Until then $q^{-1}\le K^{-1}$, so summing and using
$\E\delta_0\le C/W$ gives
\[
 \E r_K(q_{\rm stopped})
        \le\E r_K(q_1)+C(1+K^{-1})t+o(1).
\]
The actual finite core has $q_1/W$ bounded below in probability.
For the initial bond $X=B^*z^{-1}B$, functional calculus gives
\[
 q=\sum_i\frac{\eta s_i(B)^2}{|z^2-s_i(B)^2|^2}
 \ge\frac{\eta\Tr B^*B}{(|z|^2+\|B\|^2)^2}.
\]
Gaussian norm and trace concentration make this at least $aW$ with
probability tending to one, for a fixed $a>0$.
For an additional fixed layer, $L^{-1}=B^*R-zB^{-1}[I+(A-E-z)R]$. Outside a fraction $\varepsilon$ of small singular
directions of $B$, this has a fixed bound set $T_0=I+(A-E-z)R$. Then
\[
 L^{-1}=B^*R-zB^{-1}T_0,\qquad\|R\|\le\eta^{-1}.
\]
Write $B=U\Sigma V^*$ and let $P_{<h}$ project onto the singular
coordinates with $s_i(B)<h$. On
$\mathcal E=\ker(P_{<h}U^*T_0)$ one has
$\|B^{-1}T_0v\|\le h^{-1}\|T_0\|\|v\|$.
If there are at most $\varepsilon W$ such coordinates, then
$\dim\mathcal E\ge(1-\varepsilon)W$ and, on fixed norm cutoffs for
$A,B$, $\|L^{-1}v\|\le C_h\|v\|$ on $\mathcal E$.
The singular-value min--max principle implies that at least
$(1-\varepsilon)W$ eigenvalues of $LL^*$ are at least $C_h^{-2}$.
If $P_{\rm good}$ is their spectral projection, positivity yields
\[
 q_+=\Tr QLL^*\ge C_h^{-2}\Tr QP_{\rm good}
 \ge C_h^{-2}(q-\varepsilon W\|Q\|).
\]
This is a deterministic trace inequality and needs no independence
between $Q$ and $P_{\rm good}$. For $q\ge aW$, choose
$\varepsilon<a/(2\|Q\|)$, using $\|Q\|\le(4\eta)^{-1}$.
The excluded probability can be made arbitrarily small uniformly in
$W$: Gaussian QR gives
\[
 \det(B^*B)\ \overset{\rm law}=
 (b/W)^W\prod_{j=1}^W G_j,\qquad G_j\sim{\rm Gamma}(j,1)
\]
independently, and hence
$W^{-1}\E\sum_i\log^-(s_i(B)^2)\le C_b$.
Markov bounds the fraction $s_i(B)<a$ by
$C_b/|\log a^2|$ in expectation.
 For fixed $\varepsilon$, Markov bounds the probability of more than
$\varepsilon W$ such values by $C_b/(\varepsilon|\log h^2|)$.
Choose $h>0$ sufficiently small after fixing the desired exceptional
probability. The preceding deterministic inequality can then be
iterated any fixed number of times. In precise terms, for every
$\delta>0$ and fixed core length $r$, there exists $a_{r,\delta}>0$
with $\limsup_W\Prob(q_r<a_{r,\delta}W)\le\delta$.
Thus $q_1\to\infty$ in probability for the actual physical core, and
boundedness of $r_K$ gives $\E r_K(q_1)\to0$.
The frozen reciprocal is $K^{-1}$ at a crossing, proving
\eqref{eq:transmission-entrance}. Thus $g_t\to\infty$ in
probability as $t\downarrow0$.

For microscopic absorption set
$A-E-X\mapsto A-E+i\epsilon-X$,
$\epsilon=\kappa(r/W)/W^2$. The exact update becomes
\[
 Q_+=L^*QL+\epsilon C,\qquad
 C=\mathcal G_{12}^*\mathcal G_{12}.
\]
The inverse and jet estimates change by $O(\epsilon)$.
The injection has operator $L^p$ norm $O(W^{-2})$, trace norm
$O_{L^p}(W^{-1})$, and
$W^{-1}\E(\Tr C\mid X)=d_*+O(\delta_0+q/W+W^{-1})$.
It therefore adds $h(t)$ only to the drift of $P_1$; all higher
power variations and cross covariances are summable errors.
For completeness this also proves the entrance with a physical far
lead, not only with zero initial absorption. The all-length scalar
estimate with physical absorptions is
Lemma~\ref{lem:t26-passive-initial-scalar}. Let
$\widetilde q_+=\Tr QLL^*$ omit just the current injection.
Resolvent differentiation in the fresh $i\epsilon$ shift, using the
weighted inverse moments of Lemma~\ref{lem:passive-bounds}, changes
its conditional mean by at most $C\epsilon q$ and preserves the
centered bound $C_pq/\sqrt W$. Thus the same reciprocal calculation
adds at most $C\epsilon/q$ for $q\ge K$. Since
$q_+\ge\widetilde q_+$ and $r_K$ is decreasing, the injection can
only improve that upper bound. Summing the added error over $TW$
layers costs $C/(KW)$. The finite-core argument above persists with
$A-E$ replaced by $A-E+i\epsilon$, by the same singular-subspace
identity and positive injection. This proves
\eqref{eq:transmission-entrance} with bounded microscopic absorption
on the entire initial interval.
For a path without a far lead, the entrance is instead zero.
Indeed its absorption identity and Ward bound are
\[
 P_1=4\eta\sum_{k\le r}\epsilon_k
       \norm{P_kGP_{\rm port}}_{\HS}^2,\qquad
 \eta\norm{P_kGP_{\rm port}}_{\HS}^2
                        \le-\Im\Tr P_kG.
\]
Conditional block spectral averaging bounds the last expectation
by $CW$. Thus $\E P_1\le C\|\kappa\|_\infty(r/W+W^{-1})$.
For a physical regulator   $\delta/[N\rho(E)]$, the limiting
coefficient is exactly $h=4\pi\delta/\tau$.

\subsection{The finite dual and its corner}\label{subsec:s14-the-finite-dual-and-its-corner}

This subsection replaces the infinitely many transmission eigenvalues by
finitely many parameters, which enter a ratio of determinants. A duality
identity shows that, as functions of these parameters, the expected
ratios evolve by an explicit diffusion in finitely many variables, which
we call the dual. We realize the dual as a positive reweighting of
independent one-dimensional diffusions.
Lemma~\ref{lem:t26-dual-boundary} then shows that the dual leaves its
state space only when one pair of its particles reaches a corner, after
which that pair is removed. The lemma also fixes the class of processes
in which this law is unique, constructs the law, and gives a tail bound.
These are exactly the facts used in the proof of
Theorem~\ref{thm:transmission-law} in the next subsection. A kernel form
of the deletion rule is recorded in Remark~\ref{rem:t26-dual-kernel}, for
orientation only; it is not used later.

\paragraph{The dual generator.} Throughout this subsection put
$a_0(x)=x(1+x)$ and $c_0(y)=y(1-y)$. The dual acts on $p$ parameters
$x_1,\ldots,x_p\in(0,\infty)$ and $p$ parameters $y_1,\ldots,y_p\in(0,1)$,
and its generator is
\[
 \mathcal D_p=\frac{D_0}{2}
 \left[\sum_i(a_0(x_i)\partial_{x_i}^2+\beta_i^x\partial_{x_i})
       +\sum_j(c_0(y_j)\partial_{y_j}^2+\beta_j^y\partial_{y_j})\right],
\]
with the drifts
\begin{align}
 \beta_i^x&=1+2x_i+2x_i(1+x_i)
       \left[\sum_{k\ne i}\frac{1}{x_i-x_k}
                         -\sum_j\frac{1}{x_i+y_j}\right],\notag\\
 \beta_j^y&=1-2y_j+2y_j(1-y_j)
       \left[\sum_{\ell\ne j}\frac{1}{y_j-y_\ell}
                         -\sum_i\frac{1}{x_i+y_j}\right].
                                      \label{eq:ratio-dual-generator}
\end{align}
The first bracket in each drift contains a repulsion between particles of
the same species and an attraction between particles of different
species. The dual generator is related to the transmission generator by
the algebraic identity
\begin{equation}\label{eq:ratio-duality}
 \mathcal G^\kappa F=\mathcal D_pF-h(t)RF,\qquad
 R=\sum_i(x_i+y_i),
\end{equation}
where $F=F_{\mathbf x,\mathbf y}$ is the ratio \eqref{eq:transmission-ratios},
viewed on the left as a function of the transmission and on the right as
a function of the parameters, and $h$ is the absorption rate.

Here is a direct verification. Let
\[
 S_z=\sum_j\frac{t_j}{1+zt_j},\qquad
 K_z=\sum_j\frac{t_j^2}{(1+zt_j)^2}.
\]
The transmission generator applied to $\log F$ is
\[
 \frac{D_0}{2}\Bigl[\sum_i x_i(1+x_i)S_{x_i}^2+
 \sum_jy_j(1-y_j)S_{-y_j}^2\Bigr],
\]
and the bracket of $\log F$ with itself is
\[
 D_0\sum_k t_k^2(1-t_k)
 \left[\sum_i\frac{x_i}{1+x_it_k}
                    +\sum_j\frac{y_j}{1-y_jt_k}\right]^2.
\]
On the dual side, the first derivatives of $\log F$ in the parameters are
$-S_{x_i}$ and $-S_{-y_j}$, and the second derivatives are $K_{x_i}$ and
$-K_{-y_j}$. We expand the square and use the identity
\[
 S_z-S_w=(w-z)\sum_k\frac{t_k^2}{(1+zt_k)(1+wt_k)} .
\]
The terms with two distinct parameters give precisely the repulsion and
mixed drifts in \eqref{eq:ratio-dual-generator}, and the diagonal terms
give its second derivatives. The absorption differentiates the value at
zero of the integrand in \eqref{eq:transmission-ratios}, which gives the
term $-hRF$. Continuity in the trace extends the identity to a finite
measure that may have mass at zero.

\paragraph{Positivity: a reweighting of independent diffusions.} We start
from independent one-dimensional diffusions with generators
\[
 \mathcal L_X=\frac{D_0}{2}\partial_x[x(1+x)\partial_x],
 \qquad
 \mathcal L_Y=\frac{D_0}{2}\partial_y[y(1-y)\partial_y],
\]
on $[0,\infty)$ and $[0,1]$. At the endpoints their diffusion coefficients
vanish like the distance to the endpoint, which pushes the process back
inside, and they have unique laws that do not explode. We take $p$
independent copies of each and kill the system at the first collision of
two particles of the same species. The killed kernel on the ordered
chamber is the product of two determinants,
\[
 K_p^0(t;z,z')=
 \det[p_t^X(x_i,x_j')]_{i,j\le p}
 \det[p_t^Y(y_i,y_j')]_{i,j\le p}.
\]
By the Karlin--McGregor formula~\cite{KM59}, which holds for independent
copies of a one-dimensional diffusion with continuous strong Markov
paths and is applied here to each species separately, $K_p^0$ is the
transition density of the independent system killed at its first
same-species collision. In particular it is nonnegative. The kernels
$K_p^0$ and $Q_p$ below are used only in
Remark~\ref{rem:t26-dual-kernel}; the proof of
Lemma~\ref{lem:t26-dual-boundary} uses only the function $H_p$ and its
harmonicity.

We reweight by the positive function
\[
 H_p(\mathbf x,\mathbf y)=
 \det[(x_i+y_j)^{-1}]_{i,j\le p}
 =\frac{\prod_{i<k}(x_k-x_i)\prod_{j<\ell}(y_\ell-y_j)}{\prod_{i,j}(x_i+y_j)},
\]
a Cauchy determinant. Every term of the determinant is a product of
functions of disjoint pairs of coordinates, and
$(\mathcal L_X+\mathcal L_Y)(x+y)^{-1}=0$ away from the corner. Hence
$H_p$ is harmonic for the independent system, and
\[
 H_p^{-1}\mathcal L_p^0(H_pf)=\mathcal D_pf .
\]
In other words, $\mathcal D_p$ is the Doob $h$-transform of the
independent system by $H_p$. The transformed kernel
$H_p(z)^{-1}H_p(z')K_p^0(t;z,z')\,dz'$ is a positive subprobability
kernel.

\paragraph{The lost mass for one pair.} We specify the mass lost by the
transformed kernel, again only for Remark~\ref{rem:t26-dual-kernel}. For
$p=1$, the function $H=(x+y)^{-1}$ has flux
\[
 -\int_0^\varepsilon
 \frac{\varepsilon+x^2-(\varepsilon-x)^2}{\varepsilon^2}\,dx=-1
\]
through the line $x+y=\varepsilon$, and the fluxes through the other parts
of the boundary vanish. Consequently
\begin{equation}\label{eq:corner-flux}
 H(z)-P_t^0H(z)=\frac{D_0}{2}\int_0^t p_s^0(z,(0,0))\,ds,
\end{equation}
where $p_s^0(z,(0,0))=p_s^X(x,0)p_s^Y(y,0)$.

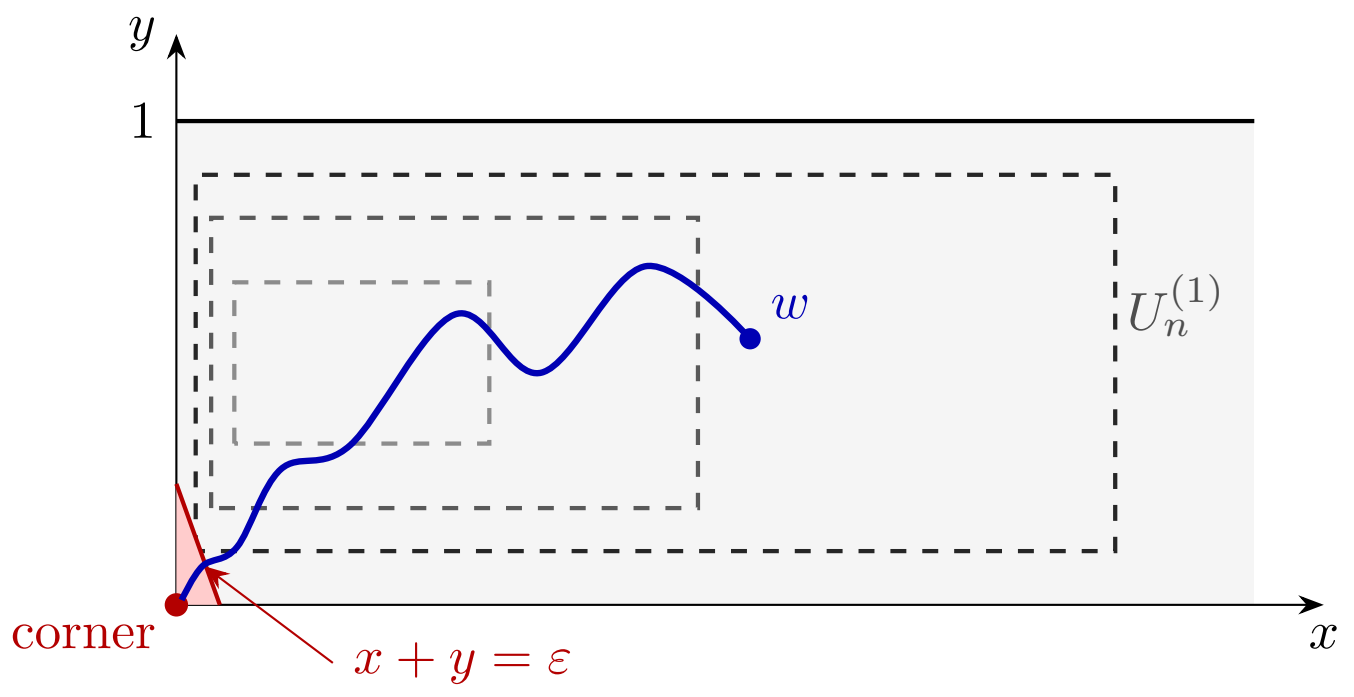
\begin{figure}[ht]
\centering
\begin{tikzpicture}[x=1.15cm,y=3.2cm,>=Stealth]
 \fill[black!4] (0,0) rectangle (6.2,1);
 \draw[->] (0,0)--(6.6,0) node[below,font=\small] {$x$};
 \draw[->] (0,0)--(0,1.18) node[left,font=\small] {$y$};
 \draw[thick] (0,1)--(6.2,1);
 \node[font=\small,left] at (0,1) {$1$};
 \foreach \n/\c in {3/black!45,5/black!65,9/black!85}{
   \draw[\c,thick,dashed] ({1/\n},{1/\n}) rectangle ({1.8*\n/3},{1-1/\n});
 }
 \node[font=\small,black!70,anchor=west] at (5.35,.62) {$U^{(1)}_n$};
 \fill[red!20] (0,0)--(.25,0)--(0,.25)--cycle;
 \draw[red!70!black,thick] (.25,0)--(0,.25);
 \draw[->,thin,red!70!black] (.9,-.12) node[right,font=\small] {$x+y=\varepsilon$} -- (.16,.08);
 \fill[red!70!black] (0,0) circle (2.2pt);
 \node[font=\small,red!70!black,below left] at (0,0) {corner};
 \draw[very thick,blue!70!black] plot[smooth,tension=.6] coordinates
   {(3.3,.55) (2.7,.7) (2.1,.48) (1.6,.6) (1.0,.33) (.6,.28) (.35,.12) (.15,.08) (.03,.01)};
 \fill[blue!70!black] (3.3,.55) circle (2pt) node[above right,font=\small] {$w$};
\end{tikzpicture}
\caption{The state space of one pair, drawn schematically; the scale in
$x$ is compressed. The $X$ coordinate lives on $(0,\infty)$ and the $Y$
coordinate on $(0,1)$. The dashed rectangles are the exhausting domains
$U^{(1)}_n$, which grow as $n$ increases. A path of the dual (blue) can
leave every $U^{(1)}_n$ in finite time only by going to the corner
$(0,0)$. The reason, proved in Step~5(a), is that $x+y$ and
$e^{-D_0t}\Lambda$ with $\Lambda=(x+y)\phi$ are nonnegative
supermartingales, so they have finite limits, while $\Lambda$ is large on
the boundary of $U^{(1)}_n$ wherever $x+y$ stays bounded below. The red
segment is the line $x+y=\varepsilon$ through which the reweighting
function $(x+y)^{-1}$ has total flux $-1$; this flux is the mass lost at
the corner, as in \eqref{eq:corner-flux}.}
\label{fig:db-one-pair}
\end{figure}

To justify \eqref{eq:corner-flux}, we use geometric coordinates. Under
$x=(\cosh r-1)/2$ and $y=(1-\cos\theta)/2$, the two scalar diffusions are
the radial parts of Brownian motion, run at speed $D_0$, on the
hyperbolic plane and on the sphere, and their radial area measures are
multiples of $dx$ and $dy$. So $p_s^X$ and $p_s^Y$ are angular averages of
smooth heat kernels. They are continuous up to the endpoints, locally
uniformly in $s>0$, and have Gaussian tails in $r$. Green's formula on the
region $\{x+y>\varepsilon,\ x<M\}$, with the flux just computed, followed
by $M\to\infty$ and then $\varepsilon\downarrow0$, gives
\eqref{eq:corner-flux}. Figure~\ref{fig:db-one-pair} shows this region.
The lifetime of the transformed process therefore has the density
$(D_0/2)(x+y)p_t^0(z,(0,0))$.

For general $p$, Remark~\ref{rem:t26-dual-kernel} computes the complete
deficit of mass. Lemma~\ref{lem:t26-dual-boundary}, in contrast,
constructs the law of the dual directly from the diffusion with generator
$\mathcal D_p$, without using kernels.

\paragraph{The dual process and its deletion rule.} A dual process
consists of $p$ particles of species $X$, at positions
$x_1<\cdots<x_p$ in $(0,\infty)$, and $p$ particles of species $Y$, at
positions $y_1<\cdots<y_p$ in $(0,1)$. Each particle moves as an
independent scalar diffusion. Killing the system when two particles of
the same species collide gives the kernel $K_p^0$. Reweighting it by the
positive function $H_p$, a Doob $h$-transform, gives $Q_p$. The
reweighting removes the collisions, but $Q_p$ still loses mass. The
following lemma shows that the diffusion with generator $\mathcal D_p$
leaves its chamber only through the paired corner, where the smallest $X$
particle and the smallest $Y$ particle reach zero together, and that the
survivors are then strictly ordered and positive. We continue by deleting
that pair and restarting with the $p-1$ surviving pairs, as in
Figure~\ref{fig:db-deletion}. The lemma fixes the class of processes in
which this law is unique, constructs it, and gives the tail bound for $R$
that is used later.

\begin{lemma}[Conservative deletion law]
\label{lem:t26-dual-boundary}
Let $D_0>0$ and let $p\ge0$ be an integer. For $q\ge1$ let
\[
 \mathcal C_q=\{0<x_1<\cdots<x_q,\ 0<y_1<\cdots<y_q<1\}
\]
be the ordered interior chamber with $q$ pairs, let
$\mathcal C_0=\{\emptyset\}$ consist of the empty state, and give
$E_p=\bigsqcup_{q=0}^p\mathcal C_q$ the topology of the disjoint union. Put
$R(\emptyset)=0$. For
$\widetilde z=(\widetilde{\mathbf x},\widetilde{\mathbf y})\in\mathcal C_{q-1}$
write
\[
 \iota\widetilde z=(0,\widetilde{\mathbf x};0,\widetilde{\mathbf y})
\]
for the paired corner with survivors $\widetilde z$: the state in which
one $X$ particle and one $Y$ particle sit at zero and the others are given
by $\widetilde z$.

Let $z\in\mathcal C_p$. A \emph{deletion process from $z$} is a filtered
probability space $(\Omega,\mathcal F,(\mathcal F_t)_{t\ge0},\Prob)$ with
an adapted process $t\mapsto(N_t,Z_t)$, where $N_t\in\{0,\dots,p\}$ counts
the surviving pairs and $Z_t\in\mathbb R^{2N_t}$ gives their positions,
with the following properties.
\begin{enumerate}[label=(D\arabic*)]
\item \emph{Start and state space.} Almost surely $N_0=p$, $Z_0=z$, $N$ is
 nonincreasing and right-continuous, and $Z_t\in\mathcal C_{N_t}$ for
 every $t\ge0$.
\item \emph{Continuity between deletions.} Almost surely $Z$ is
 continuous on every interval on which $N$ is constant.
\item \emph{Deletions at the paired corner.} Almost surely, at every jump
 time $\tau$ of $N$, we have $N_\tau=N_{\tau-}-1$ and
 $Z_t\to\iota Z_\tau$ as $t\uparrow\tau$.
\item \emph{The generator.} For $1\le q\le p$ and
 $f\in C_c^\infty(\mathcal C_q)$, with
 $\sigma_q=\inf\{t\ge0:N_t\le q\}$, the process
 \[
  M^f_t=f(Z_t)\one_{\{N_t=q\}}-f(Z_{\sigma_q})\one_{\{\sigma_q\le t\}}
        -\int_0^t(\mathcal D_qf)(Z_s)\one_{\{N_s=q\}}\,ds
 \]
 is an $(\mathcal F_t)$-martingale. Here $\mathcal D_q$ is
 \eqref{eq:ratio-dual-generator} with $q$ pairs; by (D1) and (D3),
 $N_{\sigma_q}=q$ on $\{\sigma_q<\infty\}$.
\end{enumerate}
Then the following hold.
\begin{enumerate}[label=(\alph*)]
\item \emph{Existence.} For every $z\in\mathcal C_p$ there is a deletion
 process from $z$.
\item \emph{Uniqueness.} Any two deletion processes from $z$, possibly on
 different probability spaces, have the same finite-dimensional
 distributions of $t\mapsto(N_t,Z_t)$. For $w\in\mathcal C_q$,
 $0\le q\le p$, we write $\Prob_w$ for this law with $q$ pairs.
\item \emph{Tail bound.} Every deletion process from $z$ satisfies, for
 all $K>0$ and $T\ge0$,
 \begin{equation}\label{eq:t26-dual-tail}
  \Prob\Bigl(\sup_{t\le T}R(Z_t)\ge K\Bigr)\le\frac{R(z)}{K}.
 \end{equation}
\item \emph{Restart at deletions.} For $1\le q\le p$ and every cylinder set
 $\Xi$ of paths, the map $w\mapsto\Prob_w(\Xi)$ is measurable on
 $\mathcal C_{q-1}$. Every deletion process from $z$ satisfies, for every
 $A\in\mathcal F_{\sigma_{q-1}}$ with $A\subset\{\sigma_{q-1}<\infty\}$,
 \begin{equation}\label{eq:t26-dual-restart}
  \Prob\bigl(A\cap\{(N_{\sigma_{q-1}+\cdot},Z_{\sigma_{q-1}+\cdot})\in\Xi\}\bigr)
  =\E\bigl[\one_A\,\Prob_{Z_{\sigma_{q-1}}}(\Xi)\bigr].
 \end{equation}
 That is, conditionally on $\mathcal F_{\sigma_{q-1}}$, the process after
 the deletion time $\sigma_{q-1}$ has law $\Prob_{Z_{\sigma_{q-1}}}$.
\end{enumerate}
\end{lemma}

The law $\Prob_z$ in (b) is the \emph{conservative deletion law}; we write
$\E_z$ for its expectation and $\zeta=\sigma_0$ for its extinction time,
the time at which the last pair is deleted. A deletion removes two
coordinates equal to zero, so $t\mapsto R(Z_t)$ is continuous. Hence
$\{\zeta\le t\}=\{N_t=0\}$ and $\int_0^th(t-s)R(Z_s)\,ds$ are measurable
functions of the path, and by (b) the right side of
\eqref{eq:transmission-dual-law} depends only on $(\mathbf x,\mathbf y)$.
While the process is at level $q$, and up to its exit from an open set
with compact closure in $\mathcal C_q$, a deletion process is a
continuous semimartingale with the drift and diffusion matrix of
$\mathcal D_q$ (Steps~1 and~2 of the proof). So It\^o's formula applies
there to functions that are once differentiable in time and twice in
space. These facts, together with (a), \eqref{eq:t26-dual-restart} and
\eqref{eq:t26-dual-tail}, are what the proof of
Theorem~\ref{thm:transmission-law} uses. No uniformity in $p$ is claimed.

\begin{proof}
 The dual is a diffusion with smooth
coefficients and a nondegenerate covariance inside each chamber, so it
is well defined until it reaches the boundary of the chamber (Step~1). The
whole difficulty is the boundary. Uniqueness is proved by induction on
the number of pairs, gluing the law before the first deletion to the law
after it (Step~2). The tail bound follows because $R$ has no drift
(Step~3). For existence we must know where the process can leave its
chamber. We first identify the dual, before it leaves, with the
reweighted independent system (Step~4). For one pair we show that it can
leave only through the corner and that its lifetime has no atoms
(Step~5), see Figure~\ref{fig:db-one-pair}. For several pairs we compare
the dual with independent systems of fewer pairs, and show that it leaves
only through a paired corner with ordered positive survivors (Step~6), see
Figure~\ref{fig:db-domination}. Step~7 then constructs the process by
deleting the zero pair and restarting, and the final paragraph proves the
restart property.

\emph{Preliminaries.} Throughout, $q\ge1$. Let
$\mathcal L^0_q=\sum_i\mathcal L_X^{(x_i)}+\sum_j\mathcal L_Y^{(y_j)}$ be
the generator of $2q$ independent scalar diffusions on
$O_q=(0,\infty)^q\times(0,1)^q\supset\mathcal C_q$. As shown above,
$H_q>0$ is smooth on $\mathcal C_q$, $\mathcal L^0_qH_q=0$ there, and
\begin{equation}\label{eq:t26-dual-htransform}
 \mathcal D_qf=H_q^{-1}\mathcal L^0_q(H_qf)\quad\text{on }\mathcal C_q,
 \qquad f\in C^2(\mathcal C_q).
\end{equation}
We exhaust the chamber by polytopes. For $n\ge q+2$ let
$U^{(q)}_n\subset\mathcal C_q$ be the open polytope
\[
 x_1>1/n,\quad x_{k+1}-x_k>1/n,\quad x_q<n,\qquad
 y_1>1/n,\quad y_{k+1}-y_k>1/n,\quad y_q<1-1/n .
\]
Its closure is compact, $\overline{U^{(q)}_n}\subset U^{(q)}_{n+1}$, and
$\bigcup_nU^{(q)}_n=\mathcal C_q$. We use one elementary fact about exit
times from these polytopes.
\begin{itemize}
\item[(F)] Let $w$ be continuous on $[a,b)$ with values in
 $\mathcal C_q$, where either $b=\infty$, or $b<\infty$ and $w(t)$
 converges as $t\uparrow b$ to a point outside $\mathcal C_q$. Let
 $e_n=\inf\{t\in[a,b):w(t)\notin U^{(q)}_n\}$, with $\inf\emptyset=b$.
 Then $e_n\uparrow b$, and $e_n<b$ if $b<\infty$.
\end{itemize}
Indeed, for $t<b$ the compact set $w([a,t])$ lies in some $U^{(q)}_m$,
and if $b<\infty$ the path leaves $\overline{U^{(q)}_n}$ before $b$.

\emph{Step 1: local well-posedness.} Let $b_q$ be the drift and
\[
 a_q=D_0\diag\bigl(x_i(1+x_i),\,y_j(1-y_j)\bigr)
\]
the diffusion matrix of $\mathcal D_q$, and put $\Sigma_q=a_q^{1/2}$. They
are smooth on $\mathcal C_q$, and $a_q\succ0$ there. We change them
outside $U^{(q)}_{n+1}$ into globally Lipschitz coefficients. Then the
stochastic equation
\[
 dZ_t=b_q(Z_t)\,dt+\Sigma_q(Z_t)\,dB_t,\qquad Z_0=w\in\mathcal C_q,
\]
driven by a $2q$-dimensional Brownian motion $B$, has a pathwise unique
strong solution up to the exit time $\rho_n$ from $U^{(q)}_n$, and this
solution is jointly measurable in $(w,B)$. By pathwise uniqueness, the
solutions for different $n$ agree. This defines a continuous
$\mathcal C_q$-valued process $\Phi_q(w,B)$ on $[0,\rho_\infty)$, where
$\rho_\infty=\lim_n\rho_n$. Let $\Prob^{(q)}_w$ be its law, viewed as the
law of the sequence of stopped paths $(Z^{\rho_n})_n$; this sequence
determines $\rho_\infty$ and the path on $[0,\rho_\infty)$. Conversely,
any process with the right martingale property has this law locally:
\begin{itemize}
\item[(U)] Let $U\subset\mathcal C_q$ be open with compact closure in
 $\mathcal C_q$. Let $Z'$ be a continuous adapted process with initial
 law $\mu$ on $\mathcal C_q$ and exit time $\tau'$ from $U$ (with
 $\tau'=0$ if $Z'_0\notin U$), such that
 \[
  f(Z'_{t\wedge\tau'})-f(Z'_0)-\int_0^{t\wedge\tau'}\mathcal D_qf(Z'_s)\,ds
 \]
 is a martingale for every $f\in C_c^\infty(\mathcal C_q)$. Then
 $Z'^{\tau'}$ has the law of $\Phi_q(\xi,B)$ stopped at its exit from
 $U$, where $\xi\sim\mu$ is independent of $B$.
\end{itemize}
To prove (U), take $f=\chi u$, where $\chi\in C_c^\infty(\mathcal C_q)$
equals one near $\overline U$ and $u$ is a coordinate or a product of two
coordinates. This shows that $Z'^{\tau'}$ is a continuous semimartingale
with drift $b_q$ and covariance $a_q$ up to $\tau'$. Since $a_q$ is
invertible on $\overline U$, $Z'^{\tau'}$ is a weak solution of the
stopped equation, and pathwise uniqueness gives uniqueness in law by the
Yamada--Watanabe theorem~\cite{YW71}.

\emph{Step 2: uniqueness.} We prove a slightly more general statement.
For a law $\mu$ on $\mathcal C_q$, a \emph{deletion process from $\mu$}
satisfies (D1)--(D4) with $q$ pairs, except that $Z_0$ is
$\mathcal F_0$-measurable with law $\mu$. We show by induction on $q$
that its finite-dimensional distributions depend only on $\mu$; the case
$\mu=\delta_z$ gives (b). For $q=0$ the path is constant.

Let $q\ge1$, write $\sigma=\sigma_{q-1}$ for the first deletion time, and
let $\rho_n$ be the exit time of $Z$ from $U^{(q)}_n$ before $\sigma$. By
(D1)--(D3) and (F), $\rho_n\uparrow\sigma$, $\rho_n<\sigma$ on
$\{\sigma<\infty\}$, and $N=q$ on $[0,\rho_n]$. The process $M^f$ is
bounded on $[0,t]$ by $2\|f\|_\infty+t\|\mathcal D_qf\|_\infty$, so optional
stopping of (D4) at $t\wedge\rho_n$ gives the hypothesis of (U) with
$U=U^{(q)}_n$. Hence the law of $\Pi=(Z^{\rho_n})_n$ depends only on
$\mu$; note that $Z^{\rho_n}$ is a function of $Z^{\rho_{n+1}}$. The
sequence $\Pi$ determines $\sigma$, the path on $[0,\sigma)$ and, by
(D3), $Z_\sigma$ on $\{\sigma<\infty\}$. The same test functions as in
(U) give the semimartingale property stated after the lemma, at the top
level, and at the lower levels after the shift below.

We now glue the law after the deletion. Let $T>0$ and
$A\in\mathcal F_\sigma$ with $A\subset\{\sigma\le T\}$ and $\Prob(A)>0$.
This includes every such $A\in\sigma(\Pi)$, since $Z^{\rho_n}$ is
$\mathcal F_{\rho_n}$-measurable and $\rho_n\le\sigma$. Put
\[
 N'_s=N_{\sigma+s},\qquad Z'_s=Z_{\sigma+s},\qquad
 \mathcal F'_s=\mathcal F_{\sigma+s}.
\]
Properties (D1)--(D3) pass to $(N',Z')$ with $q-1$ pairs. For $q'\le q-1$
and $f\in C_c^\infty(\mathcal C_{q'})$, the process $M^f$ vanishes on
$[0,\sigma]$, because $N>q'$ on $[0,\sigma)$ and no level is skipped. So
$M^f_{\sigma+s}$ is the process of (D4) for $(N',Z')$, and by optional
sampling at $(\sigma\wedge T)+s$ it is an $(\mathcal F'_s)$-martingale
under $\Prob(\cdot\mid A)$. Thus, under $\Prob(\cdot\mid A)$, $(N',Z')$ is
a deletion process from the law of $Z_\sigma$, and by induction its
finite-dimensional distributions are determined by the law of $\Pi$.
So $\Prob(A\cap\{(N',Z')\in\Xi\})$ is determined for all such
$A\in\sigma(\Pi)$ and all cylinder sets $\Xi$. Letting $T\to\infty$, the
$\pi$--$\lambda$ theorem determines the joint law of $\Pi$ and $(N',Z')$
on $\{\sigma<\infty\}$. Finally, $(N_t,Z_t)$ is the path before deletion
for $t<\sigma$ and $(N'_{t-\sigma},Z'_{t-\sigma})$ for $t\ge\sigma$. By
right-continuity of $(N',Z')$ in $E_{q-1}$ along dyadic times, it is a
measurable function of $(\Pi,N',Z')$, and its finite-dimensional
distributions are determined.

\emph{Step 3: the tail bound.} The function $R$ has no drift: for every
$q$, $\mathcal D_qR=0$ on $\mathcal C_q$. Indeed, the same-species drifts
contribute $2q(q-1)+2(q-1)(\sum x_i-\sum y_i)$, the mixed drifts
contribute $-2q^2-2q(\sum x_i-\sum y_i)$, and these cancel the base drift
$2q+2(\sum x_i-\sum y_i)$. Let $\tau_K=\inf\{t:R(Z_t)\ge K\}$, which is a
stopping time because $R(Z)$ is continuous, and put
$\bar\tau=T\wedge\tau_K$, so that $0\le R(Z_t)\le K$ for $t\le\bar\tau$.

Fix $q$, let $\chi_m\in C_c^\infty(\mathcal C_q)$ equal one on
$\overline{U^{(q)}_m}$, and let $e_m$ be the exit time of $Z$ from
$U^{(q)}_m$ after $\sigma_q$. By (F) on $[\sigma_q,\sigma_{q-1})$,
$e_m\uparrow\sigma_{q-1}$, and $e_m<\sigma_{q-1}$ on
$\{\sigma_{q-1}<\infty\}$. On $[\sigma_q,e_m)$ we have $N=q$ and
$\mathcal D_q(\chi_mR)(Z)=0$, so optional stopping of $M^{\chi_mR}$ at
$\bar\tau\wedge e_m$ gives
\[
 \E\bigl[(R(Z_{\bar\tau\wedge e_m})-R(Z_{\sigma_q}))
     \one_{\{\sigma_q\le\bar\tau,\ Z_{\sigma_q}\in U^{(q)}_m\}}\bigr]=0 .
\]
As $m\to\infty$, continuity of $R(Z)$ and the bound by $K$ give
\[
 \E\bigl[R(Z_{\bar\tau\wedge\sigma_{q-1}})-R(Z_{\bar\tau\wedge\sigma_q})\bigr]=0 .
\]
Summing over $q=1,\dots,p$, with $\sigma_p=0$ and $R=0$ after $\sigma_0$,
gives $\E R(Z_{\bar\tau})=R(z)$. Since $R(Z_{\bar\tau})=K$ on the event
$\{\tau_K\le T\}=\{\sup_{t\le T}R(Z_t)\ge K\}$, Markov's inequality proves
\eqref{eq:t26-dual-tail}.

\emph{Step 4: the base process and a local Doob identity.} We first check
that the independent base particles never reach the endpoints. Let
\[
 \phi_X(x)=\log\frac{(1+x)^2}{x},\qquad \phi_Y(y)=-\log(y(1-y)).
\]
Then $\mathcal L_X\phi_X=D_0/2$ and $\mathcal L_Y\phi_Y=D_0$, and both
functions are bounded below and tend to infinity at the endpoints. Hence
the scalar equations with generators $\mathcal L_X$ on $(0,\infty)$ and
$\mathcal L_Y$ on $(0,1)$, which have locally Lipschitz coefficients,
have pathwise unique solutions that never leave these intervals. Let
$\Prob^0_w$, $w\in O_q$, be the law of $2q$ independent such diffusions
started at $w$; their paths are continuous in $O_q$ for all times.

Now let $U\subset\mathcal C_q$ be open with compact closure in
$\mathcal C_q$, let $\tau_U$ be its exit time, let $t_0\ge0$, put
$T=t_0\wedge\tau_U$, and let $w\in U$. Since
$\overline U\subset U^{(q)}_n$ for some $n$, the stopped path $Z^T$ is
defined under $\Prob^{(q)}_w$. We claim that, for every event $A$ of
$Z^T$, the dual before its exit is the base process reweighted by $H_q$:
\begin{equation}\label{eq:t26-local-doob}
 \Prob^{(q)}_w(A)=\E^0_w\Bigl[\frac{H_q(Z_T)}{H_q(w)};A\Bigr].
\end{equation}
To prove this, note that under $\Prob^0_w$, It\^o's formula and
$\mathcal L^0_qH_q=0$ make $M_s=H_q(Z_{s\wedge\tau_U})/H_q(w)$ a
martingale with values in $[0,\max_{\overline U}H_q/H_q(w)]$. Let
$M_\infty=\lim_tM_t$, which exists because $M$ is a bounded martingale,
and put $\widetilde\Prob=M_\infty\Prob^0_w$, so that
$\widetilde\Prob=M_t\Prob^0_w$ on $\mathcal F_t$. For
$f\in C^2(\mathcal C_q)$, It\^o's formula and \eqref{eq:t26-dual-htransform}
make $M_sX_s$ a $\Prob^0_w$-martingale, where
\[
 X_s=f(Z_{s\wedge\tau_U})-f(w)-\int_0^{s\wedge\tau_U}\mathcal D_qf(Z_u)\,du .
\]
By Bayes' rule, $X$ is then a $\widetilde\Prob$-martingale. By (U) with
$\mu=\delta_w$, the laws of $Z^{\tau_U}$ under $\widetilde\Prob$ and under
$\Prob^{(q)}_w$ agree. Events of $Z^T$ lie in $\mathcal F_{t_0}$, where
$\widetilde\Prob=M_{t_0}\Prob^0_w$ and $M_{t_0}=H_q(Z_T)/H_q(w)$. This
proves \eqref{eq:t26-local-doob}.

\emph{Step 5: one pair.} Let $q=1$. Then $\mathcal C_1=O_1$, $H_1=1/S$
with $S=x+y$, and $\mathcal D_1S=H_1^{-1}\mathcal L^0_11=0$ by
\eqref{eq:t26-dual-htransform}. We also use the Lyapunov function
\[
 \phi=\log(1+x)-\log x-\log y-\log(1-y),\qquad \Lambda=S\phi .
\]
Since $y(1-y)\le1/4$, we have $\phi>\log4>1$. Since
$\mathcal L_X\log((1+x)/x)=0$ and $\mathcal L_Y(-\log(y(1-y)))=D_0$,
\eqref{eq:t26-dual-htransform} gives
\[
 \mathcal D_1\Lambda=S\,\mathcal L^0_1\phi=D_0S\le D_0\Lambda .
\]
Let $\rho_n$ be the exit times from $U^{(1)}_n$.

(a) \emph{Under $\Prob^{(1)}_w$, almost surely on $\{\rho_\infty<\infty\}$,
$Z_t\to(0,0)$ as $t\uparrow\rho_\infty$.} The processes $S(Z_t)$ and
$e^{-D_0t}\Lambda(Z_t)$ are nonnegative local supermartingales on
$[0,\rho_\infty)$, localized by the $\rho_n$. By the upcrossing
inequality, both have finite limits as $t\uparrow\rho_\infty$ on
$\{\rho_\infty<\infty\}$. For $0<\varepsilon\le1$, $\Lambda$ is at least
$\min(n,\varepsilon\log n)$ on $\partial U^{(1)}_n\cap\{S\ge\varepsilon\}$.
If $S(Z_t)\to s>0$, then $S(Z_{\rho_n})\ge s/2$ for large $n$, so
$\Lambda(Z_{\rho_n})\ge\min(n,\min(s/2,1)\log n)\to\infty$, which
contradicts the finite limit of $\Lambda$. Hence $S\to0$, and so
$x,y\to0$. Figure~\ref{fig:db-one-pair} illustrates the argument.

(b) \emph{The law of $\rho_\infty$ under $\Prob^{(1)}_w$ has no atoms in
$(0,\infty)$.} Apply \eqref{eq:t26-local-doob} with $U=U^{(1)}_n$,
$t_0=u$ and $A=\{\rho_n>u\}$, and let $n\to\infty$; under $\Prob^0_w$ the
exit times tend to infinity. This gives
\[
 \Prob^{(1)}_w(\rho_\infty>u)=\frac{\E^0_w[H_1(Z_u)]}{H_1(w)} .
\]
So it suffices to show that $u\mapsto\E^0_w[H_1(Z_u)]$ is continuous, for
which we need $H_1(Z_u)$ to be uniformly integrable, that is, the base
process must rarely be close to the corner. Let
$\gamma(x,y)=x/(1+x)+y\le x+y$. A direct computation, using
$x(1+x)^{-3}\le x(1+x)^{-1}$, $y(1-y)\le y$ and $1-(1+x)^{-2}\le2x/(1+x)$,
gives
\[
 \mathcal L^0_1e^{-\nu\gamma}\le\frac{D_0}{2}\,e^{-\nu\gamma}
     \bigl[(\nu^2+2\nu)\gamma-2\nu\bigr],\qquad\nu>0 .
\]
Put
\[
 \nu(s)=\frac{2e^{-D_0s}}{1-e^{-D_0s}},\qquad c(s)=(1-e^{-D_0s})^{-2},
\]
so that $\nu'=-(D_0/2)(\nu^2+2\nu)$ and $c'=-D_0\nu c$. For $s_1,T>0$, the
bounded function $\Psi(t,w')=c(s)e^{-\nu(s)\gamma(w')}$, with
$s=s_1+T-t$, satisfies $(\partial_t+\mathcal L^0_1)\Psi\le0$ on
$[0,T]\times O_1$. So $\Psi(t,Z_t)$ is a $\Prob^0_w$-supermartingale on
$[0,T]$, and
\[
 \E^0_w\bigl[c(s_1)e^{-\nu(s_1)\gamma(Z_T)}\bigr]\le c(s_1+T).
\]
Choosing $s_1$ with $\nu(s_1)=1/\eta$ gives the small-ball bound
\[
 \Prob^0_w\bigl(\gamma(Z_T)<\eta\bigr)\le e\,\frac{c(s_1+T)}{c(s_1)}
                       \le\frac{4e\,\eta^2}{(1-e^{-D_0T})^2}.
\]
Hence $\E^0_w[H_1(Z_v)^{3/2}]\le\E^0_w[\gamma(Z_v)^{-3/2}]$ is bounded for
$v$ in compact subsets of $(0,\infty)$. Since $H_1(Z_v)\to H_1(Z_u)$ almost
surely as $v\to u$, uniform integrability makes
$u\mapsto\Prob^{(1)}_w(\rho_\infty>u)$ continuous on $(0,\infty)$. Also
$\rho_\infty>0$.

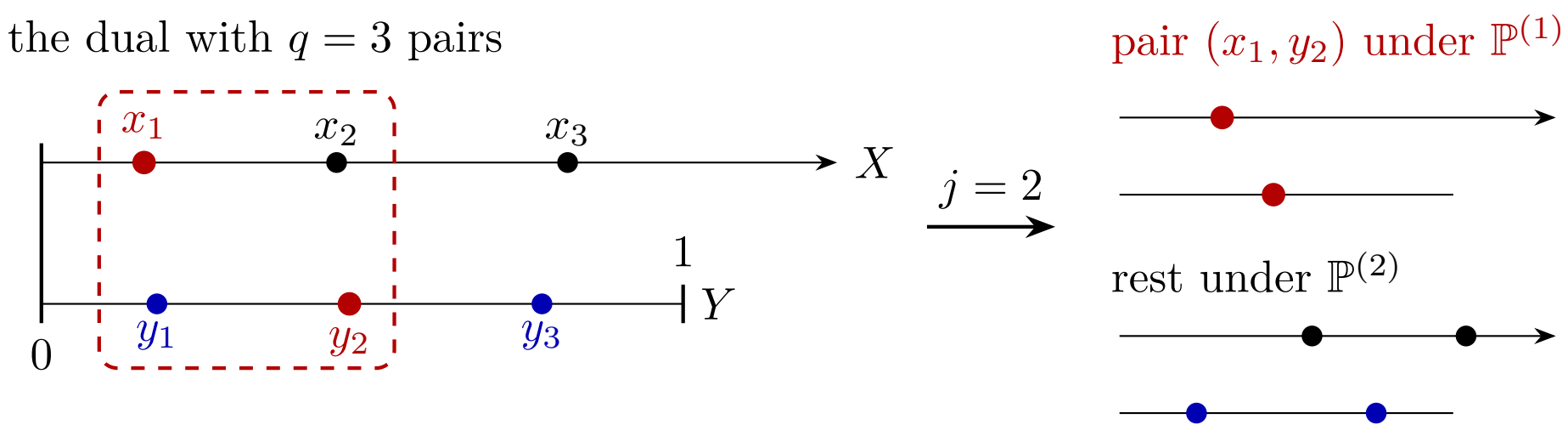
\begin{figure}[ht]
\centering
\begin{tikzpicture}[x=1cm,y=1cm,>=Stealth]
 \begin{scope}
  \node[font=\small,anchor=west] at (-.4,2.35) {the dual with $q=3$ pairs};
  \draw[->] (0,1.4)--(6.2,1.4) node[right,font=\small] {$X$};
  \draw (0,.3)--(5.0,.3) node[right,font=\small] {$Y$};
  \draw[thick] (0,.15)--(0,1.55);
  \node[font=\small,below] at (0,.15) {$0$};
  \draw[thick] (5.0,.15)--(5.0,.45) node[above,font=\small] {$1$};
  \fill[red!70!black] (.8,1.4) circle (2.6pt);
  \node[font=\small,red!70!black,above] at (.8,1.45) {$x_1$};
  \fill (2.3,1.4) circle (2.3pt) node[above,font=\small] {$x_2$};
  \fill (4.1,1.4) circle (2.3pt) node[above,font=\small] {$x_3$};
  \fill[blue!70!black] (.9,.3) circle (2.3pt) node[below,font=\small] {$y_1$};
  \fill[red!70!black] (2.4,.3) circle (2.6pt);
  \node[font=\small,red!70!black,below] at (2.4,.25) {$y_2$};
  \fill[blue!70!black] (3.9,.3) circle (2.3pt) node[below,font=\small] {$y_3$};
  \draw[red!70!black,thick,dashed,rounded corners] (.45,-.2) rectangle (2.75,1.95);
 \end{scope}
 \draw[->,thick] (6.9,.9)--(7.9,.9)
   node[midway,above,font=\small] {$j=2$};
 \begin{scope}[shift={(8.4,0)}]
  \node[font=\small,anchor=west,red!70!black] at (-.2,2.35) {pair $(x_1,y_2)$ under $\Prob^{(1)}$};
  \draw[->] (0,1.75)--(3.4,1.75);
  \draw (0,1.15)--(2.6,1.15);
  \fill[red!70!black] (.8,1.75) circle (2.6pt);
  \fill[red!70!black] (1.2,1.15) circle (2.6pt);
  \node[font=\small,anchor=west] at (-.2,.55) {rest under $\Prob^{(2)}$};
  \draw[->] (0,0.05)--(3.4,0.05);
  \draw (0,-.55)--(2.6,-.55);
  \fill (1.5,.05) circle (2.3pt);
  \fill (2.7,.05) circle (2.3pt);
  \fill[blue!70!black] (.6,-.55) circle (2.3pt);
  \fill[blue!70!black] (2.0,-.55) circle (2.3pt);
 \end{scope}
 \node[font=\small,align=center] at (5.0,-1.3)
   {$H_3\le G_1+G_2+G_3$,\qquad $G_j=H_1(x_1,y_j)\,H_2(\text{rest})$};
\end{tikzpicture}
\caption{The domination in Step~6(ii), drawn for $q=3$. The reweighting
function $H_q$ is a Cauchy determinant. Expanding it along its first row
bounds it by a sum of terms $G_j$, each of which is the product of the
one-pair function of $(x_1,y_j)$ and the $(q-1)$-pair function of the
remaining particles. Accordingly, the law of the dual is dominated, on
events before its exit, by a combination of laws $\mathcal Q^{(j)}_w$ in
which the pair $(x_1,y_j)$ (red) runs as an independent one-pair dual
and the remaining particles run as an independent $(q-1)$-pair dual, each
in its original slots. For these simpler processes the boundary behaviour
is already known: the pair dies only at the corner (Step~5(a)), and the
rest dies only at a paired corner (the induction hypothesis). A bad exit
of the $q$-pair dual would force the two independent lifetimes to
coincide, which has probability zero because the one-pair lifetime has no
atoms (Step~5(b)). Exits at a collision of two particles of the same
species, away from zero, are excluded separately in Step~6(i), because
$H_q$ is small there.}
\label{fig:db-domination}
\end{figure}

\emph{Step 6: the boundary behaviour.} We prove by induction on $q$:
\begin{itemize}
\item[(G$_q$)] Under $\Prob^{(q)}_w$, $w\in\mathcal C_q$, almost surely on
 $\{\rho_\infty<\infty\}$ the limit $z^*=\lim_{t\uparrow\rho_\infty}Z_t$
 exists and equals $\iota\widetilde z$ for some
 $\widetilde z\in\mathcal C_{q-1}$.
\end{itemize}
In words, the dual can leave its chamber only through a paired corner,
with strictly ordered positive survivors. (G$_1$) is Step~5(a). Let
$q\ge2$ and assume (G$_{q-1}$). Fix $w\in\mathcal C_q$ and $t_0>0$, let
$\rho_n$ be the exit times from $U^{(q)}_n$, and regard all laws below as
laws of the sequence $(Z^{t_0\wedge\rho_n})_n$. Let
$\mathcal A_n=\sigma(Z^{t_0\wedge\rho_n})$, which increases in $n$, and
$\mathcal A=\sigma(\bigcup_n\mathcal A_n)$, which contains all the events
below. On $\{\rho_\infty<t_0\}$ a limit $z^*$, if it exists, lies outside
$\mathcal C_q$.

(i) \emph{Collisions at positive height.} For $\eta,K,\delta>0$ let
$\mathcal B_\delta$ be the set of $w'\in\mathcal C_q$ that have some gap
$x_{k+1}-x_k$ or $y_{k+1}-y_k$ at most $\delta$, with $x_1+y_1\ge\eta$ and
$R(w')\le K$. It is relatively closed in $\mathcal C_q$, and on it
\[
 H_q\le c\delta,\qquad c=\max(1,K)^{q(q-1)}\eta^{-q^2},
\]
since in the product formula for $H_q$ one numerator factor is at most
$\delta$, the others are at most $\max(1,K)$, and every denominator factor
is at least $\eta$. Let $w\notin\mathcal B_\delta$, let $S_\delta$ be the
hitting time of $\mathcal B_\delta$, and apply \eqref{eq:t26-local-doob}
to $U'=U^{(q)}_n\setminus\mathcal B_\delta$ and to the event $A_\delta$
that $\tau_{U'}<t_0$ and $Z_{\tau_{U'}}\in\mathcal B_\delta$. On $A_\delta$,
$H_q(Z_{t_0\wedge\tau_{U'}})\le c\delta$. Moreover
$\{S_\delta<\rho_n\wedge t_0\}\subset A_\delta$, since before $\rho_n$ the
path leaves $U'$ only into $\mathcal B_\delta$. Letting $n\to\infty$ gives
\[
 \Prob^{(q)}_w\bigl(S_\delta<\rho_\infty\wedge t_0\bigr)\le\frac{c\delta}{H_q(w)} .
\]

(ii) \emph{Domination by simpler processes.} For $1\le j\le q$ and
$w'\in\mathcal C_q$, let $w'^{(j)}\in\mathcal C_{q-1}$ be obtained by
omitting $x_1$ and $y_j$, and let
\[
 G_j(w')=H_1(x_1,y_j)\,H_{q-1}(w'^{(j)}).
\]
The minor of the Cauchy matrix without row $1$ and column $j$ is the
Cauchy matrix of $w'^{(j)}$, so expanding the determinant along its first
row gives $H_q\le\sum_jG_j$ on $\mathcal C_q$; see
Figure~\ref{fig:db-domination}. Let $\mathcal Q^{(j)}_w$ be the law of the
path obtained by running the pair $(x_1,y_j)$ under $\Prob^{(1)}$ and the
remaining particles under $\Prob^{(q-1)}$, independently and each in its
original slots, up to the first of their lifetimes $\zeta_{\rm pair}$ and
$\zeta'$. The projections of $U^{(q)}_n$ to the pair and to the rest lie
in $U^{(1)}_n$ and $U^{(q-1)}_n$ (removing $y_j$ merges two gaps), so under
$\mathcal Q^{(j)}_w$ we have $\rho_n<\zeta_{\rm pair}\wedge\zeta'$ and
$\rho_\infty\le\zeta_{\rm pair}\wedge\zeta'$. Under $\Prob^0_w$ the pair
and the rest are independent base processes. So, by
\eqref{eq:t26-local-doob} for each of them and the $\pi$--$\lambda$
theorem, the product of their two bounded density martingales at $t_0$ is
the density of $\mathcal Q^{(j)}_w$ on events of the two stopped paths, and
optional stopping at $t_0\wedge\rho_n$ gives
\[
 \E^0_w\Bigl[\frac{G_j(Z_{t_0\wedge\rho_n})}{G_j(w)};A\Bigr]
 =\mathcal Q^{(j)}_w(A),\qquad A\in\mathcal A_n .
\]
With \eqref{eq:t26-local-doob} for $q$ pairs and $H_q\le\sum_jG_j$,
\[
 \Prob^{(q)}_w(A)\le\sum_jc_j\,\mathcal Q^{(j)}_w(A),\qquad
 c_j=\frac{G_j(w)}{H_q(w)},
\]
for $A$ in the algebra $\bigcup_n\mathcal A_n$, and hence on $\mathcal A$
by the monotone class theorem.

(iii) \emph{Exceptional limits are null.} Let $\Gamma$ be the event that
$\rho_\infty<t_0$ and $z^*$ either does not exist or is not of the form
$\iota\widetilde z$ with $\widetilde z\in\mathcal C_{q-1}$. Let
$\Gamma_1\subset\Gamma$ be the part on which $z^*$ exists, is finite, has
$x^*_1+y^*_1>0$ and has a zero gap, and let $\Gamma_2=\Gamma\setminus\Gamma_1$.
On $\Gamma_1$, for rational $\eta<(x^*_1+y^*_1)/2$ and $K>2R(z^*)$, the path
enters every $\mathcal B_\delta$ before $\rho_\infty$, so (i) gives
$\Prob^{(q)}_w(\Gamma_1)=0$. By (ii) and a union over rational $t_0$, it
suffices to show that $\mathcal Q^{(j)}_w(\Gamma_2)=0$ for every $j$.

Work under $\mathcal Q^{(j)}_w$ on $\{\rho_\infty<t_0\}$. On
$[0,\rho_\infty)$ the path is ordered and $x_1$ is the $x$-coordinate of
the pair. The pair is continuous in $\mathcal C_1$ and tends to $(0,0)$
at a finite lifetime, by Step~5(a). The rest is continuous in
$\mathcal C_{q-1}$ and tends, at a finite lifetime, to a paired corner
with survivors in $\mathcal C_{q-2}$, by (G$_{q-1}$). Since
$\rho_\infty\le\zeta_{\rm pair}\wedge\zeta'$, the limit $z^*$ exists and is
finite, with $y^*_q<1$. Since $z^*\notin\mathcal C_q$, we have $x^*_1=0$,
or $y^*_1=0$, or a zero gap. We consider the cases.
\begin{itemize}
\item If $x^*_1=0$, then $\rho_\infty=\zeta_{\rm pair}$, since the pair is
 positive before its lifetime. So the $y$-coordinate of the pair also
 tends to zero, and $y^*_1\le y^*_j=0$.
\item If $y^*_1=0<x^*_1$, then $\rho_\infty<\zeta_{\rm pair}$. So $j\ne1$,
 and a coordinate of the rest tends to zero. Hence $\rho_\infty=\zeta'$,
 and the smallest $x$-coordinate of the rest, $x^*_2$, is zero, which
 contradicts $x^*_2\ge x^*_1>0$.
\item If $x^*_1=y^*_1=0$ and $z^*$ is not of the form $\iota\widetilde z$
 with $\widetilde z\in\mathcal C_{q-1}$, then $\rho_\infty=\zeta_{\rm pair}$.
 If $j\ne1$, the coordinate $y_1$ of the rest tends to zero, so
 $\zeta'=\rho_\infty$. If $j=1$, the survivors are the limit of the rest,
 which lies in $\mathcal C_{q-1}$ unless $\zeta'=\rho_\infty$.
\item If $x^*_1,y^*_1>0$, then $z^*$ has a zero gap, and the path is in
 $\Gamma_1$.
\end{itemize}
Hence, up to a null set, $\Gamma_2\subset\{\zeta_{\rm pair}=\zeta'<\infty\}$:
a bad exit forces the two independent lifetimes to coincide. This event
is null, by independence and Step~5(b). This proves (G$_q$).

\emph{Step 7: existence.} We now construct the process by running the
dual, deleting the zero pair at its exit, and restarting. Let
$(\Omega,\mathcal F,(\mathcal F_t),\Prob)$, with $\mathcal F_t$ complete,
carry independent $(\mathcal F_t)$-Brownian motions $B^{(q)}$ of
dimension $2q$, $1\le q\le p$. Put $\sigma_p=0$ and $\xi_p=z$, and
proceed recursively for $q=p,\dots,1$. On $\{\sigma_q<\infty\}$, the
process $\widehat B^{(q)}=B^{(q)}_{\sigma_q+\cdot}-B^{(q)}_{\sigma_q}$ is a
Brownian motion independent of $\mathcal F_{\sigma_q}$, by the strong
Markov property, and $\xi_q\in\mathcal C_q$ is
$\mathcal F_{\sigma_q}$-measurable. Let $N_t=q$ and
\[
 Z_{\sigma_q+s}=\Phi_q(\xi_q,\widehat B^{(q)})_s
\]
for $s$ less than the lifetime $\zeta_q$ of this strong solution, which is
a stopping time of $(\mathcal F_{\sigma_q+s})_{s\ge0}$, and put
$\sigma_{q-1}=\sigma_q+\zeta_q$. So $\sigma_{q-1}$ is an
$(\mathcal F_t)$-stopping time, the segment of $Z$ on
$[\sigma_q,\sigma_{q-1})$ is continuous and adapted, and conditionally on
$\mathcal F_{\sigma_q}$ it has law $\Prob^{(q)}_{\xi_q}$. By (G$_q$),
almost surely on $\{\sigma_{q-1}<\infty\}$, $Z_t\to\iota\xi_{q-1}$ as
$t\uparrow\sigma_{q-1}$, for an $\mathcal F_{\sigma_{q-1}}$-measurable
$\xi_{q-1}\in\mathcal C_{q-1}$; on the null set where the limit fails we
take $\xi_{q-1}$ to be a fixed point of $\mathcal C_{q-1}$. Set
$Z_{\sigma_{q-1}}=\xi_{q-1}$. After $\sigma_0$ the state is empty.

Properties (D1)--(D3) hold by construction. For (D4), let
$f\in C_c^\infty(\mathcal C_q)$. On $[\sigma_q,\sigma_{q-1})$ we have
$dZ=b_q(Z)\,dt+\Sigma_q(Z)\,dB^{(q)}$, and at a finite $\sigma_{q-1}$ the
path tends to $\iota\xi_{q-1}\notin\mathcal C_q$, so $f(Z)$ and
$\nabla f(Z)$ vanish near $\sigma_{q-1}$. It\^o's formula gives, for all
$t$,
\[
 M^f_t=\int_0^t\one_{\{\sigma_q<s\le\sigma_{q-1}\}}
          \nabla f(Z_s)\cdot\Sigma_q(Z_s)\,dB^{(q)}_s .
\]
The integrand is predictable and bounded, so $M^f$ is a martingale. This
proves (a).

\emph{The restart property.} Since $\Phi_q$ is jointly measurable, the
construction from $w$ is a measurable function of $w$ and of the Brownian
motions, so $w\mapsto\Prob_w(\Xi)$ is measurable. With an
$\mathcal F_0$-measurable starting point of law $\mu$, the construction is
a deletion process from $\mu$ with law $\int\Prob_w\,\mu(dw)$, which by
Step~2 is the law of every deletion process from $\mu$. The shift
argument of Step~2 applies at $\sigma_{q-1}$ for every $q\le p$: for
$A\in\mathcal F_{\sigma_{q-1}}$ with $A\subset\{\sigma_{q-1}\le T\}$ and
$\Prob(A)>0$, the shifted process is, under $\Prob(\cdot\mid A)$, a
deletion process from the law of $Z_{\sigma_{q-1}}$, so its law is
$\E[\Prob_{Z_{\sigma_{q-1}}}(\cdot)\mid A]$. Letting $T\to\infty$ proves
\eqref{eq:t26-dual-restart}.
\end{proof}

\begin{remark}[Kernel form of the deletion]\label{rem:t26-dual-kernel}
This remark is not used later; it describes the deletion rule through
transition kernels. Let $p\ge1$, $z\in\mathcal C_p$, $H_0=1$, and let
\[
 Q_p(t;z,dz')=H_p(z)^{-1}H_p(z')K_p^0(t;z,z')\,dz'
\]
be the $H_p$-transformed killed kernel. Given the Karlin--McGregor formula
for $K_p^0$, and \eqref{eq:corner-flux} with scalar kernels continuous up
to the endpoints, $Q_p$ loses exactly the mass carried by deletions:
\begin{equation}\label{eq:t26-dual-deficit}
 1-Q_p(t;z,\mathcal C_p)
 =\int_0^t\int_{\mathcal C_{p-1}}J_p(s;z,d\widetilde z)\,ds,
\end{equation}
where
\begin{equation}\label{eq:paired-deletion}
 J_p(s;z,d\widetilde z)=
 \frac{D_0}{2}\frac{H_{p-1}(\widetilde z)}{H_p(z)}
 K_p^0(s;z;\iota\widetilde z)\,d\widetilde z
\end{equation}
is a nonnegative density with $Q_p(t)J_p(s)=J_p(t+s)$. By
\eqref{eq:t26-local-doob}, monotone convergence and the Karlin--McGregor
formula, $Q_p(t;z,\cdot)$ is the law of $Z_t$ on $\{N_t=p\}$ under
$\Prob_z$, so the left side of \eqref{eq:t26-dual-deficit} is
$\Prob_z(\sigma_{p-1}\le t)$, the probability of a first deletion by time
$t$. Heuristically, $J_p(s;z,d\widetilde z)\,ds$ is the joint law of the
first deletion time and the survivors; this is not used.
\end{remark}

\begin{proof}[Proof of Remark~\ref{rem:t26-dual-kernel}]
Steps~1 and~2 derive the remark from its two inputs, the
Karlin--McGregor formula and \eqref{eq:corner-flux}, by an exact
determinant computation. Step~3 describes the picture of paths behind
it.

\emph{Step 1: the identity for one pair.} We start with a single pair.
Put $h(x,y)=(x+y)^{-1}$, and let $h_t$ be its evolution under the two
scalar diffusions,
\[
 h_t(x,y)=\iint p_t^X(x,u)p_t^Y(y,v)h(u,v)\,du\,dv.
\]
The one-pair identity \eqref{eq:corner-flux} says exactly that $h$ is
preserved, except for a flux into the corner $(0,0)$:
\begin{equation}\label{eq:t26-rank-one-kernel-derivative}
 h_t(x,y)=h(x,y)-\frac{D_0}{2}\int_0^t
                    p_s^X(x,0)p_s^Y(y,0)\,ds.
\end{equation}
In particular $0<h_t\le h$, and $h_t\to h$ as $t\downarrow0$.

\begin{figure}[ht]
\centering
\begin{tikzpicture}[x=.78cm,y=1cm]
 \begin{scope}
  \draw[->] (0,1)--(5.4,1) node[right,font=\small] {$X\in(0,\infty)$};
  \draw (0,0)--(4,0) node[right,font=\small] {$Y\in(0,1)$};
  \draw[thick] (0,-.12)--(0,1.12);
  \node[font=\small,below] at (0,-.12) {$0$};
  \node[font=\small,below] at (4,-.05) {$1$};
  \foreach \x in {.35,1.6,2.9} \fill (\x,1) circle (2.3pt);
  \foreach \y in {.3,1.3,2.7} \fill[blue!70!black] (\y,0) circle (2.3pt);
  \draw[->,thick,red!70!black] (.35,1.25) to[bend right=30] (.05,1.1);
  \draw[->,thick,red!70!black] (.3,-.25) to[bend left=30] (.05,-.1);
  \node[font=\small,red!70!black,anchor=east] at (-.15,.5) {corner};
  \node[font=\small] at (2.3,1.6) {before: $p$ pairs};
 \end{scope}
 \draw[->,thick] (6.6,.5)--(7.8,.5) node[midway,above,font=\small] {delete};
 \begin{scope}[shift={(8.2,0)}]
  \draw[->] (0,1)--(5.4,1);
  \draw (0,0)--(4,0);
  \draw[thick] (0,-.12)--(0,1.12);
  \node[font=\small,below] at (0,-.12) {$0$};
  \node[font=\small,below] at (4,-.05) {$1$};
  \foreach \x in {1.6,2.9} \fill (\x,1) circle (2.3pt);
  \foreach \y in {1.3,2.7} \fill[blue!70!black] (\y,0) circle (2.3pt);
  \node[font=\small] at (2.3,1.6) {after: $p-1$ pairs};
 \end{scope}
\end{tikzpicture}
\caption{The deletion rule of Lemma~\ref{lem:t26-dual-boundary}. The $X$
particles (black) move on $(0,\infty)$ and the $Y$ particles (blue) on
$(0,1)$, each species strictly ordered. The transformed kernel $Q_p$ loses
mass only when the smallest particle of each species reaches zero at the
same time, the paired corner $\iota\widetilde z$. At that moment the zero
pair is deleted, and the process restarts from the $p-1$ surviving
pairs, which are strictly ordered and positive. Heuristically, the
density of the survivors at deletion is $J_p$ in \eqref{eq:paired-deletion},
and \eqref{eq:t26-dual-deficit} says that these deletions account for all
the missing mass. After at most $p$ deletions the empty state is reached,
and it is absorbing.}
\label{fig:db-deletion}
\end{figure}
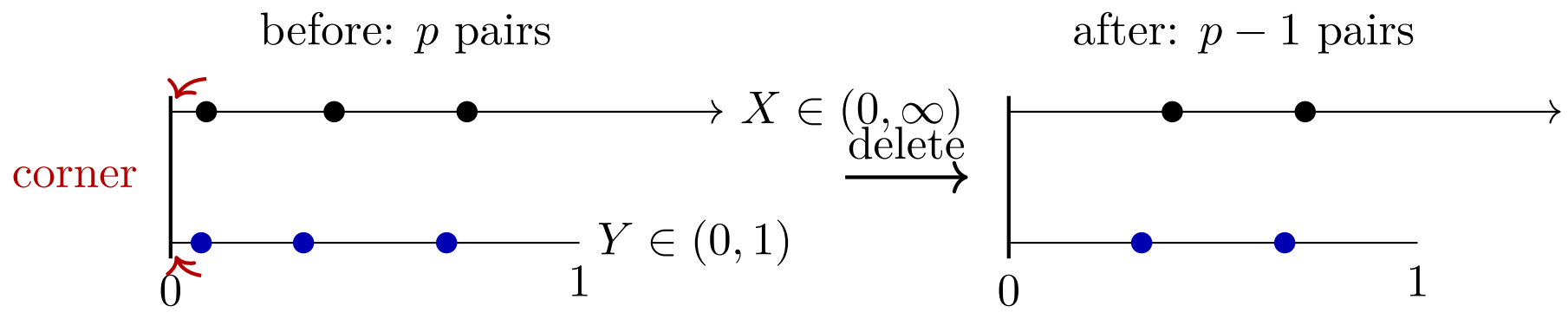

\emph{Step 2: the identity for $p$ pairs.} For $p$ pairs, $H_p$ is the
Cauchy determinant $\det[h(x_i,y_j)]_{i,j=1}^p$. We first compute the
total mass of the transformed kernel. We apply the determinant
integration identity twice, first in the ordered $u$ variables and then
in the ordered $v$ variables. This gives
\begin{equation}\label{eq:t26-integrated-cauchy-determinant}
 \int_{\mathcal C_p}K_p^0(t;z,z')H_p(z')\,dz'
                  =\det[h_t(x_i,y_j)]_{i,j=1}^p.
\end{equation}
There is no factorial factor. In each integration, the product of the two
determinants is symmetric in the integration variables, so the factor
$p!$ of the identity on the full space exactly compensates the
restriction to the ordered chamber.

Before exchanging the signed terms of the determinants, we check that all
integrals converge absolutely. We use
\[
 |H_p|\le\sum_\sigma\prod_i h(x_i,y_{\sigma(i)}),
\]
and we expand each kernel determinant in absolute value. Each resulting
integral is a product of finite values of $h_t$, and there are at most a
fixed number of permutations. This also justifies truncating the domains
and exhausting them.

Now let $C_t=(h_t(x_i,y_j))_{i,j}$, so that the right side of
\eqref{eq:t26-integrated-cauchy-determinant} is $\det C_t$. By
\eqref{eq:t26-rank-one-kernel-derivative}, its derivative is the rank-one
matrix
\[
 \frac{d}{dt}C_t=-\frac{D_0}{2}\bigl(p_t^X(x_i,0)p_t^Y(y_j,0)\bigr)_{i,j}.
\]
We differentiate the determinant by expanding in cofactors, which does
not require $C_t$ to be invertible:
\[
 -\frac{d}{dt}\det C_t
 =\frac{D_0}{2}\sum_{i,j}(-1)^{i+j}
       p_t^X(x_i,0)p_t^Y(y_j,0)\det(C_t)_{\widehat i,\widehat j}.
\]
Here $(C_t)_{\widehat i,\widehat j}$ is $C_t$ with row $i$ and column $j$
removed. We recognize the right side as the corner density: expanding the
first column of each determinant in $K_p^0(t;z;\iota\widetilde z)$ and
applying the same two determinant integrations with $p-1$ variables, the
last sum is precisely
\[
 \int_{\mathcal C_{p-1}}
       K_p^0(t;z;\iota\widetilde z)H_{p-1}(\widetilde z)
                                      \,d\widetilde z,
\]
which is nonnegative. We integrate from $0$ to $t$, use $\det C_0=H_p(z)$,
and divide by $H_p(z)$. By \eqref{eq:t26-integrated-cauchy-determinant},
the left side becomes $1-Q_p(t;z,\mathcal C_p)$ and the right side becomes
the integral of $J_p$. This proves \eqref{eq:t26-dual-deficit} globally,
including all boundary regions and regions where coordinates are large.
In particular, the deletion term accounts for the entire deficit of
mass, and no part of the boundary has been discarded. Finally, the
semigroup property of the base kernel, evaluated at the continuous corner
density and then exhausted by positive integrals, gives
$Q_p(t)J_p(s)=J_p(t+s)$.

\begin{figure}[ht]
\centering
\begin{tikzpicture}[x=1cm,y=1cm]
 \begin{scope}
  \node[font=\small,anchor=west] at (-.6,3.45) {(a) no deletion on $[0,T]$: mass $Q_p(T;z,dz')$};
  \draw[->] (0,0)--(7.4,0) node[below,font=\small] {time};
  \draw[->] (0,0)--(0,3.1) node[left,font=\small] {position};
  \draw[thick] (6.8,-.1)--(6.8,.1) node[below=3pt,font=\small] {$T$};
  \node[font=\small,below] at (0,-.05) {$0$};
  \draw[thick] plot[smooth] coordinates {(0,.6) (1.2,.9) (2.4,.55) (3.6,.8) (4.8,.5) (6.0,.75) (6.8,.65)};
  \draw[thick] plot[smooth] coordinates {(0,1.9) (1.2,1.6) (2.4,2.0) (3.6,1.7) (4.8,2.1) (6.0,1.8) (6.8,2.0)};
  \draw[thick,blue!70!black] plot[smooth] coordinates {(0,1.2) (1.2,1.35) (2.4,1.05) (3.6,1.3) (4.8,1.1) (6.0,1.3) (6.8,1.15)};
  \draw[thick,blue!70!black] plot[smooth] coordinates {(0,2.5) (1.2,2.7) (2.4,2.45) (3.6,2.75) (4.8,2.55) (6.0,2.8) (6.8,2.6)};
  \fill (0,.6) circle (1.6pt); \fill (0,1.9) circle (1.6pt);
  \fill[blue!70!black] (0,1.2) circle (1.6pt); \fill[blue!70!black] (0,2.5) circle (1.6pt);
  \node[font=\small,left] at (-.05,1.55) {$z$};
  \node[font=\small,right] at (6.85,1.55) {$z'$};
  \node[font=\small,align=left,anchor=west] at (7.6,1.6)
    {killed base bridge\\ from $z$ to $z'$:\\ no same-species\\ collision};
 \end{scope}
 \begin{scope}[yshift=-4.4cm]
  \node[font=\small,anchor=west] at (-.6,3.45) {(b) deletion at time $s\le T$: mass $J_p(s;z,d\widetilde z)\,ds$};
  \draw[->] (0,0)--(7.4,0) node[below,font=\small] {time};
  \draw[->] (0,0)--(0,3.1) node[left,font=\small] {position};
  \draw[thick] (6.8,-.1)--(6.8,.1) node[below=3pt,font=\small] {$T$};
  \draw[thick,red!70!black] (3.4,-.1)--(3.4,.1);
  \node[font=\small,below,red!70!black] at (3.4,-.08) {$s$};
  \node[font=\small,below] at (0,-.05) {$0$};
  \fill[black!6] (3.4,0) rectangle (6.8,3.0);
  \draw[thick] plot[smooth] coordinates {(0,.6) (1.0,.85) (2.0,.5) (2.8,.3) (3.4,0)};
  \draw[thick,blue!70!black] plot[smooth] coordinates {(0,1.2) (1.0,1.0) (2.0,.9) (2.8,.45) (3.4,0)};
  \draw[thick] plot[smooth] coordinates {(0,1.9) (1.2,1.6) (2.4,2.0) (3.4,1.75) (4.6,2.05) (5.8,1.8) (6.8,1.95)};
  \draw[thick,blue!70!black] plot[smooth] coordinates {(0,2.5) (1.2,2.7) (2.4,2.45) (3.4,2.65) (4.6,2.5) (5.8,2.75) (6.8,2.6)};
  \fill (0,.6) circle (1.6pt); \fill (0,1.9) circle (1.6pt);
  \fill[blue!70!black] (0,1.2) circle (1.6pt); \fill[blue!70!black] (0,2.5) circle (1.6pt);
  \fill[red!70!black] (3.4,0) circle (2.6pt);
  \draw[->,thin,red!70!black] (3.9,.42) node[right,font=\small,align=left] {paired corner $\iota\widetilde z$} -- (3.48,.06);
  \fill (3.4,1.75) circle (1.6pt);
  \fill[blue!70!black] (3.4,2.65) circle (1.6pt);
  \node[font=\small,left] at (-.05,1.55) {$z$};
  \node[font=\small,align=center] at (5.55,1.25) {$(p-1)$-pair law\\ from $\widetilde z$};
  \node[font=\small,align=left,anchor=west] at (7.6,1.6)
    {killed bridge to $\iota\widetilde z$\\ on $[0,s]$, then restart\\ from the survivors $\widetilde z$};
 \end{scope}
\end{tikzpicture}
\caption{The picture of paths in Step~3 of the proof of
Remark~\ref{rem:t26-dual-kernel}, drawn for $p=2$. The black curves are
the $X$ particles and the blue curves the $Y$ particles, each species
ordered. The two species live on $(0,\infty)$ and $(0,1)$, and are drawn on
one axis only for illustration. Panel (a) shows a path with no deletion on
$[0,T]$. It carries the mass $Q_p(T;z,dz')$ and is a killed base bridge
from $z$ to $z'$, so particles of the same species never collide. Panel
(b) shows a path whose smallest $X$ particle and smallest $Y$ particle
reach zero together at time $s$, the paired corner $\iota\widetilde z$
(red). It carries the mass $J_p(s;z,d\widetilde z)\,ds$. The zero pair is
deleted there, and on the remaining time (shaded) the law with $p-1$
pairs runs from the survivors $\widetilde z$, which are strictly positive
and distinct. By \eqref{eq:t26-dual-deficit} the masses of (a) and of
(b), integrated over $s\le T$, add up to one. The path is continuous
across the deletion, because both deleted coordinates arrive at zero.
Repeating the construction gives at most $p$ deletions on every finite
horizon.}
\label{fig:db-construction}
\end{figure}
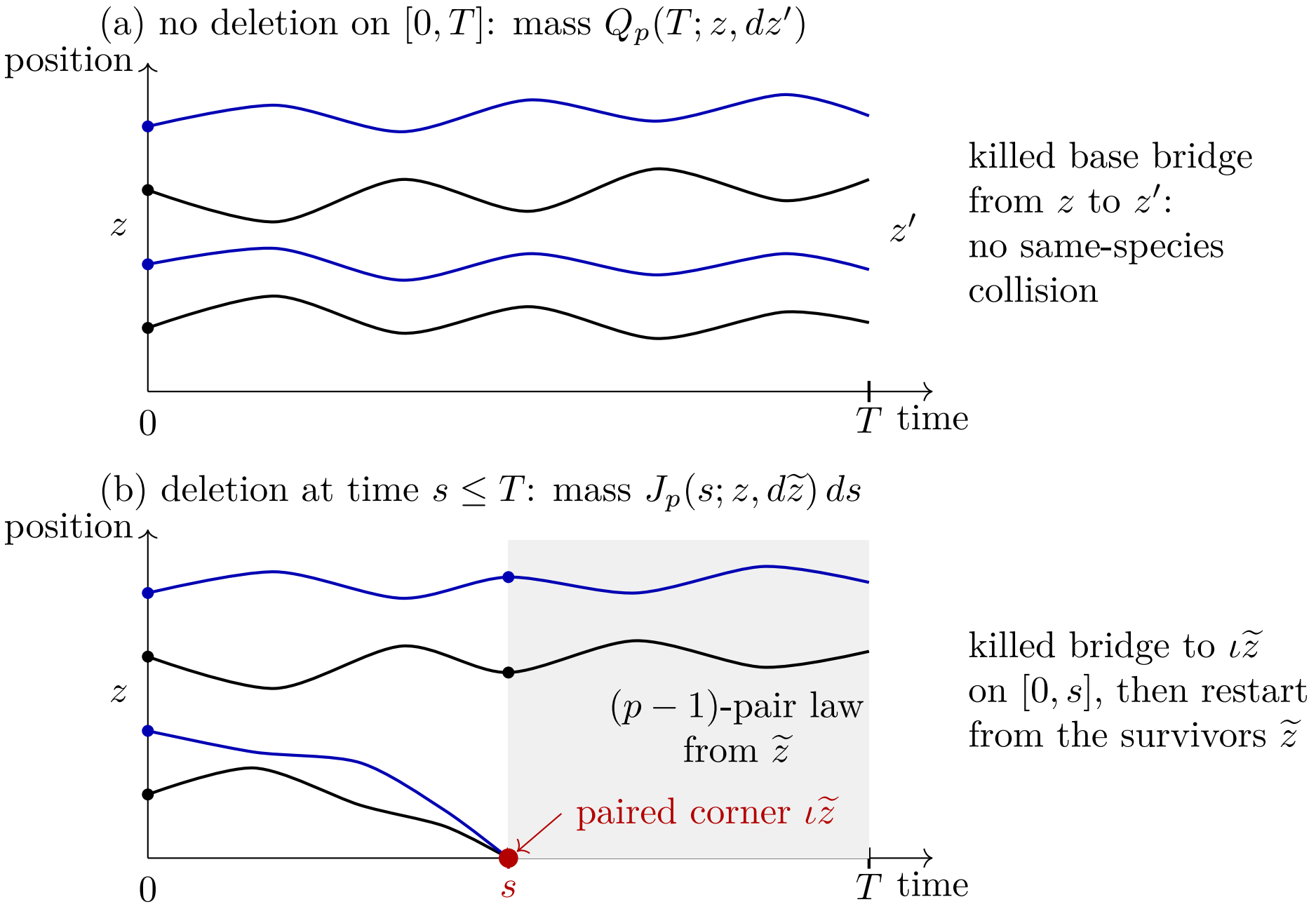

\emph{Step 3: the picture of paths.} On a horizon $T$, a path either
stays in $\mathcal C_p$, which carries the mass $Q_p(T;z,dz')$, or is
first deleted at a time $s\le T$ with survivors $\widetilde z$, which
carries the mass $J_p(s;z,d\widetilde z)\,ds$, and then continues with
$p-1$ pairs; see Figure~\ref{fig:db-construction}. By
\eqref{eq:t26-dual-deficit} these masses add up to one, and
$Q_p(t)J_p(s)=J_p(t+s)$ makes the picture consistent as $T$ changes. The
proof of Lemma~\ref{lem:t26-dual-boundary} does not use this picture.
This proves the remark. The identity $\mathcal D_pR=0$ and the tail bound
\eqref{eq:t26-dual-tail} were proved in Step~3 of the proof of
Lemma~\ref{lem:t26-dual-boundary}.
\end{proof}

The determining transforms are first
constructed for distinct interior parameters, as stated in
\eqref{eq:transmission-ratios}. The coincident parameters and the value
$y=0$ used below are limits of these bounded transforms. On sets of
bounded conductance their finite products are uniformly continuous in the
parameters, and the tail of the conductance removes that bound. Thus this
use of boundary parameters does not require an additional, freely chosen
entrance law at collisions. Let $\zeta$ be the extinction time of the
conservative deletion law.

\subsection{Determination and absence of small-channel mass}\label{subsec:s14-determination-and-absence-of-small-chann}

 We show that the finite-dimensional construction determines the
physical transmission law in Theorem~\ref{thm:transmission-law}. The
entrance fixes its initial values and the corner rule fixes its boundary
values. We solve the bounded backward equations, remove the stopping
cutoffs, and use the resulting family of determinant ratios to identify
every subsequential limit.

 \begin{proof}[Proof of the law assertion in Theorem~\ref{thm:transmission-law}]
All fixed derivatives of the ratios, on compact parameter sets with
$y_p<1$, are bounded by polynomials in $g$.
Theorem~\ref{thm:transmission-generator} and the forced extension
therefore imply \eqref{eq:ratio-duality} for every subsequential
expectation limit, with summable errors. Entrance
\eqref{eq:transmission-entrance} gives zero initial value whenever
a pair remains, because
$F_{x,y}(T)\le\exp[-(x+y)g/(1+x)]$.
At a deleted pair the ratio is one. On positive-time intervals
the error in this boundary value is bounded by
$(x_1+y_1)\E g$.

Solve the bounded backward equation recursively in the number of
pairs. Stop before collisions, before $R=K$ and before a paired
corner at distance $\varepsilon$. Diffusion duality on this compact
domain is the usual It\^o formula. Remove collision cutoffs by
noncollision, let $K\to\infty$ using $R_0/K$, and then
$\varepsilon\downarrow0$ using the boundary continuity and induction.
Nonnegative killing and $0\le F\le1$ dominate the terminal terms.
This proves \eqref{eq:transmission-dual-law} and uniqueness.

Let $y_i\downarrow0$ and permit repeated $x_i$ by continuity.
These transforms determine every mixed integer moment of the
$[0,1]$-valued variables
\[
 f_\mu(x)=\exp\!\left[-\int\frac{\log(1+xt)}t\,d\mu(t)\right],
                        \qquad x\in\mathbb Q_{>0}.
\]
Polynomial approximation on a compact cube determines their joint
law. Their map from $\mu$ is injective, since
$-\partial_x\log f_\mu(x)=\int(1+xt)^{-1}d\mu(t)$ determines
the finite measure including its mass at zero.
It is continuous on each set of bounded total mass, so its inverse
on the image is measurable. This proves one-time uniqueness.
Apply the same stopped backward calculation between two positive
times, conditionally on the earlier history, using the surviving
ratios against the earlier measure as terminal data.
It determines the transition law and hence all finite-dimensional
distributions. Tightness from the polynomial tests and conductance
moments completes process convergence.
\end{proof}

\subsection{Fractional bounds for homogeneous transmission}\label{subsec:s14-fractional-bounds-for-homogeneous-transm}
 A trace bound alone could hide mass in many almost closed channels.
For homogeneous transmission, the fractional estimates below rule out
that loss and justify finite-rank approximation. We use the negative
curvature of $\operatorname{Tr}T^p$, $0<p<1$, together with the fresh
covariance, then control the entrance and the finite-width remainders.

 Throughout the remainder of this section $\kappa\equiv0$, so the
update is $Q_+=L^*QL$ with no additive injection. We give the small-channel estimate before using finite-rank cutoffs.
For $\Phi(A)=\Tr A^p$, $0<p<1$, the exact Hessian is
\begin{equation}\label{eq:fractional-hessian}
 D^2\Phi(A)[H,H]=-\frac{p\sin(\pi p)}{\pi}
 \int_0^\infty s^{p-1}
       \Tr[(A+s)^{-1}H(A+s)^{-1}H]\,ds.
\end{equation}
Here $A\succ0$, $H$ is Hermitian and
$D^2\Phi(A)[H,H]=\partial_\tau^2\Phi(A+\tau H)|_{\tau=0}$. Indeed, with
$R=(A+s)^{-1}$, the identities
$A^p=\pi^{-1}\sin(\pi p)\int_0^\infty s^{p-1}A(A+s)^{-1}\,ds$ and
$A(A+s)^{-1}=I-sR$ give
$D^2\Phi(A)[H,H]=-2\pi^{-1}\sin(\pi p)\int_0^\infty s^p\Tr(R^2HRH)\,ds$,
differentiation under the integral being justified for
$|\tau|\norm H\le\lambda_{\min}(A)/2$, where $\norm{(A+\tau H+s)^{-1}}\le(\lambda_{\min}(A)/2+s)^{-1}$. Since
$\partial_s\Tr(RHRH)=-2\Tr(R^2HRH)$, $s^p\to0$ at $0$ and
$\Tr(RHRH)=O(s^{-2})$ at infinity, an integration by parts gives
\eqref{eq:fractional-hessian}. On the segment between
$T^{1/2}UT^{1/2}$ and $T^{1/2}\bar UT^{1/2}$ used below this applies
almost surely when $T,U,\bar U\succ0$; otherwise one regularizes as in
Lemma~\ref{lem:t22-clock-fractional-entrance}.
In an eigenbasis of the old $T$, the covariance in
Lemma~\ref{lem:passive-covariance} gives
$\E|U_{ij}-\bar U_{ij}|^2\ge c/W$ for $t_i,t_j\le1/2$,
on a small fixed jet and $g/W$ domain.
The noncancellation follows from
\begin{equation}\label{eq:jet-gap}
 |\sigma_-|\le\frac{|k_-|}{1-|k_-|d_*}
                   <\frac{k}{1-kd_*}=\sigma_+ ,
\end{equation}
which is valid at every bulk energy, at the stationary jet $s_*$ of
\eqref{eq:stationary-jet}. Indeed $|s_*|\le d_*$ gives
$|1+k_-s_*|\ge1-|k_-|d_*$, $|k_-|<c+b^2=k$
because $a^2\ne1$, $x\mapsto x/(1-xd_*)$ increases on $[0,k]$, and
$kd_*=1-b$. The gap is uniform on compact bulk sets
(Lemma~\ref{lem:source-recorded-estimates}(c)).
The fourth centered moment is $O(W^{-2})$. Uniform operator
moments let us retain this lower bound on $\{\norm U\le K\}$
for a fixed large $K$, also bounding $\norm{\bar U}$.
Along the segment between
$T^{1/2}UT^{1/2}$ and $T^{1/2}\bar UT^{1/2}$ the resolvent
dominates $(KT+s)^{-1}$. If $R\succeq R_0\succeq0$, then
$\norm{R^{1/2}HR^{1/2}}_{\HS}\ge
\norm{R_0^{1/2}HR_0^{1/2}}_{\HS}$ by contraction on both sides.
Equation~\eqref{eq:fractional-hessian} therefore bounds the
conditional concavity gap below by
\[
 \frac{c_p}{W}\sum_{t_i,t_j\le1/2}
 t_it_j\int_0^\infty\frac{s^{p-1}}{(Kt_i+s)(Kt_j+s)}\,ds
       \ge\frac{c_p'}W\left(\sum_{t_i\le1/2}t_i^p\right)^2 .
\]
For the last inequality, divide the divided difference by
$t_i^pt_j^p$; apart from a positive constant it equals
$K^{p-2}(t_i^{1-p}-t_j^{1-p})/(t_i-t_j)\ge K^{p-2}(1-p)$.
Off the retained event the concavity gap is nonnegative.
The linear term has zero mean over the whole law. Hence
\[
 \E(P_{p,+}\mid\mathcal F)
 \le(1+C_p\delta_0+C_pW^{-2})P_p
              -c_pP_p^2/W+C_pg^2/W,\qquad P_p=\Tr T^p.
\]
Multiply by a linear scalar-error weight $1+K_p\delta_{0,+}$.
The strict scalar contraction absorbs the correlated
$\delta_0P_p$ term; Cauchy--Schwarz and
$\norm{P_{p,+}}_{L^2\mid\mathcal F}\le C_pP_p$ pay the fresh
error by $C_pP_p/W$. Thus on positive-time intervals
\[
 y_{r+1}\le(1+C_p/W)y_r-c_py_r^2/W+C_p\E g_r^2/W+o(W^{-2}).
\]
The fixed high jet moments and $P_p\le W$ pay all omitted events.
The second-moment equation gives
$(\E g_t^2)'\le-D_0(\E g_t^2)^{3/2}+D_0\E g_t^2$,
so its forcing is at most $C(1+t^{-2})$.
Discrete logistic comparison, starting at $\epsilon W$ from
$P_p\le W$, proves \eqref{eq:fractional-transmission}.
For $p<q\le1$,
$\sum_{t_i\le\varepsilon}t_i^q\le\varepsilon^{q-p}P_p$.
First use $q=1$ to exclude mass at zero, then remove the cutoff
for every fractional trace both in law and in expectation.

To prove exponential decay without assuming a participation ratio,
use the Lyapunov function
$F=\sqrt g(1-P_2/(64g^2))$ on $g>0$, extended by zero at zero.
Let $r_j=P_j/g^j$. The generator gives
\[
 \frac{\mathcal GF}{D_0\sqrt g}
 =-\frac{r_2}{8}-\frac{1}{64}\left(1+\frac{15r_2^2}{8}-3r_3\right)
 +g\left[-\frac{1}{4}+\frac{r_3}{8}+
       \frac{1}{64}\left(\frac{5r_2}{4}+\frac{15r_2r_3}{8}-3r_4\right)\right].
\]
Use $0\le r_2\le1$, $r_3\le r_2^{3/2}$, $r_4\ge0$.
The coefficient of $g$ is at most $-39/512$.
For $r_2\le1/4$ the constant part is at most $-5/512$;
otherwise use $1+15x^4/8-3x^3\ge-1/8$ on $[0,1]$
to bound it by $-15/512$.
Thus $\mathcal GF\le-D_0F/128$ and
$63\sqrt g/64\le F\le\sqrt g$.
Stopping at $\varepsilon<g<K$, then removing the stops with
the conductance moments and its supermartingale property, gives
$\E\sqrt{g_t}\le(64/63)e^{-D_0(t-s)/128}\E\sqrt{g_s}$.
The uniformly bounded second moment for $t\ge1$ and H\"older give
$\E g_t\le Ce^{-D_0(t-1)/192}$.
Finally
$\E P_q\le(\E g)^{(q-p)/(1-p)}
(\E P_p)^{(1-q)/(1-p)}$ for $p<q$ proves the last assertion
of Theorem~\ref{thm:transmission-law}.
\section{The full two-input source and current}
\label{sec:full-source}
Transmission alone records how much passes through the path, but not
where the response lies inside it. To recover eigenvector profiles, we follow
sources, positive matrices that record the response to the two boundary
inputs weighted by spatial profiles, together with the current. This section
derives their joint generator.
Theorem~\ref{thm:full-source-generator} states that every joint critical
limit of finitely many sources and the current has the trace-word generator
\eqref{eq:word-generator}--\eqref{eq:word-covariance}. It keeps the exact
transported current, and it gives moment bounds for the full and tagged
clocks.
In Subsection~\ref{subsec:s15-exact-transport-and-its-ward-identity} we derive the exact transport of the two inputs and
its Ward identity. In Subsection~\ref{subsec:s15-differentiating-the-complete-source} we compute the mean contraction of
the source (Lemma~\ref{lem:complete-source-gap}) and its covariance
(Lemma~\ref{lem:full-source-covariance}). In Subsection~\ref{subsec:s15-mixed-words-and-the-closed-polynomial-ru} we show that
words containing contracting letters vanish after summation
(Lemma~\ref{lem:t26-contracting-words}), remove the stops
(Lemma~\ref{lem:t26-full-source-stops}), and complete the proof.

For a path of $j$ physical blocks let $\mathbb R=(z-\mathbb X)^{-1}$
and $\mathbb S=I-2i\eta\mathbb R$. The scattering matrix
$\mathbb S$ is unitary at real energy. Let $\mathcal U_j$ map an
incoming boundary vector in $\C^{2W}$ to the interior solution, with
the convention that the boundary amplitude is $\mathbb R$ times
that vector. Define
\[
 Y_j^q=\mathcal U_j^*\diag(q_1I_W,\ldots,q_jI_W)\mathcal U_j,\quad
 Y_{j;k}=\mathcal U_j^*P_k\mathcal U_j,\quad
 P=\begin{pmatrix}I_W&0\\0&0\end{pmatrix}.
\]
The auxiliary port masses are excluded. Write
$P'=\mathbb S^*P\mathbb S$ and $J=P'-P$.
Both $P,P'$ are projections, $\norm J\le1$ and
$g=\frac{1}{2}\Tr J^2$. A homogeneous source means a matrix
transported after its last injection.

\begin{theorem}\label{thm:full-source-generator}
For finitely many positive homogeneous sources $Q_a$ normalized
by $W$ for a tag and by $W^2$ for an interval, every joint critical
limit has the trace-word generator described by
\eqref{eq:word-generator}--\eqref{eq:word-covariance} below.
It retains the exact transported current $J$.
For every fixed $p\ge1$ and $j\le TW$,
\begin{align}
 \norm{W^{-1}\Tr Y_j^1}_{L^p}&\le C_{p,T}(j+1),&
 \E\Tr (Y_j^1)^2&\le C_TW(j+1)^3,\label{eq:full-clock-moments}\\
 \sup_{k\le j}\norm{W^{-1}\Tr Y_{j;k}}_{L^p}
                                      &\le C_{p,T}.
                                             \label{eq:full-tag-clock-moments}
\end{align}
\end{theorem}

The generator is a rule for trace observables, not an
infinite-dimensional Brownian matrix. The construction in the next
section supplies an actual probability law with these coefficients.

\subsection{Exact transport and its Ward identity}\label{subsec:s15-exact-transport-and-its-ward-identity}

 We write an update that keeps both incoming boundary amplitudes and
their interference. The recursion and Ward identity below retain the
entire old mass matrix as well as the newly added block. We derive them
by Schur inversion and use the exact conservation relation to control
the conditional drift and fluctuations.

 Using the near passive coordinates of Section~\ref{sec:transmission},
put
\[
 L=K\mathcal G_{12},\quad
 \Delta=K\mathcal G_{11}K-K,\quad
 N=\begin{pmatrix}L&\Delta\mathbb R_{12}\\0&I\end{pmatrix},
 \quad D_+=[\mathcal G_{12},\,\mathcal G_{11}K\mathbb R_{12}].
\]
Schur inversion gives
\begin{equation}\label{eq:full-clock-recursion}
 Y_{j+1}^q=N^*Y_j^qN+q_{j+1}C_+,\qquad C_+=D_+^*D_+ .
\end{equation}
The newly physical old-near block gives exactly this injection.
Let $\iota:\C^W\to\C^{2W}$ embed the near coordinates and put
\[
 F=[R-(2i\eta)^{-1}I,\mathbb R_{12}]
                       =-(2i\eta)^{-1}\iota^*\mathbb S.
\]
The full Ward identity is
\begin{equation}\label{eq:full-source-ward}
 \mathcal E:=NN^*-I=\iota\Delta F+F^*\Delta^*\iota^*,
                    \qquad FF^*=(4\eta^2)^{-1}I.
\end{equation}
For example its near diagonal entry reads
$LL^*+\Delta\mathbb R_{12}\mathbb R_{12}^*\Delta^*
=I+\Delta f+f^*\Delta^*$.
It follows by substituting the imaginary part of
$\mathcal G_{11}$ and
$-\Im X=\eta K\mathbb R_{12}\mathbb R_{12}^*K^*$.
Thus the far-lead absorption cancels before estimation.
The far incoming and outgoing amplitudes are unchanged by $N$,
which also proves the exact conservation law
\begin{equation}\label{eq:current-conservation}
 J_+=N^*JN.
\end{equation}

Both weighted orders needed below are bounded:
\[
 \Delta=K\mathcal G_{11}[(z+E)I-A+z^{-1}BB^*],\qquad
 \mathcal G_{11}K=I+\mathcal G_{11}[(z+E)I-A+z^{-1}BB^*].
\]
Lemma~\ref{lem:passive-bounds} therefore gives, on the scalar domain,
\begin{align}
 \E_j\Delta&=\Omega KM+O_{\rm op}(W^{-2}),\notag\\
 \norm{\E_j\mathcal E}&\le C(\delta_0+W^{-2}),&
 \norm{u^*(\mathcal E-\E_j\mathcal E)v}_{L^p\mid\mathcal F_j}
                &\le C_pW^{-1/2}\norm u\norm v.
                                           \label{eq:full-source-estimates}
\end{align}
The first identity follows from the weighted mean equation
$KMK-K=\Omega KM$; its residual is bounded while the weights
are still attached. Unconditionally
$\norm N_{L^p}+\norm{C_+}_{L^p}\le C_p$.

\begin{proof}[Proof of \eqref{eq:full-clock-moments}--\eqref{eq:full-tag-clock-moments}]
Write $s=W^{-1}\Tr Y$, $h=W^{-1}\Tr Y^2$.
Positivity, \eqref{eq:full-clock-recursion} and
\eqref{eq:full-source-estimates} give, for integer $p\ge2$,
\[
 \E_js_+^p\le(1+C_p\delta_0+C_p/W)s^p+C_p(s^{p-1}+1),
 \qquad
 \E_jh_+\le(1+C\delta_0+C/W)h+C(s^2+s+1).
\]
For the second estimate diagonalize $Y$ to obtain
$\E_j\Tr(Y\mathcal E Y\mathcal E)
\le C(\delta_0^2+W^{-4})\Tr Y^2+
CW^{-1}(\Tr Y)^2$.
Cross injection costs $C\Tr Y$ and
$\E_j\Tr C_+^2\le CW$.
For the first estimate the centered trace has conditional
$L^r$ norm at most $C_r(s+1)/\sqrt W$; its linear term has
zero mean and higher terms cost $C_p(s+1)^p/W$.

The passive scalar recursion has
$\delta_{0,+}\le\rho\delta_0+C/W^2+|\xi_+|$,
$\rho<1$, with $\norm{\xi_+}_{L^r\mid\mathcal F_j}\le C_r/W$
for every fixed $r$. Here $\rho<1$ is the contraction of the scalar mean
map of Subsection~\ref{subsec:s13-the-scalar-derivative-and-a-contracting-},
and the bound on $\xi_+$ is the fluctuation bound at the end of the
proof of Lemma~\ref{lem:passive-bounds}. By
\eqref{eq:passive-superexp-exit}, with $d_1$ the radius of this local
domain, the scalar exit through $TW$ layers has probability at most
$e^{-cW\log W}$.
Multiply the displayed moment recursions by
$1+K_p\delta_{0,+}$ with $K_p(1-\rho)>2C_p$.
Its predictable part absorbs $\delta_0$; its fresh part is bounded by
Cauchy--Schwarz as
$\E_j|\xi_+|s_+^p\le C_p(s^p+1)/W$ and
$\E_j|\xi_+|h_+\le C(h+1)/W$.
The resulting stopped recursions are
\[
 M_{p,j+1}\le(1+C_p/W)M_{p,j}
                  +C_pM_{p,j}^{1-1/p}+C_p,\quad
 H_{j+1}\le(1+C/W)H_j+C(M_{2,j}+\sqrt{M_{2,j}}+1).
\]
Iteration gives $M_{p,j}\le C(j+1)^p$ and $H_j\le C(j+1)^3$.
Crude unconditional fresh product moments grow at most as
$C_p^{j+1}$; H\"older and the superexponential exit probability
remove every stop. A tag starts with $C_+$ and has no further
injection, so its weighted recursion is simply
$M_{p,j+1}\le(1+C_p/W)M_{p,j}$, proving
\eqref{eq:full-tag-clock-moments}. Increasing $p$ to an integer
covers all real orders. Adding the auxiliary clock
$\mathbb R^*\mathbb R$ costs at most $2W\eta^{-2}$ in trace.
\end{proof}

The tag constant depends only on the elapsed length: the weighted
recursion starts at the birth step $k$ from $Y_{k;k}=C_+$, whose moments
are bounded uniformly in $k$, and iterates
$M_{p,j+1}\le(1+C_p/W)M_{p,j}$ over $j-k$ steps. The total length enters
only when the stops are removed, through the crude moments
$C_p^{j-k+1}$ and the exit probability $TW\,e^{-cW\log W}$ of
\eqref{eq:passive-superexp-exit}, combined by H\"older. Hence for
$R>0$, $p\ge1$ there is $C_{p,R}$, and then for every $T$ a
$W_0=W_0(p,R,T)$, such that
\begin{equation}\label{eq:tag-propagation-moments}
 \norm{W^{-1}\Tr Y_{j;k}}_{L^p}\le C_{p,R},\qquad
 W\ge W_0,\ k\le j\le TW,\ j-k\le RW,
\end{equation}
uniformly for $E$ in a compact bulk set.

\subsection{Differentiating the complete source}\label{subsec:s15-differentiating-the-complete-source}

 We identify the contracting part of the complete two-input response.
Lemma~\ref{lem:complete-source-gap} gives its mean contraction, and
Lemma~\ref{lem:full-source-covariance} supplies the joint covariance.
We differentiate the full transport identity, including the old source,
and contract the Gaussian terms while preserving their matrix order.

 \begin{lemma}\label{lem:complete-source-gap}
On the fixed scalar and first-jet domain,
\begin{equation}\label{eq:complete-source-gap}
 \E_j[N(P\mathbb S)_+N^*]
    =\theta_*P\mathbb S+
       O_{\rm op}(\delta_0+\delta_1+g/W+W^{-2}).
\end{equation}
Centered bilinear entries have conditional $L^p$ size
$C_pW^{-1/2}\norm u\norm v$.
\end{lemma}

The contraction concerns the complete transported source, including its old part. Differentiating only the newly injected source would not give this estimate.

\begin{proof}
Vary only the new near diagonal $w$, keeping the old matrix and far
parameter fixed. Normalize its incoming amplitude by
$\sqrt{\eta_w/\eta}$, where $\eta_w=(w-\bar w)/(2i)$;
the first column block of $N(w)$ acquires this factor.
Wirtinger differentiation at $w=z$ gives exactly
\begin{equation}\label{eq:source-gap-derivative}
 2i\eta\,\partial_w[N(w)N(w)^*]_{w=z}
                           =NP\mathbb S_+N^* .
\end{equation}
Indeed the two square roots contribute $NPN^*/(2i\eta)$,
whereas differentiation of the inverse contributes
$-NP\mathbb R_+N^*$. Multiplication by $2i\eta$ gives the
right side. In the Ward representation only
$\iota(\partial_w\Delta)F$ is differentiated; the adjoint inverse
is antiholomorphic. The weighted conditional mean is
$KM(w)K-K+O_{\rm op}(W^{-2})$ on a fixed $w$ disk.
Here the old $K$ is fixed, in contrast to the differentiation
in Lemma~\ref{lem:passive-source-mixing}. Explicitly
$m(w)=\langle(K-\Omega(w))^{-1}\rangle$,
$v(w)=(w-bm(w))^{-1}$ and
$\Omega(w)=c(m(w)-a)+b(v(w)-a)$. If $A=\langle M(w)^2\rangle$,
\[
 m'=A\Omega',\qquad v'=-v^2(1-bm'),\qquad
 (1-(c+b^2v^2)A)\Omega'=-bv^2.
\]
At the stationary jets, hence
$A=-s_*$, so $\Omega'=-\theta_*$.
Moreover $\partial_w(KM(w)K)=\Omega'KM(w)^2K$. Therefore
 \[
 \partial_w\E_j\Delta
   =-\theta_*I+O_{\rm op}(\delta_0+\delta_1+g/W+W^{-2}).
\]
The passive scalar trace reduction costs $O(g/W)$ as in
Lemma~\ref{lem:passive-covariance}. Cauchy's formula preserves
the weighted residual. Since $F=-(2i\eta)^{-1}\iota^*\mathbb S$,
substitution in \eqref{eq:source-gap-derivative} proves the sign
and coefficient in \eqref{eq:complete-source-gap}.
Differentiating the bounded weighted words and applying Gaussian
Poincar\'e proves the fluctuation assertion.

The projected letter needs its own identity. Triangularity gives
\[
 N(P\mathbb S_+P)N^*
 =-2i\eta\,\iota Lf_+L^*\iota^*,\qquad
 P\mathbb SP=-2i\eta\,\iota f\iota^*.
\]
The differentiated passive mean in
Lemma~\ref{lem:passive-source-mixing} consequently yields
\begin{equation}\label{eq:t26-near-source-gap}
 \E_j[N(P\mathbb S_+P)N^*]
 =\theta_*P\mathbb SP+
 O_{\rm op}(\delta_0+\delta_1+g/W+W^{-2}).
\end{equation}
Its factor $(1+\alpha q/W)^{-1}$ was absorbed using $q=g/(4\eta)$.
The same lemma supplies the centered bilinear bound. Taking adjoints
proves both conjugate gaps. This treats all four contracting letters.
 \end{proof}

\begin{lemma}\label{lem:full-source-covariance}
The centered covariance of $\mathcal E$ has crossed tensor
\begin{equation}\label{eq:full-source-covariance}
 \frac{D_0}{2W}(P\otimes P'+P'\otimes P)
 -\frac{\sigma_-}{4\eta^2W}P\mathbb S\otimes P\mathbb S
 -\frac{\bar\sigma_-}{4\eta^2W}\mathbb S^*P\otimes\mathbb S^*P
 +\mathcal R,
\end{equation}
with four-vector error bounded by
$C[W^{-2}+W^{-1}(\delta_0+\delta_1+g/W)]\prod_i\norm{v_i}$.
\end{lemma}

Both projections and their difference enter the covariance. The resulting tensor will generate the positive finite processes of the next section.

\begin{proof}The covariance we need is built from expectations
of products of two resolvents. Gaussian
integration by parts turns these expectations into a closed linear system for
two scalar traces, one at the physical block and one at the port (Step~1);
Figure~\ref{fig:fc-feedback} shows how the two traces feed into each other. We
solve this system explicitly (Step~2) and check that it is stable, because its
feedback stays below one (Step~3). Steps~4--6 turn the solution into the
covariance of two matrix entries, keeping the outside weights attached to
their inverses, and evaluate the constants. Step~7 assembles the covariance of
the full update.

\emph{Step 1: the two-resolvent system.} We use signs
$\alpha,\beta\in\{+,-\}$, with
\[
 \mathcal G_+=\mathcal G,\quad \mathcal G_-=\mathcal G^*,\qquad
 M_+=M,\quad M_-=M^*,\qquad v_+=v,\quad v_-=\bar v.
\]
For the source $T=cP_1+bP_2$, which weights the physical block by $c$
and the port block by $b$, we define
$\mathcal K^{\beta\alpha}=\E_j(\mathcal G_{\beta}T\mathcal G_{\alpha})$. The
superscript $\beta\alpha$ records the two resolvent signs. Its two
diagonal blocks satisfy
\begin{align}
 \mathcal K^{\beta\alpha}_1
  &=M_{\beta}(c+ct_1+bt_2)M_{\alpha}+O_{\rm op}(W^{-2}),\notag\\
 \mathcal K^{\beta\alpha}_2
  &=v_{\beta}(b+bt_1)v_{\alpha} I+O_{\rm op}(W^{-2}),
 \qquad t_i=W^{-1}\Tr\mathcal K^{\beta\alpha}_i.
                                      \label{eq:t26-full-replica-system}
\end{align}

To derive this system, we multiply the first inverse by its
deterministic block matrix, and we integrate each fresh supported entry
by parts. Differentiating the first inverse produces its scalar
self-energy. Differentiating the second inverse produces
$\mathcal S(\mathcal K^{\beta\alpha})\E_j\mathcal G_{\alpha}$, which is where
the traces $t_1,t_2$ of the unknown itself enter. When we separate the
scalar means, only covariances between normalized traces and fixed
inverse words remain. The interpolation contraction of
Lemma~\ref{lem:trace-matrix}, together with the passive inverse bounds,
shows that these covariances are $O_{\rm op}(W^{-2})$. We then replace
the one-resolvent means by their deterministic values and multiply by
$\operatorname{diag}(M_{\beta},v_{\beta}I)$. This gives
\eqref{eq:t26-full-replica-system}. The derivation never diagonalizes
the passive old matrix, and never commutes it with its adjoint.

\emph{Step 2: solving the system.} The system closes on the two scalars
$t_1,t_2$, as in Figure~\ref{fig:fc-feedback}. We put
\[
 A_{\beta\alpha}=W^{-1}\Tr M_{\beta}M_{\alpha},\qquad
 k_{\beta\alpha}=c+b^2v_{\beta}v_{\alpha}.
\]
Taking normalized traces in \eqref{eq:t26-full-replica-system} gives
\[
 t_1=A_{\beta\alpha}(c+ct_1+bt_2)+O(W^{-2}),\qquad
 t_2=v_{\beta}v_{\alpha}(b+bt_1)+O(W^{-2}).
\]
Substituting the second equation into the first, the bracket becomes
$c+ct_1+bt_2=k_{\beta\alpha}(1+t_1)+O(W^{-2})$. So
$1+t_1=(1-k_{\beta\alpha}A_{\beta\alpha})^{-1}+O(W^{-2})$, and therefore
\[
 \mathcal K^{\beta\alpha}_1=\sigma_{\beta\alpha}M_{\beta}M_{\alpha}+O_{\rm op}(W^{-2}),
 \qquad \sigma_{\beta\alpha}=\frac{k_{\beta\alpha}}{1-k_{\beta\alpha}A_{\beta\alpha}}.
\]

\begin{figure}[H]
\centering
\begin{tikzpicture}[x=1cm,y=1cm,>=Stealth,
  nd/.style={draw,thick,circle,minimum size=1.05cm,font=\small}]
 \node[nd] (a) at (0,0) {$t_1$};
 \node[nd] (b) at (4,0) {$t_2$};
 \draw[->,thick] (a) to[bend left=25] node[above,font=\small] {$b\,v_\beta v_\alpha$} (b);
 \draw[->,thick] (b) to[bend left=25] node[below,font=\small] {$b\,A_{\beta\alpha}$} (a);
 \draw[->,thick] (a) to[loop left,looseness=6] node[left,font=\small] {$c\,A_{\beta\alpha}$} (a);
 \node[font=\small,align=center] at (0,-1.45) {physical block};
 \node[font=\small,align=center] at (4,-1.45) {port block};
 \node[font=\small,anchor=west] at (5.6,.45) {one loop through both:};
 \node[font=\small,anchor=west] at (5.6,-.15) {$k_{\beta\alpha}A_{\beta\alpha}$,};
 \node[font=\small,anchor=west] at (5.6,-.75) {$k_{\beta\alpha}=c+b^2v_\beta v_\alpha$};
\end{tikzpicture}
\caption{The scalar feedback in the two-resolvent system
\eqref{eq:t26-full-replica-system}, used in Steps~2 and~3 of the proof
of Lemma~\ref{lem:full-source-covariance}. The trace $t_1$ of the
physical block feeds back to itself with weight $cA_{\beta\alpha}$, and
to the port trace $t_2$ with weight $bv_\beta v_\alpha$. The port feeds
back to the physical block with weight $bA_{\beta\alpha}$. Eliminating
$t_2$ combines the two routes into the single loop gain
$k_{\beta\alpha}A_{\beta\alpha}$. The system is stable when this gain
stays away from one, and then the solution is the geometric sum
$\sigma_{\beta\alpha}=k_{\beta\alpha}/(1-k_{\beta\alpha}A_{\beta\alpha})$.}
\label{fig:fc-feedback}
\end{figure}
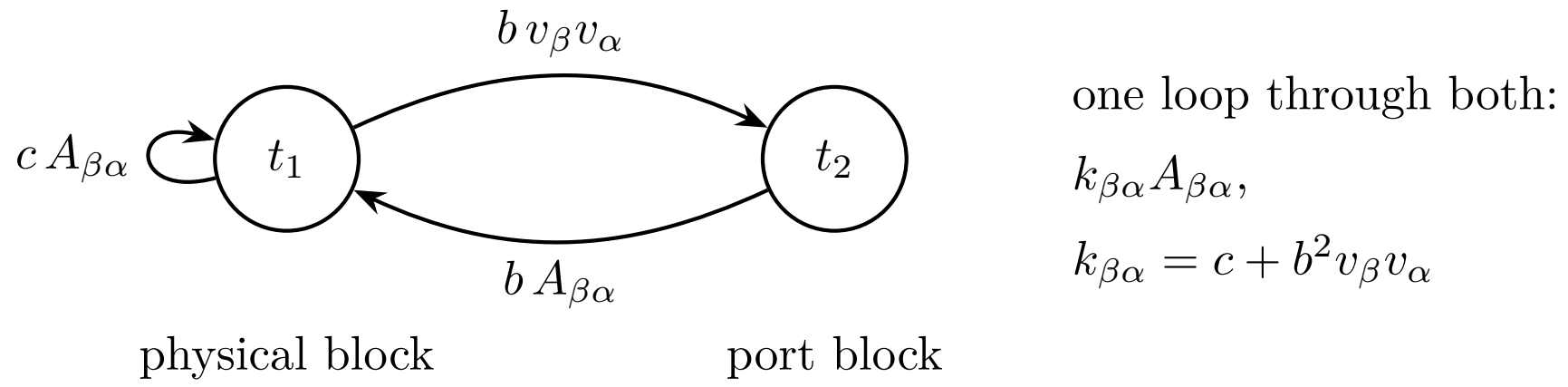

\emph{Step 3: stability.} The formula for $\sigma_{\beta\alpha}$ is useful
only if its denominator stays away from zero. For opposite signs the
denominator is $b+O(\delta_0+g/W)$, and for the same signs it is
$1+k_-s_*+O(\delta_0+\delta_1+g/W)$. To see that the systems are
stable, we use $|A_{++}|\le A_{-+}$ and $|k_{++}|\le k_{-+}$.
On the fixed small domain of the scalar errors and of $g/W$, the
advanced--retarded feedback is below $1-b/2$. Thus both two-dimensional
systems have uniformly bounded inverses.

\emph{Step 4: the entrywise covariance.} We now pass from the expected
product $\mathcal K^{\beta\alpha}$ to the covariance of two entries, and
we control the error that precedes this ladder. Consider the
integration-by-parts equation for the product
$(\mathcal G_{\alpha})_{ij}(\mathcal G_{\beta})_{kl}$. Differentiating the
second factor adds the term
\[
 W^{-1}\sum_{\gamma=1}^2\sigma_{1\gamma}
 (\mathcal G_{\beta}P_\gamma\mathcal G_{\alpha})_{kj}
                                  (\mathcal G_{\beta})_{il},
 \qquad (\sigma_{11},\sigma_{12})=(c,b).
\]
We subtract the two exact mean equations for single resolvents, and we
multiply on the left by $M_{\alpha}$. Two kinds of error remain. First,
the covariances of the centered self-energy with the product of the two
bilinear forms are $O(W^{-2})$. Indeed, in the Ornstein--Uhlenbeck
formula one of the two bilinear forms is differentiated, and its
increment is contracted with the increment of the normalized trace
before any norm is taken. Second, the product in the display can be
replaced by the product of its two means with error $O(W^{-1})$, because
each bilinear inverse word has conditional centered $L^2$ norm
$O(W^{-1/2})$. The prefactor $W^{-1}$ in the display makes this error
$O(W^{-2})$.

\emph{Step 5: the outside weights.} The words in $\Delta$ carry outside
weights. We put $U_+=V_+=K$ and $U_-=V_-=K^*$. The covariance of the
centered top words then has the crossed tensor
\begin{equation}\label{eq:t26-weighted-resolvent-tensor}
 \frac{\sigma_{\beta\alpha}}W
 (U_{\alpha}M_{\alpha}M_{\beta}V_{\beta})\otimes(U_{\beta}M_{\beta}M_{\alpha}V_{\alpha})+O(W^{-2}).
\end{equation}
All the error estimates here are four-vector estimates. The weight $K$
may be large, so we must not multiply an error by $\|K\|$ at the end.
This does not happen, because every derivative keeps its outside
weighted inverse attached, for example $K\mathcal G_{11}$ or
$\mathcal G_{11}K$. The identities preceding
\eqref{eq:full-source-estimates} bound these words and their adjoints.
The same attached orders occur when a deterministic $M_{\alpha}$ is
inserted. Thus no error was multiplied afterward by $\|K\|$.

\emph{Step 6: evaluating the weights and the coefficients.} For opposite
signs, the outside factors are, in their exact orders, $KMM^*K^*$ and
$K^*M^*MK$. For equal positive signs, both are $KM^2K$. Since
\[
 KM=MK=I+O_{\rm op}(\delta_0),
\]
all these factors equal $I+O_{\rm op}(\delta_0)$. Only the normalized
scalar words are reduced by the passive Ward identities of the preceding
replica proof. They give
\[
 A_{-+}=d_*+O(\delta_0+g/W),\qquad
 A_{++}=-s_*+O(\delta_0+\delta_1).
\]
Inserting these into the formula for $\sigma_{\beta\alpha}$, we get
$\sigma_{-+}=\sigma_++O(\delta_0+g/W)$ and
$\sigma_{++}=\sigma_-+O(\delta_0+\delta_1+g/W)$. Finally, subtracting $K$
inside $\Delta$ does not change a centered covariance, because $K$ is
fixed under the conditioning.

\emph{Step 7: assembling the covariance of $\mathcal E$.} We insert these
tensors into $\mathcal E^0=\iota\Delta^0F+F^*(\Delta^0)^*\iota^*$. The
opposite-sign terms give
\[
 \frac{\sigma_+}W\bigl[P\otimes F^*F+F^*F\otimes P\bigr].
\]
Since $F^*F=P'/(4\eta^2)$ and $\sigma_+=2\eta^2D_0$, this equals
$\frac{D_0}{2W}(P\otimes P'+P'\otimes P)$, the first two terms of
\eqref{eq:full-source-covariance}. The equal-sign terms give
$\sigma_-(\iota F)\otimes(\iota F)/W$ and its adjoint. Using
$\iota F=-(2i\eta)^{-1}P\mathbb S$, we have
$(\iota F)\otimes(\iota F)=(2i\eta)^{-2}P\mathbb S\otimes P\mathbb S
=-(4\eta^2)^{-1}P\mathbb S\otimes P\mathbb S$, which gives their minus
signs. All factors of $F$ have fixed norm, so the four-vector remainder
is the one stated in the lemma.
\end{proof}

\subsection{Mixed words and the closed polynomial rule}\label{subsec:s15-mixed-words-and-the-closed-polynomial-ru}

 To obtain a closed limiting equation, we must remove the fast
contracting response terms from products of source observables.
Lemma~\ref{lem:t26-contracting-words} gives the required summation bound,
with stopping errors handled by Lemma~\ref{lem:t26-full-source-stops}.
We classify the trace products, expand their exact one-step updates,
and pay the exit terms before identifying the polynomial rule.

 A trace word contains letters $J,Q_a$ and possibly
$P\mathbb S,\mathbb S^*P,P\mathbb SP,P\mathbb S^*P$.
Call the latter four letters contracting letters. Ordinary separators
are $I$ or $P$. Consecutive contracting letters are always separated
cyclically by at least one trace-class letter $J$ or $Q_a$.
Products of such cyclic words are admitted.

\begin{lemma}[Several contracting letters and quantitative summation]
\label{lem:t26-contracting-words}
Fix a word degree, $0<\epsilon<T$, a fixed number of source births, and
a stop $\sigma_K$ before $\sum_a\Tr Q_a+\Tr|J|>K$. Let $Z_j$ be an
admitted product of words with at least one contracting letter, on an
interval after its last tag birth $j_0\ge\epsilon W$. If the numbers of
unconjugated and conjugated contracting letters are $r$ and $\bar r$,
put $\lambda=\theta_*^r\bar\theta_*^{\bar r}$. Before the stop, the exact
update of $Z_j$ has the form
\[
 Z_{j+1}=\lambda Z_j+e_j+\xi_{j+1},\qquad
 \E_j\xi_{j+1}=0,
\]
where, for every fixed $p$,
\[
 \|1_{\{j<\sigma_K\}}e_j\|_p\le C_{K,p,\epsilon}/W,
 \qquad
 \|1_{\{j<\sigma_K\}}\xi_{j+1}\|_p
                                      \le C_{K,p,\epsilon}/\sqrt W.
\]
Continuous bounded interval injections are included. In particular,
for some $\rho<1$ depending on the degree and on the fixed parameters,
\begin{equation}\label{eq:t26-word-mixing-rate}
 \|1_{\{j<\sigma_K\}}Z_j\|_p
 \le C_{K,p,\epsilon}
       (\rho^{j-j_0}+W^{-1/2}).
\end{equation}
For any fixed ordinary trace polynomial, the total expected absolute
error in its first and second conditional increment formulas, over at
most $TW$ steps, is
\begin{equation}\label{eq:t26-word-summed-error}
 C_{K,p,\epsilon,T}(W^{-1/2}+W^{-1})=o(1).
\end{equation}
The estimate is uniform in tag locations at or after $\epsilon W$.
\end{lemma}

\begin{proof} A contracting word is the trace of a product that
contains a factor which shrinks on average at every step. We write the
one-step change of such a word as a contracting main term, a mean error and a
centered noise (Step~1). The danger is that traces of products could pick up a
factor of the dimension. Step~2 shows that this does not happen. Cutting the
trace at each random factor, as in Figure~\ref{fig:cw-cycle}, turns it into a
sum of products of scalar quantities, each of which is small. Step~3 iterates
the recursion up to the stopping time, which gives geometric decay, and
Step~4 adds up the errors over the whole path.

\emph{Step 1: the one-step update.} We substitute $N^*AN$ for each
transported letter $A=J,Q_a$, and we rotate the trace. Since
\[
 \Tr\bigl(N^*A_1N\,B_1\,N^*A_2N\,B_2\cdots\bigr)
 =\Tr\bigl(A_1\,NB_1N^*\,A_2\,NB_2N^*\cdots\bigr),
\]
each separator $B$ becomes $NBN^*$, with $B$ evaluated at the new time.
There are three kinds of separator. For a contracting $B$, the full gap
or \eqref{eq:t26-near-source-gap} gives its conditional mean as
$\theta_*B$ or $\bar\theta_*B$, with operator error
$C(\delta_0+\delta_1+g/W+W^{-2})$. This is where the factor $\lambda$
comes from. An identity separator becomes $NN^*=I+\mathcal E$. The
projected separator is treated just below. In each case we split the
new separator into its conditional mean and a centered matrix. All these
centered matrices have the bilinear $L^p$ bound
\[
 \|u^*(\text{centered matrix})v\|_p\le C_pW^{-1/2}\|u\|\|v\| .
\]

Expanding the product, each term has some number $m$ of centered
factors. The term with $m=0$ is $\lambda Z_j$ plus the errors of the
means. Terms with $m=1$ have conditional mean zero, so they go into the
martingale increment $\xi_{j+1}$. Terms with $m\ge2$ are small, as we show
next, and they go into $e_j$.

\emph{Step 2: products do not lose a dimension factor.} We now spell out
why the trace of a product of random factors is not multiplied by the
dimension. Consider a cyclic trace containing $m$ random factors
$D_1,\ldots,D_m$, and let $A_i$ be the segment between $D_i$ and
$D_{i+1}$, as in Figure~\ref{fig:cw-cycle}. We write each trace-class
segment in singular value form,
\[
 A_i=\sum_\ell s_{i\ell}u_{i\ell}v_{i\ell}^* .
\]
Then the cyclic trace becomes a sum of products of $m$ scalar bilinear
forms,
\[
 \Tr(D_1A_1D_2A_2\cdots D_mA_m)
 =\sum_{\ell_1,\ldots,\ell_m}\Bigl(\prod_is_{i\ell_i}\Bigr)
   \prod_i v_{i-1,\ell_{i-1}}^*D_iu_{i\ell_i},
\]
with indices read cyclically. The absolute sum of the coefficients is at
most $\prod_i\|A_i\|_1$. Each segment satisfies this bound with a fixed
power of $K$, because before the stop the trace norms of the letters
$Q_a$ and $J$ are at most $K$. The bounded separators $I$ and $P$, and
the contracting letters, cost only their fixed operator norms within a
segment. A segment always contains a trace-class letter, even when every
separator is chosen as a random factor.

Each scalar bilinear form is $O_{L^p}(W^{-1/2})$. So with $m$ centered
factors, conditional H\"older bounds every such sum by
$C_{K,p}W^{-m/2}$. With just one centered factor, the conditional mean
is zero. The same argument handles products of traces, by applying the
expansion in each cycle separately. Note that no free trace of $I$ or $P$
occurs, which is exactly the factor $W$ we wanted to avoid.

The noncentered factors contribute errors of size
$C_K(\delta_0+\delta_1+g/W+W^{-1})$. The moments of the scalar and jet
errors and of the conductance on $[\epsilon W,TW]$ turn this into the
stated $L^p$ bound $C/W$ for $e_j$. A new interval source has trace
$O_{L^p}(W^{-1})$ and fixed higher moments, so terms containing its
injection obey the same bound. This proves the recursion, the bound on
$e_j$, and the bound $C/\sqrt W$ on the centered part $\xi_{j+1}$.

\begin{figure}[ht]
\centering
\begin{tikzpicture}[x=1cm,y=1cm]
 \draw[thick,black!35] (0,0) circle (1.9);
 \foreach \a/\l in {90/D_1,210/D_2,330/D_3}{
   \node[draw,thick,fill=red!12,circle,minimum size=.85cm,font=\small]
     at (\a:1.9) {$\l$};
 }
 \foreach \a/\l in {150/A_1,270/A_2,30/A_3}{
   \node[draw,thick,fill=black!8,rounded corners,minimum width=1.0cm,
         minimum height=.6cm,font=\small] at (\a:1.9) {$\l$};
 }
 \node[font=\small,align=center] at (0,0) {cyclic\\ trace};
 \node[font=\small,anchor=west,align=left] at (2.9,.75)
   {$D_i$: centered fresh factors,\\ each bilinear form $O(W^{-1/2})$};
 \node[font=\small,anchor=west,align=left] at (2.9,-.55)
   {$A_i$: segments with a trace-class\\ letter, $\|A_i\|_1\le C_K$};
\end{tikzpicture}
\caption{A cyclic trace with $m=3$ centered factors, used in Step~2 of
the proof of Lemma~\ref{lem:t26-contracting-words}. Between any two
centered fresh factors $D_i$ (red) there is a segment $A_i$ (grey) that
contains a trace-class letter $Q_a$ or $J$. Expanding each segment in
singular values cuts the cycle at every $D_i$. The trace then becomes a
sum of products of $m$ scalar bilinear forms $v^*D_iu$, with coefficients
of total size at most $\prod_i\|A_i\|_1$. So the trace is
$O(W^{-m/2})$, and no free trace of an identity of dimension $W$
appears.}
\label{fig:cw-cycle}
\end{figure}
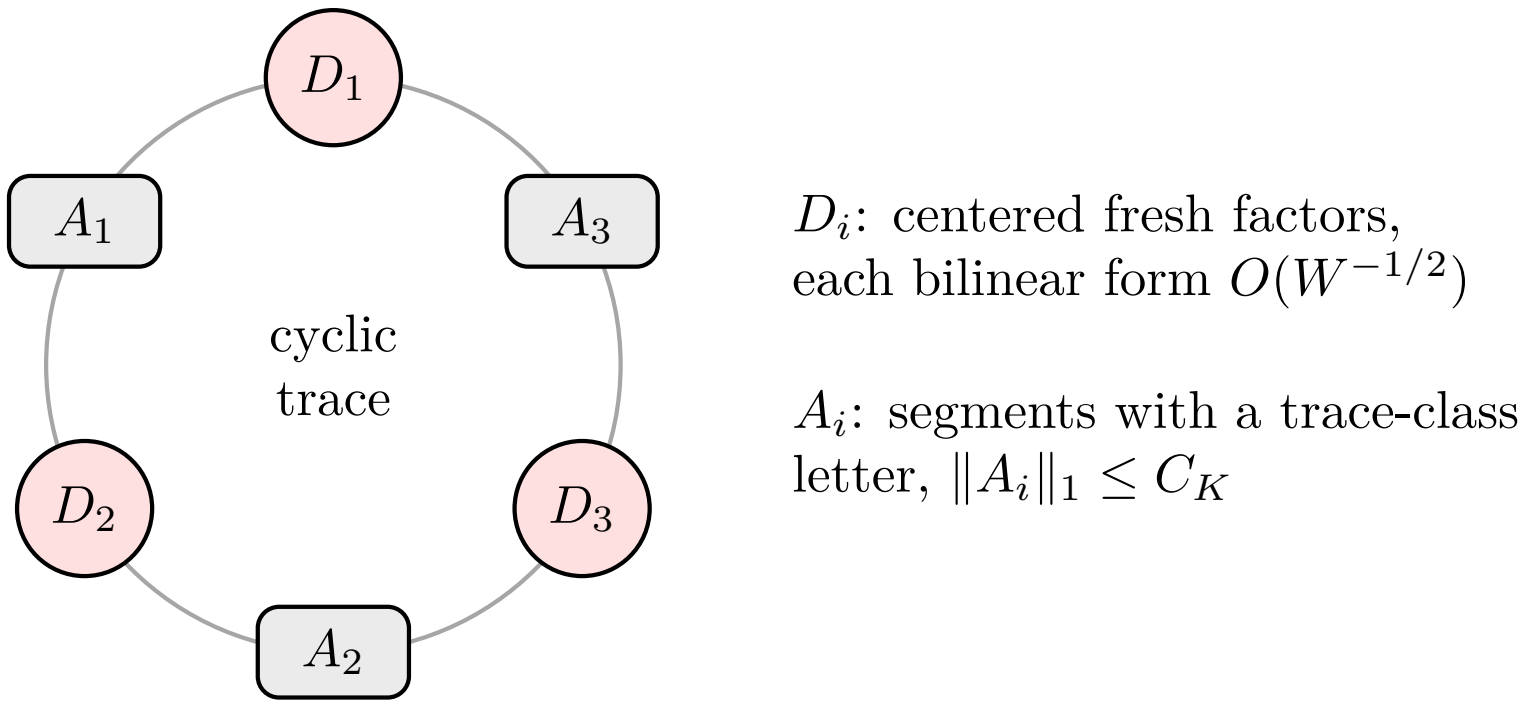

\emph{Step 3: iteration up to the stop.} We now iterate the recursion.
The survival indicator decreases in time, so multiplying by the
indicator at $j+1$ instead of $j$ can only decrease an $L^p$ norm.
Minkowski's inequality therefore gives
\[
 \|1_{\{j+1<\sigma_K\}}Z_{j+1}\|_p
 \le|\lambda|\,\|1_{\{j<\sigma_K\}}Z_j\|_p
   +\|1_{\{j<\sigma_K\}}e_j\|_p+\|1_{\{j<\sigma_K\}}\xi_{j+1}\|_p.
\]
This iteration runs on the surviving path, even though survival at
$j+1$ depends on the fresh layer. Since $|\lambda|<1$, the geometric sum
of these inequalities, started from the bounded value at $j_0$, proves
\eqref{eq:t26-word-mixing-rate}. We never use conditional centering
after conditioning on survival.

\emph{Step 4: summation of the increment errors.} Finally we sum the
errors in the first and second conditional increment formulas of a
fixed trace polynomial over at most $TW$ steps. For the first increment
we use the exact Ward mean
\[
 \E_j\mathcal E
 =-(2i\eta)^{-1}\Omega P\mathbb S+
       [-(2i\eta)^{-1}\Omega P\mathbb S]^*
                  +O_{\rm op}(\delta_0^2+W^{-2}),
\]
where $\|\Omega\|_p\le C_p/W$. Pairing its first two terms with any old
trace polynomial produces a contracting word. H\"older's inequality and
\eqref{eq:t26-word-mixing-rate} bound its expected absolute value by
$C_KW^{-1}(\rho^{j-j_0}+W^{-1/2})$.

For the second increment we use the covariance of
Lemma~\ref{lem:full-source-covariance}. The two same-sign covariance
tensors carry the prefactor $W^{-1}$. They produce either one cycle with
two contracting letters separated by trace-class segments, or two cycles
each containing such a letter. Their bound is the same. The four-vector
covariance remainder, after the singular-value expansion of Step~2, is
bounded in expectation by
\[
 C_K\{W^{-2}+W^{-1}(\E\delta_0+\E\delta_1+\E g/W)\}
 \le C_{K,\epsilon}W^{-2}.
\]
Every term with three or more centered factors has expected size
$C_KW^{-3/2}$, by the same expansion. The noncentered remainders have
size $C_KW^{-2}$.

It remains to add these bounds over the at most $TW$ steps:
\[
 \sum_{j\le TW}\Bigl[C\frac{\rho^{j-j_0}}W+\frac C{W^{3/2}}
   +\frac C{W^{3/2}}+\frac C{W^2}\Bigr]
 \le C\Bigl(\frac1W+\frac T{\sqrt W}+\frac TW\Bigr).
\]
The geometric term gives $O(W^{-1})$, and the terms of size $W^{-3/2}$
and $W^{-2}$ give $O(W^{-1/2})$ and $O(W^{-1})$. This proves
\eqref{eq:t26-word-summed-error}. At each of the finitely many tag
births we restart the geometric term. Its total sum then remains
$O(W^{-1})$, which proves the uniformity in the tag locations.
\end{proof}

For clarity, the projected separator used above has the exact form
 \[
 NPN^*=P+P\mathcal EP-
       \begin{pmatrix}C_NC_N^*&0\\0&0\end{pmatrix},
 \qquad C_N=\Delta\mathbb R_{12}.
\]
Put $A=\mathbb R_{12}\mathbb R_{12}^*$, so
$\Tr A=g/(4\eta^2)$ and $\|A\|\le\eta^{-2}$.
 The opposite-sign weighted tensor
\eqref{eq:t26-weighted-resolvent-tensor}, paired with a spectral
decomposition of $A$, gives
\[
 \E_j\Delta A\Delta^*
 =\frac{D_0g}{2W}I+R_j,
 \quad
 \|R_j\|\le C\{\delta_0^2+gW^{-2}
            +(g/W)(\delta_0+\delta_1+g/W)\}.
\]
The squared mean of $\Delta$ is included in $\delta_0^2$.
The centered quadratic expression has four-vector $L^p$ bound
$C_pg/W+C_p\delta_0\sqrt{g/W}$: diagonalize $A$ and apply the
bilinear $W^{-1/2}$ estimate to its two centered factors.
Consequently, on the trace stop and positive-time interval,
\begin{equation}\label{eq:projected-source-transport}
 NPN^*=P+P\mathcal EP-\frac{D_0g}{2W}P+\mathcal Z+\mathcal R,
\end{equation}
where $\mathcal Z$ is centered with the stated bilinear bound
and $\E\norm{\mathcal R}\le C_KW^{-2}$.
This bilinear $L^p$ size $C_K/W$,
$\E\|\mathcal R\|\le C_KW^{-2}$, and its covariance with
$\mathcal E^0$ is $O_K(W^{-3/2})$ after four-vector pairing.
This last assertion is conditional Holder followed by the scalar
moment bound. An operator bound $O(W^{-1})$ for $\mathcal Z$ is
not used. This proves all the projected-separator estimates required
in Lemma~\ref{lem:t26-contracting-words}.

Here is the resulting explicit rule. Write a cyclic word in normal form
$\Tr(A_1B_{e_1}\cdots A_rB_{e_r})$,
$A_i\in\{J,Q_a\}$, $B_0=I$, $B_1=P$. Replace each separator by
\begin{equation}\label{eq:word-generator}
 B_e+P^e\,d\mathcal B\,P^e
                  -\one_{\{e=1\}}(D_0g/2)P\,dt.
\end{equation}
Retain first drift terms and pairs of differential terms, using
\begin{align}
 \Tr(A\,d\mathcal B\,B\,d\mathcal B)
 &=\frac{D_0}{2}
       [\Tr(AP)\Tr(BP')+\Tr(AP')\Tr(BP)]\,dt,\notag\\
 \Tr(A\,d\mathcal B)\Tr(B\,d\mathcal B)
 &=\frac{D_0}{2}
       [\Tr(APBP')+\Tr(AP'BP)]\,dt.       \label{eq:word-covariance}
\end{align}
Finally expand $P'=P+J$, reduce $P^2=P$ and apply the product
rule to products of traces. The symbol $d\mathcal B$ specifies
only this finite contraction rule.
For continuous injection $q_a(t)$, add $d_*q_a(t)\,dt$ to
both $\Tr Q_a$ and $\Tr(Q_aP)$, and zero to its first-order
contribution in a word containing another trace-class letter.

The estimates of
Lemma~\ref{lem:t26-contracting-words} hold only before the stop
$\sigma_K$, the first time at which the source traces and the nuclear
norm of the current exceed $K$. To identify the limit of the source
process, we must remove this stop at the end. The following lemma provides the
three facts needed for that. First, the stop is rarely reached, i.e., its
probability is $O(1/K)$. Second, the compensators in the martingale
problem are uniformly integrable, so they converge in expectation and
not only in law. Third, the stopped processes converge together with
their exit values, provided the levels $K$ are chosen properly.

\begin{lemma}[Trace stops, exit values and compensators]
\label{lem:t26-full-source-stops}
For $0<\epsilon<T$, the physical source and current processes satisfy
\[
 \limsup_W\Prob\!\left\{
 \sup_{\epsilon W\le j\le TW}
       \Big(\sum_a\Tr Q_{a,j}+\Tr|J_j|\Big)>K\right\}
 \le C_{\epsilon,T}/K.
\]
The suprema of the source traces have every fixed positive moment on
this interval. All drift and quadratic-variation compensators of
ordinary trace polynomials in
\eqref{eq:word-generator}--\eqref{eq:word-covariance} are uniformly
integrable. The stops may be removed through increasing regular levels,
retaining their exit values.
\end{lemma}

\begin{proof}
Throughout we use the positive transport estimates, not the limiting
generator, which has not yet been identified.

\emph{Step 1: moments of the source traces.} Let $H\succeq0$ be an old
homogeneous source and $h=\Tr H$. By \eqref{eq:full-source-estimates},
one step of transport satisfies
\[
 \E_jh_+\le(1+C\delta_0+CW^{-2})h,\qquad
 \|h_+\|_{2\mid\mathcal F_j}\le Ch.
\]
So $h$ grows in mean by at most a factor $1+C/W$ per step, except for the
term $C\delta_0h$ coming from the scalar error $\delta_0$. To absorb that
term, we use the scalar-error weight from the proof of
\eqref{eq:full-clock-moments}. It makes
\[
 e^{-C(j-j_0)/W}(1+A\delta_{0,j})h_j
\]
a nonnegative supermartingale up to the scalar exit. In the conditional
calculation we keep the scalar innovation and $h_+$ together. Their
product is bounded by $Ch/W$, while the contraction of $\delta_0$,
through the condition $A(1-\rho)>C$, absorbs the old term $\delta_0h$.

For the normalized physical sources with bounded interval injection,
the same calculation applies to any integer power $p$. We use
$1+\Tr Q$ in place of $h$, so that the injection of size $O(W^{-1})$ is
included. This gives
\[
 \E_j[(1+A_p\delta_{0,+})(1+\Tr Q_+)^p]
 \le(1+C_p/W)(1+A_p\delta_0)(1+\Tr Q)^p.
\]
At each of the finitely many tag births we restart the estimate. The
traces at birth have all fixed moments, by the weighted layer bounds.
Over the at most $TW$ steps, the factor $(1+C_p/W)^{TW}$ stays bounded. So
we may apply the maximal inequality for nonnegative supermartingales to
the exponent $2p$. It gives a tail of order $u^{-2}$ for the supremum of
$(1+\Tr Q)^p$, and hence bounded $p$-th moments of the supremum of the
trace.

These bounds hold up to the original scalar exit, whose probability is
$e^{-cW\log W}$  by \eqref{eq:passive-superexp-exit}.. On the exit event we use crude moments of the fresh
products, which are at most $C_p^{TW+1}$. Their product with
$e^{-cW\log W}$ tends to zero. All these statements concern the actual
path, which we never condition on survival.

\emph{Step 2: the current.} We take $j_0=\lfloor\epsilon W\rfloor$, and we
transport the positive matrix $H_{j_0}=|J_{j_0}|$ forward by
$H_+=N^*HN$. The current is transported in the same way, by exact
current conservation. Congruence preserves the Loewner order, and the
inequalities hold at $j_0$. So
\[
 -H_j\preceq J_j\preceq H_j\qquad(j\ge j_0),
 \qquad\text{hence}\qquad \Tr|J_j|\le\Tr H_j.
\]
The carrier $H$ is positive, so Step~1 applies to it. We need only its
first moment at $j_0$. The eigenvalues of $J_{j_0}$ come in pairs
$\pm T_k^{1/2}$, where $T_k$ are the transmission eigenvalues. So the
homogeneous fractional transmission bound gives
\[
 \E\Tr H_{j_0}=2\E\Tr T_{j_0}^{1/2}\le C_\epsilon.
\]
The first-moment supermartingale of Step~1 and its maximal inequality
then give the tail $C_\epsilon/K$ for the current. This uses only the
first nuclear moment at $j_0$, and not all nuclear moments. Together
with the maximal bounds for the sources, it proves the stated tail of
the stop.

\emph{Step 3: integrability of the compensators.} The compensators are
finite combinations of ordinary cyclic words. We bound every word
without using $\Tr|J|$. A cyclic word containing a $Q_a$ is bounded by a
fixed product of source traces, because $\|J\|,\|P\|\le1$. A word made
only of $J$ and $P$ either vanishes or is bounded by $2g$. With one $J$
we use
\[
 \Tr J=0,\qquad \Tr JP=-g.
\]
With at least two copies of $J$, we put two of them in Hilbert--Schmidt
norm and all other factors in operator norm. Since $\Tr J^2=2g$, the
word is at most $2g$. Thus every coefficient in the finite contraction
rule is bounded by a fixed polynomial in $g$ and the source traces,
whatever the size of $\Tr|J|$.

The conductance $g$ has every fixed moment, uniformly at positive times.
So H\"older's inequality, together with a higher moment of the sources
and of the conductance, bounds these coefficient polynomials in
$L^{1+\delta}$ for some $\delta>0$. Jensen's inequality in time then gives
the uniform integrability of their time integrals. The quadratic
variations have the same property. In particular, removing the current
stop costs at most a positive power of $C_\epsilon/K$ in each
compensator. We never need a high nuclear moment of the current.

\emph{Step 4: tightness and continuity of stopped paths.} We now prepare
the passage to the limit at a fixed stop. On each fixed trace stop, the
centered increment of every fixed trace polynomial has conditional $L^p$
size $C_KW^{-1/2}$, and its drift has $L^p$ size $C_K/W$. These bounds
follow from the singular-value expansion in
Lemma~\ref{lem:t26-contracting-words}, where we choose higher fixed
scalar moments for its coefficient errors. The maximal bound for the
martingale part and the bound on the drift give Aldous tightness.
Moreover, away from the finitely many prescribed births,
\[
 \E\sum_j|\Delta Z_j|^p\le TW\cdot C_KW^{-p/2}=C_KW^{1-p/2}\to0
 \qquad(p>2),
\]
so the largest jump tends to zero. Every subsequential limit of the
stopped processes is therefore continuous there. We keep the full
increment at the exit, rather than discarding it.

\emph{Step 5: continuity of the nuclear norm of the current.} The stop
also involves $\Tr|J|$, which is not a trace polynomial. We show that it
has the same continuity, without assuming anything about the limiting
source law. At $j_0$ we truncate the spectral decomposition of $J_{j_0}$
to finitely many paired transmission channels. By the fractional
transmission estimate, the expected nuclear norm of the omitted part
tends to zero as the truncation grows. We transport the absolute value
of the omitted part as a positive carrier, as in Step~2. Its maximal
bound shows that the omitted nuclear norm stays uniformly small in
probability on the whole interval $[\epsilon,T]$. The retained part has
fixed rank. The entries of its signed Gram matrix are polarized traces
of rank-one sources, to which the tightness of Step~4 applies. Its
nuclear norm is a continuous function of this finite Gram matrix. So the
full nuclear norm is a limit in probability, uniformly in time, of
continuous limits, and it is itself continuous.

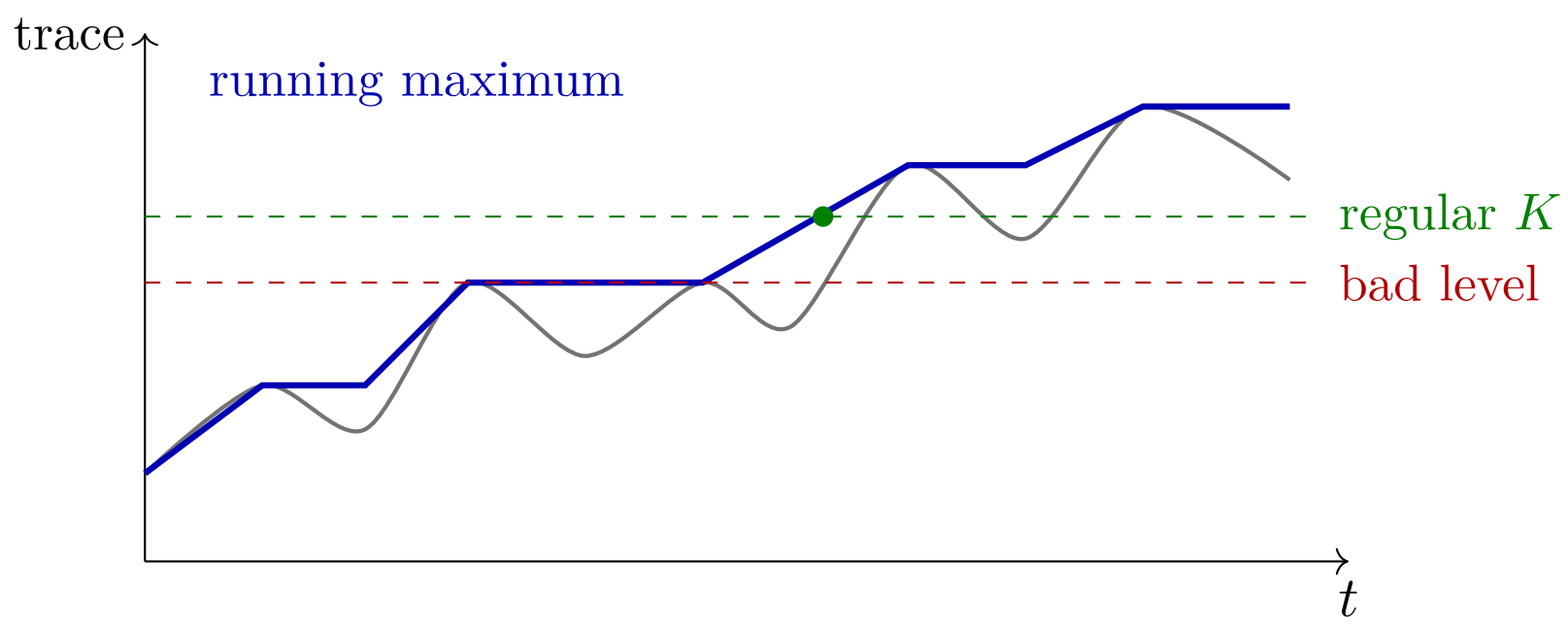
\begin{figure}[ht]
\centering
\begin{tikzpicture}[x=1cm,y=1cm]
 \draw[->] (0,0)--(8.2,0) node[below] {$t$};
 \draw[->] (0,0)--(0,3.6) node[left,font=\small] {trace};
 \draw[thick,black!55] plot[smooth] coordinates
   {(0,.6) (.8,1.2) (1.5,.9) (2.2,1.9) (3.0,1.4) (3.8,1.9) (4.4,1.6)
    (5.2,2.7) (6.0,2.2) (6.8,3.1) (7.8,2.6)};
 \draw[very thick,blue!70!black]
   (0,.6)--(.8,1.2)--(1.5,1.2)--(2.2,1.9)--(3.8,1.9)--(5.2,2.7)
   --(6.0,2.7)--(6.8,3.1)--(7.8,3.1);
 \draw[dashed,red!70!black] (0,1.9)--(8,1.9)
   node[right,font=\small,red!70!black] {bad level};
 \draw[dashed,green!50!black] (0,2.35)--(8,2.35)
   node[right,font=\small,green!50!black] {regular $K$};
 \fill[green!50!black] (4.62,2.35) circle (2pt);
 \node[font=\small,blue!70!black,anchor=west] at (.3,3.25) {running maximum};
\end{tikzpicture}
\caption{The choice of regular levels in Step~6 of the proof of
Lemma~\ref{lem:t26-full-source-stops}. The grey curve is a continuous
path and the blue curve its running maximum. At the red level the
running maximum has a plateau: the path touches the level and turns
back. A small perturbation of the path can then move the first exit time
by a macroscopic amount, so the exit map is not continuous there. Such
levels are at most countable. At a regular level, such as the green one,
the path crosses strictly, and the exit time and the exit value depend
continuously on the path. Choosing the stopping levels $K$ in this way
lets the stopped processes converge together with their exit values.}
\label{fig:fs-levels}
\end{figure}

\emph{Step 6: choosing regular levels and removing the stop.} It remains
to pass the exit time and the exit value to the limit. For a continuous
path, the running maximum has a nontrivial plateau at only countably
many levels, as in Figure~\ref{fig:fs-levels}. By Fubini's theorem we
may choose increasing deterministic levels $K$ that avoid these
exceptional levels almost surely, and that also avoid the atoms of the
terminal maximum. At such levels the map that sends a path to its first
exit time and to its exit value is continuous under uniform convergence.
This proves the passage to the limit for the stopped paths. The jumps
at the prescribed births are retained separately. Finally, the tail
estimate of Step~2 and the integrability of the compensators of Step~3
remove $K$.
\end{proof}

\begin{proof}[Completion of the proof of Theorem~\ref{thm:full-source-generator}] The theorem identifies the limit of the sources $Q_a$
and the current $J$ through the generator of their trace words. The
moment bounds \eqref{eq:full-clock-moments} and
\eqref{eq:full-tag-clock-moments} have already been proved. It remains to
show four things. The exact one-step updates produce the generator
\eqref{eq:word-generator}--\eqref{eq:word-covariance}, up to errors
and these errors vanish after summation over the path (Step~1).
The births of new sources enter correctly, both for interval injections
and for tags (Steps~2 and~3). And the stops can be removed (Step~4).
Finally we illustrate the rule on three examples (Step~5).

\emph{Step 1: the generator and the summed errors.} We expand the exact
congruences of the finite-width update. Two kinds of terms survive. The
opposite-sign terms of the covariance \eqref{eq:full-source-covariance}
give the quadratic-variation part, and the drift in
\eqref{eq:projected-source-transport} gives the drift part. Together
they give precisely \eqref{eq:word-generator}--\eqref{eq:word-covariance}.

All other terms are errors, and Lemma~\ref{lem:t26-contracting-words}
bounds their total expected absolute size by
\[
 C_{K,\epsilon,T}(W^{-1/2}+W^{-1})
\]
on each positive-time interval and before each trace stop. This bound
covers all the error terms: the Ward mean, both same-sign tensors, all
covariance remainders, and the cubic terms. Each of them contains a
contracting word or a small prefactor, so it vanishes after summation
over the at most $TW$ steps.

\emph{Step 2: the injection of an interval source.} An interval source
is born continuously along an interval, as in Figure~\ref{fig:fg-births},
and we must compute what each step adds. For this we solve the
advanced--retarded ladder with source $P_1$. As in the proof of
Lemma~\ref{lem:full-source-covariance}, its traces satisfy
\[
 t_1=A(1+ct_1+bt_2)+O(W^{-2}),\qquad t_2=b|v|^2t_1+O(W^{-2}),\qquad A=\langle MM^*\rangle.
\]
Substituting the second equation into the first gives
$t_1=A/(1-(c+b^2|v|^2)A)+O(W^{-2})$, and therefore
\[
 t_2=\frac{b|v|^2A}{1-(c+b^2|v|^2)A}+O(W^{-2})
 =d_*+O(\delta_0+g/W+O(W^{-2})).
\]
The term $O(W^{-2})$, from the covariance terms of
\eqref{eq:t26-full-replica-system} and the $O_{\rm op}(W^{-2})$ part of
$\E_j\Delta$ in \eqref{eq:full-source-estimates}, is not shown to vanish,
so it is kept; it is harmless, since $t_2$ enters each of the at most
$TW$ injection steps with the factor $W^{-1}$. The quantity $t_2$ is $W^{-1}\E_j\Tr\mathcal G_{12}^*\mathcal G_{12}$, the
expected normalized trace of the injected matrix. Its normalized trace
fluctuates by $O_{L^p}(W^{-1})$ around this mean. The far part of $C_+$
has trace $O_{L^p}(g)$, because its weighted factors stay attached and
$\Tr\mathbb R_{12}\mathbb R_{12}^*=g/(4\eta^2)$.

A bounded interval injection, normalized by $W^2$, therefore has
operator size $O_{L^p}(W^{-2})$, rank at most $W$, and near trace
\[
 \frac{d_*q}{W}+o_{L^p}(W^{-1}).
\]
So each step adds a small positive matrix whose near trace has the
deterministic value $d_*q/W$. Pairing such an injection with any other
trace-class letter costs $O_K(W^{-2})$. Over the at most $TW$ steps, the
total mixed contribution is $O_K(W^{-1})$, and it vanishes.

\emph{Step 3: the birth of a tag.} A tag is born at a single step and is
normalized by $W$. It has initial trace and near trace $d_*$, vanishing
squared trace, and vanishing overlaps with the old finite frames. A
tag produces only finitely many birth increments, and we retain them
exactly. After its birth, the contracting words that contain it obey
the uniform geometric estimate of Lemma~\ref{lem:t26-contracting-words},
started at the birth time.

\emph{Step 4: removing the stops.} The bounds of Steps~1--3 hold before
the source and current stops. Lemma~\ref{lem:t26-full-source-stops}
removes these stops, together with their exit values and every
compensator. This identifies the claimed generator in every joint
critical limit.

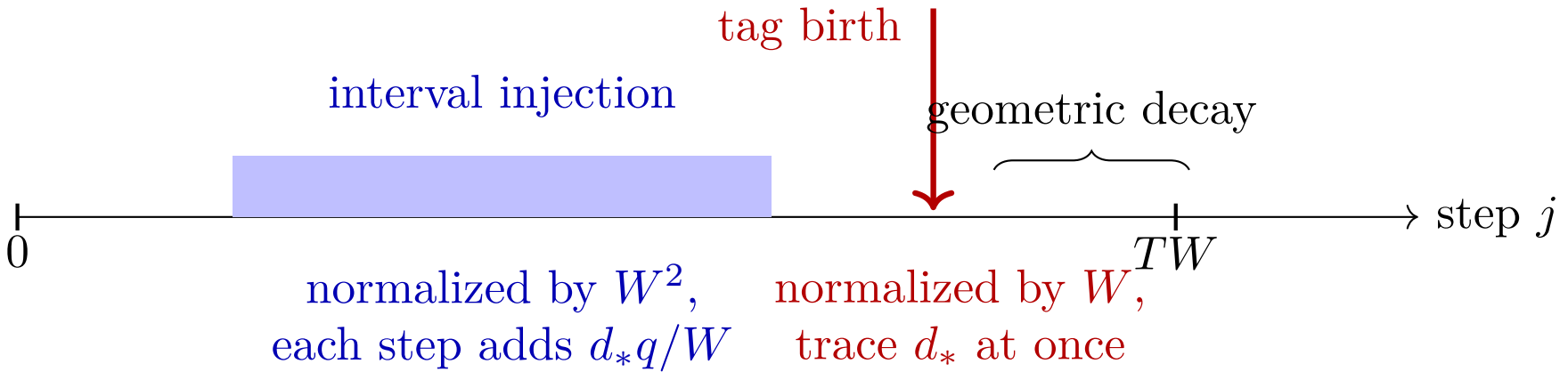
\begin{figure}[H]
\centering
\begin{tikzpicture}[x=1cm,y=1cm]
 \draw[->] (0,0)--(10.4,0) node[right,font=\small] {step $j$};
 \draw[thick] (0,-.1)--(0,.1) node[below=3pt,font=\small] {$0$};
 \draw[thick] (8.6,-.1)--(8.6,.1) node[below=3pt,font=\small] {$TW$};
 \fill[blue!25] (1.6,0) rectangle (5.6,.45);
 \node[font=\small,blue!70!black] at (3.6,.9) {interval injection};
 \node[font=\small,blue!70!black,align=center] at (3.6,-.75)
   {normalized by $W^2$,\\ each step adds $d_*q/W$};
 \draw[very thick,red!70!black,->] (6.8,1.55)--(6.8,.05);
 \node[font=\small,red!70!black,anchor=east] at (6.7,1.4) {tag birth};
 \node[font=\small,red!70!black,align=center] at (7.0,-.75)
   {normalized by $W$,\\ trace $d_*$ at once};
 \draw[decorate,decoration={brace,amplitude=4pt}] (7.25,.35)--(8.7,.35)
   node[midway,above=4pt,font=\small] {geometric decay};
\end{tikzpicture}
\caption{The two kinds of source birth in the proof of
Theorem~\ref{thm:full-source-generator}. An interval source (blue) is
injected continuously along an interval. It is normalized by $W^2$, so
each step adds a small positive matrix of operator size $O(W^{-2})$ and
near trace $d_*q/W$, and the mixed terms with other letters are of
order $W^{-2}$ per step. A tag (red) is born at a single step. It is
normalized by $W$, so it starts with trace $d_*$. Its birth increments
are kept exactly, and after the birth its contracting words decay
geometrically, by Lemma~\ref{lem:t26-contracting-words}.}
\label{fig:fg-births}
\end{figure}

\emph{Step 5: three examples of the rule.} For example, the rule gives
\begin{align*}
 d\Tr Q_a&=dM_a,\\
 d\Tr(Q_aJ)&=-\tfrac12D_0g\Tr(Q_aJ)\,dt+dM_{a,J},\\
 d\Tr(Q_aQ_b)&=\tfrac12D_0
 [\Tr(Q_aP)\Tr(Q_bP')+\Tr(Q_aP')\Tr(Q_bP)]\,dt+dM_{ab}.
\end{align*}
The first line says that the trace of a source has no drift. In the
second line the drift comes from contracting the covariance tensor
$\frac{D_0}{2W}(P\otimes P'+P'\otimes P)$ between $Q_a$ and the current,
exactly as in the third line with $J$ in place of $Q_b$. Its sign follows
from
\[
 \Tr(JP)=-g,\qquad \Tr(JP')=g.
\]
Indeed, these give
\[
 \tfrac12D_0\bigl[\Tr(Q_aP)\Tr(JP')+\Tr(Q_aP')\Tr(JP)\bigr]
 =\tfrac12D_0g\Tr\bigl(Q_a(P-P')\bigr)
 =-\tfrac12D_0g\Tr(Q_aJ),
\]
since $J=P'-P$.
The same-sign contraction of this word vanishes, because
$\Tr(JP\mathbb S)=0$. In the third line the drift is the
opposite-sign contraction of two sources. These checks also show why the
transported current cannot simply be set equal to zero: it enters the
drift of $\Tr(Q_aJ)$ with a nonzero coefficient.
\end{proof}

The estimates of Sections~\ref{sec:transmission}--\ref{sec:full-source}
are used later, also at energies that move with $W$. We record below the order
of their constants.

\begin{lemma}[Uniform constants]
\label{lem:source-recorded-estimates}
Fix $0<b<1/2$ and a compact set $I\subset(-2,2)$ of bulk energies. The
estimates listed below were proved earlier at a fixed energy; this lemma
records that their constants do not depend on the energy within $I$. More
precisely, in each item the constants, including the width threshold
$W_0$, the contraction rates and the levels of uniform integrability,
depend only on $b$, $I$ and the parameters named before them, and they are
chosen before the energy $E\in I$ and the width $W$. Each estimate holds at
every $E\in I$, with the stationary jet $s_*=s_*(E)$ of
\eqref{eq:stationary-jet}. Consequently each estimate also holds along
every deterministic sequence of energies $E_W\to E\in(-2,2)$, at any
rate, by taking for $I$ a compact neighbourhood of $E$.
\begin{enumerate}
\item[(a)] \emph{Passive one-step bounds:} For every deterministic old
 input $X$ with $-\Im X\succeq0$, the one-step mean of $Lf_+L^*$ in the
 proof of Lemma~\ref{lem:passive-source-mixing} holds with uniform
 constants. Moreover, for $p\ge1$, $W\ge W_0(p)$ and $u,v\in\C^W$, the
 centered bilinear forms satisfy
 \[
  \norm{u^*(U-\E(U\mid X))v}_{L^p\mid X}\le C_pW^{-1/2}\norm u\norm v,
  \qquad U\in\{LL^*,\ L(R_+-(2i\eta)^{-1}I)L^*\}.
 \]
\item[(b)] \emph{Covariances:} The covariance formula
 \eqref{eq:passive-covariance}, with its bound on the remainder
 $\mathcal R$ for all four vectors, holds with uniform constants. So does
 \eqref{eq:full-source-covariance}, with its four-vector error, for every
 deterministic Hermitian old two-port matrix. Consequently
 \eqref{eq:passive-source-mixing} holds with a constant
 $C_{p,\epsilon,T}$.
\item[(c)] \emph{Jet gap:} The stationary jet $s_*$ exists at every
 $E\in(-2,2)$, and the gap between the two covariance coefficients is
 uniform on $I$: $\sigma_+-|\sigma_-|\ge c_I>0$.
\item[(d)] \emph{Contracting words and stops:} The exact update and the
 estimates \eqref{eq:t26-word-mixing-rate} and
 \eqref{eq:t26-word-summed-error} hold with uniform constants, and so do
 the assertions of Lemma~\ref{lem:t26-full-source-stops}; in its tail
 bound, the supremum over $E\in I$ and over the tag locations is taken
 inside $\limsup_W$.
\item[(e)] \emph{Clock, tag and edge moments:} The clock bounds
 \eqref{eq:full-clock-moments}--\eqref{eq:full-tag-clock-moments} hold for
 $j\le TW$ with a constant $C_{p,T}$, and \eqref{eq:tag-propagation-moments}
 holds. In addition, a source injected only during a bounded interval has
 bounded trace afterwards: for $R_0>0$ and $p\ge1$ there is a constant
 $C_{p,R_0}$, and then for every $T$ a threshold $W_0(T)$, such that the
 following holds. For $W\ge W_0(T)$, $t_1\ge0$ and $0\le R\le R_0$ with
 $t_1+R+R_0\le T$, every test $0\le q\le1$ vanishing outside
 $[t_1,t_1+R]$, and every $j\le(t_1+R+R_0)W$,
 \begin{equation}\label{eq:edge-clock-finite}
  \norm{\Tr Q_j}_{L^p}\le C_{p,R_0},\qquad
  Q_j=W^{-2}Y_j^q,\quad q_i=q(i/W).
 \end{equation}
 In \eqref{eq:tag-propagation-moments} and \eqref{eq:edge-clock-finite}
 the constant is chosen before $T$ and $t_1$; only the threshold $W_0$
 depends on $T$.
\item[(f)] \emph{Tag birth:} Let $0<\epsilon<T$, and let
 $Q=W^{-1}Y_{k;k}=C_+/W$ be a tag at its birth step $k$. Its trace tends
 to $d_*$ uniformly in the birth step, and its squared trace is small: for
 $p\ge1$,
 \[
  \lim_{W\to\infty}\sup_{\epsilon W\le k\le TW}
      \norm{\Tr Q-d_*}_{L^p}=0,\qquad d_*=\frac1{1-b},
  \qquad\qquad
  \E\Tr Q^2\le C/W\quad(W\ge W_0).
 \]
\item[(g)] \emph{Finite-width fractional bound.} For homogeneous
 transport, the fractional traces of the transmission are bounded at
 finite width: for $0<p<1$ and $\epsilon>0$ there is $C$ such that for
 every $t\ge\epsilon$ there is $W_0(t)$ with
 \[
  \E\Tr T_{\lfloor Wt\rfloor}^p\le C\qquad(W\ge W_0(t)).
 \]
\item[(h)] \emph{Moving energies in the transmission law:} Let $E_W\to E$
 be deterministic, and evaluate the path at the energy $E_W$, including
 the boundary parameter $z=z_\partial(E_W)$ of \eqref{eq:matched-parameters}.
 Then, for every absorption profile $\kappa$, the transmission trace-measure
 processes converge to the limit at the fixed energy $E$, with transforms
 \eqref{eq:transmission-dual-law}, $D_0=D_0(E)$ and $h=4\pi\rho(E)\kappa$.
 For $\kappa\equiv0$, the fractional bound \eqref{eq:fractional-transmission}
 and the exponential bound of Theorem~\ref{thm:transmission-law} hold
 along $E_W$, with constants $C_p,C_q,c_q$ that depend on $b,E,p,q$ but
 not on the sequence $(E_W)$. Finally, for every $\delta_W\to0$, the
 convergence is uniform over deterministic energies $|E'-E|\le\delta_W$.
\end{enumerate}
\end{lemma}

The point of this lemma is that the physical energy is not fixed later in
the paper. It is $e(u)=E+u/(N\rho(E))$, which moves with the width $W$.
All the source, clock and transmission estimates must therefore hold with
constants that are uniform over a compact set of bulk energies, so that
they can be applied along such moving energies.

\begin{proof}
 Most items follow from one observation: the
earlier proofs use the energy only through a few explicit coefficients,
which are continuous on $(-2,2)$, and through estimates whose constants
are already uniform (Step~1). This gives (a), (b), (d) and the clock
moments in (e). The remaining items need a short additional argument
each - the jet gap (c) is an explicit inequality (Step~2), the edge
moments in (e) follow by summing tag moments (Step~3), the tag birth (f)
follows from the injection estimates (Step~4), the fractional bound (g)
from the logistic recursion (Step~5), and the moving-energy statement (h)
from tightness and uniqueness of the limit (Step~6).

\emph{Step 1: how the energy enters.} The proofs cited in the lemma use
the energy only through the following quantities: the parameters
$a,z,z_0,\eta$ of \eqref{eq:matched-parameters}; the coefficients
$D_0,\alpha,\theta_*,s_*,\sigma_\pm$; the bounds for one fresh layer from
Section~\ref{sec:layers} and Lemma~\ref{lem:passive-bounds}, which are
uniform in the deterministic old input; and
Theorem~\ref{thm:matched-traces} and
Lemma~\ref{lem:t26-passive-initial-scalar}, whose constants depend only on
$b$, $I$ and fixed parameters. The parameters and coefficients are
continuous on $(-2,2)$, and $\eta$ is bounded below on $I$. The
inequalities
\[
 |\theta_*|<1,\qquad |1+k_-s_*|\ge1-|k_-|d_*>0
\]
hold with margins that are uniform on $I$
(Subsection~\ref{subsec:s13-the-scalar-derivative-and-a-contracting-}).
Hence the coefficient bounds and the contraction rates are uniform on $I$;
these rates are $|\lambda|\le\max_I|\theta_*|$ in
Lemma~\ref{lem:t26-contracting-words} and $q_*<1$ in
Theorem~\ref{thm:matched-traces}. No step chooses a constant after the
energy, and the exit probabilities over at most $TW$ layers are bounded
uniformly on $I$. This proves (a), (b), (d) and the clock-moment
assertion of (e); the bound \eqref{eq:tag-propagation-moments} is uniform
on compact bulk sets by its derivation.

Two points deserve a separate remark. The fluctuation bound in (a) is the
last assertion of Lemma~\ref{lem:passive-bounds}, whose constants come
from the uniform inverse moments and the Gaussian $L^p$ Poincar\'e
inequality, not from $X$. In (d), Step~2 of the proof of
Lemma~\ref{lem:t26-full-source-stops} uses the fractional bound (g) at
$t=\epsilon$ and $p=1/2$.

\emph{Step 2: the jet gap (c).} Since $|a|=1$,
\[
 |k_-|^2=|c+b^2a^2|^2=k^2-cb^2|a^2-1|^2 ,
\]
and $|a^2-1|$ is bounded below on $I$, because $a^2=1$ only at the edges
of the spectrum. Hence $|k_-|\le k-c_I$. The function $x\mapsto
x/(1-xd_*)$ has derivative at least one on $[0,k]$, so
\eqref{eq:jet-gap} gives
\[
 \sigma_+-|\sigma_-|
 \ge\frac{k}{1-kd_*}-\frac{k-c_I}{1-(k-c_I)d_*}
 \ge c_I .
\]

\emph{Step 3: the edge moments in (e).} A source injected on an interval is
a sum of tags. Indeed $Y_j^q=\sum_{i\le j}q_iY_{j;i}$, so
\[
 \Tr Q_j=W^{-1}\sum_{i\le j}q_i\,W^{-1}\Tr Y_{j;i}.
\]
This sum has at most $RW+1$ nonzero terms, since $q$ vanishes outside an
interval of length $R$, and each tag has travelled for $j-i\le(R+R_0)W\le
2R_0W$ steps. Minkowski's inequality and \eqref{eq:tag-propagation-moments}
with $R=2R_0$ give \eqref{eq:edge-clock-finite}, with
$C_{p,R_0}=(R_0+1)C_{p,2R_0}$, for $W\ge W_0(p,2R_0,T)$.

\emph{Step 4: the tag birth (f).} By Step~2 of the completion of the proof
of Theorem~\ref{thm:full-source-generator}, on the scalar domain the near
part of the injection $C_+$ has normalized trace
\[
 d_*+O(\delta_0+g/W)+O_{L^p}(W^{-1}),
\]
and its far part has trace $O_{L^p}(g)$, so it contributes $O_{L^p}(g/W)$
to $\Tr Q$. Here
\[
 g=\tfrac12\Tr J^2=\Tr(I-\mathbb S_{11}^*\mathbb S_{11})=\Tr T ,
\]
since $\mathbb S_{11}=S$. For $k\ge\epsilon W$,
Lemma~\ref{lem:t26-passive-initial-scalar} gives
$\norm{\delta_0}_{L^p}\le C_{p,T}/W$, and \eqref{eq:conductance-moments}
gives $\E g^p\le C_{p,T}(1+\epsilon^{-p})$. The event of exit from the
scalar domain is paid for by $\norm{C_+}_{L^{2p}}\le C_p$, H\"older's
inequality and \eqref{eq:passive-superexp-exit}. This gives the limit,
uniformly in $k$. Finally, $C_+$ has rank at most $W$, so
$\Tr C_+^2\le W\norm{C_+}^2$, and therefore $\E\Tr Q^2\le C/W$.

\emph{Step 5: the fractional bound (g).} This is the logistic recursion for
$\E\Tr T^p$ in
Subsection~\ref{subsec:s14-fractional-bounds-for-homogeneous-transm}, run
over the last $\lfloor\epsilon W/2\rfloor$ layers before $\lfloor Wt\rfloor$.
Its inputs are the clauses following \eqref{eq:conductance-moments} and
\eqref{eq:passive-jet-moments}, which are uniform in the length of the
path. Its coefficients are uniform on $I$, and it holds at finite width.

\emph{Step 6: moving energies (h).} By (a), (e), (g) and the conductance
moments, the proofs of Theorems~\ref{thm:transmission-generator}
and~\ref{thm:transmission-law} apply along $E_W$ with constants
independent of $W$, and $D_0(E_W)\to D_0(E)$. With absorption, the forcing
coefficient is a continuous function of the energy times $\kappa$, and it
tends to $h$. Hence the processes are tight on compact intervals of
positive time, and every subsequential limit solves the martingale
problem \eqref{eq:transmission-drift}--\eqref{eq:transmission-bracket}
with the coefficient $D_0(E)$ and the entrance
\eqref{eq:transmission-entrance}. So every subsequential limit satisfies
\eqref{eq:transmission-dual-law} at the energy $E$, which determines the
law, and the whole sequence converges. The bounds hold with constants
independent of $(E_W)$, since, as before, they follow from the uniform
finite-width estimates and from the limiting generator at $E$.

Finally, the uniformity over $|E'-E|\le\delta_W$ follows by contradiction.
If the convergence of a fixed bounded test were not uniform over these
energies, we could choose, for each $W$, an energy $E_W$ with
$|E_W-E|\le\delta_W$ at which the test is nearly as far from its limit as
possible. This gives a deterministic sequence $E_W\to E$ along which the
convergence fails, contradicting what we have just proved.
\end{proof}
\section{Positive finite source laws and their physical limit}
\label{sec:positive-source-law}
The generator of Section~\ref{sec:full-source} acts on trace observables
and does not by itself define a process. This section constructs a positive
stochastic process with this generator, and identifies it with the limit of
the physical path.
Proposition~\ref{prop:finite-gram-law} states that the finite stochastic
equation \eqref{eq:finite-source-sde} has a unique nonexplosive solution, that
it preserves a pair of projections, and that its Gram law depends only on its
initial data. Theorem~\ref{thm:source-process} states that the bounded
Gaussian source transforms of the path have a unique critical limit, given by
this equation together with the transmission entrance and independent probe
births.
In Subsection~\ref{subsec:s16-the-finite-stochastic-equation} we solve the finite stochastic equation. In
Subsection~\ref{subsec:s16-a-comparison-independent-of-retained-ran} we couple constructions with different numbers of
retained channels, with a bound that does not depend on the dimension. In
Subsection~\ref{subsec:s16-physical-tails-and-closing-the-omitted-c} we compare with the actual path and control the mass lost
in omitted channels. In Subsection~\ref{subsec:s16-sampling-sources-and-identification} we pass from local sources to
sources spread over an interval, and identify the limit.

At a positive time $s$, sample the transmission sequence $(t_i)$
from Theorem~\ref{thm:transmission-law}. Its canonical pair of
projections is the orthogonal sum of
\[
 P_i=\begin{pmatrix}1&0\\0&0\end{pmatrix},\qquad
 P_i'=\begin{pmatrix}1-t_i&\sqrt{t_i(1-t_i)}\\
             \sqrt{t_i(1-t_i)}&t_i\end{pmatrix}.
\]
The difference has eigenvalues $\pm\sqrt{t_i}$.
Retain finitely many pairs. Choose a frame $V$ and a fixed
Hermitian signature matrix $\mathsf D$ with eigenvalues in
$\{-1,0,1\}$ so that $J=V\mathsf DV^*=P'-P$.
The two current columns for mode $t_i$ have norm $t_i^{1/4}$.
Probe columns have signature zero.

\begin{theorem}\label{thm:source-process}
For any finite family of positive source tests supported in
$[s,T]$, the bounded Gaussian source transforms of the actual
two-input path have a unique critical limit. It is specified by
the finite stochastic equation \eqref{eq:finite-source-sde},
the preceding transmission entrance, and independent near probe
births of squared norm $d_*m/M$ at $M$ sampled source locations.
First remove the current cutoff, then send $M\to\infty$.
The limit is positive, consistent under deleting tests or probes,
and independent of the cutoff sequence. Sources before $s$ can
then be restored by $s\downarrow0$.
\end{theorem}

Bounded Gaussian source transforms mean expectations of finite
products of $\det(I+\sum_a a_aQ_a)^{-1}$, $a_a\ge0$,
with independent final Gaussian probes for distinct factors.
They determine the corresponding positive Gaussian response measures.
The construction does not use determination by unbounded moments.

\subsection{The finite stochastic equation}\label{subsec:s16-the-finite-stochastic-equation}

 We realize the source law as an ordinary finite-dimensional stochastic
equation. Proposition~\ref{prop:finite-gram-law} proves existence,
preserves the two-projection constraint, and gives the law of all
inner products of the source vectors. We specify the Gaussian noise,
compute these inner-product equations, and verify positivity and
absence of explosion.

 In a finite ambient space   $\C^{2\Naux }$  let $P$ have rank   $\Naux $  and
let $P'=P+V\mathsf DV^*$ also be a rank-  $\Naux $  projection.
Let $Z$ have independent complex Brownian entries normalized by
$\E|dZ_{ij}|^2=2\,dt$. Define

\begin{equation}\label{eq:finite-source-sde}
 dA=\frac{\sqrt{D_0}}{2}P\,dZ\,P',\qquad
 dV=dA^*V-\frac{D_0\Naux }{4}PV\,dt.
\end{equation}

Write $M=V^*V$, $B=V^*PV$ and
\[
 C_0=B+M\mathsf DM,\quad
 C_1=B+M\mathsf DB,\quad
 C_2=B+B\mathsf DB,\quad
 g=\tfrac12\Tr(\mathsf DM\mathsf DM).
\]

\begin{proposition}\label{prop:finite-gram-law}
The finite stochastic equation \eqref{eq:finite-source-sde} has a unique
solution that exists for all times, that is, it does not explode, and
along which $P'$ remains an orthogonal projection of the same rank as
$P$. Its \emph{Gram law}, the law of the two matrices of inner products
$M=V^*V$ and $B=V^*PV$, depends only on the initial data $(M_0,B_0)$ and
the signature $\mathsf D$. The conditional means and the covariances of
the increments are
\begin{align}
 \E_t\,dM&=0,&\E_t\,dB&=-\tfrac12D_0gB\,dt,\label{eq:gram-drift}\\
 dM_{ij}\,dM_{kl}
 &=\tfrac{D_0}{2}[B_{il}(C_0)_{kj}+(C_0)_{il}B_{kj}]\,dt,\notag\\
 dM_{ij}\,dB_{kl}
 &=\tfrac{D_0}{2}[B_{il}(C_1^*)_{kj}+(C_1)_{il}B_{kj}]\,dt,\notag\\
 dB_{ij}\,dB_{kl}
 &=\tfrac{D_0}{2}[B_{il}(C_2)_{kj}+(C_2)_{il}B_{kj}]\,dt.
                                                   \label{eq:gram-covariance}
\end{align}
Here a product such as $dM_{ij}\,dB_{kl}$ denotes a bilinear covariation,
taken without complex conjugation. Finally, the Gram law does not change
if unused near or far dimensions are added to the ambient space.
\end{proposition}
Uniqueness of the Gram law holds in the following class. Call a pair
$(M,B)$ a \emph{feasible Gram state} for $\mathsf D$ if $M=V^*V$ and
$B=V^*PV$ for a frame $V$ in some finite-dimensional ambient space with
an orthogonal projection $P$ such that $P+V\mathsf DV^*$ is an
orthogonal projection of the same rank as $P$. Every continuous adapted
process with values in the feasible Gram states which solves the
martingale problem of \eqref{eq:gram-drift}--\eqref{eq:gram-covariance}
up to each stop $\Tr M\ge K$, and has the initial law of $(M_0,B_0)$,
has this Gram law.

Every feasible Gram state satisfies $M\succeq0$ and $0\preceq B\preceq M$,
since $B=(PV)^*(PV)$ and $M-B=((I-P)V)^*((I-P)V)$. The feasible states
form a closed set, so feasibility passes to limits of Gram processes.
Indeed, let $V$ have $n_V$ columns and put
$X=\operatorname{ran}V+\operatorname{ran}PV$, of dimension at most
$2n_V$. With $P'=P+V\mathsf DV^*$ we have $P'V=PV+V\mathsf DM$ and
$P'PV=PV+V\mathsf DB$, so $X$ is invariant under $P$ and $P'$, and
$P'=P$ on $X^\perp$. Hence the restrictions of $P$ and $P'$ to $X$ are orthogonal projections
of the same rank realizing the same Gram state, and, adding zero
dimensions, every feasible state is realized in $\C^{2n_V}$ by a frame
with $\|V\|_{\rm HS}^2=\Tr M$. For converging feasible states these frames
and projections have a convergent subsequence whose limit satisfies the
same closed conditions.

All coefficients in \eqref{eq:gram-drift}--\eqref{eq:gram-covariance}
are expressed through the finite Gram matrices themselves, and none of
them involves the dimension of the ambient space. This independence from
unused dimensions is what allows us to remove the cutoff on the number of
retained transmission channels later in this section.

\begin{proof}
The proof has three parts. First we compute the
drift and the covariance of the Gram matrices directly from the equation
(Steps~1 and~3), and we show that the projection constraint is preserved
and that solutions do not explode (Step~2). Second, we prove that the
Gram law is unique. The difficulty is that the Gram matrices can be
singular, and their covariance is degenerate, so standard uniqueness
results for diffusions with a nondegenerate covariance do not apply. We
therefore show that the ranks of the two Gram blocks can neither increase
(Step~5) nor decrease (Step~4), and that on each set of fixed ranks the
Gram process can be compared with the ambient solution through an
explicit global factor (Step~6). Third, we remove the stopping times used
along the way and check that padding with unused dimensions changes
nothing (Step~7).

\emph{Step 1: the drift.} The noise of the equation enters through
$d\mathcal E=dA+dA^*$, whose covariance is
$D_0(P\otimes P'+P'\otimes P)/2$. Two contractions of this noise give
the drift. The first, $dA\,dA^*=D_0\Naux P\,dt/2$, exactly cancels the
norm drift in \eqref{eq:finite-source-sde}; hence $M$ has no drift. The
second, $dA\,P\,dA^*=D_0(\Naux-g)P\,dt/2$, gives the drift of $B$ in
\eqref{eq:gram-drift}.

\emph{Step 2: the projection constraint, existence and nonexplosion.} We
check directly that $P'$ remains a projection. Let $Q'=I-P'$. The noise
and the drift of the current $J=P'-P$ are
\[
 dJ_{\rm noise}=-\frac{\sqrt{D_0}}{2}
       (P'\,dZ^*P Q'+Q'P\,dZ\,P'),\qquad
 K_J=-\frac{D_0\Naux }{4}(PJ+JP)-\frac{D_0g}{2}P'.
\]
In the decomposition $P'\oplus Q'$, the noise has zero diagonal blocks;
it only rotates $P'$, as a tangent direction should. In the $P'$ block,
the drift $-D_0gP'/2$ is cancelled by the It\^o correction
$P'(dJ)^2P'/dt=D_0gP'/2$. In the $Q'$ block, the drift
$D_0\Naux Q'PQ'/2$ equals the It\^o correction $Q'(dJ)^2Q'/dt$. These are
exactly the two diagonal blocks of the It\^o equation for the identity
$(P')^2=P'$, and its off-diagonal blocks follow from the form of the
noise. So the constraint persists.

The coefficients are smooth functions of the finitely many entries, and
they are bounded as functions of the constrained current. Additional
columns of the frame solve linear equations. This gives local existence
and pathwise uniqueness. To exclude explosion, we use that the trace
$\Tr M$ is a nonnegative martingale with fixed moment bounds on bounded
time intervals.

\emph{Step 3: the covariances.} The entries of the two Gram matrices are
\[
 dM_{ij}=v_i^*d\mathcal E\,v_j,\qquad
 dB_{ij}=v_i^*P\,d\mathcal E\,Pv_j,
\]
where $v_i$ are the columns of $V$. Pairing the noise covariance with
these frames produces four Gram matrices involving $P'$, namely
\[
 V^*P'V=C_0,\quad V^*P'PV=C_1,\quad
 V^*PP'V=C_1^*,\quad V^*PP'PV=C_2 .
\]
Inserting them into the crossed covariance of Step~1 gives each line of
\eqref{eq:gram-covariance}, including the adjoint on $C_1$. In
particular, every trace contraction agrees with the rule
\eqref{eq:word-generator}--\eqref{eq:word-covariance}.

\emph{Step 4: the ambient solution keeps its ranks.} For uniqueness we
first study the ranks of the two Gram blocks. Split the frame into its
near and far rows, $V=(L^T,F^T)^T$. Then $B=L^*L$ and $C:=M-B=F^*F$. The
near rows solve a linear stochastic equation with the coefficient acting
from the left, and the far rows solve
\[
 dF=(\sqrt{D_0}/2)\,F\,\mathsf DV^*dZ^*PV,
\]
a linear equation with the coefficient acting from the right. Each
solution can be written as a stochastic fundamental matrix, the solution
started from the identity, applied to its initial value. On bounded
localized time intervals these fundamental matrices are invertible. Hence
the ranks of $B$ and $C$ stay constant along the ambient solution.

\emph{Step 5: no solution can create rank.} A different solution of the
Gram martingale problem, that is, any process with values in the feasible Gram states and with
the drift and covariance of \eqref{eq:gram-drift}--\eqref{eq:gram-covariance}, might
in principle leave a singular state immediately. We exclude this with a
barrier function. For a positive matrix $A$, let $e_k(A)$ be the $k$th
elementary symmetric function of its eigenvalues; $e_k(A)>0$ exactly when
the rank of $A$ is at least $k$. If $A=U^*U$, the Cauchy--Binet formula
gives $e_k(A)=\|\bigwedge^k U\|_{\rm HS}^2$, where $\bigwedge^kU$ is the
$k$th exterior power of $U$, also when $A$ is singular.

We stop when $\Tr M$ exceeds a level $K$. Write the near equation as
$dL=\sum_\alpha U_\alpha L\,dB_\alpha+U_0L\,dt$ and the far equation as
$dF=\sum_\alpha F A_\alpha\,dB_\alpha$, with real Brownian motions
$B_\alpha$. For each fixed ambient dimension and number of columns, the
coefficient matrices $U_\alpha$, $A_\alpha$ and $U_0$ have bounded squared
sums. This follows directly from $P'=P+V\mathsf DV^*$ and
$\|V\|_{\rm HS}^2\le K$, without inverting $L$ or $F$. For example, the
near noise is the near block of $(\sqrt{D_0}/2)PP'dZ^*P$ multiplying
$L$, and $A_\alpha$ is the corresponding real coordinate of
$(\sqrt{D_0}/2)\mathsf DV^*dZ^*PV$.

The exterior power is multilinear, so $\bigwedge^k L$ and $\bigwedge^k F$
also solve linear stochastic equations with bounded coefficients. The
first-order term inserts one coefficient matrix in one tensor slot, and
the It\^o correction inserts two in different slots; every term still
factors through the same exterior power. Applying It\^o's formula to the
squared Hilbert--Schmidt norm therefore proves the inequalities
\begin{equation}\label{eq:t26-gram-rank-barrier}
 \mathscr G e_k(B)\le C_{K,k,\Naux }e_k(B),\qquad
 \mathscr G e_k(C)\le C_{K,k,\Naux }e_k(C),
\end{equation}
where $\mathscr G$ is the generator of the Gram process. These
inequalities involve only the drift and covariance
\eqref{eq:gram-drift}--\eqref{eq:gram-covariance}. Every possible Gram
state has a finite factorization, so the calculation applies at every
state, including states on the boundary where the rank drops, and the
inequalities hold for every solution, not only for the ambient one. If
$e_k(B_0)=0$, then positivity of $e_k$ and Gronwall's inequality for the
stopped process give $\E e_k(B_{t\wedge\sigma_K})=0$; the same holds for
$C$. Hence no rank can increase. The constants here may depend on the
padding dimension, which is harmless at this point.

\emph{Step 6: comparison on fixed ranks.} It remains to compare an
arbitrary solution with the ambient one while the ranks stay fixed. We
localize on the event that the trace is bounded and that the nonzero
eigenvalues of $B$ and $M-B$ stay above a level $\delta>0$, which we call
a spectral-gap localization. On such an event we use the global factor
\[
 V=\bigl[h_\delta(B);\,h_\delta(M-B)\bigr],\qquad
 P_0=\operatorname{diag}(I,0),
\]
where $h_\delta$ is a Lipschitz function that equals the square root on
$[\delta,\infty)$ and vanishes at zero. On the localized event this
factor reproduces the Gram data: $V^*V=M$ and $V^*P_0V=B$. The brackets of
this factor depend only on $(V^*V,V^*P_0V)$, that is, on the Gram data
alone, and these bracket identities give the comparison of the stopped
Gram laws. Two points explain why this global factor is used. First, a
single chart based on one invertible principal minor does not cover the
whole set of matrices of fixed rank. Second, the Gram covariance is
degenerate, so uniqueness cannot be obtained from ellipticity, the
nondegeneracy of the covariance. The localizations of the trace and of
the positive Gram eigenvalues are put in place before the stops are
removed in the next step.

\emph{Step 7: removing the stops and adding dimensions.} By Step~5 no
solution creates rank, so until a possible loss of rank every solution
has the same stopped law as the ambient solution, by Step~6. The ambient
solution never loses rank in finite time, because its stochastic
fundamental matrices are invertible (Step~4), and the trace localization
excludes explosion (Step~2). These facts remove the stops for every
solution. This proves uniqueness of the Gram law, also for initially
singular matrices, on each set of fixed ranks. Finally, the coefficients
in \eqref{eq:gram-drift}--\eqref{eq:gram-covariance} contain no
parameter for unused dimensions, so padding the ambient space gives the
same law.
\end{proof}

\subsection{A comparison independent of retained rank}\label{subsec:s16-a-comparison-independent-of-retained-ran}

 We need the finite construction to converge as more channels are
retained. The coupling estimate below has no dimension factor in its
drift, so nearby initial source frames stay close up to a controlled
trace cutoff. We drive both frames with the same noise, exhibit the
cancellation, and remove the cutoff using their positive trace bounds.

 Couple two models using the same $Z$ in a common padded space
and the same signature $\mathsf D$. Put
$X=\norm{V_1-V_2}_{\HS}^2$,
$R=\norm{V_1}_{\HS}^2+\norm{V_2}_{\HS}^2$ and
$K_J=P_1'-P_2'$.
The exact cross contraction is
\[
 dA_1\,dA_2^*=\tfrac{D_0}{2}\Tr(P_1'P_2')P\,dt.
\]
Since
  $\Tr(P_1'P_2')=\Naux -\norm{K_J}_{\HS}^2/2$, every factor   $\Naux $
 cancels in the mixed Gram drift:
\[
 \E_t\,d(V_1^*V_2)
       =-\tfrac{D_0}{4}\norm{K_J}_{\HS}^2V_1^*PV_2\,dt,\quad
 \E_t\,dX=\tfrac{D_0}{2}\norm{K_J}_{\HS}^2
                         \Re\Tr(V_1^*PV_2)\,dt.
\]
The elementary bounds
$\norm{K_J}_{\HS}^2\le2RX$ and
$|\Tr(V_1^*PV_2)|\le R/2$ give upper drift $D_0R^2X/2$.
Stopping at $R=K$, the maximal inequality for
$e^{-D_0K^2t/2}X_t$ yields
\begin{equation}\label{eq:rank-free-coupling}
 \Prob\{\sup_{t\le T}X_t>u\}
 \le\frac{\E R_0+2B_{\rm birth}}K
       +\frac{e^{D_0K^2T/2}\E X_0}{u}.
\end{equation}
Here $B_{\rm birth}$ is the total squared norm of finitely many
shared births. A shared birth is chosen perpendicular to both
old frames in unused near dimensions; it creates zero distance.
Equivariance and padding give the prescribed marginal law.
On $R\le K$, Schatten Cauchy--Schwarz gives
\[
 \norm{M_1-M_2}_1+\norm{B_1-B_2}_1
             \le2\sqrt{2KX},
\]
and the same bound for the current difference.
There is no dependence on the number or size of retained modes.

For nested current cutoffs,
$\E X_0=2\E\sum_{\text{omitted }i}\sqrt{t_i}\to0$ by
Theorem~\ref{thm:transmission-law}, while initial frame traces
have uniformly bounded first moment. First fix $K$, remove the
cutoff in \eqref{eq:rank-free-coupling}, then send $K\to\infty$.
This makes the finite Gram paths a Cauchy family in probability
under these couplings, in the uniform trace-norm metric.
The continuous trace-class path space is complete and separable.
It follows that their laws have a unique limit. For nonnested
cutoffs compare each with their union.
Deleting a probe does not change $P'$ because its signature is zero,
so positivity and consistency pass to this completion.

\subsection{Physical tails and closing the omitted channels}\label{subsec:s16-physical-tails-and-closing-the-omitted-c}

 We compare a finite-channel approximation with the actual path, rather
than only with another auxiliary process. The estimates below control
the current and source mass lost when omitted far channels are closed.
We dominate their signed current by a transported positive matrix,
average the closing phase, and use the small transmission tail.

 At time $s$, let $A_-,A_+$ be the paired far input and output
projections for the omitted transmission modes. Set
$J_A=A_--\mathbb S^*A_+\mathbb S$.
Initially
$\norm{J_A(s)}_1=2\sum_{i\in A}\sqrt{t_i}$.
The fixed far amplitudes give
$J_A(t)=N_{s,t}^*J_A(s)N_{s,t}$.
Transport the positive matrix $D_s=|J_A(s)|$ by the same
congruence. Exactly
$-D_t\preceq J_A(t)\preceq D_t$.
The mixed-word argument applied to this past-measurable positive
matrix, after trace normalization and a vanishing mixing interval,
makes $\Tr D_t$ a nonnegative local martingale in each critical
limit. Therefore
\begin{equation}\label{eq:physical-current-carrier}
 \Prob\{\sup_{s\le t\le T}\Tr D_t>v\}
                       \le\E\Tr D_s/v.
\end{equation}
This also holds as a critical limsup by the stopped finite-width
calculation. It requires only the first moment of the initial
nuclear tail.

Let $Q$ be a positive source born after $s$, $\gamma=\Tr Q$,
and $\ell_A=\Tr(QA_-)$. Its incoming leakage is bounded by
\[
 g_A^{\rm in}:=\Tr(P\mathbb SA_-\mathbb S^*)
 \le\Tr[A_-(I-\mathbb S^*A_+\mathbb S)A_-]
 \le\norm{J_A(t)}_1\le\Tr D_t.
\]
Here $P\preceq I-A_+$ is an output-space inequality.
Writing $N=\left(\begin{smallmatrix}L&C\\0&I\end{smallmatrix}\right)$,
$C=\Delta\mathbb R_{12}$, gives exactly
\[
 NA_-N^*-A_-=
 \begin{pmatrix}CA_-C^*&CA_-\\A_-C^*&0\end{pmatrix}.
\]
For positive old $X$, the bilinear variance of $\Delta$ gives
$\E_j\Delta X\Delta^*\preceq
C[\Tr X/W+\delta_0^2\norm X]I$.
Use $X=\mathbb R_{12}A_-\mathbb R_{12}^*$,
whose trace is $g_A^{\rm in}/(4\eta^2)$, to obtain
\begin{equation}\label{eq:source-leakage-recursion}
 \E_j\ell_{A,+}-\ell_A
 \le C\gamma[(\delta_0+W^{-2})\sqrt{g_A^{\rm in}}
                    +(W^{-1}+\delta_0^2)g_A^{\rm in}].
\end{equation}
The cross term is bounded using its actual mean; it is not
declared centered. A normalized tag's far birth mass is
$O(g_A^{\rm in}/W)$ by the weighted inverse bound.

Stop at $\gamma=K$ or $\Tr D=v$. Since
$\E\delta_0\le C_s/W$ and $\E\delta_0^2\le C_s/W^2$,
the sum of \eqref{eq:source-leakage-recursion} is at most
$C_{s,T}K(\sqrt v+v)$. The carrier exit is paid by
\eqref{eq:physical-current-carrier}; source trace exits are paid
by fixed moments. Since $\ell_A\le\gamma$, a second moment
pays their final contribution. First remove the current cutoff
at fixed $K,v$, then let $v\downarrow0$, then $K\to\infty$.
We have proved
\begin{equation}\label{eq:small-channel-leakage}
 \lim_{\text{cutoff}}\limsup_W\E\Tr(Q(T)A_-)=0 .
\end{equation}
The proof is uniform for finite positive source families and their
sampled averages, by positivity and Jensen.

At final time partition scattering inputs into omitted $A$ and
retained $B$. Write $r=S_{AA}$, $t=S_{AB}$, and let $F_A,F_B$
map these inputs to the chosen source coordinates.
Close $A$ with one independent uniform scalar phase $z$:
\[
 F_B^{\rm cl}(z)=F_B+F_Az(I-rz)^{-1}t.
\]
For the finite actual path $\norm r<1$ almost surely.
Unitarity gives $tt^*=I-rr^*$, so the Hardy-series identity
$\sum_{k\ge0}r^ktt^*(r^*)^k=I$ and phase orthogonality give
\begin{equation}\label{eq:gaussian-closing-error}
 \E_{z,h_A,h_B}
 \norm{F_B^{\rm cl}(z)h_B-(F_Ah_A+F_Bh_B)}^2
                         =2\Tr(F_AF_A^*).
\end{equation}
Both fields use the same retained Gaussian $h_B$; the removed
field $F_Ah_A$ is independent and has the same covariance as the
phase correction. Its cross term vanishes.
The function $v\mapsto e^{-\norm v^2}$ is globally Lipschitz.
Equations \eqref{eq:small-channel-leakage} and
\eqref{eq:gaussian-closing-error} therefore pay every fixed
family of positive Gaussian transforms, with the same phase for
all probes.

Closing at time $s$ or at final time eliminates the same linear
equations. In the canonical two-projection basis at time $s$,
each discarded pair becomes a unitary near reflector, while the
retained transmission eigenvalues are unchanged.
Its near reflection correction has modulus
$|tz/(1-\sqrt{1-t}\,z)|$.
For $0<p<1$, its phase mean is at most $C_pt^p$, using
$t(1+|\log t|)\le C_pt^p$.
The normalized near trace therefore changes by at most
$C_pW^{-1}\sum_{i\in A}t_i^p=o(1)$ in mean.
This checks the scalar regularity domain of the modified boundary.
After $C\log W$ source-free layers, its jets are stationary
at the $W^{-1}$ scale. The finite coordinates change by $o(1)$:
the initial scalar error is $o(1)$ and decays geometrically, and
the centered increments have total variance $O(\log W/W)$.

\subsection{Sampling sources and identification}\label{subsec:s16-sampling-sources-and-identification}

 We pass from finitely many local sources to a source spread over an
interval. The determinant comparison below has error $O(M^{-1})$ when
$M$ tag locations are sampled. We condition on the path, represent the
source by Gaussian probes, bound the sampling variance, and combine
this approximation with the finite-source law to complete its
identification.

 Write an interval source as
$Q(q)=\sum_j(q_j/W)Q_j$, where $Q_j=Y_{n;j}/W$ is a
transported tag. Put $m=\sum_jq_j/W$.
Sample $M$ independent tag locations with probabilities
$q_j/(Wm)$; at each tag inject an independent Gaussian vector
with that tag covariance and weight $m/M$.
Let $\widehat Q_M$ be the sum of their final rank-one matrices.
Condition on the actual path and an independent final Gaussian $h$.
A sample $Z$ has mean $h^*Q(q)h$ and second moment bounded by
$2m^2\E_j(h^*Q_jh)^2$.
Taylor's theorem for $e^{-x}$ and Gaussian averaging yield
\begin{equation}\label{eq:positive-source-sampling}
 \left|\E\det(I+a\widehat Q_M)^{-1}
             -\E\det(I+aQ(q))^{-1}\right|
 \le\frac{C_{\tau,q,a}}M.
\end{equation}
Indeed $\E_h(h^*Q_jh)^2\le2(\Tr Q_j)^2$, controlled by
\eqref{eq:full-tag-clock-moments}. Common samples work for a finite
family of spatial tests; independent final probes give products.
At a positive birth time the normalized injection has trace and
near trace tending to $d_*$, squared trace $O(W^{-1})$ and
operator moments $O(W^{-1})$. Gaussian bilinear concentration
makes its new column have deterministic squared norm $d_*m/M$
and zero overlap with every old finite frame and current.

\begin{proof}[Proof of Theorem~\ref{thm:source-process}]
Retain $t_i>h^2$ at time $s$.  The level $h$ is chosen off the countable set $\mathcal N_s$ defined
in the proof of Theorem~\ref{thm:complex-source-law}, at which the
entrance identification of its Step~9 is not available; the removal of
$h$ below runs through the remaining levels. For $0<p<1/2$,
$\E\Tr D_s\le C_{p,s}h^{1-2p}$; their number is tight.
At fixed $h,M$ and a finite cap on that number, close the omitted
channels as above. Subsequent layers remain independent Gaussians
conditional on the old state and phase.
The full-source calculation uses the near rank $W$, passivity and
the projection identity, and thus applies with finitely many far
open channels; unused far coordinates can be padded by decoupled
ones. Its errors are uniform at this fixed rank and probe number.
The verified scalar entrance and the deterministic birth data
therefore identify the physical finite-coordinate martingales with
the unique finite Gram law of Proposition~\ref{prop:finite-gram-law}.
In particular its three covariances are the explicit contractions
\eqref{eq:gram-covariance}, rather than a conclusion inferred from
matching only source traces.

Remove the current-count cap by the homogeneous fractional trace bound.
The physical path in this source construction has no absorption;
its source clocks record derivatives and do not inject into transmission.
 Equations \eqref{eq:small-channel-leakage} and
\eqref{eq:gaussian-closing-error} remove the physical closing error
as $h\downarrow0$, uniformly in the source sample size.
The rank-independent coupling removes $h$ in the finite positive
models at each fixed $M$.
Equation~\eqref{eq:positive-source-sampling} then proves that their
$M\to\infty$ transforms are Cauchy and coincide with the actual
ones. Sources before $s$ have expected positive trace $O(s)$
by \eqref{eq:full-tag-clock-moments}, so restoring them changes a
bounded positive transform by $O(s)$. This pays $s\downarrow0$.
All finite stages are positive and consistent, and the limits
preserve these properties. Positive spatial test families determine
the law of a finite positive response measure by their Laplace
functionals; tightness follows from the same trace bounds.
\end{proof}
\section{Complex finite sources and Gaussian replacement}
\label{sec:complex-sources}
The probe matrix of the ring needs the real part of the resolvent quadratic
forms, and not only their positive imaginary part. This section extends the
source law of Section~\ref{sec:positive-source-law} to complex probe
matrices.
Theorem~\ref{thm:complex-source-law} states that the joint law of the
complex probe matrix, the source Gram matrices and the current has a unique
critical limit. It is given by the finite source equation with one added
equation for the complex matrix.
In Subsection~\ref{subsec:s17-the-real-part-of-the-closing-error} we control the real part of the error made by
closing omitted channels, through fractional moments derived from its
imaginary part. In Subsection~\ref{subsec:s17-uniform-raw-moments-and-extended-probes} we propagate uniform raw moments and
use them to replace extended Gaussian probes by finitely many local
sources.

Let $G$ be the inverse of the full path with its two matched leads.
Fix finitely many independent Gaussian source columns
$\xi_j=\sqrt{m_j/W}\,\iota_{k_j}g_j$ on physical blocks, where
$g_j\sim{\rm CN}(0,I_W)$ and $m_j\ge0$ is fixed.
Let $\mathcal J$ collect these columns and let $V_{\rm pr}$ be
their boundary response frame. Define
\begin{equation}\label{eq:complex-source-data}
 \mathcal K=\mathcal J^*G\mathcal J,\qquad
 M_{\rm pr}=V_{\rm pr}^*V_{\rm pr},\qquad
 \Im\mathcal K=-\eta M_{\rm pr}.
\end{equation}
The last identity is the resolvent imaginary-part identity and
fixes all normalizations.

\begin{theorem}\label{thm:complex-source-law}
For fixed distinct positive continuum birth times, the actual
joint law of $\mathcal K$, the source Gram matrices and the current
has a unique critical limit. Its positive construction is
\eqref{eq:finite-source-sde}, completed as in
Theorem~\ref{thm:source-process}, with the added equation
\begin{equation}\label{eq:complex-source-sde}
 d\mathcal K=-2i\eta V_{\rm pr}^*dA\,V_{\rm pr}.
\end{equation}
At a birth of mass $m$, append diagonal entry $am$ and zero cross
entries to $\mathcal K$; the new Gram diagonal is $d_*m$.
All fixed source moments are bounded on bounded time intervals.
\end{theorem}

Since $\Im a=-\eta d_*$, the births preserve
\eqref{eq:complex-source-data}. The real-part of noise is coupled to
the same $dA$ as the Gram process; it is not an independent noise.

\begin{proof}[Proof of Theorem~\ref{thm:complex-source-law}]
The probe matrix $\mathcal K$ collects the
resolvent quadratic forms of the probes. Its imaginary part is the Gram
matrix of the probe responses, whose law was constructed in
Section~\ref{sec:positive-source-law}; its real part is new. We first
write the exact change of $\mathcal K$ when one layer is added, and see
that it is driven by the same noise as the Gram process (Step~1). The
rest of the proof makes the passage to the limit rigorous. We work under
a stop that bounds the trace of the closed frame, and we check that the
frame factors entering every increment are then bounded (Step~2). We
prove fixed fractional moments of the transmission (Step~3), which
control the shift caused by closing the omitted channels and keep the
scalar jet small along the path (Step~4). This bounds the increments of
the Gram coordinates (Step~5). At a birth of a probe, we estimate the new
entries of $\mathcal K$ directly from Schur inversion (Step~6) and sum
the birth contributions (Step~7). Together with the non-birth steps this
gives a stopped maximal moment for $\mathcal K$ (Step~8). We then
identify the entrance law of any candidate limit (Step~9) and show that
every such limit solves the martingale problem of the Gram and
complex-source generator (Step~10). Step~11 identifies the limiting
equation, the births and the brackets, and explains how the stops are
removed. Step~12 extends the rank-free coupling to the complex source.

\emph{Step 1: the exact one-step update.} Consider probes born before
the next layer. Adding the layer changes the Green matrix on the old
coordinates by Schur inversion,
$G_+|_{\rm old}-G=G\iota\Delta\iota^*G$. The boundary identity
$G^*V_{\rm ports}=GV_{\rm ports}\mathbb S^*$ then gives
\begin{equation}\label{eq:complex-source-update}
 \mathcal K_+-\mathcal K
   =V_{\rm pr}^*\iota\Delta\iota^*\mathbb S V_{\rm pr},
                       \qquad V_{{\rm pr},+}=N^*V_{\rm pr}.
\end{equation}
Put $\mathcal A=-(2i\eta)^{-1}\iota\Delta\iota^*\mathbb S$. The increment
is then $-2i\eta V_{\rm pr}^*\mathcal A V_{\rm pr}$, and
$NN^*-I=\mathcal A+\mathcal A^*$. Taking the imaginary part, the
increment of $\Im\mathcal K$ agrees exactly with the update of the Gram
matrix, as required by $\Im\mathcal K=-\eta M_{\rm pr}$.

The covariance of $\mathcal A$ with its adjoint, the opposite-sign
covariance, is $D_0P\otimes P'/(2W)$. Its covariance with itself, the
equal-sign covariance, contains the contracting letters of
Section~\ref{sec:full-source}, which disappear after summation. The
conditional mean of \eqref{eq:complex-source-update} is
\[
 O(\delta_0)V_{\rm pr}^*P\mathbb S V_{\rm pr}
 +O(\delta_0^2+W^{-2})M_{\rm pr}.
\]
The remaining steps control these terms along the whole path. The first term has $L^p$ size $C_p/W$ per step; its sum over the
order $W$ steps vanishes because $V_{\rm pr}^*P\mathbb SV_{\rm pr}$ is a
closed contracting coordinate, which decays after each restart by
\eqref{eq:closed-contracting-decay} of Step~5.

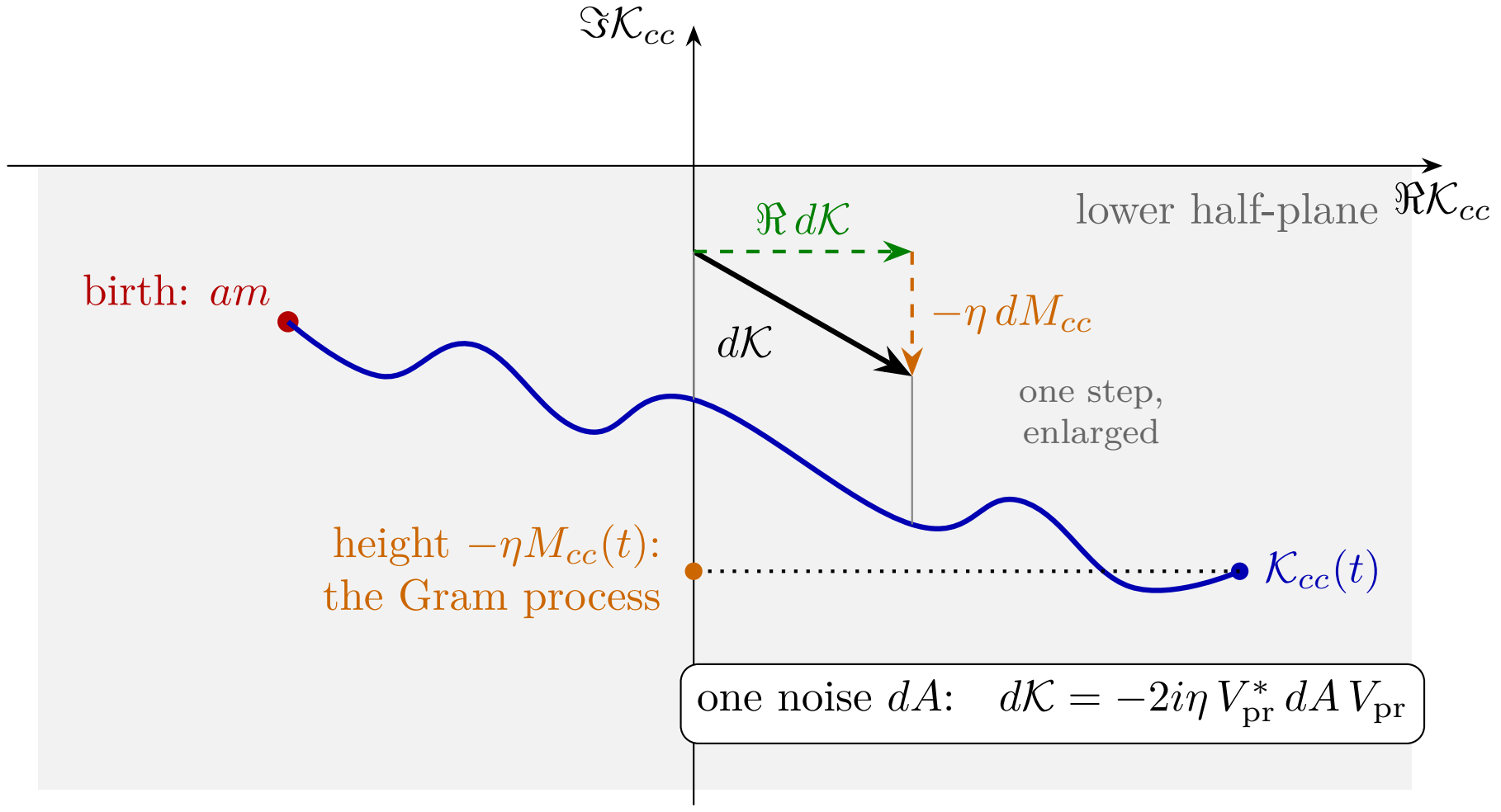
\begin{figure}[ht]
\centering
\begin{tikzpicture}[x=1.25cm,y=1.25cm,>=Stealth]
 \fill[black!5] (-4.2,0) rectangle (4.6,-4.0);
 \draw[->] (-4.4,0)--(4.8,0) node[below,font=\small] {$\Re\mathcal K_{cc}$};
 \draw[->] (0,-4.1)--(0,0.9) node[left,font=\small] {$\Im\mathcal K_{cc}$};
 \node[font=\small,black!60,anchor=east] at (4.5,-0.3) {lower half-plane};
 \fill[red!70!black] (-2.6,-1.0) circle (2.4pt);
 \node[font=\small,red!70!black,anchor=south east] at (-2.6,-1.0) {birth: $am$};
 \draw[very thick,blue!70!black] plot[smooth,tension=.7] coordinates
  {(-2.6,-1.0) (-2.0,-1.35) (-1.4,-1.15) (-0.7,-1.7) (0.0,-1.5)
   (1.4,-2.3) (2.1,-2.15) (2.8,-2.7) (3.5,-2.6)};
 \fill[blue!70!black] (3.5,-2.6) circle (2pt);
 \node[font=\small,blue!70!black,anchor=west] at (3.55,-2.6) {$\mathcal K_{cc}(t)$};
 \coordinate (P) at (0.0,-0.55);
 \coordinate (Q) at (1.4,-1.35);
 \draw[->,very thick] (P)--(Q) node[pos=.45,below left,font=\small] {$d\mathcal K$};
 \draw[->,thick,dashed,green!50!black] (P)--(1.4,-0.55)
   node[pos=.5,above,font=\small] {$\Re\,d\mathcal K$};
 \draw[->,thick,dashed,orange!80!black] (1.4,-0.55)--(Q)
   node[pos=.5,right,font=\small] {$-\eta\,dM_{cc}$};
 \draw[thin,black!50] (0.0,-1.5)--(P);
 \draw[thin,black!50] (1.4,-2.3)--(Q);
 \node[font=\scriptsize,black!60,align=center] at (2.55,-1.6) {one step,\\ enlarged};
 \draw[dotted,thick] (3.5,-2.6)--(0,-2.6);
 \fill[orange!80!black] (0,-2.6) circle (2pt);
 \node[font=\small,orange!80!black,anchor=east,align=right] at (-0.1,-2.6)
   {height $-\eta M_{cc}(t)$:\\ the Gram process};
 \node[draw,rounded corners,font=\small,align=center,fill=white] at (2.3,-3.45)
   {one noise $dA$:\quad $d\mathcal K=-2i\eta\,V_{\rm pr}^*\,dA\,V_{\rm pr}$};
\end{tikzpicture}
\caption{The mechanism behind Theorem~\ref{thm:complex-source-law}, drawn
for one diagonal entry of the probe matrix. A probe is born at the point
$am$, which lies in the lower half-plane because $\Im a<0$. Afterwards
the entry $\mathcal K_{cc}(t)$ moves in the lower half-plane (blue). Each
increment (solid arrow, shown enlarged above the path) is produced by the single noise $dA$, the same
noise that drives the source process. Its vertical component (dashed,
orange) is exactly $-\eta$ times the increment of the Gram entry
$M_{cc}$, so the height of the path is always $-\eta M_{cc}(t)$: the
imaginary part is the positive source process of
Section~\ref{sec:positive-source-law}, already constructed. The only new
information is the horizontal component (dashed, green), the real part.
It is not an independent noise; it is read off from the same $dA$, so
the real and imaginary parts are coupled and the whole complex matrix is
obtained by adding one linear equation to the Gram process
(Step~1 of the proof). The remaining steps justify this limit.}
\label{fig:complex-source}
\end{figure}

\emph{Step 2: the pre-step stop and the frame factors.}
We first fix the data and the stop. Fix $0<b<1/2$, $-2<E<2$, $0<h<1$,
finitely many masses $m_c\ge0$, and distinct birth times $s<\tau_c<T$,
with $s>0$. Write $d$ for the number of probes, $n_s=\lfloor Ws\rfloor$,
$n_T=\lfloor WT\rfloor$, and $m_*=\max_c m_c$, with $m_*=0$ if $d=0$.
The probability law $\mu_W$ used here is the product of three laws: the
law of the open Gaussian path grown layer by layer, the law of the
independent probe Gaussians, and the uniform law of the closing phase.
It is not the law of the periodic ring.

At time $n_s$ we close the transmission modes with $t_\nu\le h^2$, and we
impose a cap $K$ on the rank of the retained current. The level $h$ is taken off the set
\[
 \mathcal N_s=\{h\in(0,1):\rho\{t_k(\nu)=h^2\}>0
        \ \text{for some }k\ge0\},
\]
where, as in Step~9, $\rho$ is the law of the scalar transmission limit
$\mu(s)$ and $t_0(\nu)\ge t_1(\nu)\ge\cdots$ is the ordered transmission
sequence. Each $t_k$ has at most countably many atoms, so $\mathcal N_s$
is at most countable. At a level $h\in\mathcal N_s$ the physical cutoff
retains modes with transmission near $h^2$ at random, and the entrance
identification of Step~9 is not available at that level. For $h\notin\mathcal N_s$ the
continuity condition of Step~9 holds, and all limits $h\downarrow0$ run
through such levels. Let $V_j$ be the
full closed frame, consisting of the current columns and the probe
columns, and put $M_j=V_j^*V_j$ and $B_j=V_j^*P_nV_j$, where $P_n$ is the
projection onto the near port. For $R>0$ define
\[
 A_j=\{\text{physical current cap }K\}\cap
       \bigcap_{i=n_s}^j\{\Tr M_i<R\}.
\]
This is a \emph{pre-step event}: it is decided before step $j+1$, and the
increment through $j+1$ is kept in full, even if it makes the trace
overshoot $R$. All constants and width thresholds below are chosen for
these fixed data and for the inherited Gaussian constants, before $W$,
the time index and the coordinate; they are uniform over all admissible
orientations of the current frame. Vector norms are Euclidean, matrix
norms without a subscript are operator norms, and $L^p$ norms are taken
under $\mu_W$ unless a fresh-layer law or the open-path law is indicated.

On $A_j$, both the incoming and the outgoing near factor of every old
column have norm at most $\sqrt R$. To see this, let $\Pi$ be the closing
partial isometry, $z=e^{i\Theta}$ the closing phase and $\mathbb S_j$ the
full unitary scattering matrix, and set
\[
 H=I-z\Pi\mathbb S_j,\qquad C=H^{-1}(I-\Pi\Pi^*),\qquad Z=H^{-1}-I .
\]
Unitarity gives the energy identities
\[
 CC^*=I+Z+Z^*
 =[(I-\Pi^*\Pi)\mathbb S_j(I+Z)]^*
       [(I-\Pi^*\Pi)\mathbb S_j(I+Z)].
\]
They follow by multiplying the identity $H+H^*-HH^*=I-\Pi\Pi^*$, and its
counterpart with the factors in the opposite order, by $H^{-1}$ on both
sides. The near projection cannot increase a norm. Hence, for an old
column $x_k$ of the open frame, the two factors
\[
 u_k=P_nC^*x_k,\qquad v_k=P_n\mathbb S_j(I+Z)x_k
\]
obey the bound $\sqrt R$. The identities hold on the common set of full
measure on which all the inverses exist. Note what this gives and what it
does not: it controls the frame factors, but it is not a pointwise bound
on $\Re\mathcal K$. A stopped moment bound for $\mathcal K$ is proved in
Step~8.

\emph{Step 3: fixed fractional moments of the transmission.}
We need bounds on fractional traces of the transmission at the closing
time. Under the open path law, for every fixed $\alpha>0$ and every
finite $p>0$,
\[
 \sup_{W\ge W_0}\E_{\rm path}
       (\Tr T_{n_s}^{\alpha})^p\le H_{\alpha,p,b,E,s}<\infty.
\]
Here $T_i$ is the far transmission matrix; by unitarity of the two-port
scattering matrix, the near and far transmission spectra agree. The
higher moments are not obtained by Jensen's inequality from a first
moment; they come from a recursion, as follows.

Let $0<\alpha<1$ and let $q\ge1$ be an integer. Put $P=\Tr T^\alpha$,
$g=\Tr T$ and $\ell=\sum_{t_\nu\le1/2}t_\nu^\alpha$, the contribution of
the small transmission eigenvalues. Work on the regularity domain where
the scalar and jet errors satisfy $\delta_0,\delta_1\le d_{\rm reg}$ and
$g\le\gamma_{\rm reg}W$, for fixed $d_{\rm reg},\gamma_{\rm reg}>0$.
There the fractional tangent bound and its estimate in mean form give,
for the value $P'$ after one layer,
\[
 0\le P'\le P+\chi,\qquad
 \chi=C(\delta_0+W^{-2})P+L_\alpha-\E_{\rm fresh}L_\alpha,
 \qquad D=P+\chi-P'\ge0,\qquad
 \E_{\rm fresh}D\ge c\ell^2/W.
\]
The fluctuation term satisfies
$\|L_\alpha-\E_{\rm fresh}L_\alpha\|_k\le B_kP/\sqrt W$ for every fixed
integer $k\ge2$. We lift the one-step bound to the $q$th power by the
one-sided inequality
\[
 (P')^q\le(P+\chi)^q-2^{1-q}P^{q-1}D
               +2^{1-q}P^q\mathbf1_{\{\chi<-P/2\}}.
\]
We expand the polynomial, use the second-moment bound on the event in
the last term, use $\ell^2\ge P^2/2-4g^2$, and apply Young's inequality.
This keeps a negative multiple of $P^{q+1}/W$. The open jet and
conductance moments of order $q+1$ make irregular states have
probability $O(W^{-q-1})$, and since $P,P'\le W$ they cost only
$O(W^{-1})$ in this recursion. Consequently, with constants $c_q>0$ and
$D_q<\infty$ chosen before $W$, and $c_q\le1/(4(q+1))$,
\[
 y_{i+1}\le y_i-\frac{c_q}{W}y_i^{1+1/q}+\frac{D_q}{W},
 \qquad y_i=\E_{\rm path}P_i^q\le W^q.
\]
Iterating this recursion from $\lfloor Ws/2\rfloor$ to $n_s$, with no
assumption on the higher moments at the entrance, gives the admissible
constant
\[
 H_{\alpha,q,b,E,s}
 =\max\{1,(2D_q/c_q)^{q/(q+1)},(8q/(c_qs))^q\}.
\]
Real orders $p\le q$ follow from Lyapunov's inequality, the monotonicity
of $L^p$ norms in $p$. For $\alpha\ge1$ we use $\Tr T^\alpha\le\Tr T$.

In the phase estimate of Step~4 we use, for $p\ge1$, the choices
$q=\lceil p\rceil$ and $\alpha=1/(2q)$, so that $\alpha p\le1/2<1$. This
condition makes the phase integral in Step~4 finite; it is not a
threshold on the moments of the transmission sum. We assert no
uniformity as $\alpha\downarrow0$, $s\downarrow0$ or $p\to\infty$, or at
the endpoints of the energy or band ranges, and no exponential moment.

\emph{Step 4: the closing shift and the propagation of the scalar jet.}
Closing the small transmission modes changes the near input of the
path. We show that this change is small and that the scalar jet stays
small afterwards. Let $X_c$ and $X_o$ be the closed and open near inputs
at $n_s$, both at the physical boundary parameter $z_\partial$, and put
$R_c=(z_\partial-X_c)^{-1}$ and $R_o=(z_\partial-X_o)^{-1}$. The inputs
are not recomputed by closing at a different spectral parameter. Writing
$\mathbb R$ for the two-port resolvent at $n_s$, the near Schur
complement and the identity $\iota^*C(z)\iota=I$ give
\[
 R_o=\iota^*\mathbb R\iota=(z_\partial-X_o)^{-1},\qquad
 R_c-R_o=\iota^*\mathbb R(C(z)-I)\iota.
\]
From the canonical two-port blocks and the resolvent identity we obtain,
for $|w-z_\partial|\le\eta/2$,
\[
 \left|\frac1W\Tr\big[(w-X_c)^{-1}-(w-X_o)^{-1}\big]\right|
 \le\frac2{\eta W}\sum_{t_\nu\le h^2}
       \frac{t_\nu}{\sqrt{1-t_\nu}|1-z\sqrt{1-t_\nu}|}.
\]
Indeed, the nuclear norm of $R_c-R_o$ is at most $1/(2\eta)$ times that
sum, and in
$(I+(w-z_\partial)R_c)^{-1}(R_c-R_o)(I+(w-z_\partial)R_o)^{-1}$ the two
outer factors each have norm at most $2$.

We now bound the sum. Set $r_h=\sqrt{1-h^2}$ and
$K_{h,\alpha}=(2/r_h)^{1-\alpha}r_h^{-2\alpha}$. Each summand is bounded
both by $2/r_h$ and by $t_\nu/(r_h^2|1-z|)$; interpolating between the two
bounds, the sum is at most $K_{h,\alpha}|1-z|^{-\alpha}\Tr T_{n_s}^{\alpha}$.
The closing phase is independent of the path, and
$\int|1-e^{i\theta}|^{-u}\,d\theta/(2\pi)\le2^{-u}/(1-u)$ for $0<u<1$.
With the moments of Step~3 this gives, under the whole law and without
any stop,
\[
 \sup_{|w-z_\partial|\le\eta/2}
 \left\|\frac1W\Tr\big[(w-X_c)^{-1}-(w-X_o)^{-1}\big]\right\|_p
 \le\frac{C_{\rm sh}}W,
\]
\[
 C_{\rm sh}=\frac{2K_{h,\alpha}}\eta
       \left(\frac{2^{-\alpha p}}{1-\alpha p}\right)^{1/p}
       H_{\alpha,q,b,E,s}^{1/q},
 \qquad q=\lceil p\rceil,\quad\alpha=\frac1{2q}.
\]
The same pointwise majorant holds on the whole disc, so Cauchy's formula
bounds the shift of the first derivative with constant $2C_{\rm sh}/\eta$.

Now write $\delta_j=\delta_{0,j}^{\rm cl}$ for the scalar jet of the
closed path. Combining the shift with the open jet gives
$\|\delta_{n_s}\|_p\le S_p/W$. The passive scalar recursion has fixed
$d_0>0$, $0\le\rho<1$ and constants $C_m,Q_p$ such that
\[
 \delta_j\le d_0\ \Longrightarrow\
 \delta_{j+1}\le\rho\delta_j+C_m/W^2+\xi_{j+1},
 \qquad \|\xi_{j+1}\|_p\le Q_p/W,
 \qquad 0\le\delta_j\le1+\eta^{-1}.
\]
On the event that $\delta_{n_s}\le d_0/2$ and
$\xi_i\le(1-\rho)d_0/2$ at every later step, the jet never leaves the
contraction domain, and the geometric sum bounds its $L^p$ norm by
$[S_p+(C_m+Q_p)/(1-\rho)]/W$ for large $W$. Markov's inequality at order
$p+2$ and a union bound over at most $TW$ layers make the probability
of the complementary event at most $W^{-p}$, after a fixed width
threshold; on that event we use the deterministic bound on the jet. In
particular, at every $n_s\le j\le n_T$,
\[
 \|\delta_{0,j}^{\rm cl}\|_p\le J_p/W,
 \qquad J_p=S_p+(C_m+Q_p)/(1-\rho)+1+\eta^{-1}.
\]
Only high fixed moments and the local passive contraction are used here,
not an exponential moment of the summed jet.

\begin{figure}[ht]
\centering
\begin{tikzpicture}[>=Stealth]
 \begin{scope}[x=1cm,y=1cm]
  \draw[->,black!50] (-2.3,0)--(3.1,0) node[below,font=\small,black] {$\Re z$};
  \draw[->,black!50] (0,-2.2)--(0,2.3) node[left,font=\small,black] {$\Im z$};
  \draw[thick] (0,0) circle (1.6);
  \draw[line width=5pt,red!60,opacity=.7] (-20:1.6) arc (-20:20:1.6);
  \fill (1.6,0) circle (2pt) node[above right,font=\small] {$1$};
  \coordinate (Z) at (55:1.6);
  \fill[blue!70!black] (Z) circle (2.2pt);
  \node[font=\small,blue!70!black,anchor=south west] at (Z) {$z=e^{i\Theta}$};
  \draw[thick,blue!70!black] (Z)--(1.6,0)
    node[pos=.45,left=1pt,font=\small] {$|1-z|$};
  \node[font=\small,red!70!black,align=left,anchor=west] at (1.75,-0.9)
    {shift large\\ near $z=1$};
  \node[font=\small,align=center] at (0,-2.65)
    {uniform closing phase,\\ independent of the path};
 \end{scope}
 \begin{scope}[shift={(4.6,-2.0)},x=0.62cm,y=0.95cm]
  \fill[blue!15,domain=-3.1416:-0.035,samples=120]
    (-3.1416,0) -- plot (\x,{min(4.2,(2*abs(sin(\x r/2)))^(-0.5))}) -- (-0.035,0) -- cycle;
  \fill[blue!15,domain=0.035:3.1416,samples=120]
    (0.035,0) -- plot (\x,{min(4.2,(2*abs(sin(\x r/2)))^(-0.5))}) -- (3.1416,0) -- cycle;
  \draw[dashed,thick,red!70!black,domain=-3.1416:-0.24,samples=120]
    plot (\x,{(2*abs(sin(\x r/2)))^(-1)});
  \draw[dashed,thick,red!70!black,domain=0.24:3.1416,samples=120]
    plot (\x,{(2*abs(sin(\x r/2)))^(-1)});
  \draw[very thick,blue!70!black,domain=-3.1416:-0.035,samples=120]
    plot (\x,{min(4.2,(2*abs(sin(\x r/2)))^(-0.5))});
  \draw[very thick,blue!70!black,domain=0.035:3.1416,samples=120]
    plot (\x,{min(4.2,(2*abs(sin(\x r/2)))^(-0.5))});
  \draw[->] (-3.5,0)--(3.8,0) node[below,font=\small] {$\Theta$};
  \draw[->] (0,-0.1)--(0,4.6);
  \draw (-3.1416,-.08)--(-3.1416,.08) node[below=4pt,font=\small] {$-\pi$};
  \draw (3.1416,-.08)--(3.1416,.08) node[below=4pt,font=\small] {$\pi$};
  \node[font=\small,below] at (0,-.05) {$0$};
  \node[font=\small,blue!70!black,anchor=west] at (0.9,2.3)
    {$|1-e^{i\Theta}|^{-\alpha p}$, $\alpha p\le\frac12$};
  \node[font=\small,blue!70!black,anchor=west] at (0.9,1.85)
    {finite area};
  \node[font=\small,red!70!black,anchor=west] at (0.35,4.25)
    {$|1-e^{i\Theta}|^{-1}$: infinite area};
 \end{scope}
\end{tikzpicture}
\caption{The key estimate in Step~4 of the proof of
Theorem~\ref{thm:complex-source-law}, drawn schematically. Left: the
transmission modes with $t_\nu\le h^2$ are closed with a random phase
$z=e^{i\Theta}$, uniform on the unit circle and independent of the path.
A closed mode reflects back into the path with a factor of the form
$\bigl(1-z\sqrt{1-t_\nu}\bigr)^{-1}$, so the shift of the near input is
large only when $z$ is near $1$ (red arc). Each mode contributes at most
the smaller of a bounded term $2/r_h$ and a singular term
$t_\nu/(r_h^2|1-z|)$; interpolating between the two gives the bound
$K_{h,\alpha}\,t_\nu^{\alpha}|1-z|^{-\alpha}$. Right: averaging the
$p$th power over the phase leads to the integral of
$|1-e^{i\Theta}|^{-\alpha p}$. With $\alpha p\le1/2$ this function (solid)
blows up at $\Theta=0$ but has finite area (shaded). With the exponent one
(dashed), corresponding to using the full trace $\Tr T$ instead of a
fractional trace, the area would be infinite. This is why the proof
needs the fractional moments of $\Tr T^\alpha$ from Step~3, with
$\alpha=1/(2q)$. %
}
\label{fig:closing-shift}
\end{figure}
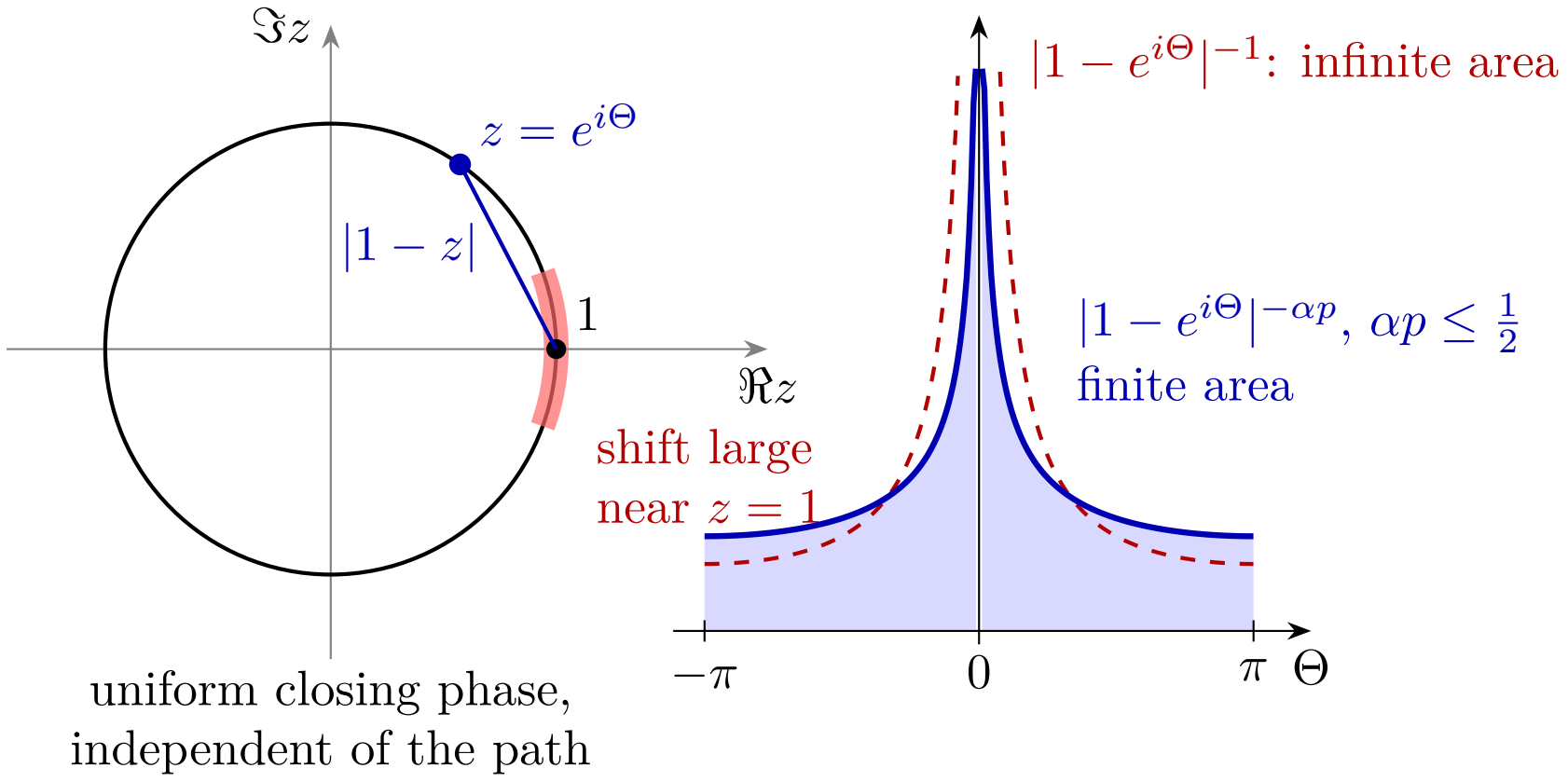

\emph{Step 5: the increments of the closed Gram coordinates, and the
closed contracting coordinates.} This step has two parts. First we bound
the one-step increments of the Gram coordinates of the closed frame
(part~(a)). Second, we show that certain contracting coordinates of the
closed frame decay after each birth (parts~(b)--(e)); this is what makes
the summed drift and covariance errors used in Step~10 tend to zero.

\emph{(a) The increments of the Gram coordinates.} We write
$\mathsf G_j=(M_j,B_j)$ for the Gram coordinates, to distinguish them from
the Green matrix. The Schur and Ward identities on the retained ports give
exactly, for any two old columns,
\[
 \Delta M_{kl}=c_\eta u_k^*\Delta v_l
                     +\overline{c_\eta u_l^*\Delta v_k},
 \qquad \Delta B_{kl}=u_k^*(U-I)u_l,
 \qquad c_\eta=-(2i\eta)^{-1},\quad U=LL^*.
\]
These identities follow from $N_cN_c^*-I=\mathcal A_c+\mathcal A_c^*$ and
$N_cP_nN_c^*=\iota U\iota^*$, together with the fact that closing and
transport commute. They hold for current columns as well as for probe
columns, and they do not use any smallness of the difference between
frames. The closed near reflection has defect trace at most $K$, because
the retained far port has dimension at most $K$. The fourth-moment
estimate for passive bilinear forms therefore gives, conditionally on the
past,
\[
 \|\Delta\mathsf G_e\|_{L^4({\rm fresh})}
 \le C_gR(1+K)(\delta_{0,j}^{\rm cl}+W^{-1/2})
 \quad\hbox{on }A_j,
\]
where $C_g$ and its width threshold depend only on the fixed $b,E$ and on
the inherited Gaussian constants. Here $u_k$ and $v_k$ have norm at most
$\sqrt R$, by Step~2. Independence of the fresh layer, Tonelli's theorem
and the jet bound of Step~4 give
\[
 \E_{\mu_W}[\mathbf1_{A_j}|\Delta\mathsf G_e|^4]
 \le\frac{C_\Gamma}{W^2},\qquad
 C_\Gamma=[C_gR(1+K)(1+J_4)]^4,
\]
with no logarithmic loss. Here an \emph{old} index is one that is not born
in the present step; the birth of a different probe in the same step is
allowed. The bound covers the current--current coordinates as well.
Consequently the stopped Lindeberg sum for the Gram coordinates, with
threshold $\epsilon$, is at most $TC_\Gamma/(\epsilon^2W)$, and an old Gram
entry at a birth step has second moment at most $\sqrt{C_\Gamma}/W$.

\emph{(b) The closed contracting coordinates, and why they must decay.}
The conditional means of the increments of $\mathcal K$ and of the Gram
coordinates contain, multiplied by a prefactor of $L^p$ size $C_p/W$, the
\emph{closed contracting coordinates}
\[
 \Upsilon^\sigma_{kl}(j)=(V_j^*\Sigma_\sigma(j)V_j)_{kl},\qquad
 \Sigma_\sigma(j)\in\{P_nU_j,\ U_j^*P_n,\ P_nU_jP_n,\ P_nU_j^*P_n\}.
\]
Here $U_j$ is the closed scatterer defined in part~(c); on the closed frame,
$P_nU_j$ acts as $P_n\mathbb S_jC_j$. For $\mathcal K$, this is the term
$O(\delta_0)V_{\rm pr}^*P\mathbb SV_{\rm pr}$ of Step~1.

These coordinates must be shown to decay. Indeed, the prefactor has size
$C_p/W$, because $\|\delta_{0,j}^{\rm cl}\|_p\le J_p/W$, but there are of
order $W$ steps. So a bound of order one on the coordinates would only give
a summed drift of order one. We show instead that they decay after each
\emph{restart}, that is, at $n_s$ and at each step that follows a birth
step. Let $r_j$ be the last restart at or before $j$, and put
$\vartheta=|\theta_*|<1$. For every $p\ge1$ there is a constant
$C<\infty$, depending on the fixed data, on $R$ and on $p$, such that for
$n_s\le j\le n_T$, all $\sigma$ and all columns $k,l$ of the closed frame,
\begin{equation}\label{eq:closed-contracting-decay}
 \|\mathbf1_{A_j}\Upsilon^\sigma_{kl}(j)\|_p
 \le C\bigl(\vartheta^{j-r_j}+W^{-1/2}\bigr).
\end{equation}
So each coordinate decays geometrically after a restart, down to the level
$W^{-1/2}$ of the fluctuations.

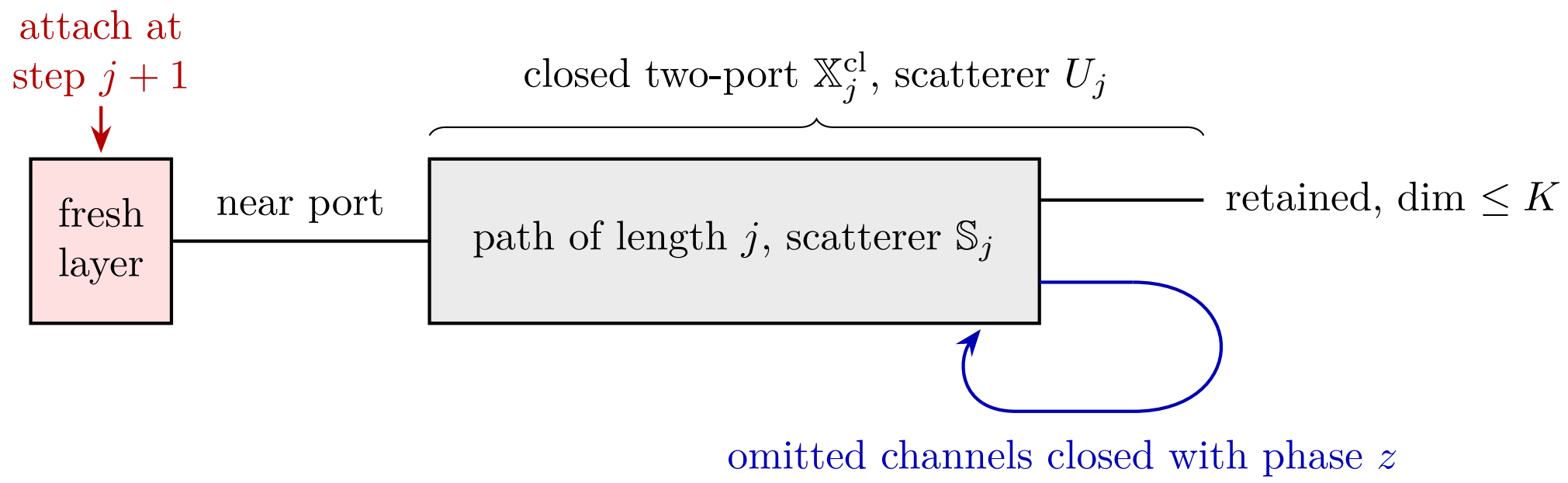
\begin{figure}[ht]
\centering
\begin{tikzpicture}[x=1cm,y=1cm,>=Stealth]
 \draw[thick,fill=red!12] (-3.4,0) rectangle (-2.2,1.4);
 \node[font=\small,align=center] at (-2.8,.7) {fresh\\ layer};
 \draw[thick,fill=black!8] (0,0) rectangle (5.2,1.4);
 \node[font=\small] at (2.6,.7) {path of length $j$, scatterer $\mathbb S_j$};
 \draw[thick] (-2.2,.7)--(0,.7);
 \node[font=\small,above] at (-1.1,.72) {near port};
 \draw[thick] (5.2,1.05)--(6.6,1.05);
 \node[font=\small,anchor=west] at (6.65,1.05) {retained, dim $\le K$};
 \draw[thick,blue!70!black] (5.2,.35)--(6.0,.35);
 \draw[thick,blue!70!black,->] (6.0,.35) .. controls (7.0,.35) and (7.0,-.75) .. (6.0,-.75)
    -- (5.0,-.75) .. controls (4.5,-.75) and (4.5,-.3) .. (4.7,-.05);
 \node[font=\small,blue!70!black] at (5.4,-1.15) {omitted channels closed with phase $z$};
 \draw[->,thick,red!70!black] (-2.8,1.85)--(-2.8,1.45);
 \node[font=\small,red!70!black,align=center] at (-2.8,2.3) {attach at\\ step $j+1$};
 \draw[decorate,decoration={brace,amplitude=4pt}] (0,1.6)--(6.6,1.6)
   node[midway,above=4pt,font=\small] {closed two-port $\mathbb X^{\rm cl}_j$, scatterer $U_j$};
\end{tikzpicture}
\caption{The closed path as a Hermitian two-port, used in part~(c) of Step~5. The
path of length $j$ (grey) has scatterer $\mathbb S_j$. At the far end, the
omitted channels are closed by feeding them back with the phase $z$
(blue), and only at most $K$ channels are retained. Together, the path
and this feedback form a single two-port with unitary scatterer $U_j$ and
a Hermitian Schur matrix $\mathbb X^{\rm cl}_j$. The next fresh layer
(red) is attached at the near end, while the feedback acts at the far
end. These two eliminations involve disjoint sets of equations, so they
commute: growing the closed two-port by one layer is the same as growing
the path and then closing it. This is why the one-layer analysis of
contracting words for a Hermitian old state applies to the closed path.}
\label{fig:closed-two-port}
\end{figure}

\emph{(c) The closed path as a path of Hermitian two-port states.} To
prove \eqref{eq:closed-contracting-decay}, we represent the closed path
as a path of Hermitian two-port states, grown by the same fresh layers as
the open path; see Figure~\ref{fig:closed-two-port}. Write $C_j,Z_j,H_j$
for the maps $C,Z,H$ of Step~2 at length $j$, let $A_-=\Pi\Pi^*$ and
$A_+=\Pi^*\Pi$ be the projections onto the omitted inputs and outputs,
and put
\[
 U_j=(I-A_+)\mathbb S_jC_j(I-A_-)+z\Pi^*A_-,\qquad
 \mathbb X_j^{\rm cl}=z_\partial I-2i\eta(I-U_j)^{-1}.
\]
We check four facts.

\emph{First, $U_j$ is unitary for $|z|=1$.} Let $A_-w=0$. Then
$u=C_jw=H_j^{-1}w$ satisfies $u=w+z\Pi\mathbb S_ju$, where the last term
lies in $\operatorname{ran}A_-$. Hence
$\|u\|^2=\|w\|^2+\|A_+\mathbb S_ju\|^2$ and, since $\mathbb S_j$ is
unitary, $\|(I-A_+)\mathbb S_ju\|=\|w\|$. As $C_jA_-=0$, the two terms of
$U_j$ are isometries, the first from $\ker A_-$ into $\ker A_+$ and the
second from $\operatorname{ran}A_-$ onto $\operatorname{ran}A_+$, and their
ranges are orthogonal. So $U_j$ is unitary.

\emph{Second, $I-U_j$ is invertible almost surely.} The entries of $U_j$
are polynomials in $z$ divided by $\det H_j$, so $\det(I-U_j)(\det H_j)^m$,
where $m$ is the dimension, is a polynomial in $z$. It is not identically
zero. Indeed, both it and $\det H_j$ are nonzero at $z=0$, where $H_j=I$ and
$U_j=(I-A_+)\mathbb S_j(I-A_-)$: if $U_jv=v$ at $z=0$, then
$\|v\|\le\|(I-A_-)v\|$ forces $A_-v=0$, then
$\|(I-A_+)\mathbb S_jv\|=\|v\|$ forces $A_+\mathbb S_jv=0$, hence
$\mathbb S_jv=v$, and $v=0$ because $I-\mathbb S_j=2i\eta\mathbb R$ is
invertible. So for each path only finitely many phases are excluded.
Since the phase is uniform and independent of the path, Tonelli's theorem
and a countable intersection over the lengths give a set of full measure
on which every $I-U_j$ is invertible.

\emph{Third, the closed state is Hermitian and reproduces the closed
data.} For a unitary $U$ with $I-U$ invertible,
$(I-U)^{-1}+((I-U)^{-1})^*=I$, so $\mathbb X_j^{\rm cl}$ is Hermitian, and
its scattering matrix $I-2i\eta(z_\partial-\mathbb X_j^{\rm cl})^{-1}$ is
$U_j$. Because $P_n\Pi^*=0$, $P_n(I-A_+)=P_n$ and $C_j(I-A_-)=C_j$, the
near rows of $U_j$ are those of the closed scatterer $\mathbb S_jC_j$.
The closed frame satisfies $A_-V_j=0$, because $C_jA_-=0$. Therefore the
closed near Schur matrix, the separators $\Sigma_\sigma(j)$ on the closed
frame, the closed current and the closed conductance are exactly those of
the Hermitian two-port $\mathbb X_j^{\rm cl}$.

\emph{Fourth, closing commutes with growth.} Attaching the fresh layer at
the near port and imposing the feedback of the omitted channels at the far
port eliminate disjoint groups of linear equations, so they can be done in
either order. On $\operatorname{ran}A_-$ there is no near wave, and the far
output stays $z\Pi^*$. Hence $\mathbb X_{j+1}^{\rm cl}$ is obtained from
$\mathbb X_j^{\rm cl}$ by the Schur growth with the fresh layer of the
step. On steps without a birth, the closed frame is transported by
$V_{j+1}=N^*V_j$, where $N$ is the transport of the old state
$\mathbb X_j^{\rm cl}$, as in \eqref{eq:complex-source-update}. The old
state is measurable before the fresh layer, and the fresh layer is
independent of the past, of the probes and of the phase.

\emph{(d) The decay estimate.} By part~(c), Steps~1 and~2 of the proof of
Lemma~\ref{lem:t26-contracting-words}, which treat one fresh layer at a
Hermitian old state, apply to the closed layer. The letter is the
rank-one matrix $V_je_le_k^*V_j^*$, which we split into two Hermitian
letters by polarization; their trace norms are at most
$\|V_je_k\|\|V_je_l\|\le\Tr M_j<R$ on $A_j$. The stop $A_j$ replaces the
stop $\sigma_K$ of that lemma, and the errors $\delta_0,\delta_1,g/W$ are
those of the closed near Schur matrix.

We bound these errors. By Step~4, $\|\delta_{0,j}^{\rm cl}\|_p\le J_p/W$.
For the first jet $\delta_{1,j}=\delta_{1,j}^{\rm cl}$, the closed near
input is passive, and the first-jet contraction behind
\eqref{eq:passive-jet-moments} gives fixed $0\le\varrho_1<1$ and constants
$C_1,C_1',Q_p'$ such that
\[
 \delta_{0,j},\delta_{1,j}\le d_0\ \Longrightarrow\
 \delta_{1,j+1}\le\varrho_1\delta_{1,j}+C_1\delta_{0,j}+C_1'/W^2
       +\xi^{(1)}_{j+1},
 \qquad \|\xi^{(1)}_{j+1}\|_p\le Q_p'/W,
\]
while the Herglotz representation gives the deterministic bound
$0\le\delta_1\le d_*+\eta^{-2}$. The open bound
\eqref{eq:passive-jet-moments} and the derivative shift $2C_{\rm sh}/\eta$
of Step~4 give $\|\delta_{1,n_s}\|_p\le S_p'/W$. We intersect the event of
Step~4 with the event that
\[
 \xi^{(1)}_i\le(1-\varrho_1)d_0/4,\qquad C_1\delta_{0,i}\le(1-\varrho_1)d_0/4
\]
at every later step. On this event the first jet stays in the
contraction domain, and the complement has probability at most $W^{-p}$,
by Markov's inequality at order $p+2$. As for $\delta_0$, this gives
\[
 \|\delta_{1,j}^{\rm cl}\|_p\le J_p'/W,\qquad
 J_p'=S_p'+\frac{C_1J_p+C_1'+Q_p'}{1-\varrho_1}+d_*+\eta^{-2}.
\]
The closed conductance is at most $K$ on the cap.

With these inputs, the one-step analysis gives, at every step without a
birth,
\[
 \Upsilon_{j+1}=\lambda_\sigma\Upsilon_j+e_j+\xi_{j+1},\qquad
 \E_j\xi_{j+1}=0,\qquad
 \|\mathbf1_{A_j}e_j\|_p\le C/W,\qquad
 \|\mathbf1_{A_j}\xi_{j+1}\|_p\le C/\sqrt W,
\]
with $\Upsilon_j=\Upsilon^\sigma_{kl}(j)$ and
$\lambda_\sigma\in\{\theta_*,\bar\theta_*\}$. So each coordinate is
multiplied by a factor of modulus $\vartheta<1$ at every step, up to a mean
error and a centered noise. At a restart, the deterministic bound
$\|\Sigma_\sigma(j)\|\le1$ gives
\[
 |\mathbf1_{A_j}\Upsilon^\sigma_{kl}(j)|\le\|V_je_k\|\|V_je_l\|<R .
\]
No step in $[r_j,j)$ is a birth step. Since $A_{j+1}\subseteq A_j$, the
iteration of Step~3 of the proof of Lemma~\ref{lem:t26-contracting-words},
started at $r_j$, gives
\[
 \|\mathbf1_{A_j}\Upsilon^\sigma_{kl}(j)\|_p
 \le R\vartheta^{j-r_j}+\frac{2C}{(1-\vartheta)\sqrt W},
\]
which is \eqref{eq:closed-contracting-decay}.

\emph{(e) The summed errors.} We now add up the contributions of the
contracting coordinates over the path. By H\"older's inequality, the mean
term of Step~1 contributes at most
\[
 C\sum_{n_s\le j<n_T}W^{-1}\bigl(\vartheta^{j-r_j}+W^{-1/2}\bigr)
 \le C\Bigl(\frac{d+1}{(1-\vartheta)W}+\frac T{\sqrt W}\Bigr)
\]
to the summed drift error of $\mathcal K$. Here the geometric terms are
summed separately after each of the at most $d+1$ restarts. The term
$O(\delta_0^2+W^{-2})M_{\rm pr}$ contributes $O(RT/W)$.

The drift errors of the Gram coordinates contain the four closed
contracting coordinates, with prefactors of $L^p$ size $C/W$, as in Step~4
of the proof of Lemma~\ref{lem:t26-contracting-words}, and they are summed
in the same way. So are the summed errors of the real covariance. Their
contracting part carries the prefactor $W^{-1}$ and products such as
$\Upsilon^\sigma_{il}(j)\Upsilon^{\sigma'}_{kj'}(j)$. Since
$|\mathbf1_{A_j}\Upsilon|<R$, H\"older's inequality and
\eqref{eq:closed-contracting-decay} bound their sum by the same quantity
$C((d+1)/((1-\vartheta)W)+T/\sqrt W)$. This proves that the summed drift
and covariance errors used in Step~10 tend to zero, for each fixed $R$.

\emph{Step 6: direct Schur estimates at a birth.}
At the birth of a probe of mass $m$, let $R_c$ now denote the passive
closed near resolvent just before this step. Schur inversion gives the
new entries of $\mathcal K$:
\[
 \begin{aligned}
 \Delta\mathcal K_{\rm old,old}&=u^*\Delta v,\\
 \mathcal K_{\rm old,new}&=u^*(I+\Delta R_c)\sqrt{m/W}\,g,\\
 \mathcal K_{\rm new,old}&=\sqrt{m/W}\,g^*(I+R_c\Delta)v,\\
 \mathcal K_{\rm new,new}&=(m/W)g^*G_b g,
       \qquad G_b=R_c+R_c\Delta R_c,
 \end{aligned}
\]
where $g$ is the Gaussian vector of the newborn probe. Before the birth,
an entry with a newborn index is zero. Passivity gives
$\|R_c\|\le\eta^{-1}$. Let $D_p$, for $p=2,4$, be the bound on the $L^p$
operator norm of $\Delta$ that is uniform over passive inputs; it is
supplied by Gaussian absorption at the fixed $b,E$.

We integrate the newborn Gaussian first, by Tonelli's theorem. This is
legitimate because its column vanishes up to the pre-step time, so the
event $A_j$, the old frames, $R_c$ and $\Delta$ do not depend on it. This
is independence in the product law $\mu_W$; it is not freshness with
respect to a filtration, since our filtration already contains all
probes.

Let $F_2$ be the constant in the passive bilinear bound
$\|u^*\Delta v\|_{L^2({\rm fresh})}\le F_2(\delta_0+W^{-1/2})\|u\|\|v\|$,
and put $V_2=\eta^{-1}+\eta^{-2}D_2$. The Gaussian identity
$\E_g|g^*Hg-\Tr H|^2=\Tr HH^*$ shows that the newborn quadratic form
fluctuates around its trace by at most $m_*V_2/\sqrt W$. The remaining
trace is estimated rather than assumed: summing the diagonal bilinear
forms of $R_c\Delta R_c$ gives
\[
 \|W^{-1}\Tr G_b-a\|_2
 \le(1+F_2\eta^{-2})J_2/W+F_2\eta^{-2}/\sqrt W.
\]
Consequently, for a newborn index $c$ and an old index $k$,
\[
 \begin{aligned}
 \|\mathbf1_{A_j}\mathcal K_{j+1,kc}\|_2,
 \ \|\mathbf1_{A_j}\mathcal K_{j+1,ck}\|_2
       &\le C_{on}/\sqrt W,\\
 \|\mathbf1_{A_j}(\mathcal K_{j+1,cc}-am_c)\|_2
       &\le C_{nn}/\sqrt W,\\
 C_{on}&=\sqrt{m_*R}(1+D_2/\eta),\\
 C_{nn}&=m_*[V_2+(1+F_2\eta^{-2})J_2+F_2\eta^{-2}].
 \end{aligned}
\]
An index that is still unborn after the step gives a zero cross entry.
We take the width threshold larger than those of $D_2$, $F_2$, $J_2$, the
Schur identities, and the separation of the floor birth steps. These are
direct estimates of the closed entries at fixed $h,K,s,R$, made before
any identification of the limit law. They use no moment of the bare
resolvent at a birth, and they do not assert an unstopped bound on that
resolvent.

\emph{Step 7: the finite birth sum.}
This step uses only the frame energy of Step~2, the exact birth
identities of Step~6, the operator moment $D_4$ and Gaussian integration
over the probes; it does not use the fractional or jet estimates. The
stopped $L^4$ norms of the increments at a birth are bounded by
$B_{oo}$ for an old/old entry, by $B_{on}/\sqrt W$ for an old/new or
new/old entry, and by $B_{nn}$ for a new/new entry, where
\[
 \begin{aligned}
 B_{oo}&=RD_4,\\
 B_{on}&=2^{1/4}\sqrt{m_*R}(1+D_4/\eta),\\
 B_{nn}&=24^{1/4}m_*(\eta^{-1}+\eta^{-2}D_4),
       \qquad B_* =\max(B_{oo},B_{on},B_{nn}).
 \end{aligned}
\]
Here we used $\E_g|w^*g|^4=2\|w\|^4$ and
$\E_g(\|g\|^2/W)^4=\prod_{r=0}^3(1+r/W)\le24$. The old/old identity
remains valid at the birth of a different probe. Let $\mathcal B_W$ be the
set of birth steps in $[n_s,n_T)$; it has at most $d$ elements. Once $W$
exceeds the fixed Schur and moment thresholds, $2/s$, and, when $d\ge2$,
$2/\min_{c\ne c'}|\tau_c-\tau_{c'}|$, distinct births occur at distinct
steps. Minkowski's inequality then gives
\[
 \E_{\mu_W}\left(\sum_{j\in\mathcal B_W}
    \mathbf1_{A_j}|\mathcal K_{j+1,kl}-\mathcal K_{j,kl}|\right)^4
 \le d^4B_*^4.
\]
In particular, this bound uses no moment of the closing correction and
no exponential moment of a summed jet.

\emph{Step 8: the stopped maximal moment of $\mathcal K$.}
For a fixed pair of probes put
$Z_j=\mathbf1_{A_j}(\mathcal K_{j+1,kl}-\mathcal K_{j,kl})$. On steps
without a birth, the passive bilinear bound and the jet estimate of
Step~4, valid under the whole law, give
\[
 \E|Z_j|^4\le C_Z/W^2,\qquad
 \|\E_jZ_j\|_4\le C_\mu/W,
 \qquad C_Z=[F_4R(J_4+1)]^4,\quad
 C_\mu=C_m'R(J_4+1).
\]
Here $F_4$ is the bilinear constant at order four and $C_m'$ is the
constant in the passive global mean bound
$\|\E_{\rm fresh}\Delta\|\le C_m'(\delta_0+W^{-1})$. The sum of the
absolute conditional means over non-birth steps has $L^4$ norm at most
$TC_\mu$. For the centered martingale differences $d_j=Z_j-\E_jZ_j$, the
maximal inequality for martingales at order four gives
\[
 \left\|\max_n\left|\sum_{j<n,\,j\notin\mathcal B_W}d_j\right|\right\|_4
 \le(17424T^2C_Z)^{1/4}.
\]
One way to see this is to use
$\E\max_n|\sum_{j<n}d_j|^4\le1089\E(\sum_j|d_j|^2)^2$, then
$\E(\sum_j|d_j|^2)^2\le N\sum_j\E|d_j|^4$ with $N\le TW$, and finally
$\E|d_j|^4\le16\E|Z_j|^4$.

All probes are unborn at $n_s$, so $\mathcal K_{n_s}=0$. The events $A_j$
decrease in $j$, and on $A_j$ all earlier stopped increments are
included. Adding the finite birth sum of Step~7, we obtain the stopped
maximal estimate, including the overshoot at the last step,
\[
 \left\|\max_{n_s\le j<n_T}\mathbf1_{A_j}
                      |\mathcal K_{j+1,kl}|\right\|_4\le P_4,
 \qquad P_4=(17424T^2C_Z)^{1/4}+TC_\mu+dB_*.
\]
If $O_4$ denotes the inherited bound on the raw open entry at each
length, then the stopped closing correction has $L^4$ norm at most
$O_4+P_4$. This is only an optional corollary; it is not an input to the
rest of the argument once the birth sum has been proved directly. The
order of dependence is: the birth sum, then the non-birth mean and
martingale estimates, then the stopped maximal moment, and only
optionally the correction moment, never the reverse. These are the
stopped controls of $\mathcal K$ used in the multiplier argument. A cap
and trace stop alone is not a pointwise bound on $\Re\mathcal K$.

\emph{Step 9: identification of the entrance law.}
This step identifies the entrance marginal of a candidate limit
conditioned on the cap; it does not assert that such a candidate exists
or that the family is tight. Write $\nu_W$ for the physical trace measure
of the transmission at time $s$, $\rho$ for the law of its scalar limit
$\mu(s)$, and $t_0(\nu)\ge t_1(\nu)\ge\cdots$ for the ordered
transmission sequence, so that $\nu=\sum_{k\ge0}t_k\delta_{t_k}$. We use
the cap $A=\{\nu:t_K(\nu)\le h^2\}$, assuming $p_A=\rho(A)>0$ and the
continuity condition $\rho\{t_k=h^2\}=0$ for every $k\ge0$. Let $A_W$ be
the physical cap, without the trace stop, and $p_W=\Prob(A_W)$. The
inherited scalar entrance and cap-convergence statements give, for every
bounded continuous real test $f$ of finite measures,
\[
 \nu_W\Rightarrow\mu(s),\qquad p_W\longrightarrow p_A,
 \qquad \E_{\rm path}[\mathbf1_{A_W}f(\nu_W)]
           \longrightarrow\int_A f(\nu)\,d\rho(\nu).
\]
Both the physical and the limiting measures are almost surely trace
measures of summable transmission sequences. Tests of either sign follow
by adding and subtracting $\|f\|_\infty$.

To pass from measures to the ordered sequence, we use an extractor that
needs no integrability assumption. For $0<a<c$, define the continuous
bounded weight on $[0,1]$ and its integral
\[
 q_{a,c}(x)=\frac{\min(1,\max(0,(x-a)/(c-a)))}{\max(x,a)},
 \qquad N_{a,c}(\nu)=\int q_{a,c}\,d\nu.
\]
For $\delta>0$, an integer $M\ge1/\delta$ and $k\ge0$, put
\[
 T_k^{\delta,M}(\nu)=\delta\sum_{j=0}^{M-1}
   \min(1,\max(0,N_{(j+1)\delta,(j+2)\delta}(\nu)-k)).
\]
This is continuous for weak convergence of finite measures. On a trace
measure, $N_{a,c}$ is the sum of the ramps
$\min(1,\max(0,(t_i-a)/(c-a)))$, so it counts, with a smooth transition,
the eigenvalues above $a$. Hence the $j$th clipped summand lies between
$\mathbf1_{\{(j+2)\delta\le t_k\}}$ and $\mathbf1_{\{(j+1)\delta\le t_k\}}$,
and therefore
\[
 0\le T_k^{\delta,M}(\nu)\le t_k(\nu)\le1,\qquad
 |T_k^{\delta,M}(\nu)-t_k(\nu)|\le2\delta.
\]
The error is uniform in the sequence, its total trace, $k$ and $W$. For a
finite set $J\subset\mathbb N$ and a bounded continuous
$f:\mathbb R^J\to\mathbb R$, uniform continuity on $[0,1]^J$ allows us to
replace $f((t_k)_{k\in J})$ by the continuous test
$f((T_k^{\delta,M})_{k\in J})$, with an error that vanishes uniformly. We
apply the cap-weighted convergence to that test and then let $\delta$
decrease to zero. Since eventually $p_W\ge p_A/2$, division gives
\[
 \frac{\E_{\rm path}[\mathbf1_{A_W}f((t_k^W)_{k\in J})]}{p_W}
 \longrightarrow
 \frac1{p_A}\int_A f((t_k(\nu))_{k\in J})\,d\rho(\nu).
\]
Now let $\mathsf Q$ be a candidate limit with the cap-conditioned
convergence of finite-dimensional distributions stipulated in the
cap-limit premise. We apply this convergence, along its subsequence, to
the entrance labels at model time zero, which is physical time $s$,
before every birth. The real test $f$ is composed with the real parts of
the coordinates, which gives a bounded continuous test on the complex
observation space. The probe and phase factors integrate out of this
observable, which depends only on the entrance. The same limit is then
$\E_{\mathsf Q}f((\widehat t_k)_{k\in J})$, so
\[
 \mathcal L_{\mathsf Q}((\widehat t_k)_{k\in J})
 =\mathcal L_{\rho(\,\cdot\mid A)}((t_k(\nu))_{k\in J}).
\]
Finite-dimensional tests, followed by uniqueness of a law on the
countable product, identify the law of the entrance sequence. The cap is
conditioned on, not discarded, and we do not assert that this entrance is
independent of the later process.

\emph{Step 10: the limit of the Gram and complex-source generator.} In
this step we show that every candidate limit solves the martingale
problem of the Gram and complex-source generator. Concretely, we justify
the limit passage for the coordinate tests of first and second order,
including their mixed and conjugate brackets. The argument has seven
parts: the setting and the test functions (a), the generator (b), the
finite-width process with its discrete compensation of the births (c), a
bound that does not depend on the auxiliary stop, together with the
exact feasibility of the finite-width Gram data (d), the Taylor expansion
at finite width (e), the passage from the discrete compensation to the
model-time compensation (f), and finally the removal of the localization
and the passage to the limit (g).

\emph{(a) The setting and the test functions.} Fix the cap-conditioned
candidate $\mathsf Q$ of Step~9 and a finite model horizon $H$. Model time
$u$ corresponds to physical time $s+u$. We apply the finite-width
estimates, which hold globally in time, at $T_H=\max(T,s+H+1)$. Write
$\widehat\mu_W=\mu_W(\,\cdot\mid A_W)$ for the law conditioned on the
physical cap. Since $p_W\to p_A>0$, every cap-weighted error estimate costs
at most a factor $2/p_A$ after conditioning, for all sufficiently large
$W$. The cap is measurable at the starting time, so conditioning does not
change the identities of the later fresh layers on that event.

Let $J_G$ and $J_K$ be the deterministic jumps of the Gram and
complex-source coordinates at the births, in model time; precisely,
\[
 J_G(u)=\sum_{c:\,\tau_c-s\le u}\Gamma_c,\qquad
 J_K(u)=\sum_{c:\,\tau_c-s\le u}\kappa_c,
\]
where $\Gamma_c$ is the deterministic jump of the Gram coordinates at the
birth of probe $c$, whose $M$-component is $d_*m_c$ at the new diagonal
entry, and $\kappa_c$ is the jump $am_c$ at the new diagonal entry of
$\mathcal K$. We subtract these jumps before forming the tests:
\[
 \begin{aligned}
 Y&=\mathsf G-J_G,\qquad Y^K=\mathcal K-J_K,
       \qquad \vartheta(Y)=\Re\Tr Y^M,\\
 A_e^\chi(u)&=\int_0^u b_e(\mathsf G(v))\chi(\vartheta(Y(v)))\,dv,\\
 x_e&=Y_e-Y_e(0)-A_e^\chi,\qquad
 x_a^K=Y_a^K-Y_a^K(0).
 \end{aligned}
\]
Here $Y^M$ is the $M$-component of $Y$, and $b_e$ is the Gram drift of
\eqref{eq:gram-drift}. So $A^\chi$ is the accumulated drift cut off by
$\chi$, and $x_e$, $x_a^K$ are the compensated coordinates. For real
compactly supported $C^3$ cutoffs $\chi$ and $\xi$, put
$q(Y^K)=\sum_a|Y_a^K|^2$. The two Gram test functions are $x_e$ and
$x_ex_f$, multiplied by $\chi(\vartheta(Y(0)))\chi(\vartheta(Y))$. The four
complex-source test functions are $x_a^K$, $x_a^Kx_b^K$, $x_a^Kx_e$ and
$x_a^K\overline{x_b^K}$, multiplied by
\[
 \chi(\vartheta(Y(0)))\xi(q(Y^K(0)))
       \chi(\vartheta(Y))\xi(q(Y^K)).
\]
We denote any of these six composite test functions by $F$, and the
augmented state by
\[
 Z=(Y,A^\chi,Y(0),Y^K,Y^K(0)).
\]

\emph{(b) The generator.} The generator $\Lambda(J_G,Z)$ applied to $F$ is
the ordinary real drift plus the second-order covariance operator,
computed from \eqref{eq:gram-drift}, \eqref{eq:gram-covariance} and the
complex and mixed brackets, which all come from the same noise $dA$ in
\eqref{eq:complex-source-sde}. In particular, $x_e$ has drift
$b_e(\mathsf G)(1-\chi)$ and the Gram covariance, while $x_a^K$ has zero
drift and zero equal-sign bracket. The product rule includes the mixed
term with $\chi'\xi'$ and both terms with one half of a second
derivative. Conjugating a Gram coordinate swaps its matrix indices, so
these bilinear brackets, together with the conjugate complex-source
brackets, determine the full real covariance; no bracket is omitted.
Define
\[
 M_F(u)=F(Z(u))-\int_0^u\Lambda(J_G(v),Z(v))\,dv.
\]
We must show that $M_F$ is a martingale under every candidate limit.

\emph{(c) The finite-width process and the discrete compensation.} The
process at finite width is defined before any auxiliary stop. Put
$n_u=\lfloor W(s+u)\rfloor$. At finite width a probe is born at a grid
step, not at its model time, so we compensate the births on the grid.
With $j_c=\lfloor W\tau_c\rfloor$, probe $c$ belongs to the closed frame at
index $j$ exactly when $j_c<j$. The aligned discrete compensators of the
births are therefore
\[
 J_{G,W}(j)=\sum_{c:\,j_c<j}\Gamma_c,\qquad
 J_{K,W}(j)=\sum_{c:\,j_c<j}\kappa_c .
\]
Since $\tau\mapsto\lfloor W\tau\rfloor$ is nondecreasing, every value of
$J_{G,W}$ and $J_{K,W}$ is a value of $J_G$ and $J_K$, respectively. With
$Y_W=\mathsf G_W-J_{G,W}$ and $Y_W^K=\mathcal K_W-J_{K,W}$, set
\[
 A_{e,W}^\chi(n)=\frac1W\sum_{n_s\le j<n}
       b_e(\mathsf G_{W,j})\chi(\vartheta(Y_{W,j})),
 \qquad
 Z_W(n)=(Y_W(n),A_W^\chi(n),Y_W(n_s),Y_W^K(n),Y_W^K(n_s)).
\]
The discretely compensated generator process is
\[
 M_{F,W}(u)=F(Z_W(n_u))-
       \frac1W\sum_{n_s\le j<n_u}\Lambda(J_{G,W}(j),Z_W(j)).
\]
Here $A_W^\chi$ denotes the drift vector, not the cap event $A_W$.

\emph{(d) A bound that does not depend on the auxiliary stop.} Choose
$R_\chi,R_\xi\ge1$ so that the supports of the cutoffs and of their first
three derivatives are contained in $[-R_\chi,R_\chi]$ and
$[-R_\xi,R_\xi]$. With $d_{\rm fr}=2K+d$ the size of the full frame, put
\[
 C_J=1+\max_{u,e}|(J_G(u))_e|,
 \qquad D_Y=2(R_\chi+d_{\rm fr}C_J)+C_J,
\]
where the maximum runs over the finitely many deterministic jump values.
Feasibility of the Gram data gives $M\succeq0$, $0\preceq B\preceq M$ and
$|\mathsf G_e|\le2\Re\Tr M$. So the trace cutoff bounds $|Y_e|$ by $D_Y$.
The cut drift is bounded on a compact set:
\[
 C_A=1+\max_{J\in\operatorname{range}J_G,\ |y_e|\le D_Y,\ e}
           |b_e(y+J)\chi(\vartheta(y))|,
 \qquad |A_e^\chi(u)|\le HC_A.
\]
Wherever a test function or its generator is nonzero, both $Y$ and $Y(0)$
lie in this box, so $|x_e|\le2D_Y+HC_A$. For the complex-source tests, the
current and initial $Y^K$ coordinates are at most $\sqrt{R_\xi}$, hence
$|x_a^K|\le2\sqrt{R_\xi}$. On the resulting compact box $\mathcal B_H$ put
\[
 C_F=\max_{Z\in\mathcal B_H}|F(Z)|,\qquad
 C_\Lambda=\max_{J\in\operatorname{range}J_G,\ Z\in\mathcal B_H}
             |\Lambda(J,Z)|,\qquad C_H=C_F+HC_\Lambda.
\]
Consequently $|M_F(u)|\le C_H$ on feasible paths for $u\le H$, also at
stopped times. Enlarging the box by using $H+1$ in place of $H$ gives the
same bound for $M_{F,W}$, uniformly for large $W$; this uses that the
discrete compensators take only values of $J_G$. These constants depend
on the fixed frame data, the cutoffs and the horizon, but not on the
auxiliary trace level $R$ introduced below. Feasibility of the candidate
limit follows from feasibility at finite width, at a countable dense set
of times off the births, and right-continuity.

\emph{Exact feasibility at finite width.} The feasibility at finite width
used here holds exactly, not only approximately. For every $j\ge n_s$, on
the set of full measure of Step~2 and for every phase at which $I-U_j$ is
invertible, $(M_j,B_j)$ is a feasible Gram state for $\mathsf D$ in the
sense of Proposition~\ref{prop:finite-gram-law}, with the boundary space
as ambient space and $P=P_n$. To see this, note that $P_nZ_j=0$, so
$P_nC_j=P_n$, and the closed scatterer $U_j$ of Step~5 satisfies
$P_nU_j=P_n\mathbb S_jC_j$. Hence the closed current is
\[
 J_j^{\rm cl}=U_j^*P_nU_j-P_n=C_j^*(\mathbb S_j^*P_n\mathbb S_j-P_n)C_j .
\]
At $n_s$, closing turns each omitted pair into a near reflector without
current and leaves the retained pairs unchanged, so
$J_{n_s}^{\rm cl}=V_{n_s}\mathsf DV_{n_s}^*$. This identity is preserved
from step to step: by the conservation law \eqref{eq:current-conservation}
for the closed layer, by the transport $V_{j+1}=N^*V_j$, and by the zero
signature of newborn probes. Therefore
\[
 P_n+V_j\mathsf DV_j^*=U_j^*P_nU_j,
\]
and the right side is an orthogonal projection of the same rank as $P_n$,
because $U_j$ is unitary. This is exactly the two-projection constraint
that defines a feasible Gram state.

\emph{Smooth extensions.} For the Taylor expansion we multiply the
composite $F$ by fixed smooth cutoffs that equal one on a neighbourhood of
these boxes, including the drift and initial coordinates. This gives an
extension of $F$ that is $C^3$ on the whole space, with bounded derivative
norms $L_i=\sup\|D^iF\|$, $i=1,2,3$, and that agrees with $F$ and its first
two derivatives on the relevant feasible states. For convergence of
expectations we likewise use bounded continuous extensions of all three
functions involved: the cut drift integrand, the endpoint test $F$, and
the generator integrand $\Lambda$. They agree with the original functions
on feasible states, so no target test is changed.

\emph{(e) The Taylor remainder and the births, at finite width.} The
estimates of this part hold at every level $R>0$, for the discretely
compensated process $M_{F,W}$. The passage to the model-time jumps is made
afterwards, in part~(f), at the levels $R\ge R_0$ of
\eqref{eq:cs-identgram-threshold}.

Let $j_R$ be the first index at or after $n_s$ with $\Re\Tr M_{W,j}\ge R$,
and let $\sigma_W^R$ be the corresponding model stopping time. The
indicator $\mathbf1_{\{j<j_R\}}$ includes the step that reaches or
overshoots $R$. Let $X_j$ collect all real Gram and complex-source
coordinates, and, summing over the steps without a birth from $n_s$ to
$n_H-1$, put
\[
 \begin{aligned}
 Q_W(R)&=\sum_{j\ {\rm nonbirth}}\E_{\widehat\mu_W}
         [\mathbf1_{\{j<j_R\}}\|\Delta X_j\|_2^2],\\
 L_W(R)&=\sum_{j\ {\rm nonbirth}}\E_{\widehat\mu_W}
         [\mathbf1_{\{j<j_R\}}\|\Delta X_j\|_2^4].
 \end{aligned}
\]
The summed errors in the conditional drift and the real covariance,
together with the $L^1$ errors at the births, have a total
$\epsilon_W(R)\to0$ for each fixed $R$. Continuity of the covariance on
feasible Gram states with trace at most $R$ gives
$Q_W(R)\le(H+1)C_R+\epsilon_W(R)$, with $C_R$ independent of $W$ and of the
values of the complex-source coordinates. The fourth-moment Lindeberg
sums of the coordinates imply $L_W(R)\to0$, since for the fixed number $q$
of real coordinates $\|\Delta X\|_2^4\le q\sum_i|\Delta X_i|^4$.

With $\mathcal E_G$ the set of complex Gram coordinates, take
$B_A=\sqrt{2|\mathcal E_G|}\,C_A$; it bounds the Euclidean norm of the
drift integrand, whose increments are therefore at most $B_A/W$. We may
then take, for the augmented state,
\[
 Q_W^{\rm aug}=2Q_W+2(H+1)B_A^2/W,\qquad
 L_W^{\rm aug}=8L_W+8(H+1)B_A^4/W^3.
\]
The Cauchy--Schwarz inequality, for counting measure in time times
$\widehat\mu_W$, bounds the summed third-order Taylor remainder by
\[
 \frac{L_3}{6}\sum_j\E_{\widehat\mu_W}
          [\mathbf1_{\{j<j_R\}}\|\Delta Z_j\|_2^3]
 \le\frac{L_3}{6}\sqrt{Q_W^{\rm aug}(R)L_W^{\rm aug}(R)}\longrightarrow0,
\]
where the sum is over steps without a birth. The other summed errors are
at most
\[
 (L_1+L_2)\epsilon_W(R)
 +L_2B_A\sqrt{(H+1)Q_W(R)/W}
 +\tfrac12L_2(H+1)B_A^2/W+o(1).
\]
The deterministic term $o(1)$ is the rounding of the endpoints, which
vanishes by uniform continuity on the fixed boxes and the bounded drift.

At a birth, we first subtract the deterministic jump and then apply the
global Lipschitz bound $L_1$ to the remaining $L^1$ birth error. We do not
Taylor-expand the prescribed jump, which is of order one. The extra drift
insertion at these at most $d$ steps costs at most $L_1dB_A/W$. No newborn
probe is declared fresh with respect to the filtration, which already
contains the phase and all probes. Telescoping over the steps, and testing
against any bounded past cylinder function $\phi_W$, we obtain for each
fixed $R$
\[
 \E_{\widehat\mu_W}\big[\phi_W
   (M_{F,W}(v\wedge\sigma_W^R)-M_{F,W}(u\wedge\sigma_W^R))\big]
 \longrightarrow0.
\]

\begin{figure}[ht]
\centering
\begin{tikzpicture}[x=1.1cm,y=1cm,>=Stealth]
 \draw[->] (-.3,0)--(10.4,0) node[below,font=\small] {model time};
 \foreach \k in {0,1.5,...,9} \draw (\k,-.08)--(\k,.08);
 \node[font=\small,below] at (3,-.1) {$\frac{j_c}W-s$};
 \node[font=\small,below] at (4.5,-.1) {$\frac{j_c+1}W-s$};
 \draw[very thick,red!70!black] (3.8,0)--(4.5,0);
 \fill[red!70!black] (3.8,0) circle (2pt);
 \node[font=\small,red!70!black,anchor=east] at (3.75,.35) {$\tau_c-s$};
 \draw[decorate,decoration={brace,amplitude=3pt}] (3.8,.25)--(4.5,.25)
   node[midway,above=3pt,font=\small,red!70!black] {$\mathcal D_W$};
 \draw[thick,blue!70!black] (-.3,1.4)--(3.8,1.4);
 \draw[thick,blue!70!black] (3.8,2.2)--(10.2,2.2);
 \fill[blue!70!black] (3.8,2.2) circle (2pt);
 \draw[blue!70!black] (3.8,1.4) circle (2pt);
 \node[font=\small,blue!70!black,anchor=west] at (10.3,2.2) {$J_G(r)$};
 \draw[thick,green!50!black,dashed] (-.3,1.25)--(4.5,1.25);
 \draw[thick,green!50!black,dashed] (4.5,2.05)--(10.2,2.05);
 \fill[green!50!black] (4.5,2.05) circle (2pt);
 \node[font=\small,green!50!black,anchor=west] at (10.3,1.25) {$J_{G,W}(n_r)$};
 \draw[thick] (6,-.75)--(6,-.95) node[below,font=\small] {$u$};
 \draw[thick] (9,-.75)--(9,-.95) node[below,font=\small] {$v$};
 \draw[decorate,decoration={brace,amplitude=4pt,mirror}] (6,-1.45)--(9,-1.45)
   node[midway,below=4pt,font=\small] {cells $[v_j^{u,v},v_{j+1}^{u,v}]$, end cells partial};
\end{tikzpicture}
\caption{The two ways of compensating a birth, used in part~(f) of Step~10. The ticks
are the grid times $j/W-s$. A probe is born at the model time $\tau_c-s$
(red dot), but at finite width it enters the closed frame only from the
next grid step $j_c+1$ on. The model-time compensator $J_G(r)$ (blue,
solid) jumps at $\tau_c-s$, while the discrete compensator $J_{G,W}(n_r)$
(green, dashed; drawn slightly lower for visibility) jumps at
$(j_c+1)/W-s$. They differ only on the short delayed-birth interval
$\mathcal D_W$ (red), of length at most $1/W$ for each birth. The
increments over $[u,v]$ are split into the clipped cells shown below the
axis; the two end cells may be partial. The fixed times $u$ and $v$ lie
outside $\mathcal D_W$ for large $W$, so the two compensations agree at
them.}
\label{fig:cs-compensation}
\end{figure}
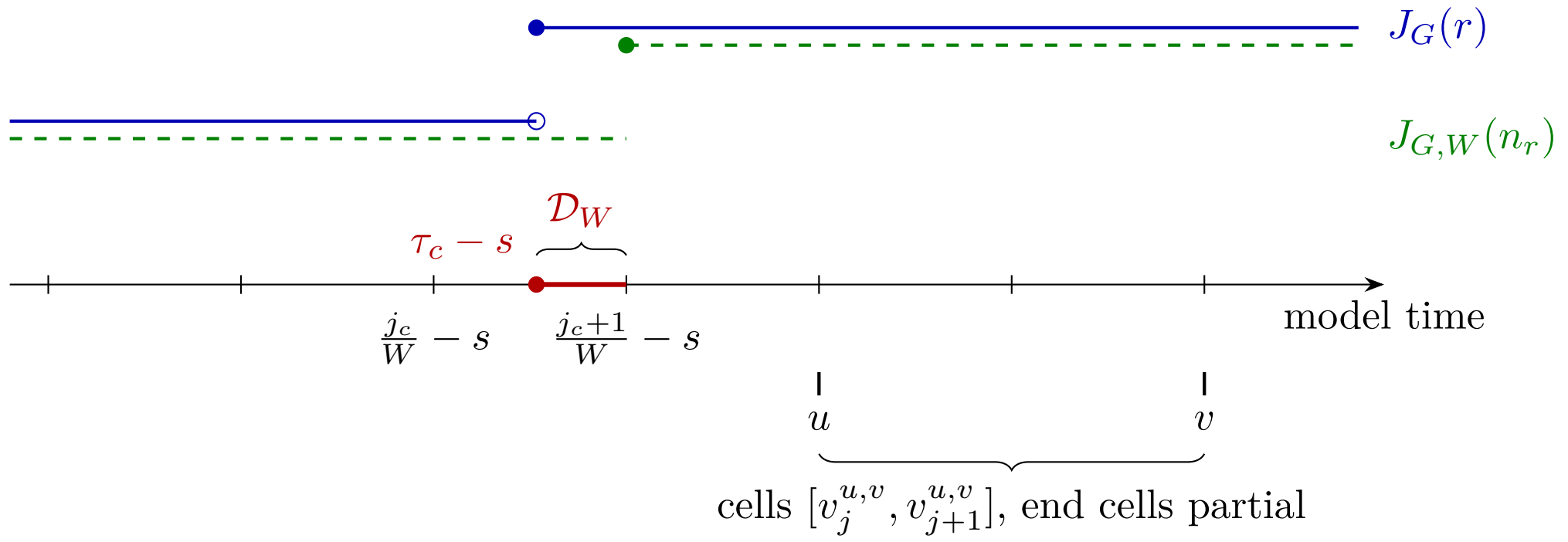

\emph{(f) Passage to the model-time jumps.} The limit is taken for the
process whose births are compensated by $J_G$ and $J_K$ at model time, not
by the discrete compensators. We now compare the two; see
Figure~\ref{fig:cs-compensation}. Put
$\mathsf G_W(r)=\mathsf G_{W,n_r}$, $\mathcal K_W(r)=\mathcal K_{W,n_r}$,
$\widetilde Y_W=\mathsf G_W-J_G$ and $\widetilde Y_W^K=\mathcal K_W-J_K$, and
\[
 \widetilde A_{e,W}^\chi(u)=\int_0^u b_e(\mathsf G_W(r))
        \chi(\vartheta(\widetilde Y_W(r)))\,dr .
\]
Let $\widetilde Z_W$ be the augmented state formed from these quantities
and the initial values, and
\[
 \widetilde M_{F,W}(u)=F(\widetilde Z_W(u))
      -\int_0^u\Lambda(J_G(r),\widetilde Z_W(r))\,dr .
\]
Fix times $0\le u\le v\le H$ off the births, and a test $\phi_W$ observed
at times at most $u$ off the births. With the clipped vertices
$v_j^{u,v}=\max(u,\min(v,j/W-s))$, the cells $[v_j^{u,v},v_{j+1}^{u,v}]$,
$n_u\le j\le n_v$, partition $[u,v]$; the two end cells may be partial. We
must show that, at suitable levels $R$, the stopped cell sum
\begin{equation}\label{eq:cs-identgram-cells}
 S_W^R=\E_{\widehat\mu_W}\Big[\phi_W\sum_{n_u\le j\le n_v}
     \mathbf1_{\{j<j_R\}}\big(\widetilde M_{F,W}(v_{j+1}^{u,v})
          -\widetilde M_{F,W}(v_j^{u,v})\big)\Big]
\end{equation}
tends to zero.

\emph{Where the two compensations differ.} They disagree only on the
delayed-birth set
\[
 \mathcal D_W=\bigcup_c\big[\tau_c-s,(j_c+1)/W-s\big),\qquad
 |\mathcal D_W|\le d/W,
\]
that is, $J_{G,W}(n_r)=J_G(r)$ and $J_{K,W}(n_r)=J_K(r)$ for
$r\notin\mathcal D_W$. The finitely many fixed times $u$, $v$ and the
observation times of $\phi_W$ have positive distance to the births, so
they lie outside $\mathcal D_W$ for all large $W$.

\emph{The drifts and the generator integrals.} Let $A_W^\sharp$ be the
integral in $r$ of the discretely compensated cut drift
$b_e(\mathsf G_W(r))\chi(\vartheta(\mathsf G_W(r)-J_{G,W}(n_r)))$. The two
drift integrands have Euclidean norm at most $B_A$ and differ only on
$\mathcal D_W$, so
\[
 \sup_{u'\le H}\|\widetilde A_W^\chi(u')-A_W^\sharp(u')\|_2\le2B_Ad/W,
\]
and the left-grid sum $A_W^\chi$ differs from $A_W^\sharp$ by at most
$2B_A/W$, coming from its two partial cells. Put $\delta_W=2B_A(d+1)/W$.
At each vertex, the values of $F$ with the same compensation but the two
different drifts differ by at most $L_1\delta_W$. Summed against the
decreasing indicators $\mathbf1_{\{j<j_R\}}$, these differences telescope
to the two ends of the interval where the process is alive, and they cost
at most $2\|\phi\|_\infty L_1\delta_W$.

Let $\omega_\Lambda$ be the modulus of continuity of the extended generator
integrand on the fixed compact box, at a fixed jump value. Outside
$\mathcal D_W$ the jump parameters agree, and only the drift and the
position within a cell change; on $\mathcal D_W$ and on the two end cells
both integrands are bounded by $C_\Lambda$. The generator integrals
therefore differ by at most
\[
 \|\phi\|_\infty\Big[H\omega_\Lambda(\delta_W+2B_A/W)
      +\frac{2C_\Lambda(d+2)}{W}\Big].
\]
These bounds tend to zero and do not depend on $R$.

\emph{The endpoint values of $F$.} It remains to compare the values of $F$
at the vertices under the two compensations. Let $E_j$ be the value of $F$
at the vertex $v_j^{u,v}$ computed from $\mathsf G_{W,j}$ and
$\mathcal K_{W,j}$ with the model-time jumps $J_G(v_j^{u,v})$,
$J_K(v_j^{u,v})$, minus its value computed with the discrete jumps
$J_{G,W}(j)$, $J_{K,W}(j)$, with the same initial state and the same drift.
For large $W$, $E_j=0$ at the two clipped endpoints, since they lie outside
$\mathcal D_W$. At an interior vertex $v_j^{u,v}=j/W-s$, the definitions
give
\[
 J_G(j/W-s)-J_{G,W}(j)=\sum_{c:\,W\tau_c=j}\Gamma_c,\qquad
 J_K(j/W-s)-J_{K,W}(j)=\sum_{c:\,W\tau_c=j}\kappa_c,
\]
so $E_j=0$ unless a birth lies exactly on that vertex. Summation by parts,
with $\mathbf1_{\{j-1<j_R\}}-\mathbf1_{\{j<j_R\}}=\mathbf1_{\{j_R=j\}}$,
gives the exact identity
\begin{equation}\label{eq:cs-identgram-residual}
 \sum_{n_u\le j\le n_v}\mathbf1_{\{j<j_R\}}(E_{j+1}-E_j)
 =\sum_{j\in\{W\tau_c\}_c\cap(n_u,n_v]\cap\mathbb N}
      \mathbf1_{\{j_R=j\}}E_j.
\end{equation}

\emph{A level at which the residual vanishes.} The term
\eqref{eq:cs-identgram-residual} has no factor $1/W$, and it need not be
small at small levels $R$, since $E_j$ compares $F$ at two states that
differ by a birth jump of order one. At large levels, however, it
vanishes identically. Put
\begin{equation}\label{eq:cs-identgram-threshold}
 B_J=\sum_cd_*m_c,\qquad R_0=R_\chi+B_J+1;
\end{equation}
any $R_0\ge R_\chi+B_J$ works, since $\chi$ and its first three
derivatives are continuous and vanish for $|x|>R_\chi$, hence also at
$\pm R_\chi$. Since $0<b<1/2$ and $m_c\ge0$, we have $d_*=1/(1-b)>0$, so
every partial sum $D$ of the Gram birth jumps satisfies
$0\le\Re\Tr D^M\le B_J$. On $\{j_R=j\}$ with $R\ge R_0$ we have
$\Re\Tr M_{W,j}\ge R$, hence $\Re\Tr(\mathsf G_{W,j}-D)^M\ge R_\chi+1$ for
both compensations. Each of the six tests contains the factor
$\chi(\vartheta(Y))$, which vanishes there, so both values of $F$ in $E_j$
vanish, whatever the complex-source and drift coordinates are. So at
these levels the process is stopped before a birth jump can be seen by
the test.

\emph{Conclusion of part~(f).} The discrete compensators match the births
in the Taylor expansion of part~(e) exactly, so
\eqref{eq:cs-identgram-residual} is the only remaining difference between
the two compensations. Together with the comparisons above and the
discrete statement of part~(e), this gives, for every $R\ge R_0$,
\begin{equation}\label{eq:cs-identgram-local}
 \E_{\widehat\mu_W}\Big[\phi_W\sum_{n_u\le j\le n_v}
     \mathbf1_{\{j<j_R\}}\big(\widetilde M_{F,W}(v_{j+1}^{u,v})
          -\widetilde M_{F,W}(v_j^{u,v})\big)\Big]\longrightarrow0.
\end{equation}
The level $R_0$ depends only on $\chi$ and on the birth data.

\emph{(g) Removing the localization and passing to the limit.} The
estimate on the tail of the trace stop, derived from the same summed
errors and Lindeberg bounds, includes the starting time $s$. For every
$\gamma>0$ it gives a level $R_{\rm tail}$ with
\[
 \limsup_W\widehat\mu_W\{\sigma_W^{R_{\rm tail}}\le H\}\le\gamma .
\]
We use the single fixed level $R_*=\max(R_{\rm tail},R_0)$. The exit
conditions are $\Re\Tr M_{W,j}\ge R_*$ and $\Re\Tr M_{W,j}\ge R_{\rm tail}$,
so $\sigma_W^{R_*}\ge\sigma_W^{R_{\rm tail}}$ on every path. Hence the tail
bound also holds at $R_*$, where \eqref{eq:cs-identgram-local} applies. On
each path, the alive cell sum in \eqref{eq:cs-identgram-local} telescopes
up to the first exit, including its overshoot, and it vanishes if the exit
precedes $u$. So it differs from $\widetilde M_{F,W}(v)-\widetilde
M_{F,W}(u)$ only on the event $\{\sigma_W^{R_*}\le H\}$. Replacing the
stopped cell sum by the unstopped increment therefore costs at most
\[
 4C_H\|\phi\|_\infty\widehat\mu_W\{\sigma_W^{R_*}\le H\}.
\]
We let $W\to\infty$ at the fixed level $R_*$, which gives
\[
 \limsup_W\big|\E_{\widehat\mu_W}\big[\phi_W(\widetilde M_{F,W}(v)
        -\widetilde M_{F,W}(u))\big]\big|\le4C_H\|\phi\|_\infty\gamma,
\]
and then we let $\gamma\downarrow0$. This is legitimate because $C_H$
depends neither on $R_{\rm tail}$ nor on $R_*$, and no limit in $W$ is
exchanged with a varying stop level. In particular, no source moment that
depends on $R$ is sent to infinity, and no functional of a stopping time
is passed through a finite-dimensional limit.

\emph{Adjoining the accumulated drift.} It remains to add the accumulated
drift $A^\chi$ to the given convergence of finite-dimensional
distributions; this does not follow from that convergence alone. Consider
a drift integrand $a(r,U)$ that is bounded and continuous in the state, and
depends on time through $J_G$, jointly continuously off the births. Let
$\pi_\delta(r)$ be the tags of a partition of mesh $\delta$, chosen off the
finite birth set, and write $I_W(t)=\int_0^t a(r,U_W(r))\,dr$ and
$I_{W,\delta}$ for the corresponding tagged sum. The Cauchy--Schwarz
inequality and Tonelli's theorem give
\[
 \E|I_W(t)-I_{W,\delta}(t)|^2
 \le H\int_0^t\E|a(r,U_W(r))-
           a(\pi_\delta(r),U_W(\pi_\delta(r)))|^2\,dr.
\]
For a fixed partition, convergence of the two-time finite-dimensional
distributions and boundedness pass the right side to the candidate limit.
The paths of the limit, after the birth jumps are subtracted, are
continuous, so the limiting integral tends to zero with the mesh, by
dominated convergence off the finite birth set. Thus
\[
 \lim_{\delta\downarrow0}\limsup_W\E|I_W-I_{W,\delta}|^2=0 .
\]
Finite tagged sums and bounded Lipschitz tests now give joint convergence
together with all accumulated cut drifts. The same reasoning applies to
finitely many vector-valued integrands.

\emph{The martingale identities.} We use the bounded continuous extensions
of $F$ and $\Lambda$ specified in part~(d). Expectations at the endpoints
converge at every time off the births. For the generator integral, we
adjoin the observation and its accumulated drift at the inner time $r$;
convergence holds for almost every $r$, and the integrand is bounded by
$\|\phi\|_\infty C_\Lambda$, so dominated convergence in $r$ passes the
integral to the limit as well. With the finite birth set
$\mathcal B=\{\tau_c-s\}$, this proves
\[
 \E_{\mathsf Q}\big[\phi((\operatorname{obs}_{i_l}(u_l))_{l<r})
                   (M_F(v)-M_F(u))\big]=0,
\]
whenever $0\le u\le v$, $u,v,u_l\notin\mathcal B$, $u_l\le u$, and
$\phi:\mathbb C^r\to\mathbb R$ is bounded and continuous. The global bound
$C_H$ also gives integrability at every time, including the births. These
are the tested martingale identities for all six coordinate tests. They
use only the stated finite-width estimates, the cap-limit premise,
feasibility and the derived stop tail. They do not use an unstopped trace
moment, an estimate of the resolvent at a birth, an exponential jet
moment, or the conclusion about the complex source that we are proving.
The later removal of a stop on the complex source uses the maximal moment
of Step~8 at a fixed Gram level.

\emph{Step 11: the limiting equation and the role of the stops.} The
product expansion preserves the strict gap of the source and the
summable remainders. Thus, at finite retained current rank, the noise
that survives in the limit is \eqref{eq:complex-source-sde}, with
brackets
\[
 d\mathcal K_{ij}\,d\mathcal K_{kl}=0,\qquad
 d\mathcal K_{ij}\,d\overline{\mathcal K_{kl}}
              =2\eta^2D_0 B_{ik}(C_0)_{lj}\,dt.
\]
All cross brackets with the Gram coordinates come from the same noise
$dA$. In particular $(d\mathcal K-d\mathcal K^*)/(2i)=-\eta\,dM_{\rm pr}$,
which is the infinitesimal form of $\Im\mathcal K=-\eta M_{\rm pr}$.

At a birth, the direct Schur estimates of Step~6 give the diagonal value
$am$ and vanishing cross entries, under the same pre-step stop, with
$L^2$ errors at most $C/\sqrt W$. In the martingale identities of Step~10 that enter the uniqueness
below, the births are compensated by the model-time jumps $J_G,J_K$.
Uniqueness of the finite law follows by adding this linear martingale
equation to the Gram uniqueness argument of
Proposition~\ref{prop:finite-gram-law}.

We stress what the stops give and what they do not. The argument above
supplies the stopped fourth-moment bound that is used later, with a
constant depending on $R$; it does not prove an unstopped fourth moment
of the trace by sending $R$ to infinity. The births are first identified
at fixed $h$, $K$ and $s$. The cap is then removed, and then
$h\downarrow0$ through levels $h\notin\mathcal N_s$, using the positive-source construction, the leakage
estimate and the coupling; only after that are the source sampling and
$s$ removed. No birth constant uniform in $h$ and no interchange of these
limits is asserted. For the bounded cutoff generator processes of
Step~10, the auxiliary trace stop is removed instead by the bound $C_H$,
which does not depend on $R$, and by the derived stop tail. This auxiliary stop is removed first, at the level
$\max(R_{\rm tail},R_0)$; it only identifies the cutoff generator
processes and is not the moment-based removal of the stops.

\emph{Step 12: the rank-free coupling for the complex source.} The
coupling \eqref{eq:rank-free-coupling} extends to the complex-source
equation. We expand the difference of the noise coefficients of the two
coupled systems at both frame factors and at $P'$. On $R\le K$, its
squared Hilbert--Schmidt norm is bounded by $C_KX$, since
$\norm{P_1'-P_2'}_{\HS}^2\le2RX$. Every side of the noise is paired with a
trace frame, so no ambient dimension appears. The Burkholder inequality
and the frame coupling then make the complex-source paths a Cauchy family
as the current cutoff is removed.

\subsection{The real part of the closing error}\label{subsec:s17-the-real-part-of-the-closing-error}

 Small positive mass in omitted channels must also control the real
part of the Green matrix. The boundary analytic estimate below gives
fractional moments of that real closing error from its imaginary part.
We prove the scalar inequality through positive harmonic
representations, then apply it to the closing-phase function and
polarize the matrix entries.

 For an analytic function $f$ on a neighborhood of the closed disk
with $f(0)=0$, and $0<p<1$,
\begin{equation}\label{eq:kolmogorov-conjugate}
 \int|f(e^{i\theta})|^p\frac{d\theta}{2\pi}
 \le C_p\left(\int|\Im f(e^{i\theta})|
                                      \frac{d\theta}{2\pi}\right)^p.
\end{equation}
To prove it split the imaginary boundary measure into positive and
negative parts. Their analytic Herglotz transforms, normalized by
zero real part at the origin, differ by $f$.
For a positive transform of mass $m$, its Cayley transform is a
Schur function vanishing at zero. The harmonic majorant of
$|(1+z)/(1-z)|^p$, integrable on the circle for $p<1$, gives
the bound $C_pm^p$. Finally use
$|x+y|^p\le|x|^p+|y|^p$.

Close omitted far channels by phase $z$ as in
Section~\ref{sec:positive-source-law}. The matrix
$\mathcal K(z)-\mathcal K(0)$ is analytic in the disk and vanishes
at zero. For its source response blocks $F_A,F_B$, let
$T_0=\Tr(F_AF_A^*+F_BF_B^*)$ and
$\ell_A=\Tr F_AF_A^*$.
Equation~\eqref{eq:complex-source-data}, the Hardy-series
calculation and Schatten Cauchy--Schwarz give
\[
 \int\norm{\Im(\mathcal K(e^{i\theta})-\mathcal K(0))}_1
                 \frac{d\theta}{2\pi}
       \le C\eta(\sqrt{T_0\ell_A}+\ell_A).
\]
One first compares the closed covariance with $F_BF_B^*$,
then adds the missing $F_AF_A^*$; the correction has mean
squared norm $\ell_A$.
Apply \eqref{eq:kolmogorov-conjugate} to finitely many quadratic
forms and polarize to recover every entry. The leakage estimate
\eqref{eq:small-channel-leakage} and the first trace moment imply
\[
 \E\norm{\mathcal K(e^{i\Theta})-\mathcal K(0)}^p\to0,
                                                   \qquad 0<p<1 .
\]
This pays the real part after phase averaging.
It is not a pointwise estimate by the imaginary part.

The Gram coordinates are positive forms, and their closing errors are
paid as in the proof of Theorem~\ref{thm:source-process}. The current
coordinate $V_{\rm pr}^*JV_{\rm pr}$, $J=\mathbb S^*P_n\mathbb S-P_n$ at
the final length, is signed and involves the outgoing side, so it needs
one more estimate. We use the notation of Step~2 at $j=n_T$, with
$\mathbb S=\mathbb S_{n_T}$, $A_-=\Pi\Pi^*$ and $A_+=\Pi^*\Pi$. Since
$P_n\Pi^*=0$ and $\mathbb S$ is unitary, $P_n\mathbb SZ^*=0$, and the energy identity
$CC^*=I+Z+Z^*$ of Step~2 gives, for a column $x$ of the open probe frame,
the closed outgoing near factor
$P_n\mathbb SC\,C^*x=P_n\mathbb Sx+f$, $f=P_n\mathbb SZx$, while the
closed incoming near factor $P_nC^*x$ enters the Gram coordinate $B$. So
the closing error of $(V_{\rm pr}^*JV_{\rm pr})_{kl}$ is minus that of
$B_{kl}$ plus a term of modulus at most
$\|x_k\|\|f_l\|+\|f_k\|\|x_l\|+\|f_k\|\|f_l\|$.

The phase average of $\|f\|^2$ is bounded by the \emph{output leakage}
$\|A_+\mathbb Sx\|^2$. Put $N_\Pi=\Pi\mathbb SA_-$. Since $\Pi=A_-\Pi$,
$Z=z(I-zN_\Pi)^{-1}\Pi\mathbb S$, and $\|N_\Pi\|=\|A_+\mathbb SA_-\|<1$
almost surely, as noted before \eqref{eq:gaussian-closing-error}; so
$z\mapsto Zx$ is analytic near the closed unit disk and vanishes at $0$.
The vector $y=Zx$ satisfies $A_-y=y$ and $y=z(w+N_\Pi y)$, $w=\Pi\mathbb
Sx$, and $\|A_+\mathbb Sy\|=\|N_\Pi y\|$. As $P_n\perp A_+$ and $\mathbb
S$ is unitary, for $|z|=1$
\[
 \|P_n\mathbb Sy\|^2\le\|y\|^2-\|N_\Pi y\|^2
 =\|w\|^2+2\Re\langle w,N_\Pi Zx\rangle ,
\]
and the last term has phase average zero. Hence, pathwise,
\begin{equation}\label{eq:output-cavity}
 \E_z\|P_n\mathbb SZx\|^2\le\|\Pi\mathbb Sx\|^2=\|A_+\mathbb Sx\|^2 .
\end{equation}

Summing over the probe columns, the phase average of $\sum_c\|f_c\|^2$
is at most $\ell^{\rm out}=\|A_+\mathbb SV_{\rm pr}\|_{\rm HS}^2$, while
\eqref{eq:small-channel-leakage} controls the input leakage
$\ell^{\rm in}=\|A_-V_{\rm pr}\|_{\rm HS}^2$. These have the same law.
Indeed, conjugating every Gaussian entry of the path and the probes,
$\omega\mapsto\bar\omega$, preserves $\mu_W$ and, the blocks being
Hermitian, transposes every Schur matrix and the full inverse $G$. With the convention
$V_{\rm pr}=V_{\rm ports}^*G^*\mathcal J$ of
\eqref{eq:complex-source-update}, the boundary identity of Step~1 gives
$V_{\rm ports}^*G\mathcal J=\mathbb SV_{\rm pr}$, and
\[
 \mathbb S(\bar\omega)=\mathbb S(\omega)^T,\qquad
 A_-(\bar\omega)=\overline{A_+(\omega)},\qquad
 V_{\rm pr}(\bar\omega)=\overline{\mathbb S(\omega)V_{\rm pr}(\omega)},
\]
where $A_+$ is the far output projection at time $s$, fixed in time, and $V_{\rm ports}$
is real. Thus $\ell^{\rm in}(\bar\omega)=\ell^{\rm out}(\omega)$.
Integrating the probes, $\E\ell^{\rm in}=\sum_cm_c\E\Tr(Q_cA_-)$, with
$Q_c$ the transported tag of the probe $c$, so
\eqref{eq:small-channel-leakage} gives
$\lim_{\rm cutoff}\limsup_W\E\ell^{\rm out}=0$. Since
$\E\|x_c\|^2=m_c\E\Tr Q_c$ is bounded by
\eqref{eq:full-tag-clock-moments}, Cauchy--Schwarz and
\eqref{eq:output-cavity} show that the expected outgoing part of the
closing error of the current tends to zero in the same order of limits.
With the Gram coordinates this pays the closing error of every
coordinate in the theorem; these are probe blocks. Entries between a
probe column and a current column of the closed frame enter the probe
current only through
$V_{\rm pr}^*J^{\rm cl}V_{\rm pr}=(M\mathsf DM)_{\rm pr}$, which is exact
because $J^{\rm cl}=V\mathsf DV^*$ (Step~10), so they need no closing
estimate. Together with finite-current admission and the extended coupling,
it completes the proof of Theorem~\ref{thm:complex-source-law}.

\end{proof}

\subsection{Uniform raw moments and extended probes}\label{subsec:s17-uniform-raw-moments-and-extended-probes}

 We justify replacing spatially extended probes by finitely many
sampled tags. Uniform raw moments control the error even for coincident
tags or births near the entrance. We propagate these moments with the
exact source update, then use Gaussian replacement and sampling to
recover the joint matrix law.

 For normalized tags, including coincident independent tags and
births arbitrarily near zero, the preceding raw update bounds give
\begin{equation}\label{eq:raw-tag-moments}
 \sup_{j,k\le TW}\E|\xi_j^*G\xi_k'|^p
                \le C_{p,T}(m_jm_k)^{p/2}.
\end{equation}
At birth this is the uniform fresh inverse estimate.
During transport the scalar error has
$L^p$ size $C_p(\rho^j+W^{-1})$ after a fixed finite entrance;
its sum through $TW$ is bounded. H\"older and Minkowski control
the means in \eqref{eq:complex-source-update}, while Burkholder
controls its centered increments using the tag trace moments.
The superexponential exit estimate \eqref{eq:passive-superexp-exit} removes stops.
For coincident tags polarize the quadratic forms with
$g\pm h,g\pm ih$. No inverse moment of a bare real Green matrix
is used.

Let $Q=\diag(q_jI_W)$, $q_j\ge0$ bounded, and
$A_Q=W^{-2}Q^{1/2}GQ^{1/2}$.
Sample $J$ with probabilities $q_j/(Wm)$,
$m=W^{-1}\sum_jq_j$, and set
$\xi=\sqrt{m/W}\,\iota_Jg$. Then
$\E(\xi\xi^*\mid G)=Q/W^2$ exactly.
Conditional Gaussian moments and Jensen applied to
\eqref{eq:raw-tag-moments} give
\begin{equation}\label{eq:extended-source-moments}
 \sup_W\E[|\Tr A_Q|^p+\norm{A_Q}_{\HS}^p]<\infty.
\end{equation}
The trace is the mean sampled quadratic form and the squared
Hilbert--Schmidt norm is the mean square independent two-tag
bilinear form. The same argument covers mixed source profiles.

We record a replacement estimate with the full quadratic
third derivative, which will also be used in the ring construction.
Let $f$ be a bounded $C^3$ test of a fixed finite collection of
complex quadratic and bilinear source entries, viewed in real
coordinates. Replace an extended Gaussian probe of covariance
$Q/W^2$ by $M^{-1/2}\sum_{a=1}^M\xi_a$.
At a single replacement, with the other probes fixed, the argument
of $f$ is $K(t)=K(0)+tA+t^2C$, and
\begin{equation}\label{eq:full-quadratic-third-derivative}
 \frac{d^3}{dt^3}f(K(t))
 =D^3f(K(t))[A+2tC,A+2tC,A+2tC]
       +6D^2f(K(t))[A+2tC,C].
\end{equation}
The uniform moments above bound its expected absolute value for
$0\le t\le M^{-1/2}$ by a constant times
\[
 \E\big[\norm A^3+\norm C^3+\norm A\norm C+\norm C^2\big].
\]
In particular the $\norm C^2$ term is retained.
The first two Taylor expectations agree, because the tag and
Gaussian replacements have equal mean, covariance and zero
pseudocovariance. Each replacement costs $CM^{-3/2}$, hence
\begin{equation}\label{eq:global-source-replacement}
 |\E f(\text{extended probes})-\E f(\text{sampled probes})|
                                      \le C_fM^{-1/2},
\end{equation}
uniformly in $W$.
Partial sums are conditionally Gaussian after fixing their tag
locations, with covariances that are convex averages of block
covariances; Jensen and \eqref{eq:raw-tag-moments} bound every
mixed moment used in this derivative estimate.

For the bounded exponential test one also has exactly
\[
 \det(I+iyA_Q)^{-1}
   =\E_\zeta e^{-iy\zeta^*G\zeta},\qquad
 \zeta\sim{\rm CN}(0,Q/W^2),\quad y\ge0,
\]
whose modulus is at most one because $\Im G\preceq0$.
If $a=v^*Gx+x^*Gv$ and $c=x^*Gx$, its third derivative is
$[-6y^2c(a+2tc)+iy^3(a+2tc)^3]e^{-iy(v+tx)^*G(v+tx)}$.
Thus the same replacement bound applies to fixed products of
these determinants and their conjugates, with independent sample
groups and the same actual path. The estimate retains the complete
determinant, including its first trace.
\section{Mixed relative eigenvectors of a fresh Gaussian matrix}
\label{sec:relative-vectors}
To identify the closing law in Part~II, we need the joint law of the
eigenvector components of the fresh collar in several boundary directions at
once. This section shows that these components are asymptotically independent
Gaussian vectors, relative to their covariance.
Theorem~\ref{thm:relative-mixed-vectors} states that, for finitely many
eigenvalue labels near $E$, the corresponding eigenvector components converge
jointly to independent standard circular Gaussian vectors, independently of
bounded tests of the eigenvalues. The error is a negative power of $W$,
uniformly in the deformation and in the directions. Combined with
Theorem~\ref{thm:external-sine}, this gives the decorated sine process used
for the collar.
In Subsection~\ref{subsec:s18-relative-estimates-along-the-short-gauss} we choose scales on which the eigenvector flow can
be localized while its covariance stays nearly constant. In
Subsection~\ref{subsec:s18-the-complex-reference-process-and-its-lo} we compare with a reference flow that keeps short
rotations and damps long ones, and prove that it localizes
(Lemma~\ref{lem:t26-relative-localization}). In Subsection~\ref{subsec:s18-the-two-generator-errors} we bound
the two errors in the generator and complete the proof.

Let $H_s=F+G_s$,   $\norm F\le \Naux ^A$, with the fixed regularity
rectangle \eqref{eq:fresh-regularity} at final variance $c$ and
center $E$. Put
$\Gamma=c|F-E-cm_c(E)|^{-2}$.
Let   $V\in\C^{\Naux \times r}$, $r$ fixed, satisfy $V^*\Gamma V=I_r$.
Give eigenvectors independent uniform phases and define
  $z_i(s)=\sqrt \Naux \,V^*u_i(s)$.

\begin{theorem}[Joint Gaussian law of eigenvector components]
\label{thm:relative-mixed-vectors}
Fix the number $r$ of directions, a number $k$ of eigenvalue labels, and
a bound on the size of the test functions and of their first two
derivatives. There are an exponent $\chi>0$, chosen after these fixed
quantities, and a set $\mathcal G$ of good paths of the full eigenvalue
process, whose complement has probability at most $\Naux^{-D}$ for any
prescribed fixed $D$, with the following property.
 
Let $i_1,\ldots,i_k$ be distinct labels whose terminal eigenvalues lie in a
fixed microscopic window around $E$, and let $g_1,\ldots,g_k$ be
independent standard circular Gaussian vectors in $\C^r$. Then for every
event $B\subset\mathcal G$ that is determined by the eigenvalue path, and
every bounded test function $f$ with two bounded derivatives within the
fixed bounds,
\begin{equation}\label{eq:relative-mixed-vectors}
 \Bigl|\E\bigl[f(z_{i_1}(c),\ldots,z_{i_k}(c))\,\one_B\bigr]
       -\E\bigl[f(g_1,\ldots,g_k)\bigr]\Prob(B)\Bigr|
 \le C\Naux^{-\chi}\Prob(B).
\end{equation}
The constant $C$ and the exponent $\chi$ (independent of $D$) are uniform over the regular
deformation class and over all normalized frames $V$.

Moreover, the decorated local eigenvalue process converges. Put
$\rho=\pi^{-1}\Im m_c(E)$ and $x_j=\Naux\rho(\lambda_j(c)-E)$, and let
the frames $V=V_\Naux$ satisfy only $V^*\Gamma V=I_r+o(1)$ as
$\Naux\to\infty$. Fix $R,L_\phi<\infty$ and a Borel function $\phi\ge0$ on
$\R\times\C^r$ that vanishes for $|x|>R$ and such that every
$\phi(x,\cdot)$ is of class $C^2$ with $C^2$ norm at most $L_\phi$. Then
\begin{equation}\label{eq:decorated-sine}
 \lim_{\Naux\to\infty}\E\exp\Bigl\{-\sum_j\phi\bigl(x_j,z_j(c)\bigr)\Bigr\}
 =\E\exp\Bigl\{-\sum_j\phi(y_j,g_j)\Bigr\},
\end{equation}
where $(y_j)$ is the unit-intensity sine process and the $g_j$ are
independent standard circular Gaussian vectors in $\C^r$, independent of
$(y_j)$. The convergence is uniform over such $\phi$, over the regular
deformation class and over frames with a common rate in
$V^*\Gamma V\to I_r$. Quantitatively, for every $\varepsilon>0$ (and fixed $r,R,L_\phi$) there are $\delta>0$ and $\Naux_0$ such that the two sides of \eqref{eq:decorated-sine} differ by at most $\varepsilon$ whenever $\Naux\ge\Naux_0$ and $\norm{V^*\Gamma V-I_r}\le\delta$, uniformly over such $\phi$ and the regular deformation class. The same holds with $x_j=\Naux\rho_\Naux(\lambda_j(c)-E)$ for any $\rho_\Naux$ with $|\rho_\Naux/\rho-1|\le\delta$, in particular whenever $\rho_\Naux/\rho\to1$. In particular the marked point process
$\sum_j\delta_{(x_j,z_j(c))}$ converges in law to the sine process with
independent standard Gaussian decorations.
\end{theorem}

The above theorem claims that near a regular energy, the components of finitely many
eigenvectors in finitely many fixed directions, normalized by the
deterministic covariance, behave like independent standard Gaussian
vectors. They are asymptotically independent of each other and of the
eigenvalues, with an error that is a small negative power of the
dimension. Independence from the eigenvalues is expressed by allowing an
arbitrary event $B$ of the eigenvalue path in
\eqref{eq:relative-mixed-vectors}; this is the same as testing against
any bounded function of the eigenvalues. The statement is formulated
through these integrals over events, with an error proportional to
$\Prob(B)$, so it does not rely on choosing a particular version of the
conditional law given the eigenvalue path. It holds on the good paths
$\mathcal G$, and the bad paths have arbitrarily small polynomial
probability, as specified in the proof.
This is a fresh full-GUE statement, not a theorem for the band matrix.

\subsection{Relative estimates along the short Gaussian flow}\label{subsec:s18-relative-estimates-along-the-short-gauss}

 We choose scales on which the eigenvector flow can be localized while
its covariance remains nearly constant. The estimates below control
individual components, averages over label intervals, and the relative
covariance error. We split off a short part of the original Gaussian
flow, apply the relative local law, and record the inequalities needed
by the reference process.

 Choose the fixed exponents

\begin{equation}\label{eq:relative-vector-scales}
\begin{aligned}
 \alpha=\frac{1}{100},\quad\theta=\frac{1}{500},\quad
 &\kappa=\frac{13}{5000},\quad\epsilon=\frac{1}{10000},\\
 h=\Naux ^{-\alpha},\quad t=\Naux ^{-\theta},\quad
 &\ell=\lfloor \Naux ^{1-\kappa}\rfloor,\\
 &L=\lceil \Naux t\Naux ^\epsilon\rceil,\quad w=\lceil \Naux \sqrt t\,\Naux ^{2\epsilon}\rceil.
\end{aligned}
\end{equation}
 
Split the same target Gaussian exactly as $H_{c-t}+G_t$.
Uniformly on $s\in[c-t,c]$, the relative coarse law
of Lemma~\ref{lem:coarse-relative}, polarized in $V$, gives
with arbitrary fixed failure power

\begin{align}
 \norm{z_i(s)}^2&\le C_r\Naux h
                   &&(|\lambda_i(s)-E|\le r_0/4),\notag\\
 \left\|\sum_{i\in I}[z_i(s)z_i(s)^*-\Gamma_i(s)]\right\|
      &\le C_r\Naux ^{1-\alpha+\epsilon},&
 \Gamma_i(s)&=V^*\frac{\Im M_s(\lambda_i(s))}{\Im m_s(\lambda_i(s))}V,\notag\\
 \norm{\Gamma_i(s)-I_r}&\le C(|\lambda_i(s)-E|+t).
                                             \label{eq:relative-interval}
\end{align}

The middle bound holds for every rank interval in that window.
For its proof integrate the relative resolvent estimate on a
rectangle with bottom height $h$. The vertical and upper sides
cost   $O(\Naux ^{2/3}\log \Naux )$  after normalization by   $\Naux $  is undone.
Positivity bounds each lower endpoint strip by   $C\Naux h$.
The scalar local law matches energy endpoints to rank endpoints.
Integration by parts transfers the smooth deterministic covariance
from its spectral integral to the rank sum.
The first bound follows directly from
$|u_i^*Vv|^2\le h\,v^*V^*\Im G_s(\lambda_i+ih)Vv$.
For the last, differentiating $sM_s^*M_s$ at fixed positive
imaginary Dyson shift bounds its quadratic-form derivatives by a
constant times those quadratic forms. This is relative stability;
it incurs no $\norm V$ factor.

For a central label $j$ the distant energy region satisfies

\begin{equation}\label{eq:relative-distant}
 \sum_{|\lambda_i-E|>r_0/4}
 \frac{z_iz_i^*}{\Naux (\lambda_i-\lambda_j)^2}\preceq C_rI_r.
\end{equation}

The left side is bounded by a fixed multiple of
$V^*\Im G_s(E+i)V$, controlled by the relative coarse law at
height one. Inside the window, for $|i-j|\ge\ell$, scalar
counting at height $h$ gives
  $|\lambda_i-\lambda_j|\ge c|i-j|/\Naux $.
The frame norm itself is only polynomial:
  $\Gamma\succeq C \Naux ^{-2\max(A,1)}I$ implies
$\norm V\le C\Naux ^{\max(A,1)}$. We use that bound only outside
the localized region.

\subsection{The complex reference process and its localization}\label{subsec:s18-the-complex-reference-process-and-its-lo}

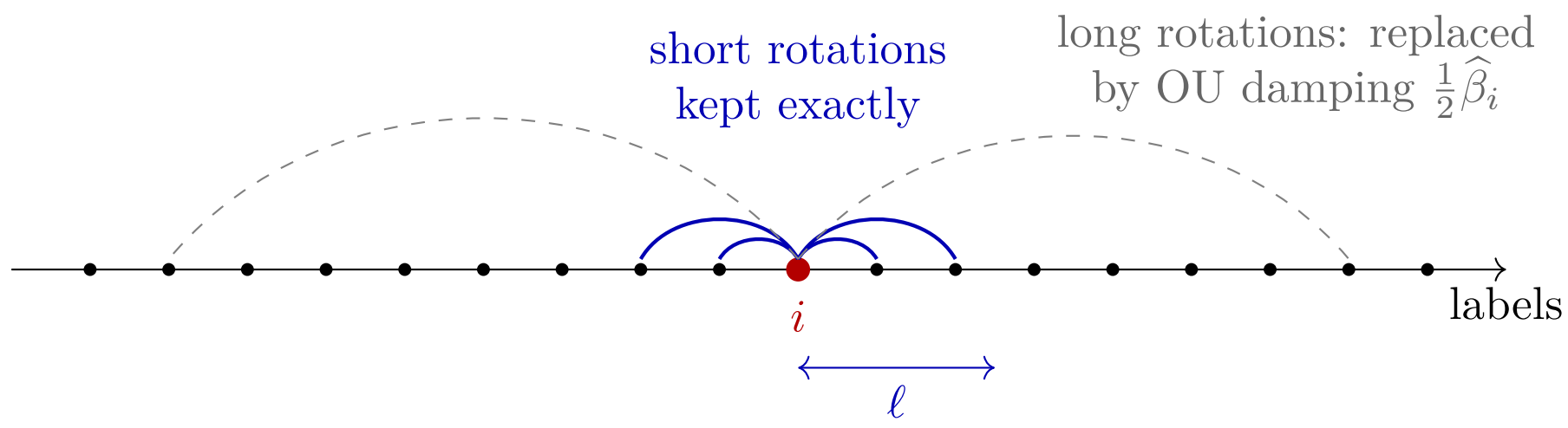
\begin{figure}[H]
\centering
\begin{tikzpicture}[x=1cm,y=1cm]
 \draw[->] (-5.6,0)--(5.8,0) node[below,font=\small] {labels};
 \foreach \x in {-5,-4.4,-3.8,-3.2,-2.6,-2.0,-1.4,-0.8,-0.2,0.4,1.0,1.6,2.2,2.8,3.4,4.0,4.6,5.2}
   \fill (\x,0) circle (1.4pt);
 \fill[red!70!black] (0.4,0) circle (2.6pt) node[below=3pt,font=\small] {$i$};
 \draw[thick,blue!70!black] (0.4,.08) to[bend left=60] (1.0,.08);
 \draw[thick,blue!70!black] (0.4,.08) to[bend left=60] (1.6,.08);
 \draw[thick,blue!70!black] (0.4,.08) to[bend right=60] (-0.2,.08);
 \draw[thick,blue!70!black] (0.4,.08) to[bend right=60] (-0.8,.08);
 \draw[dashed,black!50] (0.4,.08) to[bend left=50] (4.6,.08);
 \draw[dashed,black!50] (0.4,.08) to[bend right=50] (-4.4,.08);
 \draw[<->,blue!70!black] (0.4,-.75)--(1.9,-.75)
   node[midway,below,font=\small,blue!70!black] {$\ell$};
 \node[font=\small,blue!70!black,align=center] at (0.4,1.45)
   {short rotations\\ kept exactly};
 \node[font=\small,black!60,align=center] at (4.2,1.55)
   {long rotations: replaced\\ by OU damping $\tfrac12\widehat\beta_i$};
\end{tikzpicture}
\caption{The reference flow of Step~3. A label $i$ (red) is coupled to
every other label $j$ at rate   $c_{ij}=[\Naux (\lambda_i-\lambda_j)^2]^{-1}$.
The reference flow keeps all short rotations, with $|i-j|<\ell$ (solid
blue), exactly as in the true flow, including those outside the bulk
window. The long rotations (dashed) are replaced by an independent
circular Ornstein--Uhlenbeck process at label $i$. Its damping rate
$\frac{1}{2}\widehat\beta_i$ is half the total long coupling
  $\beta_i\asymp \Naux /\ell$, bounded below by   $b_0\Naux /\ell$. This damping makes the
reference flow contract at rate   $b_0\Naux /(2\ell)$, which is the source of
localization.}
\label{fig:rl-reference}
\end{figure}

We   show that eigenvector components  travel only a short distance in
the   eigenvalue labels during the comparison time.
Lemma~\ref{lem:t26-relative-localization} gives this localization with
the eigenvalue drift retained. We derive the complex angular flow,
replace  long rotations by  damping and noise,   and estimate the resulting
weighted propagation.  Figure~\ref{fig:rl-reference}   illustrates the
reference process.

\emph{Step 1: the eigenvalues and the angular flow.} We write $a=c-t$
for the starting variance and

\[
 c_{ij}(s)=\frac1{\Naux (\lambda_i(s)-\lambda_j(s))^2}
\]

for the coupling between the labels $i$ and $j$. The eigenvalues have the
GUE normalization

\begin{equation}\label{eq:t26-vector-dbm}
 d\lambda_i=\Naux ^{-1/2}dB_i+
              \frac{1}\Naux \sum_{j\ne i}\frac{ds}{\lambda_i-\lambda_j}.
\end{equation}

The real Brownian motions $B_i$ driving the eigenvalues are independent of
the off-diagonal Brownian motions driving the eigenvectors. To see this,
we conjugate the matrix Brownian increment by the current eigenbasis.
Its Gaussian covariance does not change, and its diagonal and
off-diagonal coordinates are orthogonal. L\'evy's characterization,
applied to their joint continuous martingale, then gives the
independence. Consequently, conditionally on the whole eigenvalue path
and on the initial matrix, the angular evolution is the unitary Dyson
flow. Conditioning only on the eigenvalue path mixes the initial angular
law, but it does not change the angular generator.

We now write this generator explicitly. For $i<j$ let $X_{ij}$ and
$Y_{ij}$ be the real vector fields
\[
 X_{ij}z_i=z_j,\quad X_{ij}z_j=-z_i,\qquad
 Y_{ij}z_i=iz_j,\quad Y_{ij}z_j=iz_i,
\]
which rotate the pair of components $(z_i,z_j)$. The angular generator
is
\[
 \frac{1}{4}\sum_{i<j}c_{ij}(X_{ij}^2+Y_{ij}^2),
\]
that is, each real rotation has Stratonovich coefficient
$\sqrt{c_{ij}/2}$. In particular, for $|z_i^v|^2$ it gives
\[
 \sum_{j\ne i}c_{ij}\bigl(|z_j^v|^2-|z_i^v|^2\bigr),
\]
where $v$ denotes any fixed component, or linear combination of
components. So mass moves from label $i$ to label $j$ at rate $c_{ij}$,
and the rate is large when $\lambda_i$ and $\lambda_j$ are close.

We make this conditional description of the angular flow precise, for
joint functions of all labels and directions. We regard $G$ as a matrix Brownian motion on $[0,\infty)$, so
that $H_s=F+G_s$ is defined for all $s\ge0$. Let $\mathscr F_v$ be the
$\sigma$-algebra generated by $G_w$, $w\le v$, and let
$\mathscr E_v=\sigma(\lambda(w):w\ge v)$ be that of the eigenvalue path
after time $v$. A function $f$ on $(\C^r)^{\Naux}$ is \emph{phase
invariant} if $f((\vartheta_iz_i)_i)=f(z)$ whenever $|\vartheta_i|=1$ for
all $i$. For a real linear vector field $A$ on $(\C^r)^{\Naux}$ we write
$A^2f(z)=D^2f(z)[Az,Az]+Df(z)[A^2z]$, the second derivative of
$f(e^{\varrho A}z)$ at $\varrho=0$, and we put
\[
 \mathcal L_{\rm true}(v)=\frac14\sum_{i<j}c_{ij}(v)\bigl(X_{ij}^2+Y_{ij}^2\bigr).
\]
For $s>0$ and $n\ge1$ let
$\tau_n=\inf\{v\ge s:\min_{i\ne j}|\lambda_i(v)-\lambda_j(v)|\le1/n\}$.
Since $\sum_i\|z_i\|^2=\Naux\|V\|_{\HS}^2$, the array
$z(v)=(z_i(v))_i$ stays in the closed ball $\mathbb B_V$ of radius
$\sqrt\Naux\|V\|_{\HS}$.

\begin{lemma}[The angular martingale problem]
\label{lem:t26-angular-martingale}
Let $s>0$ and $n\ge1$. For each $\omega$ let $f_\omega$ be a
phase-invariant $C^2$ function on $(\C^r)^{\Naux}$ such that
$\omega\mapsto f_\omega(z)$ is $\mathscr F_s\vee\mathscr E_s$-measurable for
every $z$, and $f_\omega$ and its first two derivatives are bounded on $\mathbb B_V$
uniformly in $\omega$. Then, with $f=f_\omega$ throughout, for every $u\ge s$
\begin{equation}\label{eq:t26-angular-martingale}
 \E\Bigl[f(z(u\wedge\tau_n))-f(z(s))
   -\int_s^{u\wedge\tau_n}\mathcal L_{\rm true}(v)f(z(v))\,dv\Bigr]=0 ,
\end{equation}
and the variable inside is bounded. In particular
\eqref{eq:t26-angular-martingale} holds with $f$ replaced by
$\one_{\mathcal Q}f$ for every $\mathcal Q\in\mathscr F_s\vee\mathscr E_s$.
\end{lemma}

\begin{proof}
(a) \emph{Phases and smoothness.} Let $\mathcal U$ be the open set of
Hermitian matrices with simple spectrum, with smooth spectral projections
$P_k(H)$. Near $H_0\in\mathcal U$ choose for each $k$ a standard basis
vector $y$ with $y^*P_k(H_0)y>0$; then $u_k=P_ky/(y^*P_ky)^{1/2}$ is a
smooth unit eigenvector, and other choices differ by unimodular phases,
so for phase-invariant $f$ the function
$\Phi(H)=f\bigl((\sqrt\Naux V^*u_k(H))_k\bigr)$ is well defined and $C^2$
on $\mathcal U$. For $T_\vartheta z=(\vartheta_iz_i)_i$ and
$e^{i\varphi}=\bar\vartheta_i\vartheta_j$, conjugation by $T_\vartheta$
rotates the pair $(X_{ij},Y_{ij})$ by the angle $\varphi$, which leaves
$X_{ij}^2+Y_{ij}^2$ unchanged. Hence $\mathcal L_{\rm true}(v)f$ is phase
invariant, and the variable in \eqref{eq:t26-angular-martingale} does not
depend on the choice of eigenvectors.

(b) \emph{The Laplacian identity.} We use the real inner product
$\Tr(AB)$ on Hermitian matrices. In an orthonormal basis the increments of
$G$ have covariance $\Naux^{-1}$ times the identity; this is the
normalization of \eqref{eq:t26-vector-dbm}, and the gradient of
$\lambda_k$ is $P_k$. Fix $H=\sum_k\lambda_ku_ku_k^*\in\mathcal U$ and let
$z$ be its array. We claim that
\begin{equation}\label{eq:t26-array-laplacian}
 \sum_bD^2\Phi(H)[\mathsf E_b,\mathsf E_b]=2\Naux\,\mathcal L f(z),
 \qquad D\Phi(H)[P_k]=0\quad\text{for all }k,
\end{equation}
where $\mathsf E_b$ runs over an orthonormal basis and $\mathcal L$ is the
angular generator with the couplings of $\lambda(H)$. The left side does not depend on the basis; we use $P_k=u_ku_k^*$,
$\mathsf S_{ij}=(u_iu_j^*+u_ju_i^*)/\sqrt2$ and
$\mathsf A_{ij}=i(u_ju_i^*-u_iu_j^*)/\sqrt2$ for $i<j$. For small
$|\varepsilon|$, $H+\varepsilon P_k$ has simple spectrum and the same
spectral projections as $H$, so $\Phi(H+\varepsilon P_k)=\Phi(H)$. This
gives the second identity and $D^2\Phi[P_k,P_k]=0$. Next fix $i<j$, let
$J=u_ju_i^*-u_iu_j^*$ and $\gamma(\varrho)=e^{\varrho J}He^{-\varrho J}$,
with eigenvalues $\lambda$ and eigenvectors $e^{\varrho J}u_k$, so
$\Phi(\gamma(\varrho))=f(e^{\varrho X_{ij}}z)$. Moreover
$\gamma'(0)=[J,H]=\sqrt2(\lambda_i-\lambda_j)\mathsf S_{ij}$ and
$\gamma''(0)=[J,[J,H]]=2(\lambda_i-\lambda_j)(P_j-P_i)$. Differentiating
twice at $\varrho=0$, with $D\Phi[P_i]=D\Phi[P_j]=0$,
\[
 2(\lambda_i-\lambda_j)^2D^2\Phi[\mathsf S_{ij},\mathsf S_{ij}]=X_{ij}^2f(z).
\]
With $J=i(u_ju_i^*+u_iu_j^*)$ the array is $e^{\varrho Y_{ij}}z$,
$\gamma'(0)=\sqrt2(\lambda_i-\lambda_j)\mathsf A_{ij}$ and $\gamma''(0)$
is unchanged, so
$2(\lambda_i-\lambda_j)^2D^2\Phi[\mathsf A_{ij},\mathsf A_{ij}]=Y_{ij}^2f(z)$.
Summing over the basis and using $c_{ij}=[\Naux(\lambda_i-\lambda_j)^2]^{-1}$
proves \eqref{eq:t26-array-laplacian}.

(c) \emph{The martingale.} Let first $f$ be deterministic. Since $s>0$,
$H_s=F+G_s$ has an absolutely continuous law, hence almost surely simple
spectrum, and by the noncollision of complex Dyson Brownian motion
\cite[Section~3.3, (3.10) and Remark~1]{KT} almost surely $H_v\in\mathcal U$ for all $v\ge s$.
Before $\tau_n$ the couplings are at most $n^2/\Naux$, and by first-order
perturbation theory the gradient of $\Phi$ is bounded by a constant
depending on $n$, $\Naux$, $V$ and $f$. By It\^o's formula on the open set
$\mathcal U$ and \eqref{eq:t26-array-laplacian},
\[
 M_v=\Phi(H_{v\wedge\tau_n})
     -\int_s^{v\wedge\tau_n}\mathcal L_{\rm true}(w)f(z(w))\,dw
    =\Phi(H_s)+\int_s^{v\wedge\tau_n}\Tr\bigl(\nabla\Phi(H_w)\,dG_w\bigr),
 \qquad v\ge s,
\]
is a martingale, bounded on bounded intervals. Let $\mathbf w=(B_k)_k$ be
the Brownian motions of \eqref{eq:t26-vector-dbm}, so
$\mathbf w_k(v)-\mathbf w_k(s)=\sqrt\Naux\int_s^v\Tr(P_k\,dG_w)$, and by the
second identity in \eqref{eq:t26-array-laplacian} the covariation of $M$
with every $\mathbf w_k$ vanishes.

(d) \emph{Orthogonality to the eigenvalue noise.} Let
$\mathscr G_s=\mathscr F_s\vee\sigma(\mathbf w(v)-\mathbf w(s):v\ge s)$. Take
$\mathcal Q_0\in\mathscr F_s$, times $s=\varsigma_0<\dots<\varsigma_p$,
vectors $\zeta_1,\dots,\zeta_p\in\R^{\Naux}$, let $\zeta(v)=\zeta_k$ on
$(\varsigma_{k-1},\varsigma_k]$ and $\zeta(v)=0$ for $v>\varsigma_p$, and
put
\[
 \xi_v=\one_{\mathcal Q_0}\exp\Bigl(i\int_s^v\zeta(x)\cdot d\mathbf w(x)
        +\frac12\int_s^v|\zeta(x)|^2dx\Bigr).
\]
Then $\xi$ is a bounded martingale with $d\xi=i\xi\,\zeta\cdot d\mathbf w$,
so $d\langle M,\xi\rangle=0$, $M\xi$ is a martingale, and for
$\varsigma_*\ge\max(u,\varsigma_p)$
$\E[M_u\xi_{\varsigma_*}]=\E[M_u\xi_u]=\E[M_s\xi_s]=\E[M_s\xi_{\varsigma_*}]$.
The variables $\xi_{\varsigma_*}$, constant multiples of
$\one_{\mathcal Q_0}\prod_ke^{i\zeta_k\cdot(\mathbf w(\varsigma_k)-\mathbf w(\varsigma_{k-1}))}$,
form a class closed under multiplication and conjugation generating
$\mathscr G_s$, and $M_u-M_s$ is bounded, so the monotone class theorem
gives $\E[(M_u-M_s)\xi]=0$ for every bounded
$\mathscr G_s$-measurable $\xi$.

(e) \emph{The eigenvalue path.} For $n'\ge1$, the drift in
\eqref{eq:t26-vector-dbm} is Lipschitz on the convex set of ordered vectors
with consecutive gaps at least $1/n'$; extend it to a globally Lipschitz
field $b^{(n')}$. The equation
$d\Lambda=\Naux^{-1/2}d\mathbf w+b^{(n')}(\Lambda)\,dv$, $\Lambda(s)=\lambda(s)$,
has a pathwise unique strong solution, which is a measurable function of
$\lambda(s)$ and of the increments of $\mathbf w$ after $s$, hence
$\mathscr G_s$-measurable. By pathwise uniqueness $\Lambda=\lambda$ on
$[s,\tau_{n'}]$, and $\tau_{n'}\to\infty$ almost surely by noncollision.
Thus every $\lambda(v)$, $v\ge s$, is an almost sure limit of
$\mathscr G_s$-measurable variables, and $\mathscr F_s\vee\mathscr E_s$ is
contained in the completion of $\mathscr G_s$. With (d) this gives
$\E[\one_{\mathcal Q}(M_u-M_s)]=0$ for
$\mathcal Q\in\mathscr F_s\vee\mathscr E_s$, which is
\eqref{eq:t26-angular-martingale} for $\one_{\mathcal Q}f$.

(f) \emph{Random tests.} If
$f_\omega=\sum_l\one_{\mathcal Q_l}(\omega)f^{(l)}$ for a finite partition
into sets $\mathcal Q_l\in\mathscr F_s\vee\mathscr E_s$ and deterministic
tests $f^{(l)}$, the claim follows from (e) by summation. In general, the
$2$-jet $j_\omega=(f_\omega,Df_\omega,D^2f_\omega)$ on $\mathbb B_V$ is an
$\mathscr F_s\vee\mathscr E_s$-measurable map into a separable Banach space
of continuous functions (derivatives are limits of difference quotients). Choose $\omega_1,\omega_2,\dots$ such that $(j_{\omega_l})_l$ is
dense in its range, and let $f^{[q]}_\omega=f_{\omega_l}$ for the first
$l\le q$ minimizing $\|j_\omega-j_{\omega_l}\|$. These are tests of the
previous form, and their jets converge to $j_\omega$ for every $\omega$ as
$q\to\infty$. The variable in \eqref{eq:t26-angular-martingale} is linear in
the test, depends on it only through its jet on $\mathbb B_V$, and is bounded by
$C(1+(u-s)n^2)$ times the norm of this jet, with $C$ depending on
$\Naux$, $r$ and $V$. Dominated convergence completes the proof.
\end{proof}

\emph{Step 2: a good spectral event.} To control the couplings we use an
event $\mathcal R$ whose probability is at least   $1-\Naux ^{-D}$  for any fixed
$D$. On this event, simultaneously for all $a\le s\le c$,

\begin{align}
 |\lambda_i-\lambda_j|&\ge c_1|i-j|/\Naux 
                                      &&(|i-j|\ge\ell),\notag\\
 |\lambda_i-\lambda_j|&\le C_1\ell/\Naux 
 &&(|i-j|<\ell,\quad
   \min(|\lambda_i-E|,|\lambda_j-E|)\le r_0/8).
                                                     \label{eq:t26-long-gaps}
\end{align}

The first line says that labels far apart have eigenvalues far apart.
The second is used only when the pair meets the smaller bulk interval,
whose endpoints still lie inside $r_0/4$. Here and below we may reduce
$r_0$ by a fixed factor. Note that these estimates do not bound
consecutive gaps from below, and we never need such a bound.

For completeness we explain how the counts behind \eqref{eq:t26-long-gaps}
are proved. The first line follows from a global upper count, which
comes from spectral averaging over one block and Gaussian concentration
at height $h$. Spectral averaging gives
\[
 \E\Im\av{G_s(x+ih)}\le C
\]
for every real $x$, uniformly in $F$. The Gaussian moments of the
gradient give

\[
 \|\av{G_s}-\E\av{G_s}\|_p\le C_p\Naux ^{-1}h^{-2}.
\]

A polynomial mesh in $x$ and $s$ then gives, with an arbitrary failure
power,

\[
 \#\{i:\lambda_i(s)\in I\}\le C\Naux (|I|+h)
\]

on the polynomial spectral range. Lipschitz bounds for the resolvent and
the modulus of continuity of Brownian motion fill in the mesh. The same
mesh and the local scalar law give the lower count   $c\Naux |I|$  on bulk
intervals of length comparable to   $\ell/\Naux $, since $h=o(\ell/\Naux )$. These two
counts prove \eqref{eq:t26-long-gaps}.

The counts also control the total long coupling of each label,

\[
 \beta_i:=\sum_{|i-j|\ge\ell}c_{ij}\le\frac{C\Naux }\ell,
 \qquad \beta_i\ge\frac{b_0\Naux }\ell\quad (|\lambda_i-E|\le r_0/8),
\]

for a fixed $b_0>0$. The upper bound follows from the first line of
\eqref{eq:t26-long-gaps}, since then   $c_{ij}\le \Naux /(c_1^2|i-j|^2)$, and the
sum over $|i-j|\ge\ell$ is   $O(\Naux /\ell)$. For the lower bound we use the
$\ell$ labels at rank distance between $\ell$ and $2\ell$ on one bulk side.
By the lower count, each of them has   $|\lambda_i-\lambda_j|\le C\ell/\Naux $, so
$c_{ij}\ge \Naux /(C^2\ell^2)$, and their sum is at least   $b_0\Naux /\ell$.

We also include in $\mathcal R$ the event
\[
 \sup_{a\le s\le c}\max_i|\lambda_i(s)-\lambda_i(a)|\le C\sqrt t,
\]
again with arbitrary polynomial failure probability. Gaussian moments of
the matrix norm and a dyadic subdivision of the Brownian time interval
prove $\sup_{a\le s\le c}\|H_s-H_a\|\le C\sqrt t$ with this failure
probability, and by Weyl's inequality this matrix event implies the
displayed one, so the failure bound is unchanged. Thus $\mathcal R$ is an
event of the eigenvalue path. Nothing here compares the global
counting function with a deterministic rank, and no microscopic gap
estimate is used.

\emph{Step 3: the reference flow.} We keep \emph{all} short rotations
$|i-j|<\ell$, including those outside the regular bulk window. We
replace the long rotations by the circular Ornstein--Uhlenbeck process

\[
 d\widehat z_i=-\tfrac12\widehat\beta_i\widehat z_i\,ds
                     +\sqrt{\widehat\beta_i}\,dC_i,
 \qquad \widehat\beta_i=\max(\beta_i,b_0\Naux /\ell),
 \qquad \E dC_i dC_i^*=I_r\,ds.
\]

The floor   $b_0\Naux /\ell$  only matters outside the smaller bulk window. Inside
it, $\widehat\beta_i=\beta_i$ by Step~2. The short rotation matrix is not
cut at any bulk boundary. Conditionally on $\lambda$, the Brownian
motions of the reference rotations and its Ornstein--Uhlenbeck noises
are independent of the true angular flow.

\emph{Step 4: the solution of the reference flow.} The reference flow is
linear, so we can solve it. Conditionally also on the Brownian motions
of the reference rotations, its solution is
\[
 \widehat z(u)=\mathcal M_{s,u}z(s)+\Xi_{s,u},
\]
where $\Xi$ is a circular Gaussian with
\begin{equation}\label{eq:t26-reference-covariance}
 \operatorname{Cov}(\Xi_{s,u})
       =(I-\mathcal M_{s,u}\mathcal M_{s,u}^*)\otimes I_r,
 \qquad \E\Xi_{s,u}\Xi_{s,u}^{\mathsf T}=0.
\end{equation}
Here $\mathcal M$ solves the homogeneous Stratonovich equation with the
short skew-Hermitian rotation matrices and the damping
$\frac{1}{2}\operatorname{diag}(\widehat\beta_i)$. To see
\eqref{eq:t26-reference-covariance}, we compare two matrix equations. The
covariance of the independent Ornstein--Uhlenbeck integral, and the
matrix $I-\mathcal M\mathcal M^*$, satisfy the same linear equation with
zero initial data. In both, the inhomogeneous term is
$\operatorname{diag}(\widehat\beta_i)$. Indeed, the rotations are
skew-Hermitian, so the generator $A$ of $\mathcal M$ satisfies
$A+A^*=-\operatorname{diag}(\widehat\beta_i)$, and this produces exactly the
same inhomogeneous term in the equation for $I-\mathcal M\mathcal M^*$.
So the two matrices coincide. See Lemma~\ref{lem:t26-reference-backward} detailed proofs of various properties of this reference flow.

Moreover, pathwise, the rotations preserve the norm, and only the damping
changes it:

\[
 d\|\mathcal M_{s,u}v\|^2
 =-\sum_i\widehat\beta_i|(\mathcal M_{s,u}v)_i|^2du,
 \qquad \|\mathcal M_{s,u}\|\le
                 e^{-b_0\Naux (u-s)/(2\ell)}.
\]

The bound on the norm follows from   $\widehat\beta_i\ge b_0\Naux /\ell$. This
proves the contraction of the reference flow. It also proves that the
product density $\prod_i\pi^{-r}e^{-\|z_i\|^2}$ is exactly invariant,
because \eqref{eq:t26-reference-covariance} shows that a standard
circular Gaussian initial condition is mapped to a standard circular
Gaussian.

Notice that matrix $\mathcal M_{s,c}$ is the
propagator of the reference flow from time $s$ to the final time $c$.
Its entry $(\mathcal M_{s,c})_{kj}$ measures how much of the component at
label $j$ at time $s$ reaches label $k$ at time $c$. The following lemma says that
a terminal label $k$ near the energy $E$ receives essentially nothing
from labels outside the window $J$ of $w$ labels around $E$. So
information travels at most $w$ labels during the flow. The bound is
integrated in the starting time $s$, which is the form used in the
Duhamel formula.

\begin{lemma}[Localization with the eigenvalue drift retained]
\label{lem:t26-relative-localization}
Let $k_0$ be the first rank above $E$ at time $a$, and let
$J=\{j:|j-k_0|\le w\}$. Thus $J$ is measurable at time $a$. For each fixed
$D$, outside a set of eigenvalue paths of probability   $\Naux ^{-D}$, all
terminal labels $k$ with   $|\lambda_k(c)-E|\le C/\Naux $  satisfy

\begin{equation}\label{eq:reference-finite-speed}
 \int_a^c\E_{\rm aux}\sum_{j\notin J}
             |(\mathcal M_{s,c})_{kj}|^2\,ds\le \Naux ^{-D}.
\end{equation}

The constants are uniform in the deformation class. This integrated
bound, rather than a supremum over all time intervals, is what the
Duhamel comparison requires.
\end{lemma}

\begin{proof}
We want to show that mass starting far from the
energy $E$ does not reach the labels near $E$. To measure this we use a weight
that is of order one near $E$ and decays exponentially away from it, shown in
Figure~\ref{fig:rl2-weight} (Step~1), and we follow the weighted mass of the
reference flow (Step~2). It\^o's formula shows that this weighted mass grows at most by
a bounded factor. The key is that the drift of the eigenvalues cancels the
main term coming from the rotations (Step~3). Mass that starts outside the
window has exponentially small weight, so it remains negligible near $E$
(Steps~4--5).

\emph{Step 1: the weight.} We put   $d=L/\Naux $  and $\nu=d^{-1}$. We choose a
smooth function $\psi$ that agrees with $d-\sqrt{(x-E)^2+d^2}$ on
$|x-E|\le r_0/32$, is a negative constant on $|x-E|\ge r_0/16$, and
satisfies
\[
 |\psi'|\le C,\qquad |\psi''|\le C/d.
\]
It can be chosen even and nonincreasing in $|x-E|$. Near $E$, $\psi$ is
close to zero, and away from $E$ it decreases like $-|x-E|$. We put
\[
 \phi_i(u)=e^{\nu\psi(\lambda_i(u))}.
\]
So $\phi_i$ is of order one when $\lambda_i$ is within $d$ of $E$, and it
decays like $e^{-\nu|\lambda_i-E|}$ further away. All nonconstant parts of
this weight lie strictly inside the regular window. In particular,
\[
 |\phi_i'|\le C\nu\phi_i,\qquad |\phi_i''|\le C\nu^2\phi_i.
\]

\begin{figure}[ht]
\centering
  \begin{tikzpicture}[x=1cm,y=1cm]
 \draw[->] (-5.4,0)--(5.6,0) node[below,font=\small] {energy};
 \draw[->] (0,-.1)--(0,2.6) node[left,font=\small] {$\phi$};
 \draw[very thick,blue!70!black,domain=-4.4:4.4,samples=150]
   plot (\x,{2*exp(0.35-sqrt(\x*\x+0.1225)*1.0)});
 \draw[very thick,blue!70!black] (-5.2,.03)--(-4.4,.03);
 \draw[very thick,blue!70!black] (4.4,.03)--(5.2,.03);
 \fill[red!70!black] (0,0) circle (2pt) node[below=2pt,font=\small] {$E$};
 \fill[green!50!black,opacity=.18] (-.35,0) rectangle (.35,2.5);
 \node[font=\small,green!50!black,anchor=west] at (.45,2.2) {terminal labels, $\phi_k\ge\tfrac12$};
 \draw[decorate,decoration={brace,amplitude=4pt,mirror}] (-2.4,-.35)--(2.4,-.35)
   node[midway,below=4pt,font=\small] {window $J$: energies within $Cw/\Naux $};
 \node[font=\small,anchor=west,align=left] at (2.6,1.15)
   {outside $J$:\\ $\phi_j\le e^{-c_2\Naux ^\epsilon}$};
 \draw[<->] (-.35,.9)--(.35,.9);
 \node[font=\small,anchor=east] at (-.4,.9) {$\sim d$};
\end{tikzpicture}
 \caption{The weight in the proof of Lemma~\ref{lem:t26-relative-localization},
drawn schematically. As a function of the energy, $\phi=e^{\nu\psi}$ is of
order one within distance   $d=L/\Naux $  of $E$ (green) and decays like
$e^{-\nu|x-E|}$ away from it. It is constant beyond $r_0/16$. Labels outside
the window $J$ have energies at distance at least   $cw/\Naux $  from $E$, so their
weight is at most   $e^{-c_2\Naux ^\epsilon}$. The proof shows that the weighted
mass $S$ grows at most by a bounded factor. Mass that starts outside $J$
therefore has exponentially small weight at all times, so it cannot reach
the terminal labels near $E$, where $\phi_k\ge1/2$.}
\label{fig:rl2-weight}
\end{figure}
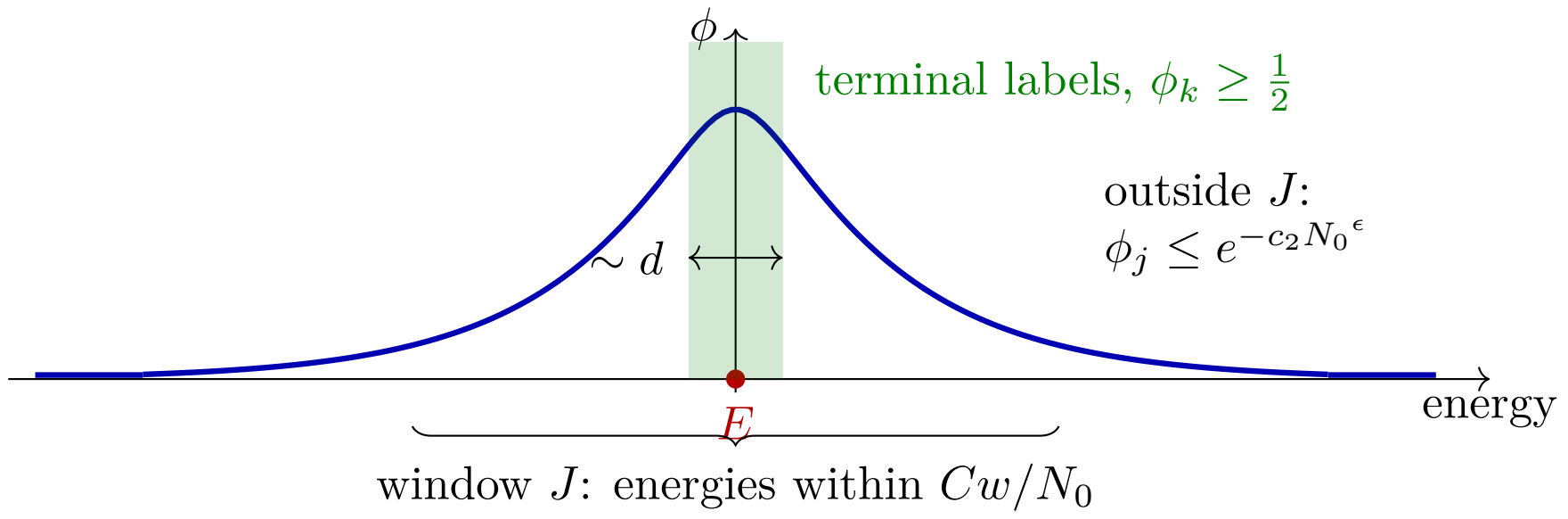

\emph{Step 2: the mass equation.} We fix the starting time $s$ and follow
the mass that starts outside $J$. For each label $i$ we put
\[
 v_i(u)=\E_{\rm aux}\sum_{j\notin J}|(\mathcal M_{s,u})_{ij}|^2 ,
\]
which is nonnegative. Conditionally on the eigenvalue path, it solves the
linear equation
\[
 \partial_uv_i=
       \sum_{0<|i-j|<\ell}c_{ij}(v_j-v_i)-\widehat\beta_iv_i,
 \qquad v_i(s)=1_{\{i\notin J\}}.
\]
The first term moves mass between labels at distance less than $\ell$,
and the second is the killing by the damping of the reference flow. So
$v$ is adapted to the eigenvalue filtration and has finite variation in
$u$.

\emph{Step 3: Itô's formula for the weighted mass.} We put
\[
 S(u)=\sum_i\phi_i(u)v_i(u),
\]
the mass weighted by $\phi$. We apply It\^o's formula to $S$, stopped
when \eqref{eq:t26-long-gaps} or the coarse count event fails. Here the coarse count event is the simultaneous event on the
stated time interval where both bounds in \eqref{eq:t26-long-gaps}
hold, $\beta_i\le C\Naux/\ell$, and
$\beta_i\ge b_0\Naux/\ell$ for $|\lambda_i-E|\le r_0/8$;
it also includes the displacement bound
$\sup_s\max_i|\lambda_i(s)-\lambda_i(a)|\le C\sqrt t$ on that interval..
The constants $c_1,C_1,b_0>0$ are chosen before $D$ and the width
threshold, and failure has outer probability at most $\Naux^{-D}$.
In particular the artificial floor is inactive in the bulk window. The
change of $S$ has three sources: the jumps of $v$, the motion of the
weights $\phi_i$ through the eigenvalues, and the killing.

The jumps of $v$ contribute $\sum_i\phi_i\sum_jc_{ij}(v_j-v_i)$, which we
symmetrize as $\sum_iv_i\sum_jc_{ij}(\phi_j-\phi_i)$. We expand
$\phi_j-\phi_i$ by Taylor's formula. For each retained pair, the linear
Taylor term

\[
 c_{ij}\phi_i'(\lambda_j-\lambda_i)=-\frac{\phi_i'}{\Naux (\lambda_i-\lambda_j)}
\]

cancels exactly the contribution   $\phi_i'/[\Naux (\lambda_i-\lambda_j)]$  of the
pair $(i,j)$ to the drift of $\phi_i$ in \eqref{eq:t26-vector-dbm}. This is
the key cancellation. For the remainder, two cases arise. If the weight
varies on the pair, then \eqref{eq:t26-long-gaps} gives
  $|\lambda_i-\lambda_j|\le C\ell/\Naux =o(d)$, so the weight is comparable at the
two points. If it does not vary, its difference vanishes. In both cases,
the Taylor remainder divided by the squared gap is at most
  $C\nu^2\phi_i/\Naux $  per retained neighbour.

The long pairs, with $|i-j|\ge\ell$, are not in the jump term, so their
part of the Dyson drift is not cancelled. Where $\phi_i'\ne0$, the first
inequality in \eqref{eq:t26-long-gaps} bounds it by

\[
 \frac{|\phi_i'|}{\Naux }
       \sum_{|i-j|\ge\ell}|\lambda_i-\lambda_j|^{-1}
 \le C\nu\log \Naux \,\phi_i.
\]

The It\^o correction from the eigenvalue noise is at most
  $C\nu^2\phi_i/\Naux $, and the killing term is nonpositive. Since there are at
most $2\ell$ retained neighbours, we conclude that

\[
 dS\le C\Bigl(\nu\log \Naux +\frac{\nu^2\ell}\Naux +\frac{\nu^2}\Naux \Bigr)S\,du+d\mathcal N,
 \qquad d\langle\mathcal N\rangle\le C\frac{\nu^2}\Naux S^2du,
\]

where $\mathcal N$ is a local martingale. The stopping makes the
martingale integrable. Positivity and Fatou's lemma then give the
Gronwall inequality for $\E S$, without any assumption on moments of
inverse gaps.

\emph{Step 4: the initial and the final weights.} We now compare the
weights at the two ends. By the upper count at time $a$ and the displacement
bound in $\mathcal R$, every $j\notin J$ satisfies

\[
 |\lambda_j(s)-E|\ge cw/\Naux  ,
\]

because the possible displacement $C\sqrt t+Ch$ is   $o(w/\Naux )$. By the lower
bulk count, the labels in $J$ have energies within   $Cw/\Naux $  of $E$ at every
time. Consequently, for $j\notin J$,

\[
 \phi_j(s)\le e^{-c_2w/(\Naux d)}\le e^{-c_2\Naux ^\epsilon}.
\]

We may multiply the initial data by the time-$s$ event on which this
holds. That event also records the counts at time $a$ and the
displacement up to time $s$, so it is measurable at time $s$, and we
apply the stopped inequality after this multiplication. The initial
weighted mass is then at most   $\Naux e^{-c_2\Naux ^\epsilon}$.

The Gronwall exponent is bounded, since

\[
 t\nu\log \Naux =\Naux ^{-\epsilon}\log \Naux ,\qquad
 t\nu^2\ell/\Naux =\Naux ^{\theta-\kappa-2\epsilon},
\]

and both tend to zero. At the terminal microscopic labels the energy is
within   $C/\Naux \ll d$  of $E$, so $\phi_k(c)\ge1/2$ for large   $\Naux $. Hence
$\sum_kv_k(c)\le2S(c)$ over these labels, and

\[
 \E\left[1_{\mathcal R}
       \sum_{k:\,|\lambda_k(c)-E|\le C/\Naux }
       \sum_{j\notin J}\E_{\rm aux}|(\mathcal M_{s,c})_{kj}|^2\right]
 \le C\Naux e^{-c_2\Naux ^\epsilon}.
\]

The sum over the terminal labels is bounded by the weighted sum itself,
so there is no reversal of the eigenvalue dynamics.

\emph{Step 5: conclusion.} We integrate in $s$ and use Tonelli's theorem
and Markov's inequality over the eigenvalue paths. Finally we adjoin the
complement of $\mathcal R$, whose probability decays faster than any
polynomial. This proves the lemma.
\end{proof}

Here we record the
properties of the reference flow that the Duhamel comparison uses in \eqref{eq:t26-reference-covariance}. The
comparison runs the true flow up to a time $s$ and the reference flow from
$s$ to the final time $c$, and it needs the expected final value of the
reference flow, as a function of its starting time and starting point, to
be smooth enough for It\^o's formula. We prove this here.
 
Fix an eigenvalue path with simple spectrum on $[a,c]$. Then the short
couplings $c_{ij}$, $0<|i-j|<\ell$, and the rates $\widehat\beta_i$ are
continuous functions of time on $[a,c]$. Explicitly, the reference flow is
the linear Stratonovich equation
\[
 d\widehat z=\sum_{i<j,\ j-i<\ell}\sqrt{c_{ij}/2}\,
   \bigl(X_{ij}\widehat z\circ dW_{ij}+Y_{ij}\widehat z\circ dW'_{ij}\bigr)
   -\tfrac12\operatorname{diag}(\widehat\beta)\widehat z\,dv
   +\operatorname{diag}(\widehat\beta)^{1/2}dC .
\]
The first sum consists of the random rotations between nearby labels,
driven by independent real Brownian motions $W_{ij},W'_{ij}$. The second
term damps each label at rate $\tfrac12\widehat\beta_i$, and the last term
adds the independent circular noises $C_i$ of Step~3. We write
$\widehat z^{s,z}$ for the solution on $[s,c]$ with $\widehat z^{s,z}(s)=z$.
Since the equation is linear, its solution splits into a part that
transports the initial value and a Gaussian part,
\[
 \widehat z^{s,z}(u)=\mathcal M_{s,u}z+\Xi_{s,u}.
\]
The generator of the reference flow at time $v$ is
\[
 \mathcal L_{\rm ref}(v)=\frac14\sum_{0<j-i<\ell}c_{ij}(v)\bigl(X_{ij}^2+Y_{ij}^2\bigr)
   +\sum_i\widehat\beta_i(v)\Bigl(-\tfrac12\,z_i\cdot\nabla_i
   +\tfrac14\Delta_i^{\rm real}\Bigr).
\]
Here $z_i\cdot\nabla_ig(z)$ is the derivative of $g$ in the direction of
the array whose only nonzero entry is $z_i$, at label $i$, and
$\Delta_i^{\rm real}$ is the Laplacian in the $2r$ real coordinates of
$z_i$. This is the normalization used in \eqref{eq:t26-floor-error}: the
circular noise has $\E|dC_{iq}|^2=ds$, which splits equally between the
real and the imaginary parts, so each real coordinate has variance
$ds/2$ and contributes $\tfrac14\widehat\beta_i$ times a second
derivative.
 
\begin{lemma}[The reference backward function]
\label{lem:t26-reference-backward}
Let $f$ be a bounded $C^2$ function on $(\C^r)^{\Naux}$ with bounded first
and second derivatives, and define the \emph{backward function}
\[
 g_s(z)=\E_{\rm aux}f(\mathcal M_{s,c}z+\Xi_{s,c}),\qquad s\in[a,c],
\]
the expected final value of $f$ under the reference flow started at time
$s$ from the point $z$. Then the following hold.
\begin{enumerate}[label=(\roman*)]
\item \emph{Bounds.} $\|D^jg_s\|_\infty\le\|D^jf\|_\infty$ for $j=0,1,2$
  and all $s$.
\item \emph{Continuity.} $g$, $Dg$ and $D^2g$ are jointly continuous on
  $[a,c]\times(\C^r)^{\Naux}$, and so is $\mathcal L_{\rm ref}(s)g_s(z)$.
\item \emph{Backward equation.} For every $z$, the map $s\mapsto g_s(z)$
  is continuously differentiable on $[a,c]$, with
  \[
   \partial_sg_s(z)=-\mathcal L_{\rm ref}(s)g_s(z).
  \]
\item \emph{Measurability in the path.} For each $s'\in[a,c]$, the map
  $(\text{path},s,z)\mapsto g_s(z)$ on $[s',c]\times(\C^r)^{\Naux}$ is
  jointly measurable with respect to $\sigma(\lambda(w):s'\le w\le c)$, on
  the set of eigenvalue paths with simple spectrum on $[s',c]$; this set
  belongs to the same $\sigma$-algebra.
\item \emph{Phase invariance.} If $f$ is invariant under multiplying each
  vector $z_i$ by a phase, then so is every $g_s$.
\end{enumerate}
\end{lemma}
 
In words: the backward function is as smooth as the test function, it
depends continuously on its starting time, it solves the backward equation
of the reference flow, it depends measurably on the eigenvalue path from
time $s'$ on, and it inherits the phase invariance of the test function.
These are exactly the properties needed to apply It\^o's formula to
$g_s$ along the true flow in the Duhamel comparison.
 
\begin{figure}[ht]
\centering
\begin{tikzpicture}[x=1cm,y=1cm,>=Stealth]
 \draw[->] (-.3,0)--(10.2,0) node[below,font=\small] {time};
 \foreach \x/\l in {0/a,3/s,6/{s'},9.4/c}
   {\draw[thick] (\x,-.1)--(\x,.1); \node[font=\small,below] at (\x,-.12) {$\l$};}
 \draw[->,thick,blue!70!black] (3,.6) to[bend left=25]
   node[midway,above,font=\small] {$\mathcal M_{s,s'},\ \Xi_{s,s'}$} (6,.6);
 \draw[->,thick,blue!70!black] (6,.6) to[bend left=25]
   node[midway,above,font=\small] {$\mathcal M_{s',c},\ \Xi_{s',c}$} (9.4,.6);
 \draw[->,thick,red!70!black] (3,-.75) to[bend right=12]
   node[midway,below,font=\small] {$\mathcal M_{s,c}=\mathcal M_{s',c}\mathcal M_{s,s'},\quad
     \Xi_{s,c}=\mathcal M_{s',c}\Xi_{s,s'}+\Xi_{s',c}$} (9.4,-.75);
 \fill[black!8] (3,.05) rectangle (6,.35);
 \fill[black!18] (6,.05) rectangle (9.4,.35);
 \node[font=\scriptsize] at (4.5,.2) {noise on $[s,s']$};
 \node[font=\scriptsize] at (7.7,.2) {noise on $[s',c]$};
\end{tikzpicture}
\caption{The flow property used in Step~3 of the proof. Running the
reference flow from $s$ to $c$ is the same as running it from $s$ to an
intermediate time $s'$ and then from $s'$ to $c$ (blue arrows); the
transport matrices compose, and the Gaussian parts add up after
transporting the first one (red). The two pieces are functions of the
noise on the disjoint intervals $[s,s']$ and $[s',c]$, so they are
independent. Averaging first over the second piece turns the final
expectation into $g_{s'}$ evaluated at the position at time $s'$, which
gives the identity $g_s(z)=\E_{\rm aux}g_{s'}(\widehat z^{s,z}(s'))$ and,
with It\^o's formula on $[s,s']$, the backward equation.}
\label{fig:reference-flow-property}
\end{figure}
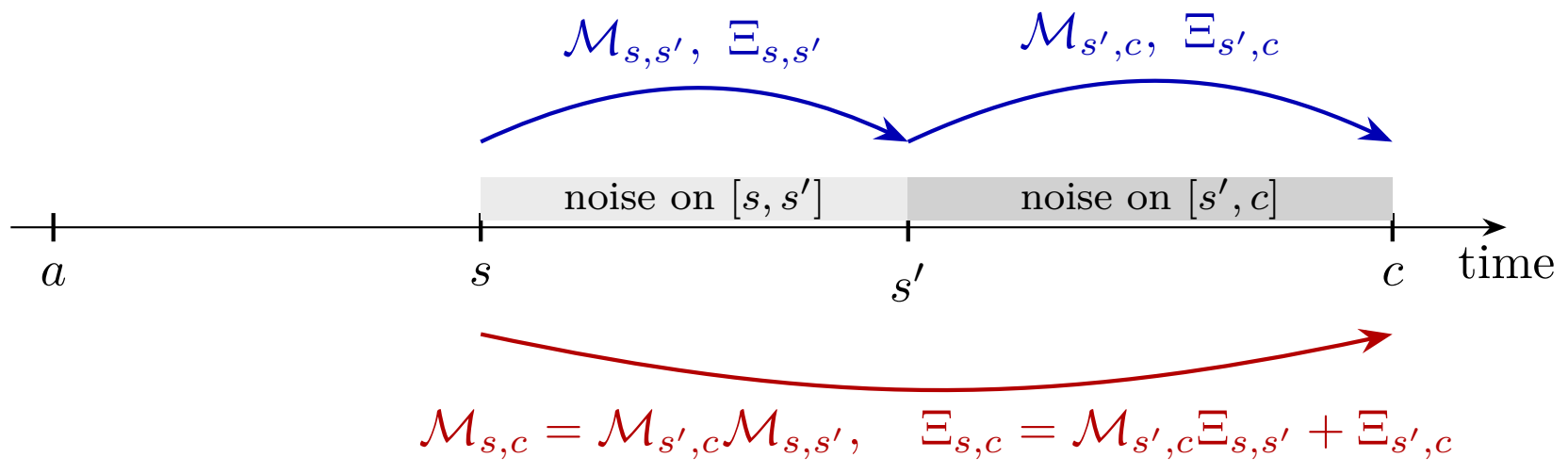
 
\begin{proof} The reference flow is a linear stochastic
equation, so we first note that it has a unique solution with all
moments (Step~0). The bounds (i) follow by differentiating under the
expectation, because the transport matrices are contractions (Step~1).
Continuity (ii) follows from a representation of all transport matrices
through one invertible fundamental matrix (Step~2). The backward equation
(iii) follows from the flow property, illustrated in
Figure~\ref{fig:reference-flow-property}, and It\^o's formula (Step~3).
Measurability (iv) follows from continuous dependence of the solution on
its coefficients (Step~4), and phase invariance (v) from a symmetry of
the noise (Step~5).
 
\emph{Step 0: existence and moments.} The equation is linear with
continuous coefficients, so it has a pathwise unique strong solution, and
for every $p$
\[
 \E\sup_{u\in[s,c]}\|\widehat z^{s,z}(u)\|^p\le C_p(1+\|z\|)^p .
\]
 
\emph{Step 1: the bounds (i).} We differentiate under the expectation
$\E_{\rm aux}$. Since $z$ enters only through $\mathcal M_{s,c}z$,
\[
 D^jg_s(z)[w_1,\dots,w_j]=\E_{\rm aux}D^jf(\mathcal M_{s,c}z+\Xi_{s,c})
 [\mathcal M_{s,c}w_1,\dots,\mathcal M_{s,c}w_j].
\]
The rotations preserve the norm and the damping only decreases it, so
$\|\mathcal M_{s,c}\|\le1$, and (i) follows.
 
\emph{Step 2: continuity (ii).} Put $\Phi_u=\mathcal M_{a,u}$, the
fundamental matrix of the flow started at time $a$. The rotation matrices
are traceless, so the Liouville formula for Stratonovich equations gives
\[
 \det\Phi_u=\exp\Bigl(-\tfrac12\int_a^u\sum_i\widehat\beta_i\Bigr)>0 ,
\]
and $\Phi_u$ is invertible. By uniqueness of the solution,
\[
 \mathcal M_{s,c}=\Phi_c\Phi_s^{-1},\qquad
 \Xi_{s,c}=\Xi_{a,c}-\mathcal M_{s,c}\Xi_{a,s}.
\]
Hence $s\mapsto(\mathcal M_{s,c},\Xi_{s,c})$ is almost surely continuous,
and the formula of Step~1, together with dominated convergence, gives the
joint continuity of $g$, $Dg$ and $D^2g$. The coefficients of
$\mathcal L_{\rm ref}$ are continuous in time, which gives the last claim.
 
\emph{Step 3: the backward equation (iii).} Let $a\le s\le s'\le c$. By
uniqueness of the solution, the flow has the composition property
\[
 \mathcal M_{s,c}=\mathcal M_{s',c}\mathcal M_{s,s'},\qquad
 \Xi_{s,c}=\mathcal M_{s',c}\Xi_{s,s'}+\Xi_{s',c},
\]
and the pairs $(\mathcal M_{s,s'},\Xi_{s,s'})$ and
$(\mathcal M_{s',c},\Xi_{s',c})$ are independent, since they are functions
of the noise on $[s,s']$ and on $[s',c]$, respectively; see
Figure~\ref{fig:reference-flow-property}. Averaging first over the second
pair gives
\[
 g_s(z)=\E_{\rm aux}g_{s'}\bigl(\widehat z^{s,z}(s')\bigr).
\]
We apply It\^o's formula to the $C^2$ function $g_{s'}$ along
$\widehat z^{s,z}$ on $[s,s']$. Its martingale part has a
square-integrable integrand, by the moment bound of Step~0, so it has
mean zero, and we obtain
\[
 g_s(z)-g_{s'}(z)=\E_{\rm aux}\int_s^{s'}
   \mathcal L_{\rm ref}(v)g_{s'}\bigl(\widehat z^{s,z}(v)\bigr)\,dv .
\]
The integrand is at most $C(1+\|\widehat z^{s,z}(v)\|^2)$, by the bounds
of Step~1 and the form of $\mathcal L_{\rm ref}$, and
\[
 \E\sup_{v\in[s,s']}\|\widehat z^{s,z}(v)-z\|^2\to0
 \qquad\text{as } s'-s\to0 .
\]
By the continuity of Step~2 and dominated convergence, the difference
quotient $(g_s(z)-g_{s'}(z))/(s'-s)$ tends to
$\mathcal L_{\rm ref}(s_0)g_{s_0}(z)$ when $s'-s\to0$ with either $s=s_0$
or $s'=s_0$ fixed. So both one-sided derivatives exist and equal
$-\mathcal L_{\rm ref}(s_0)g_{s_0}(z)$, which is continuous in $s_0$ by
Step~2. This proves (iii).
 
\emph{Step 4: measurability in the path (iv).} Fix $s'$. For $s\ge s'$,
the function $g_s$ depends on the eigenvalue path only through the
coefficients $c_{ij}$, $0<|i-j|<\ell$, and $\widehat\beta_i$ on $[s',c]$.
On paths with simple spectrum on $[s',c]$, these coefficients are
measurable functions of $(\lambda(w))_{s'\le w\le c}$, with values in
$C([s',c])$. It therefore suffices to show that $g_s(z)$ depends
continuously on the coefficients. Suppose two sets of coefficients differ
by at most $\eta\le1$, uniformly in time, and run both equations with the
same noises. It\^o's formula and Gronwall's inequality give
\[
 \E\|\widehat z^{s,z}(c)-\widetilde z^{s,z}(c)\|^2\le C\eta(1+\|z\|^2),
\]
where for the noise coefficients we used $|\sqrt x-\sqrt y|\le|x-y|^{1/2}$.
Since $g_s(z)=\E_{\rm aux}f(\widehat z^{s,z}(c))$ and $Df$ is bounded,
$g_s(z)$ is a continuous function of the coefficients for fixed $(s,z)$,
and hence a measurable function of the path. Together with the
continuity in $(s,z)$ from Step~2, this gives joint measurability.
 
\emph{Step 5: phase invariance (v).} For phases $|\vartheta_i|=1$ let
$T_\vartheta z=(\vartheta_iz_i)_i$. Conjugation by $T_\vartheta$ rotates
each pair of rotation fields $(X_{ij},Y_{ij})$ into a combination of the
two, so it maps the pair of driving noises $(W_{ij},W'_{ij})$ to a rotated
pair, which is again a pair of independent Brownian motions. It also
leaves the law of the circular noise unchanged. Hence
$(T_\vartheta^{-1}\mathcal M_{s,c}T_\vartheta,T_\vartheta^{-1}\Xi_{s,c})$ has
the same law as $(\mathcal M_{s,c},\Xi_{s,c})$, and, if $f$ is phase
invariant,
\[
 g_s(T_\vartheta z)=\E_{\rm aux}f\bigl(T_\vartheta(\mathcal M_{s,c}z+\Xi_{s,c})\bigr)=g_s(z).
\]
\end{proof}

\subsection{The two generator errors}\label{subsec:s18-the-two-generator-errors}

 We bound the difference between the true and reference eigenvector
flows and complete Theorem~\ref{thm:relative-mixed-vectors}. Expanding
the complex generator leaves a diagonal covariance error and a mixed
second-derivative error, together with the damping correction away from
the bulk. We estimate each term using interval averages and
localization, then integrate over the comparison time.

 Expanding $X_{ij}^2+Y_{ij}^2$ separates a diagonal covariance
error and a cross-Hessian error. This is the complex version of
\cite[Lemma 2.5]{BL}; here the expansion also fixes its normalization.
In real coordinates $z=x+iy$, write
$\mathscr R(H)=\left(\begin{smallmatrix}\Re H&-\Im H\\
\Im H&\Re H\end{smallmatrix}\right)$.
The diagonal error at label $j$ is one quarter of the contraction
of its real Hessian with $\mathscr R(\Delta_j)$, where

\[
 \Delta_j=\sum_{|i-j|\ge\ell}
        \frac{z_iz_i^*-I_r}{\Naux (\lambda_i-\lambda_j)^2}.
\]

The cross contribution of an unordered pair is
\[
 \frac{c_{ij}}{2}
 \{D_{ij}^2g[iz_j,iz_i]-D_{ij}^2g[z_j,z_i]\}.
\]
The first derivatives cancel exactly against the OU drift.
These are all errors in the bulk, where $\widehat\beta_j=\beta_j$.
Globally one must add the OU-floor correction
\begin{equation}\label{eq:t26-floor-error}
 \sum_j(\widehat\beta_j-\beta_j)
     \left(\tfrac12 z_j\cdot\nabla_jg-\tfrac14\Delta_j^{\rm real}g\right).
\end{equation}
It is supported outside the smaller bulk window and is estimated below;
no extra boundary-rotation term occurs.

Dyadic summation by parts using \eqref{eq:relative-interval}
gives, for $|j-k_0|\le w$,

\begin{equation}\label{eq:diagonal-relative-error}
 \norm{\Delta_j}\le C_r\left[
 \frac{\Naux ^{2-\alpha+\epsilon}}{\ell^2}
 +\frac{\Naux }{\ell}(t+w/\Naux )+\log \Naux \right].
\end{equation}

For an interval at distance $2^m\ell$, its covariance-error sum
is at most   $\Naux ^{1-\alpha+\epsilon}$  and the maximum coefficient
is   $C\Naux /(2^m\ell)^2$ ; summing the resulting geometric series
gives the first term. The deterministic covariance differs from
$I$ by   $O(t+w/\Naux )$  near the tested energy, giving the second.
Its linear variation across more distant shells sums to   $O(\log \Naux )$.
Equation~\eqref{eq:relative-distant} bounds the exterior energy
region with no polynomial norm loss.

For a bounded $C^2$ test $f$ of a fixed array of central labels,
the reference backward function has the representation
$g_s(z)=\E_{\rm aux}f((\mathcal M_{s,t}z+\Xi_{s,t})_{\rm tested})$.
Conditional on the auxiliary rotations, $\Xi$ is Gaussian with
covariance $I-\mathcal M\mathcal M^*$.
Let $T_j$ map the $j$th input vector to the tested array.
Then
$\nabla_j^2g_s=\E_{\rm aux}T_j^*(\nabla^2f)T_j$, and
$\sum_j\norm{T_j}_{\HS}^2$ is at most the fixed tested dimension,
since $\norm{\mathcal M}\le1$.
Hence
\[
 \sum_j\norm{\nabla_j^2g_s}\le C_r\norm f_{C^2}.
\]
For $j\notin J$, the sums of squared Hilbert--Schmidt norms of $T_j$
are bounded by a fixed multiple of the row tail in
\eqref{eq:reference-finite-speed}. The corresponding sum of first
norms is at most   $C\sqrt \Naux $  times the square root of this tail.
On $\mathcal R$, all long coefficients and OU rates are polynomially
bounded; also   $\|z_j\|\le\sqrt \Naux \|V\|$  deterministically.
Cauchy--Schwarz in time and the localization lemma therefore bound the
integrated exterior Hessian terms, all mixed terms with at least one
exterior index, and \eqref{eq:t26-floor-error} by   $\Naux ^{-D}$, after
choosing the lemma's failure exponent larger than the fixed frame exponent.
Inside $J$ the floor vanishes. The integrated diagonal error is at most
$C_ft$ times the right side of \eqref{eq:diagonal-relative-error}, plus
  $\Naux ^{-D}$.

For the cross error, fix the auxiliary propagator and expand a
tested Hessian entry. Its coefficient is a sum of terms of the form
$\sum_{|i-j|\ge\ell}\Theta_{ij}z_i^a z_j^b$ or its conjugated
version, with
$|\Theta_{ij}|\le c_{ij}|\mu_j^a\mu_i^b|$.
The stopping event is invariant under the independent phases of
the actual eigenvectors. Phase averaging its squared absolute value
leaves only the index pairings $(i,j)=(k,l)$ or $(l,k)$.
On the local event   $\norm{z_i}^2\le C_r\Naux ^{1-\alpha}$, so

\[
 \E_{\rm phases}\left|\sum_{i\ne j}\Theta_{ij}z_i^az_j^b\right|^2
 \le C_r\Naux ^{2-2\alpha}\sum_{i\ne j}
                  (|\Theta_{ij}|^2+|\Theta_{ij}\Theta_{ji}|).
\]

The conjugated version has fewer pairings and the same bound.
Use $2|xy|\le|x|^2+|y|^2$,
  $c_{ij}\le C\Naux /\ell^2$  and
$\sum_i|\mu_i^a|^2\le1$.
Taking square roots and integrating gives

\begin{equation}\label{eq:cross-relative-error}
 \int_{c-t}^c|\E\,\mathcal O_sg_s|\,ds
                 \le C_f\Naux ^{2-\alpha+\epsilon}t/\ell^2+\Naux ^{-D}.
\end{equation}

The auxiliary randomness is independent of the actual phases and
stopping event, so this pairing is legitimate. Hessian dependence
on the actual vectors is removed first by its uniform norm bound,
as in \cite[Lemma 2.14]{BL}.

\begin{proof}[Proof of Theorem~\ref{thm:relative-mixed-vectors}]
 We compare the true angular flow with the reference
flow of the previous paragraphs, which keeps the short rotations and
replaces the long ones by Ornstein--Uhlenbeck damping. The reference
flow is linear and contracts fast, so its terminal law is close to the
product circular Gaussian law. A Duhamel formula expresses the
difference between the two flows as a time integral of the difference of
their generators. This difference involves only the long rotations,
because the short rotations are common to both flows, as in
Figure~\ref{fig:rm-duhamel}. The proof has five steps. We fix the good
eigenvalue paths (Step~1), justify the Duhamel identity (Step~2), bound
the generator difference (Step~3), add up the errors (Step~4), and pass
to the unconditional and the point-process statements (Step~5).

\begin{figure}[H]
\centering
\begin{tikzpicture}[x=1cm,y=1cm]
 \draw[->] (0,0)--(9.4,0) node[below,font=\small] {time};
 \draw[thick] (.5,-.1)--(.5,.1) node[below=4pt,font=\small] {$a$};
 \draw[thick] (8.6,-.1)--(8.6,.1) node[below=4pt,font=\small] {$c$};
 \draw[thick] (4.6,-.1)--(4.6,.1) node[below=4pt,font=\small] {$s$};
 \draw[very thick,black] plot[smooth] coordinates
   {(.5,1.0) (1.4,1.35) (2.3,1.05) (3.2,1.45) (4.0,1.2) (4.6,1.4)};
 \draw[very thick,blue!70!black,dashed] plot[smooth] coordinates
   {(4.6,1.4) (5.5,1.1) (6.5,1.3) (7.5,1.0) (8.6,1.15)};
 \node[font=\small,anchor=south] at (2.5,1.6) {true flow on $[a,s]$};
 \node[font=\small,anchor=south,blue!70!black] at (6.6,1.45) {reference flow on $[s,c]$};
 \fill (4.6,1.4) circle (2pt);
 \node[font=\small,align=center] at (4.6,2.45)
   {switch at time $s$:\\ error rate $(\mathcal L_{\rm true}-\mathcal L_{\rm ref})g_s$};
 \node[font=\small,align=center] at (4.6,-1.15)
   {short rotations $|i-j|<\ell$ are common to both flows and cancel};
\end{tikzpicture}
\caption{The Duhamel comparison in the proof of
Theorem~\ref{thm:relative-mixed-vectors}. For each switching time $s$ we
run the true angular flow (solid) up to $s$ and the reference flow
(dashed) from $s$ to $c$. The function $g_s$ is the expected terminal
value of the reference flow started at time $s$. Moving the switching
time from $c$ back to $a$ interpolates between the true and the reference
terminal laws. The rate of change is the difference of the two
generators applied to $g_s$. Since both flows contain the same short
rotations, this difference involves only the long rotations of the true
flow and the Ornstein--Uhlenbeck part of the reference flow. Its size is
controlled by the contraction bounds on the derivatives of $g_s$.}
\label{fig:rm-duhamel}
\end{figure}
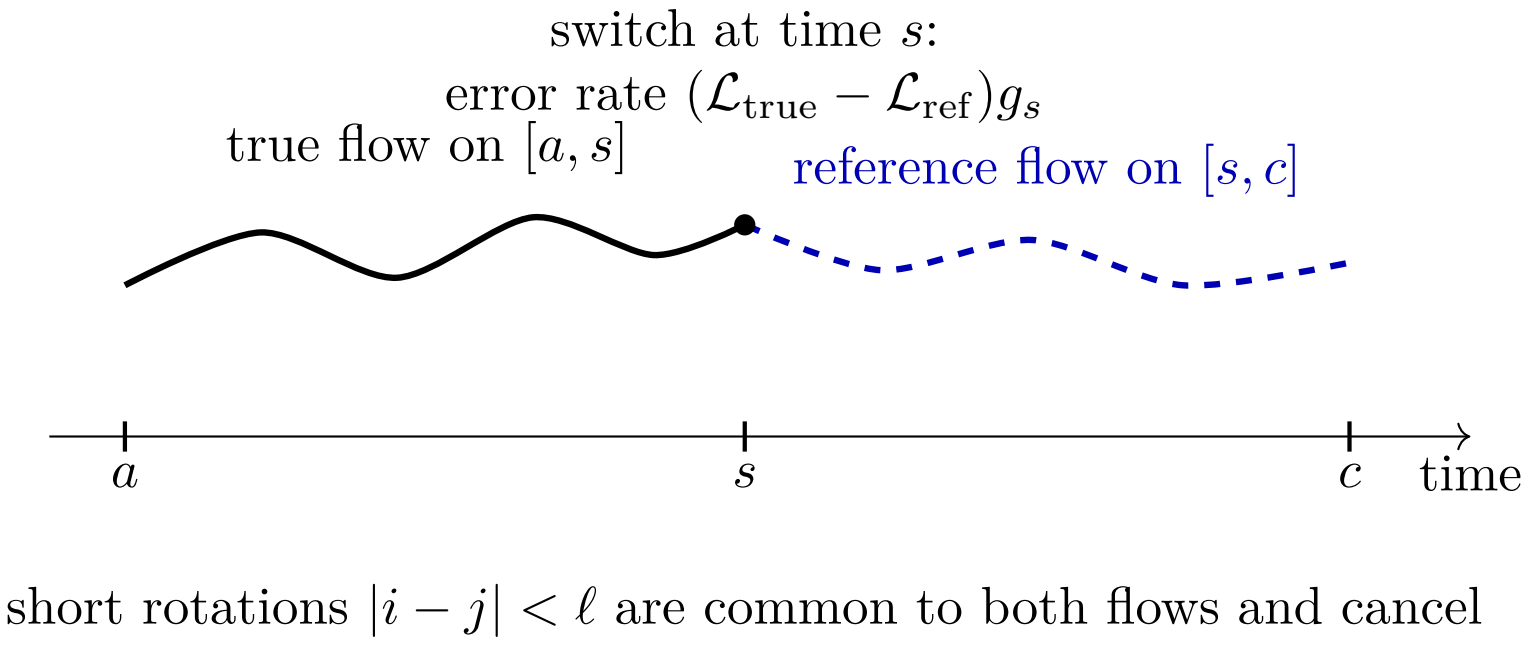

\emph{Step 1: the good eigenvalue paths.} We first specify the
conditional events and how errors are removed. Let $\mathcal A$ be the
event that \eqref{eq:relative-interval} and \eqref{eq:relative-distant}
hold throughout $[a,c]$. This event depends on the actual flow and is
invariant under the phases of the eigenvectors. By the relative coarse
law and its time mesh, its failure probability is smaller than any
power of   $\Naux $. We want to work conditionally on the eigenvalue path, so we
disintegrate. For any desired $D$, we apply Markov's inequality to
\[
 \Prob\bigl(\mathcal A^c\mid\sigma(\lambda(s):a\le s\le c)\bigr).
\]
Taking first the unconditional exponent $2D$, this excludes eigenvalue
paths of probability at most   $\Naux ^{-D}$, and on the remaining paths the
conditional failure of $\mathcal A$ is at most   $\Naux ^{-D}$. We intersect
these paths with $\mathcal R$ and with the good paths of
Lemma~\ref{lem:t26-relative-localization}. All these are measurable sets
of full eigenvalue paths. We never condition the angular evolution on
$\mathcal A$ itself.

\emph{Step 2: the Duhamel identity.} We fix such a path. We use the
backward function $g_s$ of the reference flow, that is, the expected
terminal value of the reference flow started at time $s$, and It\^o's
formula for the true angular flow. Since $g_s$
depends on the eigenvalue path after time $s$, we use It\^o's formula in
the form of Lemma~\ref{lem:t26-angular-martingale}, at restart times.
Test functions are $C^2$, with the derivative bounds of
Theorem~\ref{thm:relative-mixed-vectors}. Both sides of
\eqref{eq:relative-mixed-vectors} are unchanged if $f$ is replaced by its
average over independent uniform phases of its $k$ arguments, since the
eigenvector phases are uniform and independent of the eigenvalue path and
the moduli, and circular Gaussian vectors are phase invariant; the
average has the same bounds, so we assume $f$ phase invariant. Then
Lemma~\ref{lem:t26-reference-backward}(i)--(v) applies to $g$. Fix $n$,
let $\tau=\tau_n$ be the gap stop with $s=a$, $\mathcal Q$ an event of the
eigenvalue path, $a=s_0<\dots<s_p=c$, and $\widehat g_v=g_{v\wedge\tau}$.
As $\tau$ is determined by the eigenvalue path, by (iv) with
$s'=s_{k-1}$ the test $\one_{\mathcal Q}\one_{\{\tau\ge s_{k-1}\}}\widehat g_{s_k}$ is
admissible in Lemma~\ref{lem:t26-angular-martingale} with start
$s_{k-1}$, whose gap stop is $\tau$ on $\{\tau\ge s_{k-1}\}$. Hence
\[
 \E\bigl[\one_{\mathcal Q}\bigl(\widehat g_{s_k}(z(s_k\wedge\tau))
     -\widehat g_{s_k}(z(s_{k-1}\wedge\tau))\bigr)\bigr]
 =\E\Bigl[\one_{\mathcal Q}\int_{s_{k-1}\wedge\tau}^{s_k\wedge\tau}
      \mathcal L_{\rm true}(v)\widehat g_{s_k}(z(v))\,dv\Bigr],
\]
both sides vanishing on $\{\tau<s_{k-1}\}$, while (iii) gives
\[
 \E\bigl[\one_{\mathcal Q}(\widehat g_{s_k}-\widehat g_{s_{k-1}})(z(s_{k-1}\wedge\tau))\bigr]
 =-\E\Bigl[\one_{\mathcal Q}\int_{s_{k-1}\wedge\tau}^{s_k\wedge\tau}
      \mathcal L_{\rm ref}(v)g_v(z(s_{k-1}\wedge\tau))\,dv\Bigr].
\]
Before $\tau$ all coefficients of both generators are bounded. Summing
over $k$ and letting the mesh tend to zero by (i), (ii) and dominated
convergence,
\begin{equation}\label{eq:t26-duhamel-restart}
 \E\bigl[\one_{\mathcal Q}\,g_{c\wedge\tau}(z(c\wedge\tau))\bigr]-\E\bigl[\one_{\mathcal Q}\,g_a(z(a))\bigr]
 =\E\Bigl[\one_{\mathcal Q}\int_a^{c\wedge\tau}
   \bigl(\mathcal L_{\rm true}(v)-\mathcal L_{\rm ref}(v)\bigr)g_v(z(v))\,dv\Bigr].
\end{equation}
As $n\to\infty$, $\tau_n\to\infty$ and $g_c=f$, so the left side tends to
$\E[\one_{\mathcal Q}f(z(c))]-\E[\one_{\mathcal Q}g_a(z(a))]$. If
$\mathcal Q\subset\mathcal R$, the right side converges by dominated
convergence: the short rotations cancel in
$\mathcal L_{\rm true}-\mathcal L_{\rm ref}$, and the long couplings and
$\widehat\beta$ are bounded on $\mathcal R$ uniformly in $n$ by the first
line of \eqref{eq:t26-long-gaps} and the bound on $\beta_i$. This applies
to the events $B\subset\mathcal G$ of
Theorem~\ref{thm:relative-mixed-vectors}, since the good paths of Step~1
were intersected with the eigenvalue-path event $\mathcal R$; the
resulting integrated identity is the Duhamel identity used below. The generators of the short rotations
cancel exactly, because both flows contain them. What remains is the
difference between the long coefficients of the true flow and the
Ornstein--Uhlenbeck part of the reference flow. Three bounds make this
difference integrable. The long coefficients and the Ornstein--Uhlenbeck
floor are polynomially bounded on $\mathcal R$. The true vectors are
bounded by   $\sqrt \Naux \|V\|$. And the gradient and the Hessian of $g_s$ obey
the contraction bounds above. Thus the Duhamel identity is integrable,
and its contribution on $\mathcal A^c$ is   $O(\Naux ^{-D+C(A)})$.

The short coefficients $c_{ij}$ blow up if two eigenvalues collide, so
we check that this causes no problem. We may first stop before a gap
tends to zero, and then pass to the limit. The common short generator
cancels before this passage, and the bound we obtain contains no inverse
short gap. In fact the eigenvalues never collide. Simplicity and
noncollision of complex Dyson Brownian motion also follow from the Doob
transform of killed independent Brownian motions by the harmonic
Vandermonde $\Delta$, as in \cite[Section 3.3, (3.10) and Remark 1]{KT}.
Indeed, since $\Delta$ is harmonic and
$\partial_i\Delta/\Delta=\sum_{j\ne i}(\lambda_i-\lambda_j)^{-1}$,
\[
 \Delta^{-1}\Bigl(\frac{1}{2}\sum_i\partial_i^2\Bigr)(\Delta f)
 =\frac{1}{2}\sum_i\partial_i^2f+
 \sum_i\sum_{j\ne i}(\lambda_i-\lambda_j)^{-1}\partial_i f,
\]
and \eqref{eq:t26-vector-dbm} has exactly   $1/\Naux $  times this generator. The
square $\Delta^2$ occurs in the eigenvalue density, not as the harmonic
function. We use these facts only on the compact interval $[a,c]$ of
positive times.

\emph{Step 3: the generator difference.} On $\mathcal A$, the diagonal
calculation and the phase-pairing calculation give
\eqref{eq:diagonal-relative-error} and \eqref{eq:cross-relative-error},
conditionally on the path. One point needs care. The Hessian of $g_s$ is
evaluated at the true vectors, so it depends on their phases, and we must
check that averaging over the phases is legitimate. We first write the
generator error as a sum over the finitely many tested entries of the
Hessian. We bound each entry by $\|\nabla^2f\|_\infty$, and we apply the
Cauchy--Schwarz inequality to its coefficient. Only then do we average
the squared coefficient over the independent uniform phases of the
columns, keeping $\mathcal A$. This is exactly the pairing displayed
above. The independent reference rotations and the Ornstein--Uhlenbeck
noise are averaged last.

\emph{Step 4: adding up the errors.} The initial array has deterministic
norm at most   $\sqrt \Naux \|V\|_{\rm HS}$. By \eqref{eq:t26-reference-covariance}
and the contraction estimate, the terminal law of the reference flow is
coupled to the product circular Gaussian law with bounded-Lipschitz
error

\[
 \Naux ^{C(A,r)}e^{-b_0\Naux t/(2\ell)}.
\]

This holds on every conditioned spectral path. The total Duhamel error
tends to zero at a fixed power of   $\Naux $, because

\[
\begin{aligned}
 2\kappa-\alpha+\epsilon-\theta&=-0.0067,\\
 \kappa-\tfrac32\theta+2\epsilon&=-0.0002,\qquad
 \Naux t/\ell=\Naux ^{\kappa-\theta}=\Naux ^{0.0006}.
\end{aligned}
\]

These three exponents control, respectively, the fluctuation and cross
terms, the variation   $t(\Naux /\ell)(w/\Naux )$  of the covariance, and the
equilibration. In particular, since   $\Naux t/\ell=\Naux ^{0.0006}$, the
equilibration error   $\Naux ^{C(A,r)}e^{-b_0\Naux ^{0.0006}/2}$  is smaller than any
power of   $\Naux $. The smaller terms   $t(\Naux /\ell)t$ and $t\log \Naux $  also vanish.
The enlarged localization window costs only the elementary displacement
$\sqrt t$ from Weyl's inequality, and it still leaves positive slack in
the exponents. We choose $D$ after the fixed exponent of the deformation
and the dimension of the test. This proves the conditional
bounded-$C^2$ comparison, uniformly over all terminal microscopic labels
and all normalized frames.

\emph{Step 5: independence and the decorated sine process.} We multiply
the conditional comparison by any bounded spectral test, and we
integrate over the good paths. This proves joint independence from that
test, with the same polynomial error. No statement about unbounded
angular moments follows from this argument for bounded tests. Those
moments are treated separately in Theorem~\ref{thm:fresh-palm}.

Combining this conditional comparison with
Theorem~\ref{thm:external-sine} at the original fixed variance $c$ gives
the decorated sine process used in the collar, that is, the sine process
whose points carry independent Gaussian vectors. On a compact spectral
interval we first retain at most $k$ points, and we use their joint
bounded spectral tests. Because the comparison is uniform conditionally
on the path, the labels may be selected from the eigenvalue path. The
exponential count bound then removes the restriction to $k$ points,
retaining independent Gaussian decorations. 

In detail, fix $\varepsilon$, $\phi$ as in \eqref{eq:decorated-sine}
and an integer $k$. Let $x_j=\Naux\rho_\Naux(\lambda_j(c)-E)$ with
$|\rho_\Naux/\rho-1|\le\delta\le1/2$, and first let $V^*\Gamma V=I_r$. Let
$N_R=\#\{j:|x_j|\le R\}$, and for $1\le m\le k$ and a rank $j_0$ let
$B_{j_0,m}\subset\mathcal G$ be the (disjoint) events of the eigenvalue
path that these points are exactly $x_{j_0+1},\dots,x_{j_0+m}$; on
$\{N_R=0\}$ both exponentials below equal one. On
$B_{j_0,m}$ the sum $\sum_j\phi(x_j,z_j(c))$ involves only the $m$ selected
vectors, and $\exp\{-\sum_{l\le m}\phi(x_{j_0+l},\cdot)\}$ is a test function
whose size and first two derivatives are bounded in terms of $k$ and
$L_\phi$, with coefficients determined by the path; this is admissible
because the comparison of Step~4 holds on every conditioned spectral
path, uniformly over the terminal labels. It gives
\[
 \E\Bigl[e^{-\sum_j\phi(x_j,z_j(c))}\one_{B_{j_0,m}}\Bigr]
 =\E\Bigl[e^{-\sum_j\psi(x_j)}\one_{B_{j_0,m}}\Bigr]
  +O_{k,L_\phi}\bigl(\Naux^{-\chi}\Prob(B_{j_0,m})\bigr),\qquad
 \psi(x)=-\log\E e^{-\phi(x,g_1)} .
\]
Summing over $j_0$ and $m\le k$, and bounding the integrands by one on
$\{N_R>k\}$ and $\mathcal G^c$,
\[
 \Bigl|\E e^{-\sum_j\phi(x_j,z_j(c))}-\E e^{-\sum_j\psi(x_j)}\Bigr|
 \le C_{k,L_\phi}\Naux^{-\chi}+2\Prob(N_R>k)+2\Naux^{-D}.
\]
Since $H_c/\sqrt c$ is $F/\sqrt c$ plus a normalized GUE matrix,
Theorem~\ref{thm:external-count} with one block, $d=\Naux$,
$A_0=F/\sqrt c$ and an interval of length
$2R/(\sqrt c\,\Naux\rho_\Naux)\le4R/(\sqrt c\,\Naux\rho)$ gives
\[
 \Prob(N_R>k)\le\frac{1}{(k+1)!}\E\prod_{i=0}^{k}(N_R-i)
 \le\frac{(4CR/(\sqrt c\,\rho))^{k+1}}{(k+1)!},
\]
uniformly, because $\rho\ge v_0/\pi$ by \eqref{eq:fresh-regularity}.
As $\psi$ is Borel, $0\le\psi\le L_\phi$ and $\psi=0$ off $[-R,R]$,
$\E e^{-\sum_j\psi(x_j)}$ is the Fredholm determinant of $1-e^{-\psi}$
times the correlation kernel rescaled by $\Naux\rho_\Naux$ on $[-R,R]$
(the conjugation in \eqref{eq:t26-complex-sine-kernel} does not change
it). Theorem~\ref{thm:external-sine} with variance $c$ and
$v_*=v_0$ gives the kernel rescaled by $\Naux\rho$ within
$O(\Naux^{-1/2})$ of the sine kernel on $[-2R,2R]^2$. The kernel rescaled
by $\Naux\rho_\Naux$ is $\rho/\rho_\Naux$ times that kernel at
$(\rho/\rho_\Naux)\xi,(\rho/\rho_\Naux)\zeta$; by uniform continuity of
the sine kernel on compacts it lies within
$\omega_R(\delta)+O(\Naux^{-1/2})$ of the sine kernel on $[-R,R]^2$, where
$\omega_R(\delta)\to0$ as $\delta\to0$. By the Fredholm series and Hadamard's inequality these determinants are
Lipschitz in the kernel and in $1-e^{-\psi}$, uniformly for kernels
bounded on $[-R,R]^2$, so $\E e^{-\sum_j\psi(x_j)}$ is within
$C_R(\omega_R(\delta)+\Naux^{-1/2})$ of $\E e^{-\sum_j\psi(y_j)}$, the
right side of \eqref{eq:decorated-sine} after conditioning on the
points. If only $\norm{V^*\Gamma V-I_r}\le\delta$, apply this to
$V'=V\mathsf T^{-1}$, $\mathsf T=(V^*\Gamma V)^{1/2}$, and
$\phi_{\mathsf T}(x,z)=\phi(x,\mathsf Tz)$, whose $C^2$ norms in $z$ are at
most $3L_\phi$; since $z_j=\mathsf Tz_j'$, $\norm{\mathsf T-I}\le\delta$
and $|e^{-\phi(x,\mathsf Tg)}-e^{-\phi(x,g)}|\le L_\phi\delta|g|$, the
Laplace functionals of the decorated sine process for $\phi$ and
$\phi_{\mathsf T}$ differ by at most $C_RL_\phi\delta\sqrt r$. Now choose $k$ with
$2(4CR/(\sqrt c\,\rho))^{k+1}/(k+1)!\le\varepsilon/4$, then $\delta$ with
$C_R(\omega_R(\delta)+L_\phi\delta\sqrt r)\le\varepsilon/4$, and finally
$\Naux_0$ such that the remaining terms, with $3L_\phi$ in place of
$L_\phi$, are at most $\varepsilon/2$ for $\Naux\ge\Naux_0$. This proves
\eqref{eq:decorated-sine}.
\end{proof}

\section{Physical normalization and Gaussian reconstruction}
\label{sec:physical-reconstruction}
Parts~II and~III need one-mark statistics of the ring, that is, sums over
eigenvalues in a window of a function of their profiles. This section reduces
such a statistic to a normalized Gaussian reconstruction, by removing the
exterior weight without losing its exact normalization.
Proposition~\ref{prop:t22-exact-palm} gives the exact formulas, indexed
and unranked, with their normalization. Theorem~\ref{thm:t22-normalizer}
shows that the exterior weight converges in total variation to the original
exterior law times the energy-window weight.
Theorem~\ref{thm:t22-Gaussian-ratios} shows that bounded ratios of quadratic
forms of the reconstructed direction may be replaced by Gaussian ratios.
Theorem~\ref{thm:t22-one-mark-Gaussian} combines these into a formula for one
mark at the energy $e(u)=E+u/(N\rho(E))$.
In Subsection~\ref{subsec:s19-exact-cancellation-with-every-remaining-} we identify every factor in the energy-window
measure and cancel the coarea factors. In Subsection~\ref{subsec:s19-a-pointwise-density-lemma-for-the-needed} we compare the
density of the fresh block at one energy with the deterministic density
(Lemma~\ref{lem:t22-point-density}). In Subsection~\ref{subsec:s19-the-exact-exterior-measure-in-total-vari} we prove the
convergence of the exterior measure in total variation. In
Subsection~\ref{subsec:s19-ratios-of-arbitrary-positive-reconstruct} we compare ratios of arbitrary positive clocks. In
Subsection~\ref{subsec:s19-the-physical-one-mark-comparison} we assemble the one-mark formula.

\subsection{Exact cancellation, with every remaining factor retained}\label{subsec:s19-exact-cancellation-with-every-remaining-}
 We first identify every factor in the physical energy-window measure.
Proposition~\ref{prop:t22-exact-palm} gives the indexed and unranked
formulas and their exact normalization. We expose a fresh pair, write
the reconstruction and zero-eigenvector density, and then change to
the microscopic energy variable without discarding the determinant
or rank constraint.

 Fix an adjacent pair $S$, put   $\Naux =2W$, $s=\min(b,1-2b)$ and   $t=s/\Naux $.
Split an independent   $\operatorname{GUE}_{\Naux }(t)$  matrix $G$ from $H$
on $S$; the remaining supported Gaussian matrix is $H^{(S)}$.
For example, its two diagonal-block variances are $(c-s/2)/W$ and
its internal bond variance is $(b-s/2)/W$. All other original variances
are unchanged. This is the one-pair marginal of \eqref{eq:simultaneous-splitting}. Let $R=S^c$. At a real integration energy $e$ outside
the finite spectrum of $H_{RR}$, set
\[
 K(e)=H^{(S)}_{SS}-e-H_{SR}(H_{RR}-e)^{-1}H_{RS}.
\]
\[
 L(e)=\begin{pmatrix}I_S\\-(H_{RR}-e)^{-1}H_{RS}\end{pmatrix},
 \qquad D=L^*L=-K'(e)\succ0.
\]
Here and below the ordering of the rows of $L$ is restored to the
original block order. The real pole set has Lebesgue measure zero;
no expectation of an unweighted real inverse is used.

Let   $v\in\mathbb {CP}^{\Naux -1}$  be represented by a unit vector, $Q_v=I-vv^*$,
and let $U_v$ be an orthonormal frame for $v^\perp$. Define

\[
 a=v^*Kv,\quad b_2=v^*K^2v,\quad
 Z_\ell(v;e)=\mathbb E_{G_{\Naux -1}(t)}
  \left[\det(U_v^*KU_v+G_{\Naux -1})^2
    \mathbf1_{\{N_-(U_v^*KU_v+G_{\Naux -1})=\ell\}}\right].
\]

An unattainable sector has mass zero. With normalized projective Haar
measure   $d\nu_\Naux $, put

\begin{equation}\label{eq:t22-palm-angular}
 d\mathcal A_{\ell,e}(v)=
 \frac{t^{-(\Naux -1)}}{(\Naux -1)!\sqrt{2\pi t}}
 e^{-b_2/t+a^2/(2t)}Z_\ell(v;e)\,d\nu_\Naux (v),\qquad
 \mathcal A_e=\sum_{\ell=0}^{\Naux -1}\mathcal A_{\ell,e}.
\end{equation}

Write $Z(e)=\mathcal A_e(1)$, and
\[
 p_k(v,e)=\frac{\|P_kLv\|^2}{v^*Dv},\qquad
 q_S(v,e)=(v^*Dv)^{-1}.
\]
The exterior expectation below is the unchanged Gaussian prior of
$H^{(S)}$ and its connecting bonds.

\begin{proposition}[Complete fixed-energy and indexed Palm formula]
\label{prop:t22-exact-palm}
For every bounded measurable $f(e,p)$ with compact energy support,
\begin{align}
 \mathbb E\sum_\alpha q_{\alpha S} f(\lambda_\alpha,p_\alpha)
 &=\mathbb E_{\rm ext}\int_{\mathbb R}
           \mathcal A_e[f(e,p(v,e))]\,de,\label{eq:t22-mass-palm}\\
 \mathbb E f(\lambda_\alpha,p_\alpha)
 &=\mathbb E_{\rm ext}\int_{\mathbb R}
  \mathcal A_{\alpha-1-N_-(H_{RR}-e),e}
        [(v^*Dv)f(e,p(v,e))]\,de.\label{eq:t22-index-palm}
\end{align}
Thus pair mass cancels exactly the crossing derivative, but leaves the
regression, compressed determinant, sector and exterior measure in
\eqref{eq:t22-palm-angular}. Its full unranked mass satisfies

\begin{equation}\label{eq:t22-Z-upper}
 0\le Z(e)\le C_b\Naux 
\end{equation}

for every exterior and real $e$ (with a continuous density version).
\end{proposition}

The density records the original eigenvalue energy and ordered rank. Summation over ranks will later yield the unranked normalizer, with no conditional-probability division.

\begin{proof}
At a zero eigenvalue of $G+K(e)$ with eigenvector $v$, the corresponding
eigenvalue branch has derivative $-v^*Dv$. Schur inertia gives the
full ordered rank as $1+N_-(H_{RR}-e)+\ell$. Kac--Rice in $e$ therefore
weights the zero-eigenvector Palm measure by $v^*Dv$. The reconstructed
unit vector is $Lv/\|Lv\|$ and its pair mass is $(v^*Dv)^{-1}$, proving
the claimed cancellation and the index in the second formula.

For completeness, zero-eigenvector coarea on   $\mathbb{CP}^{\Naux -1}$ has
volume $\pi^{\Naux -1}/(\Naux -1)!$ and Jacobian
$\det(U_v^*KU_v+G_{\Naux -1})^2$. The real constraint $v^*Gv=-a$ has
density $(2\pi t)^{-1/2}e^{-a^2/(2t)}$; its   $\Naux -1$  complex transverse
constraints have joint density   $(\pi t)^{-(\Naux -1)}
 e^{-(b_2-a^2)/t}$. These independent GUE coordinates leave precisely
the displayed complementary GUE. Multiplying proves the complete
constant and exponent in \eqref{eq:t22-palm-angular}. Equivalently,
this is the local form of the coarea calculation in Section~\ref{sec:coarea}; no sector has been averaged separately. First apply the
identities to nonnegative functions, then subtract their positive and
negative parts. The integrated first identity at $f=1$ has value $2W$;
the second has value one. Finally $Z(e)$ is the one-point density at
zero of the deterministic deformation $K(e)+G$. The fresh-GUE Wegner
bound, uniform in the deformation, gives \eqref{eq:t22-Z-upper}.
The density is continuous because the independent scalar part of $G$
convolves the trace-free eigenvalue laws with a Gaussian.
\end{proof}

The same identity holds after a deterministic scalar shift of the pair.

\begin{corollary}[Mass Palm formula under a scalar pair shift]
\label{cor:t22-shifted-palm}
Let $\varsigma\in\mathbb R$ be deterministic. Let $\lambda^\varsigma_\alpha$
and $\psi^\varsigma_\alpha$ be the eigenpairs of $H+\varsigma P_S$, with
pair masses $q^\varsigma_{\alpha S}$ and block-mass vectors
$p^\varsigma_\alpha$. Let $\mathcal A^\varsigma_e$ be \eqref{eq:t22-palm-angular} with $K(e)$
replaced by $K(e)+\varsigma I_S$ in $a$, $b_2$ and $Z_\ell$. Then for every bounded measurable $f(e,p)$ with
compact energy support,
\[
 \mathbb E\sum_\alpha q^\varsigma_{\alpha S}
        f(\lambda^\varsigma_\alpha,p^\varsigma_\alpha)
 =\mathbb E_{\rm ext}\int_{\mathbb R}
           \mathcal A^\varsigma_e[f(e,p(v,e))]\,de,
\]
with the same $L(e)$, $D(e)$ and $p(v,e)$ as for $\varsigma=0$. The mass
$Z^\varsigma(e)=\mathcal A^\varsigma_e(1)$ satisfies \eqref{eq:t22-Z-upper}
for every exterior, every real $e$ and every $\varsigma$. In particular,
for every measurable $\varphi\ge0$,
\[
 \mathbb E\sum_\alpha q^\varsigma_{\alpha S}\varphi(\lambda^\varsigma_\alpha)
 \le C_b\Naux\int_{\mathbb R}\varphi(e)\,de .
\]
\end{corollary}

\begin{proof}
In the decomposition $S\oplus R$ the matrix $H+\varsigma P_S$ has the same
blocks $H_{RR}$, $H_{RS}$ and fresh $G$, and pair block
$H^{(S)}_{SS}+\varsigma I_S+G$. It is obtained by replacing the exterior
$H^{(S)}$ by $H^{(S)}+\varsigma P_S$, a deterministic change that keeps
$G$ independent of the exterior. The proof of
Proposition~\ref{prop:t22-exact-palm} is carried out at a fixed exterior,
and its steps (Schur inertia, the crossing derivative, the reconstruction
and the zero-eigenvector coarea for $G$, i.e.\ the fixed-$\omega$ form of
Theorem~\ref{thm:coarea} with Lemma~\ref{lem:eigenpair-jacobian}) use the
exterior only through $K(e)$, $L(e)$ and $D(e)=-K'(e)$, which become
$K(e)+\varsigma I_S$, $L(e)$ and $D(e)$. Applying the fixed-exterior
identity to the shifted exterior and integrating over the unchanged
exterior law proves the formula. The Wegner bound behind
\eqref{eq:t22-Z-upper} is uniform in the deformation. For the last
inequality take $f(e,p)=\varphi(e)$, truncated and restricted to a
compact energy set, use $Z^\varsigma(e)\le C_b\Naux$, and remove the
truncation by monotone convergence.
\end{proof}

Let $w\ge0$ be a fixed continuous compactly supported function with
$\int w=1$, and set   $w_W(e)=w(N\rho(E)(e-E))$. Define the finite positive
\emph{actual exterior window measure}
\begin{equation}\label{eq:t22-exterior-window}
 d\mathfrak m_W(e,\mathcal E)=\frac{n}{2}w_W(e)Z(e)\,de\,d\mathbb P_{\rm ext}(\mathcal E).
\end{equation}
Its mass is exactly $\mathbb E\sum_\alpha w_W(\lambda_\alpha)$ and tends
to one by Theorem~\ref{thm:critical-density}. It is not replaced by
the exterior prior. On a deterministic-energy exterior event $B(e)$,
\begin{equation}\label{eq:t22-exterior-admission}
 \mathfrak m_W(B)\le C_{b,E,w}\sup_{e\in\operatorname{supp}w_W}
                         \mathbb P_{\rm ext}(B(e)).
\end{equation}
Indeed   $n\Naux /N=2$  and the energy change of variables pays the apparent
factor   $n\Naux $. In particular any fixed-energy original-prior exceptional
probability tending to zero is paid under this actual window measure.

\subsection{A pointwise density lemma for the needed deformation class}\label{subsec:s19-a-pointwise-density-lemma-for-the-needed}
  We need the density of the fresh block at one energy, not an average
over an interval. Lemma~\ref{lem:t22-point-density} compares it uniformly
with the deterministic density for the required large-deformation
class. We apply the fixed-variance Gaussian kernel result, then deduce
a deterministic lower bound for the full unranked normalizer.

Let $F=F^*\in\mathbb C^{\Naux \times \Naux }$, $\|F\|\le \Naux ^A$, and let $G_c$
have GUE entry variance   $c/\Naux $, for a fixed $0<c<1$. Write
$\rho_{F,c}(x)$ for the normalized expected eigenvalue density of $F+G_c$,
and

\[
 m_{F,s}(z)=\Naux ^{-1}\operatorname{Tr}(F-z-sm_{F,s}(z))^{-1},
 \qquad r_{F,c}(x)=\pi^{-1}\Im m_{F,c}(x).
\]

Assume \eqref{eq:fresh-regularity}: $\Im m_{F,s}$ is bounded
above and below by fixed positive constants on a fixed real neighborhood
of $x$, for $0\le\Im z\le1$ and $s\in[c-s_0,c]$, with fixed $s_0>0$.
All constants and the uniformity below retain $A,c,s_0$ and those
regularity constants. The statement is valid for polynomially unbounded
$F$; a compactly supported limiting deformation measure is not assumed.

\begin{lemma}[Uniform pointwise fresh density]\label{lem:t22-point-density}
Assume the hypotheses just stated above, i.e., $F$ is Hermitian with
$\|F\|\le\Naux^A$, the Gaussian part $G_c$ has GUE entry variance
$c/\Naux$, and the Dyson solution satisfies the regularity condition
\eqref{eq:fresh-regularity} on a fixed real neighbourhood of $x$. Then
there is a deterministic sequence $\varepsilon_\Naux\to0$, depending only
on $A$, $c$, $s_0$ and the regularity constants, such that
\begin{equation}\label{eq:t22-point-density}
 |\rho_{F,c}(x)-r_{F,c}(x)|\le\varepsilon_\Naux .
\end{equation}
In words, the expected eigenvalue density of $F+G_c$ at the single point
$x$ is close to the deterministic density $r_{F,c}(x)$ given by the
free convolution, uniformly over the whole regular class.
\end{lemma}

This is a statement about one point of the density, not about an average
over an interval, and it concerns an independent full Gaussian matrix
with a fixed matrix added to it. It says nothing yet about the ring.
Passing from such a matrix to the ring still requires the normalization of
the exterior weight, which is carried out in
Theorem~\ref{thm:t22-normalizer}.

\begin{proof}The estimate is the diagonal part of the sine
kernel theorem, Theorem~\ref{thm:external-sine}. The value of the
correlation kernel on the diagonal is the expected eigenvalue density.
Before using it, we check that this density is a genuine continuous
function, so that its value at a single point makes sense (Step~1). We
then shift the energy to zero, because the theorem is stated for
energies in a fixed compact interval while $x$ need not lie in such an
interval (Step~2), and finally, we read off the estimate (Step~3).

\emph{Step 1: the density is continuous and positive.} Split the Gaussian
matrix as $G_c=G_c^0+\zeta I$, where $\zeta$ is its normalized trace and
$G_c^0$ is traceless. The scalar $\zeta$ is a real Gaussian variable
independent of $G_c^0$, and it shifts all eigenvalues by the same amount.
The expected eigenvalue distribution of $F+G_c$ is therefore the
expected distribution of $F+G_c^0$ convolved with the Gaussian density
of $\zeta$. Hence it has a continuous, strictly positive density, and
$\rho_{F,c}(x)$ is well defined at every point.

\emph{Step 2: shifting the energy to zero.} We apply
Theorem~\ref{thm:external-sine} to the shifted deformation $F-xI$, with
variance parameter $t=c$, at the energy $0$. Shifting the deformation by
$x$ shifts the whole spectrum and the Dyson solution by $x$, so the
statement at energy $0$ for $F-xI$ is the statement at energy $x$ for
$F$. This is why $x$ need not belong to the compact energy interval of
the theorem. We check its hypotheses. On the regular class the norm of
the shifted matrix still has a polynomial bound, after increasing the
fixed exponent $A$. The regularity condition gives the required fixed
positive lower bound on $\Im m_{F,c}(x)$, which is the lower bound on the
imaginary part of the Dyson solution of $F-xI$ at $0$.

\emph{Step 3: the diagonal of the kernel.} The expected counting density
of the eigenvalues is the value $K_\Naux(x,x)$ of the correlation kernel
on the diagonal. The diagonal assertion of
Theorem~\ref{thm:external-sine}, translated back to the energy $x$,
therefore gives
\[
 |\rho_{F,c}(x)-r_{F,c}(x)|\le C\Naux^{-1/2},
\]
uniformly over the regular class. This proves
\eqref{eq:t22-point-density} with $\varepsilon_\Naux=C\Naux^{-1/2}$. No
final averaging over a microscopic window is needed for the estimate; the
scalar Gaussian split of Step~1 is used only to make sure that the
density exists and is positive.
\end{proof}

 \begin{corollary}[Deterministic lower bound for the full fresh normalizer]
\label{cor:t26-deterministic-normalizer}
On the deterministic deformation class of
Lemma~\ref{lem:t22-point-density}, suppose $\Im m_{F,c}(x)\ge v_*>0$.
The complete unranked eigenpair density $Z_{F,c}(x)$, including its
Gaussian coarea constant and the sum of all ordered sectors, satisfies

\[
 Z_{F,c}(x)=\Naux \rho_{F,c}(x)\ge\frac{v_*}{2\pi}\Naux .\]
 
for all sufficiently large   $\Naux $, uniformly over this class.
\end{corollary}
\begin{proof}
Sum the fresh eigenpair coarea formula over its ordered ranks and
integrate its angular variable. The result is exactly the expected
counting density $\E\Tr\delta_x(F+G_c)$; by definition this is
  $\Naux \rho_{F,c}(x)$. Lemma~\ref{lem:t22-point-density} gives
  $\rho_{F,c}(x)\ge v_*/\pi-\varepsilon_\Naux $. Choose the common dimension
threshold so that   $\varepsilon_\Naux \le v_*/(2\pi)$.
This deterministic argument uses the point-density result, not the
averaged exterior total-variation theorem. It supplies no lower bound
for a single ordered sector.
\end{proof}

\subsection{The exact exterior measure, in total variation}\label{subsec:s19-the-exact-exterior-measure-in-total-vari}
 We show that the exact exterior weight becomes asymptotically the
original exterior probability law times the energy-window weight.
Theorem~\ref{thm:t22-normalizer} proves this in total variation. We
express the density ratio through the full fresh normalizer, apply the
pointwise density estimate on the regular set, and bound its complement
before changing energy variables.

 Return to the pair notation of the preceding pair construction, so
  $\Naux =2W$, the fresh variance parameter is $s=\min(b,1-2b)$, and
$K(e)$ is the unchanged real exterior deformation. Define a comparison
probability measure on the same energy/exterior space by

\[
 d\mathfrak p_W(u,\mathcal E)=w(u)\,du\,d\mathbb P_{\rm ext}(\mathcal E),
 \qquad e(u)=E+\frac{u}{N\rho(E)}.
\]

In these coordinates the exact physical exterior window measure is

\begin{equation}\label{eq:t22-exterior-density-ratio}
 \frac{d\mathfrak m_W}{d\mathfrak p_W}(u,\mathcal E)
             =\frac{Z(e(u))}{\Naux \rho(E)}.
\end{equation}

Here $Z$ is the full unranked angular coarea mass, including its
compressed determinant and complete Gaussian normalizer. It is not an
individual indexed sector mass.

\begin{theorem}[Unranked physical exterior normalization]
\label{thm:t22-normalizer}
For fixed $b\in(0,1/2)$, deterministic compact-bulk $E$, fixed $T$ and
fixed nonnegative continuous compactly supported $w$ with integral one,
uniformly on $4\le n\le TW$,

\begin{equation}\label{eq:t22-normalizer-L1}
 \int\left|\frac{Z(e(u))}{\Naux \rho(E)}-1\right|
                     \,d\mathfrak p_W(u,\mathcal E)\longrightarrow0.
\end{equation}

Thus the exact physical exterior window measure and this prior/energy
probability measure agree asymptotically in total variation.
In particular, for every sequence of bounded measurable exterior
observables $|F_W|\le1$, including observables depending on the whole
reconstructed clock,
\[
 \int F_W\,d\mathfrak m_W-\int F_W\,d\mathfrak p_W\longrightarrow0.
\]
\end{theorem}
For use in this theorem, fix a large constant $A$. Exclude energies
within $W^{-A}$ of an exterior pole and Gaussian bond norms exceeding
a fixed sufficiently large constant. Theorem~\ref{thm:external-count}, with factorial order one after conditioning on exterior bonds and rescaling $H_{RR}$ by $\sqrt c$, bounds
the first failure probability by   $CNW^{-A}$ ; the second is exponentially
small. The real deformation is polynomially bounded on the complement.
The residual-pair trace has every fixed moment of size $W^{-1}$ by
Theorem~\ref{thm:matched-traces}. Lemma~\ref{lem:regularity-certificate}
therefore gives the fixed energy/variance rectangle outside an arbitrarily
small polynomial exceptional probability. These estimates hold under the
unchanged exterior prior at deterministic $e$. Inequality
\eqref{eq:t22-exterior-admission} pays them in the actual energy window.

\begin{proof}
On the regular event described below, $K(e)$ has polynomial norm and satisfies
the fresh-density lemma's fixed energy/variance rectangle regularity.
Apply Lemma~\ref{lem:regularity-certificate} to
$\widetilde K=K(e)+eI$ at energy $e$; this removes the energy
already included in the definition of $K(e)$.
Then $m_{\widetilde K,s}(e)=m_{K(e),s}(0)$.
Together with the residual-pair conclusion of
Theorem~\ref{thm:matched-traces}, the lemma places this scalar root at
\[
 m_{K,s}(0)=\mu(e)+O(\epsilon_e),\qquad
 \mu(e)=(-e+i\sqrt{4-e^2})/2,
 \qquad \|\epsilon_e\|_p\le C_p/W.
\]
The certificate's bounded inverse derivatives extend its positive height
to the required fixed rectangle. All these statements have the original
exterior prior at deterministic $e$, uniformly in the allowed lengths.
The complement has arbitrarily small fixed polynomial probability,
including the real-pole norm cutoff paid by the same conditional application of Theorem~\ref{thm:external-count}.

By its exact definition   $Z(e)/\Naux =\rho_{K(e),s}(0)$.
Lemma~\ref{lem:t22-point-density} therefore gives on this event

\[
 |Z(e)/\Naux -\rho(e)|\le\varepsilon_\Naux +C\epsilon_e.
\]

On its complement use the uniform bound   $Z/\Naux \le C_b$.
The deterministic-energy expectation of the displayed difference tends
to zero, uniformly in $e$ in a compact bulk interval and $n\le TW$.
The change of variables in \eqref{eq:t22-exterior-density-ratio} uses
exactly   $n\Naux /(2N)=1$. Since $e(u)-E=O(N^{-1})$  on the support of $w$,
$\rho(e(u))/\rho(E)\to1$. Integrate the previous uniform expectation
against $w(u)du$. This proves \eqref{eq:t22-normalizer-L1} and the
bounded-observable consequence.
\end{proof}

\subsection{Ratios of arbitrary positive reconstruction clocks}\label{subsec:s19-ratios-of-arbitrary-positive-reconstruct}
 We compare reconstructed directions without assuming that their
normalizing mass matrix is well conditioned.
Theorem~\ref{thm:t22-Gaussian-ratios} replaces their bounded quadratic
ratios by Gaussian ratios uniformly over the positive clocks. We use
relative directional moments, control small denominators, and then
approximate continuous tests of the ratios.

 Let $\widehat{\mathcal A}$ be the normalized unranked zero-eigenvector
Palm law for a deterministic regular deformation, of mass $Z>0$ before
normalization.  Corollary~\ref{cor:t26-deterministic-normalizer} gives   $Z\ge c_0\Naux $  for
all sufficiently large   $\Naux $  in this class. Let $D\succeq0$ and
$B_j=B_j^*$, $-D\preceq B_j\preceq D$, $1\le j\le k$, and assume
$T=\operatorname{Tr}(\Gamma D)>0$. They are arbitrary deterministic
matrices relative to the fresh GUE; their norms and condition numbers
are not restricted. Let $g$ be standard circular complex Gaussian in
  $\mathbb C^\Naux $.

\begin{theorem}[Uniform complete ratio comparison]
\label{thm:t22-Gaussian-ratios}
For each fixed $k$ and continuous $\Phi:[-1,1]^k\to\mathbb R$,
\[
 \widehat{\mathcal A}\!\left[
  \Phi\left(\frac{v^*B_jv}{v^*Dv}:1\le j\le k\right)\right]
 -\mathbb E_g\Phi\left(
 \frac{g^*\Gamma^{1/2}B_j\Gamma^{1/2}g}{g^*\Gamma^{1/2}D\Gamma^{1/2}g}:1\le j\le k\right)
       \longrightarrow0
\]
uniformly over the preceding deformation and reconstruction class.
Zero denominators can be assigned arbitrary values. They have zero
Gaussian probability, and their possible Palm contribution is covered
by the proof's small-denominator estimate.
\end{theorem}

The denominator is an arbitrary positive reconstruction clock. Relative Palm moments and positivity allow its condition number to remain unrestricted.

\begin{proof}
For any $Y\succeq0$, average Theorem~\ref{thm:fresh-palm}
over an independent $a\sim\operatorname{CN}(0,Y)$, and sum the local
ranks. The exact conditional identity
$\mathbb E_a|a^*v|^{2r}=r!(v^*Yv)^r$ yields

\[
 \left|\widehat{\mathcal A}(\Naux v^*Yv)^r
       -\mathbb E_g(g^*\Gamma^{1/2}Y\Gamma^{1/2}g)^r\right|
 \le o_r(1)\,[\operatorname{Tr}(\Gamma Y)]^r.
\]

The moment bound for a positive Gaussian quadratic form used here is
$\mathbb E(g^*Ag)^r\le r!(\operatorname{Tr}A)^r$.
There is no dimension factor in this averaging.

Put   $X_0=\Naux v^*Dv/T$ and $X_j=\Naux v^*B_jv/T$.
The corresponding Gaussian variables are denoted $\widetilde X_j$.
To obtain all mixed moments, first replace each $B_j$ by $D+B_j\succeq0$,
then use the real polarization identity
\[
 x_1\cdots x_r=\frac{1}{r!}\sum_{S\subset\{1,\ldots,r\}}
                   (-1)^{r-|S|}\left(\sum_{j\in S}x_j\right)^r.
\]
Every sum of the positive matrices here is bounded above by $2rD$.
The preceding estimate thus compares all fixed mixed moments of
$(X_0,\ldots,X_k)$ and $(\widetilde X_0,\ldots,\widetilde X_k)$,
uniformly. In both laws $|X_j|\le X_0$, and all fixed moments of $X_0$
are uniformly bounded.

The eigenvalues $a_i$ of $\Gamma^{1/2}D\Gamma^{1/2}/T$ are nonnegative
and sum to one. Thus
\[
 \mathbb E_g e^{-s\widetilde X_0}
       =\prod_i(1+sa_i)^{-1}\le(1+s)^{-1},
 \quad
 \mathbb P_g(\widetilde X_0\le h)\le e h,
\]
by taking $s=1/h$. Also, for $0\le t<1$,
\[
 \mathbb E_g e^{t\widetilde X_0}
       =\prod_i(1-ta_i)^{-1}\le(1-t)^{-1}.
\]
The last product inequality follows by induction from
$\prod_i(1-ta_i)\ge1-t\sum_i a_i$.

For a proof of uniformity, suppose there were a sequence violating the
conclusion. Extract weakly convergent subsequences of both joint vectors.
Their higher fixed moments imply uniform integrability for each lower
moment, so the two limits have identical moments. The displayed factorial
bounds imply a finite exponential moment near zero for both limits;
with $|X_j|\le X_0$, this determines their joint laws. The Gaussian
small-ball bound, applied to open sets $\{X_0<h\}$, shows that this
common limit has no mass at $X_0=0$. Division by $X_0$ is therefore
continuous almost surely in the limit. The bounded continuous test
$\Phi$ gives the same limit on both subsequences, a contradiction.
This argument proves exactly the stated uniformity without requiring
uniform inverse moments of the clock.
\end{proof}

\subsection{The physical one-mark comparison}\label{subsec:s19-the-physical-one-mark-comparison}
 We assemble the exterior normalization and Gaussian ratio comparison
into a statement about an actual eigenvector mark.
Theorem~\ref{thm:t22-one-mark-Gaussian} identifies its microscopic
intensity by a normalized Gaussian reconstruction with a uniform
spatial rotation. We define that reconstruction, apply the two
comparisons, and remove the regularity cutoff.

 On the preceding regular exterior set define, with a new independent standard
complex Gaussian $g$,
\[
 \Psi_g=\frac{L\Gamma^{1/2}g}{\|L\Gamma^{1/2}g\|},\qquad
 \nu_g=\sum_{k=1}^n\|P_k\Psi_g\|^2\delta_{k/n}.
\]
Set $\nu_g$ to any fixed probability measure off the regular set.
Let $U_n$ be independent and uniform on $\{0,1/n,\ldots,(n-1)/n\}$,
and let $\mathsf R_x$ denote rotation of a probability measure by $x$.
The Gaussian here is a comparison variable justified by the theorem;
it is not claimed to be the exact finite-$W$ selected-vector law.

\begin{theorem}[Actual one-mark intensity via Gaussian reconstruction]
\label{thm:t22-one-mark-Gaussian}
For every continuous functional $\Phi$ on the compact space
$\mathcal P(\mathbb T)$, under the original ring law with the earlier
uniformity and actual energy marginal,
\begin{equation}\label{eq:t22-one-mark-comparison}
 \mathbb E\sum_\alpha w(u_\alpha)\Phi(\nu_\alpha)
 -\int w(u)\,\mathbb E_{{\rm ext},g,U_n}
          \Phi(\mathsf R_{U_n}\nu_g(e(u)))\,du
       \longrightarrow0.
\end{equation}
For a rotation-invariant $\Phi$, including a squared Fourier mode,
$U_n$ may be omitted.
\end{theorem}

The comparison concerns one physical mark under the stated spectral weight. A common law for several energies is obtained later from the entire determinant.

\begin{proof}
For a cylinder functional use Theorem~\ref{thm:t22-Gaussian-ratios}
with $D=L^*L$ and $B_j=L^*Q_jL$ for the real bounded block tests.
Every continuous functional on $\mathcal P(\mathbb T)$ is uniformly
approximated by such cylinder functionals: the coordinate integrals of
continuous spatial tests separate this compact space, and their
polynomial algebra contains constants. This passes the comparison to
$\Phi$ with uniform error. The original Palm bad events have vanishing
physical mass by \eqref{eq:t22-exterior-admission}. Theorem~\ref{thm:t22-normalizer} removes the full
unranked exterior weight in total variation, also for these bounded
Gaussian expectations.

For a fixed pair this proves the formula for the pair-mass weighted
one-mark intensity. Sum over all cyclic adjacent pairs. The exact
identity $\sum_S q_{\alpha S}=2$ recovers the unweighted physical
intensity. Cyclic invariance of the original law changes this pair
average into the independent uniform rotation $U_n$ of the base-pair
reconstruction. Its energy density is still $w(u)du$, with
  $e(u)=E+u/(N\rho(E))$ ; it has not been frozen or replaced by a random
selected rank. This proves \eqref{eq:t22-one-mark-comparison}.
\end{proof}
\section{Integrating the physical common port}
\label{sec:t22-common-port-profile}
The one-mark formula of Section~\ref{sec:physical-reconstruction} still
contains the three blocks that close the ring. This section integrates them
out, replacing them by their effect on the two incoming directions of the
path. The result identifies the one-mark law of the ring with that of a path
with matched leads, which is used in Part~III.
Theorem~\ref{thm:t22-common-port-ratios} compares ratios of clocks after
the common port has been integrated out, for an arbitrary old clock.
Theorem~\ref{thm:t22-physical-lead-profile} states that the one-mark intensity
of the ring equals the normalized Gaussian response of the path with matched
leads.
We first integrate the fresh coupling and compare positive clock ratios
(Lemmas~\ref{lem:t22-common-port-mean} and~\ref{lem:t26-two-positive-clocks}).
In Subsection~\ref{subsec:s20-admission-to-the-actual-physical-one-mar} we combine these with the exact coarea normalization,
pay for the mass of the removed blocks, and prove
Theorem~\ref{thm:t22-physical-lead-profile}.

Fix $0<b<1/2$, $c=1-2b$, a compact bulk energy range and $n/W\to\tau>0$.
The constants may depend on compact interior parameter sets. Put
\[
 a=(-e-i\sqrt{4-e^2})/2,\quad z_b=a^{-1}+ba,
 \quad z_0=a^{-1}+2ba=-e-ca,
 \quad d=(1-b)^{-1},\quad \gamma=(1-b)/c.
\]
Remove a distinguished physical block, denoted $0$, and its two neighbors
from the ring. The remaining $r=n-3$ blocks form the actual independent
Gaussian path. Attach its actual independent endpoint couplings to two
auxiliary $W$-dimensional ports. Let $X$ be its Hermitian $2W$ Schur
self-energy, $K=z_bI-X$, $R=K^{-1}$, and $J$ its harmonic lift from
port amplitudes to the two ports plus the interior path. Thus the
matched-lead response is
\[
 U=JR=G_{\rm lead}V,\qquad
 Y(q)=U^*Q_qU,\qquad Y=U^*U,
\]
where $V$ embeds both ports and $Q_q$ is the real diagonal block source
on these $n-1$ positions. For $0\le q\le1$, $0\preceq Y(q)\preceq Y$.
The inverse $G_{\rm lead}$ has diagonal $z_b$ at each auxiliary port,
physical diagonal $A_j-e$ in the interior and the actual Gaussian bonds.
It exists at every real energy almost surely. No real cavity pole is
removed from the law by conditioning.

The two missing neighbor diagonals are a fresh $A=\operatorname{diag}(A_1,A_2)$,
with GUE variance $c/W$. Their couplings to block $0$ form a fresh
$2W$-by-$W$ matrix $B=(B_1^T,B_2^T)^T$, with independent entries of
variance $b/W$ in both blocks. Set
\[
 \mathcal G=\begin{pmatrix}A-e-X&B\\B^*&z_0I\end{pmatrix}^{-1},
 \qquad L=K\mathcal G_{12},\qquad H=\mathcal G_{22}.
\]
All these fresh variables are independent of the old path. In this
section $L$ is rectangular, from $W$ to $2W$ dimensions.

\begin{lemma}[Complete common-port covariance]
\label{lem:t22-common-port-mean}
Write $\delta=\max_i|W^{-1}\operatorname{Tr}R_{ii}-a|$.
On its fixed small regularity domain,
\[
 \mathbb E(LL^*\mid X)=\gamma I_{2W}+O_{\rm op}(\delta+W^{-2}).
\]
For arbitrary old vectors $u,v$ and every fixed finite $p\ge2$,
\[
 \|u^*\{LL^*-\mathbb E(LL^*\mid X)\}v\|_{p\mid X}
       \le C_pW^{-1/2}\|u\|\|v\|.
\]
The fresh operator moments of $L,H$ are bounded uniformly in the old
Hermitian deformation. The assertions use the complete common-port
advanced/retarded ladder, including both physical endpoints.
\end{lemma}

The covariance comparison is relative to the matched leads of the same old path. Its scalar normalization fixes the density and the Gaussian input size.

\begin{proof}
The absorbing inverse and weighted-derivative bounds are Lemma~\ref{lem:gaussian-absorption}. Retain the outside $K,K^*$ before estimating residuals, as in Lemma~\ref{lem:passive-bounds}.
Indeed
\[
 K\mathcal G_{11}
   =I+[z_b+e-A+z_0^{-1}BB^*]\mathcal G_{11},
\]
so all required weighted factors have the same uniform moments. Gaussian
Poincare on a bilinear entry of $LL^*$, differentiating these weighted
words, gives the stated $W^{-1/2}$ bound with arbitrary old vectors.
Normalized trace--matrix contractions give $O_{\rm op}(W^{-2})$  in the
conditional mean equations.

For the source equal to the common bottom port, let $t_i$ be the two
normalized top traces of the mean $\mathcal G P_{\rm bottom}\mathcal G^*$,
and let $v_s$ be its normalized bottom trace. At the stationary scalar
point the one-resolvent means are $R$ and $a$. Put
$\mathsf A_{ij}=W^{-1}\operatorname{Tr}R_{ij}R_{ij}^*$.
The full source system is
\begin{equation}
 t=\mathsf A(ct+bv_s\mathbf1),\qquad
 v_s=1+b(t_1+t_2).                                      \label{eq:common-port-ladder}
\end{equation}
The actual equations have $O(W^{-2})$ residuals and scalar coefficients
within $O(\delta)$ of these. Because $X$ is Hermitian,
\[
 \sum_j\mathsf A_{ij}=-\frac{\Im W^{-1}\operatorname{Tr}R_{ii}}{\Im z_b}
       =d+O(\delta).
\]
At the exact row sum $d$, \eqref{eq:common-port-ladder} has
$t_1=t_2=v_s=1/c$. Its feedback matrix is
$\mathsf A(cI+b^2\mathbf1\mathbf1^T)$, with row sum
\[
 d(c+2b^2)=1-\frac{bc}{1-b}<1.
\]
This proves a uniform inverse bound in row-sum norm, and stability for
small $\delta$. The top matrix source coefficient is therefore
$ct_i+bv_s=(1-b)/c+O(\delta)$.
The one-resolvent scalar equations are
\[
 M=(-e-X-\operatorname{diag}(cm_i+bv))^{-1},\qquad
 v=(z_0-bm_1-bm_2)^{-1}.
\]
The common-port implicit-function calculation in Lemma~\ref{lem:layer-contraction}, with \eqref{eq:fresh-trace-fluctuation}, gives $m_i,v=a+O(\delta)$ for this arbitrary fixed Hermitian $X$, hence $KM=I+O_{\rm op}(\delta)$.
The weighted source equation now gives the asserted mean of $LL^*$,
including its absolute $O(W^{-2})$ residual. No off-diagonal transmission
entry was set to zero; only the exact row sum was used.
\end{proof}

\begin{theorem}[Uniform Gaussian-clock comparison across the common port]
\label{thm:t22-common-port-ratios}
Let $Y\succeq0$, $\operatorname{Tr}Y>0$, and let
$-Y\preceq Y_j\preceq Y$, $1\le j\le k$, be arbitrary old matrices.
For independent standard circular Gaussians $g\in\mathbb C^W$ and
$h\in\mathbb C^{2W}$, and each fixed continuous $\Phi:[-1,1]^k\to\mathbb R$,
\[
 \mathbb E_{A,B,g}\Phi\left(
 \frac{g^*L^*Y_jLg}{g^*L^*YLg}:1\le j\le k\right)
 -\mathbb E_h\Phi\left(
 \frac{h^*Y_jh}{h^*Yh}:1\le j\le k\right)\longrightarrow0
\]
uniformly as $W\to\infty$ and $\delta\to0$.
There is no condition-number, norm or rank assumption on these clocks.
\end{theorem}

Positive numerator and denominator clocks are compared together. The statement avoids controlling an inverse reconstruction clock before taking the limit.

\begin{proof}
Scale $\operatorname{Tr}Y=1$ and put $\mathsf U=LL^*$.
Lemma~\ref{lem:t22-common-port-mean}, in an eigenbasis of $Y$, gives
\[
 \mathbb E\|Y^{1/2}(\mathsf U-\gamma I)Y^{1/2}\|_{\rm HS}^2
       \le C(W^{-1}+\delta^2),
\]
and, for every bounded old matrix $D$,
\[
 \mathbb E|\operatorname{Tr}[D Y^{1/2}(\mathsf U-\gamma I)Y^{1/2}]|
       \le C\|D\|(W^{-1/2}+\delta)
\]
when $D$ is a bounded function of $Y$. The last assertion follows by
summing the diagonal bilinear bounds with their positive weights.
For any fixed $a\ge0$, expand $\log\det(I+aZ)$ on the line joining
$Z=\gamma Y$ to $Z=Y^{1/2}\mathsf UY^{1/2}$. All matrices on that line
are positive. The linear term is the preceding trace with
$D=(I+a\gamma Y)^{-1}$, and the remainder has absolute value at most
$a^2\|Y^{1/2}(\mathsf U-\gamma I)Y^{1/2}\|_{\rm HS}^2/2$.
Since $e^{-u}$ is 1-Lipschitz for $u\ge0$, this proves
\begin{equation}
 \mathbb E\det(I+aL^*YL)^{-1}
       -\det(I+a\gamma Y)^{-1}\longrightarrow0.          \label{eq:common-port-laplace}
\end{equation}
The same argument applies to any fixed positive combination of the
positive matrices $Y$ and $Y\pm Y_j$, after trace normalization.
Their Gaussian integral identities therefore compare the joint Laplace
transforms of all required positive quadratic forms. Uniform first
moments give tightness; a subsequence/compactness argument and uniqueness
of multivariate Laplace transforms give uniform weak comparison.

Denominators are paid before taking ratios. For any positive trace-one
$Y$, $\det(I+a\gamma Y)\ge1+a\gamma$. Use \eqref{eq:common-port-laplace} with
$a=\varepsilon^{-1/2}$ and the exponential Markov bound to obtain,
for both laws,
\[
 \mathbb P\{\hbox{denominator}\le\varepsilon\}
     \le C\sqrt\varepsilon+o_W(1)
\]
after its fixed scalar normalization. On the complementary set the
ratio map is continuous with a fixed modulus on tight compact sets.
First let $W\to\infty$, then $\varepsilon\downarrow0$.
The common scalar $\gamma$ cancels from every comparison ratio.
This proves the theorem without a dimension-dependent norm estimate.
\end{proof}

A clock is a positive quadratic form of the
eigenvector, such as $Nv^*Dv$, normalized by its expected size $T$. The following
lemma compares two such clocks, built from two positive matrices $D$
and $C$, with the same two clocks built from a Gaussian vector. It says
that the pair of clocks has asymptotically the same joint law in both
cases. It also says that their ratio has the same law, which is the form
used later. The matrix $D$ fixes the normalization, and $C$ is allowed to
be any positive matrix whose trace is at most $K$ times that of $D$. In
particular, $C$ need not be dominated by $D$, so the ratio can be large.

\begin{lemma}[Two positive clocks and trace cutoffs]
\label{lem:t26-two-positive-clocks}
Fix $K<\infty$ and $\Phi\in C_b([0,\infty))$. In the deterministic
regular class of Theorem~\ref{thm:t22-Gaussian-ratios}, take arbitrary
$D,C\succeq0$ with $T=\Tr(\Gamma D)>0$ and $\Tr(\Gamma C)\le KT$. Then the
joint laws of
\[
 \left(\frac{Nv^*Dv}{T},\frac{Nv^*Cv}{T}\right),\qquad
 \left(\frac{g^*\Gamma^{1/2}D\Gamma^{1/2}g}{T},
       \frac{g^*\Gamma^{1/2}C\Gamma^{1/2}g}{T}\right)
\]
agree asymptotically against bounded continuous tests, uniformly in
these matrices and deformations. The ratios of numerator to denominator
also agree against $\Phi$. In the actual conditional common-port law of
Lemma~\ref{lem:t22-common-port-mean}, the same assertions hold for
\[
 \left(\frac{g^*L^*DLg}{\gamma T},
       \frac{g^*L^*CLg}{\gamma T}\right),\qquad
 \left(\frac{h^*Dh}{T},\frac{h^*Ch}{T}\right),
 \quad T=\Tr D>0,\quad \Tr C\le KT,
\]
uniformly as $W\to\infty$ and $\delta\to0$. Here the old matrices are
measurable with respect to the entire old path, while all fresh
Gaussians are integrated conditionally on that path. No Loewner
domination of $C$ by $D$ is imposed.
\end{lemma}

\begin{proof}
A clock here is a positive quadratic form of a
random vector, and we want the joint law of two clocks to agree with the
law obtained when the vector is replaced by a Gaussian vector. For an
eigenvector selected at a given energy (the Palm case) we compare all
joint moments of the two clocks (Step~1). In the common-port case we
compare their joint Laplace transforms, that is, expectations of
exponentials of the clocks (Step~2). Either family determines the joint
law. To pass from the two clocks to their ratio, we must exclude a small
denominator; a Gaussian small-ball bound does this, as in
Figure~\ref{fig:tc-cutoff} (Step~3).
Step~4 shows that the comparison survives a small relative change of the
covariance, that is, it compares the Gaussian laws whenever
$(1-\omega)\Gamma_0\preceq\Gamma\preceq(1+\omega)\Gamma_0$ in the sense of
quadratic forms. There we normalize the clocks by $T_0=\Tr(\Gamma_0D)$.
If an earlier formula uses $T=\Tr(\Gamma D)$ instead, the two
normalizations differ only by the factor $T/T_0\in[1-\omega,1+\omega]$.
This rescaling is harmless. Tightness, together with uniform continuity
of the test functions on compact sets, controls it, and the ratios do not
change at all under a common rescaling.
Finally, the case of finitely many ratios is a separate consequence of
the single-ratio comparison, not an unstated part of it. We apply the
Laplace comparison for one pair of clocks to the positive combinations
$C_\lambda=\sum_j\lambda_j(D+B_j)/2\succeq0$, with $\lambda_j\ge0$. This
gives the joint Laplace transforms of the clocks $(D+B_j)/2$. The vector
of ratios takes values in a compact cube, so uniqueness of the Laplace
transform identifies its law there.

\emph{Step 1: the fresh-Gaussian comparison.} We start with the first
pair in the deterministic regular class. The positive-source moment
estimate of Theorem~\ref{thm:t22-Gaussian-ratios} applies to every
\[
 \frac{aD+bC}T\succeq0,\qquad a,b\ge0,
\]
whose trace against $\Gamma$ is at most $a+bK$. Expanding the moments of
$aX_0+bX_1$ in $a$ and $b$, polarization therefore compares every fixed
mixed moment of the two coordinates, with errors $o_r(1)(1+K)^r$. The
moments of the sum of the two coordinates are at most
$r!(1+K)^r+o_r(1)$, as for a sum of independent exponential variables
with total weight $1+K$. Tightness and the higher moments let each mixed
moment pass to any subsequential limit. These factorial bounds give a
finite exponential moment near zero in the limit. So the common moments
determine the joint law.

The comparison is uniform in the matrices and deformations. Indeed,
suppose that it is not. Then there is a sequence of matrices along which
some bounded continuous test differs by at least a fixed amount. By
tightness we may extract subsequential limits of both laws. They would
be two distinct laws with the same moments, which is impossible. This
comparison holds under the normalized \emph{unranked} Palm measure,
because the earlier sector estimates were summed before the
normalization by $Z\ge c_0N$.

\emph{Step 2: the common-port comparison.} For the second pair we use
Laplace transforms instead of moments. We put
\[
 F=\frac{aD+bC}{\gamma T},\qquad \mathsf U=LL^*.
\]
Integrating out the Gaussian $g$, the Laplace transform of the first
pair is an inverse determinant involving $F^{1/2}\mathsf UF^{1/2}$, and
the second pair corresponds to replacing $\mathsf U$ by $\gamma I$. So we
need to show that $\mathsf U$ is close to $\gamma I$, tested against $F$.
The bilinear estimates, summed in an eigenbasis of $F$, give
\[
 \E\|F^{1/2}(\mathsf U-\gamma I)F^{1/2}\|_{\rm HS}^2
 \le C(W^{-1}+\delta^2)(\Tr F)^2,
\]
\[
 \E|\Tr[(I+\gamma F)^{-1}
             F^{1/2}(\mathsf U-\gamma I)F^{1/2}]|
 \le C(W^{-1/2}+\delta)\Tr F.
\]
The first bound controls the quadratic term and the second the linear
term in the expansion of the log-determinant of a positive segment.
That expansion therefore proves the explicit joint Laplace estimate
\begin{equation}\label{eq:t26-two-clock-laplace}
 \left|\E e^{-aX_0-bX_1}-\E e^{-a\widetilde X_0-b\widetilde X_1}\right|
 \le C_{a,b,K}(W^{-1/2}+\delta),
\end{equation}
where $X$ and $\widetilde X$ denote the two common-port pairs. Uniform first
moments give tightness. Every two subsequential limits have the same
Laplace transform for $a,b\ge0$, and hence the same law. As in Step~1,
this proves the uniform joint statement in this case as well.

\emph{Step 3: small denominators.} To pass to the ratio $X_1/X_0$ we must
show that the denominator is rarely small. In both comparisons the
Gaussian denominator has a covariance of trace one. For $A\succeq0$ with
$\Tr A=1$ we have $\det(I+sA)\ge1+s$, and so
\[
 \E e^{-s\widetilde X_0}\le(1+s)^{-1},\qquad
 \Prob(\widetilde X_0\le u)\le eu.
\]
The second bound follows from the first by Chebyshev's inequality with
$s=1/u$. Consequently all limiting denominators are strictly positive,
and
\[
 \lim_{u\downarrow0}\limsup\Prob(X_0\le u)=0
\]
for both laws, uniformly in the indicated class. The numerator has first
moment at most $CK$, uniformly, so it is tight.

On the set $\{X_0\ge u,\ X_1\le R\}$, the ratio is a continuous function
on a compact set. So for this set the joint comparison of Steps~1 and~2
passes to $\Phi(X_1/X_0)$. Its complement has small probability, by the
bounds above. We first take the limit in the width, then let
$R\to\infty$, and finally let $u\downarrow0$. This proves the comparison
for every fixed $\Phi\in C_b$. We never claim inverse moments of the
denominator before the limit.

\begin{figure}[H]
\centering
\begin{tikzpicture}[x=1cm,y=1cm]
 \fill[black!10] (0,0) rectangle (.6,3.6);
 \fill[black!10] (0,3.0) rectangle (6,3.6);
 \fill[blue!12] (.6,0) rectangle (6,3.0);
 \draw[->] (0,0)--(6.4,0) node[below,font=\small] {$X_0$};
 \draw[->] (0,0)--(0,4.0) node[left,font=\small] {$X_1$};
 \draw[dashed] (.6,0)--(.6,3.6);
 \draw[dashed] (0,3.0)--(6,3.0);
 \node[font=\small,below] at (.6,0) {$u$};
 \node[font=\small,left] at (0,3.0) {$R$};
 \node[font=\small,align=center,blue!70!black] at (3.4,1.5)
   {$X_1/X_0$ continuous\\ and bounded here};
 \node[font=\small,align=center] at (3.4,3.3) {numerator large: tight};
 \draw[->,thin] (-.9,1.6) node[left,font=\small,align=right]
   {denominator\\ small: $\le eu$} -- (.3,1.6);
\end{tikzpicture}
\caption{The cutoff in Step~3 of the proof of
Lemma~\ref{lem:t26-two-positive-clocks}. On the shaded region
$\{X_0\ge u,\ X_1\le R\}$ (blue), the ratio $X_1/X_0$ is a continuous
bounded function, so the joint comparison of the two clocks passes to
it. The grey strips have small probability. Near $X_0=0$ we use
$\Prob(\widetilde X_0\le u)\le eu$, which holds because the Gaussian
denominator has a covariance of trace one. For large $X_1$ we use the
uniform first moment of the numerator. Letting first $W\to\infty$, then
$R\to\infty$ and $u\downarrow0$, removes both strips, with no inverse moments
of the denominator.}
\label{fig:tc-cutoff}
\end{figure}
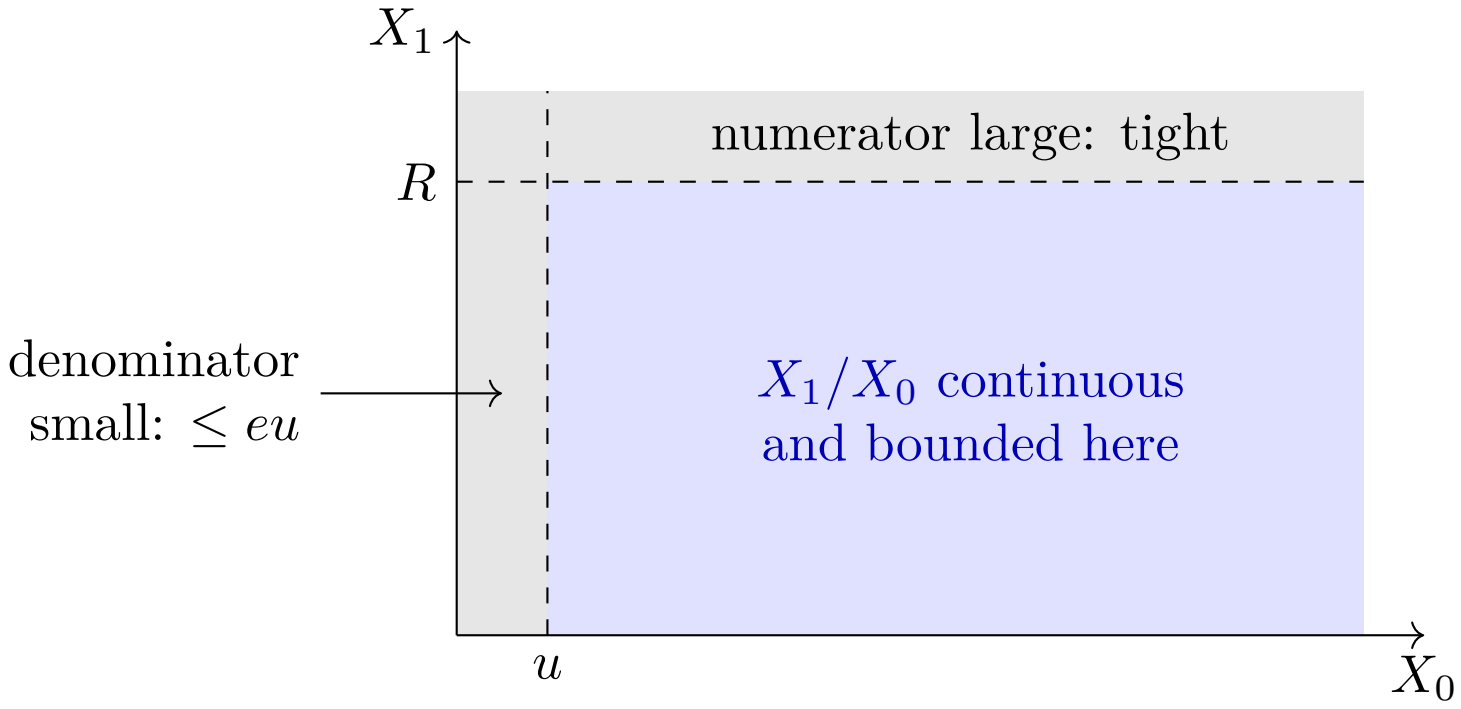

\emph{Step 4: stability under a change of covariance.} Finally, the
assertion is stable under a relative change of the covariance,
\[
 (1-\omega)\Gamma_0\preceq\Gamma\preceq(1+\omega)\Gamma_0,
 \qquad \omega\to0 .
\]
We interpolate by $\Gamma_t=\Gamma_0+t(\Gamma-\Gamma_0)$, $0\le t\le1$. For each
positive combination $F$ normalized by $\Tr(\Gamma_0D)$, Sylvester's
identity lets us write the derivative of
$\log\det(I+\Gamma_t^{1/2}F\Gamma_t^{1/2})$ as
\[
 \Tr\bigl[(I+F^{1/2}\Gamma_tF^{1/2})^{-1}
 F^{1/2}(\Gamma-\Gamma_0)F^{1/2}\bigr].
\]
Since $-\omega\Gamma_0\preceq\Gamma-\Gamma_0\preceq\omega\Gamma_0$, and the inverse
has norm at most one, this is at most $\omega\Tr(\Gamma_0F)$ in absolute
value. So the Laplace transforms change by $O(\omega)$, and the same
Laplace and small-denominator argument applies. We never commute matrix
square roots.

Two consequences of this bound are used later. First, as noted at the
start of the proof, the normalizations $T=\Tr(\Gamma D)$ of Steps~1--3
and $T_0=\Tr(\Gamma_0D)$ here satisfy $T/T_0\in[1-\omega,1+\omega]$, as
$D\succeq0$. For $\omega\le1/2$ the rescaling by $T_0/T$ moves the clocks
by at most $2\omega$ times their size; their first moments are at most
$2(1+K)$, so they are tight uniformly in the class, and the rescaling
changes every bounded continuous test by $o(1)$ as $\omega\to0$,
uniformly, without changing the ratios. With the Laplace bound and the
contradiction argument of Step~1, this proves the joint and ratio
comparisons between the Gaussian laws with covariances $\Gamma$ and
$\Gamma_0$, uniformly in the dimension and the matrices as $\omega\to0$.

Second, finitely many ratios. Let $B_1,\ldots,B_k$ be Hermitian with
$-D\preceq B_j\preceq D$, and put $P_j=(D+B_j)/2$, so that
$0\preceq P_j\preceq D$. In any one of the comparisons above, write
$X(P)$ for the clock of a matrix $P$ and $r_j=X(P_j)/X(D)\in[0,1]$ on
$\{X(D)>0\}$. For $\lambda\in[0,\infty)^k$ the matrix
$C_\lambda=\sum_j\lambda_jP_j\succeq0$ satisfies the trace condition
with $K=\sum_j\lambda_j$, because $\Tr(\Gamma P_j)\le\Tr(\Gamma D)$
(respectively $\Tr P_j\le\Tr D$), and its ratio is $\lambda\cdot r$.
The single-ratio comparison with $\Phi(x)=e^{-x}$ therefore compares
$\E e^{-\lambda\cdot r}$ under the two laws, uniformly in the matrices,
for each fixed $\lambda$. Laws on the compact cube $[0,1]^k$ are tight and determined by their
Laplace transforms on $[0,\infty)^k$, so a continuous test on the cube
with comparison error bounded below along a sequence would give two
distinct subsequential limits with the same Laplace transform. So the
ratio vector $r$, and with it $(X(B_j)/X(D))_{j\le k}=(2r_j-1)_{j\le k}$,
obeys the same uniform comparison against every continuous test on
the cube.
\end{proof}

\begin{corollary}[Uniform tags without a spectral domination assumption]
\label{cor:t26-uniform-clock-tags}
Suppose $D_j\succeq0$, $1\le j\le q$, and $\sum_jD_j=D$.
In either part of Lemma~\ref{lem:t26-two-positive-clocks}, let $J$ be
independent uniform on these $q$ indices and put $C=qD_J$.
The ratio comparison holds uniformly over the entire family
$(D_j)_{j\le q}$, even when $q$ grows with the dimension.
\end{corollary}
\begin{proof}
For the Palm part set
$a_j=q\Tr(\Gamma D_j)/\Tr(\Gamma D)$; for the common-port part
set $a_j=q\Tr D_j/\Tr D$. In both cases
$q^{-1}\sum_ja_j=1$ exactly, so the proportion of tags with $a_j>K$
is at most $1/K$. Each such tag costs at most $2\|\Phi\|_\infty$.
On the other tags the preceding lemma is uniform, giving total error
at most $o_K(1)+2\|\Phi\|_\infty/K$. Take the width limit first,
then $K\to\infty$. The argument is pointwise in the exterior and
therefore applies before or after its physical Palm weighting.
For finitely many independent tags use the union bound.
\end{proof}

 \subsection{Admission to the actual physical one-mark law}\label{subsec:s20-admission-to-the-actual-physical-one-mar}
 We finish the passage from the open path to the original ring's
one-mark law. Theorem~\ref{thm:t22-physical-lead-profile} retains the
actual microscopic energy variable and both boundary inputs. We
combine the exact coarea normalization with the common-port
comparison, remove the exceptional set, and restore spatial positions
by a uniform rotation.

 For the matched-lead path just defined and independent
$h\sim\operatorname{CN}(0,I_{2W})$, put
\[
 \nu^{\rm lead}_W(q)=\frac{h^*Y(q)h}{h^*Yh}.
\]
Give its $n-1$ blocks their original positions on the circle, with zero
mass at the removed block. Let $\mathcal R_J$ rotate positions by the
independent uniform $J\in\{0,\ldots,n-1\}$.

\begin{theorem}[Physical one-mark comparison with independent matched leads]
\label{thm:t22-physical-lead-profile}
For every continuous bounded functional $\Phi$ on probability measures
on the circle, every fixed smooth compactly supported spectral weight
$w$, and $n/W\to\tau\in(0,\infty)$,

\begin{equation}
 \mathbb E\sum_\alpha w(u_\alpha)\Phi(\nu_\alpha)
 -\int w(u)\,
       \mathbb E\Phi(\mathcal R_J\nu^{\rm lead}_W(e(u)))\,du
       \longrightarrow0,\qquad e(u)=E+\frac{u}{N\rho(E)}. \label{eq:physical-lead-profile}
\end{equation}

The expectation on the right is the actual independent path and port
prior. Both lead inputs are independent standard Gaussians. No spectral
root, uncomputed exterior determinant or local covariance $\Gamma$
remains in this one-mark comparison. The full two-lead Green response
$Y(q)$, including interference through the actual path, remains present.
The assertion concerns one-mark intensity, not joint laws of distinct
spectral points.
\end{theorem}

The resulting Gaussian inputs feed the two ends of one actual path. Their reconstructed interior fields still include the interference between those inputs.

\begin{proof} The one-block version uses $A(e)=G+K(e)$ with
$K(e)=-e-\Sigma(e)$, $D=I+T^*T=-K'(e)$ and fresh
$G\sim\GUE_W(c/W)$. At a crossing,
$\mu_k'(e)=-v^*Dv$ and $p_{\alpha0}=(v^*Dv)^{-1}$, hence
\[
 \E\sum_\alpha p_{\alpha0}g(\lambda_\alpha)
 =\E_{\rm ext}\int g(e)Z_0(e)\,de,
 \qquad Z_0(e)=W\rho_{K(e),c}(0).
\]
Cyclic symmetry gives the window prefactor $n$ and $nW/N=1$.
The one-block application of Theorem~\ref{thm:external-count} gives
$Z_0/W\le C/\sqrt c$, paying exterior admission on the exceptional
set. For the regularity certificate the residual is the complex
conjugate of $W^{-1}\Tr(z_0-\Sigma)^{-1}-a$, with centre $\bar a$.
The normalization and Gaussian-ratio argument therefore uses the
one-block Palm law with its actual exterior prior.

Theorems~\ref{thm:t22-normalizer} and~\ref{thm:t22-Gaussian-ratios}
apply also to the one-block coarea law: its fresh
variance is $c/W$, its dimension is $W$, and Theorem~\ref{thm:matched-traces} and Lemma~\ref{lem:regularity-certificate} supply its
regular deformation and polynomial tails. The normalizer proof gives
$Z/W\to\rho(e)$ in prior $L^1$; the coarea prefactor is now $n$, so
  $nW/N=1$. The all-fixed-moment argument is formulated for the general
regular fresh deformation, hence the same complete Gaussian-ratio
comparison holds here. The stationary scalar approximation changes its
covariance to $c|z_0-\Sigma|^{-2}$ with relative Loewner error $o(1)$.
Indeed $\Im z_0=c\pi\rho(e)$ stays strictly positive, and the relative
resolvent comparison is bounded by the scalar self-energy error divided
by this height. It is valid against arbitrary positive reconstruction
clocks and thus preserves the ratio comparison.

The absorbed-density identity identifies the resulting Gaussian field.
Indeed for the full inverse $G_{\rm abs}$ with imaginary diagonal
$\Im z_0=c\pi\rho(e)$ only at block $0$,
\[
 -\Im G_{\rm abs}
   =c\pi\rho(e)G_{\rm abs}P_0G_{\rm abs}^*.
\]
The real Schur lift of $c|z_0-\Sigma|^{-2}$ is consequently
$c$ times the covariance of the block column $G_{\rm abs}P_0$.
Splitting this column into the old path and block $0$ gives exactly
\[
 (J\mathcal G_{12}g,\ Hg)
      =(JR Lg,\ Hg),
\]
up to the common factor $\sqrt c$, which cancels on normalization.
This retains the original two incident bonds, the entire prior and the
actual energy. It is not a new physical matrix ensemble.

The removed-block mass is negligible in this locally weighted Palm law.
To see this without an inverse-clock bound, Theorem~\ref{thm:small-critical} and the count
moments imply
\begin{equation}
 \mathbb E\sum_\alpha |w(u_\alpha)|\sum_k p_{\alpha k}^2\longrightarrow0.
                                                               \label{eq:vanishing-block-atom}
\end{equation}
Every critical marked subsequence has nonatomic marks, so its largest
single-block mass tends to zero; compactness and the count moment bound
then give \eqref{eq:vanishing-block-atom} along the full sequence. Under the block-$0$ Palm law
the removed-block expectation is exactly
$n\mathbb E\sum_\alpha |w(u_\alpha)|p_{\alpha0}^2$, equal by cyclic
symmetry to \eqref{eq:vanishing-block-atom}. Theorems~\ref{thm:t22-normalizer} and~\ref{thm:t22-Gaussian-ratios} transfer the same bounded mass ratio to the
Gaussian field above. Removing that coordinate and renormalizing
therefore changes any continuous profile functional by $o(1)$.

For the remaining coordinates its quadratic forms are exactly
$g^*L^*Y(q)Lg$. Apply Theorem~\ref{thm:t22-common-port-ratios}
conditionally on the actual old path. Theorem~\ref{thm:matched-traces} gives $\mathbb E\delta=O(W^{-1})$,
and its arbitrary polynomial tails remove the scalar exceptional set.
The comparison is uniform in the old clocks, so their norms, poles and
condition numbers cause no extra exceptional set or normalizer loss.
This proves the comparison for finitely many continuous spatial tests;
compactness of the space of probability measures extends it to every
continuous $\Phi$.
Finally the exact identity $\sum_kp_{\alpha k}=1$ and cyclic invariance
express the unweighted physical one-mark intensity as the uniform
rotation of the block-$0$ Palm law. This proves \eqref{eq:physical-lead-profile}, with the stated
energy marginal. Signed $w$ follows by decomposition or domination by
a nonnegative compactly supported smooth weight.
\end{proof}
\section{A critical law for the complete open clock}
\label{sec:open-clock-law}
Profiles are normalized by the full clock, so we must know that this clock
is not too small. As a first step, this section identifies the critical limit
of the clock of a path with one absorbing boundary. The comparison with the
two-input clock of the ring is a separate step, carried out in
Section~\ref{sec:clock-lower-tail}.
Theorem~\ref{thm:t22-open-clock-law} states that the open clock converges
in law to the measure
$d\exp\{\sqrt{D_0}B_{\tau-x}-D_0(\tau-x)/2\}\mathbf1_{[0,\tau]}(x)\,dx$, a
geometric Brownian profile.
In Subsection~\ref{subsec:s21-the-one-ended-path} we set up the one-ended path, prove its conditional
mean and its entrance from the zero boundary
(Lemmas~\ref{lem:c01-cavity-mean}, \ref{lem:c01-zero-entrance}
and~\ref{lem:c01-open-inputs}). In Subsection~\ref{subsec:s21-the-source-gap-for-trace-class-clock-wei} we prove mixing for
the complete clock (Lemma~\ref{lem:t22-clock-mixing}). In
Subsection~\ref{subsec:s21-a-closed-determinant-laplace-equation} we derive a closed equation for the Laplace transform of
the clock (Lemma~\ref{lem:t22-det-clock}) and prove the theorem.

\subsection{The one-ended path}\label{subsec:s21-the-one-ended-path}
 We set up the one-ended path and retain the response from every
physical source block. Theorem~\ref{thm:t22-open-clock-law} gives its
limiting Brownian profile law. We establish the conditional mean,
verify entry from the actual zero-boundary initial condition, and
control the derivative and unstopped clock moments needed for the
transform equation proved below.

 Use the one-ended cavity chain
\[
 X_j=B_j^*(A_j-e-X_{j-1})^{-1}B_j,\quad X_0=0,\qquad
 a=(-e-i\sqrt{4-e^2})/2,\quad z=a^{-1}+ba,
\]
\[
 \mathcal G_j=\begin{pmatrix}A_j-e-X_{j-1}&B_j\\B_j^*&zI\end{pmatrix}^{-1},
 \quad L_j=(z-X_{j-1})(\mathcal G_j)_{12},\quad
 C_j=(\mathcal G_j)_{12}^*(\mathcal G_j)_{12}.
\]
Here $A_j$ has GUE entry variance $(1-2b)/W$, $B_j$ has independent
circular complex entries of variance $b/W$, and all layers are independent.
For a deterministic nonnegative continuous function $q$ on $[0,T]$ set
\[
 Y_0^q=0,\qquad
 Y_j^q=q(j/W)C_j+L_j^*Y_{j-1}^qL_j,
 \qquad Q_j^q=W^{-2}Y_j^q,\qquad \gamma_j^q=\operatorname{Tr}Q_j^q.
\]
The sign convention has $\operatorname{Im}z>0$. Set
\[
 d=(1-b)^{-1},\quad D_0=\frac{2}{b(4-e^2)},\qquad
 f(X)=(z-X)^{-1}-(z-\overline z)^{-1}I.
\]
All conclusions below hold for fixed $b\in(0,1/2)$, deterministic
compact-bulk $e$, $T<\infty$ and bounded continuous $q\ge0$, uniformly
for $j\le TW$. The moment orders used in the proof are fixed.

The zero boundary requires its own entrance proof. The two-matched-port
cores in Theorem~\ref{thm:matched-traces} are a different law. We give the
one-ended argument next, including the first cavity-spectral jet, the
growing moment range used to remove exits, and the initial clock.
Throughout this argument $b\in(0,1/2)$ is fixed and $e$ ranges in a fixed compact subset of
$(-2,2)$. Constants are uniform in $e$ on that set.. All finite moment assertions require sufficiently
large $W$ for the stated moment and word lengths.

\begin{figure}[ht]
\centering
\begin{tikzpicture}[x=1cm,y=1cm]
 \draw[thick,fill=black!45] (0,0) rectangle (1.6,.9);
 \node[font=\small,white] at (.8,.45) {old $X$};
 \draw[thick,fill=black!12] (2.6,0) rectangle (3.6,.9);
 \node[font=\small] at (3.1,.45) {$A$};
 \draw[very thick,dashed] (4.6,0) rectangle (5.6,.9);
 \node[font=\small] at (5.1,.45) {$z'$};
 \draw[thick] (1.6,.45)--(2.6,.45);
 \draw[thick] (3.6,.45)--(4.6,.45) node[midway,below,font=\small] {$B$};
 \node[font=\small,align=center] at (3.1,-.75) {physical block\\ trace $v$};
 \node[font=\small,align=center] at (5.1,-.75) {port\\ trace $u$};
 \node[font=\small,align=center] at (.8,-.75) {eliminated\\ old path};
 \draw[->,thick,blue!70!black] (3.3,1.05) to[out=60,in=120]
   node[above,font=\small] {$b$} (4.9,1.05);
 \draw[->,thick,blue!70!black] (4.9,.95) to[out=150,in=30]
   node[below,font=\small] {$b$} (3.3,.95);
 \draw[->,thick,blue!70!black] (2.85,1.0) .. controls (2.2,2.0) and (3.4,2.0) .. (2.95,1.05);
 \node[font=\small,blue!70!black] at (2.55,1.85) {$c$};
\end{tikzpicture}
\caption{The layer in Lemma~\ref{lem:c01-cavity-mean}. The old path
(dark) has been eliminated, and it enters only through its Schur
complement term $X$. One fresh physical block $A$, at energy $e$, is
coupled by a fresh Gaussian bond $B$ to a port with parameter $z'$
(dashed). Integrating out the fresh Gaussians couples the two
normalized traces. The physical trace $v$ feels itself with weight $c$,
the variance of $A$, and the port trace $u$ with weight $b$, the
variance of $B$. The port feels the physical block with weight $b$. This
gives the two equations \eqref{eq:c01-cavity-map}. The condition
$\Im z'\ge\eta_0$ at the port makes the system stable, for every $X$.}
\label{fig:cv-layer}
\end{figure}
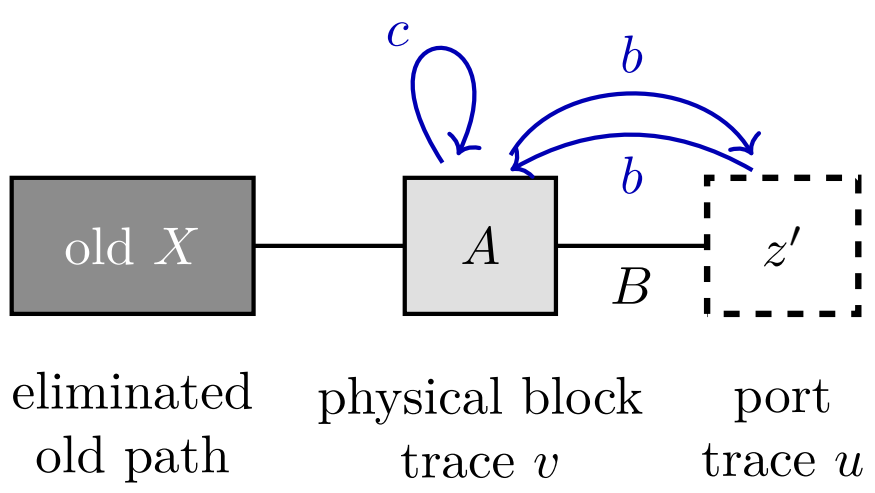

We add one fresh layer to an old path whose
Schur complement term is $X$. The layer consists of one fresh physical
block $A$, at energy $e$, coupled by a fresh bond $B$ to a port with
parameter $z'$, as in Figure~\ref{fig:cv-layer}. The new Schur
complement term is $X_+=B^*(A-e-X)^{-1}B$, and $h_+(z')$ is the
normalized trace of the new port resolvent. The following lemma says that $h_+$
is close to a deterministic number $u$, given by the two scalar
equations \eqref{eq:c01-cavity-map}. Here $v$ is the normalized trace at
the physical block, and $u$ the one at the port. The equations come
from integrating out the fresh Gaussians. The important feature is that
the lemma is global: it needs no assumption that $X$ is close to a
matched value, and its constants do not depend on $X$.

\begin{lemma}[Global one-ended conditional mean]\label{lem:c01-cavity-mean}
For an independent Hermitian $X$, write
$f(\zeta)=W^{-1}\Tr(\zeta-X)^{-1}$. For $\Im z'\ge\eta_0>0$ and
$|z'|\le M$, the equations
\begin{equation}\label{eq:c01-cavity-map}
 v=f(-e-cv-bu),\qquad u=(z'-bv)^{-1}
\end{equation}
have a unique solution with negative imaginary parts. For the actual
fresh $X_+=B^*(A-e-X)^{-1}B$ and $h_+(z')=W^{-1}\Tr(z'-X_+)^{-1}$,
\[
 |\E(h_+(z')\mid X)-u|\le CW^{-2},\qquad
 \|h_+(z')-\E(h_+(z')\mid X)\|_{p\mid X}
       \le C\sqrt p/W,\quad 2\le p\le a_0W.
\]
The same bounds hold for each fixed complex derivative in $z'$ at an
interior point of this domain, after decreasing $a_0$. The physical
energy $e$ is held fixed in these derivatives. The constants are
independent of $X$. Let $\mathcal G$ be the inverse of Lemma~\ref{lem:gaussian-absorption}
for this layer, with $E=e$ and $z'$ in place of $z$. Its physical mean
trace satisfies, with a constant independent of $X$,
\[
 |\E(W^{-1}\Tr\mathcal G_{11}\mid X)-v|\le CW^{-2}.
\]
\end{lemma}

\begin{proof}
Integrating out the fresh Gaussian block and bond
gives two approximate scalar equations, one for the average normalized trace
at the physical block and one for the average normalized trace at the port,
as in Figure~\ref{fig:cv-layer} (Step~1). We show that the exact equations
have a unique solution, and that this solution is stable, so the approximate
means are close to it. The key quantity is a gap: taking imaginary parts of
the equations gives identities (Step~2) which show that the feedback of the
system stays strictly below one, whatever the old matrix is (Step~3). This
gap gives existence and uniqueness of the solution (Step~4) and the bound on
the error of the means (Step~5). Step~6 adds the fluctuation bound and the
derivatives in the port parameter.

\emph{Step 1: the approximate mean equations.} We use the case of
Lemma~\ref{lem:gaussian-absorption} with one physical block and one
port. The bottom block of its inverse is $(z'-X_+)^{-1}$, so $h_+$ is the
normalized trace of that block. The exact Gaussian mean equations in the
proof of Lemma~\ref{lem:layer-contraction} apply in this one-block case
without any local assumption. We put $D=-e-X$, and we let $\bar v$ and
$\bar u$ be the two normalized mean traces. The port equation, the
inverse moments and the Ward identity give
\[
 (z'-b\bar v)\bar u=1+O(W^{-2}),\qquad
 |\bar u|\ge a_*>0,\qquad -\Im\bar u\ge\eta_0|\bar u|^2.
\]
Here $\bar v=\E\langle\mathcal G_{11}\rangle$ for the inverse
$\mathcal G$ of the statement. The inverted matrix has imaginary part
$\operatorname{diag}(0,\Im z'\,I)$, so
$\Im\mathcal G=-\Im z'\,\mathcal G\operatorname{diag}(0,I)\mathcal G^*$ and
\[
 -\Im\bar v=\Im z'\,\E\langle\mathcal G_{12}\mathcal G_{12}^*\rangle\ge0,
 \qquad
 \Im(D-c\bar v-b\bar u)=-c\Im\bar v-b\Im\bar u\ge b\eta_0a_*^2 .
\]
Thus $D-c\bar v-b\bar u$ has an imaginary part bounded below, so it is
invertible with bounded inverse. Multiplying the top mean equation by
this inverse and using Lemma~\ref{lem:trace-matrix} for the covariance
errors gives
\[
 \bar v=\langle(D-c\bar v-b\bar u)^{-1}\rangle+O(W^{-2}),
 \qquad \bar u=(z'-b\bar v)^{-1}+O(W^{-2}).
\]
This is \eqref{eq:c01-cavity-map} up to errors $O(W^{-2})$, since
$f(\zeta)=\langle(\zeta-X)^{-1}\rangle$.

\emph{Step 2: the imaginary-part identities.} We now study an exact
solution of \eqref{eq:c01-cavity-map}. We set
\[
 y=-\Im v,\qquad w=-\Im u,\qquad A=\langle|D-cv-bu|^{-2}\rangle.
\]
In the eigenbasis of $D$, taking imaginary parts of the two equations
gives
\[
 y=A(cy+bw),\qquad w=|u|^2(\Im z'+by).
\]
Substituting the second equation into $cy+bw$, the first becomes
$y\bigl(1-A(c+b^2|u|^2)\bigr)=Ab|u|^2\Im z'$. Dividing by $y=A(cy+bw)$, we get
the gap identity
\[
 1-A(c+b^2|u|^2)=\frac{b(\Im z')|u|^2}{cy+bw}.
\]

\emph{Step 3: the gap is bounded below.} We bound every quantity in the
gap identity independently of $D$. By the Cauchy--Schwarz inequality,
$|v|^2\le A$. Since $y=A(cy+bw)\ge Acy$, we have $A\le c^{-1}$. Since
$-b\Im v\ge0$, we have $\Im(z'-bv)\ge\eta_0$, so $|u|\le\eta_0^{-1}$. And since
$|v|\le\sqrt A\le c^{-1/2}$, we have $|z'-bv|\le M+b/\sqrt c$. Together,
\[
 |v|^2\le A\le c^{-1},\qquad
 (M+b/\sqrt c)^{-1}\le|u|\le\eta_0^{-1}.
\]
In the gap identity the numerator is at least $b\eta_0(M+b/\sqrt c)^{-2}$,
and the denominator is at most $c\cdot c^{-1/2}+b\eta_0^{-1}$, since
$y\le|v|$ and $w\le|u|$. So the gap is bounded below by a positive constant
independent of $D$:
\[
 1-A(c+b^2|u|^2)\ge
 \frac{b\eta_0(M+b/\sqrt c)^{-2}}{\sqrt c+b/\eta_0}>0 .
\]

\emph{Step 4: existence and uniqueness.} For existence, we repeat the
eigenvalues of $D$ in successively larger dimensions, and we take a
subsequence of the corresponding approximate mean pairs. These satisfy
\eqref{eq:c01-cavity-map} up to errors that vanish in the limit. The
uniform bound on the denominator in Step~3 lets us pass to the limit,
which gives an exact solution with negative imaginary parts.

For uniqueness, we take two solutions and subtract their equations. Two
Cauchy--Schwarz inequalities bound the coefficient in the resulting
comparison by the geometric mean of the quantities $A(c+b^2|u|^2)$ for
the two solutions. By Step~3, this is strictly below one, so the two
solutions coincide. The same comparison bounds $|\bar v-v|$ by $CW^{-2}$, which is the
physical-trace bound of the statement.

\emph{Step 5: the mean error.} The approximate pair $(\bar v,\bar u)$ of
Step~1 has the same gap, up to $O(W^{-2})$. To see this, we use its port
lower bound $-\Im\bar u\ge\eta_0|\bar u|^2$ in the same imaginary
identities. We do not divide by $-\Im\bar v$, which might be small. The
comparison of Step~4 between the approximate pair and the exact solution
then has a coefficient strictly below one, and errors $O(W^{-2})$. This
proves $|\E(h_+(z')\mid X)-u|\le CW^{-2}$.

\emph{Step 6: fluctuations and derivatives.} The gradient with respect
to the fresh Gaussians, together with the growing inverse moments, gives
the fluctuation bound $C\sqrt p/W$ for $p\le a_0W$, exactly as in
\eqref{eq:fresh-trace-fluctuation}. For the derivatives in $z'$, we first
note that the solution $u$ is analytic in $z'$. This follows from the
nonzero stability denominator of Step~3 and the uniqueness of Step~4,
by the implicit function theorem. We take a small fixed circle around an
interior point, inside the upper half-plane. Cauchy's formula writes
each derivative as an integral over this circle, and Minkowski's
inequality moves the bounds of Step~5 and of the fluctuation inside the
integral. This gives both error bounds for the derivatives, after
decreasing $a_0$. None of these operations differentiates an error
estimate that has not been proved.
\end{proof}

 We start a path with nothing attached at its
far end, so the Schur complement term is $X_0=0$. This is the zero
boundary. After $J$ layers we obtain $X_J$. The  following lemma has three parts.
First, the empirical eigenvalue law of $X_J$ has a deterministic limit
$\mu_J$. Second, the Stieltjes transform of $\mu_J$ at the matched test
$z$ is exponentially close to the matched value $a$: the effect of the
zero boundary dies out at rate $e^{-c_*J}$. Third, for fixed $J$ the
normalized traces concentrate, and they depend continuously on a small
deterministic initial matrix $D$. This is the entrance from which the
matched analysis of the path starts.

\begin{lemma}[The actual zero-boundary entrance]\label{lem:c01-zero-entrance}
For fixed $J$, the empirical law of $X_J$, starting at $X_0=0$, converges
in probability to a deterministic probability measure $\mu_J$. At the
matched test $z$,
\begin{equation}\label{eq:c01-finite-seed}
 |G_{\mu_J}(z)-a|\le C e^{-c_*J}.
\end{equation}
For each fixed $J$, the normalized trace and its first $z$-derivative have
centered $L^p$ norm at most $C_J\sqrt p/W$ for $2\le p\le a_JW$. Their
means and these estimates are uniform in the compact parameter sets. If
the initial matrix is instead a deterministic Hermitian $D$ with
$\|D\|\le\epsilon$, their mean changes are at most $C_J\epsilon$, with the
same concentration bound uniform in $D$.
\end{lemma}

\begin{figure}[ht]
\centering
\begin{tikzpicture}[x=1cm,y=1cm]
 \foreach \x/\l in {0/1,1.2/2,2.4/3,4.8/J}{
   \draw[thick,fill=black!12] (\x,0) circle (.32);
   \node[font=\small] at (\x,0) {$m_{\l}$};
 }
 \node[font=\small] at (3.6,0) {$\cdots$};
 \draw[very thick,dashed] (6.15,0) circle (.53);
 \node[font=\small] at (6.15,0) {$m_{J+1}$};
 \draw[thick] (.32,0)--(.88,0) (1.52,0)--(2.08,0) (2.72,0)--(3.2,0)
   (4.0,0)--(4.48,0) (5.12,0)--(5.62,0);
 \draw[thick,red!70!black] (-1.2,-.4)--(-1.2,.4);
 \node[font=\small,red!70!black,align=center] at (-1.2,-.85) {$m_0=0$\\ zero boundary};
 \draw[thick] (-1.2,0)--(-.32,0);
 \node[font=\small,align=center] at (6.15,-1.05) {port\\ $z$};
 \fill[red!70!black,opacity=.25] (0,0) circle (.45);
 \draw[->,thick,blue!70!black] (0,.65) .. controls (2.5,1.5) and (4,1.5) .. (5.8,.5)
   node[pos=.5,above,font=\small,blue!70!black] {defect decays like $e^{-c_*J}$};
\end{tikzpicture}
\caption{The chain of scalar equations \eqref{eq:c01-finite-dyson} in the
proof of Lemma~\ref{lem:c01-zero-entrance}. Sites $1,\ldots,J$ are
physical blocks, coupled to their neighbours with weight $b$ and to
themselves with weight $c$. Site $J+1$ is the port, with parameter $z$.
The matched constant solution $a$ would require a neighbour with value
$a$ to the left of site~$1$. The actual zero boundary $m_0=0$ is a defect
at site~$1$ (shaded). The difference $m-a\mathbf1$ solves a linear
system with a source at site~$1$. Because that system is uniformly
invertible and its inverse decays exponentially, the defect reaches the
port only with weight $e^{-c_*J}$, which is \eqref{eq:c01-finite-seed}.}
\label{fig:ze-chain}
\end{figure}
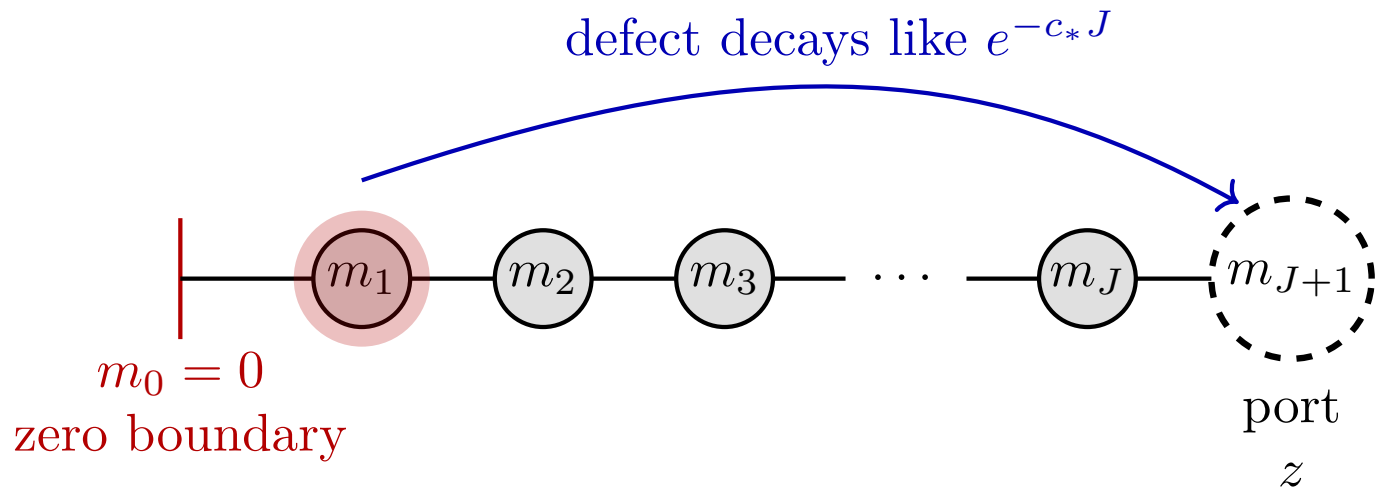

\begin{proof}
We first
show that the eigenvalue distribution of the Schur complement of a path that starts from zero boundary and crosses a fixed number of layers, converges, for
each fixed number of layers, to a deterministic limit (Step~1), described by a
finite chain of scalar equations (Step~2), drawn in
Figure~\ref{fig:ze-chain}. The zero boundary acts as a defect at the first
site of this chain. Comparing the chain with the matched constant solution,
we show that the defect is damped exponentially along the chain, using a
symmetrized form of the chain whose inverse decays exponentially
(Steps~3--7). Steps~8--10 prove the concentration and perturbation bounds by
writing each layer as a scattering matrix, a norm-preserving matrix that maps
incoming to outgoing waves, so that no large real inverse ever has to be
estimated.

\emph{Step 1: the limiting measures.} We obtain the limiting measures
from the global conditional map of Lemma~\ref{lem:c01-cavity-mean}, and
not from an assumed free equilibrium. We first need tightness of the
empirical laws. On the event $\|B\|\le K_0$, singular-value inequalities
and the count bound for deformed Gaussian matrices give
\begin{align*}
 \#\{|\lambda(X_+)|>R\}
 &\le\#\{|\lambda(A-e-X)|<K_0^2/R\},\\
 \Prob\{\#\{|\lambda(X_+)|>R\}\ge\varepsilon W\}
 &\le Ce^{-cW}+(C/(R\varepsilon))^{\lceil\varepsilon W\rceil}.
\end{align*}
The first line holds because $X_+=B^*(A-e-X)^{-1}B$, so a large
eigenvalue of $X_+$ needs a small eigenvalue of $A-e-X$. This gives
tightness, uniformly in the old Hermitian input.

The solution of \eqref{eq:c01-cavity-map} depends continuously on its
input probability measure, in the weak topology. Indeed, its means and
inverse denominators are uniformly bounded, and every subsequential
limit solves the same equations, whose solution is unique. We start
with $\mu_0=\delta_0$. Conditional concentration, together with a
countable set of Stieltjes tests, then proves the limit for each fixed
$J$ by induction.

\emph{Step 2: the finite chain of equations.} Nesting the scalar
equations of \eqref{eq:c01-cavity-map} over the $J$ layers gives numbers
$m_1,\ldots,m_{J+1}$ with negative imaginary parts that satisfy
\begin{equation}\label{eq:c01-finite-dyson}
 \begin{aligned}
 m_\ell^{-1}&=-e-cm_\ell-bm_{\ell-1}-bm_{\ell+1},
       &&1\le\ell\le J,\quad m_0=0,\\
 m_{J+1}^{-1}&=z-bm_J,\qquad&G_{\mu_J}(z)&=m_{J+1}.
 \end{aligned}
\end{equation}
For example, the outer $v$ is the preceding Stieltjes transform
evaluated at $-e-cv-bu$, and expanding that transform recursively gives
exactly the preceding row of this system. The same recursion selects
the solution with negative imaginary parts.

\emph{Step 3: comparison with the matched solution.} We now prove the
decay estimate \eqref{eq:c01-finite-seed}, globally. We write the system
in matrix form. Let $S$ be the tridiagonal matrix with diagonal entries
$c$ on the first $J$ sites, last diagonal entry zero, and adjacent
entries $b$. Then \eqref{eq:c01-finite-dyson} reads
\[
 m_\ell^{-1}=d_\ell-(Sm)_\ell,\qquad d=(-e,\ldots,-e,z).
\]
The constant vector $a\mathbf1$ solves the same system with $d_1$
replaced by $-e-ba$. Indeed, on an interior row this is
$a^{-1}=-e-(c+2b)a=-e-a$, on the first row the replacement supplies the
missing left neighbour, and on the last row it is $z-ba=a^{-1}$. Both
diagonal vectors have nonnegative imaginary parts. With
$y_\ell=-\Im m_\ell>0$ and $y_a=-\Im a$, taking imaginary parts gives
\[
 y_\ell=|m_\ell|^2(\Im d_\ell+(Sy)_\ell),\qquad
 y_a=\Im d_\ell^{\rm ref}+(S(y_a\mathbf1))_\ell,
\]
where the second identity uses $|a|=1$.

\emph{Step 4: bounds on $m$.} The self-interaction gives
$|m_\ell|\le c^{-1/2}$ for $\ell\le J$, since $y_\ell\ge|m_\ell|^2cy_\ell$. The port
gives $|m_{J+1}|\le(\Im z)^{-1}$, since $\Im(z-bm_J)\ge\Im z$. So the
denominators in \eqref{eq:c01-finite-dyson} are bounded above, and hence
\[
 0<r\le|m_\ell|\le R<\infty
\]
uniformly in $J$.

\emph{Step 5: the symmetrized matrix.} We put $R_m=\diag(\sqrt{|m_\ell|})$
and $A=R_mSR_m$. The two imaginary equations give $(Sy)_\ell\le y_\ell/|m_\ell|^2$
and $(S(y_a\mathbf1))_\ell\le y_a$. By the Cauchy--Schwarz inequality,
\[
 \diag(|m_\ell|)S\sqrt{yy_a}\le\sqrt{yy_a}
\]
entrywise. So the positive vector $R_m^{-1}\sqrt{yy_a}$ satisfies
$A(R_m^{-1}\sqrt{yy_a})\le R_m^{-1}\sqrt{yy_a}$. Since $A$ is real,
symmetric and nonnegative, this subinvariant vector shows that
$\|A\|\le1$. The first $J$ diagonal entries of $A$ are $c|m_\ell|\ge d_*=cr>0$.

We also need $A$ to stay away from $-1$. Splitting $A$ into its diagonal
and its nonnegative off-diagonal part, we have, for real $v$,
\[
 v^TAv\ge-|v|^TA|v|+2d_*\sum_{\ell\le J}|v_\ell|^2 .
\]
Suppose a unit eigenvector $v$ has eigenvalue $\lambda\le-1/2$. Since the
last diagonal entry is zero, its last row reads
$\lambda v_{J+1}=A_{J+1,J}v_J$, with $|A_{J+1,J}|\le1$. So
$|v_{J+1}|\le2|v_J|$, and at least one fifth of its mass lies in the first
$J$ sites. The inequality then gives $\lambda\ge-1+2d_*/5$. Consequently
\[
 \operatorname{Spec}A\subset[-1+\delta_*,1],\qquad
 \delta_*=\min(1/2,2d_*/5)>0 .
\]

\emph{Step 6: uniform invertibility and exponential decay.} We put
$\Theta=\diag(am_\ell/|m_\ell|)$, a diagonal unitary. The argument of $a$
stays away from $0$ and $-\pi$, and $\arg m_\ell\in(-\pi,0)$, so
$|1-\Theta_{\ell\ell}|\ge\delta>0$. Let $v$ be a unit complex vector and
$\epsilon=\|(I-\Theta A)v\|\le1$. Since $\Theta$ is unitary,
$\|Av\|\ge1-\epsilon$, and by the spectral bound of Step~5,
\[
 1-\|Av\|^2\le2\epsilon,\qquad
 \|(I-A)v\|^2\le\frac{2}{\delta_*}(1-\|Av\|^2)
                  \le4\epsilon/\delta_* .
\]
It follows that
\[
 \delta\le\|(I-\Theta)v\|\le\epsilon+2\sqrt{\epsilon/\delta_*},
\]
so $\epsilon$ is bounded below by a positive constant. This is a uniform
lower bound on the smallest singular value of $I-\Theta A$.

To get decay, we conjugate by $\diag(e^{t\ell})$. This changes the
tridiagonal matrix $I-\Theta A$ in norm by at most $2(e^{|t|}-1)$, since
only the off-diagonal entries are multiplied by $e^{\pm t}$. So the bound
on the inverse persists for small fixed $|t|$. This gives exponential
decay of the entries of the inverse, in both directions.

\emph{Step 7: the seed estimate.} Finally, the exact difference
$\Delta=m-a\mathbf1$ solves
\[
 (I-\diag(am_\ell)S)\Delta=-bm_1a^2e_1,\qquad
 I-\diag(am_\ell)S=R_m(I-\Theta A)R_m^{-1}.
\]
The first identity follows by subtracting the reference system, since
the two diagonal vectors differ only by $ba$ in the first entry. The
second follows from $\diag(am_\ell)=\Theta R_m^2$. The source sits at
site~$1$, so by Step~6 the last coordinate of $\Delta$ is at most
$Ce^{-c_*J}$. This proves \eqref{eq:c01-finite-seed}. In particular, the
argument does not assume that $z^{-1}$ is already close to $a$.

\emph{Step 8: one layer as a scattering matrix.} We now justify the
concentration and perturbation statements for the finite chain. We
avoid estimating real cavity inverses, which may be large. Instead we
describe each layer by its bounded unitary scattering matrix
\[
 \mathcal S=
 \begin{pmatrix}-I+2iR_0&-2R_0B\\-2B^*R_0&I-2iB^*R_0B\end{pmatrix},
 \quad R_0=(A-e+i(I+BB^*))^{-1},\quad \|R_0\|\le1.
\]
The inverse identity verifies directly that $\mathcal S$ is unitary.
Every derivative of its blocks has a fixed polynomial bound in
$1+\|A\|+\|B\|$.

To glue the layer to an independent unitary old reflector $U$, we need
the inverse of $I-Ur$, where $r=-I+2iR_0$. We put
\[
 F=i(I+U)^{-1}(I-U),\qquad D_0=A-e+iBB^*,\qquad M_0=D_0+F.
\]
Direct multiplication gives
\[
 (I-Ur)^{-1}
 =\tfrac12(D_0+iI)\{M_0^{-1}(I+iD_0)-iI\}.
\]
Lemma~\ref{lem:gaussian-absorption}, with its port equal to $iI$, bounds
all factors in $L^p$ uniformly in the Hermitian $F$ for $p\le aW$.
Continuity and Fatou's lemma cover the excluded eigenvalue $-1$ of $U$,
where $F$ is not defined. Thus the gluing estimate holds conditionally on
any independent old reflector. It does not apply to a reflector selected
from the fresh layer.

\emph{Step 9: concentration and perturbation after $J$ layers.} Starting
with $U_0=I$, which corresponds to $X_0=0$, the full reflector after $J$
layers determines the port resolvent by
\[
 (z'-X_J)^{-1}
 =(I+U_J)[(z'+i)I+(z'-i)U_J]^{-1},
\]
because $X=i(U-I)(I+U)^{-1}$. For $z'$ in a fixed compact subset of the
upper half-plane, this map and all its fixed derivatives on the unitary
group are bounded. Its Gaussian differential is a sum of a fixed number
of products of derivatives of scattering matrices and of gluing
inverses. Conditional bounds followed by H\"older's inequality control
their $L^p$ norms for $p\le a_JW$. We do not claim any independence
between the factors. So the gradient of the normalized trace has norm at
most $C_J/\sqrt W$ in $L^p$. The Gaussian covariance and the $L^p$
gradient inequality give the concentration bound $C_J\sqrt p/W$.
Cauchy's formula gives the estimate for the first $z'$-derivative.

For an initial deterministic $D$, the reflector is
$U_0(D)=(iI+D)(iI-D)^{-1}$. Along the path $tD$, $0\le t\le1$, its
derivative is bounded by $C\|D\|$. Differentiating the same fixed
products proves that the mean changes by at most $C_J\|D\|$. The
Gaussian bound is uniform in these deterministic $D$.

\emph{Step 10: uniformity in the parameters.} Differentiation in the
energy at fixed $J$ obeys the same finite-product bounds. We couple the
blocks as $A=\sqrt c\,A^{(0)}$ and $B=\sqrt b\,B^{(0)}$. Then
differentiation in $b$ also obeys these bounds, on compact subsets of
$(0,1/2)$. A finite grid in the parameters upgrades the convergence of
the means for fixed $J$ to uniform convergence on the compact parameter
sets. This completes every asserted uniformity.
\end{proof}

 We grow a path from one end, starting from a
small deterministic matrix $D$. The following lemma gives three inputs for the
later analysis. First, the normalized port trace $h_j(z)$ is within
$O(W^{-1})$ of the matched value $a$, once $j$ exceeds $j_0=(\log W)^2$
steps. Second, its derivative $h_j'(z)$ is within $O(W^{-1})$ of a
deterministic number $s_*$, the first jet. Third, the clocks
$\gamma_j^q$, which record how much source the path has accumulated, grow
at most linearly in $j$. The clock bounds hold without any stopping
time.

\begin{lemma}[One-ended first jet and unstopped clocks]
\label{lem:c01-open-inputs}
Let $h_j(\zeta)=W^{-1}\Tr(\zeta-X_j)^{-1}$ and $k=c+b^2a^2$. There is a
unique $s_*$ arising from the deterministic finite-chain limits, with
$|s_*|\le(1-b)^{-1}$ and
\[
 ks_*^2+(1+ca^2)s_*+a^2=0.
\]
Write $\delta_{0,j}=|h_j(z)-a|$ and $\delta_{1,j}=|h_j'(z)-s_*|$. For each
fixed finite moment family there is a small $\epsilon>0$ such that,
uniformly in deterministic Hermitian initial $X_0=D$ with
$\|D\|\le\epsilon$, and with $Y_0^q=0$, the following hold. For
$j_0=\lceil(\log W)^2\rceil$,
\begin{equation}\label{eq:t22-clock-input-moments}
 \|\delta_{0,j}\|_p+\|\delta_{1,j}\|_p\le C_p/W
       \quad(j_0\le j\le TW),\qquad
 \|\gamma_j^q\|_p\le C_{p,T,q}(j+1)/W
       \quad(0\le j\le TW).
\end{equation}
The scalar and jet bound also holds as $C\sqrt p/W$ for $2\le p\le a_*W$,
with a sufficiently small fixed $a_*>0$. The probability of an exit from
the scalar domain after a fixed entrance is at most $e^{-c_*W\log W}$.
The clock bounds are unstopped.
\end{lemma}

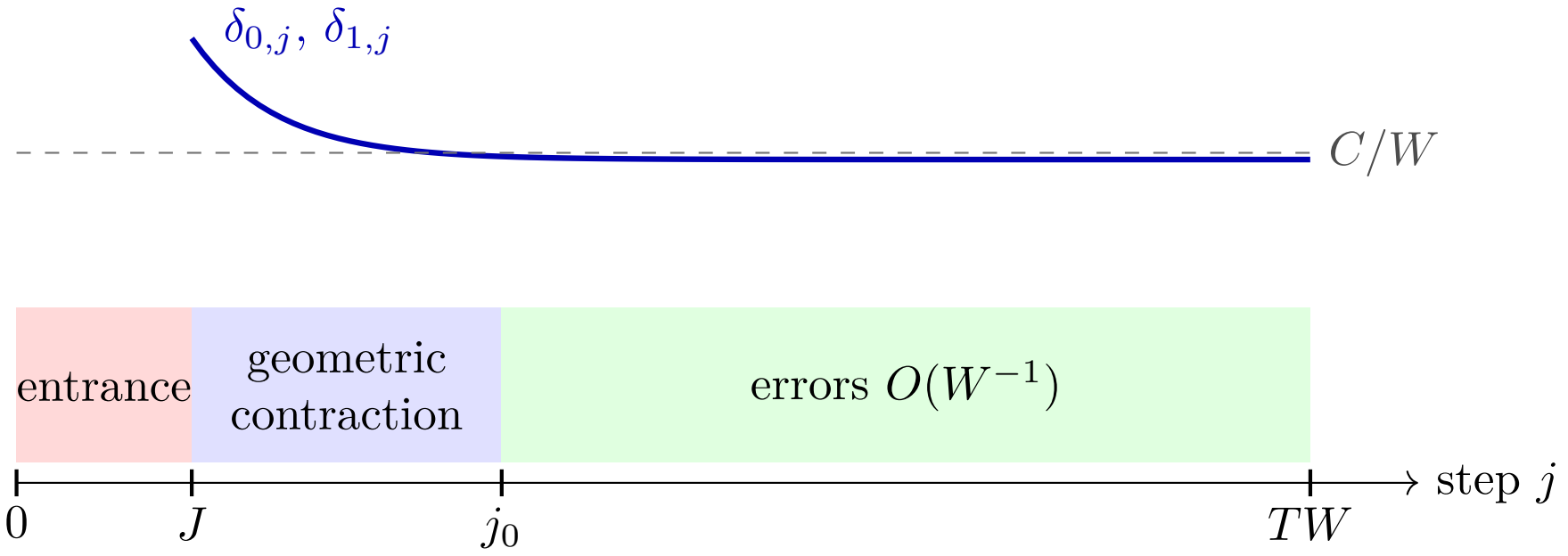
\begin{figure}[ht]
\centering
\begin{tikzpicture}[x=1cm,y=1cm]
 \draw[->] (0,0)--(10.4,0) node[right,font=\small] {step $j$};
 \foreach \x/\l in {0/0,1.3/J,3.6/j_0,9.6/TW}{
   \draw[thick] (\x,-.1)--(\x,.1) node[below=4pt,font=\small] {$\l$};
 }
 \fill[red!15] (0,.15) rectangle (1.3,1.3);
 \fill[blue!12] (1.3,.15) rectangle (3.6,1.3);
 \fill[green!12] (3.6,.15) rectangle (9.6,1.3);
 \node[font=\small,align=center] at (.65,.72) {entrance};
 \node[font=\small,align=center] at (2.45,.72) {geometric\\ contraction};
 \node[font=\small,align=center] at (6.6,.72) {errors $O(W^{-1})$};
 \draw[very thick,blue!70!black,domain=1.3:9.6,samples=100]
   plot (\x,{2.4+0.9*exp(-1.6*(\x-1.3))});
 \node[font=\small,blue!70!black,anchor=west] at (1.4,3.35) {$\delta_{0,j}$, $\delta_{1,j}$};
 \draw[dashed,black!55] (0,2.45)--(9.6,2.45)
   node[right,font=\small,black!70] {$C/W$};
\end{tikzpicture}
\caption{The three phases in the proof of Lemma~\ref{lem:c01-open-inputs},
drawn schematically. During the entrance, the first $J$ steps from the
zero boundary, Lemma~\ref{lem:c01-zero-entrance} brings the scalar error
within a small fixed radius. From step $J$ on, each layer contracts the
scalar error by $q_*<1$ and the first-jet error by $q_1^2<1$, up to fresh
errors of size $\sqrt p/W$. After $j_0=\lceil(\log W)^2\rceil$ steps, the
initial error has been multiplied by $q^{(\log W)^2}$, which is smaller
than every power of $W$. So from $j_0$ to $TW$ only the fresh errors
remain, and $\delta_{0,j}+\delta_{1,j}=O(W^{-1})$ in $L^p$.}
\label{fig:oi-timeline}
\end{figure}

\begin{proof}
The proof follows the three phases of
Figure~\ref{fig:oi-timeline}. In a fixed entrance phase, the previous lemma
brings the boundary trace close to its matched value (Step~1). After that,
each layer shrinks the error of the trace by a fixed factor and adds a small
fresh error, so the error stays of the order of the inverse block size
(Step~2). The same contraction holds for the derivative of the trace in the
port parameter, the first jet, whose limit solves a fixed-point equation
(Steps~3--5). Finally we bound the clocks, which record the accumulated
source, through a recursion for their moments with a weight that absorbs the
scalar error (Steps~6--7).

\emph{Step 1: the entrance.} We choose the local scalar radius as small
as needed below. By Lemma~\ref{lem:c01-zero-entrance} we choose a fixed
$J$ whose zero-boundary mean is within one eighth of that radius. We
choose $\epsilon$ so small that the sensitivity $C_J\epsilon$ to the initial
matrix is also within one eighth. The concentration for the finite
chain, with $p=\lfloor a_JW\rfloor$, makes the probability that the
entrance fails at most $(C/\sqrt W)^{\lfloor a_JW\rfloor}$.

\emph{Step 2: the scalar recursion.} From step $J$ on, we use the local
comparison for one block, whose coefficient is
\eqref{eq:strict-contraction}, together with
Lemma~\ref{lem:c01-cavity-mean}. On the scalar domain they give
\[
 \delta_{0,j+1}\le q_*\delta_{0,j}+C/W^2+|\xi_{j+1}|,\qquad
 \|\xi_{j+1}\|_{p\mid\mathcal F_j}\le C\sqrt p/W.
\]
So each step multiplies the old error by $q_*<1$ and adds a fresh error
$\xi_{j+1}$. To leave the scalar domain for the first time, a step needs
a fresh error of fixed size. The moment bound on $\xi_{j+1}$ and a union
bound over the at most $TW$ steps therefore give an exit probability at
most $e^{-c_*W\log W}$, uniformly in $D$. Before the exit, Minkowski's
inequality and the geometric sum give
\[
 \|\delta_{0,j}\|_p\le Cq_*^{j-J}+C\sqrt p/W.
\]
The trace is deterministically bounded by $(\Im z)^{-1}$. So after
reducing $a_*$, the exit event contributes only $O(W^{-1})$ for
$p\le a_*W$.

\emph{Step 3: the equation for the first jet.} For the derivative in the
port parameter, we differentiate \eqref{eq:c01-cavity-map}. With
$\zeta=-e-cv-bu$ we have $v'=f'(\zeta)\zeta'$, $\zeta'=-cv'-bu'$, and
$u'=-u^2(1-bv')$. Solving these three linear equations gives
\begin{equation}\label{eq:c01-first-jet}
 u'=-u^2\frac{1+cf'(\zeta)}{1+(c+b^2u^2)f'(\zeta)},\qquad
 \zeta'=\frac{bu^2}{1+(c+b^2u^2)f'(\zeta)}.
\end{equation}
On the scalar domain, $\zeta=z+O(\delta_0)$ and
$|f'(z)|\le(1-b)^{-1}+O(\delta_0)$. The derivatives of Stieltjes
transforms at a fixed height are bounded deterministically. Since also
$u=a+O(\delta_0)$, the outgoing first jet is therefore
$F(f'(z))+O(\delta_0)$, where
\[
 F(s)=-a^2\frac{1+cs}{1+ks},\qquad k=c+b^2a^2 .
\]

\emph{Step 4: the jet map is a contraction.} A direct computation gives
$F'(s)=a^4b^2/(1+ks)^2$, so, since $|a|=1$,
\[
 |F'(s)|=\frac{b^2}{|1+ks|^2}\le q_1^2<1
       \qquad\bigl(|s|\le(1-b)^{-1}+\epsilon_1\bigr).
\]
To see the bound, note that $|1+ks|\ge(1-b-|k|)/(1-b)$ for
$|s|\le(1-b)^{-1}$. So it follows from $b(1-b)/(1-b-|k|)<1$. This holds
uniformly in the bulk, because $|k|<c+b^2=(1-b)^2$ whenever $a^2\ne1$. We
choose $\epsilon_1$ small. The disk is convex, so $F$ is Lipschitz with
constant $q_1^2$ between any two of its points. We do not need $F$ to map
the whole disk into itself.

The first jets of the deterministic measures $\mu_J$ follow this
contraction, with exponentially small scalar errors by
\eqref{eq:c01-finite-seed}. Their successive differences are therefore
summable, so they converge to a limit $s_*$. The contraction gives the
uniqueness of $s_*$ in the stated disk. Moreover, $s_*$ is a fixed point,
$s_*=F(s_*)$, and multiplying out $s_*(1+ks_*)=-a^2(1+cs_*)$ gives the
quadratic equation in the statement.

\emph{Step 5: the actual first jet.} The conditional derivative errors
of Lemma~\ref{lem:c01-cavity-mean} give the same geometric recursion for
the actual first jet, now with the forcing term
$C\delta_{0,j}+C\sqrt p/W$. Its deterministic bound pays for the exits.
Since $q^{(\log W)^2}$ is smaller than every power of $W$, the initial
errors are negligible after $j_0$ steps. This proves the scalar and
first-jet assertions in \eqref{eq:t22-clock-input-moments}.

\emph{Step 6: a recursion for the clock moments.} For the clocks we put
$u_j=W^{-1}\Tr Y_j^q$. The fresh identities of
Lemma~\ref{lem:weighted-layer} apply to one block with the present
independent Hermitian past. In particular, $\|L\|$ and $\|C_+\|$ have
uniform conditional fixed moments. For integer $p\ge2$ this gives
\[
 \E(u_+^p\mid\mathcal F_j)
 \le(1+C_p\delta_{0,j}+C_p/W)u_j^p
                    +C_{p,q}(u_j^{p-1}+1).
\]
To see this, we expand $u_+^p$ about the conditional mean. The linear
centered term vanishes. The terms of degree at least two cost
$C_p(1+u_j)^p/W$, by \eqref{eq:clock-noise}. Crucially, the
Cauchy--Schwarz inequality gives
\[
 \E(|\xi_+|u_+^p\mid\mathcal F_j)\le C_p(1+u_j^p)/W
\]
without assuming independence between the scalar error and the clock.

The term $C_p\delta_{0,j}u_j^p$ could accumulate over many steps, so we
absorb it with a weight. We choose $K_p(1-q_*)>2C_p$, and we shrink the
scalar radius if necessary. Multiplying by $1+K_p\delta_{0,j+1}$, the
contraction of $\delta_0$ in Step~2 absorbs the predictable scalar error.
This gives
\[
 M_{j+1}\le(1+C_p/W)M_j+C_pM_j^{1-1/p}+C_p,\qquad
 M_j=\E[\mathbf1_{\rm no\ exit}(1+K_p\delta_{0,j})u_j^p].
\]

\emph{Step 7: iteration and the unstopped bound.} The finitely many
entrance steps have bounded clock moments, by the unrestricted fresh
norm bounds. Iterating the recursion of Step~6 then gives
$M_j\le C_{p,T,q}(j+1)^p$, since $M_j^{1/p}$ grows by at most a constant
per step. Without stopping, the same fresh bounds give a crude $2p$-th
moment at most $C_p^{j+1}$. By H\"older's inequality and the exit
probability $e^{-c_*W\log W}$, these pay the entire exceptional
contribution up to $TW$, so the clock bound holds without stopping.
Real fixed $p\ge1$ follow by increasing to an integer. Positivity
handles a general bounded $q\ge0$. In particular,
$\E\gamma_{j_0}^q\le Cj_0/W=o(1)$, which is the bound on the initial
clock required by the determinant equation.
\end{proof}

These are estimates under the complete one-ended prior, conditional
only on an independent initial $D$ when present. No eigenvector,
spectral root, matched far absorber or assumption of scalar equilibrium
is imposed.

Let $g_W\in\mathbb C^W$ be standard circular complex Gaussian independent
of the entire chain. Write $Y_{n;k}$ for the positive contribution of
the injection $C_k$ transported to time $n$, so that
$Y_n^q=\sum_{k\le n}q(k/W)Y_{n;k}$. Define the finite positive measure
\[
 \Sigma_W=\frac{1}{W^2}\sum_{k=1}^n(g_W^*Y_{n;k}g_W)\delta_{k/W}.
\]

\begin{theorem}[Actual open-source Brownian profile law]
\label{thm:t22-open-clock-law}
If $n/W\to\tau\in(0,T]$, then $\Sigma_W$ converges in law, weakly as a
finite measure on $[0,T]$, to
\begin{equation}\label{eq:t22-open-Brownian-measure}
 \Sigma_\tau(dx)=d\exp\left\{\sqrt{D_0}B_{\tau-x}
                    -\frac{D_0}{2}(\tau-x)\right\}
                    \mathbf1_{[0,\tau]}(x)\,dx,
\end{equation}
where $B$ is standard real Brownian motion. In particular the normalized
open-source profile converges to $\Sigma_\tau/\Sigma_\tau([0,\tau])$.
All fixed positive moments of bounded source integrals converge.

This limit has a strictly positive, continuous density. Every real
continuous nonconstant test on $[0,\tau]$ has strictly positive variance
under its normalized profile. After rescaling the interval to $[0,1]$,
that profile's mean-square deviation from Lebesgue measure is $O(\tau)$
as $\tau\downarrow0$. These are assertions about the stated open-source
law, not about the target's original-ring spectral marks.
\end{theorem}

The law is specified by a positive Brownian exponential functional. This is a probe law for the open clock; physical ring normalization is handled separately.

\subsection{The source gap for trace-class clock weights}\label{subsec:s21-the-source-gap-for-trace-class-clock-wei}
 We must mix the entire positive clock, rather than a fixed list of
vectors chosen in advance. Lemma~\ref{lem:t22-clock-mixing} bounds its
weighted scalar error by $O(W^{-1/2})$. We condition on the complete
past, diagonalize the old clock, and iterate the source contraction
while retaining its trace weights and centered fluctuations.

 The following extension keeps the whole positive clock. It does not
substitute a future-selected frame into a finite-rank theorem.
For $Q=Q_j^q$ put
\[
 Z_j=\operatorname{Tr}(Qf(X_j)),\qquad
 S_j=\operatorname{Tr}(Qf(X_j)Qf(X_j)^*)
        =\|Q^{1/2}f(X_j)Q^{1/2}\|_{\rm HS}^2.
\]
\begin{lemma}[Complete clock source mixing]\label{lem:t22-clock-mixing}
For every fixed finite $p\ge1$, uniformly on $j_0\le j\le TW$,
\begin{equation}\label{eq:t22-clock-mixing}
 \|Z_j\|_p+\|\sqrt{S_j}\|_p\le C_{p,T,q}W^{-1/2}.
\end{equation}
\end{lemma}

The weighted scalar error contracts along the fresh layers. The retained clock factor ensures that this estimate controls the full source response.

\begin{proof}
Condition on the complete past, diagonalize its positive $Q$, and put
$T_+=L f(X_+)L^*$. Lemma~\ref{lem:passive-source-mixing}, at zero absorption, gives
\[
 \mathbb E(T_+\mid\mathcal F)=\theta f(X)+R,
 \quad |\theta|<1,\quad
 \|R\|\le C(\delta_0+\delta_1+W^{-2}),
\]
and every centered bilinear entry has conditional $L^p$ norm
$CW^{-1/2}\|u\|\|v\|$, uniformly in the past, by
Lemma~\ref{lem:source-recorded-estimates}(a). If the eigenvalues of $Q$ are $q_i\ge0$,
Minkowski applied to
$\sum_{ij}q_iq_j|(T_+-\mathbb ET_+)_{ij}|^2$ consequently gives
\[
 \left\|\|Q^{1/2}(T_+-\mathbb ET_+)Q^{1/2}\|_{\rm HS}
                      \right\|_{p\mid\mathcal F}
 \le C_pW^{-1/2}\operatorname{Tr}Q.
\]
The same bound holds for its trace. It is independent of the number
of nonzero eigenvalues of $Q$.

Without the injection, $Z_+=\operatorname{Tr}(QT_+)$ and
$\sqrt{S_+}=\|Q^{1/2}T_+Q^{1/2}\|_{\rm HS}$. The mean contributes
at most $|\theta|\sqrt S+C(\delta_0+\delta_1+W^{-2})\gamma$.
The injection is $B=q((j+1)/W)C_+/W^2$. Its contribution to $Z_+$ is
$\operatorname{Tr}(Bf(X_+))$, of modulus at most
$\|q\|_\infty\|f\|\,\|C_+\|/W$. It therefore has conditional $L^p$ size
$C_p/W$, because $\|C_+\|$ has uniform conditional fixed moments
(Step~6 of the proof of Lemma~\ref{lem:c01-open-inputs}). For $A=L^*QL$ the exact identity
\[
 S(A+B,f)=S(A,f)+S(B,f)+
        \operatorname{Tr}(AfBf^*)+\operatorname{Tr}(BfAf^*)
\]
has nonnegative cross terms bounded by
$2\|B\|\|f\|^2\operatorname{Tr}A$. Thus the extra conditional $L^p$
cost in $\sqrt{S_+}$ is at most
$C_p(W^{-1}\sqrt\gamma+W^{-3/2})$.
Use \eqref{eq:t22-clock-input-moments} and H\"older on the products
with the jet errors. Taking full $L^p$ norms yields geometric recursions
with coefficient $|\theta|<1$ and forcing $C_pW^{-1/2}$.
At $j_0$, the trivial bounds $|Z|,\sqrt S\le C\gamma$ cost
$C_pj_0/W=o(W^{-1/2})$.
Iteration proves \eqref{eq:t22-clock-mixing}. The scalar exceptional
sets are removed with the unstopped moments and arbitrary failure powers;
$f$ itself has a deterministic bounded operator norm.
\end{proof}

\subsection{A closed determinant-Laplace equation}\label{subsec:s21-a-closed-determinant-laplace-equation}
 We reduce the matrix-valued clock law to a scalar transform equation.
Lemma~\ref{lem:t22-det-clock} gives the evolution of
$\mathbb E\det(I+aQ)^{-1}$ with a summable error. We differentiate this
determinant, use the complete-clock mixing estimate to close its
one-step expansion, and identify the limiting equation by a Brownian
representation.

 For a positive matrix $Q$ and $a>0$ define
\[
 F_a(Q)=\det(I+aQ)^{-1},\qquad A_a=aQ(I+aQ)^{-1}.
\]
Then $0<F_a\le1$, $\operatorname{Tr}A_a\le a\gamma$, and
\begin{equation}\label{eq:t22-det-derivatives}
 a^2\partial_a^2F_a
 =F_a\left[(\operatorname{Tr}A_a)^2+\operatorname{Tr}A_a^2\right].
\end{equation}

\begin{lemma}[Critical determinant evolution]\label{lem:t22-det-clock}
Uniformly for $a$ in any compact subset of $(0,\infty)$,
\begin{equation}\label{eq:t22-det-clock-step}
 \mathbb E F_a(Q_{j+1}^q)-\mathbb E F_a(Q_j^q)
 =\frac{1}W\left[\frac{D_0a^2}{2}\partial_a^2
              \mathbb E F_a(Q_j^q)
       -d q((j+1)/W)a\mathbb E F_a(Q_j^q)\right]+r_{W,j}(a),
\end{equation}
where $\sum_{j=j_0}^{\lfloor TW\rfloor}\sup_a|r_{W,j}(a)|\to0$.
The derivatives are ordinary derivatives of the displayed finite
expectation. No determinant in this statement is a spectral coarea
weight; it is the exact Laplace transform of the independent probe.
\end{lemma}

The determinant test is bounded and its derivatives have explicit trace bounds. These properties allow the discrete generator remainder to be summed before removing stops.

\begin{proof}
We display the contractions, because an absolute $O(W^{-1})$ error per
layer would not be sufficient at criticality. Write $U=LL^*-I$ and
$U^0=U-\mathbb E(U\mid\mathcal F)$. The one-ended layer is the fresh layer preceding \eqref{eq:passive-mean},
with Hermitian $X$, the same $z$ and the same $f(X)$, and with zero
absorption, so that the absorption terms in $q/W$ there vanish. Then
\eqref{eq:passive-cancelled-mean} gives, with $\omega=-\Omega$, where
$\Omega=O(\delta_0)$ on the stable branch,
\[
 \mathbb E(U\mid\mathcal F)
   =-\omega f(X)-\overline\omega f(X)^*
                     +O(\delta_0^2+W^{-2}),\qquad |\omega|\le C\delta_0.
\]
Lemma~\ref{lem:passive-covariance} at zero absorption, where its $T$ and $q$ vanish, with the Hermitian imaginary parts taken,
gives for any positive past-measurable $A$,
\begin{align}
 \mathbb E\operatorname{Tr}(AU^0AU^0\mid\mathcal F)
 &=\frac{D_0}{W}(\operatorname{Tr}A)^2+
 O\!\left(\frac{|\operatorname{Tr}(Af)|^2}{W}
       +\left[\frac{\delta_0+\delta_1}{W}+\frac{1}{W^2}\right](\operatorname{Tr}A)^2\right),
                                                        \label{eq:t22-cov-trace-square}\\
 \mathbb E[(\operatorname{Tr}AU^0)^2\mid\mathcal F]
 &=\frac{D_0}{W}\operatorname{Tr}A^2+
 O\!\left(\frac{\operatorname{Tr}(AfAf^*)}{W}
       +\left[\frac{\delta_0+\delta_1}{W}+\frac{1}{W^2}\right](\operatorname{Tr}A)^2\right).
                                                        \label{eq:t22-cov-square-trace}
\end{align}
These follow by diagonalizing $A$ and summing the four-vector identities;
the error sums to $(\operatorname{Tr}A)^2$, without a dimension factor.
The opposite-sign coefficient is $D_0=\sigma_+/(2(\operatorname{Im}z)^2)$.
The same-sign contractions have $H_-=-2i\operatorname{Im}z,f(X)$,
which yields the two displayed source errors. Their bounded scalar
coefficient need not have a limit. Perturbing $H_+$ and the
opposite-sign scalar coefficient costs $O(\delta_0)$.

For $A=A_a$, positivity and $A_a\preceq aQ$ give
\[
 \operatorname{Tr}(A_afA_af^*)\le a^2 S,
 \qquad
 |\operatorname{Tr}(A_af)|
 \le a|Z|+a^2\gamma\sqrt S.
\]
The second inequality uses $A_a=aQ-a^2Q^2(I+aQ)^{-1}$ and
Hilbert--Schmidt Cauchy--Schwarz. Lemma~\ref{lem:t22-clock-mixing},
its higher fixed moments and \eqref{eq:t22-clock-input-moments}
therefore make every source error in the two covariance formulas
summable to $o(1)$. The Ward mean contributes
$\mathbb E|\omega\operatorname{Tr}(A_af)|=O(W^{-3/2})$.
Its quadratic terms and its $O(\delta_0^2+W^{-2})$ remainder are
also summable to $o(1)$. The $\delta_1/W$ remainders of
Lemma~\ref{lem:passive-covariance} contribute at most
$Ca^2\delta_1\gamma^2/W$ per layer; by
\eqref{eq:t22-clock-input-moments} and H\"older their sum over
$j_0\le j\le TW$ is $O(W^{-1})$.

Here is a uniform Taylor-remainder justification. The determinant lemma
expresses the transported ratio as
\[
 \frac{F_a(L^*QL)}{F_a(Q)}
   =\det(I+A_a^{1/2}UA_a^{1/2})^{-1}.
\]
Centered bilinear moments and Minkowski give conditional $L^p$ bounds
$C_pW^{-1/2}\operatorname{Tr}A_a$ for both the Hilbert--Schmidt norm
of the centered matrix and its trace. The mean is bounded by
$C\delta_0\operatorname{Tr}A_a$. On the event where its full
Hilbert--Schmidt norm is at most $1/2$ and its trace has modulus at
most one, the logarithm and exponential Taylor expansions through
second order have remainder bounded by
\[
 C F_a(Q)\left(\|A_a^{1/2}UA_a^{1/2}\|_{\rm HS}^3
                      +|\operatorname{Tr}(A_aU)|^3\right).
\]
The exponential factor is bounded on this event. Its mean is
$O(W^{-3/2})$ using the clock moments. Off this event the two original
Laplace functions are bounded by one. Markov at a fixed $p>6$, then
H\"older with the clock and jet moments, bounds the omitted contribution
and the omitted Taylor polynomial by $o(W^{-1})$, summably.
This argument requires no exponential moment of $\gamma$.
The surviving second-order terms are one half of
\eqref{eq:t22-cov-trace-square} plus one half of
\eqref{eq:t22-cov-square-trace}, multiplied by $F_a(Q)$.
Equation~\eqref{eq:t22-det-derivatives} identifies their sum.

Finally restore the positive injection $B=q C_+/W^2$.
For positive $B$ the Laplace function's second directional derivative
is at most $2a^2(\operatorname{Tr}B)^2$. Its first derivative at
$L^*QL$ is
$-aF_a(L^*QL)\operatorname{Tr}[(I+aL^*QL)^{-1}B]$.
Replacing the inverse by $I$ costs at most
$Ca^2W^{-2}\|C_+\|\operatorname{Tr}(L^*QL)$ in mean.
Replacing $F_a(L^*QL)$ by $F_a(Q)$ in this term costs $O(W^{-3/2})$
by the same first-order bounds and fresh moments. The proved
injection mean is
$\mathbb E\operatorname{Tr}C_+/W=d+O(\delta_0+W^{-2})$. For the injection mean write $\mathcal G=\mathcal G_{j+1}$, so that
$C_+=\mathcal G_{12}^*\mathcal G_{12}$. The inverted matrix has
imaginary part $\operatorname{diag}(0,\eta I)$, $\eta=\Im z$, so the
exact Ward identity
$\Im\mathcal G=-\mathcal G\operatorname{diag}(0,\eta I)\mathcal G^*$
gives $\operatorname{Tr}C_+=-\eta^{-1}\Im\operatorname{Tr}\mathcal G_{11}$.
Lemma~\ref{lem:c01-cavity-mean} gives
$\mathbb E(W^{-1}\operatorname{Tr}\mathcal G_{11}\mid\mathcal F)=v+O(W^{-2})$,
where $(v,u)$ solves \eqref{eq:c01-cavity-map} with $z'=z$ and $X=X_j$.
With $a$ the matched value of
Subsection~\ref{subsec:s21-the-one-ended-path} (not the Laplace variable),
since $-e-(c+b)a=z$ and $(z-ba)^{-1}=a$, $(a,a)$ solves these equations up to the error $h_j(z)-a$, and the
stable comparison gives $v=a+O(\delta_0)$. Finally $|a|=1$ gives
$\eta=\Im(\bar a+ba)=-(1-b)\Im a$, so $-\Im a/\eta=d$. Hence
\[
 \mathbb E(\operatorname{Tr}C_+/W\mid\mathcal F)=d+O(\delta_0+W^{-2}).
\]
Thus the injection gives exactly the displayed killing term, with
summable remainder. All constants above are uniform on compact
$a$-intervals. This proves the lemma.
\end{proof}

\begin{proof}[Proof of Theorem~\ref{thm:t22-open-clock-law}]
Let $f_W(t,a)=\mathbb E F_a(Q_{\lfloor Wt\rfloor}^q)$, interpolated
linearly in $t$. The clock moments, the determinant derivative formulas,
and Lemma~\ref{lem:t22-det-clock} give equicontinuity on compact
$(t,a)$ domains, up to an error vanishing uniformly in $t$.
At the entrance, $|1-f_W(j_0/W,a)|\le a\mathbb E\gamma_{j_0}^q=o(1)$.
Every subsequential limit is bounded between zero and one and solves
\begin{equation}\label{eq:t22-clock-pde}
 \partial_t f=\frac{D_0a^2}{2}\partial_a^2 f-dq(t)af,
                    \qquad f(0,a)=1.
\end{equation}
For precision, multiply the discrete equation by a smooth compactly
supported function of $a$, integrate twice by parts, and sum in time.
Uniform moment bounds control the derivative terms; the summed remainder
tends to zero. This gives the weak equation. The operator is uniformly
parabolic on compact $a$-intervals, so its bounded weak solution has the
usual interior regularity. The bounded solution on $a>0$ is unique:
stop the geometric Brownian diffusion at $[K^{-1},K]$, apply the ordinary
Feynman--Kac formula there, and let $K\to\infty$. Its exit probability
tends to zero on each finite time interval, while the nonnegative
killing factor and $0\le f\le1$ bound all boundary terms.
Consequently
\begin{equation}\label{eq:t22-clock-laplace}
 f(t,a)=\mathbb E\exp\left\{-ad\int_0^t q(t-s)
                  e^{\sqrt{D_0}B_s-D_0s/2}\,ds\right\}.
\end{equation}
No interchange of a critical limit with an unbounded observable entered
this determination.

The exact complex Gaussian integral gives
$\mathbb E_{g_W}e^{-a g_W^*Q_n^qg_W}=F_a(Q_n^q)$.
Thus \eqref{eq:t22-clock-laplace} is the Laplace functional of the
actual measure $\Sigma_W$, and it has the unique limit
\eqref{eq:t22-open-Brownian-measure}. Tightness follows from the bounded
first moments of total mass on the compact interval. Positive continuous
tests determine finite positive measures, so every subsequence has the
same law. One may equivalently apply the scalar Laplace limit to every
positive linear combination of any finite family of nonnegative tests.
The limit's total mass is strictly positive almost surely, so normalization
is continuous at its sample measures.

The conditional Gaussian identity
$\mathbb E_g(g^*Qg)^p\le p!(\operatorname{Tr}Q)^p$ for integer $p$
and the higher clock moments supply uniform integrability of every fixed
power; hence moments of bounded source integrals converge.
Brownian paths are continuous and their exponentials positive, proving
the density assertion. After reversing the spatial interval, Brownian measure has full support
on continuous paths starting at zero. For a nonconstant test $\phi$, the functional
$h\mapsto\int\phi e^h/\int e^h$ is nonconstant even on that path space:
use $h=\epsilon(\phi-\phi(0))$ and differentiate at zero, obtaining
the strictly positive Lebesgue variance of $\phi$. Neighborhoods of two
such paths have positive probability. This proves strict variance.
The small-$\tau$ claim follows from the $L^p$ supremum bounds and
exponential moments of $\sqrt{D_0\tau}B_u-D_0\tau u/2$ on $[0,1]$,
by expanding the normalized exponential or applying the mean value theorem.
\end{proof}

The comparison in Section~\ref{sec:clock-lower-tail} starts the one-ended
chain from a small deterministic cavity; the law is the same, uniformly
in that cavity.

\begin{corollary}[Open-clock law from a small initial cavity]
\label{cor:t22-open-clock-initial}
Fix $T<\infty$, and let $\epsilon>0$ be given by
Lemma~\ref{lem:c01-open-inputs} for the finitely many moments used in the
proofs of Lemmas~\ref{lem:t22-clock-mixing} and~\ref{lem:t22-det-clock} and
of Theorem~\ref{thm:t22-open-clock-law}. For deterministic Hermitian $D$
with $\norm D\le\epsilon$, let $\Sigma^D_W$ be the measure $\Sigma_W$ of the
one-ended chain started from $X_0=D$, with $Y^q_0=0$. If
$n_W/W\to\tau\in(0,T]$ and $n_W\le TW$, then for every bounded continuous
$q\ge0$ on $[0,T]$ and every bounded continuous $\phi$,
\[
 \lim_{W\to\infty}\ \sup_{D=D^*,\ \norm D\le\epsilon}
 \Bigl|\mathbb E\phi\Bigl(\int q\,d\Sigma^D_W\Bigr)
       -\mathbb E\phi\Bigl(\int q\,d\Sigma_\tau\Bigr)\Bigr|=0 .
\]
For $q=1$ the limit variable is
$\Sigma_\tau([0,\tau])=d\int_0^\tau e^{\sqrt{D_0}B_s-D_0s/2}\,ds$.
\end{corollary}

\begin{proof}
Lemma~\ref{lem:c01-open-inputs} gives \eqref{eq:t22-clock-input-moments},
the exit bound and the unstopped clock moments with constants that do not
depend on $D$. Every other estimate in the proofs of
Lemmas~\ref{lem:t22-clock-mixing} and~\ref{lem:t22-det-clock} is a
fresh-layer estimate uniform in the past. Write
$f^D_W$ and $r^D_{W,j}$ for $f_W$ and the remainder in
\eqref{eq:t22-det-clock-step} for the chain started from $D$. Then
$\sup_D\sum_{j=j_0}^{\lfloor TW\rfloor}\sup_a|r^D_{W,j}(a)|\to0$, the
equicontinuity in the proof of Theorem~\ref{thm:t22-open-clock-law} is
uniform in $D$, and
$\sup_D|1-f^D_W(j_0/W,a)|\le a\sup_D\mathbb E\gamma^q_{j_0}\to0$.
If the assertion fails, there are $\delta>0$ and matrices $D_W$ along a
subsequence with difference at least $\delta$. Along it $f^{D_W}_W$ has
the same equicontinuity and summed remainder, so every further
subsequential limit solves
\eqref{eq:t22-clock-pde} with $f(0,a)=1$. By uniqueness it equals
\eqref{eq:t22-clock-laplace}, and equicontinuity in $t$ passes from
$\lfloor W\tau\rfloor$ to $n_W$. Since
$\mathbb E e^{-a\int q\,d\Sigma^{D_W}_W}=\mathbb EF_a(Q^q_{n_W})$, the
Laplace transforms of $\int q\,d\Sigma^{D_W}_W\ge0$ converge for every
$a>0$ to that of $\int q\,d\Sigma_\tau$, so these variables converge in
law, a contradiction.
\end{proof}
\section{Lower tails of the full clock}
\label{sec:clock-lower-tail}
When we normalize a profile we divide by the full two-input clock, so we
need to know that it is not too small. This section proves this for every
critical limit. It is used for local mass flatness in Part~III.
Theorem~\ref{thm:t22-negative-clock} states that the full clock is tight,
that every limit is strictly positive, and that its negative moments of every
order are bounded by $C\tau^{-p}$. No inverse moment is claimed before the
limit.
We resolve the incoming directions by a phase average
(Lemma~\ref{lem:t22-phase-resolution}). We then compare an interval of
boundary conditions with the one-ended clock of
Section~\ref{sec:open-clock-law} (Lemma~\ref{lem:t22-robin-entrance}), and
turn its Laplace bound into lower-tail estimates. Finally we record a
small-ball bound (Lemma~\ref{lem:t26-clock-small-ball}).

 Here $Y$ denotes the full two-input interior clock of Section~\ref{sec:full-source}. All constants are uniform for $E$ in a fixed compact bulk set, at
fixed $b$, and $0<\tau\le T_0$, and the statements of this section hold
along deterministic $E_W\to E\in(-2,2)$, with limiting objects at $E$.
Their finite-width inputs are the open-clock law of
Section~\ref{sec:open-clock-law} (Lemmas~\ref{lem:c01-cavity-mean},
\ref{lem:c01-zero-entrance}, \ref{lem:c01-open-inputs},
\ref{lem:t22-clock-mixing} and~\ref{lem:t22-det-clock},
Theorem~\ref{thm:t22-open-clock-law} and
Corollary~\ref{cor:t22-open-clock-initial}),
\eqref{eq:full-clock-moments}, and the exact phase resolution of
Lemma~\ref{lem:t22-phase-resolution}. The estimate \eqref{eq:full-clock-moments} and the passive inputs
of Lemmas~\ref{lem:passive-source-mixing}
and~\ref{lem:passive-covariance} are uniform in the energy by
Lemma~\ref{lem:source-recorded-estimates}. The open-clock law has
constants uniform for $e$ in a compact set, and its limit depends
continuously on $e$ through $D_0(e)$; so along $E_W$
Corollary~\ref{cor:t22-open-clock-initial} holds with the limit at $E$,
by the contradiction argument of its proof. The phase interval of
Lemma~\ref{lem:t22-robin-entrance} is chosen at $E$ with a
$W$-independent margin, so that it serves for $E_W$ for large $W$.

\begin{lemma}[Exact scalar phase resolution]
\label{lem:t22-phase-resolution}
Let the full unitary scattering matrix of a finite path of $n\geq 1$ be
$S=\left(\begin{smallmatrix}r&t'\\t&r'\end{smallmatrix}\right)$.
Close the far port by the scalar unitary $e^{i\theta}I$ and put
\[
 \Gamma_\theta=(I-e^{i\theta}r')^{-1}e^{i\theta}t,
 \qquad A_\theta=\binom{I}{\Gamma_\theta}.
\]
For the actual finite path, almost surely,
\begin{equation}\label{eq:t22-phase-resolution}
 \int_0^{2\pi}A_\theta A_\theta^*\frac{d\theta}{2\pi}=I_{2W}.
\end{equation}
If $Y$ is its full two-input interior clock, the one-input interior
clock after this closure is $A_\theta^*YA_\theta$.
\end{lemma}

This identity averages the full incoming space exactly, including interference between the two ports. Invertibility of the finite Gaussian bonds is used before any limit.

\begin{proof}
All Gaussian bonds are invertible almost surely. Successive block
inversion then makes the finite-path transmission $t$ invertible;
unitarity gives $tt^*=I-r'r'^*\succ0$, so $\|r'\|<1$.
The uniformly convergent geometric expansion gives
$\int\Gamma_\theta=0$ and
\[
 \int\Gamma_\theta\Gamma_\theta^*\frac{d\theta}{2\pi}
 =\sum_{j\ge0}(r')^j(I-r'r'^*)(r'^*)^j=I.
\]
The finite partial sum is $I-(r')^{N+1}(r'^*)^{N+1}$, which tends
to $I$. This proves the identity, without assuming that $r'$ is normal.
The scattering equations say the full incoming vector is
$A_\theta g$ when the near incoming vector is $g$. The internal
solution is consequently $UA_\theta g$, proving the clock identity.
\end{proof}

\begin{lemma}[A controlled Robin entrance]
\label{lem:t22-robin-entrance}
For each fixed $b,E$ there is a positive-measure phase interval
 $I=I_{b,E}$, independent of $W$, such that uniformly in $\theta\in I$
the Gaussian-probe interior
clock after this closure, divided by $W^2$, converges to
\[
 X_\tau=d\int_0^\tau
       \exp\{\sqrt{D_0}B_s-D_0s/2\}\,ds,
 \quad d=(1-b)^{-1},\quad D_0=2/[b(4-E^2)].
\]
The same assertion holds with any prescribed bounded continuous
nonnegative spatial source. The interval may be chosen separately for
a fixed set of moment bounds used in the argument; its measure has a
uniform positive lower bound for $E$ in a compact bulk set, at fixed $b$.
\end{lemma}

Only a fixed positive-measure interval of boundary phases is needed. Its uniformity permits comparison with the full clock without controlling every real boundary condition.

\begin{proof}
Write $z=\Re z+i\eta$ for the matched auxiliary diagonal. Its scattering
relation is outgoing $=$ incoming $-2i\eta$ times the auxiliary
amplitude. The condition incoming $=e^{i\theta}$ outgoing therefore
replaces this diagonal by the real number
\[
 \lambda_\theta=z+\frac{2i\eta e^{i\theta}}{1-e^{i\theta}}
              =\Re z-\eta\cot(\theta/2).
\]
Choose $I$ with $\lambda_\theta\in[\Lambda,2\Lambda]$ for a sufficiently
large fixed $\Lambda$. Eliminating the far auxiliary coordinate leaves
the initial Hermitian cavity $X_0=B_0^*\lambda_\theta^{-1}B_0$ and
zero initial interior clock. The remaining fresh layers and near
matched absorber are exactly those of Section~\ref{sec:open-clock-law}.

The clocks are identified exactly. Let $\psi$ be the solution of the
closed path with near incoming vector $g$, let $\psi_j$ be its amplitude on
the $j$th physical block and $p$ its far auxiliary amplitude. With the
diagonal $z$ restored, $\psi$ solves the open equations with far incoming
amplitude $(z-\lambda_\theta)p$, and its far outgoing amplitude is
$(z-\lambda_\theta)p-2i\eta p=(\bar z-\lambda_\theta)p$. The choice of
$\lambda_\theta$ is the closure condition
$e^{i\theta}(\bar z-\lambda_\theta)=z-\lambda_\theta$ of
Lemma~\ref{lem:t22-phase-resolution}, which determines the full incoming
vector uniquely as $A_\theta g$, since $\norm{r'}<1$. Hence $(\psi_j)_j$ is the interior part of $UA_\theta g$, and
\[
 W^{-2}g^*A_\theta^*YA_\theta g=W^{-2}\sum_{j=1}^n\norm{\psi_j}^2 .
\]
On the other hand, eliminating $p$ leaves the one-ended chain started from
$X_0$. Its new block amplitude is $(\mathcal G_n)_{12}g$, and $L_ng$ is the
incoming vector of the chain of the first $n-1$ blocks. The recursion
$Y^q_j=q(j/W)C_j+L_j^*Y^q_{j-1}L_j$ with
$C_j=(\mathcal G_j)_{12}^*(\mathcal G_j)_{12}$ therefore gives, by
induction, $g^*Y^q_ng=\sum_jq(j/W)\norm{\psi_j}^2$. So the closed clock divided by $W^2$ equals the open clock $g^*Q^1_ng$
of Section~\ref{sec:open-clock-law}, with no factor depending on $\eta$
or $\lambda_\theta$, also with a spatial weight $q$. In particular the constant $d=(1-b)^{-1}$ in $X_\tau$ is the constant of
\eqref{eq:t22-open-Brownian-measure}, which comes from the injection mean
$\mathbb E\operatorname{Tr}C_+/W=d+O(\delta_0+W^{-2})$ in the proof of
Lemma~\ref{lem:t22-det-clock}, and $X_\tau=\Sigma_\tau([0,\tau])$.

Lemma~\ref{lem:c01-zero-entrance} supplies a fixed entrance length $J$
and its uniform $C_J\|D\|$ perturbation estimate, including the first
cavity-spectral jet. Lemma~\ref{lem:c01-open-inputs} then gives the
scalar and jet errors $C_p/W$ after $\lceil\log^2W\rceil$, the
superexponential exit bound, and the unstopped clock moments uniformly
over $\|D\|\le\epsilon$. The discarded clock has expected scaled
mass $O(\log^2W/W)$. The clock source differentiates only the interior
physical diagonal: $\lambda_\theta$, $B_0$ and thus $D$ are held fixed.
Consequently the determinant equation in
 Theorem~\ref{thm:t22-open-clock-law} has the same initial value $1$
and the same coefficients, uniformly for these $D$. This is
Corollary~\ref{cor:t22-open-clock-initial}, which also covers the spatial
weights $q$.

The initial bond is independent of the remaining layers. Choose a fixed
Gaussian norm cutoff $K$ and then $\Lambda\ge K^2/\epsilon$. On
$\|B_0\|\le K$ all these deterministic-input conclusions apply,
uniformly in $\theta\in I$.
Since $d\theta/d\lambda=2\eta/((\lambda-\Re z)^2+\eta^2)$,
the interval corresponding to $[\Lambda,2\Lambda]$ has measure bounded
below uniformly on the compact parameter sets. The complementary event has probability
tending to zero exponentially in $W$. Bounded Gaussian Laplace
functions cost at most this probability there. Hence their limit is
Theorem~\ref{thm:t22-open-clock-law}'s unchanged Brownian functional. More precisely, $B_0$ is independent of the remaining layers, so the
closed clock law is the mixture over $B_0$ of the open laws started from
$X_0(B_0)$. For bounded continuous $\phi$ its distance from
$\mathbb E\phi(X_\tau)$ is therefore at most the supremum in
Corollary~\ref{cor:t22-open-clock-initial} plus
$2\norm\phi_\infty\mathbb P(\norm{B_0}>K)$, uniformly in $\theta\in I$. More precisely, $B_0$ is independent of the remaining layers, so the
closed clock law is the mixture over $B_0$ of the open laws started from
$X_0(B_0)$. For bounded continuous $\phi$ its distance from
$\mathbb E\phi(X_\tau)$ is therefore at most the supremum in
Corollary~\ref{cor:t22-open-clock-initial} plus
$2\norm\phi_\infty\mathbb P(\norm{B_0}>K)$, uniformly in $\theta\in I$. 
\end{proof}

\begin{theorem}[Negative moments of every limiting full clock]
\label{thm:t22-negative-clock}
Let $h\sim\mathrm{CN}(0,I_{2W})$ be independent of the actual path,
and $T_W=W^{-2}h^*Yh$. The family is tight. Every subsequential limit
$T$ is strictly positive and, for every fixed $p>0$,
\begin{equation}\label{eq:t22-negative-clock}
 \mathbb E T^{-p}\le C_{p,T_0}\tau^{-p}.
\end{equation}
This is an estimate on limits; uniform prelimit inverse moments are
not asserted.
\end{theorem}

The limiting clock can therefore normalize a Gaussian response without creating a zero denominator. The order of limits avoids a stronger finite-width inverse-moment estimate.

\begin{proof}
Tightness follows from \eqref{eq:full-clock-moments} and Gaussian conditioning. For any $s>0$ and positive
$Q=Y/W^2$, set
$F_Q(C)=\det(I+sQ^{1/2}CQ^{1/2})^{-1}$, $C\succeq0$.
It is decreasing in Loewner order and convex: it is exactly
$\mathbb E_g\exp(-s g^*Q^{1/2}CQ^{1/2}g)$, an average of convex
exponentials of linear functionals. Let $a=|I|/(2\pi)>0$. By
\eqref{eq:t22-phase-resolution}, $\int_I A_\theta A_\theta^*d\theta/(2\pi)
\preceq I$. Monotonicity and Jensen give the pointwise inequality
\[
 F_Q(I)\le\frac{1}a\int_I F_Q(aA_\theta A_\theta^*)\frac{d\theta}{2\pi}.
\]
Take the path expectation, then a critical subsequence. Bounded Laplace
convergence and Lemma~\ref{lem:t22-robin-entrance} show
\[
 \mathbb E e^{-sT}\le\mathbb E e^{-saX_\tau}\quad(s>0).
\]
Letting $s\to\infty$ first proves $\mathbb P(T=0)=0$.
Integrating this inequality against $s^{p-1}ds/\Gamma(p)$ gives
$\mathbb ET^{-p}\le a^{-p}\mathbb EX_\tau^{-p}$ by Tonelli.
Finally
\[
 X_\tau\ge d\tau\exp\{-\sqrt{D_0}\sup_{s\le\tau}|B_s|-D_0\tau/2\},
\]
and the Brownian supremum has all fixed exponential moments uniformly
for $\tau\le T_0$. This proves the stated bound.
\end{proof}
 \begin{lemma}[Rescaled limiting clock and trace small balls]
\label{lem:t26-clock-small-ball}
In the actual Gaussian path law of Theorem~\ref{thm:t22-negative-clock},
let $n/W=\tau_W\to\tau\in(0,T_0]$, and, for an independent standard
$h\in\C^{2W}$, put
\[
 A_W=\frac{\Tr Y}{nW},\qquad
 M_W=\frac{h^*Yh}{nW}=\frac{T_W}{\tau_W}.
\]
Every subsequential limit $M$ is strictly positive and satisfies
$\E M^{-p}\le C_{p,T_0}$ for each fixed $p>0$.
For $0<\varepsilon<1$,
\begin{equation}\label{eq:t26-clock-small-ball}
 \limsup_{W\to\infty,\,n/W\to\tau}\Prob(A_W\le\varepsilon)
 \le\frac{C_{p,T_0}\varepsilon^{p/2}}{1-\sqrt\varepsilon}.
\end{equation}
The constants are uniform in $\tau\in(0,T_0]$ and, at fixed $b$, for
$E$ in a fixed compact bulk set. The width limit precedes the small-ball limit.
No finite-width negative moment is used.
\end{lemma}
\begin{proof}
Since $T_W=\tau_WM_W$, tightness and Slutsky's theorem identify
$T=\tau M$ along every convergent subsequence. Consequently
$\E M^{-p}=\tau^p\E T^{-p}\le C_{p,T_0}$, and $M>0$ almost surely.
Conditional on the entire path (and thus on $Y$),
$\E_hM_W=A_W$. Markov's inequality gives the finite-width statement
\[
 \one_{\{A_W\le\varepsilon\}}
 \Prob_h(M_W\le\sqrt\varepsilon\mid Y)
 \ge(1-\sqrt\varepsilon)\one_{\{A_W\le\varepsilon\}}.
\]
Integrating yields
$(1-\sqrt\varepsilon)\Prob(A_W\le\varepsilon)
\le\Prob(M_W\le\sqrt\varepsilon)$.
Take a subsequence realizing the limsup and a further subsequence along
which $M_W$ converges. The closed-set part of Portmanteau, followed by
the limiting inverse-moment bound, gives
\[
 \limsup_W\Prob(M_W\le\sqrt\varepsilon)
 \le\Prob(M\le\sqrt\varepsilon)
 \le\varepsilon^{p/2}\E M^{-p}.
\]
This proves the claim, including the stated order of limits.
\end{proof}
\section{The weighted physical closing distribution}
\label{sec:collar-density}
To close the path into a ring, we need the law of the fresh three-block
collar as seen from its boundary. This section bounds the weighted boundary
kernel of the collar and identifies its mean, even when the old real Schur
complement is large. It prepares the closing law of
Section~\ref{sec:t22-analytic-collar}.
Theorem~\ref{thm:t22-collar-kernel} bounds the expected positive boundary
kernel of the collar. Theorem~\ref{thm:t22-exact-boundary-mean} shows that its
mean tends to the identity, so that no spectral mass escapes to infinity.
Corollary~\ref{cor:t22-physical-one-channel-Cauchy} shows that single boundary
entries have the standard Cauchy law.
We first represent dissipation by a real Cauchy average
(Lemma~\ref{lem:t22-dilation-variance}). In Subsection~\ref{subsec:s23-weighted-fresh-spectral-densities} we prove
spectral density bounds that survive multiplication by large exterior factors
(Lemmas~\ref{lem:global-weighted-density}
and~\ref{lem:t26-relative-matrix-density}). In Subsection~\ref{subsec:s23-an-exact-identity-for-the-boundary-kerne} we derive
an exact identity for the boundary kernel and prove
Theorem~\ref{thm:t22-collar-kernel}. In Subsection~\ref{subsec:s23-the-limiting-mean-and-scalar-boundary-la} we identify the
limiting mean and the scalar boundary law.

\begin{lemma}[Dilation and its exact matrix variance]
\label{lem:t22-dilation-variance}
Let $X=X_h-iY$ with $X_h=X_h^*$ and $Y\succeq0$, and let
$z=\xi+i\eta$, $\eta>0$. For a standard real Cauchy variable $t$ put
\[
 R=(z-X_h+iY)^{-1},\qquad R_t=(z-X_h-tY)^{-1}.
\]
Then, in both matrix orders,
\begin{align*}
 \mathbb E_t R_t&=R,\\
 \mathbb E_t(R_t-R)^*(R_t-R)&=\eta^{-1}R^*YR,\\
 \mathbb E_t(R_t-R)(R_t-R)^*&=\eta^{-1}RYR^*.
\end{align*}
If the dimension is at most $c_0W$, $\|D\|\le1$, and
$g=4\eta\operatorname{Tr}R^*YR$, then for every fixed $k\ge1$,
\[
 \mathbb E_t\left|W^{-1}\operatorname{Tr}D(R_t^k-R^k)\right|^2
       \le C_{k,c_0,\eta}\,g/W.
\]
Consequently a fixed finite list of normalized resolvent jets changes
by at most $\epsilon$ outside a phase set of Cauchy mass at most
$C g/(W\epsilon^2)$.
\end{lemma}

A scalar Cauchy average represents the dissipative boundary parameter exactly. The weighted estimate controls the exceptional real phases of that representation.

\begin{proof}
The matrix resolvent is bounded and analytic for $\Im t<0$;
its norm is at most $\eta^{-1}$. The scalar Cauchy integral evaluates
it at $t=-i$. This proves the first identity, also for singular $Y$
by a positive regularization or bounded analytic continuation.
For Hermitian $X_h+tY$, its Ward identity is
$-\Im R_t=\eta R_t^*R_t=\eta R_tR_t^*$.
For $R$, it is
$-\Im R=R^*(\eta I+Y)R=R(\eta I+Y)R^*$.
Taking the expectation and subtracting the squared mean proves both
variance formulas. For the jet estimate use
\[
 R_t^k-R^k=\sum_{j=0}^{k-1}R_t^j(R_t-R)R^{k-1-j},
\]
so its Hilbert--Schmidt norm is at most
$k\eta^{1-k}\|R_t-R\|_{\rm HS}$. Trace Cauchy--Schwarz and the
matrix variance identity give the assertion. Ordinary derivatives in
$z$ are factorial multiples of these powers. A union bound and
Chebyshev finish the proof. No normality of $X$ was used.
\end{proof}

\subsection{Weighted fresh spectral densities}\label{subsec:s23-weighted-fresh-spectral-densities}

 We need a spectral density bound that survives multiplication by the
large exterior factors from Schur elimination.
Lemma~\ref{lem:global-weighted-density} gives such a bound for arbitrary
deformations, with a sharper version on the regular class. We average
the full fresh matrix first, then sum the indexed sectors and compare
the resulting weighted matrix density.

 For a Hermitian deformation $F$, let $\Sigma_F(e)$ be the expected
matrix spectral density of $F+G$, where $G$ is   $\Naux $ -dimensional GUE
with entry variance   $c/\Naux $. Thus $\int_I\Sigma_F(e)\,de
=\E\one_I(F+G)$. A bound on this matrix after multiplying by $F-e$
will retain the weights arising from Schur complementation.

\begin{lemma}\label{lem:global-weighted-density}
For every Hermitian $F$, real $e$, and fixed $\eta>0$,

\begin{equation}\label{eq:global-weighted-density}
 (F-e+i\eta)\Sigma_F(e)(F-e-i\eta)\preceq C_{c,\eta}\sqrt \Naux  I.
\end{equation}

On the fixed regular polynomial deformation class
\eqref{eq:fresh-regularity}, the stronger bound $CI$ holds.
The first assertion has no norm bound on $F$. Moreover $\Sigma_F(e)\preceq C_c\sqrt\Naux I$. More generally, let
$\Naux=kW$ and let $G'=\bigoplus_{i=1}^kG_i$ have independent blocks
$G_i\sim\GUE_W(c_i/W)$ with $c_i\ge c_0>0$. For every Hermitian $F$,
possibly random and independent of $G'$, the expected matrix spectral
density $\Sigma'_F(e)$ of $F+G'$ satisfies
\begin{equation}\label{eq:k-block-density}
 \Sigma'_F(e)\preceq C_{k,c_0}\sqrt W I.
\end{equation}
\end{lemma}

The bound permits arbitrary deformations before a regularity event is imposed. Its dimension dependence is used only to pay the exceptional set.

\begin{proof}
For a deterministic unit vector $v$, conjugate $G$ so that $v$ is
its first coordinate and condition on every entry except $G_{11}$.
Rank-one spectral averaging says
$\int\langle v,\delta_e(H+t vv^*)v\rangle\,dt=1$.
It follows either by Schur inversion and integrating the scalar
Poisson kernel, or first at positive height and then by positivity.
The density of $G_{11}$ is bounded by   $C\sqrt \Naux $, giving
$\Sigma_F(e)\preceq C\sqrt \Naux  I$.
For $k$ independent Gaussian diagonal blocks of size $W$, the same argument proves \eqref{eq:k-block-density}: split $v$ into its block components and
use $|\sum_i v_i|_D^2\le k\sum_i|v_i|_D^2$ for a positive matrix $D$.
All off-diagonal entries can be conditioned on.

At an eigenvalue $e$, multiplication by $F-e+i\eta$ equals
multiplication by $-G+i\eta$. Therefore
\[
 (F-e+i\eta)\Sigma_F(e)(F-e-i\eta)
       =\E[(G-i\eta)\delta_e(F+G)(G+i\eta)].
\]
Write $G=G_1+G_2$ with independent variance-$c/2$ halves.
Positivity bounds the quadratic form by three terms, with outside
vectors $G_1v,G_2v,i\eta v$. Conditional averaging over $G_2$
bounds the first by   $C\sqrt \Naux \|G_1v\|^2$ ; then
$\E G_1^2=(c/2)I$. Exchange the halves for the second term.
This proves \eqref{eq:global-weighted-density}, including its
pointwise version by scalar Gaussian smoothing.

On the regular class, Theorem~\ref{thm:fresh-palm} and
Lemma~ \ref{lem:t22-point-density} give
$\Sigma_F(e)\preceq C\Gamma_F$, with
$\Gamma_F=c|F-e-cm_c(e)|^{-2}$.
The scalar self-energy and its reciprocal imaginary part are bounded.
Diagonalizing $F$ gives
\[
 \sup_{\lambda\in\R}
 \frac{c((\lambda-e)^2+\eta^2)}{|\lambda-e-cm_c(e)|^2}\le C,
\]
which proves the regular bound without discarding the weights.
\end{proof}

\begin{lemma}[Complete sector sum and relative matrix density]
\label{lem:t26-relative-matrix-density}
On the fixed regular polynomial deformation class, uniformly for $x$
in a smaller fixed regularity interval,

\begin{equation}\label{eq:t26-relative-matrix-density}
 -\epsilon_\Naux \Gamma_F(x)\preceq
 \Sigma_F(x)-r_{F,c}(x)\Gamma_F(x)
 \preceq\epsilon_\Naux \Gamma_F(x),\qquad \epsilon_\Naux \longrightarrow0.
\end{equation}

This is an operator inequality with no dimension factor from
polarization. It can be multiplied on both sides by any deterministic
rectangular matrix, regardless of its norm.
\end{lemma}

\begin{proof}
The exact scalar Gaussian decomposition gives, as a matrix density,

\[
 \Sigma_F(x)=\sum_{\ell=1}^\Naux \E[w_{\ell,x}u_\ell u_\ell^*],
 \qquad \sum_{\ell=1}^\Naux Z_{\ell,x}=\Naux \rho_{F,c}(x).
\]

Both identities use all ordered sectors and their original scalar
variance   $c/\Naux ^2$. For an arbitrary deterministic complex vector $q$,
sum the $r=1$ assertion of Theorem~\ref{thm:fresh-palm} and divide
by   $\Naux $. It gives

\[
 |q^*(\Sigma_F(x)-\rho_{F,c}(x)\Gamma_F(x))q|
 \le C\{\Naux ^{-\kappa_1}\rho_{F,c}(x)+\Naux ^{-K}\}
                                      q^*\Gamma_F(x)q.
\]

In particular the sum of additive sector errors is   $\Naux ^{-K}$  after
this division. Lemma~\ref{lem:t22-point-density} bounds the scalar
density and replaces it by $r_{F,c}(x)$ with uniform error tending
to zero. The resulting bound holds for every $q$, with the same
constant. Since the difference is Hermitian, this is exactly
\eqref{eq:t26-relative-matrix-density}. It is not necessary to take
a union over a vector net or to bound individual matrix entries.
Complex polarization, if used to recover entries, is performed only
after this quadratic-form inequality. Congruence preserves the two
Loewner bounds, proving the last assertion.
\end{proof}

 \subsection{An exact identity for the boundary kernel}\label{subsec:s23-an-exact-identity-for-the-boundary-kerne}

 We express the closing kernel in a form that keeps its exterior
weights attached. The exact identity below leads to the uniform collar
density estimate and Theorem~\ref{thm:t22-collar-kernel}. We factor the
boundary inverse, write the actual three-block collar, and apply the
weighted fresh-density bounds before averaging the remaining Gaussian
variables.

 For a contraction $V$ define its positive Poisson kernel by
$\mathcal P(V)=(I-V)^{-1}+(I-V^*)^{-1}-I$.
For Hermitian $K_c,H_o$, $Y\succ0$, and $\eta>0$, put
\[
 Z=H_o+i(\eta I+Y),\quad r=I-2i\eta Z^{-1},\quad
 U=I-2i\eta(K_c+i\eta)^{-1},\quad R_c=(K_c+H_o+iY)^{-1}.
\]
Then, with the displayed matrix orders,
\begin{equation}\label{eq:weighted-boundary-identity}
 \mathcal P(Ur)
 =\eta^{-1}(I-ZR_c)Y(I-R_c^*Z^*)
 =\eta^{-1}(K_c-i\eta)R_cYR_c^*(K_c+i\eta).
\end{equation}
Indeed
$I-Ur=2i\eta(K_c+i\eta)^{-1}(K_c+H_o+iY)Z^{-1}$.
Invert this factorization and use $-\Im R_c=R_cYR_c^*$.
These identities also explain why one cannot multiply an unweighted
small error by an unbounded exterior reaction after averaging.

Let $c=1-2b>0$, $z_b=\xi+i\eta$, $z_0=-e-ca$, and
$a=(-e-i\sqrt{4-e^2})/2$ as in Section~\ref{sec:t22-common-port-profile}. Constants retain fixed compact
interior $b,e$ sets. A collar consists of three consecutive physical
blocks: two fresh end diagonals $A=\operatorname{diag}(A_1,A_2)$,
a fresh middle diagonal $G$, all of GUE variance $c/W$, and the two
fresh bonds $B=(B_1^T,B_2^T)^T$ of variance $b/W$.
All these variables are mutually independent. Its terminating reaction
and scattering matrix, on the two outer ports, are
\[
 K_c=A-e-\xi-B(G-e)^{-1}B^*,\qquad
 U=I-2i\eta(K_c+i\eta)^{-1}.
\]
This is the actual finite Gaussian collar, not a full GUE substitution.

\begin{lemma}[Uniform weighted collar density]
\label{lem:t22-collar-density}
For any Hermitian $2W$ matrix $X$, let $\Psi_X$ be the expected
top-top matrix spectral density at zero of
\[
 \mathcal H_X=\begin{pmatrix}A-e-X&B\\B^*&G-e\end{pmatrix}.
\]
With $K_X=z_b-X$, one has
\begin{equation}
 K_X\Psi_XK_X^*\preceq C\sqrt W I.                              \label{eq:collar-global-density}
\end{equation}
There are fixed $d_*,C>0$ such that the stronger bound $CI$ holds
whenever
\begin{equation}
 \max_{i=1,2}\left|W^{-1}\operatorname{Tr}
                    [(z_b-X)^{-1}]_{ii}-a\right|\le d_* . \label{eq:collar-regularity}
\end{equation}
Neither assertion requires a norm bound on $X$.
\end{lemma}

The three-block collar keeps the weights generated by the old response. Regularity improves the density bound on the overwhelmingly likely part of that law.

\begin{proof}
At a zero eigenvector $(v,w)$ of $\mathcal H_X$,
$K_Xv=(z_b+e-A)v-Bw$. Thus the weighted spectral projector
in \eqref{eq:collar-global-density} is obtained by multiplying the full projector by the fresh
row $(z_b+e-A,-B)$. Split the full fresh Gaussian collar into two independent halves
and its outside row into the two fresh rows plus the bounded deterministic
part. Conditional on one half, the other has the three-block matrix
spectral-density bound $C\sqrt W I$ proved in
Lemma~\ref{lem:global-weighted-density}. Each outside fresh row has
expected row square at most $CI_{2W}$. Positivity and averaging over the
other half therefore give \eqref{eq:collar-global-density}, uniformly in $X$.
This keeps the complete finite matrix before taking the expectation.

Here are the regular and exceptional parts of \eqref{eq:collar-regularity}.
Conditional on $A,B$, put $D=A-e-X$, $T=-D^{-1}B$,
$F=-e-B^*D^{-1}B$, and $H=(z_0-B^*D^{-1}B)^{-1}$.
The top density, averaged first over $G$, is exactly
$T\Sigma_F(0)T^*$. This follows by Schur inversion and spectral
averaging in the fresh $G$; equivalently first give its bottom block
a positive height and then take the density boundary limit.
The effective spectral parameter's positive derivative cancels in
that conditional density. No physical energy-clock weight is dropped.

The common-port calculation of Lemma~\ref{lem:layer-contraction}, \eqref{eq:fresh-trace-fluctuation}, and then Lemma~\ref{lem:regularity-certificate} give the fixed fresh regularity rectangle for
$F$ outside an arbitrarily small polynomial-probability set when
$d_*$ is sufficiently small. The required polynomial norm of $F$
also has such a cutoff: the two-block deformed Gaussian Wegner bound
gives $\mathbb P(\|D^{-1}\|>W^M)\le CW^{1-M}$ uniformly in $X$,
and $B$ has all fixed operator moments. This Wegner bound is Theorem~\ref{thm:external-count} with two fresh blocks. For the certificate at $\mu=\bar a$, take the complex conjugate of the lower-half-plane trace residual; Hermiticity gives $c|z_0-\Sigma|^{-2}=c|\bar z_0-\Sigma|^{-2}$. The exponent $M$ is chosen
after the desired exceptional probability. No bound on $\|X\|$ enters.

On this event Theorems~\ref{thm:fresh-palm} and Lemma~\ref{lem:t22-point-density} gives
\[
 \Sigma_F(0)\preceq Cc\,HH^*,\qquad
 K_XT HH^*T^*K_X^*=LL^*,\quad L=K_X\mathcal G_{12},
\]
with the absorbed $\mathcal G$ of Section~\ref{sec:t22-common-port-profile}. Its uniform fresh operator
moments and conditional mean bound therefore give $CI$.
To pay the exceptional set, apply Lemma~\ref{lem:global-weighted-density} to the bottom full GUE, with outside factor
$F-c\Re a+i\Im z_0=F-ca$. The harmless bounded real shift is included
in the same spectral-projector proof. It gives, for every $A,B$,
\[
 K_XT\Sigma_F(0)T^*K_X^*\preceq CW\,LL^*.
\]
Lemma~\ref{lem:t22-common-port-mean}'s uniform higher operator moments of $L$ and H\"older make its
exceptional expectation $o(1)$, by choosing the polynomial power first.
This proves \eqref{eq:collar-regularity} with the complete normalizers and both incident bonds.
\end{proof}

\begin{theorem}[Actual collar Poisson-kernel bound]
\label{thm:t22-collar-kernel}
Let $r=I-2i\eta(z_b-X)^{-1}$ with $X=X_h-iY$, $Y\succeq0$,
be independent of the fresh collar. Put
$g=\operatorname{Tr}(I-r^*r)$ and
$d(X)=\max_i|W^{-1}\operatorname{Tr}[(z_b-X)^{-1}]_{ii}-a|$.
For $d(X)\le d_*/2$,
\begin{equation}
 \mathbb E_{\rm collar}\mathcal P(Ur)
       \preceq C(1+g)I,\qquad
 \mathcal P(V)=(I-V)^{-1}+(I-V^*)^{-1}-I.                 \label{eq:collar-poisson-bound}
\end{equation}
A universal $C\sqrt W I$ bound holds without the scalar condition.
At unitary boundary points these are bounds for the limiting expected
kernel under a radial regulator; the unregularized kernel can be zero
almost surely and must not replace that boundary measure.
The same estimates hold if every physical collar diagonal has a common
positive height $\nu\le\nu_0$ and a sufficiently small common real shift.
\end{theorem}

The averaged boundary inverse belongs to the physical collar distribution. This estimate pays closing errors without replacing the collar by a uniform scattering phase.

\begin{proof}
For $Y\succ0$, apply \eqref{eq:weighted-boundary-identity} to this actual $K_c$.
Scalar Cauchy dilation gives, with $t$ standard Cauchy and
$X_t=X_h-tY$, $K_t=z_b-X_t$,
\begin{equation}
 \mathbb E_{\rm collar}\mathcal P(Ur)
       =\frac{\pi}{\eta}\mathbb E_t[K_t\Psi_{X_t}K_t^*].  \label{eq:collar-dilation-mean}
\end{equation}
This is also obtained by moving the terminating factors through the
top spectral projector before eliminating the fresh middle block.
That order is why the unbounded Schur load $B(G-e)^{-1}B^*$ creates
no unmatched weight. Lemma~\ref{lem:t22-dilation-variance} shows that \eqref{eq:collar-regularity} fails
on Cauchy mass at most $Cg/W$ when $d(X)\le d_*/2$.
Apply \eqref{eq:collar-regularity} on the remaining phases and \eqref{eq:collar-global-density} on the exceptional ones.
Their sum is $C(1+g)I$. Formula \eqref{eq:collar-global-density} alone gives the universal bound.
For semidefinite $Y$, use radial contraction and the limiting expected
kernel, retaining the expectation before removing the regulator.

For clarity, the small-height assertion does not integrate an
unbounded second moment against a Cauchy tail. For each fixed vector,
the collar-averaged kernel is a positive harmonic function of its
common complex diagonal shift in the upper half-plane. The preceding
boundary bound holds on a fixed real interval about zero: the
common-port scalar equations and one-point certificate are stable
on that interval. At height one its value is bounded uniformly in
$r,W$. Indeed the finite collar reaction has imaginary part at least
$I$ and all fixed operator moments; the Cayley defect then bounds
$\|(I-U_1r)^{-1}\|$ by a polynomial in those fresh norms.
In the positive harmonic representation, the height-$\nu$ Poisson
kernel outside the fixed real interval is at most $C\nu$ times the
height-one kernel. The possible positive linear harmonic term obeys
the same bound. Inside the interval use \eqref{eq:collar-poisson-bound}, or the universal
$C\sqrt W$ bound. This proves the claimed extension. The physical height
  $\delta/(N\rho(E))$  is included.
\end{proof}

\subsection{The limiting mean and scalar boundary law}\label{subsec:s23-the-limiting-mean-and-scalar-boundary-la}

  We fix the normalization of the closing law  in each boundary
direction.   Theorem~\ref{thm:t22-exact-boundary-mean} shows that the mean
positive boundary kernel tends to the identity, and the following
corollary identifies a scalar Cauchy law. We compare the weighted
spectral densities, bound the exceptional contributions, and then use
the analytic boundary representation.

\begin{theorem}[The physical mean becomes the identity]
\label{thm:t22-exact-boundary-mean}
Let $r$ be independent of the actual collar, $d(X)\to0$ and $g=O(1)$
in the preceding notation. Uniformly for deterministic such old responses,
\begin{equation}
 \mathbb E_{\rm collar}\mathcal P(Ur)=I+o_{\rm op}(1).     \label{eq:collar-mean-identity}
\end{equation}
More explicitly, for a deterministic sequence $\varepsilon_W\to0$,
depending only on the compact parameter ranges and fixed fresh
cutoffs, the error is at most
\begin{equation}
 C\{d(X)+\varepsilon_W+W^{-1/8}+gW^{-1/4}\}.              \label{eq:collar-mean-error}
\end{equation}
Boundary expectations are interpreted as in Theorem~\ref{thm:t22-collar-kernel}.
\end{theorem}

The limiting mean includes the endpoint factors and its exact
normalization. It determines the distant imaginary mass in the matrix
closing argument.

\begin{proof} The mean of the boundary kernel is a weighted
spectral density of the collar. We first treat a Hermitian old response.
Conditionally on the outer collar blocks, we compare this weighted density
with its deterministic value through a relative comparison that keeps all
large factors attached to the inverses that control them (Steps~1--3). We pay
for the exceptional event with a crude bound (Step~4), and taking expectations
gives the identity matrix in the limit (Step~5). Step~6 passes to absorbing
old responses by writing them as averages of Hermitian ones over a Cauchy
dilation, an average over real boundary conditions with Cauchy weights, as in
Figure~\ref{fig:bm-dilation}.

\emph{Step 1: the scalar error of the collar.} We first prove the
uniform statement for a Hermitian old response, including small real
shifts $s$ of the middle spectrum. We use the conditional matrices
$D,T,F,H$ of Lemma~\ref{lem:t22-collar-density}, and we set
\[
 e_{AB}=|W^{-1}\Tr H-a| ,
\]
the scalar error of the collar given its fresh blocks $A,B$.
Conditionally on the deterministic old $X$, the common-port scalar
comparison and concentration give
\[
 \|e_{AB}\|_{L^p(A,B)}\le C_p(d(X)+W^{-1}).
\]

\emph{Step 2: a good event and the relative density.} We choose fixed
large exponents $M$ and $D_0$, after the exceptional power we need. The
event
\[
 \|D^{-1}\|\le W^M,\qquad \|B\|\le C_B,\qquad e_{AB}\le e_*
\]
has a complement of probability $O(W^{-D_0})$, uniformly for $d(X)$
small enough. This follows from three bounds: the two-block Wegner
bound $\Prob(\|D^{-1}\|>W^M)\le CW^{1-M}$, the Gaussian tail of the bond,
and the fluctuation of the scalar error about its conditional mean,
which is $O(d(X))$. On this event, $F$ has a fixed polynomial bound on
its norm.

We apply the regularity certificate to the conjugate scalar equation of
$F$. It gives regularity and, uniformly for $|s|\le s_0$,
\[
 |r_{F,c}(s)-\rho(E)|\le C(e_{AB}+|s|),\qquad
 (1-C(e_{AB}+|s|))cHH^*\preceq\Gamma_F(s)
 \preceq(1+C(e_{AB}+|s|))cHH^*.
\]
Here $s_0$ and $e_*$ are fixed small numbers. The second inequality can
also be seen directly in the spectral basis of the Hermitian $F$. There
the shifts in its two denominators are $O(e_{AB}+|s|)$, and their
imaginary parts are bounded below.

\emph{Step 3: attaching the outside factors.} On this event we apply
Lemma~\ref{lem:t26-relative-matrix-density}. It compares the spectral
density matrix $\Sigma_F(s)$ with $\rho(E)\Gamma_F(s)$ in the relative
Loewner sense. We then take the congruence by $K_XT$, which attaches the
outside factors. With $L=K_XTH$, this gives
\begin{equation}\label{eq:t26-sharp-conditional-collar}
 -C(e_{AB}+|s|+\epsilon_W)LL^*\preceq
 K_XT\Sigma_F(s)T^*K_X^*-\rho(E)cLL^*
 \preceq C(e_{AB}+|s|+\epsilon_W)LL^*.
\end{equation}
All sectors and all complex directions are included before these
unbounded outside factors are attached. A congruence preserves the
Loewner order, so the comparison survives it. Thus the error contains
neither a small sector mass nor $\|T\|$.

\emph{Step 4: the exceptional event.} We now control the complement of
the good event, in the same order. For all $A,B$ and $|s|\le s_0$, the
global weighted bound gives
\[
 K_XT\Sigma_F(s)T^*K_X^*\preceq C\sqrt W\,LL^*.
\]
Indeed, we multiply the global bound for $\Sigma_F(s)$ by $F-ca$. The two
weighting polynomials have a uniformly bounded relative ratio on the
real spectrum of $F$, since the height is fixed and $|s|\le s_0$. The
congruence by $K_XT(F-ca)^{-1}$ then gives the display. The operator
moments of $L$ are uniformly bounded, so by H\"older's inequality its
expectation on the exceptional event is
\[
 O\bigl(W^{1/2-D_0(1-1/p)}\bigr)
\]
for any fixed $p>1$. We choose $D_0$ so that this is $o(W^{-1})$. The same
argument pays for the part of $LL^*$ itself that we removed.

\emph{Step 5: taking expectations.} We take expectations in
\eqref{eq:t26-sharp-conditional-collar}. By H\"older's inequality and the
moments of the scalar error and of $L$,
\[
 \E[e_{AB}\|LL^*\|]\le C(d(X)+W^{-1}).
\]
Lemma~\ref{lem:t22-common-port-mean} gives the mean of $LL^*$, which is
close to $\gamma I$. Together, these prove the uniform operator
asymptotic, with the shift included,
\begin{equation}\label{eq:t26-sharp-collar-density}
 \E_{A,B}[K_XT\Sigma_F(s)T^*K_X^*]
 =\rho(E)c\gamma I+
 O_{\rm op}(d(X)+\epsilon_W+W^{-1}+|s|).
\end{equation}
At $s=0$ the left side is exactly $K_X\Psi_XK_X^*$, by the complete Schur
spectral-density identity. Since $c\gamma=1-b$ and
$\eta=(1-b)\pi\rho(E)$, its leading value is
\[
 \rho(E)(1-b)I=(\eta/\pi)I.
\]
This proof is uniform in the old Hermitian response, whose norm was never
cut off.

\begin{figure}[H]
\centering
\begin{tikzpicture}[x=1cm,y=1cm]
 \draw[line width=5pt,blue!35] (0,0) circle (1.6);
 \draw[line width=6pt,red!70!black] (70:1.6) arc (70:84:1.6);
 \draw[line width=6pt,red!70!black] (205:1.6) arc (205:213:1.6);
 \draw[line width=6pt,red!70!black] (320:1.6) arc (320:326:1.6);
 \node[font=\small,align=center] at (0,0) {dilation\\ phases};
 \node[font=\small,anchor=west,align=left,blue!70!black] at (2.2,.55)
   {retained: old traces within $W^{-1/8}$\\
    of undilated values, use \eqref{eq:t26-sharp-collar-density}};
 \node[font=\small,anchor=west,align=left,red!70!black] at (2.2,-.65)
   {omitted: mass $\le CgW^{-3/4}$,\\
    crude bound $C\sqrt W$, total $CgW^{-1/4}$};
\end{tikzpicture}
\caption{The passive Cauchy dilation in the last step of the proof of
Theorem~\ref{thm:t22-exact-boundary-mean}. The passive old response is
written as an average over dilation phases of Hermitian responses. On
most phases (blue), the two normalized old resolvent traces stay within
$W^{-1/8}$ of their undilated values, so the Hermitian asymptotic
\eqref{eq:t26-sharp-collar-density} applies with error $O(W^{-1/8})$. The
remaining phases (red) have total mass at most $CgW^{-3/4}$, by
Lemma~\ref{lem:t22-dilation-variance}. On them we use only the crude
bound $C\sqrt W$ of Step~4, so they contribute at most $CgW^{-1/4}$.
These are the last two terms in \eqref{eq:collar-mean-error}.}
\label{fig:bm-dilation}
\end{figure}
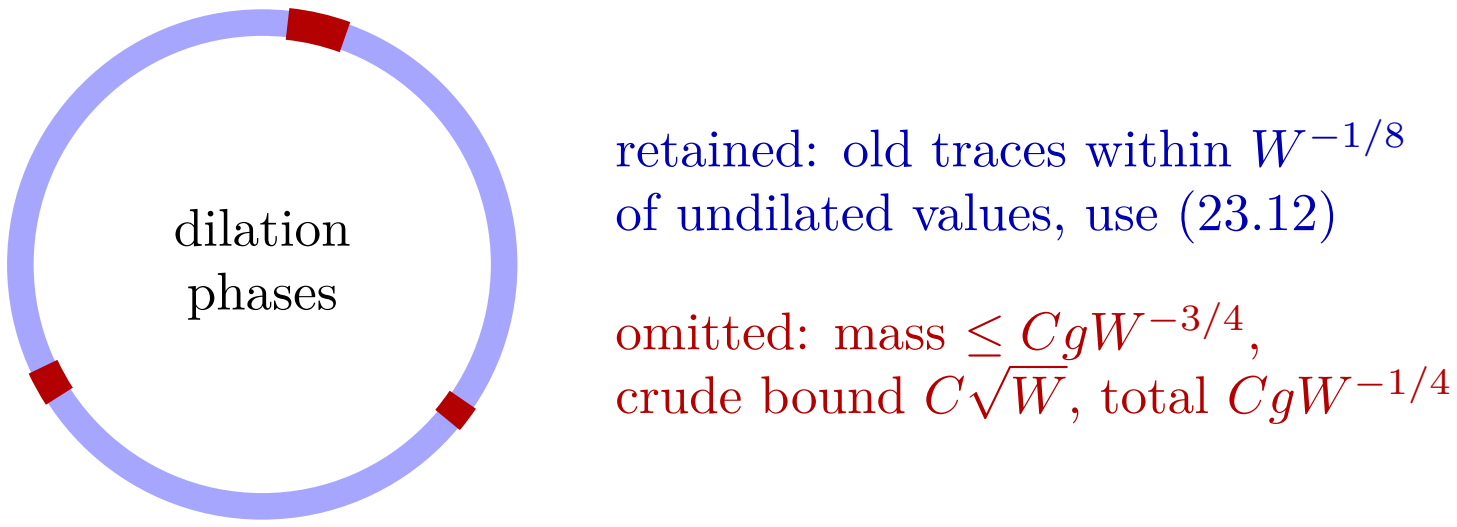

\emph{Step 6: passive responses by Cauchy dilation.} It remains to pass
from Hermitian to passive old responses. We write a passive response as
an average of Hermitian ones over a Cauchy dilation, as in
Figure~\ref{fig:bm-dilation}. We retain the phases for which the two
normalized old resolvent traces differ from their undilated values by at
most $W^{-1/8}$. Lemma~\ref{lem:t22-dilation-variance} bounds the mass of
the omitted phases by $CgW^{-3/4}$. On the retained phases we apply the
regular asymptotic of Step~5. On the omitted phases, the crude bound
$C\sqrt W$ of Step~4 makes their contribution at most
$CgW^{-3/4}\cdot C\sqrt W=CgW^{-1/4}$. Equation~\eqref{eq:collar-dilation-mean},
including its factor $\pi/\eta$, now gives \eqref{eq:collar-mean-error}.
Finally we remove a radial regulator at fixed $W$, as before.
\end{proof}

\begin{corollary}[An actual conditional one-channel Cauchy law]
\label{cor:t22-physical-one-channel-Cauchy}
Suppose the complete old two-port response $r$ is unitary and
$d(X)\to0$. Let $v$ be any old measurable unit vector, independent of
the fresh collar. Then under that collar law
\begin{equation}
 v^*(I-Ur)^{-1}v\ \Longrightarrow\ \frac{1}{2}+\frac{i}{2}\,\mathcal C,
                                                               \label{eq:scalar-collar-cauchy}
\end{equation}
uniformly in these choices, where $\mathcal C$ is standard real Cauchy.
The corresponding circle phase law converges in total variation.
No independence between several directions or several energies is asserted.
\end{corollary}

This conclusion identifies one quadratic form at a time. The joint matrix law still requires the decorated spectrum and the Herglotz tail argument.

\begin{proof}
For the lossless product $Ur$, the inverse satisfies $R+R^*=I$.
Thus $q=v^*Rv$ has real part $1/2$, and
$\zeta=1-q^{-1}$ lies on the unit circle almost surely.
Use the old rank-one damped twist
$r(I-vv^*+\rho e^{i\theta}vv^*)$, $0\le\rho<1$.
It changes the scalar traces by $O(W^{-1})$ and has conductance
$1-\rho^2\le1$. Writing $w=\rho e^{i\theta}$, the rank-one inverse identity gives
\[
 q_w=v^*(I-Ur(I-vv^*+wvv^*))^{-1}v
   =\frac{q}{w+(1-w)q}=\frac{1}{1-w\zeta}.
\]
Hence its $v$-kernel $2\Re q_w-1$ is exactly the scalar
Poisson kernel $(1-|w\zeta|^2)/|1-w\zeta|^2$.
Theorem~\ref{thm:t22-exact-boundary-mean} gives its expectation
$1+o(1)$ uniformly in $\rho,\theta$.
Letting $\rho\uparrow1$ shows that the law of $\zeta$ has density
between $1-o(1)$ and $1+o(1)$ relative to uniform circle measure:
integrate the Poisson estimates against continuous nonnegative tests
before taking that limit. This proves total variation convergence.
Finally $1/(1-e^{i\theta})=1/2+(i/2)\cot(\theta/2)$ proves \eqref{eq:scalar-collar-cauchy}.
\end{proof}
\section{The common matrix closing law}
\label{sec:t22-analytic-collar}
This section identifies one matrix-valued law for all retained boundary
directions at once. A single common closing matrix is needed so that all
probes see the same closure, which keeps their correlations.
Theorem~\ref{thm:matrix-closing} states that the closing matrix of the
actual collar converges to the canonical matrix $\mathsf C_r(0)$, whose poles
form a sine process and whose residues are independent Gaussian vectors. The
convergence is stable against bounded tests of the old path.
In Subsection~\ref{subsec:s24-the-positive-canonical-matrix-series} we construct the canonical matrix series. In
Subsection~\ref{subsec:s24-identification-of-the-distant-spectral-t} we show that, for a Herglotz function, convergence of the
local poles, a correct mean and the scalar Cauchy boundary laws determine the
limit (Lemma~\ref{lem:t22-Herglotz-tail}). In Subsection~\ref{subsec:s24-verification-for-the-actual-collar} we verify
these three conditions for the actual Gaussian collar.

We  use the actual three-block Gaussian collar of Section~\ref{sec:collar-density}, independent of a
lossless old two-port response $\mathfrak r$. The old Hermitian reaction is $X$,
$K=z_b-X$, $H_o=\Re z_b-X$, and $\mathfrak r=I-2i\eta K^{-1}$.
Assume $d(X)\to0$. All statements are uniform over these old matrices
and old orthonormal frames $J:\mathbb C^r\to\mathbb C^{2W}$ of fixed
rank $r$. Constants depend on the fixed $b$ and are uniform for $E$ in a compact bulk set. The mixed local array of Theorem~\ref{thm:relative-mixed-vectors}
identifies the poles and residues. Positivity determines the imaginary
tail; scalar boundary laws and bounded analytic means determine the
remaining real constant. These separate steps identify one common matrix
law for every finite collection of boundary directions.

\subsection{The positive canonical matrix series}\label{subsec:s24-the-positive-canonical-matrix-series}
 We construct the candidate closing law explicitly.
Theorem~\ref{thm:matrix-closing} identifies the actual finite collection
of boundary responses with this limit. We form a symmetrically cut-off
series of positive Gaussian residues at sine-process points, prove its
analytic convergence, and state the conditional physical comparison.

 Let $\{x_j\}$ be the unit-intensity sine process on $\mathbb R$ and
let $g_j$ be independent standard circular Gaussian vectors in
$\mathbb C^r$, independent of this process. Define
\begin{equation}
 \mathsf C_r(z)=\frac{1}{\pi}\lim_{R\to\infty}
                    \sum_{|x_j|\le R}\frac{g_jg_j^*}{x_j-z},
       \qquad \Im z>0.                                  \label{eq:canonical-cauchy}
\end{equation}
The cutoff is symmetric and retains its deterministic real center. The boundary matrix used below uses the same symmetric cutoff
and independent sine/Gaussian decoration law:
\[
 \mathsf C_r(0)=\frac1\pi\lim_{R\to\infty}
       \sum_{|x_j|\le R}\frac{g_jg_j^*}{x_j},
 \qquad\hbox{entrywise in probability}.
\]
This is not an undefined substitution into the formula for $\Im z>0$.
The limit exists locally uniformly in probability and in $L^2$ at
each nonreal point. Indeed split $g_jg_j^*=I+(g_jg_j^*-I)$.
The centered decoration tail has squared $L^2$ norm at most
$C_r\int_{|x|>R}|x-z|^{-2}dx$. For the scalar sine part, its
orthogonal-projection kernel gives
$\operatorname{Var}\sum f(x_j)\le\int|f|^2$; the same tail bound
applies, and its symmetric mean tends to $i\pi$.
The derivative versions, or the mean-value inequality on small circles,
give compact analytic convergence. Positivity of the imaginary part
passes from the finite sums. In particular
\begin{equation}
 \mathbb E\Im\mathsf C_r(z)=I_r,\qquad
 \mathbb E\|\mathsf C_r(iY)-iI_r\|_{\rm HS}^2\le C_r/Y.
                                                               \label{eq:canonical-cauchy-moments}
\end{equation}
Boundary values away from the locally finite poles are obtained by
first isolating a finite pole interval and then taking the same tail
limit. At $z=0$ no pole occurs almost surely. Unitary invariance and
compression consistency follow from the Gaussian decorations.
These facts are used only for the sine kernel, not for the non-Hermitian
biorthogonal kernel of deformed GUE.

\begin{theorem}[Complete conditional mixed physical boundary law]
\label{thm:matrix-closing}
Under the actual three-block collar law and uniformly over the old
lossless response and fixed frames specified above,
\begin{equation}
 -iJ^*\{2(I-U\mathfrak r)^{-1}-I\}J\ \Longrightarrow\ \mathsf C_r(0).
                                                               \label{eq:matrix-closing-law}
\end{equation}
The convergence is stable against bounded old-measurable tests. The
limiting matrix is independent of the old path. Its positive
construction is \eqref{eq:canonical-cauchy}, with a single common sine process and common
Gaussian decorations for all entries. All finite-rank compressions
and finite families of middle microscopic shifts converge jointly.
\end{theorem}

This joint matrix law is stronger than separate scalar Cauchy marginals.
Its common sine locations and Gaussian decorations retain all correlations
needed when several spatial sources are closed at once.

\subsection{Identification of the distant spectral tail}\label{subsec:s24-identification-of-the-distant-spectral-t}
 Convergence of nearby poles does not exclude extra mass or a real
constant coming from faraway eigenvalues.
Lemma~\ref{lem:t22-Herglotz-tail} shows that the positive mean and the
scalar Cauchy boundary laws determine these remaining terms. We use
positivity to control the imaginary tail, bounded analytic means to
fix the real constant, and then recover real boundary values.

 \begin{lemma}[Positive mass and bounded analytic means determine the tail]
\label{lem:t22-Herglotz-tail}
Let $F_N$ be random $r$-by-$r$ Herglotz functions of the form
$C_N+\int(x-z)^{-1}\mu_N(dx)$, where $C_N=C_N^*$ and $\mu_N$ is a
finite atomic positive matrix measure (so $F_N$ is rational and has
no continuous representing part). Suppose:
\begin{enumerate}
\item On every compact real interval their pole configurations and
matrix residues converge to those of \eqref{eq:canonical-cauchy}, including multiplicities.
\item For each $Y>0$, $\mathbb E\Im F_N(iY)\to I_r$.
\item For every real $x$ and deterministic $v\in\mathbb C^r$,
$v^*F_N(x)v$ converges to a real Cauchy variable of center zero and
scale $\|v\|^2$. Boundary values are taken outside the poles.
\end{enumerate}
Then $F_N\Rightarrow\mathsf C_r$ locally uniformly in the upper
half-plane and jointly at every finite set of real points. No
stationarity of $F_N$ or independence of its remote tail is assumed.
\end{lemma}

Local decorated poles alone leave a possible real constant and escaped positive mass. The two additional hypotheses determine both, including any correlation with the local poles.

\begin{proof}
First $F_N(i)$ is tight. Its imaginary part is tight by positivity
and assumption 2. For a scalar quadratic form, compare $F_N(i)$
with $F_N(0)$. Outside $(-a,a)$ the difference kernel obeys
\[
 \left|\frac{1}{x-i}-\frac{1}x\right|
       \le \sqrt{1+a^{-2}}\,\frac{1}{x^2+1}.
\]
The event that there is a pole in $(-a,a)$ has limiting probability
$O(a)$ by assumption 1. On its complement the preceding inequality,
assumption 2 and the Cauchy tightness in assumption 3 apply. First
choose $a$, then the tightness threshold. Polarization in finitely
many fixed vectors proves matrix tightness. The Herglotz representation
then gives tightness in the analytic-function topology: positive
harmonic functions satisfy Harnack bounds, and their real parts are
controlled by their value at $i$ and the Cauchy estimates.

In any joint subsequential limit, the representing spectral measure
is the local measure in \eqref{eq:canonical-cauchy}. This follows by distributional Stieltjes
inversion against compact tests; positivity precludes cancellation.
The only mass which can escape to infinity in this representation is
a nonnegative linear term $Az$, $A\succeq0$. Thus
\[
 F(z)=\mathsf C_r(z)+C+Az,\qquad C=C^*.
\]
Fatou and assumption 2 at $iY$, together with \eqref{eq:canonical-cauchy-moments}, give
$\mathbb E[\Im\mathsf C_r(iY)+YA]\preceq I_r$.
The first expectation equals $I_r$, so $A=0$ almost surely.
Equivalently one may saturate the positive imaginary mass with compact
spectral cutoffs. No principal-value cancellation was used here.

It remains to determine the possibly correlated real constant $C$.
For $t\ge0$, $\exp(itv^*F_N(z)v)$ is bounded analytic by one in the
upper half-plane. Its value is the Poisson integral of its boundary
values. Assumption 3 and dominated convergence, first on a fixed
real interval and then on its Poisson tail, imply
\begin{equation}
 \lim_N\mathbb E e^{itv^*F_N(z)v}=e^{-t\|v\|^2}
                 \quad(\Im z>0).                       \label{eq:bounded-analytic-mean}
\end{equation}
The same identity holds for the analytic limit. Now let $z=iY$ and
$Y\to\infty$. Equation \eqref{eq:canonical-cauchy-moments}, and the bound on the exponentials,
show that \eqref{eq:bounded-analytic-mean} becomes
$e^{-t\|v\|^2}\mathbb E e^{itv^*Cv}=e^{-t\|v\|^2}$.
Thus $v^*Cv=0$ almost surely, and polarization gives $C=0$.
In particular no independence between $C$ and the local decorated
process was used in this argument.

For real boundary points, remove an interval containing the point
and all nearby poles. The near part converges by assumption 1.
The far part is analytic across a smaller interval and converges
there by its Herglotz representation and analytic convergence off
the line. The probability of a pole arbitrarily close to the fixed
point tends to zero. This proves joint boundary convergence. The same
argument allows real arguments tending to the chosen fixed points.
\end{proof}

For the scalar Herglotz topology and boundary passage, see
Aizenman--Warzel~\cite[Sections 3 and 6]{AW}. The extra
argument \eqref{eq:bounded-analytic-mean} supplies the large-height constant from
bounded analytic means; it is not an application of their stationary
scalar Cauchy theorem to a matrix law.

\subsection{Verification for the actual collar}\label{subsec:s24-verification-for-the-actual-collar}
\begin{figure}[H]
\centering
\begin{tikzpicture}[x=.85cm,y=1cm]
 \draw[->] (-5.8,0)--(6.0,0) node[below,font=\small] {$x=W\rho(E)\lambda$};
 \foreach \x/\h in {-.9/.7,-.35/1.0,.25/.55,.8/.9,1.4/.6,-1.5/.8}
   \draw[very thick,blue!70!black] (\x,0)--(\x,\h);
 \foreach \x/\h in {-5.2/.4,-4.4/.5,-3.7/.35,3.6/.45,4.3/.3,5.1/.5}
   \draw[thick,black!55] (\x,0)--(\x,\h);
 \fill[blue!10,opacity=.6] (-2,0) rectangle (2,1.3);
 \draw[dashed,red!70!black,domain=-5.6:5.6,samples=120]
   plot (\x,{2.0*0.8/(\x*\x+0.64)});
 \node[font=\small,blue!70!black,align=center] at (0,-.75)
   {local atoms: sine process,\\ residues $g_jg_j^*/\pi$ (hyp.~1)};
 \node[font=\small,black!65,align=center] at (-4.3,-.75)
   {distant spectrum\\ (hyp.~2)};
 \node[font=\small,black!65,align=center] at (4.3,-.75)
   {distant spectrum\\ (hyp.~2)};
 \node[font=\small,red!70!black,anchor=west] at (1.1,2.25)
   {Poisson kernel at height $Y$};
\end{tikzpicture}
\caption{The representing measure of the collar function
\eqref{eq:physical-collar-function}, drawn schematically in the rescaled
variable $x=W\rho(E)\lambda$. Its atoms are the eigenvalues of $Q$, with
matrix residues $W\rho(E)V^*u_ju_j^*V$. Near zero (blue), these atoms
converge to the unit-intensity sine process with residues $g_jg_j^*/\pi$,
which is hypothesis~1. Far from zero (grey), the atoms are not
controlled individually. Hypothesis~2 asks only that, seen through the
Poisson kernel at height $Y$ (dashed), the total mean mass is $I_r$, so
none of it escapes to infinity. Hypothesis~3 concerns the boundary
values of single quadratic forms, which must have the standard Cauchy
law.}
\label{fig:cf-measure}
\end{figure}
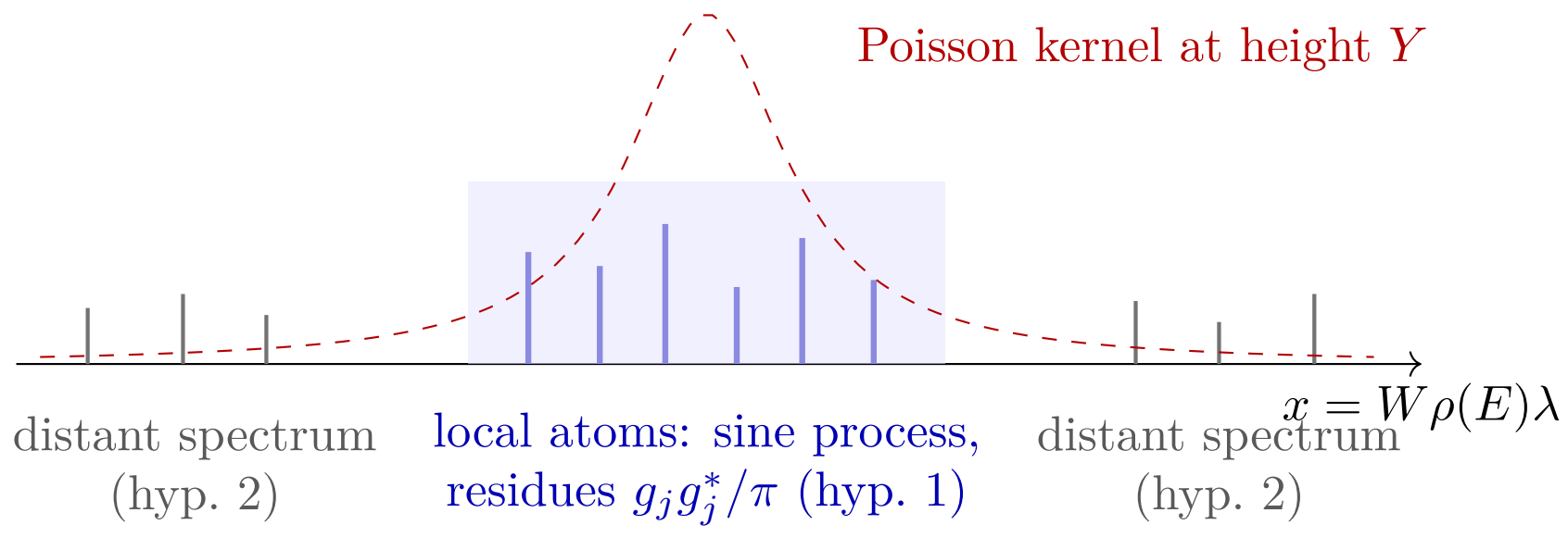

 The closing matrix is identified through
Lemma~\ref{lem:t22-Herglotz-tail}, which has three hypotheses: the local
poles and residues converge (hypothesis~1), no spectral mass escapes to
infinity (hypothesis~2), and single entries have the Cauchy boundary law
(hypothesis~3). We first write the physical collar as an exact Herglotz
function, and then we verify the three hypotheses in turn.
Figure~\ref{fig:cf-measure} shows which part of its representing measure
each hypothesis concerns.

\emph{The collar function.} We keep the actual collar variables of
Section~\ref{sec:collar-density}, and we set
\[
 D=A-E-X,\quad T=-D^{-1}B,\quad F=-E-B^*D^{-1}B,
 \quad Q=G+F,\quad V=\eta^{-1/2}T^*K^*J.
\]
The real inverses exist almost surely. For a shift only in the fresh
middle diagonal, we put
\begin{equation}
 F_W(z)=\eta^{-1}J^*(KD^{-1}K^*-H_o)J
          +V^*(Q-z/(W\rho(E)))^{-1}V.                  \label{eq:physical-collar-function}
\end{equation}
Its representing measure is the finite atomic spectral measure of $Q$,
with the positive matrix residues displayed. So it satisfies the
atomicity hypothesis of Lemma~\ref{lem:t22-Herglotz-tail}. It is an
exact Herglotz function, and its physical base value is
\begin{equation}
 F_W(0)=-iJ^*\{2(I-U\mathfrak r)^{-1}-I\}J.                     \label{eq:physical-collar-value}
\end{equation}

To verify \eqref{eq:physical-collar-value} with the factors in the right
order, we write $K_c$ for the collar reaction. Then
\[
 (K_c+i\eta)(I-U\mathfrak r)(H_o+i\eta)=2i\eta(K_c+H_o),
\]
and hence
\[
 (I-U\mathfrak r)^{-1}=\frac{K-K(K_c+H_o)^{-1}K^*}{2i\eta}.
\]
The Schur expansion $(K_c+H_o)^{-1}=D^{-1}+TQ^{-1}T^*$ now gives
\eqref{eq:physical-collar-function} and \eqref{eq:physical-collar-value}.
Every endpoint and every deterministic real normalizer is still present.

\emph{Hypothesis 1: the local atoms.} Conditionally on $A$ and $B$,
Theorem~\ref{thm:matched-traces} and Lemma~\ref{lem:regularity-certificate}
give the regularity of $F$ and a polynomial cutoff for it, with
arbitrarily small exceptional probability, as in
Lemma~\ref{lem:t22-collar-density}. The weighted frame $V$ also has a
polynomial cutoff, even when $\|X\|$ is unbounded, because
$KD^{-1}=I+(z_b+E-A)D^{-1}$. We let $H=(F-ca)^{-1}$, so that $L=KTH$.
Theorem~\ref{thm:relative-mixed-vectors} gives the complete decorated
local eigenvalue process of $Q$. Its projected covariance satisfies
\begin{equation}
 V^*\Gamma_F V=\frac{c}{\eta} J^*LL^*J+o_{\rm op}(1)
       \longrightarrow\frac{1}{\pi\rho(E)}I_r.          \label{eq:collar-projection-covariance}
\end{equation}
Here Lemma~\ref{lem:t22-common-port-mean} pays both the mean and the
bilinear fluctuations of $LL^*$. The identities $c\gamma=1-b$ and
$\eta=(1-b)\pi\rho(E)$ fix the constant, since
$(c/\eta)\gamma=(1-b)/\eta=1/(\pi\rho(E))$.

We must also fix the intensity of the eigenvalues in the scale used here.
The same common-port scalar estimate gives
$W^{-1}\Tr(F-ca)^{-1}=a+o(1)$ outside the exceptional set already paid
for. We apply Lemma~\ref{lem:regularity-certificate} to the conjugate
equation, with deformation $\widetilde F=F+EI$ at energy $E$. Indeed
$\widetilde F-E-c\bar a=F-c\bar a$, and
$m_{\widetilde F,c}(E)=m_{F,c}(0)$. So the positive Stieltjes branch
satisfies $m_{F,c}(0)=\bar a+o(1)$. Consequently, its free-convolution
density is
\[
 r_{F,c}(0)=\pi^{-1}\Im m_{F,c}(0)=\rho(E)+o(1).
\]
Theorem~\ref{thm:relative-mixed-vectors} therefore gives sine locations of
unit intensity in the scale $W\rho(E)\lambda_j(Q)$, rather than in an
unspecified local density scale. The residues of
\eqref{eq:physical-collar-function} at $x_j=W\rho(E)\lambda_j(Q)$ are
$W\rho(E)V^*u_ju_j^*V$. By \eqref{eq:collar-projection-covariance}, their
limits are precisely $g_jg_j^*/\pi$. The conditional convergence is
uniform, so it remains valid after averaging over the actual $A,B$, and
independently of bounded tests of the old path. This proves
hypothesis~1 of the lemma.

Hypotheses~2 and~3 use the collar with only its middle diagonal shifted.
For the weighted density this shift is included in
\eqref{eq:t26-sharp-collar-density}; Theorem~\ref{thm:t22-exact-boundary-mean}
is stated for $s=0$, and the common shift in
Theorem~\ref{thm:t22-collar-kernel} also moves the end diagonals. With
$e=E$ as in Section~\ref{sec:collar-density}, the shifted collar is the
unshifted one at a shifted energy: for $e_s=e+s$ and $X_s=X-sI$,
\begin{equation}\label{eq:middle-shift-matrix}
 \begin{pmatrix}A-e-X&B\\B^*&G-e-s\end{pmatrix}
 =\begin{pmatrix}A-e_s-X_s&B\\B^*&G-e_s\end{pmatrix},
\end{equation}
so only the boundary reference $z_b$ and the outside weight change. The appropriate version of the result is included in Lemma~\ref{lem:middle-shift-collar} which will be appearing at the end of this subsection.

\emph{Hypothesis 2: the mean of the whole spectrum.} For clarity, we
note that hypothesis~2 includes the full distant spectrum, not only the
local window. Equation~\eqref{eq:t26-sharp-collar-density}, with $V$ as
in \eqref{eq:physical-collar-function}, gives uniformly for $|s|\le s_0$
\begin{equation}
 \mathbb E_{A,B,G}V^*\delta_s(Q)V
 =\pi^{-1}I_r+
 O_{\rm op}(d(X)+\epsilon_W+W^{-1}+|s|).
                                  \label{eq:shifted-weighted-density}
\end{equation}
The factor is $\eta^{-1}\rho(E)c\gamma=1/\pi$. We obtain it by testing the
full Loewner inequality in every vector $Jv$, so it is uniform for the
old orthonormal frame and for all its polarizations. Its proof used the
complete sum over ordered sectors, the relative covariance and the
exceptional weighted bound. So no unweighted error is multiplied by
$\|T\|^2$. The shift of the middle eigenvalue is already covered by the
fixed interval in \eqref{eq:t26-sharp-conditional-collar}, with the same
$H$ and $L$ and an unchanged Gaussian conditioning.

Equation~\eqref{eq:shifted-weighted-density} controls the spectrum near zero.
For the distant spectrum we use a deterministic comparison at a fixed
positive height:
\begin{equation}
 \|(Q-i)^{-1}V\|_{\rm HS}
 \le (1+|i-ca|+\|G\|)\|HV\|_{\rm HS},\qquad
 \|HV\|_{\rm HS}^2=\eta^{-1}\operatorname{Tr}J^*LL^*J.
                                                               \label{eq:distant-collar-bound}
\end{equation}
It follows by inserting $(F-ca)H=I$ to the right of $(Q-i)^{-1}$. Indeed
$(Q-i)^{-1}(F-ca)=I+(Q-i)^{-1}(i-ca-G)$, and $\|(Q-i)^{-1}\|\le1$ since
$Q$ is Hermitian. The identity for $\|HV\|_{\rm HS}^2$ uses that
$H=(F-ca)^{-1}$ is normal, being a function of the Hermitian $F$, so
$H^*H=HH^*$. Lemma~\ref{lem:t22-common-port-mean} and the fresh Gaussian
norm moments therefore bound the expectation of
$V^*\Im(Q-i)^{-1}V$ by $C_rI$ in quadratic forms.

Consequently the eigenvalues with $|\lambda(Q)|>s_0$ contribute at most
$C_{r,s_0}Y/W$ to $\mathbb E\Im F_W(iY)$. On $|\lambda|\le s_0$, we integrate
\eqref{eq:shifted-weighted-density} against the Poisson kernel of height
$Y/(W\rho(E))$. We first let $W\to\infty$, and then $s_0\downarrow0$. The
answer is $I_r$. This proves hypothesis~2, including the polynomially
distant deformation and the complete normalizer.

\emph{Hypothesis 3: the Cauchy boundary law.} For hypothesis~3 we repeat
the rank-one argument of Corollary~\ref{cor:t22-physical-one-channel-Cauchy},
with the middle diagonal shifted by $x/(W\rho(E))$. By
\eqref{eq:shifted-weighted-density} and the same exceptional-set
estimate for the Cauchy dilation, \eqref{eq:collar-mean-error} changes by
at most $C|x|/W$ for each fixed $x$. Precisely, Lemma~\ref{lem:middle-shift-collar} with
$s=x/(W\rho(E))$ gives this bound for the twisted old responses below,
which have $d(X')\le d(X)+C/W$ and $g\le1$. Each twist $r'$ with $\rho<1$
has a passive Cayley inverse $X'$: if $r'y=y$, then $r(y-cv)=y$ with
$c=(1-\rho e^{i\theta})v^*y$, so unitarity of $r$ gives
$(1-\rho^2)|v^*y|^2=0$, hence $ry=y$ and $y=0$ since $I-r$ is invertible;
and $I-r'r'^*=(1-\rho^2)rvv^*r^*\succeq0$.

The argument uses rank-one twists of the old response. We write them as
\[
 r'=r(I-vv^*+\rho e^{i\theta}vv^*),\qquad \|v\|=1 .
\]
Since $(z_b-X')^{-1}=(I-r')/(2i\eta)$, the two normalized block traces
change by at most $C/W$, uniformly in $\rho$ and $\theta$. Moreover
$\Tr(I-r'^*r')=1-\rho^2\le1$. The Cayley inverse $X'$ may be unbounded, but
the sharp density bound and the global exceptional bounds were uniform in
exactly that variable. Singular Cayley inverses are handled by radial
limits. Thus the shifted kernel estimate is uniform in every twist that
we use to approximate the boundary measure. We integrate its scalar
Poisson kernels against nonnegative continuous test functions on the
circle, before letting $\rho\uparrow1$. This gives convergence of the law of
the phase in total variation.

Hence every quadratic form of \eqref{eq:physical-collar-function} in a
unit vector has the standard Cauchy boundary limit. The Poisson-kernel
estimate is uniform up to the radial boundary, exactly as in
Theorem~\ref{thm:t22-exact-boundary-mean}, so no singular boundary mass is
discarded. General vectors $v$ follow by homogeneity. This proves
hypothesis~3.

\begin{proof}[Proof of Theorem~\ref{thm:matrix-closing}]
Apply the preceding lemma after \eqref{eq:physical-collar-function}--\eqref{eq:distant-collar-bound}. The sequential uniformity
of all three estimates proves uniform bounded-continuous comparison
by contradiction along a violating old deterministic sequence.
Conditioning then gives the stable statement. Enlarging the old frame
handles any finite collection of compressions simultaneously.
\end{proof}

\begin{lemma}[Middle-only shift of the collar mean]
\label{lem:middle-shift-collar}
For real $s$ let $U_s$ be the scattering matrix of the actual collar whose
middle diagonal $G-e$ is replaced by $G-e-s$, the end diagonals $A-e$ and
$z_b$ being unchanged:
\[
 K_c^{(s)}=A-e-\xi-B(G-e-s)^{-1}B^*,\qquad
 U_s=I-2i\eta(K_c^{(s)}+i\eta)^{-1}.
\]
Let $X=X_h-iY$, $Y\succeq0$, be deterministic,
$r=I-2i\eta(z_b-X)^{-1}$ and $g=\Tr(I-r^*r)$. Fix $g_0\ge1$. There are $s_0,d_1>0$
and $C$ such that, uniformly for $|s|\le s_0$, $d(X)\le d_1$ and
$g\le g_0$,
\[
 \|\E_{\rm collar}\mathcal P(U_sr)-I\|
 \le C\{d(X)+\varepsilon_W+W^{-1/8}+gW^{-1/4}+|s|\}.
\]
Boundary expectations are interpreted as in
Theorem~\ref{thm:t22-collar-kernel}. In particular, for
$s=x/(W\rho(E))$ with fixed $x$, the error \eqref{eq:collar-mean-error}
increases by at most $C|x|/W$.
\end{lemma}

\begin{proof}
Put $z_s=z_b(e_s)=\xi_s+i\eta_s$ and $\Delta=z_s+s-z_b$. As $z_b$, $a$
and $\rho$ are smooth in the energy,
$|\Delta|+|\eta_s-\eta|+|a(e_s)-a(e)|\le C|s|$. For small $s_0$, $e_s$
stays in a slightly enlarged compact interior energy set, on which we
take the constants of \eqref{eq:t26-sharp-collar-density} and of
Theorem~\ref{thm:t22-exact-boundary-mean} and $\eta,\eta_s$ are bounded
below. With $d_*$ as in Theorem~\ref{thm:t22-collar-kernel},
\eqref{eq:collar-mean-error} is used on $d(X)\le d_*/2$, with constants
depending on a bound for $g$.

With $Z=z_b-X$, $Z-Z^*=2i(\eta+Y)$ gives $I-rr^*=4\eta Z^{-1}YZ^{-*}$,
so $\|r\|<1$ if and only if $Y\succ0$. First let $Y\succ0$. We prove the bound for $d(X)\le2d_1$, $g\le2g_0$. Put $Z_s=z_s-X_s$, so that $Z=Z_s-\Delta$,
and $X_s$ also has imaginary part $-Y$. As $\Im Z\succeq\eta$ and
$\Im Z_s\succeq\eta_s$, we have $\|Z^{-1}\|\le\eta^{-1}$ and
$\|Z_s^{-1}\|\le\eta_s^{-1}$. Let $d_s$, $r_s=I-2i\eta_sZ_s^{-1}$ and $g_s=\Tr(I-r_s^*r_s)$ be the
scalar distance, old scattering matrix and conductance for $X_s$ at
$e_s$. From
$Z_s^{-1}-Z^{-1}=-\Delta Z_s^{-1}Z^{-1}$, every diagonal entry of
$Z_s^{-1}-Z^{-1}$ is at most $|\Delta|/(\eta\eta_s)$ in modulus, so
\begin{equation}\label{eq:middle-shift-distance}
 d_s\le d(X)+C|s|.
\end{equation}
Likewise $I-r_sr_s^*=4\eta_sZ_s^{-1}YZ_s^{-*}$, and since
$Z_s^{-1}=(I-\Delta Z_s^{-1})Z^{-1}$,
\begin{equation}\label{eq:middle-shift-conductance}
 g_s=4\eta_s\Tr(Z_s^{-1}YZ_s^{-*})
 \le\frac{\eta_s}{\eta}\Bigl(1+\frac{|\Delta|}{\eta_s}\Bigr)^2g\le Cg .
\end{equation}
We choose $d_1$ and $s_0$ with $2d_1+Cs_0\le d_*/2$. Then $d_s\le d_*/2$
and $g_s\le Cg_0$.

Put $\mathcal N_s=A-e-X-B(G-e-s)^{-1}B^*$. Then
$K_c^{(s)}+(\xi-X_h)+iY=\mathcal N_s$. The reaction
$\widetilde K_c=A-e_s-\xi_s-B(G-e_s)^{-1}B^*$ of the unshifted collar at
energy $e_s$, with scattering matrix
$\widetilde U=I-2i\eta_s(\widetilde K_c+i\eta_s)^{-1}$, satisfies
$\widetilde K_c+(\xi_s-X_h+s)+iY=\mathcal N_s$ as well, by
\eqref{eq:middle-shift-matrix}. Apply
\eqref{eq:weighted-boundary-identity} to both, with $H_o=\xi-X_h$ and
$H_o=\xi_s-X_h+s$ respectively. With $\Pi=I-Z\mathcal N_s^{-1}$ and
$\Pi_s=I-Z_s\mathcal N_s^{-1}$, this gives
\[
 \eta\,\mathcal P(U_sr)=\Pi Y\Pi^*,\qquad
 \eta_s\,\mathcal P(\widetilde Ur_s)=\Pi_sY\Pi_s^* .
\]
Hence $\Pi=\Pi_s+\Delta\mathcal N_s^{-1}$, and
$\mathcal N_s^{-1}=Z_s^{-1}(I-\Pi_s)$. For a vector $q$ put
$\tilde q=Z_s^{-*}q$, $y_1=Y^{1/2}\Pi_s^*q$ and
$y_2=Y^{1/2}\mathcal N_s^{-*}q=Y^{1/2}\tilde q-Y^{1/2}\Pi_s^*\tilde q$.
Then, for every collar sample,
\[
 \bigl|\eta\,q^*\mathcal P(U_sr)q-\eta_s\,q^*\mathcal P(\widetilde Ur_s)q\bigr|
 =\bigl|\|y_1+\bar\Delta y_2\|^2-\|y_1\|^2\bigr|
 \le|\Delta|(\|y_1\|^2+\|y_2\|^2)+|\Delta|^2\|y_2\|^2,
\]
where $\|y_1\|^2=\eta_s\,q^*\mathcal P(\widetilde Ur_s)q$ and, since
$Z_s^{-1}YZ_s^{-*}\preceq(4\eta_s)^{-1}I$,
\[
 \|y_2\|^2\le2q^*Z_s^{-1}YZ_s^{-*}q
 +2\eta_s\,\tilde q^*\mathcal P(\widetilde Ur_s)\tilde q
 \le\frac{\|q\|^2}{2\eta_s}
 +2\eta_s\,\tilde q^*\mathcal P(\widetilde Ur_s)\tilde q .
\]
In particular $\mathcal P(U_sr)$ is integrable whenever
$\mathcal P(\widetilde Ur_s)$ is. The collar law does not depend on the
energy, so Theorem~\ref{thm:t22-exact-boundary-mean} applies at $e_s$ to
the deterministic old response $X_s$. By
\eqref{eq:middle-shift-distance} and \eqref{eq:middle-shift-conductance},
\[
 \|\E_{\rm collar}\mathcal P(\widetilde Ur_s)-I\|\le\omega_s
 :=C\{d(X)+|s|+\varepsilon_W+W^{-1/8}+gW^{-1/4}\}.
\]
Take expectations above and use $\|\tilde q\|\le\|q\|/\eta_s$. Then
$\E\|y_1\|^2\le\eta_s(1+\omega_s)\|q\|^2$,
$\E\|y_2\|^2\le C(1+\omega_s)\|q\|^2$, and
$|\eta_s\E q^*\mathcal P(\widetilde Ur_s)q-\eta\|q\|^2|
\le(\eta_s\omega_s+|\eta_s-\eta|)\|q\|^2$. As $|\Delta|\le C|s|\le C$,
\[
 \bigl|q^*(\E_{\rm collar}\mathcal P(U_sr)-I)q\bigr|
 \le C(|s|+\omega_s)\|q\|^2 .
\]
The matrix $\E_{\rm collar}\mathcal P(U_sr)-I$ is Hermitian, so this
proves the bound when $Y\succ0$.

For $Y\succeq0$ and $0<\varsigma<1$, $\varsigma U_sr=U_s(\varsigma r)$,
where $\varsigma r$ is the old scattering matrix of the passive response
$z_b-2i\eta(I-\varsigma r)^{-1}$, with strictly positive $Y$ as
$\|\varsigma r\|<1$; as
$\varsigma\uparrow1$ at fixed $W$, $(I-\varsigma r)/(2i\eta)\to(I-r)/(2i\eta)$
and $\Tr(I-\varsigma^2r^*r)\to g$; so for $\varsigma$ near one the scalar
distance and conductance of $\varsigma r$ are at most $2d_1$, $2g_0$ and
tend to those of $r$, and the preceding case passes to the limiting
expected kernel under the radial regulator. Finally, for $s=x/(W\rho(E))$ we have $|s|\le s_0$
once $W$ is large, and the extra term is $C|x|/(W\rho(E))\le C|x|/W$.
\end{proof}
\section{The entire determinant and the joint marked law}
\label{sec:entire-law}
\label{sec:t22-full-critical}
This section completes Part~II and proves Theorem~\ref{thm:critical}. It
assembles the source law of the path and the closing law of the collar into
the full critical spectral limit.
Theorem~\ref{thm:entire-law} states that the determinant ratio
\eqref{eq:entire-determinant} converges in law, locally uniformly, to a unique
random entire function $D$, given by an explicit positive construction.
Theorem~\ref{thm:joint-recovery} states that the marked point process of
rescaled eigenvalues and profiles converges to a unique limit, with simple
points and unit intensity. The limit is recovered from the zeros of $D$ and
the ratios of its partial derivatives at those zeros.
In Subsection~\ref{subsec:s25-global-gaussian-probes-of-the-actual-mat} we build Gaussian probes spread over the path
(Lemma~\ref{lem:t22-global-Gaussian-source}). In Subsection~\ref{subsec:s25-inserting-the-actual-physical-closing-la} we
insert the closing law through Woodbury's formula. In Subsection~\ref{subsec:s25-restoring-sources-on-all-physical-blocks} we
restore the sources on the three collar blocks. In Subsection~\ref{subsec:s25-determination-and-compactness} we
prove that the limit is unique and tight, which proves
Theorem~\ref{thm:entire-law}. In Subsection~\ref{subsec:s25-recovering-all-energies-and-marks-jointl} we recover all
eigenvalues and marks and prove Theorem~\ref{thm:joint-recovery}.

 Return to the original ring, with $n/W\to\tau\in(0,\infty)$,
$0<b<1/2$ fixed and $|E|<2$ fixed. All limits here first send $W$
to infinity at fixed $\tau$. The constants may be uniform, at fixed $b$, for $E$ in compact
bulk sets and $\tau$ in compact positive intervals. The only
random matrix ensemble in the comparison is the original
Gaussian block ring. The preceding positive source process and common matrix closing law
now determine the complete determinant quotient. Its zeros and source
derivatives recover all local eigenvalues and profile marks jointly.

For bounded nonnegative continuous functions $q_0=1,q_1,\ldots,q_k$
on $\T$, let $Q_a=\sum_jq_a(j/n)P_j$ and define
\begin{equation}
 D_W(z_0,\ldots,z_k)
   =\frac{\det(H-E+W^{-2}\sum_{a=0}^k z_aQ_a)}{\det(H-E)},\qquad Q_0=I.               \label{eq:entire-determinant}
\end{equation}
The denominator is nonzero almost surely and $D_W(0)=1$.

\begin{theorem}[Full original-ring extended determinant law]
\label{thm:entire-law}
The random entire function \eqref{eq:entire-determinant} converges in law, locally uniformly
on $\mathbb C^{k+1}$, to a unique random entire function $D$.
Its law is specified by the positive construction \eqref{eq:global-probe-replacement}--\eqref{eq:limiting-closed-source}:
for any finite list $y^{(j)}\in[0,\infty)^{k+1}$ and signs
$\sigma_j\in\{1,-1\}$, take one global probe with profile
$\sum_a y_a^{(j)}q_a$ per factor. Construct their common matrix
$\mathcal T$ by \eqref{eq:limiting-closed-source}. Then
\begin{equation}
 \mathbb E\prod_{j:\sigma_j=1}D(iy^{(j)})^{-1}
           \prod_{j:\sigma_j=-1}\overline{D(iy^{(j)})^{-1}}
     =\mathbb E\exp\{-i\sum_j\sigma_j\mathcal T_{jj}\}.
                                                               \label{eq:determining-bounded-moments}
\end{equation}
These mixed moments uniquely determine the law of $D$ among laws of random entire functions on $\mathbb C^{k+1}$ with
$|D(iy)|\ge1$ almost surely for every $y$ in the closed positive cone
$[0,\infty)^{k+1}$. The determination concerns $D$ itself, including its
dependence on the first-trace variable $z_0$; no exponential factor in
$z_0$ is divided out. The limit lies in this class, since
$|D_W(iy)|\ge1$ at every finite width. The limits are consistent in the
profile family: deleting a profile $q_a$ with $a\ge1$ replaces $D$ by
its restriction to $\{z_a=0\}$.
\end{theorem}

The determining variables are bounded inverse determinants on the positive
imaginary cone. Retaining the first trace identifies exactly the physical
quotient in \eqref{eq:entire-determinant}, without a separate restoration
of a random exponential factor.

\subsection{Global Gaussian probes of the actual matched path}\label{subsec:s25-global-gaussian-probes-of-the-actual-mat}
 We build the source matrices for profiles spread across the physical
path. Lemma~\ref{lem:t22-global-Gaussian-source} gives their unique
critical limit from finitely many random local source births. We
sample the birth locations, use Gaussian replacement to pass to
extended probes, and retain the joint real and imaginary parts.

 Theorem~\ref{thm:complex-source-law} is
proved at the terminal lengths $n_T=\lfloor WT\rfloor$. In applications
the path need not have exactly this length. In
Lemma~\ref{lem:t22-global-Gaussian-source} the physical path is what
remains of the ring after the collar is removed, so its number of blocks
is fixed by the ring and only its ratio to $W$ converges to the critical
time. In general it differs from $\lfloor WT\rfloor$ by an amount that is
small compared with $W$. The next lemma shows that such a change of the
terminal length on the scale $o(W)$ does not change the limit law of the
probe matrix. This is what allows Lemma~\ref{lem:t22-global-Gaussian-source}
to use the finite-source limit of Theorem~\ref{thm:complex-source-law} at
the actual length of the path.

We grow the path layer by layer. Step $k\ge0$ adds block $k+1$ in the
labelling of Section~\ref{sec:full-source}, so the probe with birth time
$\tau_c$ is born at step $\lfloor W\tau_c\rfloor$. We write
$\mathcal K^{(r)}$ for the probe matrix \eqref{eq:complex-source-data}
after $r$ layers.

\begin{lemma}[Terminal lengths]\label{lem:t22-terminal-lengths}
Fix $0<b<1/2$, $-2<E<2$, masses $m_c\ge0$ and distinct birth times
$0<\tau_c<T$. Let $W_j\to\infty$, and let $r_j$ be integers with
$r_j/W_j\to T$. Then, under the product of the open path law and the law
of the probe Gaussians, the law of $\mathcal K^{(r_j)}$ at width $W_j$
converges to the critical limit of Theorem~\ref{thm:complex-source-law}
at terminal time $T$.
\end{lemma}

\begin{figure}[ht]
\centering
\begin{tikzpicture}[x=1cm,y=1cm,>=Stealth]
 \draw[->] (-.2,0)--(12.0,0) node[right,font=\small] {layer $n$};
 \foreach \x in {1.4,3.1,4.6}{
   \fill[red!70!black] (\x,0) circle (2.3pt);}
 \node[font=\small,red!70!black,above] at (3.1,.15) {births at $\lfloor W\tau_c\rfloor$};
 \draw[thick] (6.2,-.1)--(6.2,.1) node[below=6pt,font=\small] {$(T-\delta)W$};
 \draw[thick] (10.4,-.1)--(10.4,.1) node[below=6pt,font=\small] {$(T+\delta)W$};
 \fill[blue!15] (8.0,-.12) rectangle (8.7,.12);
 \draw[thick,blue!70!black] (8.0,-.25)--(8.0,.25) node[above,font=\small] {$R_j$};
 \draw[thick,blue!70!black] (8.7,-.25)--(8.7,.25) node[above,font=\small] {$r_j$};
 \draw[decorate,decoration={brace,amplitude=4pt,mirror}] (8.0,-.7)--(8.7,-.7)
   node[midway,below=4pt,font=\small,blue!70!black] {$N_j=o(W)$ layers};
 \draw[decorate,decoration={brace,amplitude=4pt}] (6.2,.75)--(10.4,.75)
   node[midway,above=4pt,font=\small] {window: no births};
\end{tikzpicture}
\caption{The two terminal lengths in the proof of
Lemma~\ref{lem:t22-terminal-lengths}. All probes are born (red) before the
window $[(T-\delta)W,(T+\delta)W]$. The reference length
$R_j=\lfloor W_jT\rfloor$ and the actual length $r_j$ both lie in this
window, and they are separated by only $N_j=o(W)$ layers (blue). Each of
these layers changes the probe matrix by a conditional mean of order
$1/W$ and a centered fluctuation of variance of order $1/W$, so the total
change over $o(W)$ layers tends to zero.}
\label{fig:terminal-lengths}
\end{figure}
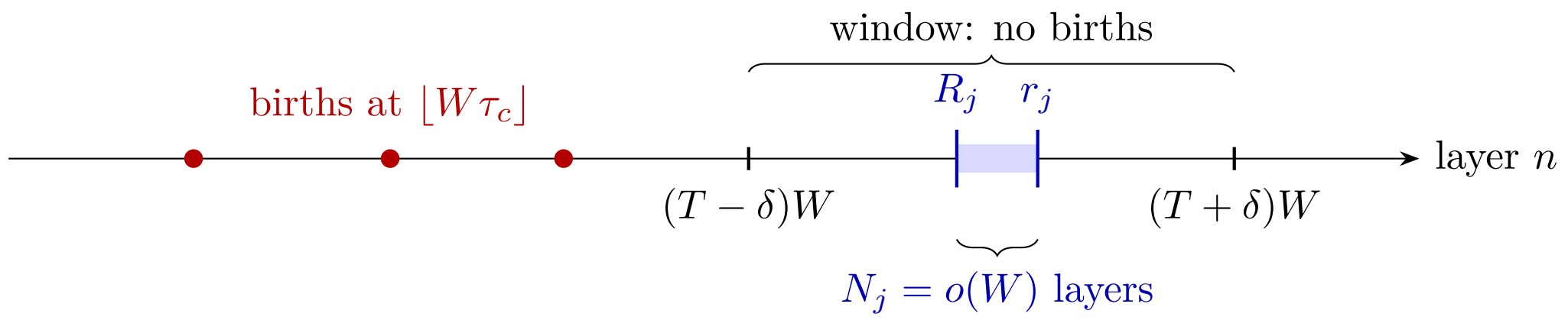

\begin{proof}
\emph{Plan of the proof.} Theorem~\ref{thm:complex-source-law} gives the
limit at the reference lengths $R_j=\lfloor W_jT\rfloor$, so it suffices to
show that the probe matrix changes little between the lengths $R_j$ and
$r_j$ (Step~1). These lengths are separated by $o(W)$ layers that contain
no birth (Step~2), see Figure~\ref{fig:terminal-lengths}. On each such
layer, the increment of the probe matrix has a conditional mean of order
$1/W$ and a fluctuation with conditional variance of order $1/W$
(Steps~3 and~4). Summing over $o(W)$ layers shows that the total change
tends to zero in probability (Step~5).

\emph{Step 1: reduction.} Put $R_j=\lfloor W_jT\rfloor$.
Theorem~\ref{thm:complex-source-law} gives the limit along $W_j$ at the
lengths $R_j$. Both $\mathcal K^{(r_j)}$ and $\mathcal K^{(R_j)}$ are
defined on the same probability space, the product of the open path law
and the probe law. So it suffices to show that
$\mathcal K^{(r_j)}-\mathcal K^{(R_j)}\to0$ in probability: since bounded
Lipschitz tests determine convergence in law, the two laws then have the
same limit.

\emph{Step 2: the window between the two lengths.} Let
$\delta=(T-\max_c\tau_c)/2$. For large $j$, both lengths lie in
$((T-\delta/2)W_j,(T+\delta/2)W_j)$, and all probes are born before step
$(T-\delta)W_j$. So the
\[
 N_j=|r_j-R_j|=o(W_j)
\]
layers between the two lengths are not birth steps, and they start at
lengths $(T-\delta)W_j\le n\le(T+\delta)W_j$.

\emph{Step 3: the one-step increment.} Let $\mathcal F_n$ be the
$\sigma$-algebra generated by the first $n$ layers and by all probe
Gaussians. Then $\mathcal K^{(n)}$ is $\mathcal F_n$-measurable, and the
next layer is independent of $\mathcal F_n$. At a step without a birth,
\eqref{eq:complex-source-update} writes the increment of an entry as
$u_k^*\Delta v_l$, with the $\mathcal F_n$-measurable vectors
\[
 u_k=\iota^*V_{\rm pr}e_k,\qquad v_l=\iota^*\mathbb SV_{\rm pr}e_l .
\]
Since $\mathbb S$ is unitary, their norms are at most $\|V_{\rm pr}e_k\|$
and $\|V_{\rm pr}e_l\|$. These norms are controlled by the probe matrix
itself, through the imaginary-part identity:
\[
 \|V_{\rm pr}e_k\|^2=(M_{\rm pr})_{kk}=-\eta^{-1}\Im\mathcal K_{kk}.
\]
By \eqref{eq:raw-tag-moments}, with $T+\delta$ in place of $T$, these norms
have all fixed moments bounded, uniformly in $W$ and over the window.

\emph{Step 4: conditional mean and variance of one layer.} We use the
passive mean and bilinear bounds, with the constants $C_m'$ and $F_2$ of
Steps~8 and~6 of the proof of Theorem~\ref{thm:complex-source-law}. They
hold conditionally on the past, uniformly over passive inputs, and we
check them in two cases.

On the scalar domain $\delta_0\le d_0$, they follow from
Lemma~\ref{lem:passive-bounds} and the expansion of $\E_j\Delta$ in
\eqref{eq:full-source-estimates}. There the scalar error $\Omega=O(\delta_0)$
multiplies the bounded weighted inverse $KM$, and the Gaussian $L^p$
Poincar\'e argument of that lemma bounds the centered bilinear forms of
$\Delta$ by $CW^{-1/2}\|u\|\|v\|$.

Off the scalar domain, $\delta_0>d_0$, they follow with the constants
$D_2/d_0$, because then
\[
 \|\E_{\rm fresh}\Delta\|\le D_2\le\frac{D_2}{d_0}\,\delta_0,
 \qquad
 \|u^*\Delta v\|_{L^2({\rm fresh})}\le D_2\|u\|\|v\|
 \le\frac{D_2}{d_0}\,\delta_0\|u\|\|v\|.
\]

At the open near input these bounds give
\[
 \begin{aligned}
 |\E[\Delta\mathcal K_{kl}\mid\mathcal F_n]|
   &\le C_m'(\delta_0+W^{-1})\|V_{\rm pr}e_k\|\|V_{\rm pr}e_l\|,\\
 \E[|\Delta\mathcal K_{kl}|^2\mid\mathcal F_n]
   &\le F_2^2(\delta_0+W^{-1/2})^2\|V_{\rm pr}e_k\|^2\|V_{\rm pr}e_l\|^2,
 \end{aligned}
\]
where $\delta_0$ is the scalar jet of the open path at length $n$. By
H\"older's inequality, Lemma~\ref{lem:t26-passive-initial-scalar} with $T$
replaced by $T+\delta$ (or \eqref{eq:passive-jet-moments}), and the
moment bounds of Step~3, there are constants $C$ and $W_0$ such that, for
$W\ge W_0$, every $n$ in the window and every real coordinate $x$ of
$\mathcal K$ (the real or the imaginary part of an entry),
\[
 \E\bigl|\E[x_{n+1}-x_n\mid\mathcal F_n]\bigr|\le\frac CW,\qquad
 \E d_n^2\le\frac CW,\qquad
 d_n=x_{n+1}-x_n-\E[x_{n+1}-x_n\mid\mathcal F_n].
\]

\emph{Step 5: summation over the window.} Between the two lengths, the
difference of a real coordinate $x$ is the sum of $N_j$ conditional means
and of the martingale differences $d_n$, which are orthogonal in $L^2$. We
bound the first sum by Markov's inequality and the second by Chebyshev's
inequality, each at the level $\epsilon/2$. This gives, for every
$\epsilon>0$,
\[
 \Prob\bigl(|x_{r_j}-x_{R_j}|\ge\epsilon\bigr)
   \le\Big(\frac2\epsilon+\frac4{\epsilon^2}\Big)\frac{CN_j}{W_j}
   \longrightarrow0,
\]
since $N_j=o(W_j)$. A union bound over the finitely many real
coordinates proves that $\mathcal K^{(r_j)}-\mathcal K^{(R_j)}\to0$ in
probability, and Step~1 concludes the proof.
\end{proof}

 Fix finitely many bounded nonnegative continuous block profiles $q_a$.
For the moment omit the three actual collar blocks from each source.
Let $\zeta_a$ be independent circular Gaussian vectors on the remaining
physical path, with covariance $Q_a/W^2$, independent of the path.
Set
\[
 \mathcal K_W=(\zeta_a^*G_{\rm lead}\zeta_b)_{a,b\le p},
 \qquad B_W=-\eta^{-1}\Im\mathcal K_W\succeq0.
\]
The imaginary part of a matrix is $(K-K^*)/(2i)$.

\begin{lemma}[Global Gaussian source matrices from finite births]
\label{lem:t22-global-Gaussian-source}
The joint law of $\mathcal K_W,B_W$ has a unique critical limit
$\mathcal K,B$. Its positive explicit construction is the following.
For each probe $a$, put $m_a=\tau\int_\mathbb T q_a$ and sample
$M$ independent times in $[0,\tau]$ with density $q_a(t/\tau)/m_a$.
Run the complex finite-source equation of Theorem~\ref{thm:complex-source-law} with a mass-$m_a$ Gaussian tag at
each sampled time. At the terminal time sum its source matrix entries
within the probe groups, with coefficient $1/M$. Send $M\to\infty$.
Zero-mass groups are identically zero. All groups use the same current
noise. The convergence determines the full complex matrix law, not
only its positive determinant transforms.
\end{lemma}

Sampling produces ordinary finite source evolutions with shared noise. The error uniform in width allows the sampling cutoff to be removed after the critical limit.

\begin{proof}
At finite width, a sampled block tag $\xi$ has covariance $Q_a/W^2$
exactly, as in Section~\ref{sec:complex-sources}. Replace each $\zeta_a$ by
$M^{-1/2}\sum_{j=1}^M\xi_{aj}$, independently within its group.
For a bounded $C^3$ test of the finitely many real and imaginary
entries of the resulting source matrix, Lindeberg replacement costs
\begin{equation}
 C_{p,q,\tau}\|f\|_{C^3}/\sqrt M,                     \label{eq:global-probe-replacement}
\end{equation}
uniformly in $W$. Indeed a single replacement changes the matrix by
$tA+t^2C$. Its first two Taylor terms depend only on the matched
mean, covariance and zero pseudo-covariance. The exact third derivative is
\[
 D^3f(K(t))[A+2tC,A+2tC,A+2tC]
       +6D^2f(K(t))[A+2tC,C].
\]
For $0\le t\le1$ its absolute value is bounded by
$C\|f\|_{C^3}(\|A\|^3+\|C\|^3+\|A\|\|C\|+\|C\|^2)$.
In particular the quadratic $\|C\|^2$ term is present.
Equations~\eqref{eq:raw-tag-moments} and~\eqref{eq:extended-source-moments}
bound all these terms uniformly for every intermediate covariance mixture.
There are $pM$ replacements, each at $t=M^{-1/2}$.
No real-axis inverse moment of the closed ring is used.

For fixed $M$, Theorem~\ref{thm:complex-source-law} identifies the
finite tag matrix after two reductions, since the sampled tags have
$W$-dependent masses and blocks. The three collar blocks are the ring
blocks $n-1,0,1$. Grow the remaining path of $r=n-3$ blocks
$2,\dots,n-2$ from ring block $2$, with steps numbered from $0$ as in
Lemma~\ref{lem:t22-terminal-lengths}: step $k<r$ adds ring block $k+2$,
which is block $k+1$ in the labelling of Section~\ref{sec:full-source}.
Put $m_{a,W}=W^{-1}\sum_{k<r}q_a((k+2)/n)$. A sampled tag of
group $a$ has mass $m_{a,W}$ and sits at step $k$ with probability
$p_{a,W}(k)=q_a((k+2)/n)/(Wm_{a,W})$. If $m_a=0$, then $q_a\equiv0$ and
group $a$ vanishes before and after the limit. Let $m_a>0$; Riemann
sums and $r/W\to\tau$ give $m_{a,W}\to m_a$.

\emph{Masses.} Put $c_{a,W}=m_a/m_{a,W}\to1$ and
$C_W=\diag(c_{a,W})$. Tags of mass $m_a$ at the same random blocks have
covariance $c_{a,W}Q_a/W^2$. They therefore replace the probes
$c_{a,W}^{1/2}\zeta_a$, whose source matrix is
$C_W^{1/2}\mathcal K_WC_W^{1/2}$, and
\eqref{eq:global-probe-replacement} holds for the profiles
$c_{a,W}q_a$ with one constant for all large $W$. As $\mathcal K_W$ is tight by \eqref{eq:extended-source-moments},
$C_W^{1/2}\mathcal K_WC_W^{1/2}-\mathcal K_W\to0$ in probability, so all
tags may be given the fixed masses $m_a$.

\emph{Positions.} A time $t$ with density $q_a(t/\tau)/m_a$ on
$[0,\tau]$ places its tag at step $\lfloor Wt\rfloor$ if
$\lfloor Wt\rfloor<r$, and nowhere otherwise; let $\pi_{a,W}$ be the
law of this placement. For $k<r$ and $k/W\le t\le(k+1)/W$,
\[
 \Big|\frac t\tau-\frac{k+2}n\Big|
 \le\frac1{W\tau}+\frac r{n\tau}\Big|\frac nW-\tau\Big|+\frac2n,
\]
which tends to zero uniformly. The steps $k\ge W\tau-1$ number $o(W)$,
and a tag is placed nowhere with probability at most
$\|q_a\|_\infty(\tau-r/W)_+/m_a$. Uniform continuity of $q_a$ and
$m_{a,W}\to m_a$ therefore give
$\sum_k|\pi_{a,W}(k)-p_{a,W}(k)|\to0$. Coupling each of the $pM$ independent discrete steps with a sampled
time so that they differ with probability at most
$\sum_k|\pi_{a,W}(k)-p_{a,W}(k)|$, with the same path and tag Gaussians,
the two tag matrices differ at fixed $M$ only on an event of probability
$o(1)$.

\emph{Identification.} Almost surely the sampled times $t_{aj}$ are
distinct and lie in $(0,\tau)$. Conditionally on them, the tags are
those of Theorem~\ref{thm:complex-source-law}, with masses $m_a$ and
birth times $t_{aj}<\tau$, on the path of length $r=n-3$. Here
$r/W\to\tau$, but in general $r\ne\lfloor W\tau\rfloor$; by
Lemma~\ref{lem:t22-terminal-lengths}, for almost every family of times
the conditional law of the tag matrix converges to the critical limit
with these births at terminal time $\tau$. Grouping with coefficient $1/M$ is continuous, and Fubini's theorem and
dominated convergence in the times (the conditional limit is measurable
in them, as an almost everywhere limit) give, for every bounded
continuous test, the convergence of the $M$-sampled matrix to the
construction of the lemma at this $M$. Bound \eqref{eq:global-probe-replacement} and the same uniform source moments prove tightness
and uniqueness for all bounded smooth matrix tests. Approximation on
tight compact sets extends this to bounded continuous tests.
The exact imaginary-part identity passes to the limit.
\end{proof}

For specificity use a finite ambient space   $\C^{2\Naux }$, with
$\rank P=\rank P'=\Naux $, and the signature matrix $\mathsf D$
of Section~\ref{sec:positive-source-law}. Here   $\Naux $  is the freely padded
finite projection rank; it is unrelated to the physical dimension
  $N=nW$. The finite SDE is

\begin{equation}
 dA=\frac{\sqrt{D_0}}{2}P\,dZ\,P',\quad
 dV=dA^*V-\frac{D_0\Naux }{4}PV\,dt,\quad
 P'=P+V\mathsf DV^*,\quad
 d\mathcal K_{\rm tag}=-2i\eta V_{\rm pr}^*dA V_{\rm pr},
                                                               \label{eq:global-probe-sde}
\end{equation}
 
where $\mathbb E|dZ_{ij}|^2=2dt$. Probe columns have signature zero;
$\mathsf D$ is not the entire determinant $D$.
A mass-$m$ birth is a new perpendicular near vector of squared norm
$dm$, $d=(1-b)^{-1}$, and a new source diagonal $am$ with zero cross
entries. Theorem~\ref{thm:transmission-law} specifies the positive transmission entrance at a
source-free time $s>0$. The canonical two-projection truncations
and coupling \eqref{eq:rank-free-coupling} remove its current cutoff; Theorems~\ref{thm:source-process} and~\ref{thm:complex-source-law}
pay the actual comparison and let $s$ decrease below all sampled
births. Thus \eqref{eq:global-probe-sde} means a specified convergent positive construction,
not a formal infinite-dimensional stochastic equation.

\subsection{Inserting the actual physical closing law}\label{subsec:s25-inserting-the-actual-physical-closing-la}
 We now insert the physical collar without losing any source-dependent
determinant factor. The identity below expresses the complete closed
source matrix through the open-path matrix and the boundary inverse.
We derive it by Woodbury's formula, apply the common matrix closing
law, and obtain the limiting closed-source distribution.

 The following identity retains the full determinant quotient. For a
matrix of old physical source columns $\mathcal J$, put
$K=\mathcal J^*G_{\rm lead}\mathcal J$ and
$V=U_{\rm lead}^*\mathcal J$. Write $S$ for the lossless old scattering
matrix, $U$ for the independent physical collar and
$R_c=(I-US)^{-1}$. By Woodbury's identity and the periodic
closing determinant,
\begin{align}
 &\frac{\det(H_{\rm ring}-E+\mathcal J C\mathcal J^*)}{\det(H_{\rm ring}-E)}\\
 &\quad=\det(I+CK)\,
       \frac{\det(I-US_C)}{\det(I-US)}
       =\det\{I+C[K-2i\eta V^*(R_c-I)V]\}.             \label{eq:complete-closing-determinant}
\end{align}
The source has zero coordinates in the collar. All source-free Schur
factors cancel between numerator and denominator; no factor depending
on $C$ is discarded. The formulas initially hold where the indicated
inverses exist and extend as polynomial identities.
Consequently the exact closed physical source Green matrix is
\begin{equation}
 K^{\rm ring}=\Re K+\eta V^*\mathcal C_W V,
 \qquad \mathcal C_W=-i(2R_c-I),\qquad V^*V=-\Im K/\eta.
                                                               \label{eq:closed-source-matrix}
\end{equation}
Here $\Re K=(K+K^*)/2$. In particular the right side is Hermitian.

Apply this to the fixed number $p$ of global Gaussian probes.
Take an old orthonormal frame spanning their $V$ columns, padded to
rank $p$. This frame is independent of the fresh collar. Theorem~\ref{thm:matrix-closing} applies
uniformly, and Equation~\eqref{eq:extended-source-moments} controls $K$ and $V^*V$ on tight compact sets.
The limiting physical source matrix is therefore explicitly
\begin{equation}\label{eq:limiting-closed-source}
 \mathcal T=\Re\mathcal K+
                   \eta B^{1/2}\mathsf C_p(0)B^{1/2},
 \qquad B=-\Im\mathcal K/\eta,                                                \end{equation}
where $\mathsf C_p$ is independent of the path source matrix.
Compression consistency of \eqref{eq:canonical-cauchy} covers singular $B$; no lower bound
on its smallest eigenvalue is required. One common $\mathsf C_p$
is used for all probes in a mixed moment. Its entries must not be
replaced by independent scalar Cauchy variables.

\subsection{Restoring sources on all physical blocks}\label{subsec:s25-restoring-sources-on-all-physical-blocks}
 The probe construction omitted three collar blocks; we must restore
their sources under the original ring law. The estimate below shows
that their contribution vanishes in the critical limit. We exclude a
short spectral interval around the observation energy, bound the
remaining inverse traces by counting, and use cyclic symmetry to pay
for the collar mass.

 This step pays the finite collar under the original ring law. Let
$\Pi$ project onto its three physical blocks and $R=(H_{
\rm ring}-E)^{-1}$. Write $\Lambda=nW\rho(E)$ and
$\mathcal E_\epsilon=\{\min_\alpha|u_\alpha|>\epsilon\}$.
Theorem~\ref{thm:counting} gives
$\mathbb P(\mathcal E_\epsilon^c)\le C\epsilon$.
Its global intensity bound and total count $nW$ imply
\begin{equation}
 \mathbb E[\operatorname{Tr}|R|\,\mathbf1_{\mathcal E_\epsilon}]
 \le C\Lambda\{1+\log((nW)/\epsilon)\},\qquad
 \mathbb E[\operatorname{Tr}R^2\,\mathbf1_{\mathcal E_\epsilon}]
 \le C\Lambda^2/\epsilon.                              \label{eq:ring-truncated-inverses}
\end{equation}
For the first inequality split the unfolded integral at $nW$; beyond
that point use its total count rather than integrate a constant
density over an infinite logarithmic tail. No operator-norm event
or real principal-value cancellation is needed.

Cyclic invariance, including that of $\mathcal E_\epsilon$, gives
exactly $\mathbb E\operatorname{Tr}(\Pi f(R))\mathbf1_{\mathcal E_\epsilon}
=(3/n)\mathbb E\operatorname{Tr}f(R)\mathbf1_{\mathcal E_\epsilon}$
for $f(R)=|R|,R^2$. If a full global probe is split as
$\zeta=\zeta_o+\zeta_c$ on the old path and collar, conditional
Gaussian integration and \eqref{eq:ring-truncated-inverses} therefore yield
\begin{equation}
 \mathbb E|\zeta_c^*R\zeta_c|\mathbf1_{\mathcal E_\epsilon}
       \le C\frac{1+\log((nW)/\epsilon)}{n},\qquad
 \mathbb E|\zeta_o^*R\zeta_c|^2\mathbf1_{\mathcal E_\epsilon}
       \le C/(n\epsilon).                              \label{eq:collar-source-removal}
\end{equation}
The same bounds treat the finitely many mixed entries and different
bounded profiles, by conditioning and Cauchy--Schwarz. At fixed
$\epsilon$ they vanish as $W\to\infty$; then send $\epsilon$ to zero.
Thus \eqref{eq:limiting-closed-source} also identifies the entire joint matrix of probes with
covariances $Q_a/W^2$ on \emph{all} original ring blocks. This pays
both the omitted source mass and its interference with the old path.

\subsection{Determination and compactness}\label{subsec:s25-determination-and-compactness}

 We finish Theorem~\ref{thm:entire-law} by proving both uniqueness and
tightness of the determinant limit. Bounded inverse-determinant moments
determine its law, while Gaussian probe estimates control the entire
functions on compact sets. We establish these two facts separately,
then identify every subsequential limit.

 \begin{proof}[Proof of Theorem~\ref{thm:entire-law}]
At finite width the Gaussian integral gives \eqref{eq:determining-bounded-moments} exactly with the
physical probe matrix, and every inverse determinant on the positive
imaginary cone lies in the closed unit disk. Equation \eqref{eq:limiting-closed-source}, with
\eqref{eq:collar-source-removal}, gives its limiting joint mixed moments. Polynomials in real
and imaginary parts are dense in continuous functions on the compact
polydisk, so these moments determine their laws. This is bounded
moment determination; no Cauchy or lognormal moments are asserted.

We also need tight entire functions, not just their inverse values.
Here it follows from the physical probe law without unbounded moments.
If $A=A^*$ and $g$ is standard circular Gaussian, put
$X=g^*Ag$ and $s_A^2=(\operatorname{Tr}A)^2+\operatorname{Tr}A^2$.
Gaussian moments give
\begin{equation}
 \mathbb E_g X^2=s_A^2,\quad
 \mathbb E_g X^4\le72s_A^4,\quad
 \mathbb P_g\{|X|\ge s_A/\sqrt2\}\ge1/288.             \label{eq:probe-tightness-converse}
\end{equation}
The last inequality is Paley--Zygmund for $X^2$. Therefore tightness
of the actual one-probe variable $X$ implies tightness of both
$\operatorname{Tr}A$ and $\|A\|_{\rm HS}$ in its underlying matrix law.
Apply this to $A=W^{-2}Q_*^{1/2}RQ_*^{1/2}$ for
$Q_*=\sum Q_a$, and to the finitely many $Q_a$ separately. Their
one-probe variables are tight by \eqref{eq:limiting-closed-source}--\eqref{eq:collar-source-removal}.

Using Sylvester's identity and the commuting bounded diagonal
coefficients $Q_*^{-1/2}Q_aQ_*^{-1/2}$, the determinant bound
\[
 |\det(I+B_0)|\le
 \exp\{\Re\operatorname{Tr}B_0+\tfrac12\|B_0\|_{\rm HS}^2\}
\]
now gives uniform bounds in probability on each compact source set.
Montel gives tightness in the Polish space of entire functions.
Every limit has $D(0)=1$ and $|D(iy)|\ge1$ on the positive imaginary
cone. A countable dense set of its inverse values determines the
entire function by the multivariable identity theorem, applied
successively in each coordinate. The injective evaluation map has a
Borel inverse on its image. Hence the already determined bounded
mixed moments identify every entire subsequential law. This proves
full convergence and retains the complete first trace term. For consistency, deleting $q_a$, $a\ge1$, sets $z_a=0$ in
\eqref{eq:entire-determinant} at every finite width, and restriction to
a coordinate subspace is continuous for locally uniform convergence.
\end{proof}

\subsection{Recovering all energies and marks jointly}\label{subsec:s25-recovering-all-energies-and-marks-jointl}

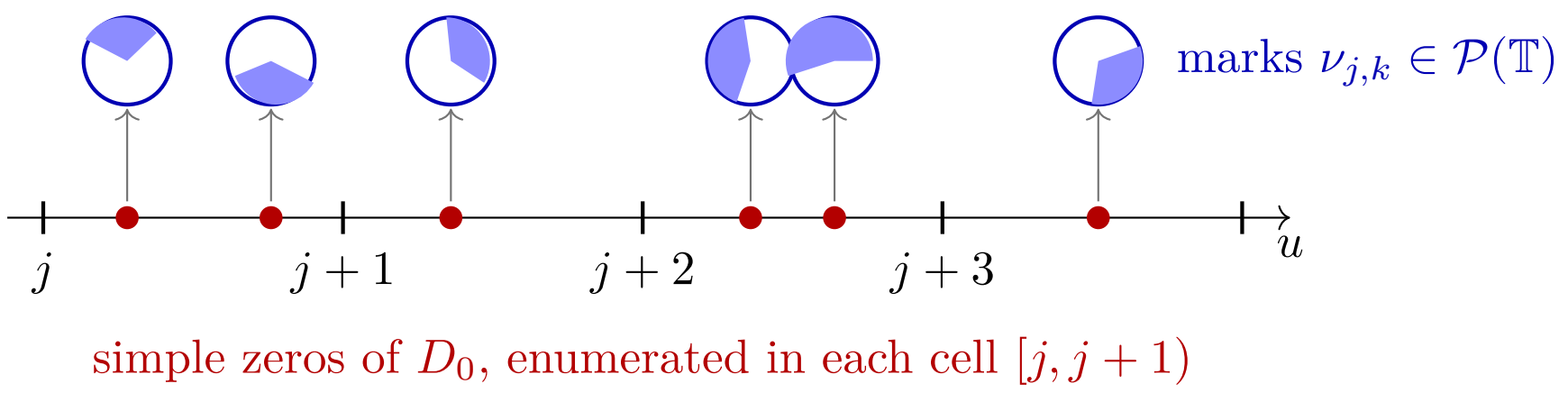
\begin{figure}[H]
\centering
\begin{tikzpicture}[x=.88cm,y=1cm]
 \draw[->] (-.3,0)--(10.4,0) node[below,font=\small] {$u$};
 \foreach \x in {0,2.5,5,7.5,10}{
   \draw[thick] (\x,-.12)--(\x,.12);
 }
 \node[font=\small,below] at (0,-.12) {$j$};
 \node[font=\small,below] at (2.5,-.12) {$j+1$};
 \node[font=\small,below] at (5,-.12) {$j+2$};
 \node[font=\small,below] at (7.5,-.12) {$j+3$};
 \foreach \x in {0.7,1.9,3.4,5.9,6.6,8.8}{
   \fill[red!70!black] (\x,0) circle (2.4pt);
   \draw[thick,blue!70!black] (\x,1.15) circle (.32cm);
   \draw[->,thin,black!55] (\x,.12)--(\x,.8);
 }
 \fill[blue!45] (0.7,1.15) ++(40:.32) arc (40:150:.32cm) -- (0.7,1.15) -- cycle;
 \fill[blue!45] (1.9,1.15) ++(200:.32) arc (200:330:.32cm) -- (1.9,1.15) -- cycle;
 \fill[blue!45] (3.4,1.15) ++(-30:.32) arc (-30:90:.32cm) -- (3.4,1.15) -- cycle;
 \fill[blue!45] (5.9,1.15) ++(100:.32) arc (100:250:.32cm) -- (5.9,1.15) -- cycle;
 \fill[blue!45] (6.6,1.15) ++(0:.32) arc (0:200:.32cm) -- (6.6,1.15) -- cycle;
 \fill[blue!45] (8.8,1.15) ++(260:.32) arc (260:380:.32cm) -- (8.8,1.15) -- cycle;
 \node[font=\small,red!70!black,align=center] at (5,-1.05)
   {simple zeros of $D_0$, enumerated in each cell $[j,j+1)$};
 \node[font=\small,blue!70!black,anchor=west] at (9.3,1.15)
   {marks $\nu_{j,k}\in\mathcal P(\T)$};
\end{tikzpicture}
\caption{The extraction in the proof of Theorem~\ref{thm:joint-recovery},
drawn schematically. The real zeros of $D_0(z)=D(z,0,\ldots,0)$ (red),
in the unfolded variable $u=-\tau\rho(E)z$, are the eigenvalue locations.
They are simple, and in each unit cell $[j,j+1)$ they are enumerated in
increasing order. At each zero, the ratios of partial derivatives of $D$
give the integrals $\nu_{j,k}(q_a)$ of the mark against the profiles
$q_a$. The proof shows that these numbers are, simultaneously for all
zeros, the values of a probability measure on the circle (blue), which is
then recovered by the Riesz representation.}
\label{fig:jr-extraction}
\end{figure}

Theorem~\ref{thm:entire-law} identifies
the law of the random entire function $D$ for each finite family of
profiles. The point of the next theorem is to turn this into the law
of the whole marked point process. The zeros of $D$ on the first axis
give the eigenvalue locations, and its derivatives at those zeros give
the values of the marks on the profiles, as in
Figure~\ref{fig:jr-extraction}. Three issues need care. The families of
profiles must fit together into one object. The zeros must be enumerated
measurably. And the numbers $\nu_\alpha(q_a)$ must come from a genuine
probability measure, for every zero simultaneously.
\begin{theorem}[Full critical crossover, original ring]
\label{thm:joint-recovery}
For the original complex Gaussian block ring, both parities, every
fixed $0<b<1/2$, $|E|<2$ and $n/W\to\tau\in(0,\infty)$,
\[
 \Xi_W\Rightarrow\Xi_{\tau,b,E}
\]
as a locally finite marked point process, with compact spectral tests
and weak probability-measure marks on the circle. The limit is unique,
has simple spectral locations and unit spectral intensity. It is
explicitly constructed by the global source law of Lemma~\ref{lem:t22-global-Gaussian-source}, realized through \eqref{eq:global-probe-sde}, the independent physical
matrix closing \eqref{eq:canonical-cauchy}, \eqref{eq:limiting-closed-source}, the determining hierarchy \eqref{eq:determining-bounded-moments}, and
the following root/derivative rule. If $z_\alpha$ is a real zero of
$z\mapsto D(z,0,\ldots,0)$, then
\begin{equation}
 u_\alpha=-\tau\rho(E)z_\alpha,\qquad
 \nu_\alpha(q_a)=
       \frac{\partial_{z_a}D(z_\alpha,0,\ldots,0)}{\partial_{z_0}D(z_\alpha,0,\ldots,0)}.       \label{eq:root-mark-extractor}
\end{equation}
Use a countable dense collection of smooth nonnegative profiles to
recover each probability measure. \end{theorem}

\begin{proof} At finite width, the zeros of the determinant ratio
are exactly the rescaled eigenvalues, and ratios of its derivatives at those
zeros are the mark integrals (Step~1). We keep this description in the limit as well. We first combine the laws for growing families of
profiles into one consistent law (Steps~2--3), and then list the zeros of the
limit in a measurable way (Step~4), as in Figure~\ref{fig:jr-extraction}. Next
we show that the derivative ratios at each zero define a positive functional
of total mass one, and hence a probability measure on the circle
(Steps~5--6). Step~7 checks that the resulting marked process is the limit of
the finite-width processes.

\emph{Step 1: the finite-width identities.} At finite width, the zeros of
\eqref{eq:entire-determinant} on the $z_0$ axis are exactly
$-W^2(\lambda_\alpha-E)$. So their unfolded values are
$-\tau_W\rho(E)z_\alpha$, with $\tau_W=n/W$. At each simple zero,
differentiating the determinant along a rank-one spectral projector gives
exactly the second identity in \eqref{eq:root-mark-extractor}, with the
physical block mass $\langle\psi_\alpha,Q_a\psi_\alpha\rangle$. This is the
complete original determinant. So it introduces neither a spectral Palm
weight nor a new energy marginal.

\emph{Step 2: a countable family of profiles.} We now give the
simultaneous construction, including its null sets and its
measurability. Let $\mathcal V$ be the countable rational vector space of
real trigonometric polynomials with rational coefficients, which
contains $1$. We enumerate its nonnegative cone as $q_0=1,q_1,q_2,\ldots$.
This cone is uniformly dense in the nonnegative continuous functions:
we approximate by a rational trigonometric polynomial and add a small
rational constant to make it nonnegative. The rational differences of
the cone span $\mathcal V$. For the initial family $q_0,\ldots,q_m$, let
$\mu_m$ be the law of Theorem~\ref{thm:entire-law} on the Polish space
$\mathcal H_m=\operatorname{Hol}(\C^{m+1})$, with the topology of local
uniform convergence.

\emph{Step 3: a common projective family.} The different families of
profiles must give one consistent object. The restriction
\[
 \pi_m F(z_0,\ldots,z_m)=F(z_0,\ldots,z_m,0)
\]
is continuous. At finite width we have exactly
$\pi_mD_W^{(m+1)}=D_W^{(m)}$, because setting $z_{m+1}=0$ removes the
profile $q_{m+1}$. So the convergence gives $(\pi_m)_*\mu_{m+1}=\mu_m$.

There is a unique probability law of a nested sequence $(D^{(m)})_{m\ge0}$
with these marginals, in which the first $m+1$ coordinates are the
restrictions of the last one. To see existence directly, we take the
actual finite-width sequences in $\prod_m\mathcal H_m$. Marginal tightness
implies product tightness: we choose compact sets $K_m$ whose complements
have probability at most $\epsilon2^{-m-1}$, and we use the compactness of
$\prod_mK_m$. Prokhorov's theorem gives a subsequential product law. The
restriction identities are closed conditions, so they persist in the
limit. The prescribed finite-coordinate laws make the product law unique.
This also proves simultaneous convergence along the original width
sequence. All algebraic identities induced by rational relations among
the $q_a$ are identities between continuous restrictions along linear
source maps, so they pass to this common law.

\emph{Step 4: the zeros and their measurable enumeration.} We write
$D_0(z)=D^{(0)}(z)$, so that $D_0(0)=1$. By Hurwitz's theorem, $D_0$ has no
nonreal zeros, and its real zero measure, with multiplicities, agrees
with every subsequential spectral limit under $u=-\tau\rho(E)z$. The count
estimate and its consequence for close pairs, \eqref{eq:close-pairs},
exclude multiple limiting locations on each compact interval. The bound
on the first-moment density excludes any fixed rational location. Thus,
on one event of probability one, the zeros are simple and locally finite,
and no rational energy is a zero.

The real zero counting measure is a Borel function of $D_0$. Indeed, on
compact contours that avoid the zeros, its counts are given by the
integrals of the argument principle, and a countable set of rational
contours exhausts the complement of the discrete zero set. Equivalently,
by Hurwitz's theorem the map is continuous near simple zeros away from
the endpoints of intervals. In each half-open energy cell $[j,j+1)$,
$j\in\mathbb Z$, we enumerate its finitely many zeros in increasing order.
The $k$-th location is
\[
 u_{j,k}=\inf\{r\in\mathbb Q\cap(j,j+1):
                         N_D([j,r])\ge k\},
 \qquad 1\le k\le N_D([j,j+1)).
\]
This is a countable measurable enumeration, and cells with no zeros
contribute no point. We put $z_{j,k}=-u_{j,k}/(\tau\rho(E))$.

\emph{Step 5: positive functionals of mass one.} For each enumerated zero
we define
\[
 \ell_{j,k}(q_a)=
 \frac{\partial_{z_a}D^{(m)}(z_{j,k},0,\ldots,0)}{\partial_{z_0}D^{(m)}(z_{j,k},0,\ldots,0)},\qquad m\ge a.
\]
By the restriction consistency of Step~3, the definition does not depend
on $m$. The denominator is nonzero because the zero is simple.

We must check that $\ell_{j,k}$ behaves like integration against a
probability measure, on one event for all zeros. For this we take a
joint subsequential limit of this countable entire family and of the
original marked processes. The latter are tight, since their mark space
$\mathcal P(\T)$ is compact and their local counts are tight. We use a
Skorokhod realization in these Polish spaces. For every compact energy
interval with no zero at its boundary, Hurwitz's theorem matches its
finitely many roots, and Cauchy's formula gives convergence of all
derivatives at them. At finite width, the extractor is exactly
integration against a probability measure, by Step~1. Consequently,
simultaneously for every zero and every rational relation in the
countable space $\mathcal V$,
\[
 \ell_{j,k}(1)=1,\qquad \ell_{j,k}(q)\ge0\ \ (q\ge0),\qquad
 |\ell_{j,k}(q)-\ell_{j,k}(\widetilde q)|
                    \le\|q-\widetilde q\|_\infty,
\]
and additivity and rational homogeneity hold. This verification is a
statement about the support of the entire-family law, which is already
unique. It introduces no choice of subsequence into the final
construction.

\emph{Step 6: from functionals to marks.} We extend $\ell_{j,k}$ first to
the rational differences $\mathcal V$, and then continuously to
$C(\T;\R)$. The Lipschitz bound of Step~5 makes this a well-defined real
linear functional of norm one, and positivity persists by uniform
approximation. The Riesz representation theorem then gives a unique
$\nu_{j,k}\in\mathcal P(\T)$. This representation is measurable. Indeed,
the map $\nu\mapsto(\nu(q_a))_{a\ge0}$ is a continuous injection from the
compact space $\mathcal P(\T)$ into a countable product of real lines, so
its inverse on its image is continuous. The coordinates here are Borel
ratios of derivatives. We assign an arbitrary fixed measure off the
admissible set of probability one, and we define
\[
 \Xi=\sum_{j\in\mathbb Z}\sum_{k=1}^{N_D([j,j+1))}
                             \delta_{(u_{j,k},\nu_{j,k})}.
\]
This is a measurable, locally finite marked point process.

\emph{Step 7: compatibility with vague convergence.} Finally we show that
$\Xi$ is the limit of the finite-width processes. In the joint
realization of Step~5, each compact interval with no boundary point has
finitely many matched roots. Convergence of every coordinate $q_a$, and
the compactness of $\mathcal P(\T)$, imply weak convergence of each matched
mark, since this family separates measures. So every continuous
compactly supported test on $\R\times\mathcal P(\T)$ has a convergent sum
on that interval. We exhaust the real line by such intervals. Hence every
joint subsequential marked limit equals the measurable construction
above, whose law is unique. This proves convergence of the full marked
process, and not only of its finite lists of spatial integrals.

Unit spectral intensity and stationarity follow from
Theorem~\ref{thm:critical-density} and Proposition~\ref{prop:stationarity},
with the uniform integrability given by the count moments. The same
count moments give uniform integrability of the weighted quadratic
profile statistic.

This proves Theorem~\ref{thm:critical}. The variance and the spectral
endpoints are established in the following sections from this unique
law.
\end{proof}

The common source noise and common matrix closing retain the correlations
between distinct energies and profiles. The construction does not imply
independence of marks or independence between marks and eigenvalues.
\section{The smoothed spectral and profile hierarchy}
\label{sec:smoothed-hierarchy}
This section expresses smoothed eigenvalue and profile observables directly
through the limiting random entire function. It gives explicit formulas for
the joint moments of the critical law, which are not needed elsewhere in the
proof.
Theorem~\ref{thm:smoothed-hierarchy} states that the joint moments of the
smoothed observables converge at every fixed order, with limits given by
source derivatives of $D$.
We write the source derivatives with the physical normalization, use the
counting estimates for uniform integrability, and remove the smoothing when we
recover the marked observables.

For $\delta>0$ let $p_\delta(x)=\delta/[\pi(x^2+\delta^2)]$.
For a compactly supported continuous $f:\R\to\R$ and a continuous
real $q:\T\to\R$, put
\begin{equation}\label{eq:smoothed-observable}
 \zeta_{W,\delta}(f\otimes q)
 =\int_\R f(u)\sum_\alpha p_\delta(u-u_\alpha)\nu_\alpha(q)\,du.
\end{equation}
Extend this definition linearly to finite sums of product tests.
Every real spatial test is the difference of two nonnegative ones,
so the determinant construction applies to the required source family.

\begin{theorem}\label{thm:smoothed-hierarchy}
Fix $p<\infty$ a positive integer $p$, $\delta_j>0$, compactly supported continuous $f_j$,
and continuous real $q_j$. Construct $D$ using the source family
$1,q_1,\ldots,q_p$, interpreting signed sources by positive
decomposition. Define
\begin{equation}\label{eq:limiting-smoothed-derivative}
 Z_{\delta,j}(D)
 =\frac{1}{\pi\tau\rho(E)}
   \int_\R f_j(u)\Im\!
    \left[
     \frac{\partial_{z_{q_j}}D}{D}
       \left(-\frac{u+i\delta_j}{\tau\rho(E)},0,\ldots,0\right)
    \right]du.
\end{equation}
Then the joint smoothed observables converge in law and in all fixed
mixed moments:
\begin{equation}\label{eq:smoothed-hierarchy-law}
 \lim_{W\to\infty}\E\prod_{j=1}^p
       \zeta_{W,\delta_j}(f_j\otimes q_j)
       =\E\prod_{j=1}^p Z_{\delta,j}(D).
\end{equation}
The same assertion holds for arbitrary compactly supported continuous
tests on $\R\times\T$, by uniform approximation. These moments and their
positive consistent construction determine the smoothed hierarchy
of $\Xi_{\tau,b,E}$. Taking the smoothing parameters to zero recovers
the unsmoothed observables and the joint marked law.
\end{theorem}

The order $p$ and the positive smoothing heights are fixed before
the critical limit. No moment of a real Cauchy boundary value is
used. In particular, the assertion follows from the entire determinant
law and the counting bounds, rather than from a formal differentiation
of an averaged log determinant.

\begin{proof}
Write   $\Lambda=N\rho(E)$  and $\tau_W=n/W$.
For $w\notin\R$, set $z_0(w)=-w/(\tau_W\rho(E))$.
The exact finite determinant identity gives
\begin{align}
 \frac{\partial_{z_q}D_W}{D_W}(z_0(w),0)
 &=W^{-2}\Tr\!\left[Q_q
       (H-E+z_0(w)W^{-2})^{-1}\right]\notag\\
 &=\tau_W\rho(E)\sum_\alpha
                   \frac{\nu_\alpha(q)}{u_\alpha-w}.
                                      \label{eq:finite-source-derivative}
\end{align}
Indeed the spectral denominator in the first line is
$(u_\alpha-w)/\Lambda$, and $\Lambda/W^2=\tau_W\rho(E)$.
For $w=u+i\delta$ its reciprocal has positive imaginary part.
Multiplication by $f(u)$ and integration therefore prove
\eqref{eq:limiting-smoothed-derivative} at finite width with
$D_W,\tau_W$ in place of $D,\tau$.

The restriction $z\mapsto D_W(z,0)$ has only real zeros.
Local uniform convergence and $D(0)=1$ imply, by Hurwitz's theorem,
that $D(z,0)$ has no nonreal zero. On every compact set a positive
distance from the real axis, local uniform convergence of the entire
functions also gives convergence of all first source derivatives.
Division by $D$ is continuous there. Thus the finite integrals
in \eqref{eq:limiting-smoothed-derivative} converge jointly in law.

It remains to justify the expectations. Set
$S_{W,\delta}(u)=\sum_\alpha p_\delta(u-u_\alpha)$.
The fixed-regulator count bound \eqref{eq:smoothed-moments} gives,
for every fixed integer $r$ and $\delta_*=\min_j\delta_j>0$,
\[
 \E\prod_{\ell=1}^r S_{W,\delta_{j_\ell}}(u_\ell)
       \le r!C_{\delta_*}^{\,r}.
\]
Since
$|\nu_\alpha(q)|\le\|q\|_\infty$, Tonelli yields
\begin{equation}\label{eq:hierarchy-product-bound}
 \E\prod_{j=1}^p
       |\zeta_{W,\delta_j}(f_j\otimes q_j)|
 \le p!C_{\delta_*}^{\,p}
                     \prod_{j=1}^p\|f_j\|_1\|q_j\|_\infty.
\end{equation}
Apply the same estimate to the $2p$ factors obtained by repeating
each test twice. This bounds the second moment of the product
in \eqref{eq:smoothed-hierarchy-law}, uniformly in $W$.
The products are uniformly integrable, so convergence in law proves
\eqref{eq:smoothed-hierarchy-law}.

For a continuous test $F$ supported in $[-R,R]\times\T$,
finite sums of products $f(u)q(r)$ approximate $F$ uniformly
on a common compact set. The error in its smoothed observable is
bounded by the uniform approximation error times
$\int_{-R-1}^{R+1}S_{W,\delta}(u)\,du$.
The corresponding fixed moments are bounded by
\eqref{eq:hierarchy-product-bound}. This proves the assertion for
general $F$, including mixed products of finitely many such tests.

Positivity for nonnegative tests and consistency under removal of
tests hold at every finite width and pass to the stated limits.
Finally Lemma~\ref{lem:smoothing} pays the distant spectral tails
and uniformly removes the Poisson smoothing against compact
continuous tests. Applied first at fixed height and then as the
height decreases to zero, it identifies the hierarchy with that of
the marked process recovered in Theorem~\ref{thm:joint-recovery}.
The exponential count bound \eqref{eq:count-exponential} determines
the finite compact-window laws from these moments. An exhaustion of
the spectral line and a countable dense family of spatial tests
then determine the full marked law.
\end{proof}

 \part{}
\section{A signed response and positive profile variance}
\label{sec:t22-signed-response}
This section begins Part~III. It shows that the critical eigenvector
profiles keep nonzero spatial fluctuations, which proves the positivity of
the variance in Theorem~\ref{thm:endpoints}.
Theorem~\ref{thm:t22-response-variance} bounds the square of a signed
response, the response of the transmission diffusion to an oscillating spatial
weight, by the corresponding Fourier variance of the profiles.
Proposition~\ref{prop:t22-nonzero-response} shows that this response is
nonzero in every Fourier mode. Together they show that the profile variance
is strictly positive.
We first justify differentiation under the physical path law
(Lemma~\ref{lem:t22-pass-derivative}) and derive the variance bound. We then
prove the variance assertions of Theorem~\ref{thm:endpoints}. In
Subsection~\ref{subsec:s27-the-scalar-lifetime-bridge} we represent the signed response through a scalar
diffusion bridge, and show that it is positive
(Lemmas~\ref{lem:t22-lifetime-bridge} and~\ref{lem:t22-bridge-occupation}).

Use the actual path law of Section~\ref{sec:t22-common-port-profile},
with $b$ fixed and $E$ in a fixed compact bulk set; let $j/W\to\tau\in(0,\infty)$.
For a bounded real continuous spatial test $q$ on $[0,\tau]$, put
\[
 A=W^{-2}U^*P_{\rm interior}U,\quad
 B_q=W^{-2}U^*Q_qU,
 \qquad U=G_0V,
\]
where the source is on physical interior blocks only. For $q\ge0$ insert
$i\kappa Q_q/W^2$ in their diagonal, and denote the resulting scattering
matrix by $\mathbb S_\kappa=I-2i\eta V^*G_\kappa V$. Write
\[
 T_L(\kappa)=I-S_{11}(\kappa)^*S_{11}(\kappa),\quad
 T_R(\kappa)=I-S_{22}(\kappa)^*S_{22}(\kappa),
\]
\[
 F_L(\kappa)=\frac{\det(I-yT_L(\kappa))}{\det(I+xT_L(\kappa))},
 \quad F_R(\kappa)=\frac{\det(I-yT_R(\kappa))}{\det(I+xT_R(\kappa))},
 \quad x\ge0,\quad0\le y<1.
\]
All source derivatives below are right derivatives along $q\ge0$.
For signed $q$, their difference along $q_+$ and $q_-$ defines the
linear derivative $\partial_q$ at zero. No negative absorber is used.

\begin{lemma}[Uniform passage of the physical derivative]
\label{lem:t22-pass-derivative}
For every fixed bounded nonnegative $q$ and $0<s<1$, the expected
right difference quotients of $F_L+F_R$ differ from their expected
zero-source derivatives by $O(\kappa^{s/2}+\kappa)$, uniformly in $W$.
Their critical derivative is
\begin{align}
 \lim_{W\to\infty}\mathbb E\partial_q(F_L+F_R)
 =-4\pi\rho(E)\,\mathbb E_{x,y}\left[
 \mathbf1_{\{\sigma_0\le\tau\}}
 \int_0^\tau\{q(s)+q(\tau-s)\}R_s\,ds\right].       \label{eq:t22-signed-limit}
\end{align}
Here $(X_s,Y_s)$, $R_s=X_s+Y_s$ and $\sigma_0$ are precisely the $p=1$ diffusion of Section~\ref{sec:transmission-law},
with absorption at its corner. The convergence is
uniform for deterministic energies in the actual bounded unfolded window. This means uniformity for $E'=E+u/(nW\rho(E))$, $|u|\le U$, with fixed $U$ where the matched lead is evaluated at $E'$ as well.  These energies satisfy $|E'-E|\le\delta_W$ with $\delta_W\to0$,
uniformly in $|u|\le U$, so Lemma~\ref{lem:source-recorded-estimates}
applies along them.
\end{lemma}

The derivative is passed through the physical limit with a uniform difference-quotient error. No unbounded real-axis inverse moment is assumed.

\begin{proof}
Let $M=Q_q^{1/2}G_0V$ and $K_q=Q_q^{1/2}G_0Q_q^{1/2}$.
The resolvent and Ward identities give exactly
\begin{equation}
 \mathbb S_\kappa-\mathbb S_0
 =-\frac{2\eta\kappa}{W^2}\mathbb S_0M^*
      (I+i\kappa K_q/W^2)^{-1}M,\qquad
 \partial_q\mathbb S=-2\eta\mathbb S_0B_q.             \label{eq:signed-source-update}
\end{equation}
Since $-\Im K_q\succeq0$, the middle inverse is a contraction. Therefore
\begin{equation}
 \left\|\frac{\mathbb S_\kappa-\mathbb S_0}{\kappa}
                   +2\eta\mathbb S_0B_q\right\|_1
 \le2\eta\operatorname{Tr}B_q
          \min(2,\kappa\|K_q\|/W^2).                \label{eq:signed-source-error}
\end{equation}
Here $\|\cdot\|_1$ is trace norm; $B_q\succeq0$ in this paragraph.

A fractional norm estimate pays this error without a microscopic local
law. Schur out the two fixed auxiliary leads, conditional on their
Gaussian couplings. The remaining inverse is
$(H_{\rm path}+D-e+iC)^{-1}$, with $D=D^*$ and $C\succeq0$ independent
of the physical Gaussian diagonal blocks. For real $t$, the
uniform deformed-block Wegner estimate \eqref{eq:external-count} gives

\[
 \mathbb P\{\operatorname{dist}(e,\operatorname{spec}
        (H_{\rm path}+D+tC))<a\}\le\min(1,CN a).
\]

Integrating this tail yields a uniform fractional inverse norm bound
  $C_sN^s$  for $0<s<1$. The norm to power $s$ of the analytic inverse
is subharmonic in the upper $t$ half-plane. Apply its scalar Poisson
inequality at $t=i$ and then Tonelli; the same bound follows with
$iC$. To justify possible real poles, first add a scalar positive height,
use the identical Wegner bound, and pass by Fatou. Schur reconstruction
costs only fixed moments of the endpoint-coupling norms. Consequently
\begin{equation}
 \mathbb E(\|G_0\|/W^2)^s\le C_{s,\tau},\qquad 0<s<1.          \label{eq:path-fractional-inverse}
\end{equation}
This bound holds under the actual path law, uniformly in the microscopic
energy window. It does not assert an inverse moment of a reconstruction
clock or of a smallest Ginibre singular value.

Use $\min(2,a)^2\le C_s a^s$, \eqref{eq:path-fractional-inverse}, Cauchy--Schwarz and
\eqref{eq:full-clock-moments}. The expected right side of
\eqref{eq:signed-source-error} is $O(\kappa^{s/2})$. The reflection defects are contractions
for every $\kappa\ge0$. Their first differences are bounded in trace
norm by $C\kappa\operatorname{Tr}B_q$ and their quadratic correction
by $C\kappa^2(\operatorname{Tr}B_q)^2$.
On the convex set $0\preceq T\preceq I$, the first and second
variations of $F_{x,y}$ are bounded respectively by
$C_{x,y}\|D\|_1$ and $C_{x,y}\|D\|_1^2$. This follows by
differentiating its log determinant, using $y<1$, and $0\le F\le1$.
The full-clock second moment therefore bounds the remaining expected
difference-quotient error by $C\kappa$.

The bounded-transform assertion of Theorem~\ref{thm:transmission-law}
supplies the limiting value at each fixed positive $\kappa$; no forced
fractional or no-small-channel assertion is used.
Apply it once from each end; reversing the layer order replaces $q(s)$
by $q(\tau-s)$. Its killing is $4\pi\rho(E)\kappa q$.
The right derivative of its bounded exponential is obtained by
monotone convergence, since
$\mathbb E\int_0^\tau R_sds=\tau(x+y)<\infty$.
First pass $W\to\infty$ at fixed $\kappa$, then let
$\kappa\downarrow0$ using the uniform error just proved. Subtract the
two positive-source derivatives for signed $q$. This proves
\eqref{eq:t22-signed-limit}. The actual unfolded-energy uniformity of
all inputs allows integration against any fixed compact spectral weight. These inputs are \eqref{eq:full-clock-moments} at $E'$ and the
bounded-transform assertion of Theorem~\ref{thm:transmission-law} along
$E'$, with the limit at $E$, uniformly in $|u|\le U$
(Lemma~\ref{lem:source-recorded-estimates}(e),(h)), and
\eqref{eq:path-fractional-inverse}, which is already uniform in the window.
\end{proof}

\begin{theorem}[Physical Fourier variance bounds the signed response]
\label{thm:t22-response-variance}
For fixed $m\in\mathbb Z\setminus\{0\}$ let
\[
 V_{W,m}=\mathbb E\sum_\alpha w(u_\alpha)
       |\nu_\alpha(e^{2\pi im\cdot})|^2,
 \qquad w\ge0,\quad w\in C_c^\infty,\quad\int w=1.
\]
Use $q_m(s)=\cos(2\pi ms/\tau)$ and define the finite explicit number
\[
 J_{\tau,m}(x,y)=\mathbb E_{x,y}\left[
 \mathbf1_{\{\sigma_0\le\tau\}}
          \int_0^\tau q_m(s)R_sds\right].
\]
For every fixed $x\ge0$, $0\le y<1$, there is a finite constant
$C_{x,y,b,E,\tau}$ such that
\begin{equation}
 (8\pi\rho(E))^2|J_{\tau,m}(x,y)|^2
       \le C_{x,y,b,E,\tau}\liminf_{W\to\infty}V_{W,m}.   \label{eq:response-variance-bound}
\end{equation}
Thus any nonzero value of this explicit scalar response proves a
strict positive physical Fourier-variance liminf. The bridge argument below will make this response nonzero.
\end{theorem}

A nonzero signed source response forces positive physical profile variance. The constants are finite at each fixed choice of the dual variables.

\begin{proof}
At zero absorption $\mathbb S=\mathbb S_0$ is unitary. Let $P$ be
its near-port orthogonal projection and $P'=\mathbb S^*P\mathbb S$.
The left and right transmission defects have identical spectra, so
$F_L(0)=F_R(0)=F$. Put
\[
 r(t)=\frac{y}{1-yt}+\frac{x}{1+xt},\qquad r_0=x+y,
 \qquad
 C=P'r(T_L)P+(I-P')r(T_R)(I-P),
\]
where each function is embedded on its port subspace. Differentiating
\eqref{eq:signed-source-update} and both determinant ratios gives the exact identity
\begin{equation}
 \partial_q(F_L+F_R)=-4\eta F\Re\operatorname{Tr}(B_qC).       \label{eq:signed-ratio-derivative}
\end{equation}
Writing $g=\operatorname{Tr}T_L$, the projection identity gives
$\|P-P'\|_{\rm HS}^2=2g$. The constant part of $C-r_0I$ is
$r_0(P'-P)(2P-I)$; the rest has Hilbert--Schmidt norm at most
$C_{x,y}(\operatorname{Tr}T_L^2+\operatorname{Tr}T_R^2)^{1/2}$.
Since $0\preceq T_{L,R}\preceq I$,
\begin{equation}
 \|C-r_0I\|_{\rm HS}\le C_{x,y}\sqrt g,
 \quad
 |\partial_q(F_L+F_R)|
       \le C(|\operatorname{Tr}B_q|+\sqrt g\|B_q\|_{\rm HS}).   \label{eq:signed-response-trace}
\end{equation}
The scalar part retains the signed trace needed for this estimate.

Apply the following positive Gaussian quadratic-form inequality
to the complete matched-lead clock. For a standard complex Gaussian $h$,
\begin{equation}
 \mathbb E_h\left|\frac{h^*B_qh}{h^*Ah}\right|^2
 \ge\frac{|\operatorname{Tr}B_q|^2+\operatorname{Tr}B_q^2}{24(\operatorname{Tr}A)^2}.                      \label{eq:gaussian-variance-lower}
\end{equation}
To verify this inequality, scale $\Tr A=1$, put $X=h^*B_qh$,
$Y=h^*Ah$, and $M=\E_h|X|^2=|\Tr B_q|^2+\Tr B_q^2$.
Wick's formula for $\E_h|X|^2Y^2$ sums over the $24$ permutations
of $B_q,B_q^*,A,A$. Each cycle product is bounded in modulus by $M$:
a cycle containing both $B_q,B_q^*$ is bounded by $\|B_q\|_{\HS}^2$;
separate cycles are bounded by $|\Tr B_q|^2$, by $\|B_q\|_{\HS}^2$,
or by their arithmetic mean. Here $\Tr A=1$ implies
$\|A\|,\|A\|_{\HS}\le1$.
Thus $\E_h|X|^2Y^2\le24M$, and
\[
 M^2\le\E_h(|X|^2/Y^2)\,\E_h(|X|^2Y^2)
       \le24M\E_h(|X|^2/Y^2).
\]
Divide when $M>0$; the zero case is immediate. Restoring the trace
scale proves the displayed inequality.

The source is on interior coordinates. The auxiliary port masses in Section~\ref{sec:t22-common-port-profile} may be
removed: their physical Palm expectation is bounded by a fixed multiple
of $\mathbb E\sum_\alpha |w(u_\alpha)|\sum_kp_{\alpha k}^2=o(1)$,
as in \eqref{eq:vanishing-block-atom}, and Theorem~\ref{thm:t22-common-port-ratios} transfers that bound.
Alternatively retain the port clock in \eqref{eq:gaussian-variance-lower}; its added trace is
$O(W^{-1})$ after division by $W^2$, and its bounded mass ratio vanishes
by the same argument. The cosine square is bounded by the squared
complex Fourier mode for every Gaussian profile.

Cauchy--Schwarz in the actual path and energy measure, \eqref{eq:signed-response-trace}--\eqref{eq:gaussian-variance-lower},
and Theorem~\ref{thm:t22-physical-lead-profile} give
\[
 \left|\int w(u)\mathbb E\partial_{q_m}(F_L+F_R)(e(u))du\right|
 \le C\sqrt{V_{W,m}+o(1)}.
\]
The multiplier costs
$\mathbb E[(\operatorname{Tr}A)^2(1+g)]$, bounded by the full-clock
fourth moment and the positive-time conductance second moment in Section~\ref{sec:transmission}.
No unbounded comparison observable is admitted by total variation alone.
The negligible endpoint/source-location changes are controlled by
$\|q_m^{(W)}-q_m\|_\infty\operatorname{Tr}A$.
Finally $q_m(\tau-s)=q_m(s)$, so Lemma~\ref{lem:t22-pass-derivative}
identifies the left-hand limit with $-8\pi\rho(E)J_{\tau,m}(x,y)$.
Taking a subsequence realizing the variance liminf proves \eqref{eq:response-variance-bound}.
\end{proof}

\subsection{The scalar lifetime bridge}\label{subsec:s27-the-scalar-lifetime-bridge}

 We find a positive representation of the signed response so that its
Fourier contribution cannot vanish identically.
Lemma~\ref{lem:t22-lifetime-bridge} identifies the dual diffusion's
lifetime and conditional path law. We compute the corner-hitting
density, study where the bridge spends its time, and use this to prove
Proposition~\ref{prop:t22-nonzero-response}.

 \begin{lemma}[Corner-lifetime bridge]
\label{lem:t22-lifetime-bridge}
Let $\mathcal D_0$ be the independent Legendre--Jacobi base generator in Section~\ref{sec:transmission-law},
with transition density $p_t^0(z,z')$ relative to $dx'\,dy'$, and let
$H(x,y)=(x+y)^{-1}$. For $z=(x,y)\ne(0,0)$ the lifetime of its
absorbing $H$ transform has density
\begin{equation}
 \mathbb P_z(\sigma_0\in dt)
       =\frac{D_0}{2}(x+y)p_t^0(z,(0,0))\,dt.          \label{eq:corner-lifetime-density}
\end{equation}
It is strictly positive for every $t>0$. Conditional on $\sigma_0=t$,
the transformed path up to $t$ has exactly the base bridge from $z$ to
$(0,0)$ of length $t$. In particular its two coordinates are independent
base bridges with these endpoints.
\end{lemma}

The dual diffusion reaches its absorbing corner at a random time. Its lifetime density gives a positive integral representation for the signed response.

\begin{proof}
The base generator is in divergence form for Lebesgue measure:
\[
 \mathcal D_0=\frac{D_0}{2}\{\partial_x[x(1+x)\partial_x]
                 +\partial_y[y(1-y)\partial_y]\}.
\]
Away from the corner it annihilates $H$. The outward flux of
$(x(1+x)H_x,y(1-y)H_y)$ through $x+y=\epsilon$ is
\[
 -\int_0^\epsilon
       \frac{\epsilon+x^2-(\epsilon-x)^2}{\epsilon^2}\,dx=-1.
\]
The other edges have zero flux. Thus, distributionally on the closed
strip with its natural boundary convention,
$-\mathcal D_0H=(D_0/2)\delta_{(0,0)}$.
Convolving with the conservative base heat kernel gives
\begin{equation}
 H(z)-P_t^0H(z)=\frac{D_0}{2}\int_0^t p_s^0(z,(0,0))\,ds.     \label{eq:corner-heat-identity}
\end{equation}
This can also be obtained by smoothing $H$ outside $x+y<\epsilon$,
using the computed flux, and then sending $\epsilon\downarrow0$.
The corner singularity is locally integrable for $dx\,dy$.
A cutoff at large $x$ causes no contribution: the one-dimensional
hyperbolic radial representation below has Gaussian tails in $\log x$
at each positive time. These tails dominate the cutoff derivatives.

The base coordinate $X=(\cosh r-1)/2$ is the radial part of hyperbolic
Brownian motion in dimension two with generator
$(D_0/2)(\partial_r^2+\coth r\,\partial_r)$.
The coordinate $Y=(1-\cos\theta)/2$ is the radial part of spherical
Brownian motion in dimension two with generator
$(D_0/2)(\partial_\theta^2+\cot\theta\,\partial_\theta)$.
Their invariant radial area measures become constant multiples of
$dx$ and $dy$. Their kernels are positive for positive times, including
the continuous density values at the poles. For the hyperbolic factor
this follows also from the bounded-potential Bessel comparison below;
for the spherical factor it is positivity of the smooth elliptic heat
kernel on the connected sphere. Independence proves positivity in \eqref{eq:corner-lifetime-density}.

The survival probability of the $H$ transform is
$P_t^0H(z)/H(z)$. Equation \eqref{eq:corner-heat-identity} gives \eqref{eq:corner-lifetime-density}. Finally multiply any
finite string of its killed transition densities by the remaining
lifetime density. The factors $H(z_i)$ telescope. After dividing by
\eqref{eq:corner-lifetime-density}, the result is
\[
 \frac{p_{t_1}^0(z,z_1)p_{t_2-t_1}^0(z_1,z_2)\cdots
                  p_{t-t_k}^0(z_k,(0,0))}{p_t^0(z,(0,0))}\,dz_1\cdots dz_k,
\]
which is precisely the base bridge. Its continuous version follows
from the diffusion bridge construction, or its endpoint conditioning
limit. The product base kernel proves coordinate independence.
\end{proof}

\begin{lemma}[Weighted occupation is concentrated near the initial end]
\label{lem:t22-bridge-occupation}
Fix $\tau>0$. Start the base bridge at
$x=(\cosh r_0-1)/2$, $y=0$, end it at $(0,0)$, and let its duration be
$0<T\le\tau$. For all sufficiently large $r_0$, uniformly in $T$,
\begin{equation}
 \mathbb E^{\rm br}_{r_0,T}\int_0^T(X_s+Y_s)ds
       \ge c_\tau e^{r_0}T/r_0.                       \label{eq:bridge-occupation-lower}
\end{equation}
For every fixed $0<\delta<\tau$, with the late integral set to zero
when $T<\delta$,
\begin{equation}
 \frac{\mathbb E^{\rm br}_{r_0,T}\int_\delta^T(X_s+Y_s)ds}{\mathbb E^{\rm br}_{r_0,T}\int_0^T(X_s+Y_s)ds}
 \le C_\tau\{e^{-\delta r_0/\tau}+r_0e^{-r_0}\}.       \label{eq:bridge-late-occupation}
\end{equation}
Constants may depend on $D_0\tau$, and are independent of $T,r_0$ in
this range. These are expectations under bridges, not pathwise claims.
\end{lemma}

These bounds distinguish early and late occupation along a lifetime bridge. They prevent cancellation of every nonzero spatial Fourier mode.

\begin{proof}
Let $L_B=\frac{1}{2}(\partial_r^2+r^{-1}\partial_r)$ be the planar
Bessel generator. Set
\[
 h(r)=(\sinh r/r)^{1/2},\quad h(0)=1,\qquad
 V(r)=\frac{L_Bh}{h}
     =\frac{1}{8}+\frac{1}{8}\left(\frac{1}{r^2}-\frac{1}{\sinh^2r}\right).
\]
The potential is continuous and bounded on $[0,\infty)$:
$V(0)=1/6$ and $V(r)\to1/8$ at infinity. Direct conjugation gives
\begin{equation}
 h^{-1}(L_B-V)(hf)=\tfrac12(f''+\coth r\,f').          \label{eq:bessel-conjugation}
\end{equation}
The hyperbolic radial bridge from $r_0$ to zero is therefore the planar
Bessel bridge with density proportional to
$\exp\{-D_0\int_0^T V(r_s)ds\}$. The unnormalized transition formula
also contains $h(r_T)/h(r_0)$; this endpoint factor cancels only after
normalizing the bridge. Taking $r_T\downarrow0$ is legitimate since
$h(r_T)\to1$ and the potential is bounded. The resulting normalized
density is bounded above and below by positive constants depending
only on $D_0\tau\operatorname{osc}V$.

The planar Bessel bridge is the radius of a planar Brownian bridge
from $r_0e_1$ to zero, with covariance coefficient $D_0$.
At time $s$ this vector has mean $(1-s/T)r_0e_1$ and covariance
$D_0s(1-s/T)I$. The function
$v\mapsto(\cosh\|v\|-1)/2$ is convex. Jensen therefore gives
\[
 \mathbb E X_s\ge(\cosh((1-s/T)r_0)-1)/2
\]
under this Bessel bridge. Integrating on $0\le s\le T/r_0$ and using
the lower bridge density bound proves \eqref{eq:bridge-occupation-lower}, for $r_0\ge2$ enlarged by
a fixed constant if necessary.
Conversely a centered planar Gaussian with covariance at most
$D_0\tau I$ has a bounded exponential moment of its radius. Hence
\[
 \mathbb E^{\rm br}_{r_0,T}X_s
       \le C_\tau e^{(1-s/T)r_0}.
\]
The Jacobi coordinate always lies in $[0,1]$. For $T\ge\delta$,
\[
 \mathbb E^{\rm br}_{r_0,T}\int_\delta^T(X_s+Y_s)ds
 \le C_\tau e^{r_0}(T/r_0)e^{-\delta r_0/T}+T.
\]
Divide by \eqref{eq:bridge-occupation-lower}, and use $T\le\tau$ to obtain \eqref{eq:bridge-late-occupation}. If $T<\delta$
the asserted upper bound is immediate. No estimate deteriorates as
$T\downarrow0$.
\end{proof}

\begin{proposition}[Nonzero response in every Fourier mode]
\label{prop:t22-nonzero-response}
For every $\tau>0$ and every integer $m\ne0$, there exists a finite
$x>0$ such that $J_{\tau,m}(x,0)>0$, with $J$ as in Theorem~\ref{thm:t22-response-variance}.
\end{proposition}

Every nonzero Fourier frequency has at least one detecting choice of dual variables. Together with the Gaussian variance bound, this gives positivity for every nonconstant smooth real test.

\begin{proof}
Set $\delta=\tau/(6|m|)$. Then
$q_m(s)=\cos(2\pi ms/\tau)\ge1/2$ on $[0,\delta]$, and $q_m\ge-1$
everywhere. Thus for each lifetime $0<T\le\tau$, its conditional
expected weighted occupation is at least
\[
 \frac{1}{2}\mathbb E^{\rm br}_{r_0,T}\int_0^T R_sds
       -\frac{3}{2}\mathbb E^{\rm br}_{r_0,T}\int_\delta^T R_sds.
\]
Choose a sufficiently large finite $r_0$ in \eqref{eq:bridge-late-occupation} to make the ratio
strictly smaller than $1/3$, uniformly in $T$. This lower bound is then
strictly positive for every $T\in(0,\tau]$.
Lemma~\ref{lem:t22-lifetime-bridge} mixes these conditional expectations
against a strictly positive lifetime density. Their mixture is exactly
$J_{\tau,m}((\cosh r_0-1)/2,0)$, since the transformed process stays
zero after its lifetime. Absolute integrability is paid by the true
martingale identity $\mathbb E R_s=R_0$:
$\mathbb E\int_0^\tau R_sds=\tau R_0<\infty$.
This proves the strict inequality.
\end{proof}

\begin{proof}[Proof of the variance assertions in Theorem~\ref{thm:endpoints}]
Theorem~\ref{thm:joint-recovery} and the fixed count moments give
convergence of the bounded quadratic mark statistic on each compact
spectral window. Its limiting weighted intensity is translation
invariant in the spectral coordinate by Proposition~\ref{prop:stationarity}.
It is therefore a constant times Lebesgue measure, so a weight $w$
of integral one yields the same variance for every such $w$.

For each $m\ne0$, choose a finite $x$ from
Proposition~\ref{prop:t22-nonzero-response}. Equation
\eqref{eq:response-variance-bound} gives positive variance in that mode.
Spatial rotation invariance diagonalizes the Fourier covariance:
\[
 V_{\tau,b,E}(\phi)=\sum_{m\ne0}|\widehat\phi_m|^2v_m,
 \qquad
 v_m=\E\sum_{(u,\nu)\in\Xi_{\tau,b,E}}w(u)|\nu(e^{2\pi im\cdot})|^2>0.
\]
Smoothness gives absolute summability and justifies expansion under
expectation. Every nonconstant real smooth $\phi$ has a nonzero
coefficient at some $m\ne0$, proving strict positivity. Finally
Theorem~\ref{thm:small-critical} gives
$V_{\tau,b,E}(\phi)\le C_\epsilon\|\phi-\bar\phi\|_\infty^2
\tau^{1-\epsilon}$ for $0<\tau\le1$. This proves the vanishing
small-$\tau$ limit with the required order of limits.
\end{proof}
\section{Local mass flatness and the sine endpoint}
\label{sec:sine-endpoint}
This section proves the sine endpoint of Theorem~\ref{thm:endpoints}. The
key input is that, when $\tau$ is small, the mass of an eigenvector on a
single pair of blocks is nearly constant in the microscopic window.
Theorem~\ref{thm:t22-tagged-convex-order} states that a tagged response is
dominated, in the convex order, by geometric Brownian motion.
Proposition~\ref{prop:t26-microscopic-tag-comparison} identifies the mass on
a block with the ratio of a tagged response to the average response of the
open path. Theorem~\ref{thm:t22-local-mass-flat} combines them to show that
the scaled mass on a pair differs from one by $O(\sqrt\tau)$ in $L^1$.
Theorem~\ref{thm:t22-sine-endpoint} then states that the local spectrum is
within $C\tau^{1/8}$ of the sine process.
In Subsection~\ref{subsec:s28-a-closed-inequality-for-a-positive-tagge} we prove the convex-order comparison for a tagged
source. In Subsection~\ref{subsec:s28-local-mass-flatness-under-the-actual-mic} we transfer it to the mass of a block in the
actual microscopic window, using the lower tail of the full clock, and prove
local mass flatness. In Subsection~\ref{subsec:s28-the-scalar-pair-flow} we follow the eigenvalues along a
scalar flow on one pair and compare them with the crossing values, which have
sine statistics (Lemmas~\ref{lem:t22-crossing-sine}
and~\ref{lem:t22-flow-admission}).

Fix $T_0<\infty$; constants below may depend on $b,E,T_0$.

\begin{theorem}[Original-ring sine endpoint, for every critical subsequence]
\label{thm:t22-sine-endpoint}
For every nonnegative smooth compactly supported $f$ on $\mathbb R$,
\begin{equation}\label{eq:t22-sine-endpoint}
 \limsup_{W\to\infty,\ n/W\to\tau}
 \left|\mathbb E e^{-\sum_\alpha f(u_\alpha)}
              -\mathbb E_{\rm sine}e^{-\sum_{u\in\mathcal S}f(u)}\right|
       \le C_f\tau^{1/8},\qquad 0<\tau\le\min(1,T_0).
\end{equation}
Consequently $\mathcal X_{\tau,b,E}$ converges to the unit-intensity
sine process as $\tau\downarrow0$, after the critical limit.
\end{theorem}

The estimate concerns the original ring and gives a rate in the critical
parameter after taking the width limit. The proof first controls the
physical pair mass, then compares exact crossing locations with energies.

\subsection{A closed inequality for a positive tagged source}\label{subsec:s28-a-closed-inequality-for-a-positive-tagge}
 We bound the fluctuations created after a source is inserted at one
block. Theorem~\ref{thm:t22-tagged-convex-order} compares convex
polynomials of its response with geometric Brownian motion and gives
small-time bounds after averaging over tags. We derive the drift and
covariance, use positivity to bound the noise, and integrate the
resulting comparison along the path.

 Let $Q_{n;k}=Y_{n;k}/W$ and, for an independent standard full Gaussian
$h$, put $Z_{n;k}=h^*Q_{n;k}h$. These are local unnormalized intensities.
A tag is a source injected at a single
block $k$, at time $s=k/W$, and $Z_{n;k}$ is its response measured at the
end of the path, at time $\tau=n/W$. The theorem says that $Z_{n;k}$ is
dominated, in the convex order, by the value at time $\tau-s$ of a
geometric Brownian motion started at $d$. So it has mean $d$, and its
fluctuations around $d$ are at most those of geometric Brownian motion
over the remaining time. The second part averages this over all tags, and
shows that the centered moments are of order $\tau^p$. This is what makes
the tagged response flat when $\tau$ is small.

\begin{theorem}[Convex polynomial comparison for the tagged response]
\label{thm:t22-tagged-convex-order}
If $k/W\to s>0$ and $n/W\to\tau\ge s$, every limiting law $Z$ of
$Z_{n;k}$ satisfies, for every convex polynomial $\varphi$ on
$[0,\infty)$,
\[
 \mathbb E\varphi(Z)\le
 \mathbb E\varphi\left(d e^{\sqrt{D_0}B_{\tau-s}
                                      -D_0(\tau-s)/2}\right).
\]
Its mean equals $d$. All fixed moments converge along such a limiting
subsequence. In particular, for every positive integer $p$,
\begin{equation}\label{eq:t22-tagged-flat}
 \limsup_W\frac{1}n\sum_{k=1}^n
             \mathbb E|Z_{n;k}-d|^{2p}
                 \le C_{p,T_0}\tau^p.
\end{equation}
The same bound holds for
$M_W=n^{-1}\sum_kZ_{n;k}=T_W/\tau_W$ in place of the averaged tag. All assertions hold along deterministic $E_W\to E$, with $D_0=D_0(E)$.
\end{theorem}

The comparison gives upper bounds for every fixed even centered moment.
It will be combined with the lower tail of the full clock to normalize
local intensities.

\begin{proof} After the source is injected, the tagged response
evolves with no further input. We compute the one-step change of a polynomial
function of the response, and show that, up to errors that vanish after
summation (Steps~1--3), its drift is at most that of a geometric Brownian
motion, the positive process whose logarithm is a Brownian motion with
drift. For convex test functions, this comparison of drifts gives the
comparison in convex order, by a backward argument from the end of the path
(Step~4); Figure~\ref{fig:tco-gbm} illustrates it. Step~5 identifies the value
of the response at the injection, and Step~6 averages the bound over the
location of the tag.

\emph{Step 1: the errors vanish after summation.} After its single
injection, the tagged clock has no additive source.
Equation~\eqref{eq:full-tag-clock-moments} gives all its fixed trace
moments. For $k\ge\epsilon W$ the source-gap argument applies to this
homogeneous clock. Lemma~\ref{lem:t26-contracting-words} gives,
uniformly in the tag,
\[
 \|Z_j^{\rm word}\|_p\le C_{p,\epsilon}(\rho^{j-k}+W^{-1/2})
\]
for every contracting word created by the Wick expansion below. There
is no need to discard an initial mixing segment, because its geometric
contribution to the sum is $O(W^{-1})$.

\emph{Step 2: the one-step change of a polynomial test.} For a polynomial
$\varphi$ we put
\[
 \Psi_\varphi(Q)=\mathbb E_h\varphi(h^*Qh),
\]
where $h$ is a standard Gaussian vector. After one step the clock becomes
$Q_+=N^*QN$. The Gaussian law of $h$, and the equality of the nonzero
spectra of the two Gram matrices $N^*QN$ and $Q^{1/2}NN^*Q^{1/2}$, let us
write its transported value with the same $h$ as
\[
 \varphi\bigl(h^*Q^{1/2}(I+\mathcal E)Q^{1/2}h\bigr),
 \qquad NN^*=I+\mathcal E.
\]
We put $v=Q^{1/2}h$ and $X=\|v\|^2$. A Taylor expansion in $\mathcal E$
and the covariance \eqref{eq:full-source-covariance} give, with summed
error $o(1)$,
\begin{equation}
 \mathbb E\Psi_\varphi(Q_+)-\mathbb E\Psi_\varphi(Q)
 =\frac{D_0}{2W}\mathbb E\mathbb E_h
       [\varphi''(X)(v^*Pv)(v^*P'v)]+r_j.
                                                               \label{eq:tagged-convex-generator}
\end{equation}
The right side is the opposite-sign term of the covariance, paired with
the second derivative of $\varphi$.

\emph{Step 3: the remainder $r_j$.} We now check that every other term is
summable. The Gaussian Wick contraction of the Ward term multiplied by
$\varphi'(X)$ contains traces $\operatorname{Tr}(Q^rPS)$. These are bounded
by the trace and Hilbert--Schmidt errors of the source times fixed powers
of $\operatorname{Tr}Q$. The equal-sign terms with $\varphi''(X)$ have the
same bound, with two such insertions. Their fixed higher moments
therefore make them summable. The four-vector remainder in
\eqref{eq:tagged-convex-generator} is at most $C_{p,\epsilon}W^{-2}$ after
Gaussian averaging. The centered terms of third and higher order cost
$C_{p,\epsilon}W^{-3/2}$.

More explicitly, for a polynomial of fixed degree, the total error over
$k\le j\le n$ is bounded by
\[
 \sum_{j=k}^{n-1}\E|r_j|
 \le C_{p,\epsilon,T_0}\sum_{j=k}^{n-1}
 \left[W^{-1}\rho^{j-k}+W^{-3/2}+W^{-2}\right]
 \le C_{p,\epsilon,T_0}(W^{-1/2}+W^{-1}).
\]
Here the geometric term sums to $O(W^{-1})$, and the other two terms, over
at most $T_0W$ steps, give $O(W^{-1/2})$ and $O(W^{-1})$. We first use this
estimate before a trace stop, and then remove the stop by
Lemma~\ref{lem:t26-full-source-stops}. The higher tagged moments pay for
every polynomial factor. In the Wick expansion the same-sign terms are
products of traces $\Tr(Q^rP\mathbb S)$, or cycles such as
$\Tr(Q^rP\mathbb SQ^sP\mathbb S)$, with $r,s\ge1$. So the required
separation by trace-class letters is present even when $\varphi''$ is
constant. These claims use the tagged moments, not an approximation of
$N$ in operator norm.

\emph{Step 4: comparison with geometric Brownian motion.} Let $\varphi$ be
convex, so $\varphi''\ge0$. Since $P$ and $P'$ are projections,
\[
 0\le(v^*Pv)(v^*P'v)\le X^2 .
\]
So, in time $t=j/W$, the drift in \eqref{eq:tagged-convex-generator} is at
most that of the generator $(D_0/2)x^2\partial_x^2$, which is the generator
of geometric Brownian motion. To make the comparison rigorous, we use
the backward polynomial
\[
 u(t,x)=\mathbb E\varphi\bigl(xe^{\sqrt{D_0}B_t-D_0t/2}\bigr)
\]
in \eqref{eq:tagged-convex-generator}. It solves the backward equation
$\partial_tu=(D_0/2)x^2\partial_x^2u$, and it stays convex in $x$, since it
is an average of convex functions of positive multiples of $x$. So the
inequality of the generators applies to it at every step. Summing the
resulting backward inequality from the tag to the end of the path proves
the convex polynomial bound. For linear $\varphi$ the second derivative
vanishes, so the drift is zero up to the summed error, and the mean is
preserved.

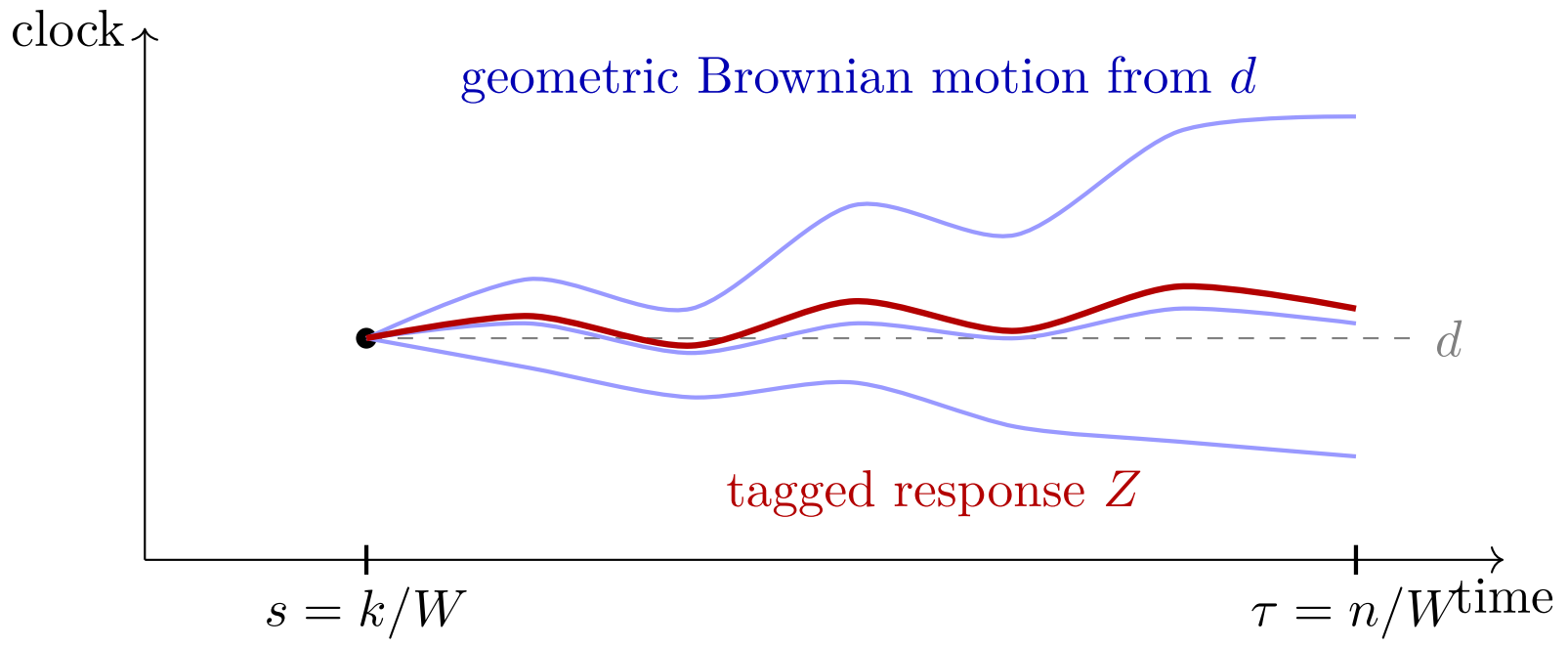
\begin{figure}[H]
\centering
\begin{tikzpicture}[x=1cm,y=1cm]
 \draw[->] (0,0)--(9.2,0) node[below,font=\small] {time};
 \draw[->] (0,0)--(0,3.6) node[left,font=\small] {clock};
 \draw[thick] (1.5,-.1)--(1.5,.1) node[below=4pt,font=\small] {$s=k/W$};
 \draw[thick] (8.2,-.1)--(8.2,.1) node[below=4pt,font=\small] {$\tau=n/W$};
 \draw[dashed,black!50] (1.5,1.5)--(8.6,1.5) node[right,font=\small] {$d$};
 \fill (1.5,1.5) circle (2pt);
 \draw[thick,blue!40] plot[smooth] coordinates
   {(1.5,1.5) (2.6,1.9) (3.7,1.7) (4.8,2.4) (5.9,2.2) (7.0,2.9) (8.2,3.0)};
 \draw[thick,blue!40] plot[smooth] coordinates
   {(1.5,1.5) (2.6,1.3) (3.7,1.1) (4.8,1.2) (5.9,.9) (7.0,.8) (8.2,.7)};
 \draw[thick,blue!40] plot[smooth] coordinates
   {(1.5,1.5) (2.6,1.6) (3.7,1.4) (4.8,1.6) (5.9,1.5) (7.0,1.7) (8.2,1.6)};
 \draw[very thick,red!70!black] plot[smooth] coordinates
   {(1.5,1.5) (2.6,1.65) (3.7,1.45) (4.8,1.75) (5.9,1.55) (7.0,1.85) (8.2,1.7)};
 \node[font=\small,blue!70!black,anchor=west] at (2.0,3.25)
   {geometric Brownian motion from $d$};
 \node[font=\small,red!70!black,anchor=west] at (3.8,.45)
   {tagged response $Z$};
\end{tikzpicture}
\caption{The comparison in Step~4 of the proof of
Theorem~\ref{thm:t22-tagged-convex-order}, drawn schematically. The tag is
injected at time $s=k/W$ with response $d$. From then on, the one-step
drift of a convex test is at most that of the geometric Brownian motion
$de^{\sqrt{D_0}B_t-D_0t/2}$ (blue sample paths), because
$(v^*Pv)(v^*P'v)\le X^2$. So the tagged response (red) is dominated in the
convex order by the geometric Brownian motion at time $\tau-s$. Both have
mean $d$, and the centered moments of the response are at most those of
the geometric Brownian motion, which are of order $(\tau-s)^p$.}
\label{fig:tco-gbm}
\end{figure}

\emph{Step 5: the initial value.} At the injection, the source is
$Q=C_+/W$. Its trace tends to $d$ in every fixed moment, and
$\mathbb E\operatorname{Tr}Q^2=O(W^{-1})$, uniformly in
$\epsilon W\le k\le T_0W$; this is
Lemma~\ref{lem:source-recorded-estimates}(f), with $d=d_*$. Indeed, the near part of $C_+$ has
normalized trace $d+O(\delta+g/W)$, by the injection ladder. Its far part
has trace $O_{L^p}(g)$, from $\|\mathbb R_{12}\|_{\rm HS}^2=g/(4\eta^2)$ and
the weighted fresh bounds. Gaussian trace concentration supplies the
fluctuations, and the operator moments of $C_+$ bound its squared trace
norm as stated. Consequently its Gaussian response $h^*Qh$ tends to the
constant $d$, since its variance is controlled by
$\mathbb E\operatorname{Tr}Q^2\to0$. This supplies the initial value of the
comparison, including the case $\tau=s$. All higher moments are uniformly
integrable, by \eqref{eq:full-tag-clock-moments}, so all fixed moments
converge along a limiting subsequence.

\emph{Step 6: the averaged bound.} For $\varphi(x)=(x-d)^{2p}$, the geometric
Brownian bound is
\[
 \E\bigl(de^{\sqrt{D_0}B_{\tau-s}-D_0(\tau-s)/2}-d\bigr)^{2p}
 \le C_{p,T_0}(\tau-s)^p.
\]
This bound is uniform for $s\ge\epsilon$. Otherwise, a contradicting
sequence would have a convergent tag time in that compact interval, and
the argument just given would apply to it. The tags with $k<\epsilon W$
are a fraction $\epsilon W/n$ of all tags, so by
\eqref{eq:full-tag-clock-moments} they contribute at most
$C_{p,T_0}\epsilon/\tau$. We let $W\to\infty$, and then $\epsilon\downarrow0$.
Since $\tau-s\le\tau$, this proves \eqref{eq:t22-tagged-flat}. Finally,
Jensen's inequality for the spatial average gives
$|M_W-d|^{2p}\le n^{-1}\sum_k|Z_{n;k}-d|^{2p}$, which proves the assertion
for $M_W$.

\emph{Moving energies.} Steps~1--6 use finite-width estimates only
through Lemmas~\ref{lem:full-source-covariance},
\ref{lem:t26-contracting-words} and~\ref{lem:t26-full-source-stops},
\eqref{eq:full-tag-clock-moments} and the tag-birth value of Step~5,
which hold along $E_W\to E$ by
Lemma~\ref{lem:source-recorded-estimates}(b),(d)--(f). The summed errors
of Steps~1--3 vanish as before. In Step~4 we use the backward polynomial
built with the drift coefficient $D_0(E_W)$; for polynomial $\varphi$ its
coefficients are continuous in $D_0$, and $D_0(E_W)\to D_0(E)$. Hence
every limiting law along $E_W$ satisfies the comparison with
$D_0=D_0(E)$, and Steps~5 and~6 are unchanged.
\end{proof}

\subsection{Local mass flatness under the actual microscopic window law}\label{subsec:s28-local-mass-flatness-under-the-actual-mic}
Let $w\ge0$ be smooth and compactly supported, with $\int w=1$.
  We transfer the tagged path estimates to the mass of a block in the
actual microscopic spectral window.
Proposition~\ref{prop:t26-microscopic-tag-comparison} identifies the
normalized path response, and Theorem~\ref{thm:t22-local-mass-flat}
shows that this mass becomes flat as the critical ratio decreases.
We remove the collar, compare the two positive clocks, and combine
tagged fluctuations with the lower tail of the total clock.

\begin{proposition}[Physical microscopic tag comparison]
\label{prop:t26-microscopic-tag-comparison}
Fix $n/W\to\tau\in(0,T_0]$, a fixed compact-bulk energy $E$, a smooth
compactly supported $w$, and $\Phi\in C_b([0,\infty))$. Let $r=n-3$, let
$Y^{\rm int}=\sum_{j=1}^rY_{r;j}$ be the interior clock of the actual
independent $r$-block path with two matched auxiliary leads, and let $K$
be independent and uniform on $\{1,\ldots,r\}$. For an independent
standard circular $h\in\C^{2W}$, put
\[
 Z_{r;K}=W^{-1}h^*Y_{r;K}h,\qquad
 M_r=(rW)^{-1}h^*Y^{\rm int}h.
\]
For every physical block $k$ of the original ring,

\begin{equation}\label{eq:t26-microscopic-tag-comparison}
 \E\sum_\alpha w(u_\alpha)\Phi(np_{\alpha k})
 -\int w(u)\E_{{\rm path},h,K}
       \Phi(Z_{r;K}/M_r)(e(u))\,du\longrightarrow0,
 \quad e(u)=E+\frac{u}{N\rho(E)}.
\end{equation}

The convergence is uniform in the physical tag, and the path expectation
has its original prior. The critical limit is taken before any trace or
small-denominator cutoff is removed.
\end{proposition}

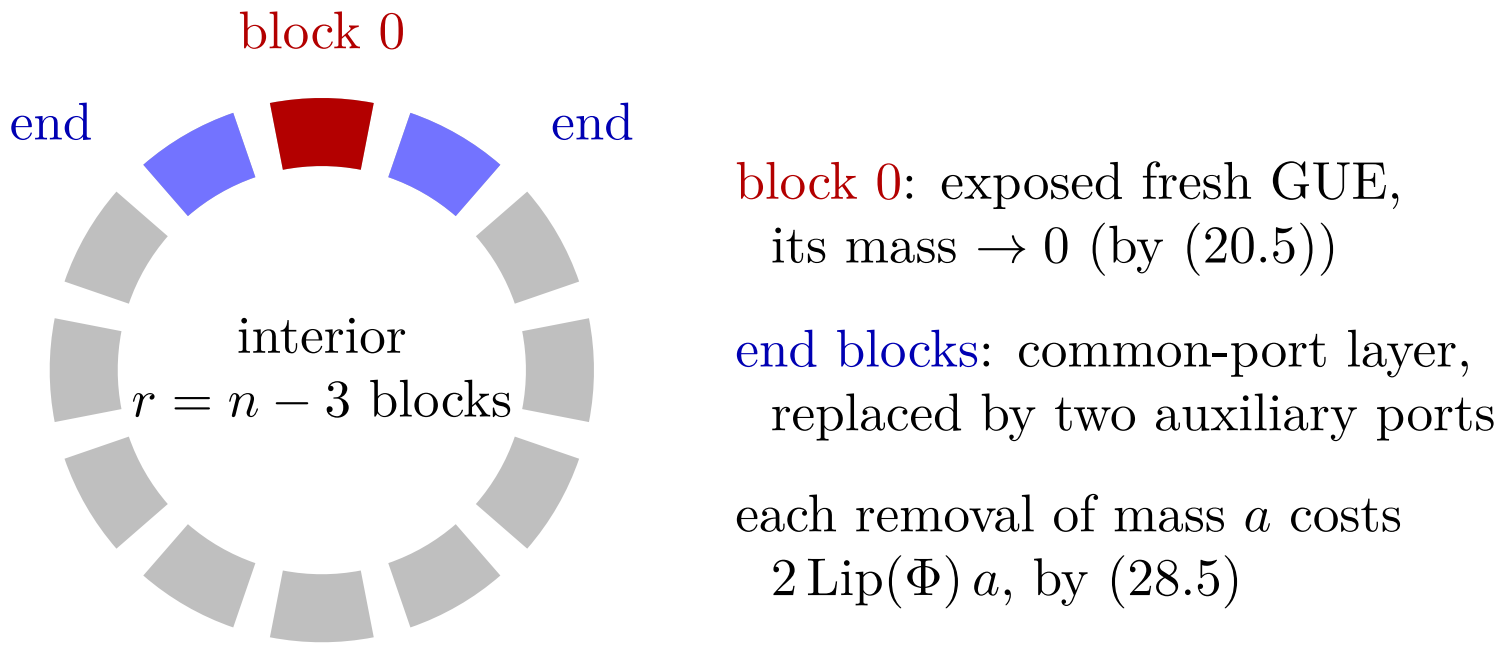
\begin{figure}[ht]
\centering
\begin{tikzpicture}[x=1cm,y=1cm]
 \foreach \a in {0,30,...,330}{
   \fill[black!25] (\a-11:1.8) arc (\a-11:\a+11:1.8) -- (\a+11:1.35) arc (\a+11:\a-11:1.35) -- cycle;
 }
 \fill[red!70!black] (90-11:1.8) arc (90-11:90+11:1.8) -- (90+11:1.35) arc (90+11:90-11:1.35) -- cycle;
 \fill[blue!55] (60-11:1.8) arc (60-11:60+11:1.8) -- (60+11:1.35) arc (60+11:60-11:1.35) -- cycle;
 \fill[blue!55] (120-11:1.8) arc (120-11:120+11:1.8) -- (120+11:1.35) arc (120+11:120-11:1.35) -- cycle;
 \node[font=\small,red!70!black] at (90:2.25) {block $0$};
 \node[font=\small,blue!70!black,anchor=west] at (50:2.15) {end};
 \node[font=\small,blue!70!black,anchor=east] at (130:2.15) {end};
 \node[font=\small,align=center] at (0,0) {interior\\ $r=n-3$ blocks};
 \node[font=\small,anchor=west,align=left] at (2.6,1.0)
   {{\color{red!70!black}block $0$}: exposed fresh GUE,\\
    \ \ its mass $\to0$ (by \eqref{eq:vanishing-block-atom})};
 \node[font=\small,anchor=west,align=left] at (2.6,-.1)
   {{\color{blue!70!black}end blocks}: common-port layer,\\
    \ \ replaced by two auxiliary ports};
 \node[font=\small,anchor=west,align=left] at (2.6,-1.2)
   {each removal of mass $a$ costs\\
    \ \ $2\operatorname{Lip}(\Phi)\,a$, by \eqref{eq:t26-tag-deletion}};
\end{tikzpicture}
\caption{The positions in the proof of
Proposition~\ref{prop:t26-microscopic-tag-comparison}. The ring has $n$
blocks. Block $0$ (red) is exposed as a fresh GUE matrix, and the coarea
formula is applied there. The two blocks next to it (blue), with their
bonds to block $0$, form the common-port layer, and in the final
comparison they are replaced by two matched auxiliary leads. What remains
is a path of $r=n-3$ interior blocks. A uniformly chosen tag lands among
them with probability $r/n\to1$. Removing block $0$ and the two ports
changes the averaged test only by a multiple of the mass they carry, by
the deletion identity \eqref{eq:t26-tag-deletion}, and that mass tends to
zero.}
\label{fig:mt-positions}
\end{figure}

\begin{proof} We want to compare the scaled mass of an eigenvector
on one block with a ratio of two clocks of the open path. Because the ring is
invariant under rotations, it suffices to treat a test averaged over all
blocks (Step~1). We apply the exact formula on an exposed block and replace
the eigenvector by a Gaussian reconstruction, using the two-clock comparison
(Steps~2--3). We then remove the exposed block and the two auxiliary ports,
which carry little mass; a simple identity shows that deleting positions of
small total mass changes the averaged test only by a multiple of that mass
(Steps~4--5). Figure~\ref{fig:mt-positions} shows the positions involved.
Steps~6--8 pass from a uniform tag to an interior tag and to general bounded
tests. The lead-profile inputs include the Gaussian-field identification
in the proof of Theorem~\ref{thm:t22-physical-lead-profile},
\eqref{eq:vanishing-block-atom}, and cyclic invariance of the original
ring law. These are proof-internal inputs, not consequences of the
displayed profile formula alone. The separate covariance-closeness
input is not used below.

\emph{Step 1: reduction to an averaged, rotation-invariant test.} It
suffices first to take a nonnegative $w$ of integral one and a bounded
Lipschitz $\Phi$. By cyclic invariance, the left side of
\eqref{eq:t26-microscopic-tag-comparison} does not depend on $k$. So it
equals its average over all $n$ tags,
\[
 \E\sum_\alpha w(u_\alpha)\,\frac1n\sum_{j=1}^n\Phi(np_{\alpha j}),
\]
and this averaged functional is invariant under cyclic rotation. Its
expectation in the Palm law of the mass on block $0$, with density
$np_{\alpha0}$, is exactly its unweighted spectral expectation, because
$\sum_jp_{\alpha j}=1$. So the one-block coarea formula used in
Theorem~\ref{thm:t22-physical-lead-profile} applies to this observable,
which depends on $W$. We do not need continuity in the weak topology of
profiles.

\emph{Step 2: the clock ratio on the exterior.} We condition on the
complete exterior of the fresh $W$-dimensional GUE block. Let
$\mathcal T$ be its physical Schur lift, and put
\[
 D_j=\mathcal T^*P_j\mathcal T,\qquad
 D=\sum_jD_j=\mathcal T^*\mathcal T.
\]
Then the scaled mass on block $j$ is $nv^*D_jv/v^*Dv$, a ratio of two
positive clocks. Corollary~\ref{cor:t26-uniform-clock-tags} compares the
averaged test $\Phi(nv^*D_jv/v^*Dv)$ with its Gaussian counterpart. The
exceptional tag bound $1/K_0$ holds for each exterior, including its
covariance $\Gamma$, and no comparison of a selected small sector mass is
involved. The good event of the deformation is exactly the one in the
proof of the one-block coarea formula. Its complement has vanishing
physical mass, and the full unranked normalizer satisfies
$Z/W\to\rho(e)$ in prior $L^1$.

The energy Jacobian and the coarea factor combine to

\[
 \frac{nW}{N\rho(E)}=\frac1{\rho(E)},
\]

since   $N=nW$. Hence this step leaves exactly $w(u)\,du$ times the actual
exterior prior in the limit. Total variation allows these bounded
observables to depend on $W$, on the tag, and on all the exterior clocks.

\emph{Step 3: the Gaussian field and the exposed block.} The relative
change of covariance in that proof is legitimate for these sources, by
the last paragraph of Lemma~\ref{lem:t26-two-positive-clocks}. The
resulting field is exactly $(ULg,Hg)$, up to the common factor $\sqrt c$,
with $U=JR$. All incident bonds and fresh diagonals keep their original
independent laws.

The mass of this field on block $0$ tends to zero in expectation. Indeed,
the exact Palm expectation of the physical mass is
\[
 n\E\sum_\alpha w(u_\alpha)p_{\alpha0}^2
 =\E\sum_\alpha w(u_\alpha)\sum_jp_{\alpha j}^2=o(1),
\]
by cyclic invariance and \eqref{eq:vanishing-block-atom}. The dominated
ratio comparison already proved transfers this bounded mass to the
Gaussian field.

\emph{Step 4: removing positions costs only their mass.} We now give the
quantitative argument for removing positions at the microscopic scale.
Let $p$ be any probability vector on $n$ sites. We remove coordinates of
total mass $a<1$ and renormalize the rest to $\widehat p$, assigning zero
to the removed sites. Then
\begin{equation}\label{eq:t26-tag-deletion}
 \frac{1}n\sum_{j=1}^n|np_j-n\widehat p_j|
       =\sum_j|p_j-\widehat p_j|=2a.
\end{equation}
Indeed, the removed sites contribute $a$, and on the remaining sites
$\widehat p_j-p_j=p_ja/(1-a)$, which also sums to $a$. Consequently,
removing block $0$ changes the averaged test by at most
$2\operatorname{Lip}(\Phi)a$, despite the factor $n$ inside $\Phi$.

The remaining positive clocks are exactly $L^*Y_jL$, where $Y_j=U^*P_jU$
on the $n-1$ old and auxiliary positions, and $\sum_jY_j=Y^{\rm all}=U^*U$.
We apply the common-port part of
Corollary~\ref{cor:t26-uniform-clock-tags}, with a zero clock at block
$0$, conditionally on the entire old path. The fixed scalar regularity
error has $\E\delta=O(W^{-1})$, and its exceptional probability is smaller
than any power of $W$. So the uniform conditional comparison integrates
under the prior, and it replaces $Lg$ by $\sqrt\gamma h$. The scalar
$\gamma$ cancels from every ratio.

\emph{Step 5: removing the auxiliary ports.} We remove the auxiliary
positions only after this comparison. Their clock is $Y^{\rm port}=R^*R$,
with $\Tr Y^{\rm port}\le2W/(\Im z_b)^2$. Thus
\[
 \E h^*Y^{\rm port}h/W^2\le C/W.
\]
The interior clock is not small. Every subsequential limit of
$h^*Y^{\rm int}h/W^2$ is strictly positive, by
Theorem~\ref{thm:t22-negative-clock}, since $r/W\to\tau$. Let $a_{\rm port}$
be the normalized mass on the ports. If $a_{\rm port}>v$ while the
interior clock exceeds $t$, then the port clock exceeds $vt$. So, for any
$v,t>0$, Markov's inequality gives
\[
 \Prob(a_{\rm port}>v)
 \le\Prob\{h^*Y^{\rm int}h/W^2\le t\}
             +\frac{C}{Wvt}.
\]
We take the width limit, and then $t\downarrow0$. Since $a_{\rm port}$ is
bounded, this gives $\E a_{\rm port}\to0$. Equation~\eqref{eq:t26-tag-deletion}
then pays for removing both ports and renormalizing.

\emph{Step 6: from a uniform tag to an interior tag.} A uniform tag now
lands on one of the $r$ interior positions with probability
$r/n\to1$. The other three positions cost $O(\|\Phi\|_\infty/n)$. On an
interior position, the argument of $\Phi$ is $(n/r)Z_{r;K}/M_r$. For every
sample,
\[
 \frac1r\sum_KZ_{r;K}/M_r=1,
\]
because $\sum_KZ_{r;K}=W^{-1}h^*Y^{\rm int}h=rM_r$. So changing the factor
$n/r$ to one costs at most $\operatorname{Lip}(\Phi)|n/r-1|$. This proves
\eqref{eq:t26-microscopic-tag-comparison} for Lipschitz tests.

\emph{Step 7: general tests.} For a general bounded continuous $\Phi$, we
approximate it by a Lipschitz function on $[0,R_0]$. Every normalized mass
vector in this proof has average scaled mass one. So, by Markov's
inequality, the proportion of tags beyond $R_0$ is at most $1/R_0$,
pointwise, for every exterior and every Gaussian sample. This removes the
tails uniformly. We take the width limit, and then let the approximation
error and $R_0^{-1}$ tend to zero. A signed $w$ follows by domination. All
the source estimates are uniform for $e=e(u)$ on its compact window, so
the actual energy marginal remains.

\emph{Step 8: the fixed-tag statement with cutoffs.} For clarity, we note
that the stronger statement for a fixed tag, with cutoffs, also follows
from the joint lemma. On the event
\[
 A_W=\frac{\Tr Y^{\rm int}}{rW}\in[\epsilon,K_0],\qquad
 \frac{\Tr Y_{r;j}}W\le K_0,
\]
the trace ratio of $(Y^{\rm int},rY_{r;j})$ is at most $K_0/\epsilon$,
uniformly in $j$. Equations
\eqref{eq:full-clock-moments}--\eqref{eq:full-tag-clock-moments} pay for
the upper cutoffs, at cost $C_pK_0^{-p}$, uniformly in the tag.
Lemma~\ref{lem:t26-clock-small-ball} pays for the lower cutoff, at cost
$C_p\epsilon^{p/2}/(1-\sqrt\epsilon)$ in the critical limsup. On the cutoff
event, the common-port estimate is \eqref{eq:t26-two-clock-laplace}, and
the corresponding Palm estimate is the fixed-moment estimate with trace
ratio $K_0/\epsilon$. The denominator cutoff in the joint lemma is removed
after the width limit, followed by $K_0\to\infty$ and $\epsilon\downarrow0$.
The uniform-tag argument above additionally controls the conditional
trace ratios of the fresh exterior before this cutoff on the old path is
imposed. It is therefore a complete comparison from Palm to prior, with
all the stated conditioning, rather than an inference from weak
convergence of profiles.
\end{proof}

 An eigenvector spread evenly over the
ring would put mass $1/n$ on each block, so its scaled mass $np_{\alpha k}$
would equal one. The next theorem says that, for the eigenvalues near $E$, the
scaled mass on a single block, or on a single pair, differs from one by
$O(\sqrt\tau)$ in $L^1$. So when $\tau$ is small, the eigenvectors are
locally flat, even at the scale of one block. This is the input that the
sine endpoint needs.

\begin{theorem}[Strong local mass flatness]
\label{thm:t22-local-mass-flat}
For the original periodic matrix, every fixed adjacent pair $S$ and
$0<\tau\le T_0$ satisfy
\begin{equation}\label{eq:t22-local-mass-flat}
 \limsup_{W\to\infty,\ n/W\to\tau}
 \mathbb E\sum_\alpha w(u_\alpha)
       \left|\frac{n}{2}\|P_S\psi_\alpha\|^2-1\right|
       \le C_{b,E,T_0}\sqrt\tau.
\end{equation}
The same statement holds with $n\|P_k\psi_\alpha\|^2$ for one block.
\end{theorem}

\begin{proof} By the previous proposition, the scaled block mass
behaves like the ratio of a tagged response to the average response of the
open path. When the critical ratio is small, both have the same mean and small
fluctuations, by the convex-order comparison, and the denominator is not
small, by the lower tail of the clock; so the ratio is close to one
(Steps~1--2). Step~3 transfers this to the physical block mass. The deviation
from one is measured by an unbounded function, so Step~4 splits it, using the
identity of Figure~\ref{fig:lm-identity}, into a linear part, known exactly,
and a bounded part, to which the comparison applies. Step~5 passes from one
block to a pair.

\emph{Step 1: the ratio $Z/M$ is close to one.} We choose a uniform
independent tag $K\in\{1,\ldots,n\}$ in the actual Gaussian path, and we
take a joint subsequential limit $(Z,M)$ of $(Z_{n;K},M_W)$.
Theorem~\ref{thm:t22-tagged-convex-order} supplies all positive moments,
with
\[
 \mathbb E|Z-d|^{2p}+\mathbb E|M-d|^{2p}\le C_p\tau^p.
\]
Theorem~\ref{thm:t22-negative-clock}, with the rescaling in
Lemma~\ref{lem:t26-clock-small-ball}, gives $M>0$ and
$\E M^{-p}\le C_{p,T_0}$. In particular, since
$\E|Z-M|^2\le2(\E|Z-d|^2+\E|M-d|^2)\le C\tau$, the Cauchy--Schwarz
inequality gives
\[
 \mathbb E|Z/M-1|
 \le(\mathbb E|Z-M|^2)^{1/2}(\mathbb E M^{-2})^{1/2}
 \le C\sqrt\tau.
\]

\emph{Step 2: the limiting ratio has mean exactly one.} For every bounded
continuous $h_0$, we have the exact identity
\[
 \mathbb E[Z_{n;K}h_0(M_W)]=\mathbb E[M_Wh_0(M_W)],
\]
because the average of $Z_{n;K}$ over the uniform tag $K$ is $M_W$. The
positive moments are uniformly integrable, so this identity passes to the
limit. It says that $\mathbb E(Z\mid M)=M$, and then Tonelli's theorem
gives $\mathbb E(Z/M)=1$. This argument does not assume any inverse
moments before the limit.

\emph{Step 3: transfer to the physical block mass.} We apply
Proposition~\ref{prop:t26-microscopic-tag-comparison} with $r=n-3$. Its
interior path has $r/W\to\tau$, so the argument of Steps~1 and~2 applies
to it, with $r$ in place of $n$. The proposition compares every bounded
continuous function of the physical scaled block mass with the same
function of $Z/M$, and it retains the actual energy marginal. The
positive and negative limiting moment bounds are uniform for the
energies $e(u)$ in the compact window.  Indeed $e(u)\to E$ uniformly for $u$ in the window, and these are
the moving-energy forms of Theorem~\ref{thm:t22-tagged-convex-order}
and of Section~\ref{sec:clock-lower-tail}. So every joint subsequential limit,
together with $u$, satisfies the same bound $\E|Z/M-1|\le C\sqrt\tau$.

\begin{figure}[ht]
\centering
\begin{tikzpicture}[x=1.3cm,y=1.1cm]
 \draw[->] (-.1,0)--(3.4,0) node[below,font=\small] {$x$};
 \draw[->] (0,-.1)--(0,3.1);
 \draw[thick,black] (0,1)--(1,0)--(3.1,2.1);
 \draw[thick,dashed,blue!70!black] (0,1)--(2.1,3.1);
 \draw[thick,red!70!black] (0,0)--(1,2)--(3.1,2);
 \draw[dotted] (1,0)--(1,2);
 \node[font=\small,below] at (1,0) {$1$};
 \node[font=\small,anchor=south west] at (3.1,2.15) {$|x-1|$};
 \node[font=\small,blue!70!black,anchor=south east] at (2.05,3.05) {$x+1$};
 \node[font=\small,red!70!black,anchor=north west] at (3.1,1.95) {$2\min(x,1)$};
\end{tikzpicture}
\caption{The identity $|x-1|=x+1-2\min(x,1)$ for $x\ge0$, used in Step~4 of
the proof of Theorem~\ref{thm:t22-local-mass-flat}. The comparison of
Proposition~\ref{prop:t26-microscopic-tag-comparison} applies only to
bounded continuous tests, and $|x-1|$ is unbounded. The identity splits it
into the linear part $x+1$ (dashed), whose expectation is known exactly by
cyclic invariance, and the bounded part $2\min(x,1)$ (red), to which the
comparison applies. So no uniform integrability of the scaled masses is
needed before the limit.}
\label{fig:lm-identity}
\end{figure}
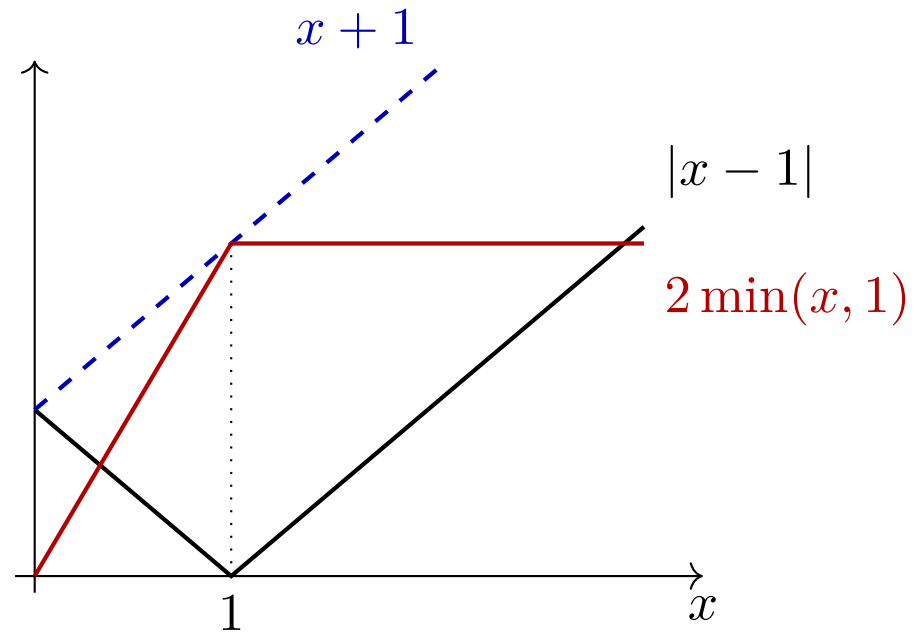

\emph{Step 4: the $L^1$ bound for one block.} The physical statistic in
\eqref{eq:t22-local-mass-flat} involves the unbounded function $|x-1|$.
Cyclic invariance gives exactly
\[
 \mathbb E\sum_\alpha w(u_\alpha)np_{\alpha k}=\mathbb E\sum_\alpha w(u_\alpha).
\]
For $x\ge0$ we use the identity
\[
 |x-1|=x+1-2\min(x,1),
\]
shown in Figure~\ref{fig:lm-identity}. The linear part $x+1$ is handled by
the exact identity above, together with
Theorem~\ref{thm:critical-density} for the limit of
$\E\sum_\alpha w(u_\alpha)$. The bounded part $\min(x,1)$ is handled by the
comparison of Step~3. For the limiting ratio, the same identity gives
$\E|Z/M-1|=2-2\E\min(Z/M,1)$, since $\E(Z/M)=1$ by Step~2. So the limit of
the physical statistic is $\int w(u)\,\E|Z/M-1|\,du\le C\sqrt\tau\int w$.
This proves the assertion for a single block in $L^1$, without any
unproved uniform integrability of the scaled masses before the limit.

\emph{Step 5: the pair, and general $w$.} For the adjacent pair
$S=\{i,i+1\}$, we have $\frac{n}{2}\|P_S\psi_\alpha\|^2=\frac{1}{2}(np_{\alpha i}+np_{\alpha,i+1})$.
So the deviation from one is at most the average of the deviations of
the two blocks, and the pair bound follows from the bound for one block.
A general nonnegative compactly supported $w$ follows by scaling.
\end{proof}

\subsection{The scalar pair flow}\label{subsec:s28-the-scalar-pair-flow}

 We compare eigenvalue locations with the parameters at which a scalar
perturbation of one pair produces a crossing at the observation energy.
Lemma~\ref{lem:t22-crossing-sine} gives sine statistics for the
crossings. We then use local mass flatness to show that their velocities
are nearly one, control the error along the whole flow, and finish the
physical sine endpoint.

 Set $\rho=\rho(E)$ and   $\Naux =2W$. For a fixed adjacent pair $S$, consider

\[
 H(t)=H+\frac{t}{\Naux \rho}P_S,\qquad
 u_\alpha(t)=N\rho(\lambda_\alpha(H(t))-E).
\]

Its exact velocity is
\begin{equation}\label{eq:t22-pair-flow-velocity}
 \dot u_\alpha(t)=\frac{n}{2}\|P_S\psi_\alpha(t)\|^2=:v_\alpha(t)>0.
\end{equation}
Strict positivity follows from the invertible bonds and the eigenvector
equation; a vector vanishing on two adjacent blocks vanishes everywhere.
Eigenvalues are simple along this one-parameter flow almost surely;
alternatively all flow estimates below hold on the almost-everywhere
simple set, which is sufficient for absolute continuity and integration.

\begin{lemma}[The crossing process has a sine limit]
\label{lem:t22-crossing-sine}
For the scalar pair flow $H(t)=H+t(\Naux\rho)^{-1}P_S$, call a real number
$t_j$ a \emph{crossing parameter} if the fixed energy $E$ is an eigenvalue
of $H(t_j)$, and let $\mathcal C_W=\sum_j\delta_{t_j}$ be the point process
of crossing parameters. Then, for every fixed critical ratio, $\mathcal
C_W$ converges in law to the sine process of unit intensity. Note that
$\mathcal C_W$ is a point process in the perturbation parameter $t$ of
the pair, not yet in the energy variable.
\end{lemma}

The crossing parameters are, after rescaling, the
eigenvalues of a fresh deformed Gaussian matrix on the pair. So the sine
kernel theorem for a fresh deformed GUE matrix applies directly in the
parameter $t$. The comparison with the physical energies is made only
afterwards, in the proof of Theorem~\ref{thm:t22-sine-endpoint}, where
local mass flatness makes the eigenvalues move at almost unit speed in
$t$.

\begin{proof}
 We first fix the good event for the exterior,
the part of the ring outside the pair, on which the fresh Gaussian block
is regular (Step~1). Conditionally on the exterior, we show that the
crossing parameters are exactly a rescaling of the eigenvalues of the
Schur complement on the pair, which is a full GUE matrix plus an
independent deformation (Step~2). The sine kernel theorem,
Theorem~\ref{thm:external-sine}, then gives convergence of all fixed
correlation functions (Step~3). Finally, uniform bounds on the number of
eigenvalues in an interval pay for the exceptional exterior event and
identify the limit as the sine process (Step~4).

\emph{Step 1: the regular exterior event.} Write $R=S^c$ for the
exterior of the pair, and recall from
Section~\ref{sec:physical-reconstruction} the deformation
\[
 K(E)=H^{(S)}_{SS}-E-H_{SR}(H_{RR}-E)^{-1}H_{RS},
\]
where $H^{(S)}$ is the matrix with the fresh Gaussian block $G$ removed.
We use the regular exterior event constructed in the proof of
Theorem~\ref{thm:t22-normalizer}, not merely the total-variation
conclusion of that theorem. Schematically, this event is
\begin{equation}\label{eq:crossing-regular-event}
\begin{aligned}
 \mathcal R_W={}&
 \Bigl\{\operatorname{dist}\bigl(E,\operatorname{Spec}H_{RR}\bigr)\ge W^{-A}\Bigr\}
 \cap\Bigl\{\text{bond norms}\le K_0\Bigr\}\\
 &\cap\Bigl\{K(E)\ \text{satisfies \eqref{eq:fresh-regularity}}\Bigr\}.
\end{aligned}
\end{equation}
The three conditions are a cutoff away from the real poles of the
exterior, a cutoff on the norms of the Gaussian bonds, and the fixed
regularity rectangle of the fresh Gaussian block. They are formulated
under the original prior law of the exterior, at the deterministic
energy $E$, and we use them as inputs here; their complement has
probability tending to zero. The real inverse $(H_{RR}-E)^{-1}$ exists
almost surely. The counting estimate we use conditionally on the
exterior is Theorem~\ref{thm:external-count}.

\emph{Step 2: crossing parameters as eigenvalues.} We use the fresh-pair
split of Section~\ref{sec:physical-reconstruction}. Conditionally on the
exterior, the Schur complement of $H-E$ on the pair is $G+K(E)$, where $G$
is a full $\Naux$-dimensional GUE matrix with variance parameter
$s=\min(b,1-2b)$, independent of $K(E)$. The flow adds the multiple
$t(\Naux\rho)^{-1}$ of the identity on the pair, and the Schur
determinant formula gives
\[
 \det\bigl(H(t)-E\bigr)
 =\det\bigl(H_{RR}-E\bigr)\,
  \det\Bigl(G+K(E)+\frac{t}{\Naux\rho}\,I\Bigr).
\]
The first factor does not vanish, because almost surely the exterior has
no eigenvalue exactly at the deterministic energy $E$. Hence
\[
\begin{aligned}
 E\in\operatorname{Spec}H(t)
 &\iff
 0\in\operatorname{Spec}\Bigl(G+K(E)+\frac{t}{\Naux\rho}\,I\Bigr)\\
 &\iff
 -\frac{t}{\Naux\rho}\in\operatorname{Spec}\bigl(G+K(E)\bigr).
\end{aligned}
\]
If $\mu_j$ are the eigenvalues of $G+K(E)$, the crossing parameters are
therefore exactly
\begin{equation}
 t_j=-\Naux \rho\mu_j.                                      \label{eq:pair-crossing-locations}
\end{equation}
In particular, no term coming from a pole of the exterior is missing in
this identity.

\emph{Step 3: convergence of the correlation functions.} Put
\[
 \mathcal M_W=\sum_j\delta_{\Naux\rho\mu_j},
\]
the rescaled eigenvalue process of the fresh block. On the regular event
$\mathcal R_W$, the deformation $K(E)$ has a polynomially bounded norm,
its free convolution with the fresh Gaussian block is regular on a fixed
rectangle of energies and variances, and its density at zero satisfies
$r_{K(E),s}(0)\to\rho$. We apply Theorem~\ref{thm:external-sine} at every
fixed correlation order $k$. Its hypotheses hold: the variance parameter
$s$ is fixed and positive, which is exactly the scope of that theorem,
proved directly from the complex kernel. The theorem gives, uniformly for
$\xi,\zeta$ in compact sets,
\[
 \frac{e^{-h(\xi)}\,K_\Naux\bigl(\xi/(\Naux\rho),\zeta/(\Naux\rho)\bigr)\,
       e^{h(\zeta)}}{\Naux\rho}
 \longrightarrow K_{\rm sine}(\xi-\zeta),
 \qquad K_{\rm sine}(x)=\frac{\sin\pi x}{\pi x},
\]
for a real conjugation function $h$. Since the local density tends to
$\rho$, the rescaling by the density tends to one uniformly on compact
sets of coordinates. The conjugation factors cancel in each correlation
determinant, because
\[
 \det\Bigl[e^{-h(\xi_i)}K(\xi_i,\xi_l)e^{h(\xi_l)}\Bigr]_{i,l=1}^k
 =\det\bigl[K(\xi_i,\xi_l)\bigr]_{i,l=1}^k .
\]
Consequently, for every fixed $k$ and every smooth compactly supported
$\varphi$ on $\mathbb R^k$,
\begin{equation}\label{eq:crossing-correlations}
\begin{aligned}
 &\E\Bigl[\one_{\mathcal R_W}\!\!\sum_{j_1,\ldots,j_k\ \rm distinct}
        \varphi\bigl(\Naux\rho\mu_{j_1},\ldots,\Naux\rho\mu_{j_k}\bigr)\Bigr]\\
 &\qquad\longrightarrow
 \int_{\mathbb R^k}\varphi(x)\,
   \det\bigl[K_{\rm sine}(x_i-x_l)\bigr]_{i,l=1}^k\,dx .
\end{aligned}
\end{equation}

\emph{Step 4: counting bounds and identification of the limit.} The
counting bound of Theorem~\ref{thm:external-count}, for the deformed GUE
block of variance $s$, holds uniformly in the deformation. For a bounded
interval $I$ and the count $\mathcal M_W(I)$, it gives
\[
 \E\bigl[(\mathcal M_W(I))_r\bigr]\le(C|I|)^r
 \quad(r\ge1),
\]
and hence
\[
 \E\,e^{\lambda\mathcal M_W(I)}\le\exp\{C|I|(e^\lambda-1)\}
 \qquad(\lambda\ge0),
\]
where $(x)_r=x(x-1)\cdots(x-r+1)$ is the falling factorial; the second
bound follows from the first by the binomial expansion used in the proof
of Theorem~\ref{thm:counting}. These uniform moments do three things.
First, they pay for the exceptional exterior event: for a statistic $X_W$
built from finitely many counts,
\[
 \bigl|\E\bigl[X_W\one_{\mathcal R_W^c}\bigr]\bigr|
 \le\bigl(\E X_W^2\bigr)^{1/2}\,\Prob(\mathcal R_W^c)^{1/2}
 \longrightarrow0,
\]
so the indicator in \eqref{eq:crossing-correlations} can be removed.
Second, they give tightness of the point processes $\mathcal M_W$.
Third, they make the limiting family of correlation functions determine
the law. The correlation functions on compact sets therefore identify the
limit of $\mathcal M_W$ as the complex Hermitian sine process.

Finally, by \eqref{eq:pair-crossing-locations} the crossing process is
the reflection of $\mathcal M_W$, that is,
$\mathcal C_W=\sum_j\delta_{-\Naux\rho\mu_j}$. The sine kernel is even, so
\[
 \det\bigl[K_{\rm sine}(-x_i+x_l)\bigr]_{i,l=1}^k
 =\det\bigl[K_{\rm sine}(x_i-x_l)\bigr]_{i,l=1}^k ,
\]
and the sine process is invariant under reflection. Hence $\mathcal C_W$
has the same sine limit. The lemma is proved under the original law of
the exterior; the crossing parameters are not reweighted by the
selection of an eigenvalue of the periodic ring.
\end{proof}

\begin{lemma}[Admission of the velocity error along the flow]
\label{lem:t22-flow-admission}
For every fixed $T,K<\infty$ and smooth compactly supported
$\chi\ge0$, equal to one on $[-K,K]$,
\begin{equation}\label{eq:t22-flow-admission}
 \limsup_W\mathbb E\int_{-T}^T\sum_\alpha\chi(u_\alpha(t))
                           |v_\alpha(t)-1|\,dt
                  \le C_{T,\chi}\tau^{1/4}.
\end{equation}
No density of a selected eigenvalue or unproved inverse clock moment
is used in this change of measure.
\end{lemma}

The likelihood-weighted speed error is measured under the original matrix law. Pair-mass coarea cancellation prevents a factor equal to the number of blocks.

\begin{proof}
Only the   $\Naux $  real diagonal entries in $S$ are shifted. Their original
variance is $c/W$, $c=1-2b$. The likelihood $R_t$ of the shifted matrix
law relative to the original law obeys exactly
\begin{equation}
 R_t^2=\exp\{t^2/(2c\rho^2)\}\,R_{2t},
 \qquad \mathrm{KL}(H(t)\Vert H)=t^2/(4c\rho^2).        \label{eq:pair-shift-likelihood}
\end{equation}
Put $x_\alpha=n\|P_S\psi_\alpha\|^2/2$. Cauchy--Schwarz with respect
to the nonnegative measure $d\mathbb P\sum_\alpha
\chi(u_\alpha)|x_\alpha-1|$ gives
\begin{align*}
 \mathbb E_t\sum_\alpha\chi(u_\alpha)|x_\alpha-1|
 &\le\left(\mathbb E\sum_\alpha\chi(u_\alpha)|x_\alpha-1|\right)^{1/2}
 \left(\mathbb E\sum_\alpha\chi(u_\alpha)|x_\alpha-1|R_t^2\right)^{1/2}.
\end{align*}
The first factor is at most $C\tau^{1/4}$ in limsup by
Theorem~\ref{thm:t22-local-mass-flat}. For the second use
$|x-1|\le x+1$ and \eqref{eq:pair-shift-likelihood}. Counts under any deterministic diagonal
shift satisfy the bound \eqref{eq:external-count}. The mass-weighted
term is also uniformly bounded: \eqref{eq:t22-mass-palm}
cancels the crossing derivative and its deformed fresh-GUE density is
at most   $C_b\Naux $, regardless of the scalar shift. Multiplication by $n/2$
and the energy-window Jacobian   $1/(N\rho)$  cancels the dimension,
since   $n\Naux /(2N)=1$. Hence
\[
 \sup_{|t|\le T,W}\mathbb E_{2t}\sum_\alpha
                       \chi(u_\alpha)(x_\alpha+1)<\infty.
\]
This bounds the second factor uniformly. All estimates are uniform
on the compact parameter interval, so integration proves
\eqref{eq:t22-flow-admission}. One may first keep a finite-$W$ error
in the local mass bound and then let $W$ tend to infinity; no exchange
with an uncontrolled supremum in energy is needed.
\end{proof}

\begin{proof}[Proof of Theorem~\ref{thm:t22-sine-endpoint}]
Choose $R$ containing the support of $f$, $T>R+3$, and then $K>R+T+5$.
Let $A_W$ be the nonnegative integrated velocity error in
\eqref{eq:t22-flow-admission}, with $\chi=1$ on $[-K,K]$.
On $A_W\le\epsilon<1$, every eigenvalue branch starting in
$[-R-1,R+1]$ remains within the chosen larger window until times
$\pm T$ and satisfies
\[
 |u_\alpha(t)-u_\alpha(0)-t|\le\epsilon.
\]
Indeed this follows by integrating its velocity error until a first
exit; the inequality itself precludes that exit. Strict monotonicity
then gives exactly one zero, with
$|t_\alpha+u_\alpha(0)|\le\epsilon$.
Conversely a branch crossing at $|t_j|\le R+1$ satisfies the same
estimate back to time zero. These correspondences are injective,
because each increasing branch crosses zero at most once.

Thus the two locally relevant point configurations
$\{u_\alpha(0)\}$ and $\{-t_j\}$ are paired within $\epsilon$.
Extending $f$ by zero, its Lipschitz constant bounds the difference
of the two Laplace functionals on this event by $C_f\epsilon$ times
the sum of their counts in $[-R-2,R+2]$. Both expected counts are
uniformly bounded: use original counting for the first process and
the conditional fresh-GUE count bound for the second. On the
complement each Laplace functional lies in $[0,1]$, so Markov gives
\[
 \left|\mathbb E e^{-\sum f(u_\alpha(0))}
          -\mathbb E e^{-\sum f(-t_j)}\right|
 \le C_f\epsilon+\epsilon^{-1}\mathbb EA_W.
\]
Take the critical limsup, apply Lemmas
\ref{lem:t22-crossing-sine} and \ref{lem:t22-flow-admission}, and
choose $\epsilon=\tau^{1/8}$. This proves
\eqref{eq:t22-sine-endpoint}. These compactly supported Laplace
functionals determine locally finite point-process laws. The
count bounds \eqref{eq:factorial-count} give tightness and extend from smooth tests to continuous
ones. The sine process is simple, stationary, and has factorial density
\[
 \det\left[\frac{\sin\pi(u_i-u_j)}{\pi(u_i-u_j)}\right]_{i,j=1}^p,
\]
with diagonal one. The unfolding in \eqref{eq:pair-crossing-locations} and
\eqref{eq:t22-pair-flow-velocity} fixes its intensity at one.
\end{proof}

\paragraph{Why the comparison works.}
At a fixed real energy the fresh full Gaussian pair supplies a sine
process in a scalar perturbation parameter. Strong local mass flatness
makes that parameter flow through the physical energy window at nearly
unit speed. The pair-mass coarea cancellation controls the likelihood
weighted speed error, so Gaussian change of measure costs no factor of
$n$ or $W$. This uses an actual periodic matrix throughout. It neither
replaces the ring by an open chain at the spectral level nor infers
spectral universality from the earlier external-supremum QUE bound alone.
\section{Fractional Schatten decay of remote complete clocks}
\label{sec:t22-remote-schatten}

The Poisson endpoint is proved by cutting the ring and comparing it with
smaller rings. For this we must know that an absorber placed far inside a
long path is barely seen at its boundaries. This section quantifies that
decay. Theorem~\ref{thm:t22-near-Schatten-decay} shows that the fractional
trace of a clock decays exponentially in the time after its last
injection. Theorem~\ref{thm:t22-full-remote-Schatten} shows that the
complete response of a source at distance $R$ from both ends of the path
has fractional trace of order $e^{-cR}$. Theorem~\ref{thm:t22-remote-higher}
extends this to products of fractional traces, when their total degree is
below one. The argument proceeds in four stages. We first collect the
properties of one fresh layer that drive the propagation of a clock
(Lemma~\ref{lem:t22-one-port-inputs}), and bound the clock while and
shortly after its source is injected (Lemma~\ref{lem:edge-clock-moments}).
We then show that the fractional trace is finite after a positive travel
time (Lemma~\ref{lem:t22-clock-fractional-entrance}), use the concavity of
fractional powers to prove strict decay, and finally combine the two
endpoints of the path. In
Subsection~\ref{subsec:s29-higher-fractional-moments} we prove the bounds
for higher fractional moments.

 Throughout this section the law is that of the
actual independent Gaussian path with its two matched leads, at a
deterministic bulk energy. The constants depend on the fixed $b$ and are
uniform for $E$ in a compact bulk set; for the endpoint assertions here,
$(b,E)$ are fixed before the constants. All statements take $W\to\infty$
at each fixed finite path length, measured in units of $W$, before that
length is increased. We set $D=D_0$. The sources below are on physical
blocks and are normalized by $W^2$. The reference path has no physical
absorption: the positive clocks below are derivatives of that unabsorbed
path with respect to a source, and their injections end before the
homogeneous fractional decay is applied.

\begin{figure}[ht]
\centering
\begin{tikzpicture}[x=1.35cm,y=1cm,>=Stealth]
 \draw[->] (-.2,0)--(11.3,0) node[below,font=\small] {scaled time};
 \draw[thick] (0,-.1)--(0,.1) node[below=5pt,font=\small] {$0$};
 \fill[black!8] (0,.15) rectangle (3.2,.95);
 \node[font=\scriptsize,align=center] at (1.6,.55) {old path,\\ about $t_1W$ blocks};
 \fill[red!20] (3.2,.15) rectangle (5.0,.95);
 \node[font=\scriptsize,red!70!black] at (4.1,.55) {injection};
 \fill[blue!12] (5.0,.15) rectangle (6.8,.95);
 \node[font=\scriptsize,blue!70!black] at (5.9,.55) {transport};
 \fill[green!12] (6.8,.15) rectangle (11.0,.95);
 \node[font=\scriptsize,green!40!black,align=center] at (8.9,.55)
   {homogeneous propagation\\ $A_+=L^*AL$};
 \foreach \x/\l in {3.2/t_1,5.0/{t_1+R},6.8/{t_1+R+R_0}}
   {\draw[thick] (\x,-.1)--(\x,.1); \node[font=\small,below] at (\x,-.12) {$\l$};}
 \draw[decorate,decoration={brace,amplitude=4pt}] (3.2,1.15)--(6.8,1.15)
   node[midway,above=4pt,font=\small] {bounded trace (Lemma~\ref{lem:edge-clock-moments})};
 \draw[decorate,decoration={brace,amplitude=4pt}] (6.85,1.15)--(11.0,1.15)
   node[midway,above=4pt,font=\small] {fractional decay of $H_p$};
\end{tikzpicture}
\caption{The life of a clock, drawn schematically. The old path (grey)
can be arbitrarily long. The source is injected during an interval of
length $R\le R_0$ (red); during the injection and for a further time
$R_0$ (blue), Lemma~\ref{lem:edge-clock-moments} bounds all moments of its
trace, with constants that do not depend on the length of the old path.
After the last injection, each fresh layer acts on the near-near block
$A$ by the congruence $A_+=L^*AL$ (green), and the fractional trace
$H_p=\Tr A^p$ decays exponentially
(Theorem~\ref{thm:t22-near-Schatten-decay}). The properties of the
single layer that drive this propagation are collected in
Lemma~\ref{lem:t22-one-port-inputs}.}
\label{fig:clock-timeline}
\end{figure}
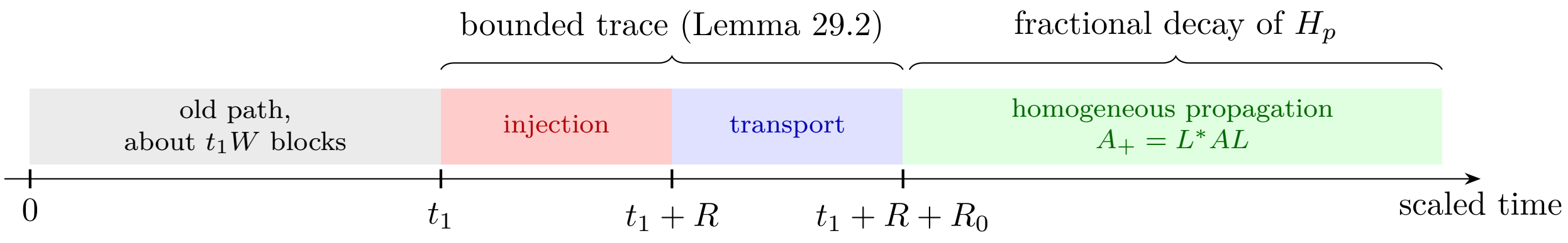

\paragraph{The clock after its last injection.} Consider a positive source
on an interval of scaled length at most one, and let $A$ be the near-near
block of its clock after the last injection; see
Figure~\ref{fig:clock-timeline}. During the subsequent propagation, each
fresh layer acts on $A$ exactly by congruence,
\[
 A_+=L^*AL .
\]
We write
\[
 \gamma=\operatorname{Tr}A,\qquad H_p=\operatorname{Tr}A^p\quad(0<p<1),
\]
and we let $T$ be the current transmission matrix, with conductance
$g=\operatorname{Tr}T$. Note that a fractional moment of the trace $\gamma$
alone does not bound $H_p$: a matrix of small trace can still have a large
fractional trace if its mass is spread over many small eigenvalues. This
is why $H_p$ is controlled directly below.

\paragraph{The fresh one-port layer.} The propagation $A_+=L^*AL$ uses the
fresh one-port layer of Section~\ref{sec:transmission}, attached at the
near end of the path. The next lemma collects the three properties of this
layer that we use: bounds on the mean of $LL^*$, moment bounds on $L$, and
the contraction of the scalar error together with a weight that absorbs
it.

\begin{lemma}[One-port passive inputs]\label{lem:t22-one-port-inputs}
Let $X$ be the near-end Schur matrix of the old path, with its far port
eliminated as in Section~\ref{sec:transmission}, so that $-\Im X\succeq0$.
Let $K=z-X$ and $L=K\mathcal G_{12}$ for the fresh one-port layer of
\eqref{eq:passive-update}, let
\[
 \delta_0=|W^{-1}\Tr(z-X)^{-1}-a|
\]
be the scalar error of the old input, and let $\mathbb E_j$ integrate the
fresh layer at fixed $X$. The constants depend only on $b$ and on a
compact energy set $I$, they are chosen before $W$ and $X$, and the
assertions hold for every such $X$ and every $E\in I$.
\begin{enumerate}
\item[(a)] \emph{The mean of $LL^*$.} For $W\ge W_0$,
 \[
  \mathbb E_jLL^*\preceq(1+C\delta_0+CW^{-2})I,
 \]
 and $\mathbb E_jLL^*\succeq I/2$ if moreover $\delta_0$ and $g/W$ are at
 most a small fixed constant.
\item[(b)] \emph{Moments of $L$.} $\mathbb E_j\norm L^p\le C_p$ for $p\ge1$
 and $W\ge W_0(p)$. Consequently, for $A\succeq0$ and $0<p<1$,
 \[
  \Tr(L^*AL)^p\le\norm L^{2p}\Tr A^p,
  \qquad
  \mathbb P_j(\norm{LL^*}>K')\le C_4K'^{-2}.
 \]
\item[(c)] \emph{Contraction of the scalar error, and a weight.} Let
 $\delta_{0,+}$ be the scalar error of the new Schur matrix $X_+$, and
 \[
  \xi_+=W^{-1}\Tr(z-X_+)^{-1}-\mathbb E_jW^{-1}\Tr(z-X_+)^{-1}
 \]
 its fluctuation. There are $d_0>0$, $q_*<1$, $C$ and $W_0$ such that, on
 $\delta_0\le d_0$,
 \[
  \delta_{0,+}\le q_*\delta_0+CW^{-2}+|\xi_+|,
 \]
 and $\norm{\xi_+}_{L^r\mid X}\le C_r/W$ for fixed $r\ge1$ and
 $W\ge W_0(r)$. Consequently, let $C_*\ge0$, choose $K_p$ with
 $K_p(1-q_*)\ge2C_*$, and reduce $d_0$ so that $K_pq_*d_0\le1$. Then the
 weight $w=1+K_p\delta_0$ satisfies, on $\delta_0\le d_0$,
 \[
  w_+\le\tilde w+K_p|\xi_+|,\qquad
  \tilde w=1+K_pq_*\delta_0+K_pCW^{-2},\qquad
  \tilde w(1+C_*\delta_0+C_*W^{-2})\le w(1+C'W^{-2}).
 \]
 Here $\tilde w$ is measurable with respect to the past,
 $1\le\tilde w\le3$ for large $W$, and $C'$ depends only on $C_*$, $K_p$
 and $C$.
\end{enumerate}
\end{lemma}

The weight in (c) is the tool that lets a growth factor of the form
$1+C_*\delta_0$, which comes from the scalar error, be absorbed by the
contraction of that error, so that it costs only $1+C'W^{-2}$ per layer.

\begin{proof}
\emph{(a) The mean of $LL^*$.} On a fixed small scalar domain, the
expansion \eqref{eq:passive-cancelled-mean} applies, with $|\Omega|\le
C\delta_0$ and
\[
 \norm f\le\norm R+(2\eta)^{-1}\le\frac{3}{2\eta},
\]
since $\Im K\succeq\eta I$. Since $1+\alpha q/W\ge1$ and $q\le W/(4\eta)$,
the latter because $T=4\eta Q\preceq I$, the remainder in that expansion
is $O(\delta_0+W^{-2})$. The numerator lies between
$(1-3|\Omega|/\eta)I$ and $(1+3|\Omega|/\eta)I$, and dividing by
$1+\alpha q/W\ge1$ keeps the upper bound. This gives the upper bound.

For the lower bound, let $\delta_0\le\varepsilon$ and $g/W\le\varepsilon$.
Then $q/W=g/(4\eta W)\le\varepsilon/(4\eta)$, and the remainder is at least
$-C(\varepsilon+W^{-2})I$. Hence
\[
 \mathbb E_jLL^*\succeq\Bigl(\frac{1-C\varepsilon}{1+\alpha\varepsilon/(4\eta)}
   -C\varepsilon-CW^{-2}\Bigr)I\succeq\frac I2
\]
for $\varepsilon$ small and $W\ge W_0$.

Off the scalar domain, part~(b) gives
$\norm{\mathbb E_jLL^*}\le C_2\le C_2d^{-1}\delta_0$, where $d$ is the
radius of the domain, so the upper bound holds there too.

\emph{(b) Moments of $L$.} The moment bound is the first assertion of
Lemma~\ref{lem:passive-bounds}. The fresh layer there is the present one,
and its proof uses only the passive inverse moments and the Gaussian
moments, uniformly in $X$. The trace inequality holds because $L^*AL$ and
$A^{1/2}LL^*A^{1/2}$ have the same nonzero eigenvalues, and
$A^{1/2}LL^*A^{1/2}\preceq\norm L^2A$; since $x\mapsto x^p$ is monotone,
$\Tr(A^{1/2}LL^*A^{1/2})^p\le\norm L^{2p}\Tr A^p$. The tail bound is
Markov's inequality applied to $\norm L^4$.

\emph{(c) The scalar error and the weight.} The contraction and the
fluctuation bound are the passive scalar recursion of Step~3 of the proof
of Lemma~\ref{lem:t26-passive-initial-scalar}, which is also used in the
proof of \eqref{eq:full-clock-moments}; its rate is
\[
 q_*=\frac{b}{1-|k_-|d_*}+O(d_0)<1 .
\]
Multiplying the contraction by $K_p$ and adding one gives
$w_+\le\tilde w+K_p|\xi_+|$. For the last inequality we expand:
\[
 \tilde w(1+C_*\delta_0+C_*W^{-2})-1
 \le(K_pq_*+C_*+K_pq_*C_*\delta_0)\delta_0+C'W^{-2}.
\]
Since $K_pq_*d_0\le1$, the coefficient of $\delta_0$ satisfies
$K_pq_*+C_*+K_pq_*C_*\delta_0\le K_pq_*+2C_*\le K_p$, the last step by
$K_p(1-q_*)\ge2C_*$. As $w\ge1$, this proves
$\tilde w(1+C_*\delta_0+C_*W^{-2})\le w(1+C'W^{-2})$.
\end{proof}

\paragraph{The clock during and after injection.} Before the homogeneous
decay can start, we must control the clock while its source is being
injected and for a bounded time afterwards, with constants that do not
depend on how long the old path is. This uniformity is needed for the
spectral cutting in Section~\ref{sec:t22-positive-log-poisson}, where the
source sits far inside a long path.

\begin{lemma}[Uniform moments after a bounded injection interval]
\label{lem:edge-clock-moments}
Consider a clock started from zero, with injection profile $0\le q\le M$,
injected during a scaled interval of length at most $R_0$, and then
propagated for at most another $R_0$. Then, for every fixed $p>0$,
\[
 \limsup_W\norm{\Tr Q}_{L^p}\le MC_{p,R_0}
\]
in the critical limsup, and the same bound holds for $\Tr A$. The
constants do not depend on the length of the preceding actual path.
\end{lemma}

The bound depends only on the duration of the injection and of the
subsequent transport. It is uniform in the length of the old middle part
of the path, which is what the spectral cutting requires.

\begin{proof}
Let the test satisfy $0\le q\le M$ and vanish outside the injection
interval $[t_1,t_1+R]$, with $R\le R_0$; the preceding actual path has
about $t_1W$ blocks. Each tag clock $Y_{j;i}$ is positive semidefinite, so
the clock grows monotonically with the injection profile, and the
near-near block is part of the full clock:
\[
 0\le\Tr A\le\Tr Q_j\le M\Tr Q_j^{(1)},
\]
where $Q_j^{(1)}$ is the clock with the test $\mathbf 1_{[t_1,t_1+R]}$.
By Lemma~\ref{lem:source-recorded-estimates}(e) with $T=t_1+R+R_0$,
\[
 \norm{\Tr Q_j}_{L^p}\le MC_{p,R_0}
 \qquad\text{for } p\ge1,\ W\ge W_0(T),\ j\le(t_1+R+R_0)W,
\]
which covers both the injection and the subsequent propagation. Lower
moments, $0<p<1$, follow by Jensen's inequality. The constant $C_{p,R_0}$
is chosen before $T$ and $t_1$; only the threshold $W_0(T)$ depends on
them. Letting $W\to\infty$ at fixed lengths therefore gives the bound in
the critical limsup, without any dependence on the length of the old
path.
\end{proof}

\begin{lemma}[Finite fractional entrance after positive travel]
\label{lem:t22-clock-fractional-entrance}
For every $0<p<1$ and $\epsilon>0$, the actual near source satisfies
\[
 \limsup_W\mathbb E H_p(\epsilon)\le C_{p,\epsilon},
\]
where time zero denotes its last injection. The constant is independent
of the preceding actual middle length. On each subsequent fixed
homogeneous time interval these expectations are bounded uniformly in
$W$, with a constant permitted to depend on that interval.
For $p<q<1$, no vanishing eigenvalues carry hidden $q$th power:
\begin{equation}
 \mathbb E\sum_{\lambda_i(A)\le a}\lambda_i(A)^q
                       \le a^{q-p}\mathbb E H_p.
                                                               \label{eq:small-clock-channels}
\end{equation}
\end{lemma}

The boundary starts used here are the $p=1$ entrance faces $x=0<y$ or $y=0<x$, not the paired corner $(0,0)$, where the corner convention is $J(0,0)=0$. No boundary-start extension for $p\ge2$ is asserted. A bound on the trace alone would allow mass in arbitrarily many small eigenvalues. The fractional estimate rules out that loss before homogeneous transport is used.

\begin{proof}
 We follow the fractional trace $H_p=\Tr A^p$
layer by layer, at finite width. One fresh layer acts on $A$ by the
congruence $A_+=L^*AL$. Since $x\mapsto x^p$ is concave, the random
fluctuation of $LL^*$ decreases the expected value of $H_p$, and the
decrease is quadratic in $H_p$. To make this quantitative we keep only the
eigendirections of $A$ in which the transmission is small (Step~1), show
that $LL^*$ fluctuates by a definite amount in these directions (Step~2),
and turn this into a concavity gap of order $H_p^2/W$ (Step~3). This
gives a one-step inequality of logistic type (Step~4). A weight absorbs
the scalar error (Step~5), and the exceptional events and the moment
inputs are handled in Step~6. A discrete logistic comparison then brings
$H_p$ down from its crude initial size to a bounded value after any
positive travel time (Step~7), as in Figure~\ref{fig:clock-logistic}.
Step~8 proves the pointwise bound \eqref{eq:small-clock-channels}.
 
\emph{Step 1: the setting and the retained directions.} We work at finite
$W$, on a domain of positive time on which the scalar and jet errors are
small and $g/W$ is sufficiently small, as in the proof of the fractional
transmission bound in Section~\ref{sec:transmission-law}. In an
eigenbasis of $A$, with eigenvalues $\lambda_i$, we retain the indices
with $T_{ii}\le1/2$, that is, the directions in which the transmission is
at most one half. Since $\Tr T=g$, at most $2g$ indices are discarded.
 
\emph{Step 2: a lower bound on the fluctuation of $LL^*$.} The gap
$\sigma_+>|\sigma_-|$ between the opposite-sign and the equal-sign
covariance coefficients, \eqref{eq:jet-gap} in
Lemma~\ref{lem:source-recorded-estimates}(c), gives, for retained $i,j$,
\begin{equation}
 \mathbb E_j|(LL^*-\mathbb E_jLL^*)_{ij}|^2\ge c/W.
                                                               \label{eq:clock-covariance-gap}
\end{equation}
This holds in the eigenbasis of $A$, and not only in an eigenbasis of the
transmission. To see this, write $r$ for the old reflection and
$S=r^*r=I-T$. Lemma~\ref{lem:passive-covariance} gives the variance of the
entry, up to its stated remainder, as
\[
 \frac{\sigma_+}{4\eta^2W}(S_{ii}+S_{jj})
 -\frac{1}{2\eta^2W}\Re(\sigma_-r_{ii}r_{jj}),
\]
because the fast source there is $f=-(2i\eta)^{-1}r$. The diagonal
inequality $|r_{ii}|\le\sqrt{S_{ii}}$ holds in every orthonormal basis, so
$|r_{ii}r_{jj}|\le\sqrt{S_{ii}S_{jj}}\le(S_{ii}+S_{jj})/2$, and the
expression is at least
\[
 \frac{(\sigma_+-|\sigma_-|)(S_{ii}+S_{jj})}{4\eta^2W}.
\]
Both retained diagonal entries satisfy $S_{ii},S_{jj}\ge1/2$, and the
fixed small errors in the scalar, the jet and $g/W$ preserve half of this
strict lower bound.
 
We also restrict to an event on which $LL^*$ is bounded. The fourth
bilinear moment and the operator moments of the fresh layer, from
Lemma~\ref{lem:t22-one-port-inputs}(b), let us retain the event
$\|LL^*\|\le K$, for a fixed $K$, while keeping half of the lower bound
\eqref{eq:clock-covariance-gap}.
 
\emph{Step 3: the concavity gap.} We apply the integral formula
\eqref{eq:fractional-hessian} for the Hessian of $\Tr B^p$, after a
positive regularization, to $B=A^{1/2}LL^*A^{1/2}$. On the event
$\|LL^*\|\le K$, the segment between $B$ and its conditional mean is
bounded by $KA$. The concavity gap, that is, the amount by which the
expectation of $\Tr B^p$ falls below its value at the mean, is therefore
at least a constant times
\[
 \frac{1}W\sum_{i,j\ {\rm retained}}\lambda_i\lambda_j
     \int_0^\infty
       \frac{s^{p-1}\,ds}{(K\lambda_i+s)(K\lambda_j+s)}
 \ \ge\ \frac{c_p}{W\gamma^p}
               \left(\sum_{i\ {\rm retained}}\lambda_i^p\right)^2.
\]
For the last inequality we use $\lambda_i\le\gamma$ and the divided
difference of $x^{1-p}$, whose derivative is at least $(1-p)\gamma^{-p}$
on $(0,\gamma]$. The case of coinciding eigenvalues and of zero
eigenvalues follows by positive regularization.
 
The discarded indices do not spoil this bound. Their power sum is at most
$(2g)^{1-p}\gamma^p$, by concavity of $x\mapsto x^p$ applied to at most
$2g$ eigenvalues of total mass at most $\gamma$.
 
\emph{Step 4: the one-step inequality.} The passive mean bound
$\mathbb E_jLL^*\preceq(1+C\delta_0+CW^{-2})I$ of
Lemma~\ref{lem:t22-one-port-inputs}(a) controls the value at the mean.
Together with Step~3 it gives
\begin{equation}
 \mathbb E_j H_{p,+}\le(1+C_p\delta_0+C_pW^{-2})H_p
 -\frac{c_p}{W}\frac{H_p^2}{\gamma^p}
 +\frac{C_p}{W}g^{2(1-p)}\gamma^p.                    \label{eq:clock-fractional-logistic}
\end{equation}
The quotient $H_p^2/\gamma^p$ is set to zero when $A=0$. The middle term
is the gain from concavity, and the last term pays for the discarded
directions.
 
\emph{Step 5: absorbing the scalar error by a weight.} The factor
$1+C_p\delta_0$ in \eqref{eq:clock-fractional-logistic} could accumulate
over many layers, so we absorb it with the scalar weight $1+K_p\delta_0$
used in Section~\ref{sec:transmission-law}. Its predictable contraction
absorbs the term $C_p\delta_0H_p$, and its centered fluctuation, which is
correlated with $H_{p,+}$, costs only $C_pH_p/W$, since
$H_{p,+}\le\|L\|^{2p}H_p$. The negative term is kept.
 
Precisely, let $w=1+K_p\delta_0$, and let $\tilde w$ be as in
Lemma~\ref{lem:t22-one-port-inputs}(c) with $C_*=C_p$. By part~(b) of that
lemma, $\mathbb E_jH_{p,+}^2\le C_pH_p^2$. Then
\eqref{eq:clock-fractional-logistic} gives, on the scalar domain,
\[
 \mathbb E_j[w_+H_{p,+}]\le w(1+C_p'W^{-2})H_p
   -\frac{c_p}{W}\frac{H_p^2}{\gamma^p}
   +\frac{C_p'}{W}\bigl(g^{2(1-p)}\gamma^p+H_p\bigr).
\]
Indeed, $w_+\le\tilde w+K_p|\xi_+|$; the term with $\tilde w$ is handled by
\eqref{eq:clock-fractional-logistic} and the last inequality of
Lemma~\ref{lem:t22-one-port-inputs}(c), and the term with $K_p|\xi_+|$ by
the Cauchy--Schwarz inequality, $\norm{\xi_+}_{L^2\mid X}\le C/W$, and the
bound on $\mathbb E_jH_{p,+}^2$.
 
\emph{Step 6: exceptional events and moment inputs.} The exceptional sets,
where the scalar error or the conductance leaves its domain, are paid for
with arbitrary polynomial failure powers, together with the crude bound
$H_p\le W^{1-p}\gamma^p$. All fixed moments of the source trace on the
bounded intervals used here are uniform in the length of the preceding
middle part of the path, by Lemma~\ref{lem:edge-clock-moments} and
\eqref{eq:edge-clock-finite}, whose constant is chosen before the total
length. The moments of the conductance are uniform at absolute times
bounded away from zero, by Theorem~\ref{thm:transmission-generator}. We
may start the comparison of Step~7 after half of the prescribed positive
travel time, so these moment hypotheses do not require any control at the
entrance time zero.
 
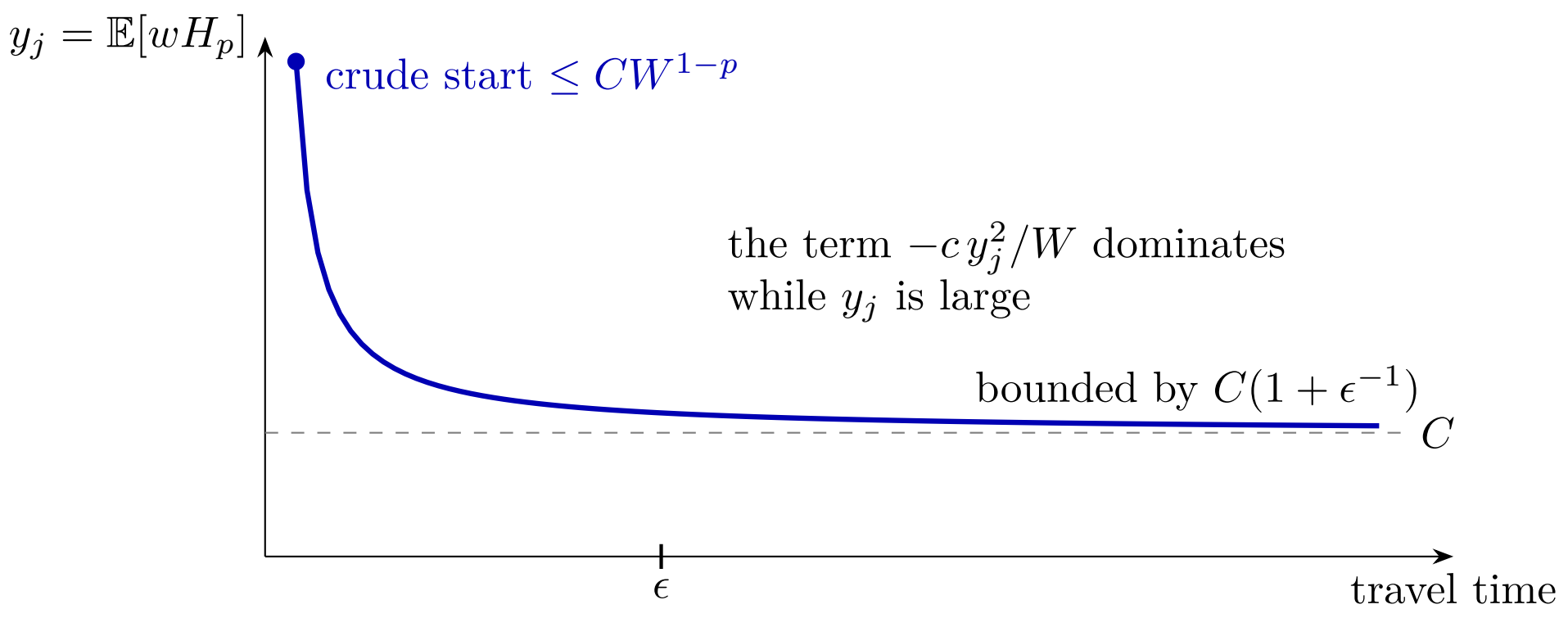
\begin{figure}[ht]
\centering
\begin{tikzpicture}[x=1cm,y=1cm,>=Stealth]
 \draw[->] (0,0)--(9.6,0) node[below,font=\small] {travel time};
 \draw[->] (0,0)--(0,4.2) node[left,font=\small] {$y_j=\mathbb E[w H_p]$};
 \draw[dashed,black!50] (0,1.0)--(9.2,1.0) node[right,font=\small,black] {$C$};
 \draw[very thick,blue!70!black,domain=0.25:9,samples=100]
   plot (\x,{1.0+3.0/(1+6*(\x-0.25))});
 \fill[blue!70!black] (0.25,4.0) circle (2pt);
 \node[font=\small,blue!70!black,anchor=west] at (0.35,3.9) {crude start $\le CW^{1-p}$};
 \draw[thick] (3.2,-.1)--(3.2,.1) node[below=4pt,font=\small] {$\epsilon$};
 \node[font=\small,align=left,anchor=west] at (3.6,2.3)
   {the term $-c\,y_j^2/W$ dominates\\ while $y_j$ is large};
 \node[font=\small,anchor=west] at (5.6,1.35) {bounded by $C(1+\epsilon^{-1})$};
\end{tikzpicture}
\caption{The discrete logistic comparison in Step~7, drawn schematically.
The weighted fractional trace $y_j$ starts from the crude bound
$CW^{1-p}$, which grows with the width. The recursion
$y_{j+1}\le y_j+W^{-1}(Cy_j-cy_j^2+C)$ has a negative quadratic term,
which dominates while $y_j$ is large and brings $y_j$ down, on the time
scale of $W$ layers, to a level that depends only on the elapsed travel
time. After any positive travel time $\epsilon$ the bound is
$C(1+\epsilon^{-1})$, uniformly in $W$ and in the length of the old path.}
\label{fig:clock-logistic}
\end{figure}
 
\emph{Step 7: the logistic comparison.} Let $y_j$ be the weighted
expectation. By the Cauchy--Schwarz inequality,
\[
 \mathbb E\frac{H_p^2}{\gamma^p}\ge\frac{(\mathbb EH_p)^2}{\mathbb E\gamma^p}.
\]
The denominator and the forcing in \eqref{eq:clock-fractional-logistic}
are bounded uniformly on the fixed comparison interval, by Step~6. Since
the weight lies between $1$ and $3$, this gives the recursion
\[
 y_{j+1}\le y_j+W^{-1}(C y_j-c y_j^2+C)+o(W^{-2}).
\]
Its initial value is at most $CW^{1-p}$. As in
Figure~\ref{fig:clock-logistic}, discrete logistic comparison gives the
bound $C(1+\epsilon^{-1})$ at any positive elapsed time $\epsilon$, with
constants that also contain the bounded source moments. Dropping the
negative term then propagates a finite bound over any later fixed time
interval. This proves the assertions at finite $W$, and hence the
assertions in the critical limsup.
 
\emph{Step 8: the small-channel bound.} Formula \eqref{eq:small-clock-channels}
holds pointwise, before taking expectations: for an eigenvalue
$\lambda_i\le a$ and $p<q$, $\lambda_i^q=\lambda_i^{q-p}\lambda_i^p\le
a^{q-p}\lambda_i^p$, and summing over these eigenvalues gives
$\sum_{\lambda_i\le a}\lambda_i^q\le a^{q-p}H_p$. It shows that the
small eigenvalues of $A$ carry little $q$th power, so it excludes any
mass of the trace hidden in many vanishing eigenvalues, and it allows
fractional powers to pass to the limit along every critical subsequence
of trace observables.
\end{proof}

\begin{theorem}[Strict homogeneous Schatten decay]
\label{thm:t22-near-Schatten-decay}
For clocks of Lemma~\ref{lem:t22-clock-fractional-entrance},
at every time $s>0$ after their last injection, every critical limit
satisfies
\begin{equation}
 \mathbb E H_p(s+t)\le
 e^{-D p(1-p)t/2}\mathbb E H_p(s),\qquad t\ge0.
                                                               \label{eq:clock-homogeneous-decay}
\end{equation}
For clocks of the preceding lemma, the finite fractional entrance at
one unit of travel and \eqref{eq:small-clock-channels} therefore imply the same exponential
bound for actual critical limsups at all later travel times, with
constants independent of the total length.
\end{theorem}

After injection ends, the clock evolves homogeneously. The negative fractional drift gives decay on the physical length scale used in the Poisson argument.

\begin{proof}
The projected generator \eqref{eq:word-generator}--\eqref{eq:word-covariance} has mean $-DgA/2$ and covariance
$D(I\otimes S+S\otimes I)/2$, $S=I-T$. In an eigenbasis of $A$,
$\mathbb E|dA_{ij}|^2=(D/2)\lambda_i\lambda_j(S_{ii}+S_{jj})dt$.
For a smooth concave scalar function $f$ with $f(0)=0$, its trace
generator is consequently
\begin{align}
 \mathcal G\operatorname{Tr}f(A)
 ={}&-\frac{D}{2}g\sum_i\lambda_i f'(\lambda_i)
       +\frac{D}{2}\sum_i\lambda_i^2S_{ii}f''(\lambda_i)\notag\\
 &+\frac{D}{2}\sum_{i<j}\lambda_i\lambda_j(S_{ii}+S_{jj})
       \frac{f'(\lambda_i)-f'(\lambda_j)}{\lambda_i-\lambda_j}.
                                                               \label{eq:t22-Schatten-generator}
\end{align}
The last summands are nonpositive. This follows by pairing the projected covariance with the first
and second differentials of the trace.
For the regularized functions used next it follows by polynomial
approximation through three derivatives on a trace stop; take the
approximating polynomial to vanish at zero. The trace error is then
bounded by its derivative error times $\gamma$, with no rank factor.
The fixed clock/current moments remove the stops.

Formally taking $f(x)=x^p$ gives
\[
 \mathcal GH_p\le-\frac{Dp}{2}\sum_i\lambda_i^p
                  [g+(1-p)S_{ii}]
             \le-\frac{Dp(1-p)}{2}H_p,
\]
since $T_{ii}\le g$. To justify the singularity, choose $0<p_0<p$
and use $f_\epsilon(x)=(x+\epsilon)^p-\epsilon^p$.
The first two terms of \eqref{eq:t22-Schatten-generator} are at most
\[
 -\frac{Dp(1-p)}{2}\sum_i\lambda_i^2
                                  (\lambda_i+\epsilon)^{p-2}.
\]
Indeed, with $u=\lambda_i/(\lambda_i+\epsilon)$,
$gu+(1-p)u^2S_{ii}\ge(1-p)u^2$. Moreover, for all $x\ge0$,
\begin{equation}
 f_\epsilon(x)-x^2(x+\epsilon)^{p-2}
                 \le C_{p,p_0}\epsilon^{p-p_0}x^{p_0}.
                                                               \label{eq:clock-power-regularization}
\end{equation}
Separate $x\le\epsilon$ and $x\ge\epsilon$ to verify this scalar bound.
The entrance lemma bounds $\mathbb E H_{p_0}$ on each fixed later
interval. Integrating the regularized generator, then taking
$\epsilon\downarrow0$, removes its error by \eqref{eq:clock-power-regularization}.
Monotone convergence passes $f_\epsilon$ to $x^p$, including the
integrated negative term; \eqref{eq:small-clock-channels} removes any zero-eigenvalue dust.
Gronwall proves \eqref{eq:clock-homogeneous-decay}. The statement is deliberately restricted to
the entrance class, which supplies the smaller-power bound.
No assertion for an arbitrary abstract clock with only finite $H_p$
is needed.

The constants at one unit of travel are independent of the earlier
middle and source location. Applying the exact limiting decay on
each subsequent finite interval therefore leaves constants independent
of its final length. The finite-$W$ smaller-power and higher-trace
bounds pass $H_p$ expectations to critical limits by first cutting
off small and large eigenvalues. This proves the stated actual limsup
version without interchanging $W$ with an unbounded time.
\end{proof}

\begin{theorem}[Complete remote-clock Schatten tail]
\label{thm:t22-full-remote-Schatten}
Let $Q_C=W^{-2}U^*P_CU$ be the complete two-input matched-lead clock
of a union $C$ of physical intervals, with $0\le P_C\le I$, all at
scaled distance at least $R\ge2$ from both endpoints of a path of
length $\tau\ge2R$. Then, for every fixed $0<p<1$,
\begin{equation}
 \limsup_{W\to\infty,\ n/W\to\tau}
                \mathbb E\operatorname{Tr}Q_C^p
                    \le C_p e^{-c_pR}.
                                                               \label{eq:remote-clock-tail}
\end{equation}
There are no hidden small clock channels in this estimate.
In particular, for each fixed $a\ge0$,
\begin{equation}
 \limsup_W\mathbb E\log\det(I+aQ_C)
                  \le C_p a^p e^{-c_pR}.             \label{eq:remote-clock-log-bound}
\end{equation}
\end{theorem}

The estimate includes both incoming ports of the actual path. It bounds the error made by changing a boundary far from the supported source.

\begin{proof}
For a unit source interval, the near-near block is covered by the
preceding two results, with travel equal to its distance from the near
end. Reverse the actual independent path to treat the far-far block.
For the complete positive clock, pinching and trace concavity give
\[
 \operatorname{Tr}Q^p
       \le\operatorname{Tr}Q_{\rm near}^p
                            +\operatorname{Tr}Q_{\rm far}^p.
\]
Thus the full interference is retained before this inequality.
Partition $C$ into intervals of length at most one. For positive
matrices and $0<p<1$, $\operatorname{Tr}(A+B)^p\le
\operatorname{Tr}A^p+\operatorname{Tr}B^p$.
One direct proof uses the integral representation of
$\operatorname{Tr}A^p$ and
$(A+B+t)^{-1}\preceq(A+t)^{-1},(B+t)^{-1}$ after splitting
$\operatorname{Tr}(A+B)(A+B+t)^{-1}$ into its two positive trace terms.
The two sums of exponential travel bounds are at most
$C_pe^{-c_pR}$ independently of the number of intervals. This proves
\eqref{eq:remote-clock-tail}. Finally $\log(1+x)\le C_px^p$ for $x\ge0$ proves \eqref{eq:remote-clock-log-bound}.

The estimate retains both boundary inputs. The weighted collar
bound will convert it into a comparison of closed spectral determinants.
\end{proof}

\subsection{Higher fractional moments}\label{subsec:s29-higher-fractional-moments}

 The closing argument also needs products of fractional boundary
errors. Theorem~\ref{thm:t22-remote-higher} bounds these higher moments
when their total degree is below one. We apply the same negative
curvature to a power of the fractional trace, control its entrance,
and sum the contributions from source intervals at increasing
distances.

 \begin{theorem}[Fractional trace moments below total degree one]
\label{thm:t22-remote-higher}
If $0<p<1$, $r\ge1$ and $pr<1$, then for every source $0\le q\le1$
at distance at least $R\ge2$ from both path endpoints,
\begin{equation}
 \limsup_W\mathbb E\{\operatorname{Tr}Q(q)^p\}^{r}
                  \le C_{p,r}e^{-c_{p,r}R}.                 \label{eq:remote-clock-higher}
\end{equation}
The constants are independent of the fixed total critical length.
For a near clock propagated homogeneously after its positive-time
entrance, the limiting bound is more precisely
\begin{equation}
 \mathbb E H_p(s+t)^r
 \le e^{-D_0rp(1-rp)t/2}\mathbb E H_p(s)^r.                \label{eq:clock-higher-decay}
\end{equation}
These assertions are for actual critical limits and limsups, not a
uniform prelimit theorem on lengths increasing with $W$. Proposition~\ref{prop:t22-remote-higher-window} gives the bound
\eqref{eq:remote-clock-higher}, with some rate $c>0$ in place of
$c_{p,r}$, at energies within $O(W^{-2})$ of $E$ and lengths
$n_W/W\to\tau$; it does not assert \eqref{eq:clock-higher-decay} there.
\end{theorem}

The restriction on the total fractional degree pays products of boundary errors. No higher-degree inverse or trace estimate is inferred.

\begin{proof}
Here are the additional finite-width entrance details. For positive $A$
write $H=\operatorname{Tr}A^p$ and $\gamma=\operatorname{Tr}A$.
In an eigenbasis of $A$, the trace Hessian and Cauchy--Schwarz give
\[
 (DH[A_1])^2\le\frac{p}{1-p}H\,(-D^2H[A_1,A_1]).
\]
Indeed the first differential uses only diagonal entries of $A_1$;
the diagonal part of the negative Hessian is
$p(1-p)\sum_i\lambda_i^{p-2}|(A_1)_{ii}|^2$, and its other terms
are nonnegative. Consequently
\begin{equation}
 D^2(H^r)[A_1,A_1]
 \le\frac{r(1-pr)}{1-p}H^{r-1}D^2H[A_1,A_1]\le0.          \label{eq:fractional-trace-concavity}
\end{equation}
No assumption on the smallest eigenvalue is used after positive
regularization.

On the preceding fixed small scalar/jet and $g/W$ domain, its passive mean
also satisfies $\mathbb E_jLL^*\succeq I/2$. Applying Lemma~\ref{lem:t22-one-port-inputs}(a), retain the same event
$\|LL^*\|\le K$. On the first half of the segment from the conditional
mean to the new matrix, its clock lies between $A/4$ and $K'A$.
Thus \eqref{eq:fractional-trace-concavity} and \eqref{eq:clock-covariance-gap} retain a Jensen
deficit at least
\[
 \frac{c_{p,r}}{W\gamma^p}H^{r-1}
       (H-C_pg^{1-p}\gamma^p)_+^2.
\]
Using $(x-a)_+^2x^{r-1}\ge c_rx^{r+1}-C_ra^{r+1}$ yields
\begin{equation}
 \mathbb E_jH_+^r\le(1+C\delta_0+CW^{-2})H^r
 -\frac{c}W\frac{H^{r+1}}{\gamma^p}
 +\frac{C}W g^{(r+1)(1-p)}\gamma^{pr}.                    \label{eq:higher-clock-logistic}
\end{equation}
The zero clock is interpreted by continuity.

Use the same linear scalar weight as above.  Lemma~\ref{lem:t22-one-port-inputs}(c)
applies with $C_*$ the constant of the first term of
\eqref{eq:higher-clock-logistic}. Its correlation costs
$CH^r/W$, since $H_+^r\le\|L\|^{2pr}H^r$.
All exceptional sets are paid with arbitrary fixed polynomial powers
and $H^r\le W^{r(1-p)}\gamma^{pr}$. The fixed finite-edge moments
of Lemma~\ref{lem:edge-clock-moments} bound the forcing in \eqref{eq:higher-clock-logistic} during each entrance interval,
uniformly in the preceding middle length. H\"older gives
\[
 \mathbb E(H^{r+1}/\gamma^p)
 \ge\frac{(\mathbb EH^r)^{(r+1)/r}}{(\mathbb E\gamma^{pr})^{1/r}}.
\]
The weighted expectation therefore obeys a discrete comparison
$y_+\le y+W^{-1}(Cy-cy^{1+1/r}+C)+o(W^{-2})$.
It comes down from its initial polynomial bound after every positive
travel time. As in Lemma~\ref{lem:t22-clock-fractional-entrance}, start the scalar argument after half that time;
no conductance entrance assumption at time zero is made.

In the limiting generator, Theorem~\ref{thm:t22-near-Schatten-decay} gives
$\mathcal GH_p\le-D_0p(1-p)H_p/2$. Its bracket is exactly
\[
 \Gamma(H_p,H_p)=D_0p^2\operatorname{Tr}(A^{2p}S)
                   \le D_0p^2H_p^2,\qquad S=I-T.
\]
The chain rule proves
$\mathcal G H_p^r\le-D_0rp(1-rp)H_p^r/2$.
For singular powers, use the smaller-power entrance and a slightly
larger moment $r'>r$ with $pr'<1$. They pay respectively the small
eigenvalues and uniform integrability. the smooth regularization \eqref{eq:clock-power-regularization},
then trace stopping, passes this chain rule and Gronwall to the limit.
This proves \eqref{eq:clock-higher-decay}.

Reverse the path for its far clock, pinch the complete positive clock,
and partition the source into intervals of length at most one.
The deterministic inequality $H_p(Q_1+Q_2)\le H_p(Q_1)+H_p(Q_2)$
and Minkowski in $L^r$ sum the exponentially decreasing travel bounds.
The sums are geometric independently of the total length.
This proves \eqref{eq:remote-clock-higher}, retaining the full two-input interference.
\end{proof}

We
close this section by putting its results in the exact form in which the
Poisson endpoint uses them. The theorems above show that an absorber deep
inside a long path is barely felt at its ends: the fractional trace of its
clock decays exponentially in the distance to the boundary. They were
proved at the fixed energy $E$, at the lengths $\lfloor Wt\rfloor$, and for
general sources $0\le q\le1$. The comparison of the ring with smaller
rings in Section~\ref{sec:t22-positive-log-poisson} needs three changes.
First, the spectral statistics there are evaluated at energies that move
with $W$, at distance of order $W^{-2}$ from $E$, which is the size of an
eigenvalue spacing of the ring. Second, the pieces of the ring have
lengths that are only asymptotically proportional to $W$, not exactly
equal to $\lfloor Wt\rfloor$. Third, the absorber sits on a union of
blocks, so its source is an indicator. The next proposition shows that all
the bounds of this section survive these changes, with the same kind of
constants. The idea is that a shift of the energy by $O(W^{-2})$ can be
undone by a small deterministic change of the Gaussian path, whose density
is bounded, so the bounds can be transferred from $E$ to the moving
energy.
 
\begin{proposition}[Window-uniform fractional clock bounds]
\label{prop:t22-remote-higher-window}
Fix $b$ and $E$ as in this section, and a window constant $K\ge0$. We call
deterministic energies $e_W$ \emph{window energies} if
$|e_W-E|\le KW^{-2}$ for all large $W$. At such an energy the whole path is
evaluated at $e_W$: its physical diagonals are $A-e_W$, both matched leads
carry the parameter $z_\partial(e_W)$ of \eqref{eq:matched-parameters}, and
the Gaussian law is unchanged. In (a) and (b) the clock is evaluated at
arbitrary lengths $m_W$ with $m_W/W\to t_0+s$, where $t_0$ is the scaled
position of its last injection and $s$ the travel time. In each assertion
the constants are chosen after $b$, $E$, $K$ and the exponents, and, in
(a), after $\epsilon$ and $s_1\le s_2$; they are chosen before the source,
the length of the preceding middle part, the travel time, the total
length, the energies and the lengths.
\begin{enumerate}
\item[(a)]  For $0<p<1$ and $\epsilon>0$, the clocks of
 Lemma~\ref{lem:t22-clock-fractional-entrance} satisfy
 $\limsup_W\mathbb EH_p(\epsilon)\le C_{p,\epsilon}$ along window energies
 and such lengths. Moreover, for $0<s_1\le s_2$ there are $C$ and, for each
 source, $W_0$ such that $\mathbb EH_p(s)\le C$ at the lengths
 $\lfloor W(t_0+s)\rfloor$, for $W\ge W_0$, $s\in[s_1,s_2]$ and every energy
 $e'$ with $|e'-E|\le KW^{-2}$.
\item[(b)]  For $0<p<1$, $r\ge1$ and $pr<1$ there are $C$ and
 $c>0$ such that these clocks satisfy
 \[
  \limsup_W\mathbb EH_p(s)^r\le Ce^{-cs}
 \]
 for every travel time $s\ge1$, along window energies and such lengths.
\item[(c)]  For $0<p<1$, $r\ge1$ and $pr<1$ there are
 $C$ and $c>0$ such that the following holds for every $\tau>0$, every
 $R\ge2$, all window energies $e_W$, all lengths $n_W$ with
 $n_W/W\to\tau$, and every source $0\le q\le1$ that vanishes within scaled
 distance $R$ of both endpoints:
 \begin{equation}\label{eq:remote-clock-higher-window}
  \limsup_W\mathbb E\{\operatorname{Tr}Q(q)^p\}^r\le Ce^{-cR},
 \end{equation}
 where $Q(q)$ is the complete two-input clock at energy $e_W$ and length
 $n_W$. Since $C$ and $c$ do not depend on the energies, the same bound
 holds for
 $\limsup_W\sup_{|e-E|\le KW^{-2}}\mathbb E\{\operatorname{Tr}Q(q)^p\}^r$.
\end{enumerate}
In (c) the source may be the indicator of any union of blocks whose
distances $R_1,R_2$ to the two ends are both at least $R$.
\end{proposition}
 
\begin{figure}[ht]
\centering
\begin{tikzpicture}[x=.85cm,y=1cm,>=Stealth]
 \begin{scope}[shift={(0,0)}]
  \node[font=\small,anchor=west] at (-.2,2.05) {path $\pi$ at energy $e$};
  \draw[thick,dashed] (0,.4) rectangle (.8,1.4);
  \node[font=\small] at (.4,.9) {$z'$};
  \foreach \x in {1.8,2.8,3.8,4.8}
    \draw[thick,fill=black!8] (\x,.4) rectangle ++(.8,1);
  \node[font=\scriptsize,black,align=center] at (3.7,-.05) {diagonals $A-e$};
  \draw[thick,dashed] (6.6,.4) rectangle (7.4,1.4);
  \node[font=\small] at (7.0,.9) {$z'$};
  \draw[thick,black] (.8,.9)--(1.8,.9) node[midway,above,font=\scriptsize] {$B_\partial$};
  \draw[thick,black] (5.6,.9)--(6.6,.9) node[midway,above,font=\scriptsize] {$B_\partial$};
 \end{scope}
 \draw[<->,thick,red!70!black] (7.9,.9)--(9.1,.9) node[midway,above,font=\small] {$T$};
 \begin{scope}[shift={(9.6,0)}]
  \node[font=\small,anchor=west] at (-.2,2.05) {path $T\pi$ at energy $E$};
  \draw[thick,dashed] (0,.4) rectangle (.8,1.4);
  \node[font=\small] at (.4,.9) {$z$};
  \foreach \x in {1.8,2.8,3.8,4.8}
    \draw[thick,fill=red!10] (\x,.4) rectangle ++(.8,1);
  \node[font=\scriptsize,red!70!black,align=center] at (3.7,-.05) {diagonals shifted by\\ $-\delta I+\alpha B_\partial B_\partial^*$};
  \draw[thick,dashed] (6.6,.4) rectangle (7.4,1.4);
  \node[font=\small] at (7.0,.9) {$z$};
  \draw[thick,red!70!black] (.8,.9)--(1.8,.9) node[midway,above,font=\scriptsize] {$\lambda B_\partial$};
  \draw[thick,red!70!black] (5.6,.9)--(6.6,.9) node[midway,above,font=\scriptsize] {$\lambda B_\partial$};
 \end{scope}
\end{tikzpicture}
\caption{The change of energy in the proof of
Proposition~\ref{prop:t22-remote-higher-window}. Left: a path sample $\pi$
evaluated at the window energy $e=E+\delta$, with both leads at the
parameter $z'=z_\partial(e)$. Right: the transformed sample $T\pi$,
evaluated at the fixed energy $E$ with leads at $z=z_\partial(E)$. It has
the same internal bonds, its physical diagonals are shifted by the small
Hermitian matrix $-\delta I+\alpha B_\partial B_\partial^*$, and its lead
bonds are multiplied by $\lambda$ (red). After the leads are eliminated,
the two interior matrices coincide, so every clock of $\pi$ at $e$ equals
$|\varpi|^2$ times the corresponding clock of $T\pi$ at $E$. Because
$\delta$, $\alpha$ and $\lambda^2-1$ are all of order $W^{-2}$, the law of
$T\pi$ has a density with bounded moments relative to the law of $\pi$,
and the bounds at $E$ transfer to $e$.}
\label{fig:window-energy-map}
\end{figure}
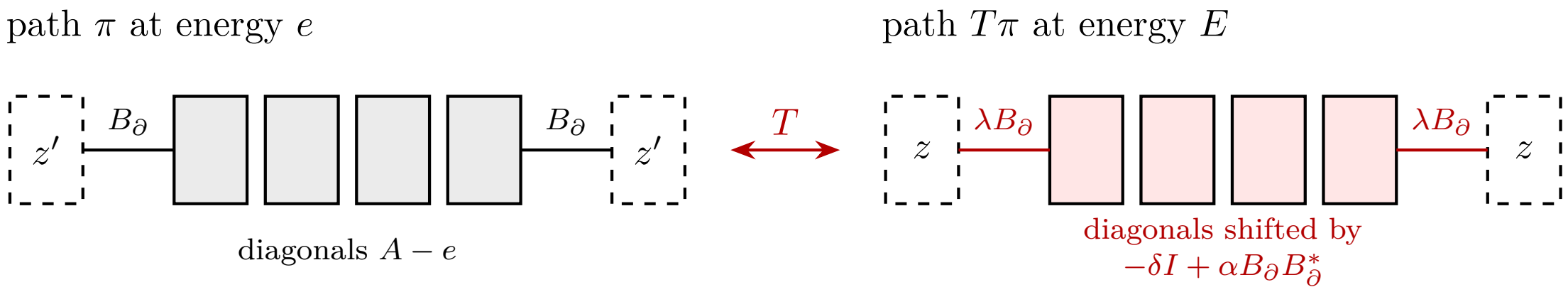
 
\begin{proof}
\emph{Plan of the proof.} The changes of length and of source are mild:
the finite-width bounds of this section are uniform over travel
intervals, and lengths that differ by $o(W)$ give the same critical limits
(Step~1). The change of energy is handled by an exact identity. A shift of
the energy by $\delta=O(W^{-2})$ is the same as a small deterministic
change of the path at the fixed energy $E$, see
Figure~\ref{fig:window-energy-map} (Step~2). This change has a density with
bounded moments (Step~3), so H\"older's inequality transfers every
fractional moment from $E$ to the window energy, at the cost of a
slightly larger exponent (Step~4). Step~5 applies this transfer to prove
(a), (b) and (c).
 
\emph{Step 1: lengths.} At the energy $E$, every finite-width bound in the
proofs of Lemma~\ref{lem:t22-clock-fractional-entrance} and
Theorem~\ref{thm:t22-remote-higher} is uniform over a travel interval. By
Theorem~\ref{thm:full-source-generator}, the trace words are tight with
continuous limits. Hence, since $|m_W-\lfloor Wt\rfloor|=o(W)$, their
values at the lengths $m_W$ and $\lfloor Wt\rfloor$ have the same critical
limits along every subsequence. The passage to critical limits in those
proofs, with its cutoffs by \eqref{eq:small-clock-channels} and by the
higher trace moments, therefore applies at the lengths $m_W$. So all the
assertions hold at $e_W=E$ along such lengths, with the same constants.
 
\emph{Step 2: the change of energy as a change of the path.} Let
$|e-E|\le KW^{-2}$, and put
\[
 \delta=e-E,\qquad z=z_\partial(E),\qquad z'=z_\partial(e).
\]
Let $B_\partial$ be the coupling of the path to its two leads; the matrix
$B_\partial B_\partial^*$ acts only on the end blocks. Since $\Im(1/z)$ and
$\Im(1/z')$ are negative, we may put
\[
 \lambda=\Bigl(\frac{\Im(1/z')}{\Im(1/z)}\Bigr)^{1/2},\qquad
 \alpha=\frac{\lambda^2}{z}-\frac1{z'} .
\]
The number $\alpha$ is real, because $\Im(\lambda^2/z)=\lambda^2\Im(1/z)=
\Im(1/z')$. For a path sample $\pi$, let $T\pi$ be the sample with the same
internal bonds, with physical diagonals shifted by
$-\delta I+\alpha B_\partial B_\partial^*$, and with lead bonds multiplied by
$\lambda$. After the leads are eliminated, the interior matrix of $T\pi$ at
energy $E$ and lead parameter $z$ is
\[
 H_{\rm path}-\delta+\alpha B_\partial B_\partial^*-E
   -\frac{\lambda^2}{z}B_\partial B_\partial^*
 =H_{\rm path}-e-\frac1{z'}B_\partial B_\partial^*,
\]
which is exactly the interior matrix of $\pi$ at energy $e$ and lead
parameter $z'$. Hence the maps from boundary amplitudes to the interior
satisfy
\[
 \mathcal U_e(\pi)=\varpi\,\mathcal U_E(T\pi),\qquad \varpi=\frac{z}{z'\lambda},
\]
and every clock of $\pi$ at energy $e$ is $|\varpi|^2$ times the
corresponding clock of $T\pi$ at $E$. Moreover, since $z_\partial$ is
smooth in the bulk and $|\delta|\le KW^{-2}$,
\[
 |\lambda^2-1|+|\alpha|+\bigl||\varpi|^2-1\bigr|\le C_KW^{-2}.
\]
 
\emph{Step 3: the density of the change.} Under $T$, the lead bonds, of
law $\operatorname{Gin}(b/W)$, are rescaled by $\lambda$. Given them, the
physical diagonals, whose Gaussian entries have variance $(1-2b)/W$,
receive the deterministic Hermitian shift
$-\delta I+\alpha B_\partial B_\partial^*$. Let
\[
 \mathcal E=\{\norm{B_\partial}_{\HS}^2\le8W\};
\]
its complement has probability at most $e^{-c_1W^2}$. On $\mathcal E$, for a
path of $m\le TW$ blocks, the density of the law of $T\pi$ with respect to
the law of $\pi$ has moments of each fixed order $a'$ bounded by
\[
 \exp\bigl\{C_{a'}\bigl(mW^2\delta^2
   +\alpha^2W\norm{B_\partial B_\partial^*}_{\HS}^2\bigr)/(1-2b)\bigr\}
 \bigl(1+C_{a'}|\lambda^2-1|^2\bigr)^{2W^2}\le C_{K,a'},
\]
for $W\ge W_0(T)$. The first factor comes from the Gaussian shift of the
diagonals and the second from the rescaling of the $2W^2$ entries of the
lead bonds. The bound holds because
\[
 mW^2\delta^2\le\frac{TK^2}{W},\qquad
 \alpha^2W\norm{B_\partial B_\partial^*}_{\HS}^2\le\frac{64C_K^2}W,
 \qquad |\lambda^2-1|^2W^2\le\frac{C_K^2}{W^2} .
\]
 
\emph{Step 4: the transfer inequality.} Choose $1<a<1/(pr)$ and $a'$ with
$1/a+1/a'=1$. By Step~2, $H_p^r$ at energy $e$ equals $|\varpi|^{2pr}$ times
$H_p^r$ of $T\pi$ at $E$. H\"older's inequality on $\mathcal E$, with the
density bound of Step~3, and the Cauchy--Schwarz inequality off
$\mathcal E$ give
\[
 \mathbb E_eH_p^r\le C_{K,a'}|\varpi|^{2pr}
   \bigl(\mathbb E_EH_p^{ra}\bigr)^{1/a}
   +\bigl(\mathbb E_eH_p^{2r}\bigr)^{1/2}e^{-c_1W^2/2},
\]
where $\mathbb E_e$ refers to the path at energy $e$. The second term is
negligible. Indeed $H_p\le W^{1-p}\gamma^p$, and $\gamma$ is at most the
normalized trace of the complete clock, so \eqref{eq:full-clock-moments},
whose constants are uniform for energies in a compact bulk set
(Lemma~\ref{lem:source-recorded-estimates}(e)), bounds $\mathbb E_eH_p^{2r}$
by a power of $W$. Thus there is a constant $C_1$, depending only on $b$,
$E$, $K$, $p$, $r$ and $a$, and chosen before $T$, such that
\begin{equation}\label{eq:window-transfer}
 \mathbb E_eH_p^r\le C_1\bigl(\mathbb E_EH_p^{ra}\bigr)^{1/a}+W^{-1}
\end{equation}
for $W\ge W_0(T)$, every $|e-E|\le KW^{-2}$, every length at most $TW$ and
every source $0\le q\le1$. The price of moving the energy is only the
slightly larger exponent $ra$, which still satisfies $pra<1$.
 
\emph{Step 5: conclusion.} For (a), we apply \eqref{eq:window-transfer}
with $r=1$. At the energy $E$ and the exponent $a$, with $pa<1$, the
finite-width higher recursion in the proof of
Theorem~\ref{thm:t22-remote-higher} comes down from its initial polynomial
bound after every positive travel time, and it stays bounded on later
fixed intervals, also at the lengths $m_W$ by Step~1.
 
For (b), the proof of Theorem~\ref{thm:t22-remote-higher} gives, from
\eqref{eq:clock-higher-decay} and the entrance bound (a),
\[
 \limsup_W\mathbb E_EH_p(s)^{ra}\le C'e^{-c's}
\]
at the lengths $m_W$. Then \eqref{eq:window-transfer} gives (b), with
$C=C_1(C')^{1/a}$ and $c=c'/a$.
 
For (c), we follow the proof of Theorem~\ref{thm:t22-remote-higher}: the
reversal of the path, the pinching inequality, the partition of the source
into intervals of length at most one, and Minkowski's inequality. These
steps are deterministic at each width, so they apply at the energy $e_W$,
and the reversed path has the same law, energy and leads. An interval at
scaled distance $R+k$ from an end contributes at most $Ce^{-c(R+k)}$, by
(b), and the sum over $k$ is geometric, independently of $\tau$. This also
covers indicator sources of any union of blocks with $R_1,R_2\ge R$.
Finally, for the supremum over energies, we choose at each $W$ an energy
$e$ with $|e-E|\le KW^{-2}$ at which the expectation is nearly maximal. The
resulting sequence consists of window energies, so the bound just proved
applies to it.
\end{proof}
 
\section{Positive log absorption and the original-ring Poisson endpoint}
\label{sec:t22-positive-log-poisson}

In this section we prove the Poisson endpoint of Theorem~\ref{thm:endpoints}.
We compare the original ring with a family of independent smaller rings
of the same total dimension, and we show that the comparison error tends
to zero when the rings are long. The comparison is carried out on
determinants with a finite positive absorber. We never use a complex
source disk or a pointwise bound on a differentiated clock. All Gaussian
variables have the original law with fixed $b$. Limits in $W$ are taken
at fixed $\tau$, and only afterwards do we let $\tau\to\infty$.

Throughout we write

\[
 \Lambda=N\rho(E)
\]

for the physical unfolding factor, so that $u_\alpha=\Lambda(\lambda_\alpha-E)$.

\begin{theorem}[Original-ring Poisson endpoint]
\label{thm:t22-Poisson-endpoint}
For every nonnegative smooth compactly supported $f$,
\begin{equation}
 \lim_{\tau\to\infty}\ \limsup_{W\to\infty,\ n/W\to\tau}
 \left|\mathbb E\exp\{-\operatorname{Tr}f(\Lambda(H-E))\}
       -\exp\{-\int_{\mathbb R}(1-e^{-f(x)})\,dx\}\right|=0.
                                                               \label{eq:poisson-laplace-limit}
\end{equation}
Thus $\mathcal X_{\tau,b,E}$ tends to the unit-intensity Poisson
process as $\tau\to\infty$, after the critical limit, for every
fixed $0<b<1/2$ and bulk $E$.
\end{theorem}

The right side of \eqref{eq:poisson-laplace-limit} is the Laplace
functional of the Poisson process of intensity one. Since Laplace
functionals determine the law of a point process, the theorem identifies
the limit completely, including the independence of counts in disjoint
intervals.

\paragraph{Strategy.} We prove the theorem in four steps, which
correspond to the four subsections below.
\begin{enumerate}
\item We bound fractional moments of boundary inverses of a collar by a
 phase-averaging argument (Lemma~\ref{lem:t22-phase-fractional-inverse}).
\item We use these bounds to compare the real logarithms of two closing
 determinants, and we apply the comparison when the change is caused by
 an absorber far from the collar
 (Lemma~\ref{lem:t22-unitary-log-comparison} and
 \eqref{eq:single-collar-log-error}).
\item We show that a smoothed positive logarithm of the determinant
 recovers the spectral statistic $\Tr f(\Lambda(H-E))$
 (Lemma~\ref{lem:t22-finite-log-recovery}).
\item We rewire the ring into independent rings of the same total
 volume, sum the collar errors, and compute the Laplace functional of
 the independent rings (proof of Theorem~\ref{thm:t22-Poisson-endpoint}).
\end{enumerate}

\subsection{A fractional boundary inverse from rank-one phase tests}\label{subsec:s30-a-fractional-boundary-inverse-from-rank-}

We  bound the potentially large inverse that appears when one collar
closes the path. Lemma~\ref{lem:t22-phase-fractional-inverse} gives a
projected fractional inverse estimate, which yields
Lemma~\ref{lem:t22-unitary-log-comparison} for real log-determinants.
We first compress to the actual ports, use rank-one phase tests, and
then integrate a controlled change of boundary condition.

We  begin with the closing determinant of one collar. Let a unitary $V$ act
on a finite port space, and let $P$ be the projection onto the two
$W$-dimensional ports of the collar. We put $\widehat U=U\oplus0$, where
$U$ is the collar response on $\operatorname{Ran}P$ and $0$ acts on the
complement. Writing $V$ in the decomposition
$\operatorname{Ran}P\oplus\operatorname{Ran}P^\perp$ as
\[
 V=\begin{pmatrix}r&B\\ *&*\end{pmatrix},
 \qquad r=PVP|_{\operatorname{Ran}P},
\]
we have $\widehat UV=\begin{pmatrix}Ur&UB\\0&0\end{pmatrix}$, and hence
\begin{equation}\label{eq:collar-determinant-compression}
 \det(I-\widehat UV)=\det(I-Ur).
\end{equation}
The complement of $\operatorname{Ran}P$ may contain arbitrarily many
unused lead or dilation coordinates; by
\eqref{eq:collar-determinant-compression} they do not change the
determinant.

\begin{lemma}[Projected fractional inverse]
\label{lem:t22-phase-fractional-inverse}
Suppose the scalar condition of Theorem~\ref{thm:t22-collar-kernel}
holds for all rank-one twists $V\exp(i\theta vv^*)$, with a common
kernel bound $M I$. Then for $0<p<1$ and every unit vector $v$,
\begin{equation}
 \mathbb E_U\left|v^*(I-\widehat UV)^{-1}v\right|^p
                         \le C_p M.                     \label{eq:phase-fractional-inverse}
\end{equation}
The expectation is over the actual fresh collar. Boundary inverses
are obtained by radial limits, not by replacing an expected boundary
Poisson measure with its almost-sure zero value.
\end{lemma}

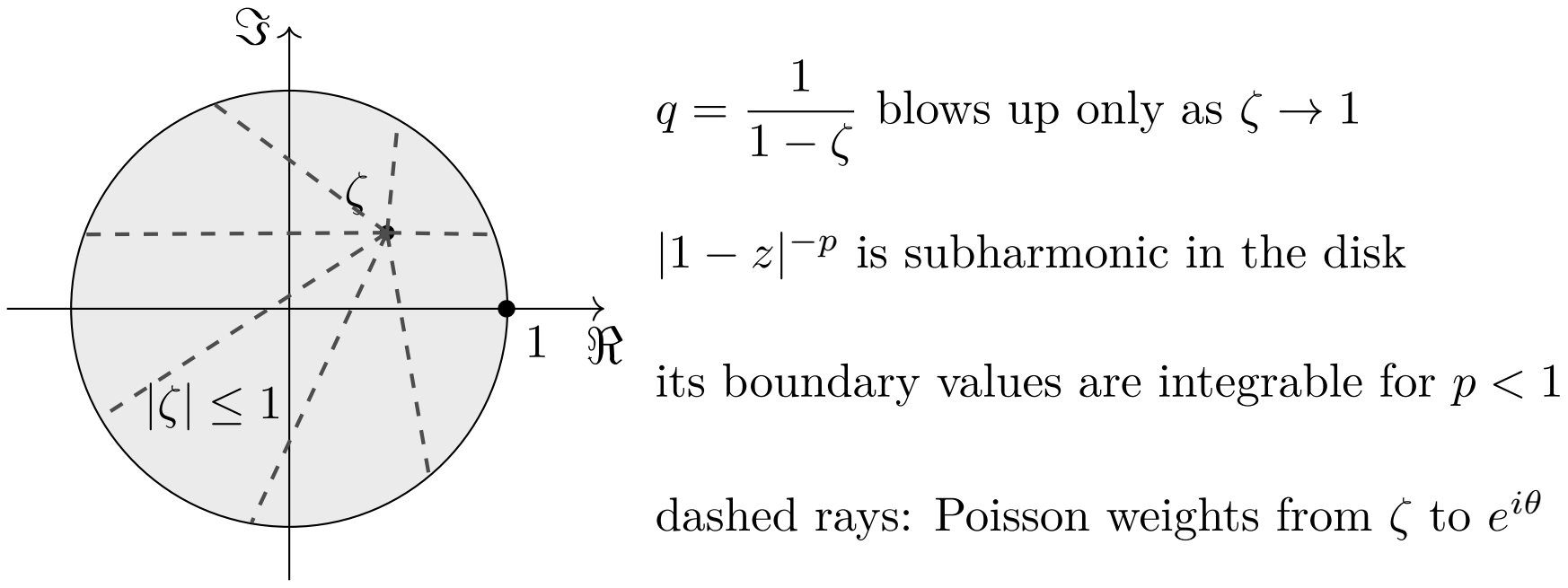
\begin{figure}[ht]
\centering
\begin{tikzpicture}[x=1.6cm,y=1.6cm]
 \draw[thick] (0,0) circle (1);
 \fill[black!8] (0,0) circle (1);
 \draw[->] (-1.3,0)--(1.45,0) node[below] {$\Re$};
 \draw[->] (0,-1.25)--(0,1.3) node[left] {$\Im$};
 \fill (.45,.35) circle (.035) node[above left,font=\small] {$\zeta$};
 \fill (1,0) circle (.04);
 \node[below right,font=\small] at (1,0) {$1$};
 \foreach \a in {20,60,110,160,210,260,310}{
  \draw[black!70,dashed,thick] (.45,.35)--({cos(\a)},{sin(\a)});
 }
 \node[font=\small] at (-.35,-.45) {$|\zeta|\le1$};
 \node[font=\small,anchor=west] at (1.6,.9) {$q=\dfrac1{1-\zeta}$ blows up only as $\zeta\to1$};
 \node[font=\small,anchor=west] at (1.6,.25) {$|1-z|^{-p}$ is subharmonic in the disk};
 \node[font=\small,anchor=west] at (1.6,-.35) {its boundary values are integrable for $p<1$};
 \node[font=\small,anchor=west] at (1.6,-.95) {dashed rays: Poisson weights from $\zeta$ to $e^{i\theta}$};
\end{tikzpicture}
\caption{The scalar step in the proof of
Lemma~\ref{lem:t22-phase-fractional-inverse}. A rank-one twist of angle
$\theta$ turns the quantity $q$ into $q_\theta=1/(1-\zeta e^{i\theta})$,
where the point $\zeta$ lies in the closed unit disk (shaded). The
inverse $q$ is large only when $\zeta$ is close to the boundary point $1$.
Since $|1-z|^{-p}$ is subharmonic, its value at $\zeta$ is at most its
Poisson average over the circle (dashed rays). For $p<1$ the boundary
function $|1-e^{i\theta}|^{-p}$ is integrable, and the Poisson weights are
exactly the real parts $2\Re q_\theta-1$, which the collar hypothesis
controls in expectation. In this way a pointwise inverse, which may be
huge, is bounded through an averaged kernel, which is not.}
\label{fig:phase-twist}
\end{figure}

The idea of the proof is that a rank-one twist moves the scalar
$q=v^*(I-\widehat UV)^{-1}v$ along a circle. The average of $|q|^p$ over
the twist can be written through the Poisson kernel of the twisted
matrix, and the hypothesis bounds the expected Poisson kernel.
Figure~\ref{fig:phase-twist} illustrates the scalar part of the argument.

\begin{proof}
We first record that the full Poisson kernel is bounded by its collar
part. In the decomposition above, the Poisson kernel of $\widehat UV$
has the block form
\[
 \begin{pmatrix}
 \mathcal P(Ur)&(I-Ur)^{-1}UB\\
 B^*U^*(I-r^*U^*)^{-1}&I
 \end{pmatrix}.
\]
This matrix is positive semidefinite. Taking the Schur complement of the
lower right identity block, we see that it is bounded above by
\begin{equation}\label{eq:full-kernel-upper}
 2\begin{pmatrix}\mathcal P(Ur)&0\\0&I\end{pmatrix}.
\end{equation}
Hence the expected full kernel is bounded by twice the block-diagonal
matrix formed by the expected collar kernel and the identity,
independently of the dimension of the complement.

\emph{Step 1: a radial regularization and the sign of $\Re q$.} At first we
replace $\widehat U$ by $t\widehat U$ with $0<t<1$, so that
$X=t\widehat UV$ is a strict contraction. Put $q=v^*(I-X)^{-1}v$. We use
the identity
\[
 (I-X)^{-1}=\frac{1}{2}\Bigl[I+(I+X)(I-X)^{-1}\Bigr],
\]
together with the formula for the real part of a Cayley transform,
\begin{equation}\label{eq:cayley-real-part}
 \Re\,(I+X)(I-X)^{-1}=(I-X)^{-*}(I-X^*X)(I-X)^{-1}\succeq0 .
\end{equation}
The right side of \eqref{eq:cayley-real-part} is the Poisson kernel of
$X$. Taking the quadratic form in $v$ gives
\begin{equation}\label{eq:re-q-half}
 \Re q=\frac{1}{2}+\frac{1}{2}\,v^*\Re\bigl[(I+X)(I-X)^{-1}\bigr]v\ \ge\ \frac{1}{2}.
\end{equation}

\emph{Step 2: the rank-one twist moves $q$ on a circle.} Fix $\theta$ and
consider the twisted matrix $X_\theta=X(I+(e^{i\theta}-1)vv^*)$. Since
$I-X_\theta=(I-X)-(e^{i\theta}-1)Xvv^*$, the Sherman--Morrison formula
gives
\[
 q_\theta:=v^*(I-X_\theta)^{-1}v
 =\frac{q}{1-(e^{i\theta}-1)\,v^*(I-X)^{-1}Xv}.
\]
Because $(I-X)^{-1}X=(I-X)^{-1}-I$, the quadratic form in the denominator
equals $q-1$. Simplifying, we obtain
\begin{equation}\label{eq:rank-one-circle}
 q_\theta=\frac{1}{1-\zeta e^{i\theta}},
 \qquad \zeta=1-q^{-1}.
\end{equation}
The condition $\Re q\ge1/2$ in \eqref{eq:re-q-half} is equivalent to
$|q^{-1}-1|\le1$, that is, to $|\zeta|\le1$. So $\zeta$ lies in the closed
unit disk of Figure~\ref{fig:phase-twist}, and $q=q_0=(1-\zeta)^{-1}$.

\emph{Step 3: Poisson weights are kernel entries.} A direct computation
from \eqref{eq:rank-one-circle} gives
\begin{equation}\label{eq:scalar-poisson-weight}
 2\Re q_\theta-1
 =\Re\frac{1+\zeta e^{i\theta}}{1-\zeta e^{i\theta}}
 =\frac{1-|\zeta|^2}{|1-\zeta e^{i\theta}|^2},
\end{equation}
which is the Poisson kernel of the disk at $\zeta$, evaluated at the
boundary point $e^{-i\theta}$. On the other hand, applying
\eqref{eq:re-q-half} to $X_\theta$ shows that $2\Re q_\theta-1$ is the
quadratic form in $v$ of the Poisson kernel of $X_\theta$. By the
hypothesis and \eqref{eq:full-kernel-upper}, its expectation over the
collar is at most $CM$ for every $\theta$. In the kernel argument of this step, we first replace the upper bound $M$ by $\max\{M,1\}$, retaining the notation $M$ for the enlarged bound. This pays $2\|P^\perp v\|^2$ in \eqref{eq:full-kernel-upper}. In the application $M=C(1+G)$, so the stated conclusion is unchanged.

\emph{Step 4: the subharmonic majorant.} The function
$z\mapsto|1-z|^{-p}=\exp\{-p\log|1-z|\}$ is subharmonic in the open disk,
because $\log|1-z|$ is harmonic there. For $p<1$ its boundary values
$|1-e^{i\theta}|^{-p}$ are integrable. Hence, as in
Figure~\ref{fig:phase-twist}, the value at $\zeta$ is bounded by its
Poisson average:
\[
 |q|^p=|1-\zeta|^{-p}
 \le\int_0^{2\pi}\frac{1-|\zeta|^2}{|1-\zeta e^{i\theta}|^2}\,
       |1-e^{-i\theta}|^{-p}\,\frac{d\theta}{2\pi}.
\]
We insert \eqref{eq:scalar-poisson-weight}, take the expectation over the
collar, and use the bound $CM$ from Step~3 for each $\theta$:
\begin{equation}\label{eq:fractional-inverse-radial}
 \E_U|q|^p
 \le\int_0^{2\pi}\E_U\bigl[2\Re q_\theta-1\bigr]\,
       |1-e^{-i\theta}|^{-p}\,\frac{d\theta}{2\pi}
 \le CM\int_0^{2\pi}|1-e^{i\theta}|^{-p}\,\frac{d\theta}{2\pi}
 =C_pM.
\end{equation}

\emph{Step 5: removing the radial factor.} Here we use the smaller scalar margin $d_1=d_*/4$. The estimate $d(tr)\le d(r)+(1-t)/(2\eta)$ keeps the radial regularization inside an admissible collar threshold for $t$ sufficiently close to one. The radial factor $t$ increases
the conductance of the collar by at most $2W(1-t^2)$, so the scalar
condition of Theorem~\ref{thm:t22-collar-kernel} still holds for $t$ close
to one. The bound \eqref{eq:fractional-inverse-radial} is therefore uniform
in such $t$. We let $t\uparrow1$ at fixed $W$; by Fatou's lemma,
\[
 \E_U\bigl|v^*(I-\widehat UV)^{-1}v\bigr|^p
 \le\liminf_{t\uparrow1}\E_U\bigl|v^*(I-t\widehat UV)^{-1}v\bigr|^p
 \le C_pM.
\]
This proves \eqref{eq:phase-fractional-inverse}.
\end{proof}

We now use the lemma to compare two closing determinants. The comparison
concerns the real logarithm of the modulus, so no choice of phase is
needed along the way.

\begin{lemma}[A real log-determinant comparison]
\label{lem:t22-unitary-log-comparison}
Let $V_0,V_1$ be unitary and independent of the collar. Put
\[
 \sigma=\operatorname{Tr}|V_1-V_0|^p,\quad 0<p<1,\qquad
 g_0=\operatorname{Tr}(I-r_0^*r_0).
\]
Assume $d(r_0)\le d_*/8$ in the notation of the collar theorem and
$\sigma\le c_*W$, for a sufficiently small fixed $c_*>0$.
Then
\begin{equation}
 \mathbb E_U\left|
 \log|\det(I-\widehat UV_1)|
 -\log|\det(I-\widehat UV_0)|\right|
               \le C_p(1+g_0+\sigma)\sigma.              \label{eq:unitary-log-comparison}
\end{equation}
The statement concerns the real logarithm of the modulus and has
no branch or winding ambiguity.
\end{lemma}

\begin{proof} We pass from the first unitary matrix to the second
by changing one eigenphase at a time, so that consecutive matrices differ by a
rank-one twist, a change acting in a single direction (Step~1). Along this
interpolation the hypotheses of the fractional inverse bound of
Lemma~\ref{lem:t22-phase-fractional-inverse} remain valid (Step~2). For one
twist, the ratio of the two determinants is an explicit function of a single
quadratic form of the inverse (Step~3), and the logarithm of its modulus is
controlled by fractional moments of that quadratic form, which the previous
lemma bounds (Step~4). Step~5 adds up the contributions of all eigenphases.

\emph{Step 1: an interpolation by rank-one twists.} The matrix $V_0^*V_1$
is unitary. We diagonalize it as
$V_0^*V_1=\sum_j\zeta_jw_jw_j^*$ with $|\zeta_j|=1$ and an orthonormal
basis $(w_j)$. Since $|V_1-V_0|=|V_0^*V_1-I|$, we have
\begin{equation}\label{eq:sigma-phases}
 \sigma=\sum_j|\zeta_j-1|^p.
\end{equation}
We change the eigenphases one at a time. For $k\ge0$ put
\[
 V^{(k)}=V_0\prod_{j\le k}\bigl(I+(\zeta_j-1)w_jw_j^*\bigr),
\]
so that $V^{(0)}=V_0$, the last matrix is $V_1$, and consecutive matrices
differ by one rank-one twist in the direction $w_k$.

\emph{Step 2: the hypotheses hold along the interpolation.} Each $V^{(k)}$,
and each rank-one twist of it, differs from $V_0$ in trace norm by at most
$C_p\sigma+2$. Compression to the collar ports preserves this bound. Since
$(z_b-X)^{-1}=(I-r)/(2i\eta)$, the normalized scalar traces change by at
most $C(\sigma+1)/W$, and the conductance satisfies
$g\le g_0+C(\sigma+1)$. With $\sigma\le c_*W$ the scalar condition of
Theorem~\ref{thm:t22-collar-kernel} holds uniformly, and
Lemma~\ref{lem:t22-phase-fractional-inverse} applies to all these
matrices with
\begin{equation}\label{eq:interpolation-M}
 M=C(1+g_0+\sigma).
\end{equation}

\emph{Step 3: one determinant ratio.} Fix one step, write $V=V^{(k-1)}$,
$\zeta=\zeta_k$, $w=w_k$, and $X=\widehat UV$. Then
$I-\widehat UV^{(k)}=(I-X)-(\zeta-1)Xww^*$. The matrix determinant lemma
and the identity $(I-X)^{-1}X=(I-X)^{-1}-I$ give
\begin{equation}\label{eq:one-phase-ratio}
 z:=\frac{\det(I-\widehat UV^{(k)})}{\det(I-\widehat UV^{(k-1)})}
 =1-(\zeta-1)(q-1),\qquad q=w^*(I-X)^{-1}w.
\end{equation}
Exchanging the roles of the two endpoints, the reciprocal ratio has the
same form with the phase $\bar\zeta$ and the quadratic form $q'$ of the
other endpoint:
\[
 |z^{-1}-1|=|\zeta-1|\,|q'-1| .
\]

\emph{Step 4: the real logarithm.} For every nonzero complex number $z$ we
use the elementary inequality
\begin{equation}\label{eq:log-two-sided}
 \bigl|\log|z|\bigr|\le C_p\bigl(|z-1|^p+|z^{-1}-1|^p\bigr).
\end{equation}
Indeed, if $|z|\ge1$ then $\log|z|\le\log(1+|z-1|)\le C_p|z-1|^p$, and if
$|z|<1$ the same argument applies to $z^{-1}$. Combining
\eqref{eq:one-phase-ratio}--\eqref{eq:log-two-sided} with
Lemma~\ref{lem:t22-phase-fractional-inverse} at both endpoints, we get for
one step
\[
 \E_U\bigl|\log|z|\bigr|
 \le C_p|\zeta-1|^p\bigl(1+\E_U|q|^p+\E_U|q'|^p\bigr)
 \le C_p|\zeta-1|^p\,(1+M).
\]

\emph{Step 5: summation.} The logarithm of the full ratio is the sum of the
logarithms of the one-step ratios. Summing the last bound over $k$ and
using \eqref{eq:sigma-phases} and \eqref{eq:interpolation-M}, we obtain
\[
 \E_U\bigl|\log|\det(I-\widehat UV_1)|-\log|\det(I-\widehat UV_0)|\bigr|
 \le C_p(1+M)\sum_j|\zeta_j-1|^p
 \le C_p(1+g_0+\sigma)\sigma.
\]
Only the determinants at the endpoints of each step enter. If a zero is
crossed along a phase path, no branch of the logarithm is needed. This
proves \eqref{eq:unitary-log-comparison}.
\end{proof}

\subsection{Applying the log comparison to a remote positive source}\label{subsec:s30-applying-the-log-comparison-to-a-remote-}

\begin{figure}[H]
\centering
\begin{tikzpicture}[x=1cm,y=1cm]
 \foreach \k in {0,...,13}{
  \pgfmathsetmacro{\x}{\k*0.75}
  \ifnum\k<3
   \draw[fill=black!35] (\x,0) rectangle ++(.6,.6);
  \else\ifnum\k>10
   \draw[fill=black!35] (\x,0) rectangle ++(.6,.6);
  \else\ifnum\k<5
   \draw[fill=black!15] (\x,0) rectangle ++(.6,.6);
  \else\ifnum\k>8
   \draw[fill=black!15] (\x,0) rectangle ++(.6,.6);
  \else
   \draw[pattern=north east lines] (\x,0) rectangle ++(.6,.6);
  \fi\fi\fi\fi
 }
 \draw[decorate,decoration={brace,amplitude=4pt}] (0,.75)--(2.1,.75)
   node[midway,above=4pt,font=\small] {collar};
 \draw[decorate,decoration={brace,amplitude=4pt}] (8.25,.75)--(10.35,.75)
   node[midway,above=4pt,font=\small] {collar};
 \draw[decorate,decoration={brace,amplitude=4pt,mirror}] (2.25,-.15)--(3.6,-.15)
   node[midway,below=4pt,font=\small] {buffer $RW$};
 \draw[decorate,decoration={brace,amplitude=4pt,mirror}] (6.75,-.15)--(8.1,-.15)
   node[midway,below=4pt,font=\small] {buffer $RW$};
 \draw[decorate,decoration={brace,amplitude=4pt}] (3.75,.75)--(6.6,.75)
   node[midway,above=4pt,font=\small] {absorber $i\kappa P_C/W^2$ on $C$};
\end{tikzpicture}
\caption{The geometry used in
\eqref{eq:remote-scattering-source}--\eqref{eq:single-collar-log-error}.
Dark blocks form collars; the hatched blocks form the set $C$, which
carries the positive absorber; the light blocks are buffers of length
$RW$ that separate $C$ from every collar. The absorber changes the
scattering matrix seen by a collar only through waves that cross a buffer.
Theorem~\ref{thm:t22-remote-higher} shows that this effect decays like
$e^{-cR}$, and Lemma~\ref{lem:t22-unitary-log-comparison} turns the small
change of the scattering matrix into a small change of the closing
log-determinant.}
\label{fig:collar-buffer}
\end{figure}
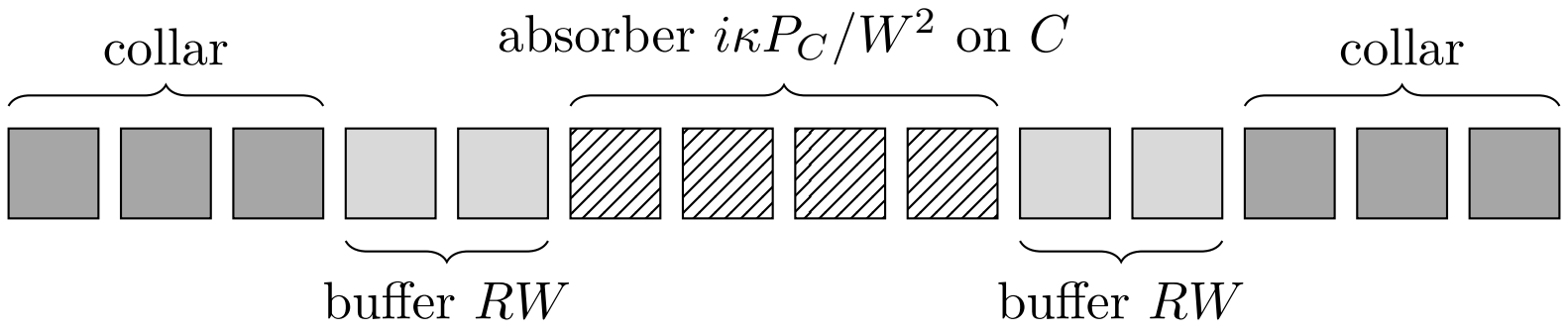

We   apply the log-determinant comparison to a positive absorber  far
from every collar. The   resulting error decreases exponentially with
the buffer length. We place the source in the core  shown in
Figure~\ref{fig:collar-buffer}, realize the absorbing response by a
unitary extension, and combine remote-clock decay with the fractional
boundary estimate, including the exceptional events.

Consider a disjoint union of at most $m$ actual Gaussian paths, with all
their endpoint matched leads included among the input ports. Without
absorption its complete scattering matrix $S_0$ is unitary. Let
\[
 Q_C=W^{-2}U_{\rm lead}^*P_CU_{\rm lead},
\]
where $C$ is a set of blocks at distance at least $R\ge2$ from every
endpoint, as in Figure~\ref{fig:collar-buffer}. We add the absorber
$i\kappa P_C/W^2$ with $\kappa\ge0$ and write $S_\kappa$ for the resulting
scattering matrix. The exact resolvent and Ward identities give
\begin{equation}
 \operatorname{Tr}|S_\kappa-S_0|^p
       \le(2\eta\kappa)^p\operatorname{Tr}Q_C^p,\qquad
 0\preceq I-S_\kappa^*S_\kappa\preceq4\eta\kappa Q_C.       \label{eq:remote-scattering-source}
\end{equation}
The second inequality uses the same contraction
$(I+i\kappa K_C/W^2)^{-1}$ as in \eqref{eq:signed-source-update}, before
either Gram product is taken. Thus the change of the scattering matrix is
controlled by fractional traces of $Q_C$, evaluated on the unabsorbed
reference path $S_0$. The finite absorber in this comparison is not a
persistent injection into the transmission estimates.

The matrix $S_\kappa$ is a contraction but not unitary, while
Lemma~\ref{lem:t22-unitary-log-comparison} is stated for unitary
matrices. We therefore pass to the unitary dilation
\[
 \mathcal V(S)=
 \begin{pmatrix}
 S&(I-SS^*)^{1/2}\\ (I-S^*S)^{1/2}&-S^*
 \end{pmatrix}.
\]
At $S=S_0$ the off-diagonal blocks vanish. The two defect matrices
$I-SS^*$ and $I-S^*S$ have the same nonzero spectrum. By subadditivity of
$\Tr|\cdot|^p$ for $0<p<1$ and by \eqref{eq:remote-scattering-source},
\begin{equation}
 \sigma:=\operatorname{Tr}|\mathcal V(S_\kappa)-\mathcal V(S_0)|^p
 \le C_p\{\kappa^p\operatorname{Tr}Q_C^p+
                  \kappa^{p/2}\operatorname{Tr}Q_C^{p/2}\}.       \label{eq:source-unitary-dilation}
\end{equation}
The off-diagonal blocks contribute the term with $p/2$, since
$\Tr(I-S_\kappa^*S_\kappa)^{p/2}\le(4\eta\kappa)^{p/2}\Tr Q_C^{p/2}$ by
\eqref{eq:remote-scattering-source}. In particular, we do not need a
fractional Schatten estimate for the difference of two arbitrary matrix
square roots.

We now choose a fixed exponent, for example $p=1/8$. Then $p/2$ and $p$ are
both below one, and Theorem~\ref{thm:t22-remote-higher}, in the window form
\eqref{eq:remote-clock-higher-window} of
Proposition~\ref{prop:t22-remote-higher-window}(c), gives, uniformly
for bounded $\kappa$,
\begin{equation}
 \limsup_W\mathbb E(\sigma+\sigma^2)
                       \le C_\kappa m^2 e^{-cR}.         \label{eq:dilation-fractional-moments}
\end{equation}
The factor $m^2$ accounts for the at most $m$ components and the square.

The comparison \eqref{eq:unitary-log-comparison} also contains the product
$g_0\sigma$. For a selected collar, the compressed conductance $g_0$
without absorption is zero if no other lead remains. Otherwise it is the
sum of at most two transmission conductances of paths of length at least
$R$. Theorem~\ref{thm:transmission-generator} gives $\E g_0^2\le C$,
independently of the total length. By Cauchy--Schwarz,
\begin{equation}\label{eq:g0-sigma}
 \E[g_0\sigma]\le(\E g_0^2)^{1/2}(\E\sigma^2)^{1/2}
 \le C\,(\E\sigma^2)^{1/2}.
\end{equation}
Combined with \eqref{eq:dilation-fractional-moments}, the right side of
\eqref{eq:unitary-log-comparison} is small in expectation, after
taking the square root of $m^2e^{-cR}$ in \eqref{eq:g0-sigma}. The bound is applied at every energy $e=E+x/\Lambda$ with $x$ in the
bounded support of $f$, to paths whose lengths $n_i$ satisfy
$n_i/W\to a_i\tau$, and to the indicator source of $C$, whose distances to
the two ends of a path may differ but are at least $R$. Since
$\Lambda=N\rho(E)$ and $N/W^2=n/W\to\tau$, these energies satisfy
$|e-E|\le KW^{-2}$ for large $W$, with $K=2\sup_{\operatorname{supp}f}|x|/(\tau\rho(E))$.
This is the window of Proposition~\ref{prop:t22-remote-higher-window}.

It remains to handle the events on which
Lemma~\ref{lem:t22-unitary-log-comparison} does not apply. The scalar
exceptional set has arbitrarily small polynomial probability by the actual
passive trace estimates. The event $\sigma>c_*W$ has probability tending to
zero by \eqref{eq:dilation-fractional-moments} and Markov's inequality. To
remove these events without creating an unbounded log error, we write each
source-centered real determinant as
\begin{equation}
 B_\kappa(e)=\log\left|
 \frac{\det(H_{\rm eff}-e+i\kappa P_C/W^2)}{\det(H_{\rm eff}-e)}\right|\ge0.                  \label{eq:positive-log-increment}
\end{equation}
For a matched path, eliminating its fixed auxiliary ports gives
\[
 H_{\rm eff}=H_{\rm path}+D+iY,\qquad D=D^*,\quad Y\succeq0,
\]
where $D$ and $Y$ depend only on endpoint bonds and are independent of the
physical Gaussian diagonals. Lemma~\ref{lem:t26-passive-log}, with
$a=\kappa/W^2$ and $Q=P_C$, gives the exact Cauchy integral
\[
 B_\kappa(e)=\int_\R B_{\kappa/W^2}(H_{\rm path}+D+tY;P_C)
                         \frac{dt}{\pi(1+t^2)}.
\]
Its proof supplies the scalar regularization and the order of integration.
Each integrand is nonnegative and bounded by its full-absorber version,
whose spectral kernel is explicit:
\begin{equation}\label{eq:full-absorber-kernel}
 \log\left|\frac{\det(H-e+i\kappa W^{-2})}{\det(H-e)}\right|
 =\sum_\lambda\frac{1}{2}\log\Bigl(1+\frac{\kappa^2}{W^4(\lambda-e)^2}\Bigr).
\end{equation}
Lemma~\ref{lem:t26-singular-log-moments} bounds all its fixed moments
 uniformly in $D+tY$ and  dimensions at most   $N$, since
$N\kappa/W^2=\kappa\tau_W\le K$. Lemma~\ref{lem:t26-passive-log}
preserves the bound for matched passive components, and the Hermitian
case applies to closed components. Explicitly, for any exceptional event
$\mathcal E_W$ and fixed $r>1$,
\[
 \E\one_{\mathcal E_W}
 |B_\kappa^{\rm closed}-B_\kappa^{\rm matched}|
 \le C_{r,K,c}\Prob(\mathcal E_W)^{1-1/r}=o(1).
\]
The uniformity in deterministic $e$ includes every bounded unfolded
window, so Tonelli applies before spectral integration.

We conclude that one actual three-block collar can be replaced by its two
matched leads in the source-centered real logarithm, at the cost
\begin{equation}
 \limsup_W\mathbb E|B_\kappa^{\rm closed}(e)
                           -B_\kappa^{\rm matched}(e)|
                         \le C_\kappa m^2e^{-c'R}.       \label{eq:single-collar-log-error}
\end{equation}
Here and below the physical source satisfies $\kappa\tau_W\le K$, so the
exceptional-set log moments are uniform in the total length; the
uniformity includes the bounded window of unfolded energies.

The reason why only the last factor matters is the exact factorization
\[
 \frac{\det H_\kappa^{\rm closed}}{\det H_\kappa^{\rm matched}}
 =\det(G-e)\,
   \frac{\det(K_c+i\eta)}{(2i\eta)^{2W}}\,
                    \det(I-U r_\kappa).
\]
Every factor before the last determinant is independent of the source
$\kappa$. It therefore cancels exactly in
\eqref{eq:positive-log-increment}, and neither its size nor its degree is
discarded. Only $\det(I-Ur_\kappa)$ remains, and this is the quantity
compared in Lemma~\ref{lem:t22-unitary-log-comparison}.

\subsection{Recovering spectral statistics from a positive logarithm}\label{subsec:s30-recovering-spectral-statistics-from-a-po}

  We recover ordinary spectral statistics from the positive
 log-determinants   that the closing argument controls.
Lemma~\ref{lem:t22-finite-log-recovery} bounds the error by the absorber
scale and the fraction of omitted blocks. We first compare with a full
absorber, then use monotonicity and projected spectral averaging to
pay for restricting it to the core.

Let $H$ be an actual self-adjoint original ring, or a direct sum of
original rings, and let $f\ge0$ be smooth with compact support. For a
deterministic union $C$ of core blocks we define
\begin{equation}
 \mathcal R_\delta^C(H;f)=\frac{1}{\pi\delta}\int f(x)
 \log\left|\frac{\det(H-E-x/\Lambda+i\delta P_C/\Lambda)}{\det(H-E-x/\Lambda)}\right|dx,
                 \qquad\delta>0.                       \label{eq:positive-log-statistic}
\end{equation}
This quantity is nonnegative, and it retains the complete determinant.

\begin{lemma}[Physical finite-log recovery]
\label{lem:t22-finite-log-recovery}
If $B=C^c$ is a union of blocks, then
\begin{equation}
 \mathbb E\left|\mathcal R_\delta^C(H;f)
        -\operatorname{Tr}f(\Lambda(H-E))\right|
 \le C_f\{\delta(1+|\log\delta|)+|B|/n\}.                 \label{eq:finite-log-recovery}
\end{equation}
The constants are uniform over the original ring and its independent
ring partition, with the same global unfolding.
\end{lemma}

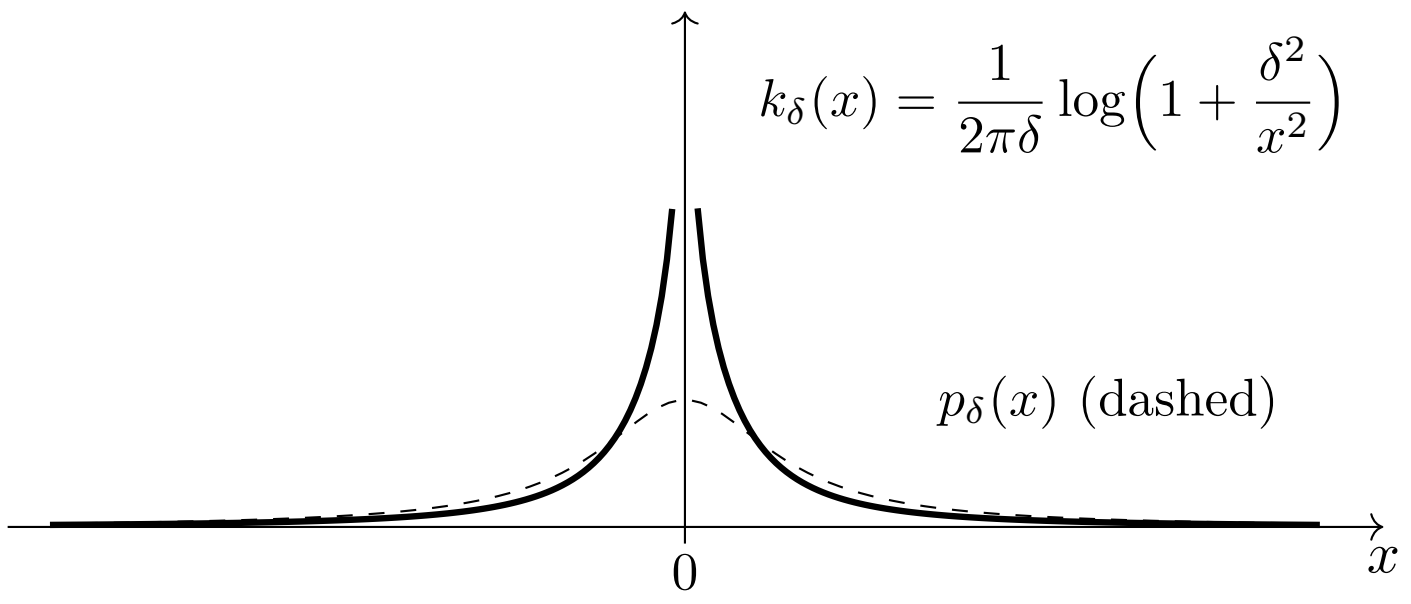
\begin{figure}[H]
\centering
\begin{tikzpicture}[x=1.4cm,y=1.1cm]
 \draw[->] (-3.2,0)--(3.3,0) node[below] {$x$};
 \draw[->] (0,-.1)--(0,3.1);
 \draw[very thick,domain=-3:-0.06,samples=120] plot (\x,{0.45*ln(1+0.25/(\x*\x))});
 \draw[very thick,domain=0.06:3,samples=120] plot (\x,{0.45*ln(1+0.25/(\x*\x))});
 \draw[dashed,domain=-3:3,samples=120] plot (\x,{1.2*0.5/(3.14159*(\x*\x+0.25))});
 \node[font=\small,anchor=west] at (.25,2.6) {$k_\delta(x)=\dfrac1{2\pi\delta}\log\Bigl(1+\dfrac{\delta^2}{x^2}\Bigr)$};
 \node[font=\small,anchor=west] at (1.1,.75) {$p_\delta(x)$ (dashed)};
 \node[font=\small,below] at (0,-.05) {$0$};
\end{tikzpicture}
\caption{The kernel $k_\delta$ of \eqref{eq:log-kernel} (solid) and the
Cauchy kernel $p_\delta$ (dashed), drawn for one value of $\delta$; the
curves are schematic. The kernel $k_\delta$ is the average of the Cauchy
kernels $p_t$ over $0<t<\delta$. It has a logarithmic singularity at $0$
but total mass one, so it is an approximate identity at scale $\delta$.
In the proof of Lemma~\ref{lem:t22-finite-log-recovery}, the positive
logarithm \eqref{eq:positive-log-statistic} with absorber on all blocks
equals the spectral statistic of the smoothed test $k_\delta*f$, and the
counting estimates bound the smoothing error.}
\label{fig:log-kernel}
\end{figure}

The first error term comes from smoothing the spectral statistic at scale
$\delta$. The second comes from the blocks in $B$, which carry no absorber;
there are $|B|$ such blocks among $n$. Figure~\ref{fig:log-kernel} shows the
smoothing kernel.

\begin{proof} When the absorber acts on every block, the positive
logarithm is exactly a sum over the eigenvalues of a smoothed version of the
test function, with the smoothing kernel of Figure~\ref{fig:log-kernel}, so the
counting estimates bound its distance to the spectral statistic (Step~1). In
general the absorber acts only on the core. The logarithm increases with the
absorber (Step~2), so the partial version lies below the full one, and their
difference is the effect of adding absorption on the remaining blocks, which
a projected spectral-averaging bound makes proportional to the fraction of
those blocks (Step~3).

\emph{Step 1: absorber on all blocks.} We first take $C$ to be the whole
space, so that $P_C=I$. Each eigenvalue $\lambda_\alpha$ of $H$
contributes one factor to both determinants in
\eqref{eq:positive-log-statistic}. With $u_\alpha=\Lambda(\lambda_\alpha-E)$
we obtain
\[
 \frac{1}{\pi\delta}\log\left|
 \frac{\det(H-E-x/\Lambda+i\delta/\Lambda)}{\det(H-E-x/\Lambda)}\right|
 =\sum_\alpha\frac{1}{2\pi\delta}
   \log\Bigl(1+\frac{\delta^2}{(u_\alpha-x)^2}\Bigr).
\]
Integrating against $f$, this gives
$\mathcal R_\delta^I(H;f)=\sum_\alpha(k_\delta*f)(u_\alpha)$, with the kernel
\begin{equation}\label{eq:log-kernel}
 k_\delta(x)=\frac{1}{2\pi\delta}\log\Bigl(1+\frac{\delta^2}{x^2}\Bigr)
 =\frac{1}{\delta}\int_0^\delta p_t(x)\,dt,
 \qquad p_t(x)=\frac{t}{\pi(x^2+t^2)}.
\end{equation}
The second equality follows from
$\int_0^\delta t\,dt/(x^2+t^2)=\frac{1}{2}\log(1+\delta^2/x^2)$. Since each
$p_t$ has mass one, so does $k_\delta$; see Figure~\ref{fig:log-kernel}.
Use Theorem~\ref{thm:external-count}, conditional on the bonds,
for either the ring or the direct sum of rings, and the Tonelli
estimate in the proof of Lemma~\ref{lem:smoothing}:
\[
 \int p_t(s)\min\{a|s|,b\}\,ds
 \le\frac{at}{\pi}\log(1+t^{-2})+\frac{2bt}{\pi}.
\]
Averaging $0<t<\delta$ gives $C_f\delta(1+|\log\delta|)$ for the
expected difference between $\sum_\alpha(k_\delta*f)(u_\alpha)$ and
$\sum_\alpha f(u_\alpha)$. The rate comes from that proof estimate.
This is the first error in \eqref{eq:finite-log-recovery}.

\emph{Step 2: monotonicity in the absorber.} Let $A$ be Hermitian, $P\succeq0$,
and $G_s=(A+isP)^{-1}$. Then
\begin{equation}\label{eq:log-monotone}
 \frac{d}{ds}\log|\det(A+isP)|=\Re\Tr(iPG_s)=-\Im\Tr(PG_s)\ \ge0,
\end{equation}
because $\Im G_s=-sG_s^*PG_s\preceq0$. Increasing a positive absorber
therefore increases the real logarithm of the determinant. Since
$P_C\preceq I$, we obtain
\[
 0\le\mathcal R_\delta^C(H;f)\le\mathcal R_\delta^I(H;f).
\]

\emph{Step 3: the missing blocks.} The difference
$\mathcal R_\delta^I-\mathcal R_\delta^C$ is obtained by increasing the
absorber on $B$ from $0$ to $\delta/\Lambda$, with the absorber on $C$
already present. By \eqref{eq:log-monotone},
\[
 \mathcal R_\delta^I-\mathcal R_\delta^C
 =\frac{1}{\pi\delta}\int f(x)\int_0^{\delta/\Lambda}
   \bigl[-\Im\Tr(P_BG_s(x))\bigr]\,ds\,dx,
\]
where $G_s(x)$ is the resolvent with absorber $\delta/\Lambda$ on $C$ and
$s$ on $B$. Lemma~\ref{lem:t26-projected-wegner}, conditioned on every variable
except the retained diagonal block, gives
$\E[-\Im\Tr P_BG_s(x)]\le CW|B|$ at every intermediate value.
Its hypothesis holds because these absorbers are deterministic; the
original or rewired bonds are independent of all diagonal GUE blocks.
One may first insert a scalar regulator, use the nonnegative integrands
and Tonelli, and remove that regulator as in
Lemma~\ref{lem:t26-passive-log}. Hence
\[
 \E\bigl[\mathcal R_\delta^I-\mathcal R_\delta^C\bigr]
 \le\frac{\norm f_1}{\pi\delta}\cdot\frac{\delta}{\Lambda}\cdot CW|B|
 =\frac{C\norm f_1}{\pi\rho(E)}\cdot\frac{|B|}n,
\]
where we used $\Lambda=nW\rho(E)$. This is the second error in
\eqref{eq:finite-log-recovery}. The two steps together prove the lemma. We
did not use that the partial absorber is diagonal in the eigenbasis of
$H$, which it is not.
\end{proof}

\subsection{Same-volume rewiring and the Poisson endpoint}\label{subsec:s30-same-volume-rewiring-and-the-poisson-end}

\begin{figure}[ht]
\centering
\begin{tikzpicture}[x=1cm,y=1cm]
 \begin{scope}
  \draw[thick] (0,0) circle (1.6);
  \foreach \a in {90,210,330}{
   \fill[black!40] (\a-7:1.6) arc (\a-7:\a+7:1.6) -- (\a+7:1.35) arc (\a+7:\a-7:1.35) -- cycle;
  }
  \node[font=\small] at (0,0) {original ring};
  \node[font=\small] at (0,-2.05) {cyclic pairing};
 \end{scope}
 \begin{scope}[shift={(6,0)}]
  \foreach \a/\b in {90/210,210/330,330/450}{
   \pgfmathsetmacro{\m}{(\a+\b)/2}
   \draw[thick] (\a+9:1.6) arc (\a+9:\b-9:1.6);
   \draw[thick] (\a+9:1.6) .. controls (\m:.55) .. (\b-9:1.6);
  }
  \foreach \a in {90,210,330}{
   \fill[black!40] (\a-7:1.6) arc (\a-7:\a+7:1.6) -- (\a+7:1.35) arc (\a+7:\a-7:1.35) -- cycle;
  }
  \node[font=\small] at (0,-2.05) {independent-ring pairing};
 \end{scope}
 \draw[->,thick] (2.2,0)--(3.8,0) node[midway,above,font=\small] {rewire};
\end{tikzpicture}
\caption{The same-volume rewiring for $m=3$. Dark arcs are the collars.
Left: with the cyclic pairing of the bonds at the collars, the blocks form
the original ring. Right: pairing both ends of each complementary path
through one collar closes each path into a ring of its own. Every block
keeps its diagonal entries and every bond keeps its law, so the two
models have the same total dimension and are coupled on one probability
space. The proof of Theorem~\ref{thm:t22-Poisson-endpoint} compares the
two models one collar at a time, using
\eqref{eq:single-collar-log-error} and the buffers of
Figure~\ref{fig:collar-buffer}.}
\label{fig:rewiring}
\end{figure}
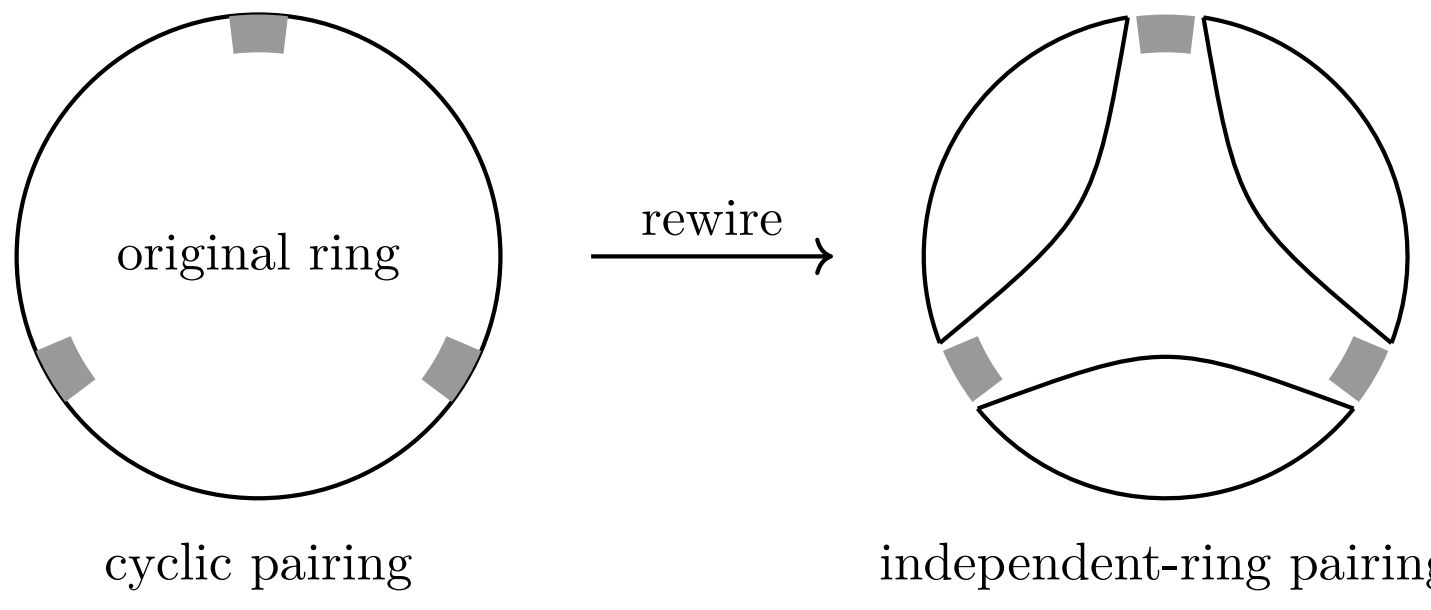

We   complete the Poisson endpoint by comparing the  ring with
independent   smaller rings of the same total dimension. The sequential
collar estimate controls the rewiring error without changing the
Gaussian prior. We couple the two geometries in
 Figure~\ref{fig:rewiring}  , remove collars one at a time, recover the
spectral statistics, and use independence and the counting bounds to
compute the limiting Poisson Laplace functional.

We choose $m$ disjoint three-block collars around the ring. Their
complement consists of $m$ independent actual Gaussian paths. The two
incident Gaussian bonds of each collar block can be paired in two ways.
The original cyclic pairing gives the original ring. Pairing both ends of
each complementary path through one collar gives $m$ independent original
rings, as on the right of Figure~\ref{fig:rewiring}. All diagonal blocks
and bonds keep their original distributions, and each physical block
belongs to exactly one component in either model. So this is a coupling
of two matrices of the same dimension; no buffer or spectrum is
duplicated.

Let $C$ consist of all blocks except the collars and a buffer of $RW$
blocks at each end of every collar. Assume that each complementary path
has at least $4RW$ physical blocks (integer roundings do not affect the
estimates). The following lemma supplies the stage-by-stage law needed
for the telescope.

In the proof of the Poisson endpoint, we
compare the original ring with independent smaller rings by removing
the $m$ collars one at a time and then putting them back in the other
pairing. Each single removal is controlled by the one-collar estimates,
which need three inputs. The fresh collar must have its original
Gaussian law, independently of the rest. The old response it sees must
be regular. And the determinant must factor so that only the closing
determinant depends on the source. The lemma checks these three inputs
at every stage, under the original prior. It also checks that the two
pairings end at the same fully matched model, so that the two sequences
of removals can be compared.

\begin{lemma}[Sequential collar removal with the actual Gaussian prior]
\label{lem:t26-sequential-collars}
Fix a deterministic order of the $m$ collars. In either pairing, remove
one collar at a time by deleting its middle block and replacing its two
outer physical diagonals by the matched scalar diagonal $zI_W$. Keep its
two exterior incident bonds. At each step:
\begin{enumerate}[label=(\roman*)]
\item The selected fresh collar consists of its three GUE diagonals and
its two internal Ginibre bonds, independent of the entire old network.
Every old component incident to it is an actual Gaussian path, and
previously removed outer coordinates are matched auxiliary ports.
\item The compressed zero-source response $r_0$ satisfies the collar
scalar certificate, meaning $d(X)\le d_*/2$ in the Cayley representation of Theorem~\ref{thm:t22-collar-kernel}, except on probability $O_{m,D}(W^{-D})$ for every
fixed $D$, at each deterministic compact-bulk energy. Its defect
$g_0=\Tr(I-r_0^*r_0)$ has $\E g_0^2\le C$, uniformly in total length.
\item The determinant factorization has exactly the source-independent
prefactors in \eqref{eq:t26-stage-factorization} below. The fully matched
final models for the two pairings agree up to a deterministic
permutation of auxiliary coordinates.
\end{enumerate}
The width limit is taken for fixed $m,R,\tau$. Statements about events of
the old data are under the original prior, without conditioning on the
success of earlier stages.
\end{lemma}

\begin{figure}[ht]
\centering
\begin{tikzpicture}[x=1cm,y=1cm]
 \begin{scope}
  \node[font=\small,anchor=west] at (-.3,1.25) {before removal};
  \foreach \x/\f in {0/black!35,.85/black!35,1.7/black!12,2.55/red!35,3.4/black!12,4.25/black!35,5.1/black!35}
    \draw[thick,fill=\f] (\x,0) rectangle ++(.7,.7);
  \draw[decorate,decoration={brace,amplitude=4pt,mirror}] (1.7,-.15)--(4.1,-.15)
    node[midway,below=4pt,font=\small] {collar: $A$, $G$, $A$};
  \node[font=\small] at (2.9,.35) {$G$};
 \end{scope}
 \begin{scope}[yshift=-2.6cm]
  \node[font=\small,anchor=west] at (-.3,1.25) {after removal};
  \foreach \x/\f in {0/black!35,.85/black!35,4.25/black!35,5.1/black!35}
    \draw[thick,fill=\f] (\x,0) rectangle ++(.7,.7);
  \draw[very thick,dashed] (1.7,0) rectangle ++(.7,.7);
  \draw[very thick,dashed] (3.4,0) rectangle ++(.7,.7);
  \node[font=\small] at (2.05,.35) {$zI$};
  \node[font=\small] at (3.75,.35) {$zI$};
  \node[font=\small,align=center] at (2.9,-.55) {middle block deleted,\\ outer diagonals matched};
 \end{scope}
 \node[font=\small,anchor=west,align=left] at (6.2,.35)
   {old paths (dark) keep\\ their incident bonds};
 \node[font=\small,anchor=west,align=left] at (6.2,-2.25)
   {two matched ports $zI_W$\\ (one-port reflections)};
\end{tikzpicture}
\caption{One stage of the sequential removal in
Lemma~\ref{lem:t26-sequential-collars}. Before the stage (top), the
selected collar consists of two outer physical blocks $A$ (light) and a
middle block $G$ (red), with two internal bonds, between old Gaussian
paths (dark). The stage deletes the middle block and replaces the two
outer diagonals by the matched scalar diagonal $zI_W$ (dashed), keeping
the two exterior bonds to the old paths. The old paths are therefore
terminated by matched ports, which is why their response $r_0$ is
regular (item~(ii)). The determinant ratio between the two pictures is
the factorization \eqref{eq:t26-stage-factorization}, in which only the
closing determinant $\det(I-Ur_\kappa)$ depends on the source.}
\label{fig:sc-stage}
\end{figure}
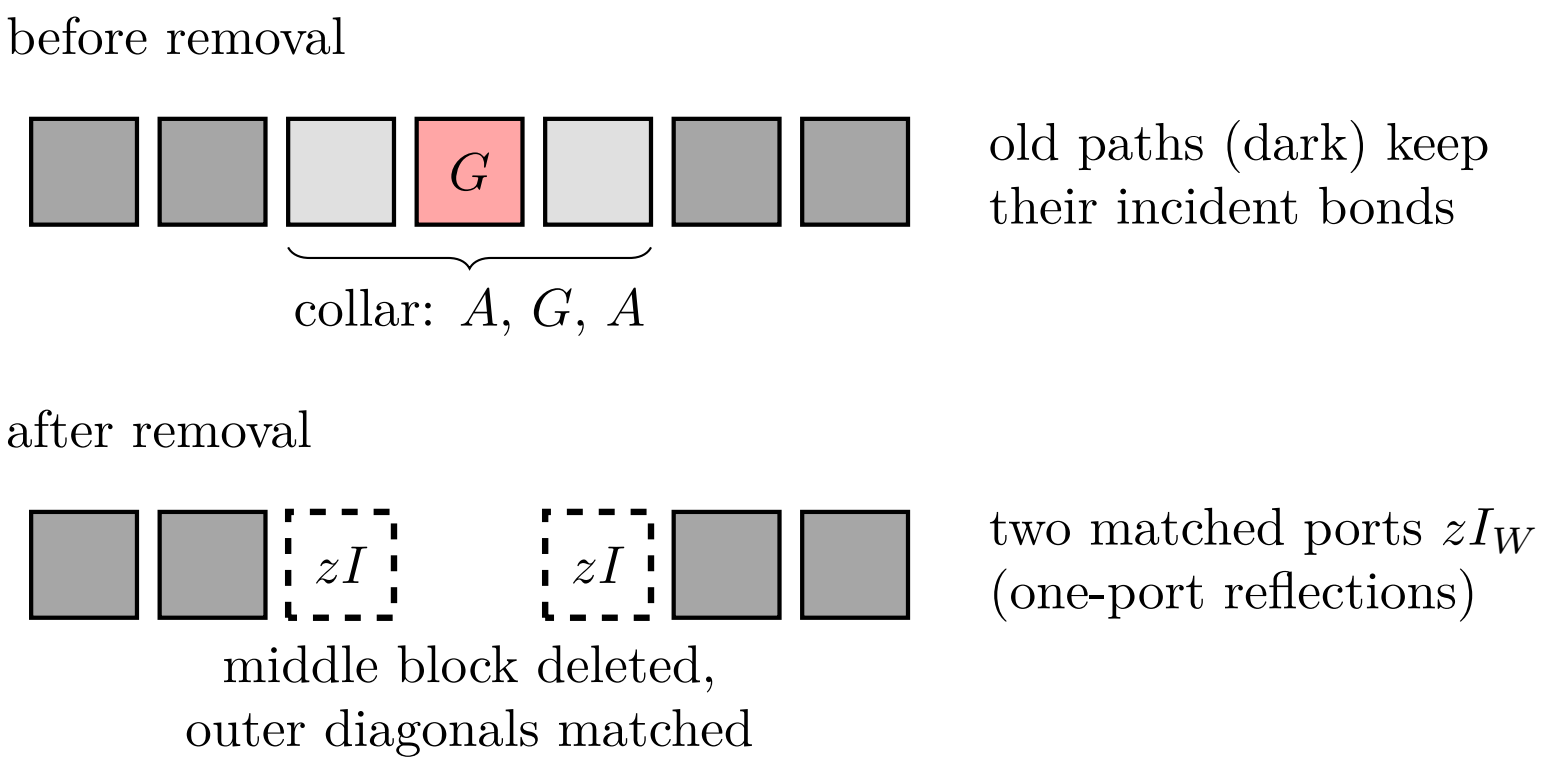

\begin{proof} We remove the collars one at a time, as in
Figure~\ref{fig:sc-stage}, and check three things at each stage. First, the
Gaussian variables split into independent groups, so the selected collar keeps
its original law, independently of everything else (Steps~1--2). Second, the
response of the old network seen by the collar is regular, because it is built
from actual Gaussian paths that end at matched ports (Steps~3--4). Third, the
ratio of the determinants before and after the removal factors exactly, with
only one factor depending on the source (Step~5). Step~6 shows that the two
ways of pairing the collars lead to the same final model.

\emph{Step 1: the independent variables.} We partition the independent
Gaussian variables into three groups: the disjoint collar sets, each with
three diagonals and two internal bonds, the sets of the complementary
paths, and the exterior incident bonds. We assign each incident bond to
its endpoint on a path. The cyclic pairing and the box pairing differ
only in the collar to which the other endpoint of each such bond is
attached. In particular, each final box uses a disjoint subset of the
independent variables, and it has exactly the law of an original ring,
including its two incident bonds and its block parity.

\emph{Step 2: item (i).} For the cyclic pairing, once the first collar is
removed, the remaining physical graph is a path. Each later removal
splits paths into paths. The collars not yet removed are ordinary
physical blocks inside these paths. Conditionally on every variable
outside the selected collar set, its five fresh blocks keep their
independent Gaussian law. No earlier comparison changes this probability
measure. Indeed, the telescoping sum is a sum of random variables on the
original joint space, and each summand is averaged over its own selected
collar. For the box pairing, removing a collar opens its box into the
complementary path, and the other boxes are disconnected and cancel from
the selected determinant ratio. This proves (i).

\emph{Step 3: the shape of the old response.} We now study the response
$r_0$ that the selected collar sees. We remove the selected collar when we
form its old response, and we temporarily expose all the remaining matched
ports. At zero source, the complete scattering matrix of the incident
paths is unitary. There are two cases, according to how the two selected
ports are connected. If they are connected to each other, with no other
port on their component, then $r_0$ is their full unitary two-port
response, so $g_0=0$. Otherwise, each selected port is the near end of one
of two disjoint paths, each ending at a previously matched port. Then
$r_0$ is the direct sum of their two one-port reflections, and $g_0$ is the
sum of their transmission conductances. Each such path contains at least
one whole complementary path, so its length is at least $4RW$.

\emph{Step 4: item (ii).} In the first case, $(I-r_0)/(2i\eta)$ is the full
two-port resolvent of an actual Gaussian path. In the second case, its
two diagonal blocks are exactly the near diagonal blocks of the two
corresponding two-port resolvents. In both cases,
Theorem~\ref{thm:matched-traces} and
Lemma~\ref{lem:t26-passive-initial-scalar} give
\[
 \left\|\max_{i=1,2}
 \left|\frac{1}W\Tr\left[\frac{I-r_0}{2i\eta}\right]_{ii}-a\right|
 \right\|_p\le C_{p,\tau}/W.
\]
A sufficiently high fixed moment, with Markov's inequality, pays for the
prescribed power of the scalar failure, and a union bound pays for the
finitely many stages. These are estimates under the old law itself, not
after selection by earlier good events. The conductance bound follows
from \eqref{eq:conductance-moments}, applied at lengths at least $4RW$.
The matrix $r_0$ is a contraction, and it has the passive reaction
required by Theorem~\ref{thm:t22-collar-kernel}. The inverse defining
that reaction exists almost surely, by the Schur formulas for a finite
path. This proves (ii).

\emph{Step 5: the determinant factorization.} For (iii), let $M_\kappa$ be
the old interior matrix. It includes all the previously matched diagonals
and the source $i\kappa P_C/W^2$. Let $C_o$ contain its bonds to the
selected outer coordinates. At a fixed real energy we put
$X_\kappa=C_o^*M_\kappa^{-1}C_o$, the Schur complement term of the old
network at the outer coordinates. Eliminating first the old interior and
then the middle block $G$, Schur inversion gives
\[
 \det H_\kappa^{\rm before}
 =\det M_\kappa\det(G-e)
   \det(A-e-X_\kappa-B(G-e)^{-1}B^*),
\]
while, after the outer diagonals are replaced by $zI_W$,
\[
 \det H_\kappa^{\rm after}=\det M_\kappa\det(z-X_\kappa).
\]
We put
\[
 K_c=A-e-\xi-B(G-e)^{-1}B^*,\qquad
 U=I-2i\eta(K_c+i\eta)^{-1},\qquad
 r_\kappa=I-2i\eta(z-X_\kappa)^{-1},
\]
with $z=\xi+i\eta$. Multiplying out, we obtain the exact identity
\[
 I-Ur_\kappa=2i\eta(K_c+i\eta)^{-1}
             (K_c+\xi-X_\kappa)(z-X_\kappa)^{-1}.
\]
Since $K_c+\xi-X_\kappa=A-e-X_\kappa-B(G-e)^{-1}B^*$, dividing the two
determinants proves
\begin{equation}\label{eq:t26-stage-factorization}
 \frac{\det H_\kappa^{\rm before}}{\det H_\kappa^{\rm after}}
 =\det(G-e)\frac{\det(K_c+i\eta)}{(2i\eta)^{2W}}
                   \det(I-Ur_\kappa).
\end{equation}
The source vanishes on the collar and on its buffers. So the two
prefactors $\det(G-e)$ and $\det(K_c+i\eta)/(2i\eta)^{2W}$ do not depend on
$\kappa$, and they cancel when we subtract the logarithm at zero source.
All the real inverses displayed exist almost surely at a fixed energy.
Alternatively, we may prove the identity with positive regulators and
then take the boundary values of the determinants.

\emph{Step 6: the two pairings end at the same model.} When every collar
has been removed, each complementary path keeps its original two incident
bonds and its two scalar matched endpoints. These endpoints have the same
diagonal $zI_W$ in both pairings. Renaming them by a deterministic
permutation identifies the fully matched matrices in the two pairings,
including their source matrices. No physical variable is duplicated. The
factors of the middle blocks removed at intermediate steps did not depend
on the source, and they have already cancelled exactly. This proves (iii).
\end{proof}

At each stage all currently exposed ports are at distance at least $RW$
from $C$. The old scattering comparison therefore uses at most $m$
actual paths with the remote-source estimate
\eqref{eq:dilation-fractional-moments}. Its energies, lengths and
sources are those of Proposition~\ref{prop:t22-remote-higher-window}(c),
as noted after \eqref{eq:dilation-fractional-moments}. Disconnected closed components
cancel from that stage and are omitted. The preceding lemma supplies
precisely the independence, scalar certificate and bounded $g_0$ needed
in \eqref{eq:single-collar-log-error}. The exceptional sets are removed
by Lemma~\ref{lem:t26-passive-log}; their dependence on other stages
is irrelevant to the sum of expectations.
For the physical source $\kappa=\delta/(\tau_W\rho(E))$,
 $\kappa\tau_W=\delta/\rho(E)$, so the log-moment constants are uniform
for $\tau\ge1$, and $\kappa$ is bounded there. Telescope each pairing
to the common fully matched model, add the at most $2m$ stage errors,
and integrate against $f$. We obtain
\begin{equation}
 \limsup_W\mathbb E|
 \mathcal R_\delta^C(H^{\rm ring};f)
       -\mathcal R_\delta^C(H^{\rm boxes};f)|
              \le C_{f,\delta}m^3e^{-cR}.                \label{eq:same-volume-log-error}
\end{equation}
The constants are uniform for $\tau\ge1$. They come from
\eqref{eq:source-unitary-dilation}--\eqref{eq:dilation-fractional-moments},
the fixed collar estimates, and the global microscopic log moments; none
of them grows exponentially in $\tau$.

\paragraph{The Laplace functional of independent rings.} Let the
independent rings have $n_i$ blocks, and put $a_i=n_i/n$, so that
$\sum_ia_i=1$. Write
\[
 Y_i=\Tr f(\Lambda(H_i-E)),\qquad h=1-e^{-f}\in[0,1).
\]
We first compare $1-e^{-Y_i}$ with the linear statistic of $h$. For a
configuration with points $x_1,\ldots,x_k$ in the support of $f$, we have
$e^{-Y_i}=\prod_j(1-h(x_j))$. Expanding the product and using the
Bonferroni inequalities,
\begin{equation}\label{eq:bonferroni-ring}
 \Bigl|1-\prod_{j=1}^k(1-h(x_j))-\sum_{j=1}^kh(x_j)\Bigr|
 \le\sum_{j<l}h(x_j)h(x_l)\le\frac{(k)_2}{2}.
\end{equation}
The factorial count bound \eqref{eq:factorial-count}, applied to the ring
$H_i$ of dimension $n_iW$ and an interval of length $|\operatorname{supp}f|/\Lambda$,
gives $\E(k)_2\le(Kn_iW|\operatorname{supp}f|/\Lambda)^2\le C_fa_i^2$.
Taking expectations in \eqref{eq:bonferroni-ring},
\begin{equation}\label{eq:one-ring-linear}
 \bigl|\E(1-e^{-Y_i})-\E\Tr h(\Lambda(H_i-E))\bigr|\le C_fa_i^2 .
\end{equation}
For each fixed partition, Theorem~\ref{thm:critical-density} shows that the
critical mean of the linear statistic is $a_i\int h$, since the ring $H_i$
carries the fraction $a_i$ of the global unfolding. By independence of the
rings,
\[
 \E e^{-\sum_iY_i}=\prod_i\E e^{-Y_i}
 =\prod_i\Bigl(1-a_i\int h+\epsilon_i\Bigr),
 \qquad \limsup_W|\epsilon_i|\le C_fa_i^2 .
\]
We compare this product with
$\exp\{-\int h\}=\prod_ie^{-a_i\int h}$. All factors have modulus at most
one, so the difference of the two products is at most the sum of the
differences of their factors. Since $|1-x-e^{-x}|\le x^2/2$ for $x\ge0$,
\begin{equation}\label{eq:independent-ring-laplace}
 \limsup_W\left|\E e^{-\sum_iY_i}
          -\exp\Bigl\{-\int(1-e^{-f})\Bigr\}\right|
       \le C_f\sum_i a_i^2.
\end{equation}
Each partition is fixed before $W\to\infty$; we do not use an estimate
that is uniform in a box length growing simultaneously with $W$.

\begin{proof}[Proof of Theorem~\ref{thm:t22-Poisson-endpoint}]
The proof assembles the estimates of this section.
All statistics involved are positive, so two Laplace functionals differ by at
most the expected difference of the statistics. This lets us pass from the
original ring to the independent rings of Figure~\ref{fig:rewiring} in three
moves: recover the positive logarithm, rewire the collars, and recover the
spectral statistic again (Step~1). Step~2 collects the errors of these moves
and of the Laplace functional of the independent rings, Step~3 chooses the
number of rings and the buffer length so that every error vanishes, and
Step~4 takes the limits in the right order.

\emph{Step 1: reduction to positive statistics.} For $a,b\ge0$ we have
$|e^{-a}-e^{-b}|\le|a-b|$. Both the positive logarithm
\eqref{eq:positive-log-statistic} and the spectral statistic are
nonnegative, so their Laplace functionals differ by at most the expected
absolute difference of the statistics. We pass from the original ring to
the independent rings in three moves:
\[
 \Tr f(\Lambda(H^{\rm ring}-E))
 \ \xrightarrow{\ \text{Lemma~\ref{lem:t22-finite-log-recovery}}\ }\
 \mathcal R_\delta^C(H^{\rm ring};f)
 \ \xrightarrow{\ \eqref{eq:same-volume-log-error}\ }\
 \mathcal R_\delta^C(H^{\rm boxes};f)
 \ \xrightarrow{\ \text{Lemma~\ref{lem:t22-finite-log-recovery}}\ }\
 \sum_iY_i.
\]

\emph{Step 2: collecting the errors.} The set $B=C^c$ consists of the $m$
collars and their buffers, so $|B|\le m(3+2RW)$ and
$|B|/n\le CmR/\tau$ for large $W$. Combining
\eqref{eq:finite-log-recovery}, \eqref{eq:same-volume-log-error} and
\eqref{eq:independent-ring-laplace}, we obtain
\begin{equation}\label{eq:poisson-error-sum}
 \limsup_W\Bigl|\E e^{-\Tr f(\Lambda(H-E))}
   -e^{-\int(1-e^{-f})}\Bigr|
 \le C_f\delta(1+|\log\delta|)
 +C_f\frac{mR}{\tau}+C_{f,\delta}m^3e^{-cR}
 +C_f\sum_i a_i^2 .
\end{equation}

\emph{Step 3: the choice of parameters.} We take $m$ of order $\sqrt\tau$,
nearly equal box lengths, and $R=A\log(2+\tau)$ with $A>2/c$. Then
\[
 \frac{mR}{\tau}\asymp\frac{\log\tau}{\sqrt\tau},\qquad
 m^3e^{-cR}\asymp\tau^{3/2}(2+\tau)^{-cA},\qquad
 \sum_ia_i^2\asymp\frac{1}m\asymp\frac{1}{\sqrt\tau},
\]
and all three terms tend to zero as $\tau\to\infty$, since $cA>2>3/2$. The
complementary paths have length about $\tau/m\asymp\sqrt\tau$, which is at
least $4R$ for large $\tau$, as required in
\eqref{eq:same-volume-log-error}.

\emph{Step 4: the order of limits.} We first take the critical limit
$W\to\infty$ for the fixed partition. Then we let $\tau\to\infty$ at fixed
$\delta$; by Step~3 the last three terms of \eqref{eq:poisson-error-sum}
vanish. Finally we let $\delta\downarrow0$, which removes the first term.
This proves \eqref{eq:poisson-laplace-limit}. The count bounds
\eqref{eq:factorial-count} extend this Laplace characterization to the
simple locally finite point process and its counts in disjoint intervals.
All volume and unfolding constants are unchanged.
\end{proof}

The comparison uses only finite positive absorbers, the actual weighted
closing distribution, and the decay of remote clocks. Together with
Theorem~\ref{thm:t22-sine-endpoint} and the proof of positive variance,
it completes the proof of Theorem~\ref{thm:endpoints}.

\appendix

\section{AxiomProver Lean certificate and AI declaration}\label{app:lean}

\paragraph{Lean certificate.}
Theorems~\ref{thm:que}, \ref{thm:critical} and~\ref{thm:endpoints} have been formalized
in Lean~4 using AxiomProver. Each formal statement contains the model, the normalizations
and the quantifiers of the corresponding theorem. The proofs assume no unproved
hypothesis, and the only axioms used are the three standard axioms of Lean's foundations
(\texttt{propext}, \texttt{Classical.choice}, \texttt{Quot.sound}). 

The development formalizes the arguments of this paper in about $1.0$ million lines of
Lean. %
It also formalizes,
in about $0.2$ million further lines, the results from other papers that the proofs use,
so that none of them is assumed:
\begin{itemize}
  \item the eigenvalue counting estimate of Peled, Schenker, Shamis and
        Sodin~\cite{PSSS} (Theorem~\ref{thm:external-count}), together with the Gaussian
        regularization estimates of Aizenman, Peled, Schenker, Shamis and Sodin~\cite{APSSS} used in
        its proof;
  \item the local laws and eigenstate thermalization of Cipolloni, Erd\H{o}s and
        Henheik~\cite{CEH} (Theorems~\ref{thm:external-local} and~\ref{thm:external-eth}),
        for the GUE-deformed model used here;
  \item the band local law of Erd\H{o}s and Riabov~\cite{ER}
        (Theorem~\ref{thm:external-band-local}), in the short-ring setting used here;
  \item Johansson's double-contour formula~\cite{JohanssonKernel}, used in
        Theorem~\ref{thm:external-sine};
  \item Pr\'ekopa's theorem on marginals of log-concave functions~\cite{Prekopa}.
\end{itemize}
Classical tools used along the way, such as the Weyl integration formula, the Karlin–McGregor formula and the noncollision of Dyson Brownian motion, are likewise either proved in the formalization or replaced by direct arguments.
Our treatment of conditional eigenvector moments follows the methods of
Benigni~\cite{Benigni} and of Bourgade, Huang and Yau~\cite{BHY}, recalled in
Theorems~\ref{thm:external-emf} and~\ref{thm:external-finite-speed}. The
proofs use complex analogues of these results, established in this paper
along the same lines and included in the formalization. The formal statements, proof sources and instructions for
checking the proofs are available at:
\begin{center}
\url{https://github.com/AxiomMath/BandMatrixQUE}.
\end{center}

\paragraph{Statement on the use of AI.}
The results of this paper were obtained through a human--AI collaboration using
AxiomProver and publicly available AI models. AxiomProver was used to
autoformalize in Lean~4 the content of this paper, including the cited results listed
above. The core ideas behind the overall architecture of the results were provided by the authors. The writing of the proofs in this manuscript were drafted with LLMs and
carefully edited by the authors, with the most extensive human writing and
revision in the introduction. The authors take full responsibility for the mathematical statements, exposition and conclusions presented here. The correctness of Theorems~\ref{thm:que}, \ref{thm:critical} and~\ref{thm:endpoints}
is certified by the Lean formalization described above.

\section*{Acknowledgement}
We are grateful to Carina Hong, Ken Ono and Shubho Sengupta for all their support and encouragement throughout this project, to Jesse Thorner for many useful discussions, as well as to the whole team of mathematicians and engineers at Axiom Math, whose work on AxiomProver made the Lean formalization of this paper possible. GR's research was supported by Axiom Math. PG's research was partially supported by NSF career award DMS 2542899. PG thanks Kevin Yang for explaining previous works on random band matrices.

\newcommand{\etalchar}[1]{$^{#1}$}

\end{document}